%% file: m5-main.tex
\documentclass{gtbook}

\usepackage{amsmath, amssymb, amscd}  

\newtheorem{theorem}{Theorem}[chapter]
\newtheorem{lemma}[theorem]{Lemma}
\newtheorem{cor}{Corollary}[theorem]
\newtheorem*{add}{Addendum}
\newtheorem*{thm}{Theorem}
\newtheorem*{lem}{Lemma}

\renewcommand{\geq}{\geqslant}
\renewcommand{\leq}{\leqslant}

\theoremstyle{definition}

\newtheoremstyle{numbered}{14pt plus6.3pt minus6.3pt}{7.4pt plus3pt minus3pt}%
{\sl}{}{\bf}{}{0.75em}{\thmname{#1}\thmnumber{ #2}\thmnote{\bf\ #3}}

\theoremstyle{numbered}
\newtheorem*{thmm}{Theorem}

\newcommand{\andeq}{\text{and}\qquad\qquad}
\newcommand{\eqand}{\qquad\qquad\phantom{\text{and}}}

\begin{document}

\frontmatter
\thispagestyle{empty}
\title{Four-manifolds, geometries and knots}

\author{J.A.Hillman}
\address{The University of Sydney}
\email{jonathan.hillman@sydney.edu.au}

\makebooktitle

\thispagestyle{empty}

\bgroup
\gtm\small\nl
ISSN 1464-8997 (on-line) 1464-8989 (printed)\nl
Volume 5 (2002)\nl
Four-manifolds, geometries and knots,
by J.A.Hillman\nl
Published 9 December 2002\nl
Revised 29 March 2007 and 14 November 2022\nl
\copyright\ \gtp\nl
All commercial rights reserved\nl
\medskip

\vfill \hrule \smallskip \small \parskip =0.07truein

\gt \ is a fully refereed international journal dealing with all
aspects of geometry and topology. It is intended to provide free
electronic dissemination of high quality research.  The \gtm\ series
is intended to provide a similar forum for conference proceedings and
research monographs.  Printed copy is available. Full details,
including pricing, are on the journal web site (see below).

Submission must be made electronically. For full instructions visit
the journal web site at:\nl {\tt http://www.maths.warwick.ac.uk/gt/}\nl or
choose your nearest mirror site from the EMIS network:\quad {\tt
http://www.emis.de/}\nl or use anonymous ftp to:\quad {\tt
ftp.maths.warwick.ac.uk}

{\it Academic Editorial Board:}\nl Joan Birman, Martin Bridson, Gunnar
Carlsson, Ralph Cohen, Simon Donaldson, Bill Dwyer, Yasha Eliashberg,
Steve Ferry, Ron Fintushel, Mike Freedman, David Gabai, Tom
Goodwillie, Cameron Gordon, Vaughan Jones, Rob Kirby, Frances Kirwan,
Dieter Kotschick, Peter Kronheimer, Wolfgang Metzler, Haynes Miller,
John Morgan, Shigeyuki Morita, Tom Mrowka, Walter Neumann,
Jean-Pierrre Otal, Ron Stern, Gang Tian

{\it Managing Editors:}\nl Colin Rourke, Brian Sanderson

Geometry and Topology Publications
\nl Mathematics Institute\nl University of Warwick\nl Coventry, CV4 7AL, UK\nl
Email:\quad {\tt gt@maths.warwick.ac.uk}\qquad Fax:\quad +44-1203-524182 

{\it For printed copy see:}\nl
{\tt http://www.maths.warwick.ac.uk/gt/gtp-subscription.html}

To order use the on-line order form:\nl
{\tt http://www.maths.warwick.ac.uk/gt/gtp-online-order-form.html}

Or write to Geometry and Topology Orders at the above address or:\nl
Email:\quad {\tt gtorders@maths.warwick.ac.uk}\qquad Fax:\quad +44-1203-524182

\egroup
\np

\include{m5-toc}

\include{m5-p}

\include{m5-a}

\mainmatter

\include{m5-1}

\include{m5-2}

\include{m5-3}

\include{m5-4}

\include{m5-5}

\include{m5-6}

\include{m5-7}

\include{m5-8}

\include{m5-9}

\include{m5-10}

\include{m5-11}

\include{m5-12}

\include{m5-13}
\include{m5-14}

\include{m5-15}

\include{m5-16}

\include{m5-17}

\include{m5-18}

\backmatter

\include{m5-bib}

\include{m5-i}

\end{document}

%% file: m5-toc.tex


\bgroup

\def\contentsline#1#2#3{{#2} \dotfill #3}
\def\numberline#1{#1\qua}

\parskip 2pt plus 2pt

\markboth{Contents}{Contents}

\cl{\Large\bf Contents}\bigskip\bigskip

\leftskip 10pt

Contents \dotfill (iii)

Preface \dotfill (ix)

Acknowledgement \dotfill (xiii)

Changes in later revisions (2007, 2014 and 2022) \dotfill (xiv)

\medskip
\leftskip 0pt

\contentsline {part}{\large\bf Part I : Manifolds and $PD$-complexes}{1}

\medskip

\contentsline {chapter}{\numberline {Chapter 1:}Group theoretic preliminaries}{3}

\medskip\leftskip 10pt

\contentsline {section}{\numberline {1.1}Group theoretic notation and terminology}{3}

\contentsline {section}{\numberline {1.2}Matrix groups}{5}

\contentsline {section}{\numberline {1.3}The Hirsch-Plotkin radical}{7}

\contentsline {section}{\numberline {1.4}Amenable groups}{8}

\contentsline {section}{\numberline {1.5}Hirsch length}{10}

\contentsline {section}{\numberline {1.6}Modules and finiteness conditions}{13}

\contentsline {section}{\numberline {1.7}Ends and cohomology with free coefficients}{16}

\contentsline {section}{\numberline {1.8}Poincar\'e duality groups}{21}

\contentsline {section}{\numberline {1.9}Hilbert modules}{23}

\medskip\leftskip 0pt

\contentsline {chapter}{\numberline {Chapter 2:}2-Complexes and $PD_3$-complexes}{25}

\medskip\leftskip 10pt

\contentsline {section}{\numberline {2.1}Notation}{25}

\contentsline {section}{\numberline {2.2}$L^2$-Betti numbers}{26}

\contentsline {section}{\numberline {2.3}2-Complexes and finitely presentable groups}{28}

\contentsline {section}{\numberline {2.4}Poincar\'e duality}{32}

\contentsline {section}{\numberline {2.5}$PD_3$-complexes}{33}

\contentsline {section}{\numberline {2.6}The spherical cases}{35}

\contentsline {section}{\numberline {2.7}$PD_3 $-groups}{37}

\contentsline {section}{\numberline {2.8}Subgroups of $PD_3$-groups and 3-manifold groups}{42}

\contentsline {section}{\numberline {2.9}$\pi _2(P)$ as a $\mathbb {Z}[\pi ]$-module}{44}

\medskip\leftskip 0pt

\contentsline {chapter}{\numberline {Chapter 3:}Homotopy invariants of $PD_4$-complexes}{47}

\medskip\leftskip 10pt

\contentsline {section}{\numberline {3.1}Homotopy equivalence and asphericity}{47}

\contentsline {section}{\numberline {3.2}Finitely dominated covering spaces}{53}

\contentsline {section}{\numberline {3.3}Minimizing the Euler characteristic}{57}

\contentsline {section}{\numberline {3.4}Euler Characteristic 0}{62}

\contentsline {section}{\numberline {3.5}The intersection pairing}{66}

\medskip\leftskip 0pt

\contentsline {chapter}{\numberline {Chapter 4:}Mapping tori and circle bundles}{69}

\medskip\leftskip 10pt

\contentsline {section}{\numberline {4.1}$PD_r$-covers of $PD_n$-spaces}{70}

\contentsline {section}{\numberline {4.2}Novikov rings and Ranicki's criterion}{73}

\contentsline {section}{\numberline {4.3}Infinite cyclic covers}{76}

\contentsline {section}{\numberline {4.4}The case $n=4$}{78}

\contentsline {section}{\numberline {4.5}Products}{82}

\contentsline {section}{\numberline {4.6}Ascendant subgroups}{83}

\contentsline {section}{\numberline {4.7}Circle bundles}{84}

\medskip\leftskip 0pt

\contentsline {chapter}{\numberline {Chapter 5:}Surface bundles}{89}

\medskip\leftskip 10pt

\contentsline {section}{\numberline {5.1}Some general results}{89}

\contentsline {section}{\numberline {5.2}Bundles with base and fibre aspherical surfaces}{91}

\contentsline {section}{\numberline {5.3}Bundles with aspherical base and fibre $S^2$ or $RP^2$}{97}

\contentsline {section}{\numberline {5.4}Bundles over $S^2 $}{104}

\contentsline {section}{\numberline {5.5}Bundles over $RP^2$}{106}

\contentsline {section}{\numberline {5.6}Bundles over $RP^2$ with $\partial =0$}{108}

\contentsline {section}{\numberline {5.7}Sections of surface bundles}{110}

\medskip\leftskip 0pt

\contentsline {chapter}{\numberline {Chapter 6:}Simple homotopy type and surgery}{111}

\medskip\leftskip 10pt

\contentsline {section}{\numberline {6.1}The Whitehead group}{112}

\contentsline {section}{\numberline {6.2}The $s$-cobordism structure set}{116}

\contentsline {section}{\numberline {6.3}Stabilization and $h$-cobordism}{122}

\contentsline {section}{\numberline {6.4}
Manifolds with $\pi _1$ elementary amenable and $\chi =0$}{123}

\contentsline {section}{\numberline {6.5}Bundles over aspherical surfaces}{126}

\medskip\leftskip 0pt

\contentsline {part}{\large\bf Part II : 4-dimensional Geometries}{129}

\medskip

\contentsline {chapter}{\numberline {Chapter 7:}Geometries and decompositions}{131}

\medskip\leftskip 10pt

\contentsline {section}{\numberline {7.1}Geometries}{132}

\contentsline {section}{\numberline {7.2}Infranilmanifolds}{133}

\contentsline {section}{\numberline {7.3}Infrasolvmanifolds}{135}

\contentsline {section}{\numberline {7.4}Orbifold bundles}{138}

\contentsline {section}{\numberline {7.5}Geometric decompositions}{139}

\contentsline {section}{\numberline {7.6}
Realization of virtual bundle groups}{143}

\contentsline {section}{\numberline {7.7}Seifert fibrations}{145}

\contentsline {section}{\numberline {7.8}
Complex surfaces and related structures}{148}

\medskip\leftskip 0pt

\contentsline {chapter}{\numberline {Chapter 8:}Solvable Lie geometries}{151}

\medskip\leftskip 10pt

\contentsline {section}{\numberline {8.1}The characterization}{151}

\contentsline {section}{\numberline {8.2}Flat 3-manifold groups and their automorphisms}{153}

\contentsline {section}{\numberline {8.3}Flat 4-manifold groups with infinite abelianization}{157}

\contentsline {section}{\numberline {8.4}Flat 4-manifold groups with finite abelianization}{161}

\contentsline {section}{\numberline {8.5}Distinguishing between the geometries}{164}

\contentsline {section}{\numberline {8.6}Mapping tori of self homeomorphisms of $\mathbb {E}^3 $-manifolds}{166}

\contentsline {section}{\numberline {8.7}Mapping tori of self homeomorphisms of $\mathbb {N}il^3$-manifolds}{168}

\contentsline {section}{\numberline {8.8}Mapping tori of self homeomorphisms of $\mathbb {S}ol^3 $-manifolds}{172}

\contentsline {section}{\numberline {8.9}Realization and classification}{174}

\contentsline {section}{\numberline {8.10}Diffeomorphism}{177}

\medskip\leftskip 0pt

\contentsline {chapter}{\numberline {Chapter 9:}The other aspherical geometries}{179}

\medskip\leftskip 10pt

\contentsline {section}{\numberline {9.1}Aspherical Seifert fibred 4-manifolds}{179}

\contentsline {section}{\numberline {9.2}The Seifert geometries: 
$\mathbb {H}^2\times \mathbb {E}^2$ and $\widetilde{\mathbb{SL}}\times\mathbb{E}^1$}{182}

\contentsline {section}{\numberline {9.3}$\mathbb {H}^3\times \mathbb {E}^1$-manifolds}{185}

\contentsline {section}{\numberline {9.4}Mapping tori}{186}

\contentsline {section}{\numberline {9.5}The semisimple geometries: 
$\mathbb {H}^2\times \mathbb {H}^2$, 
$\mathbb {H}^4$ and $\mathbb {H}^2(\mathbb {C})$}{188}

\contentsline {section}{\numberline {9.6}Miscellany}{193}

\medskip\leftskip 0pt

\contentsline {chapter}{\numberline {Chapter 10:}Manifolds covered by $S^2 \times\mathbb{R}^2$}{195}

\medskip\leftskip 10pt

\contentsline {section}{\numberline {10.1}Fundamental groups}{195}

\contentsline {section}{\numberline {10.2}The first $k$-invariant}{196}

\contentsline {section}{\numberline {10.3}Homotopy type}{199}

\contentsline {section}{\numberline {10.4}Bundle spaces are geometric}{203}

\contentsline {section}{\numberline {10.5}Fundamental groups of 
$\mathbb {S}^2\times \mathbb {E}^2$-manifolds}{208}

\contentsline {section}{\numberline {10.6}
Homotopy types of $\mathbb {S}^2\times \mathbb {E}^2$-manifolds}{210}

\contentsline {section}{\numberline {10.7}
Some remarks on the homeomorphism types}{213}

\contentsline {section}{\numberline {10.8}Minimal models}{213}

\medskip\leftskip 0pt

\contentsline {chapter}{\numberline {Chapter 11:}Manifolds covered by 
$S^3\times\mathbb{R}$}{217}

\medskip\leftskip 10pt

\contentsline {section}{\numberline {11.1}Invariants for the homotopy type}{217}

\contentsline {section}{\numberline {11.2}The action of $\pi /F$ on $F$}{219}

\contentsline {section}{\numberline {11.3}Extensions of $D$}{223}

\contentsline {section}{\numberline {11.4}$\mathbb {S}^3 \times \mathbb {E}^1 $-manifolds}{224}

\contentsline {section}{\numberline {11.5}Realization of the invariants}{226}

\contentsline {section}{\numberline {11.6}$T$- and $Kb$-bundles over $RP^2$ 
with $\partial \not =0$}{229}

\contentsline {section}{\numberline {11.7}Some remarks on the homeomorphism types}{232}

\medskip\leftskip 0pt

\contentsline {chapter}{\numberline {Chapter 12:}Geometries with compact models}{233}

\medskip\leftskip 10pt

\contentsline {section}{\numberline {12.1}The geometries $\mathbb {S}^4$ and $\mathbb {CP}^2$}{234}

\contentsline {section}{\numberline {12.2}The geometry $\mathbb {S}^2\times \mathbb {S}^2$}{235}

\contentsline {section}{\numberline {12.3}Bundle spaces}{236}

\contentsline {section}{\numberline {12.4}Cohomology and Stiefel-Whitney classes}{238}

\contentsline {section}{\numberline {12.5}The action of $\pi $ on $\pi _2 (M)$}{239}

\contentsline {section}{\numberline {12.6}Homotopy type}{241}

\contentsline {section}{\numberline {12.7}Surgery}{244}

\contentsline {section}{\numberline {12.8}A smooth fake version 
of $S^2\times{S^2}/\langle\tau(I,-I)\rangle$?}{247}

\medskip\leftskip 0pt

\contentsline {chapter}{\numberline {Chapter 13:}Geometric decompositions of bundle spaces}{249}

\medskip\leftskip 10pt

\contentsline {section}{\numberline {13.1}Mapping tori}{249}

\contentsline {section}{\numberline {13.2}Surface bundles and geometries}{254}

\contentsline {section}{\numberline {13.3}Geometric decompositions of 
Seifert fibred $4$-manifolds}{257}

\contentsline {section}{\numberline {13.4}Complex surfaces and fibrations}{259}

\contentsline {section}{\numberline {13.5}$S^1$-Actions and foliations by circles}{263}

\contentsline {section}{\numberline {13.6}Symplectic structures}{265}

\medskip\leftskip 0pt

\contentsline {part}{\large\bf Part III : 2-Knots}{267}

\medskip

\contentsline {chapter}{\numberline {Chapter 14:}Knots and links}{269}

\medskip\leftskip 10pt

\contentsline {section}{\numberline {14.1}Knots}{269}

\contentsline {section}{\numberline {14.2}Covering spaces}{271}

\contentsline {section}{\numberline {14.3}Sums, factorization and satellites}{272}

\contentsline {section}{\numberline {14.4}Spinning, twist spinning
and deform spinning}{273}

\contentsline {section}{\numberline {14.5}Ribbon and slice knots}{274}

\contentsline {section}{\numberline {14.6}The Kervaire conditions}{276}

\contentsline {section}{\numberline {14.7}Weight elements, classes and orbits}{277}

\contentsline {section}{\numberline {14.8}The commutator subgroup}{279}

\contentsline {section}{\numberline {14.9}
Deficiency and geometric dimension}{281}

\contentsline {section}{\numberline {14.10}Asphericity}{283}

\contentsline {section}{\numberline {14.11}Links}{284}

\contentsline {section}{\numberline {14.12}Link groups}{288}

\contentsline {section}{\numberline {14.13}Homology spheres}{290}

\medskip\leftskip 0pt

\contentsline {chapter}{\numberline {Chapter 15:}Restrained normal subgroups}{293}

\medskip\leftskip 10pt

\contentsline {section}{\numberline {15.1}Finite commutator subgroup}{293}

\contentsline {section}{\numberline {15.2}Finite normal subgroups}{296}

\contentsline {section}{\numberline {15.3}The group $\Phi $}{297}

\contentsline {section}{\numberline {15.4}Nilpotent normal subgroups}{299}

\contentsline {section}{\numberline {15.5}Almost coherent, restrained and locally virtually indicable}{303}

\contentsline {section}{\numberline {15.6}A Tits alternative for
 2-knot groups}{306}

\contentsline {section}{\numberline {15.7}Abelian HNN bases}{307}

\contentsline {section}{\numberline {15.8}Locally-finite normal 
subgroups}{309}

\medskip\leftskip 0pt

\contentsline {chapter}{\numberline {Chapter 16:}Abelian normal subgroups of rank $\geq 2$}{311}

\medskip\leftskip 10pt

\contentsline {section}{\numberline {16.1}The Brieskorn manifolds $M(p,q,r)$}{311}

\contentsline {section}{\numberline {16.2}Rank 2 subgroups}{312}

\contentsline {section}{\numberline {16.3}Twist spins of torus knots}{314}

\contentsline {section}{\numberline {16.4}Solvable $PD_4$-groups}{319}

\medskip\leftskip 0pt

\contentsline {chapter}{\numberline {Chapter 17:}Knot manifolds and geometries}{327}

\medskip\leftskip 10pt

\contentsline {section}{\numberline {17.1}Homotopy classification of $M(K)$}{327}

\contentsline {section}{\numberline {17.2}Surgery}{328}

\contentsline {section}{\numberline {17.3}The aspherical cases}{329}

\contentsline {section}{\numberline {17.4}Quasifibres and minimal Seifert hypersurfaces}{330}

\contentsline {section}{\numberline {17.5}The spherical cases}{332}

\contentsline {section}{\numberline {17.6}Finite geometric dimension 2}{333}

\contentsline {section}{\numberline {17.7}Geometric 2-knot manifolds}{336}

\contentsline {section}{\numberline {17.8}Complex surfaces and 2-knot manifolds}{338}

\medskip\leftskip 0pt

\contentsline {chapter}{\numberline {Chapter 18:}Reflexivity}{341}

\medskip\leftskip 10pt

\contentsline {section}{\numberline {18.1}Sections of the mapping torus}{341}

\contentsline {section}{\numberline {18.2}Reflexivity for fibred 2-knots}{342}

\contentsline {section}{\numberline {18.3}Cappell-Shaneson knots}{345}

\contentsline {section}{\numberline {18.4}The Hantzsche-Wendt flat 3-manifold}{348}

\contentsline {section}{\numberline {18.5} 2-knots with group $G(\pm)$}{349}

\contentsline {section}{\numberline {18.6}$\mathbb {N}il^3$-fibred knots}{352}

\contentsline {section}{\numberline {18.7}Other geometrically fibred knots}{357}

\contentsline {section}{\numberline {18.8}Smooth 2-knots in $S^4$}{362}

\medskip

Bibliography \dotfill 365\par
Index \dotfill 392

\egroup


%% file: m5-p.tex
\chapter*{Preface}

Every closed surface admits a geometry of constant curvature, 
and may be classified topologically either by its
fundamental group or by its Euler characteristic and orientation character.
Closed 3-manifolds have decompositions into geometric pieces, 
and are determined up to homeomorphism by invariants associated 
with the fundamental group (whereas the Euler characteristic is always 0).
In dimension 4 the Euler characteristic and fundamental group are largely 
independent, and the class of closed 4-manifolds which admit a geometric 
decomposition is rather restricted. 
For instance, there are only 11 such manifolds with finite fundamental group.
On the other hand, many complex surfaces admit geometric structures, 
as do all the manifolds arising from surgery on twist spun simple knots.

The goal of this book is to characterize algebraically the closed 4-manifolds 
that fibre nontrivially or admit geometries, or which are obtained by surgery 
on 2-knots, and to provide a reference for the topology of such manifolds and 
knots.
In many cases the Euler characteristic, fundamental group and Stiefel-Whitney 
classes together form a complete system of invariants for the homotopy type of 
such manifolds, 
and the possible values of the invariants can be described explicitly.
If the fundamental group is elementary amenable we may use topological surgery 
to obtain classifications up to homeomorphism.
Surgery techniques also work well ``stably" in dimension 4 
(i.e., modulo connected sums with copies of $S^2\times S^2$).
However, in our situation the fundamental group may have nonabelian free 
subgroups and the Euler characteristic is usually the minimal possible for 
the group, and it is not known whether $s$-cobordisms between such 4-manifolds 
are always topologically products.
Our strongest results are characterizations of infrasolvmanifolds (up to 
homeomorphism) and aspherical manifolds which fibre over a surface or which 
admit a geometry of rank $>1$ (up to TOP $s$-cobordism).
As a consequence 2-knots whose groups are poly-$Z$ are determined up to Gluck 
reconstruction and change of orientations by their groups alone.

We shall now outline the chapters in somewhat greater detail.
The first chapter is purely algebraic; here we summarize the relevant group 
theory and present the notions of amenable group,
Hirsch length of an elementary amenable group, 
finiteness conditions,
criteria for the vanishing of cohomology of a 
group with coefficients in a free module,
Poincar\'e duality groups, and Hilbert modules over the 
von Neumann algebra of a group. 
The rest of the book may be divided into three parts: general results on 
homotopy and surgery (Chapters 2-6), geometries and geometric decompositions 
(Chapters 7-13), and 2-knots (Chapters 14-18).            
                           
Some of the later arguments are applied in microcosm to 2-complexes and 
$PD_3$-complexes in Chapter 2, which presents equivariant cohomology, 
$L^2$-Betti numbers and Poincar\'e duality.
Chapter 3 gives general criteria for two closed 4-manifolds to be homotopy 
equivalent, and we show that a closed 4-manifold $M$ is aspherical if and only 
if $\pi_1(M)$ is a $PD_4$-group of type $FF$ and $\chi(M)=\chi(\pi)$.  
We show that if the universal cover of a closed 4-manifold 
is finitely dominated 
then it is contractible or homotopy equivalent to $S^2$ or $S^3$ or the 
fundamental group is finite. 
We also consider at length the relationship between fundamental group and 
Euler characteristic for closed 4-manifolds.
In Chapter 4 we show that a closed 4-manifold $M$ fibres homotopically over 
$S^1$ with fibre a $PD_3$-complex if and only if $\chi(M)=0$ and 
$\pi_1(M)$ is an extension of $\mathbb{Z}$ by 
a finitely presentable normal subgroup. 
(There remains the problem of recognizing which $PD_3$-complexes are
homotopy equivalent to 3-manifolds).
The dual problem of characterizing the total spaces of $S^1$-bundles over 
3-dimensional bases seems more difficult.
We give a criterion that applies
under some restrictions on the fundamental group.
In Chapter 5 we characterize the homotopy types of total spaces of surface 
bundles.
(Our results are incomplete if the base is $RP^2$).
In particular, a closed 4-manifold $M$ is simple homotopy equivalent to 
the total space of an $F$-bundle over $B$ (where $B$ and $F$ are closed 
surfaces and $B$ is aspherical) if and only if $\chi(M)=\chi(B)\chi(F)$ and 
$\pi_1(M)$ is an extension of $\pi_1(B)$ by a normal subgroup isomorphic to 
$\pi_1 (F)$. (The extension should split if $F=RP^2$).
Any such extension is the fundamental group of such a bundle space; 
the bundle is determined by the extension of groups in the aspherical cases 
and by the group and Stiefel-Whitney classes if the fibre is $S^2 $ or $RP^2$. 
This characterization is improved in Chapter 6, which considers Whitehead 
groups and obstructions to constructing $s$-cobordisms via surgery.

The next seven chapters consider geometries and geometric decompositions.
Chapter 7 introduces the 4-dimensional geometries and demonstrates the 
limitations of geometric methods in this dimension. 
It also gives a brief outline of the connections between geometries, 
Seifert fibrations and complex surfaces.
In Chapter 8 we show that a closed 4-manifold $M$ is homeomorphic to an 
infrasolvmanifold if and only if $\chi(M)=0$ and $\pi_1 (M)$ has a locally 
nilpotent normal subgroup of Hirsch length at least 3, 
and two such manifolds are homeomorphic if and only if their fundamental groups 
are isomorphic. 
Moreover $\pi_1 (M)$ is then a torsion free virtually poly-$Z$ group of Hirsch 
length 4 and every such group is the fundamental group of an infrasolvmanifold.
We also consider in detail the question of when such a manifold is the mapping 
torus of a self homeomorphism of a 3-manifold, and give a direct and 
elementary derivation of the fundamental groups of flat 4-manifolds.
At the end of this chapter we show that all orientable 4-dimensional
infrasolvmanifolds are determined up to diffeomorphism by their fundamental
groups. (The corresponding result in other dimensions was known).

Chapters 9-12 consider the remaining 4-dimensional geometries, grouped 
according to whether the model is homeomorphic to $\mathbb{R}^4$, 
$S^2\times\mathbb{R}^2$, $S^3 \times\mathbb{R}$ or is compact.
Aspherical geometric 4-manifolds are determined up to $s$-cobordism by
their homotopy type. However there are only partial characterizations of 
the groups arising as fundamental groups of 
$\mathbb{H}^2\times\mathbb{H}^2$-manifolds, 
while very little is known about $\mathbb{H}^4$- 
or $\mathbb{H}^2(\mathbb{C})$-manifolds.
We show that the homotopy types of manifolds covered by $S^2\times\mathbb{R}^2$ 
are determined up to finite ambiguity by their fundamental groups.
If the fundamental group is torsion free such a manifold is $s$-cobordant to 
the total space of an $S^2$-bundle over an aspherical surface.
The homotopy types of manifolds covered by $S^3\times\mathbb{R}$ 
are determined by the fundamental group and first nonzero $k$-invariant; 
much is known about the possible fundamental groups, 
but less is known about which $k$-invariants are realized.
Moreover, although the fundamental groups are all ``good", so that in principle 
surgery may be used to give a classification up to homeomorphism,
the problem of computing surgery obstructions seems very difficult.
We conclude the geometric section of the book in Chapter 13 by considering
geometric decompositions of 4-manifolds which are also mapping tori or
total spaces of surface bundles, and we characterize the complex surfaces 
which fibre over $S^1$ or over a closed orientable 2-manifold.

The final five chapters are on 2-knots.
Chapter 14 is an overview of knot theory; in particular it is shown how the 
classification of higher-dimensional knots may be largely reduced to the 
classification of knot manifolds.
The knot exterior is determined by the knot manifold and 
the conjugacy class of a normal generator for the knot group, 
and at most two knots share a given exterior.
An essential step is to characterize 2-knot groups.
Kervaire gave homological conditions which characterize high dimensional knot 
groups and which 2-knot groups must satisfy, and showed that any high 
dimensional knot group with a presentation of deficiency 1 is a 2-knot group.
Bridging the gap between the homological and combinatorial conditions appears 
to be a delicate task.
In Chapter 15 we investigate 2-knot groups with infinite
normal subgroups which have no noncyclic free subgroups.                                  
We show that under mild coherence hypotheses such 2-knot groups
usually have nontrivial abelian normal subgroups,
and we determine all 2-knot groups with finite commutator subgroup.
In Chapter 16 we show that if there is an abelian normal subgroup of rank $>1$ 
then the knot manifold is either $s$-cobordant to a 
$\widetilde{\mathbb{SL}}\times\mathbb{E}^1$-manifold or
is homeomorphic to an infrasolvmanifold.
In Chapter 17 we characterize the closed 4-manifolds obtained by surgery on 
certain 2-knots, and show that just eight of the 4-dimensional geometries
are realised by knot manifolds.
We also consider when the knot manifold admits a complex structure.
The final chapter considers when a fibred 2-knot with geometric fibre
is determined by its exterior.
We settle this question when the monodromy has finite order
or when the fibre is $\mathbb{R}^3/\mathbb{Z}^3$ 
or is a coset space of the Lie group $Nil^3$.

This book arose out of two earlier books of mine, on ``{\it 2-Knots and their 
Groups}" and ``{\it The Algebraic Characterization of Geometric 4-Manifolds}",
published by Cambridge University Press for the Australian Mathematical Society
and for the London Mathematical Society, respectively.
About a quarter of the present text has been taken from these 
books.\footnote{See the following Acknowledgment for a
summary of the textual borrowings.}
However the arguments have been improved in many cases, notably in
using Bowditch's homological criterion for virtual surface groups
to streamline the results on surface bundles,
using $L^2$-methods instead of localization,
completing the characterization of mapping tori, 
relaxing the hypotheses on torsion 
or on abelian normal subgroups in the fundamental group
and in deriving the results on 2-knot groups from the work on 4-manifolds.
The main tools used here beyond what can be found in {\it Algebraic Topology} 
\cite{[Sp]} are cohomology of groups, equivariant Poincar\'e duality
and (to a lesser extent) $L^2$-(co)homology.
Our references for these are the books {\it Homological Dimension of Discrete 
Groups} \cite{[Bi]}, {\it Surgery on Compact Manifolds} \cite{[Wl]}
and {\it $L^2$-Invariants: Theory and Applications to Geometry and $K$-Theory}
\cite{[Lu]}, respectively.
We also use properties of 3-manifolds (for the construction of examples)
and calculations of Whitehead groups and surgery obstructions.

This work has been supported in part by ARC small grants, 
enabling visits by Steve Plotnick, Mike Dyer, Charles Thomas and Fang Fuquan. 
I would like to thank them all for their advice, and in particular 
Steve Plotnick for the collaboration reported in Chapter 18. 
I would also like to thank Robert Bieri, Robin Cobb, Peter Linnell 
and Steve Wilson for their collaboration, and Warren Dicks, William Dunbar, 
Ross Geoghegan, F.T.Farrell, Ian Hambleton, Derek Holt, K.F.Lai, Eamonn O'Brien,
Peter Scott and Shmuel Weinberger for their correspondance and advice on aspects
of this work.

\rightline{\small\it Jonathan Hillman}

%% file: m5-a.tex
\chapter*{Acknowledgment}

I wish to thank Cambridge University Press for their permission to use
material from my earlier books [H1] and [H2].
The textual borrowings in each Chapter [of the 2002 version] are outlined below.

\smallskip
\noindent 1. \S1, Lemmas 1.7 and 1.10 and Theorem 1.11, 
\S6 (up to the discussion of $\chi(\pi)$),
the first paragraph of \S7 and Theorem 1.15 are from [H2:Chapter I].
(Lemma 1.1 is from [H1]). \S3 is from [H2:Chapter VI].

\noindent 2. 
\S1, most of \S4 and part of \S5 are from [H2:Chapter II and Appendix].

\noindent 3. 
Lemma 3.1, Theorems 3.2, 3.7-3.9 and Corollaries 3.9.1-3.9.3 of
Theorem 3.12 are from [H2:Chapter II]. (Theorem 3.9 has been improved).

\noindent 4.
The statements of Corollaries 4.5.1-4.5.3,
Corollary 4.5.4 and most of \S7 are from [H2:Chapter III].
(Theorem 11 and the subsequent discussion have been improved).
 
\noindent 5. Part of Lemma 5.15 and \S4-\S5 are from [H2:Chapter IV]. 
(Theorem 5.19 and Lemmas 5.21 and 5.22 have been improved).

\noindent 6.
\S1 (excepting Theorem 6.1), Theorem 6.12 and the proof of Theorem 6.14 are 
from [H2:Chapter V].

\noindent 8.
Part of Theorem 8.1, \S6, most of \S7 and \S8 are from [H2:Chapter VI].

\noindent 9. Theorems 9.1, 9.2 and 9.7 are from [H2:Chapter VI], 
with improvements. 

\noindent 10. 
Theorems 10.10-10.12 and \S6 are largely from [H2:Chapter VII]. 
(Theorem 10.10 has been improved).

\noindent 11.
Lemma 11.3, \S3 and the first three paragraphs of \S5 are 
from [H2:Chapter VIII].
\S6 is from [H2:Chapter IV].

\noindent 12. The introduction, \S1-\S3, \S5, most of
\S6 (from Lemma 12.5 onwards) and \S7 are from [H2:Chapter IX], 
with improvements (particularly in \S7). 

\noindent 14. 
\S1-\S5 are from [H1:Chapter I].
\S6 and \S7 are from [H1:Chapter II]. 

\noindent 16.
Most of \S3 is from [H1:Chapter V].(Theorem 16.4 is new
and Theorems 16.5 and 16.6 have been improved).

\noindent 17.
Lemma 2 and Theorem 7 are from [H1:Chapter VIII],
while Corollary 17.6.1 is from [H1:Chapter VII].
The first two paragraphs of \S8 and Lemma 17.12 are from [H2:Chapter X].
 
\newpage
{\bf [Added in later revisions]}

In 2007 some of the material was improved, particularly as regards 

(a) {\sl finiteness conditions} (Chapters 3 and 4); and 

(b) {\sl (aspherical) Seifert fibred 4-manifolds} (Chapters 7 and 9).

Some results on the equivariant intersection pairing and the notion of
(strongly) minimal $PD_4$-complex were added as new sections \S3.5 and 
(now) \S10.8. 

Further improvements were made in 2014,
particularly as regards 

(c) {\sl non-aspherical geometric 4-manifolds} (Chapters 10--12).

In Chapter 10 we have shown that every closed 4-manifold 
with universal cover $\simeq{S^2}$ is homotopy equivalent to 
either the total space of an $S^2$-orbifold bundle 
over an aspherical 2-orbifold or to the total space 
of an $RP^2$-bundle over an aspherical surface.
(Most such $S^2$-orbifold bundle spaces are geometric \cite{[Hi13]}.)
In Chapter 11 we have repaired an oversight (Milnor's groups $P_{48r}''$),
and summarized work by Davis and Weinberger \cite{[DW07]} 
on mapping tori of orientation-reversing self homotopy
equivalences of lens spaces.

We have cited Perelman's work on geometrization to simplify some statements, 
and have used the Virtual Fibration Theorem of Agol 
\cite{[Ag13]} in Chapters 9 and 13, 
in connection with the geometry $\mathbb{H}^3\times\mathbb{E}^1$.
We have used the notion of orbifold more widely,
improved the discussion of surgery in Chapter 6
and tightened some of the results on 2-knots.
In particular, there is a new family of 2-knots with 
torsion free, solvable groups, which was overlooked before.
(See Theorem 16.15.) 
Six pages have been added to Chapter 18.
In these we settle the questions of reflexivity, 
amphicheirality and invertibility
for the other such knots.

The classification of 4-dimensional infrasolvmanifold groups 
is now essentially known \cite{[LT15],[Hi18]}.
The discussion of quotients of $S^2\times{S^2}$ 
by $Z/4Z$ (in Chapter 12) was rewritten in 2021.
Chapter 15 has also been revised extensively,
and a section on smooth 2-knots in $S^4$
has been added to Chapter 18.

The errors and typos discovered to date (14 November 2022) have been corrected.

\smallskip
I would like to thank J.F. Davis and S.Weinberger for permission to
include a summary of \cite{[DW07]},
I.Hambleton, M.Kemp, D.H.Kochloukova and S.K.Roushon for their collaboration 
in relation to some of the improvements recorded here,
and J.G.Ratcliffe for alerting me to some gaps in \S4 of Chapter 8.

%% file: m5-1.tex
\part{Manifolds and $PD$-complexes}

\chapter{Group theoretic preliminaries}

The key algebraic idea used in this book is to study the homology groups
of covering spaces as modules over the group ring of the group of covering 
transformations.
In this chapter we shall summarize the relevant notions from group theory,
in particular, the Hirsch-Plotkin radical, amenable groups, Hirsch length, 
finiteness conditions, the connection between ends and the vanishing 
of cohomology with coefficients in a free module, 
Poincar\'e duality groups and Hilbert modules.

Our principal references for group theory are \cite{[Bi]}, \cite{[DD]} and \cite{[Ro]}.

\section{Group theoretic notation and terminology}

We shall write $\mathbb{Z}$ for the ring of integers,
and also the free (abelian) group of rank 1.
We may otherwise write $Z$ for an infinite cyclic group 
with no specified generator.
We shall also identify the units $\mathbb{Z}^\times=\{\pm1\}$,
the field $\mathbb{F}_2$ and $Z/2Z$, when convenient.
Let $F(r)$ be the free group of rank $r$.

Let $G$ be a group. 
Then $G'$ and $\zeta G$ denote the commutator subgroup and 
centre of $G$, respectively.
The {\it outer automorphism group} of $G$ is $Out(G)=Aut(G)/Inn(G)$,
where $Inn(G)\cong G/\zeta G$ is the subgroup of $Aut(G)$ consisting of conjugations by elements of $G$.
If $H$ is a subgroup of $G$ let $N_G (H)$ and $C_G (H)$ 
denote the normalizer and centralizer of $H$ in $G$, respectively.
The subgroup $H$ is a {\it characteristic} subgroup of $G$ 
if it is preserved under all automorphisms of $G$.
In particular, $I(G)=\{ g\in G\mid \exists n>0,\medspace g^n\in G'\}$ 
is a characteristic subgroup of $G$, 
and the quotient $G/I(G)$ is a torsion-free abelian group 
of rank $\beta_1 (G)$.
A group $G$ is {\it indicable} if there is an epimorphism $p:G\to\mathbb{Z}$, 
or if $G=1$.
If $S$ is a subset of $G$ then $\langle{S}\rangle$ and 
$\langle\langle S\rangle\rangle_G$ (or just $\langle\langle S\rangle\rangle$)
are the subgroup generated by $S$ and the {\it normal closure} of 
$S$ in $G$ (the intersection of the
normal subgroups of $G$ which contain $S$), respectively.

If $P$ and $Q$ are classes of groups let $PQ$ denote the class of 
(``$P$ {\it by} $Q$")
groups $G$ which have a normal subgroup $H$ in $P$ such that the quotient $G/H$ 
is in $Q$, and let $\ell P$ denote the class of (``{\it locally} $P$") groups 
such that each finitely generated subgroup is in the class $P$.
In particular, if $F$ is the class of finite groups $\ell F$ is the class of 
{\it locally finite} groups. 
In any group the union of all the locally-finite normal 
subgroups is the unique maximal locally-finite normal subgroup. 
Clearly there are no nontrivial homomorphisms from such a group 
to a torsion-free group.
Let {\it poly}-$P$ be the class of groups with a finite composition
series such that each subquotient is in $P$. 
Thus if $Ab$ is the class of abelian groups
poly-$Ab$ is the class of solvable groups.
                                             
Let $P$ be a class of groups which is closed under taking subgroups
of finite index.
A group is {\it virtually} $P$ if it has a subgroup of finite index in $P$.
Let $vP$ be the class of groups which are virtually $P$.
Thus a {\it virtually poly-}$Z$ group is one which has a subgroup of finite
index with a composition series whose factors are all infinite cyclic.
The number of infinite cyclic factors is independent of the choice of
finite index subgroup or composition series, 
and is called the {\it Hirsch length} of the group.
We shall also say that a space virtually has some property if it has a finite
regular covering space with that property.

If $p:G\to Q$ is an epimorphism with kernel $N$ we shall say that 
$G$ is an {\it extension of $Q=G/N$ by the normal subgroup $N$}.                 The action of $G$ on $N$ by conjugation determines a homomorphism 
from $G$ to $Aut(N)$ with kernel $C_G(N)$ and hence a homomorphism 
from $G/N$ to $Out(N)=Aut(N)/Inn(N)$.
If $G/N\cong\mathbb{Z}$ the extension splits: 
a choice of element $t$ in $G$ which 
projects to a generator of $G/N$ determines a right inverse to $p$.              Let $\theta$ be the automorphism of $N$ determined by conjugation by $t$ in $G$.
Then $G$ is isomorphic to the semidirect product $N\rtimes_\theta\mathbb{Z}$.
Every automorphism of $N$ arises in this way, and automorphisms whose images 
in $Out(N)$ are conjugate determine isomorphic semidirect products.
In particular, $G\cong N\times\mathbb{Z}$ if $\theta$ is an inner automorphism.

\begin{lemma}
Let $\theta$ and $\phi$ automorphisms of a group $G$
such that $H_1(\theta;\mathbb{Q})-1$ and $H_1(\phi;\mathbb{Q})-1$ are automorphisms of 
$H_1(G;\mathbb{Q})=(G/G')\otimes \mathbb{Q}$.
Then the semidirect products $\pi_\theta=G\rtimes_\theta\mathbb{Z}$ and 
$\pi_\phi=G\rtimes_\phi\mathbb{Z}$ are isomorphic
if and only if $\theta$ is conjugate to $\phi$ or $\phi^{-1}$ in $Out(G)$.
\end{lemma}

\begin{proof}
Let $t$ and $u$ be fixed elements of $\pi_\theta$ and $\pi_\phi$, respectively,
which map to 1 in $\mathbb{Z}$.
Since $H_1(\pi_\theta;\mathbb{Q})\cong H_1(\pi_\phi;\mathbb{Q})\cong Q$ the image
of $G$ in each group is characteristic. Hence an isomorphism $h:\pi_\theta\to\pi_\phi$
induces an isomorphism $e:\mathbb{Z}\to\mathbb{Z}$ of the quotients, 
for some $e=\pm 1$, and so $h(t)=u^eg$ for some $g$ in $G$.  
Therefore $h(\theta(h^{-1}(j))))=h(th^{-1}(j)t^{-1})=u^egjg^{-1}u^{-e}=\phi^e(gjg^{-1})$  for all $j$ in $G$.
Thus $\theta$ is conjugate to $\phi^e$ in $Out(G)$.

Conversely, if $\theta$ and $\phi^e$ are conjugate in $Out(G)$
there is an $f$ in $Aut(G)$ and a $g$ in $G$ such that
$\theta(j)=f^{-1}\phi^ef(gjg^{-1})$ for all $j$ in $G$.
Hence $F(j)=f(j)$ for all $j$ in $G$ and $F(t)=u^ef(g)$ defines an isomorphism $F:\pi_\theta\to\pi_\phi$.
\end{proof}

A subgroup $K$ of a group $G$ is {\it ascendant} if there is 
an increasing sequence of subgroups $N_\alpha$, 
indexed by ordinals $\leq\beth$, such that $N_0=K$, 
$N_\alpha$ is normal in $N_{\alpha+1}$ if $\alpha<\beth$,
$N_\beta=\cup_{\alpha<\beta}N_\alpha$ for all limit ordinals
$\beta\leq\beth$ and $N_\beth=G$.
If $\beth$ is finite $K$ is {\it subnormal} in $G$.
Such ascendant series are well suited to arguments by transfinite induction.

\section{Matrix groups}

In this section we shall recall some useful facts about matrices over $\mathbb{Z}$.                       

\begin{lemma}
Let $p$ be an odd prime. 
Then the kernel of the reduction modulo $(p)$ homomorphism from $SL(n,\mathbb{Z})$ 
to $SL(n,\mathbb{F}_p)$ is torsion-free.
\end{lemma}

\begin{proof}
This follows easily from the observation that if 
$A$ is an integral matrix and $k=p^vq$ with $q$ not divisible by $p$
then $(I+p^rA)^k\equiv I+kp^rA$ {\it mod} $(p^{2r+v})$, and $kp^r\not\equiv0$ 
{\it mod} $(p^{2r+v})$ if $r\geq1$. 
\end{proof}

Similarly, the kernel of reduction {\it mod} (4) is torsion-free.

Since $SL(n,\mathbb{F}_p)$ has order $(\Pi_{j=0}^{j=n-1} (p^n-p^j))/(p-1)$, 
it follows that the order of any finite subgroup of $SL(n,\mathbb{Z})$ must 
divide the highest common factor of these numbers, as $p$ varies over all odd primes.
In particular, finite subgroups of $SL(2,\mathbb{Z})$ have order dividing 24, and so are solvable.                            

Let $A=\left(\begin{smallmatrix}
0&1\\ -1&0
\end{smallmatrix}\right)$, 
$B=\left(\begin{smallmatrix}
0&-1\\ 1&1
\end{smallmatrix}\right)$ 
and $R=\left(\begin{smallmatrix}
0&1\\ 1&0
\end{smallmatrix}\right)$. 
Then $A^2=B^3=-I$ and $A^4=B^6=I$.                                   
The matrices $A$ and $R$ generate a dihedral group of order 8,
while $B$ and $R$ generate a dihedral group of order 12. 

\begin{theorem}
Let $G$ be a nontrivial finite subgroup of $GL(2,\mathbb{Z})$.
Then $G$ is conjugate to one of the cyclic groups generated by $A$, 
$A^2=-I$, $B$, $B^2$, $R$ or $RA$, 
or to one of the dihedral groups generated by 
$\{ A,R\}$, $\{-I, R\}$, $\{ A^2,RA\}$, $\{ B,R\}$, $\{ B^2,R\}$ or 
$\{ B^2,RB\}$.
If $G\not=\langle-I_2\rangle$ then $N_{GL(2,\mathbb{Z})}(G)$ is finite.
\end{theorem}

\begin{proof}  
If $M\in GL(2,\mathbb{Z})$ has finite order 
then its characteristic polynomial has cyclotomic factors.
If the characteristic polynomial is $(X\pm1)^2$ then $M=\mp I$.
(This uses the finite order of $M$.)
If the characteristic polynomial is $X^2-1$ then $M$ is conjugate to $R$ or $RA$.
If the characteristic polynomial is $X^2+1$, $X^2-X+1$ or $X^2+X+1$
then it is irreducible, and the corresponding ring of algebraic numbers is a PID. 
Since any $\mathbb{Z}$-torsion-free module over such a ring is free 
it follows easily that $M$ is conjugate to $A$, $B$ or $B^2$.

The normalizers in $SL(2,\mathbb{Z})$ of the subgroups generated by $A$, 
$B$ or $B^2$ are easily seen to be finite cyclic. 
Since $G\cap SL(2,\mathbb{Z})$ is solvable it must be cyclic also.
As it has index at most 2 in $G$ the rest of the theorem follows easily.
\end{proof}

Although the 12 groups listed in the theorem represent distinct conjugacy classes in $GL(2,\mathbb{Z})$,
some of these conjugacy classes coalesce in $GL(2,\mathbb{R})$. 
(For instance, $R$ and $RA$ are conjugate in $GL(2,\mathbb{Z}[\frac12 ])$.)

\begin{cor}
Let $G$ be a locally finite subgroup of 
$GL(2,\mathbb{Q})$. 
Then $G$ is finite, and is conjugate to one of the above subgroups of 
$GL(2,\mathbb{Z})$.
\end{cor}

\begin{proof}
Let $L$ be a finitely generated subgroup of rank 2 in $\mathbb{Q}^2$.
If $G$ is finite then $\cup_{g\in G} gL$ is finitely generated, $G$-invariant
and of rank 2,
and so $G$ is conjugate to a subgroup of $GL(2,\mathbb{Z})$.
In general, as the finite subgroups of $G$ have bounded order 
$G$ must be finite.
\end{proof}

Theorem 1.3 also follows from the fact that 
$PSL(2,\mathbb{Z})=SL(2,\mathbb{Z})/\langle\pm I\rangle$ is a free product 
$(Z/2Z)*(Z/3Z)$, generated by the images of $A$ and $B$.                        (In fact $\langle A,B\mid A^2=B^3\! ,\medspace A^4=1\rangle$ 
is a presentation for $SL(2,\mathbb{Z})$.)
Moreover, $SL(2,\mathbb{Z})'\cong PSL(2,\mathbb{Z})'$ is freely
generated by the images of
$ABA^{-1}B^{-1}=\left(\begin{smallmatrix}
2&1\\1&1
\end{smallmatrix}\right)$
and $A^{-1}B^{-1}AB=\left(\begin{smallmatrix}
1&1\\1&2
\end{smallmatrix}\right)$,
while the abelianizations are generated by the images of 
$AB=\left(\begin{smallmatrix}
1&0\\ 1&1
\end{smallmatrix}\right)$ \cite[\S6.2]{[Ro]}.

The groups arising as extensions of such groups $G$ by $\mathbb{Z}^2$ are 
the flat $2$-orbifold groups, or $2$-dimensional crystallographic groups.
In three cases $H^2(G;\mathbb{Z}^2)\not=0$, and there are in fact
17 isomorphism classes of such groups.

Let $\Lambda=\mathbb{Z}[t,t^{-1}]$ be the ring of integral Laurent polynomials.
The next theorem is a special case of a classical result of Latimer and
MacDuffee.

\begin{theorem}
There is a $1$-$1$ correspondance between conjugacy classes of matrices in
$GL(n,\mathbb{Z})$ with irreducible characteristic polynomial $\Delta(t)$
and isomorphism classes of ideals in $\Lambda/(\Delta(t))$.
The set of such ideal classes is finite.
\end{theorem}

\begin{proof}
Let $A\in GL(n,\mathbb{Z})$ have characteristic polynomial $\Delta(t)$
and let $R=\Lambda/(\Delta(t))$.
As $\Delta(A)=0$, by the Cayley-Hamilton Theorem, we may define
an $R$-module $M_A$ with underlying abelian group $\mathbb{Z}^n$ by
$t.z=A(z)$ for all $z\in\mathbb{Z}^n$.
As $R$ is a domain and has rank $n$ as an abelian group,
$M_A$ is torsion-free and of rank 1 as an $R$-module, 
and so is isomorphic to an ideal of $R$.
Conversely every $R$-ideal arises in this way.
The isomorphism of abelian groups underlying an $R$-isomorphism 
between two such modules $M_A$ and $M_B$ determines a matrix $C\in
GL(n,\mathbb{Z})$ such that $CA=BC$.
The final assertion follows from the Jordan-Zassenhaus Theorem.
\end{proof}

\section{The Hirsch-Plotkin radical}

The {\it Hirsch-Plotkin radical} $\sqrt G$ of a group $G$ is its maximal locally-nilpotent normal subgroup;
in a virtually poly-$Z$ group every subgroup is finitely generated, and so $\sqrt G$ is then the
maximal nilpotent normal subgroup.                       
If $H$ is normal in $G$ then $\sqrt H$ is normal in $G$ also,
since it is a characteristic subgroup of $H$, 
and in particular it is a subgroup of $\sqrt G$.

For each natural number $q\geq 1$ let $\Gamma_q $ be the group with 
presentation 
\begin{equation*}
\langle x,y,z\mid xz=zx,\medspace yz=zy,\medspace xy=z^q yx\rangle.
\end{equation*}
Every such group $\Gamma_q $ is torsion-free and nilpotent of 
Hirsch length 3. 

\begin{theorem}
Let $G$ be a finitely generated torsion-free nilpotent 
group of Hirsch length $h(G)\leq 4$. Then either
\begin{enumerate}
\item $G$ is free abelian; or

\item $h(G)=3$ and $G\cong\Gamma_q$ for some $q\geq 1$; or

\item $h(G)=4$, $\zeta G\cong\mathbb{Z}^2 $ and 
$G\cong\Gamma_q \times\mathbb{Z}$ for some $q\geq 1$; or

\item $h(G)=4$, $\zeta G\cong\mathbb{Z}$ and $G/\zeta G\cong\Gamma_q $ 
for some $q\geq 1$.
\end{enumerate}

\noindent In the latter case $G$ has characteristic subgroups which 
are free abelian of rank $1$, $2$ and $3$.
In all cases $G$ is an extension of $\mathbb{Z}$ by 
a free abelian normal subgroup.
\end{theorem}

\begin{proof} 
The centre $\zeta G$ is nontrivial and the quotient 
$G/\zeta G$ is again torsion-free nilpotent
\cite[Proposition 5.2.19]{[Ro]}.
We may assume that $G$ is not abelian, and hence that $G/\zeta G$ is not cyclic.
Hence $h(G/\zeta G)\geq 2$, so $h(G)\geq 3$ and $1\leq h(\zeta G)\leq h(G)-2$.
In all cases $\zeta G$ is free abelian.

If $h(G)=3$ then $\zeta G\cong\mathbb{Z}$ and $G/\zeta G\cong\mathbb{Z}^2 $.
On choosing elements $x$ and $y$ representing a basis of $G/\zeta G$ and 
$z$ generating $\zeta G$ we quickly find that $G$ is isomorphic 
to one of the groups $\Gamma_q $, and thus is an extension of 
$\mathbb{Z}$ by $\mathbb{Z}^2 $.

If $h(G)=4$ and $\zeta G\cong\mathbb{Z}^2 $ then $G/\zeta G\cong\mathbb{Z}^2$,
so $G'\leq\zeta G$.
Since $G$ may be generated by elements $x,y,t$ and $u$ where $x$ and $y$ 
represent a basis of $G/\zeta G$ and $t$ and $u$ are central it follows easily 
that $G'$ is infinite cyclic.
Therefore $\zeta G$ is not contained in $G'$ and $G$ has an infinite cyclic 
direct factor.
Hence $G\cong\Gamma_q\times\mathbb{Z} $, for some $q\geq 1$, 
and thus is an extension of $\mathbb{Z}$ by $\mathbb{Z}^3 $. 

The remaining possibility is that $h(G)=4$ and $\zeta G\cong\mathbb{Z}$.
In this case $G/\zeta G$ is torsion-free nilpotent of Hirsch length 3.
If $G/\zeta G$ were abelian $G'$ would also be infinite cyclic, 
and the pairing from $G/\zeta G\times G/\zeta G$ into $G'$ 
defined by the commutator would be nondegenerate and skewsymmetric.
But there are no such pairings on free abelian groups of odd rank.
Therefore $G/\zeta G\cong\Gamma_q $, for some $q\geq 1$.

Let $\zeta_2 G$ be the preimage in $G$ of $\zeta(G/\zeta G)$. Then
$\zeta_2 G\cong\mathbb{Z}^2$ and is a characteristic subgroup of $G$, 
so $C_G (\zeta_2 G)$ is also characteristic in $G$.
The quotient $G/\zeta_2 G$ acts by conjugation on $\zeta_2 G$.
Since $Aut(\mathbb{Z}^2)=GL(2,\mathbb{Z})$ is virtually free and 
$G/\zeta_2 G\cong\Gamma_q /\zeta\Gamma_q \cong\mathbb{Z}^2 $
and since $\zeta_2 G\not=\zeta G$ it follows that $h(C_G (\zeta_2 G))=3$. 
Since $C_G (\zeta_2 G)$ is nilpotent and has centre of rank $\geq 2$ it is abelian,
and so $C_G (\zeta_2 G)\cong\mathbb{Z}^3$. 
The preimage in $G$ of the torsion subgroup of $G/C_G (\zeta_2 G)$ is 
torsion-free, nilpotent of Hirsch length 3 and virtually abelian and hence is abelian. 
Therefore $G/C_G (\zeta_2 G)\cong\mathbb{Z}$. 
\end{proof}

\begin{theorem}
Let $\pi$ be a torsion-free virtually poly-$Z$ group of Hirsch length $4$.
Then $h(\sqrt\pi)\geq 3$.
\end{theorem}
                                                                      
\begin{proof}
Let $S$ be a solvable normal subgroup of finite index in $\pi$.
Then the lowest nontrivial term of the derived series of $S$ is an abelian 
subgroup which is characteristic in $S$ and so normal in $\pi$. 
Hence $\sqrt\pi\not=1$.
If $h(\sqrt\pi)\leq 2$ then $\sqrt\pi\cong\mathbb{Z}$ or $\mathbb{Z}^2$.
Suppose $\pi$ has an infinite cyclic normal subgroup $A$.
On replacing $\pi$ by a normal subgroup $\sigma$ of finite index we may assume 
that $A$ is central and that $\sigma/A$ is poly-$Z$. 
Let $B$ be the preimage in $\sigma$ of a nontrivial abelian normal subgroup of 
$\sigma/A$. 
Then $B$ is nilpotent (since $A$ is central and $B/A$ is abelian) and $h(B)>1$ 
(since $B/A\not=1$ and $\sigma/A$ is torsion-free). 
Hence $h(\sqrt\pi)\geq h(\sqrt\sigma)>1$.        

If $\pi$ has a normal subgroup $N\cong\mathbb{Z}^2$ then 
$Aut(N)\cong{GL(2,\mathbb{Z})}$ 
is virtually free, and so the kernel of 
the natural map from $\pi$ to $Aut(N)$ is nontrivial.
Hence $h(C_\pi (N))\geq 3$.
Since $h(\pi/N)=2$ the quotient $\pi/N$ is virtually abelian, 
and so $C_\pi (N)$ is virtually nilpotent.

In all cases we must have $h(\sqrt\pi)\geq 3$. 
\end{proof}

\section{Amenable groups}

The class of {\it amenable} groups arose first in connection 
with the Banach-Tarski paradox.
A group is amenable if it admits an invariant mean 
for bounded $\mathbb{C}$-valued functions \cite{[Pi]}.
There is a more geometric characterization of finitely presentable 
amenable groups that is more convenient for our purposes. 
Let $X$ be a finite cell-complex with universal cover $\widetilde X$.
Then $\widetilde X$ is an increasing union of finite subcomplexes
$X_j\subseteq X_{j+1}\subseteq \widetilde X=\cup_{n\geq1} X_n$ such that 
$X_j$ is the union of $N_j<\infty$ translates of some fundamental 
domain $D$ for $G=\pi_1(X)$.
Let $N_j'$ be the number of translates of $D$ which meet the frontier of 
$X_j$ in $\widetilde X$.
The sequence $\{ X_j\}$ is a {\it F\o lner exhaustion} for $\widetilde X$ 
if $\lim (N_j'/N_j)=0$,
and $\pi_1(X)$ is amenable if and only if $\widetilde X$ has a F\o lner exhaustion.
This class contains all finite groups and $\mathbb{Z}$, and
is closed under the operations of extension, increasing union,
and under the formation of sub- and quotient groups. 
(However nonabelian free groups are not amenable.)

The subclass $EG$ generated from finite groups and $\mathbb{Z}$ 
by the operations of extension and increasing union is the class 
of {\it elementary amenable} groups.
We may construct this class as follows.
Let $U_0 =1$ and $U_1 $ be the class of finitely generated virtually 
abelian groups.
If $U_\alpha$ has been defined for some ordinal $\alpha$ 
let $U_{\alpha+1}=(\ell U_\alpha)U_1$ and if $U_\alpha $ has been defined 
for all ordinals less than some limit ordinal
$\beta$ let $U_\beta=\cup_{\alpha<\beta} U_\alpha$.
Let $\kappa$ be the first uncountable ordinal.
Then $EG=\ell U_\kappa$.

This class is well adapted to arguments by transfinite induction on 
the ordinal $\alpha(G)=\min\{\alpha|G\in U_\alpha\} $.         
It is closed under extension (in fact 
$U_\alpha U_\beta \subseteq U_{\alpha+\beta}$)
and increasing union, and under the formation of sub- and quotient groups.
As $U_\kappa$ contains every countable elementary amenable group,
$U_\lambda =\ell U_\kappa=EG$ if $\lambda>\kappa$.
Torsion groups in $EG$ are locally finite and elementary amenable free groups
are cyclic.
Every locally-finite by virtually solvable group is elementary amenable;
however this inclusion is proper.

For example, let $\mathbb{Z}^\infty$ be the free abelian group with basis 
$\{ x_i\mid i\in\mathbb{Z}\}$ and let $G$ be the subgroup 
of $Aut(\mathbb{Z}^\infty)$ generated by $\{ e_i\mid i\in\mathbb{Z}\}$,
where $e_i(x_i)=x_i+x_{i+1}$ and $e_i(x_j)=x_j$ if $j\not= i$.
Then $G$ is the increasing union of subgroups isomorphic to groups of upper 
triangular matrices, and so is locally nilpotent. 
However it has no nontrivial abelian normal subgroups.
If we let $\phi$ be the automorphism of $G$ defined by $\phi(e_i)=e_{i+1}$
for all $i$ then $G\rtimes_\phi\mathbb{Z}$ is a finitely generated 
torsion-free elementary amenable group which is not virtually solvable.

It can be shown (using the F\o lner condition) that finitely generated 
groups of subexponential growth are amenable.
The class $SG$ generated from such groups by extensions and 
increasing unions contains $EG$ (since finite groups and finitely generated
abelian groups have polynomial growth), and
is the largest class of groups over which topological surgery 
techniques are known to work in dimension 4 \cite{[FT95]}.
There is a finitely presentable group in $SG$ which is not elementary amenable 
\cite{[Gr98]},
and a finitely presentable amenable group which is not in $SG$ \cite{[BV05]}.

A group is {\it restrained} if it has no noncyclic free subgroup.
Amenable groups are restrained, 
but there are finitely presentable restrained groups which are
not amenable \cite{[OS02],[LM16]}.
There are also infinite finitely generated torsion groups \cite[\S14.2]{[Ro]}. 
These are restrained, but are not elementary amenable.
No known example is also finitely presentable.                   

\section{Hirsch length}

In this section we shall use transfinite induction to extend the notion of 
Hirsch length (as a measure of the size of a solvable group) to elementary 
amenable groups, and to establish the basic properties of this invariant.

\begin{lemma}
Let $G$ be a finitely generated infinite elementary amenable group.
Then $G$ has normal subgroups $K<H$ such that $G/H$ is finite, 
$H/K$ is free abelian of positive rank and the action of $G/H$ 
on $H/K$ by conjugation is effective.
\end{lemma}

\begin{proof}
We may show that $G$ has a normal subgroup $K$ such that $G/K$
is an infinite virtually abelian group, by transfinite induction on $\alpha(G)$.
We may assume that $G/K$ has no nontrivial finite normal subgroup.
If $H$ is a subgroup of $G$ which contains $K$ and is such that $H/K$
is a maximal abelian normal subgroup of $G/K$ then $H$ and $K$ satisfy 
the above conditions. 
\end{proof}

In particular, finitely generated infinite elementary amenable groups are 
virtually indicable.

If $G$ is in $U_1 $ let $h(G)$ be the rank of an abelian subgroup of finite 
index in $G$.
If $h(G)$ has been defined for all $G$ in $U_\alpha$ and $H$ is in 
$\ell U_\alpha$ let 
\begin{equation*}
h(H)=\mathrm{l.u.b.}\{ h(F)|F\leq H,\medspace F\in U_\alpha\}.
\end{equation*}
Finally, if $G$ is in $U_{\alpha+1} $, so has a normal subgroup $H$ in 
$\ell U_\alpha$ with $G/H$ in $U_1 $, let $h(G)=h(H)+h(G/H)$.

\begin{theorem}
Let $G$ be an elementary amenable group. Then
\begin{enumerate}
\item{} $h(G)$ is well defined;

\item{} If $H$ is a subgroup of $G$ then $h(H)\leq h(G)$;

\item{} $h(G)=\mathrm{l.u.b.}\{ h(F)\mid F~is~a~finitely~generated~subgroup~of~G\}$;

\item{} if $H$ is a normal subgroup of $G$ then $h(G)=h(H)+h(G/H)$.
\end{enumerate}
\end{theorem}

\begin{proof}
We shall prove all four assertions simultaneously by induction on $\alpha(G)$.
They are clearly true when $\alpha(G)=1$.
Suppose that they hold for all groups in $U_\alpha$ and that
$\alpha(G)=\alpha+1$.
If $G$ is in $\ell{U}_\alpha$ so is any subgroup, and (1) and (2) are immediate,
while (3) follows since it holds for groups in $U_\alpha$ and since each
finitely generated subgroup of $G$ is in $U_\alpha$.
To prove (4) we may assume that $h(H)$ is finite,
for otherwise $h(G)=h(H)+h(G/H)=\infty$, by (2).
Therefore by (3) there is a finitely generated subgroup $J\leq H$ with
$h(J)=h(H)$.
Given a finitely generated subgroup $Q$ of $G/H$ we may choose a finitely
generated subgroup $F$ of $G$ containing $J$ and whose image in $G/H$ is $Q$.
Since $F$ is finitely generated it is in $U_\alpha$ and so $h(F)=h(H)+h(Q)$.
Taking least upper bounds over all such $Q$ we have $h(G)\geq h(H)+h(G/H)$.
On the other hand if $F$ is any $U_\alpha$-subgroup of $G$ then $h(F)=h(F\cap
H)+h(FH/H)$, since (4) holds for $F$, and so $h(G)\leq h(H)+h(G/H)$.
Thus (4) holds for $G$ also.

Now suppose that $G$ is not in $\ell{U}_\alpha$, but has a normal subgroup 
$K$ in $\ell{U}_\alpha$ such that $G/K$ is in $U_1$.
If $K_1$ is another such subgroup then (4) holds for $K$ and $K_1$ by the
hypothesis of induction and so $h(K)=h(K\cap K_1)+h(KK_1/K)$.
Since we also have $h(G/K)=h(G/KK_1)+h(KK_1/K)$ and
$h(G/K_1)=h(G/KK_1)+h(KK_1/K_1)$ it follows that
$h(K_1)+h(G/K_1)=h(K)+h(G/K)$ and so $h(G)$ is well defined.
Property (2) follows easily, as any subgroup of $G$ is an extension of a 
subgroup of $G/K$ by a subgroup of $K$.
Property (3) holds for $K$ by the hypothesis of induction.
Therefore if $h(K)$ is finite $K$ has a finitely generated subgroup $J$ with
$h(J)=h(K)$.
Since $G/K$ is finitely generated there is a finitely generated subgroup $F$ of
$G$ containing $J$ and such that $FK/K=G/K$.
Clearly $h(F)=h(G)$.
If $h(K)$ is infinite then for every $n\geq0$ there is a finitely generated
subgroup $J_n$ of $K$ with $h(J_n)\geq n$.
In either case, (3) also holds for $G$.
If $H$ is a normal subgroup of $G$ then $H$ and $G/H$ are also in
$U_{\alpha+1}$, while $H\cap K$ and $KH/H=K/H\cap K$ are in 
$\ell{U}_\alpha$ and
$HK/K=H/H\cap K$ and $G/HK$ are in $U_1$.
Therefore 
\begin{equation*}
\begin{split}
h(H)+h(G/H) &=h(H\cap K)+h(HK/K)+h(HK/H)+h(G/HK) \\
&=h(H\cap K)+h(HK/H)+h(HK/K)+h(G/HK).\\
\end{split}
\end{equation*}
Since $K$ is in $\ell{U}_\alpha$ and $G/K$ is in $U_1$ this sum gives
$h(G)=h(K)+h(G/K)$ and so (4) holds for $G$.
This completes the inductive step.
\end{proof}

Let $\Lambda(G)$ be the maximal locally-finite normal subgroup of $G$.

\begin{theorem}
There are functions $d$ and $M$ from 
$\mathbb{Z}_{\geq0}$ to
$\mathbb{Z}_{\geq0}$ 
such that if $G$ is an elementary amenable group of Hirsch
length at most $h$ and $\Lambda(G)$ is its maximal locally finite normal
subgroup then
$G/\Lambda(G)$ has a maximal solvable normal subgroup of derived length 
$\leq{d(h)}$ and index $\leq{M(h)}$.
\end{theorem}

\begin{proof}
We argue by induction on $h$.
Since an elementary amenable group has Hirsch length 0 if and only if it is
locally finite we may set $d(0)=0$ and $M(0)=1$.
Assume that the result is true for all such groups with Hirsch length 
$\leq{h}$ and that $G$ is an elementary amenable group with $h(G)=h+1$.

Suppose first that $G$ is finitely generated. 
Then by Lemma 1.7 there are normal subgroups $K<H$ in $G$ such that $G/H$ is
finite, $H/K$ is free abelian of rank $r\geq1$ and the action of $G/H$ on $H/K$
by conjugation is effective.
(Note that $r=h(G/K)\leq h(G)=h+1$.)
Since the kernel of the natural map from $GL(r,\mathbb{Z})$ to
$GL(r,\mathbb{F}_3)$ is torsion-free, by Lemma 1.2,
we see that $G/H$ embeds in $GL(r,\mathbb{F}_3)$ 
and so has order $\leq3^{r^2}$.
Since $h(K)=h(G)-r\leq h$ the inductive hypothesis applies for $K$, 
so it has a normal subgroup $L$ containing $\Lambda(K)$ 
and of index $\leq{M(h)}$ such that $L/\Lambda(K)$ has derived length 
$\leq{d(h)}$ and is the maximal solvable normal subgroup of $K/\Lambda(K)$.
As $\Lambda(K)$ and $L$ are characteristic in $K$ they are normal in $G$.
(In particular, $\Lambda(K)=K\cap\Lambda(G)$.)
The centralizer of $K/L$ in $H/L$ is a normal solvable subgroup of $G/L$ with
index $\leq[K:L]![G:H]$ and derived length $\leq2$.
Set $M(h+1)=M(h)!3^{(h+1)^2}$ and $d(h+1)=M(h+1)+2+d(h)$.
Then $G.\Lambda(G)$ has a maximal solvable normal subgroup of index 
$\leq{M(h+1)}$ and derived length $\leq{d(h+1)}$ 
(since it contains the preimage of the centralizer of $K/L$ in $H/L$).

In general, let $\{ G_i\mid i\in I\}$ be the set of 
finitely generated subgroups of $G$.
By the above argument $G_i$ has a normal subgroup $H_i$ containing $\Lambda(G_i)$ and such that
$H_i/\Lambda(G_i)$ is a maximal normal solvable subgroup of $G_i/\Lambda(G_i)$
and has derived length $\leq{d(h+1)}$ and index $\leq{M(h+1)}$.
Let $N=\max\{[G_i:H_i]\mid i\in I\}$ and choose $\alpha\in I$ 
such that $[G_\alpha:H_\alpha]=N$.
If $G_i\geq G_\alpha$ then $H_i\cap G_\alpha\leq H_\alpha$.
Since $[G_\alpha:H_\alpha]\leq [G_\alpha:H_i\cap G_\alpha]=
[H_iG_\alpha:H_i]\leq [G_i:H_i]$ we have $[G_i:H_i]=N$ and $H_i\geq H_\alpha$.
It follows easily that if $G_\alpha\leq G_i\leq G_j$ then $H_i\leq H_j$.

Set $J=\{ i\in I\mid H_\alpha\leq H_i\}$ and $H=\cup_{i\in J}H_i$.
If $x,y\in H$ and $g\in G$ then there are indices $i,k$ and $k\in J$ such that 
$x\in H_i$, $y\in H_j$ and $g\in G_k$.
Choose $l\in J$ such that $G_l$ contains $G_i\cup G_j\cup G_k$.
Then $xy^{-1}$ and $gxg^{-1}$ are in $H_l\leq H$, and so $H$ is a normal
subgroup of $G$.
Moreover if $x_1,\dots,x_N$ is a set of coset representatives for $H_\alpha$ in
$G_\alpha$ then it remains a set of coset representatives for $H$ in $G$, 
and so $[G:H]=N$.

Let $D_i$ be the $d(h+1)$th derived subgroup of $H_i$.
Then $D_i$ is a locally-finite normal subgroup of $G_i$ and so,
by an argument similar to that of the above paragraph $\cup_{i\in J}D_i$ is a
locally-finite normal subgroup of $G$.
Since the $d(h+1)$th derived subgroup of $H$ is contained in 
$\cup_{i\in J}D_i$ (as each iterated commutator involves only 
finitely many elements of $H$) it follows that 
$H\Lambda(G)/\Lambda(G)\cong H/H\cap\Lambda(G)$ is solvable and of derived
length $\leq{d(h+1)}$.
\end{proof}

The above result is from \cite{[HL92]}.
The argument can be simplified to some extent if $G$ is countable and
torsion-free.
(In fact a virtually solvable group of finite Hirsch length and with no 
nontrivial locally-finite normal subgroup must be countable 
\cite[Lemma 7.9]{[Bi]}.)

\begin{lemma}
Let $G$ be an elementary amenable group. If $h(G)=\infty$ 
then for every $k>0 $ there is a subgroup $H$ of $G$ with $k<h(H)<\infty$.
\end{lemma}
                                                            
\begin{proof}
We shall argue by induction on $\alpha(G)$. 
The result is vacuously true if $\alpha (G)=1$.
Suppose that it is true for all groups in $U_\alpha$ and $G$ is in 
$\ell U_\alpha$.
Since $h(G)=$ l.u.b.$\{ h(F)|F\leq G,\medspace F\in U_\alpha\}$ either there 
is a subgroup $F$ of $G$ in $U_\alpha$ with $h(F)=\infty$, in which case the 
result is true by the inductive hypothesis, or $h(G)$ is 
the least upper bound of a set of natural numbers and the result is true.
If $G$ is in $U_{\alpha+1} $ then it has a normal subgroup $N$ which is in 
$\ell U_\alpha$ with quotient $G/N$ in $U_1 $. 
But then $h(N)=h(G)=\infty$ and so $N$ has such a subgroup. 
\end{proof}

\begin{theorem}
Let $G$ be an elementary amenable group of finite cohomological dimension. 
Then $h(G)\leq c.d.G$ and $G$ is virtually solvable.
\end{theorem}
                                                      
\begin{proof}
Since $c.d.G<\infty$ the group $G$ is torsion-free.
Let $H$ be a subgroup of finite Hirsch length.
Then $H$ is virtually solvable and $c.d.H\leq c.d.G$ so $h(H)\leq c.d.G$.
The theorem now follows from Theorem 1.9 and Lemma 1.10.
\end{proof}

\section{Modules and finiteness conditions}

Let $G$ be a group and $w:G\to\mathbb{Z}^\times$ a homomorphism,
and let $R$ be a commutative ring.
Then $\bar g=w(g)g^{-1} $ defines an anti-involution on $R[G]$.
If $L$ is a left $R[G]$-module $\overline{L}$ shall denote the {\it conjugate} right $R[G]$-module 
with the same underlying $R$-module and $R[G]$-action                    
given by $l.g=\bar g.l$, for all $l\in L$ and $g\in G$.
(We shall also use the overline to denote the conjugate of a right $R[G]$-module.)
The conjugate of a free left (right) module is a free right (left) module 
of the same rank.

Let $\mathbb{Z}^w $ denote the $G$-module 
with underlying abelian group $\mathbb{Z}$ and $G$-action given by 
$g.n=w(g)n$ for all $g$ in $G$ and $n$ in $\mathbb{Z}$. 
                                                             
\begin{lemma}
{\rm[Wl65]}\qua 
Let $G$ and $H$ be groups such that $G$ is finitely presentable 
and there are homomorphisms $j:H\to G$ and $\rho:G\to H$ with $\rho j=id_H $.
Then $H$ is also finitely presentable.
\end{lemma}
                      
\begin{proof}
Since $G$ is finitely presentable there is an epimorphism $p:F\to G$
from a free group $F(X)$ with a finite basis $X$ onto $G$, with kernel the 
normal closure of a finite set of relators $R$.
We may choose elements $w_x $ in $F(X)$ such that $j\rho p(x)=p(w_x )$, for all $x$ in $X$.
Then $\rho$ factors through the group $K$ with presentation 
$\langle X\mid R, x^{-1} w_x ,\forall x\in X\rangle$, say $\rho=vu$.
Now $uj$ is clearly onto, while $vuj=\rho j=id_H$, and so $v$ and $uj$ are mutually inverse
isomomorphisms. Therefore $H\cong K$ is finitely presentable. 
\end{proof}

A group $G$ is $FP_n$ if the augmentation $\mathbb{Z}[G]$-module $\mathbb{Z}$ 
has a projective resolution
which is finitely generated in degrees $\leq n$, and it is $FP$ if it has finite cohomological
dimension and is $FP_n $ for $n=c.d.G$. 
It is $FF$ if moreover $\mathbb{Z}$ has a finite resolution
consisting of finitely generated free $\mathbb{Z}[G]$-modules.
``Finitely generated" is equivalent to $FP_1 $, while ``finitely presentable" implies $FP_2 $. 
Groups which are $FP_2 $ are also said to be {\it almost finitely presentable}.
(There are $FP$ groups which are not finitely presentable \cite{[BB97]}.)
An elementary amenable group $G$ is $FP_\infty$ if and only if 
it is virtually $FP$, 
and is then virtually constructible and solvable of finite Hirsch length 
\cite{[Kr93]}.

If the augmentation $\mathbb{Q}[\pi]$-module $\mathbb{Q}$ has a finite resolution $F_* $ by
finitely generated projective modules then 
$\chi(\pi)=\Sigma (-1)^i dim_{\mathbb{Q}}(\mathbb{Q}\otimes_\pi F_i)$ is independent 
of the resolution.
(If $\pi$ is the fundamental group of an aspherical finite complex $K$ then
$\chi(\pi)=\chi(K)$.) 
We may extend this definition to groups $\sigma$ which have a subgroup $\pi$ 
of finite index with such a resolution by setting 
$\chi(\sigma)=\chi(\pi)/[\sigma:\pi]$.
(It is not hard to see that this is well defined.)

Let $P$ be a finitely generated projective $\mathbb{Z}[\pi]$-module.
Then $P$ is a direct summand of $\mathbb{Z}[\pi]^r$, for some $r\geq0$,
and so is the image of some idempotent $r\times r$-matrix $M$ with entries 
in $\mathbb{Z}[\pi]$.
The {\it Kaplansky rank} $\kappa(P)$ is the coefficient of $1\in\pi$ 
in the trace of $M$.
It depends only on $P$ and is strictly positive if $P\not=0$.
The group $\pi$ satisfies the {\it Weak Bass Conjecture} if 
$\kappa(P)= dim_{\mathbb{Q}}\mathbb{Q}\otimes_\pi P$.
There is also a {\it Strong\/} Bass Conjecture, 
which we shall not formulate here, although it is invoked in Theorem 3.4.
Both conjectures have been confirmed for linear groups, 
residually finite groups, solvable groups, 
groups of cohomological dimension $\leq2$ over $\mathbb{Q}$
and $PD_3$-groups. 
(See \cite{[Ec01]} for further details.)

The following result from \cite{[BS78]} shall be useful.

\begin{theorem}
[Bieri-Strebel] 
Let $G$ be an $FP_2$ group with $G/G'$ infinite.
Then $G$ is an HNN extension with finitely generated base 
and associated subgroups.
\end{theorem}

\begin{proof}
(Sketch -- We shall assume that $G$ is finitely presentable.)
Let $h:F(m)\to G$ be an epimorphism, and let $g_i=h(x_i)$ for $1\leq i\leq m$.
We may assume that $g_m$ has infinite order 
modulo  $\langle\langle g_i\mid 1\leq i<m\rangle\rangle$.
Since $G$ is finitely presentable the kernel of $h$ is the normal closure 
of finitely many relators, of weight 0 in the letter $x_m$.
Each such relator is a product of powers of conjugates of the generators 
$\{ x_i\mid 1\leq i<m\}$ by powers of $x_m$.
Thus we may assume the relators are contained in 
$\langle{ x_m^jx_ix_m^{-j}}\mid 1\leq i<m,\, -p\leq j\leq p\rangle$, 
for some sufficiently large $p$.
Let $U= \langle{ g_m^jg_ig_m^{-j}}\mid 1\leq i<m,\, -p\leq j<p\rangle$,
and let $V=g_mUg_m^{-1}$.
Let $B$ be the subgroup of $G$ generated by $U\cup V$ and let $\tilde G$ 
be the HNN extension with base $B$ and associated subgroups $U$ and $V$ 
presented by 
$\tilde G=\langle B,s\mid sus^{-1}=\tau(u)\medspace \forall u\in U\rangle$,
where $\tau:U\to V$ is the isomorphism determined by conjugation by $g_m$ 
in $G$.
There are obvious epimorphisms $\xi:F(m)\to\tilde G$ and 
$\psi:\tilde G\to G$ with composite $h$. 
It is easy to see that $\mathrm{Ker}(h)\leq\mathrm{Ker}(\xi)$ and so 
$\tilde G\cong G$.
\end{proof}

An HNN extension is restrained if and only if 
it is ascending and the base is restrained.

A ring $R$ is {\it weakly finite} if every onto endomorphism of $R^n$
is an isomorphism, for all $n\geq0$.
(In \cite{[H2]} the term ``SIBN ring" was used instead.)
Finitely generated stably free modules over weakly finite rings have 
well defined ranks, and the rank is strictly positive if the module is nonzero.
Skew fields are weakly finite, as are subrings of weakly finite rings.
If $G$ is a group its complex group algebra $\mathbb{C}[G]$ is weakly finite,
by a result of Kaplansky.
(See \cite{[Ro84]} for a proof.)

A ring $R$ is {\it (regular) coherent} if every finitely presentable left
$R$-module has a (finite) resolution by finitely generated projective $R$-modules,
and is {\it (regular) noetherian} if moreover every finitely generated $R$-module is finitely presentable.
A group $G$ is regular coherent or regular noetherian if $R[G]$ is 
regular coherent or regular noetherian (respectively) for any regular noetherian ring $R$.                
It is coherent as a {\it group} if its finitely generated subgroups 
are finitely presentable.

\begin{lemma}
If $G$ is a group such that $\mathbb{Z}[G]$ is coherent then every finitely
generated subgroup of $G$ is $FP_\infty$.
\end{lemma}

\begin{proof} 
Let $H$ be a subgroup of $G$.               
Since $\mathbb{Z}[H]\leq \mathbb{Z}[G]$ is a faithfully flat ring extension
a left $\mathbb{Z}[H]$-module is finitely generated over $\mathbb{Z}[H]$ 
if and only if the induced module $\mathbb{Z}[G]\otimes_HM$ 
is finitely generated over $\mathbb{Z}[G]$.
Hence $M$ is $FP_n$ over $\mathbb{Z}[H]$ if and only if 
$\mathbb{Z}[G]\otimes_HM$ is $FP_n$ over $\mathbb{Z}[G]$,
by induction on $n$.

If $H$ is finitely generated then the augmentation $\mathbb{Z}[H]$-module 
$\mathbb{Z}$ is finitely presentable over $\mathbb{Z}[H]$. 
Hence $\mathbb{Z}[G]\otimes_H\mathbb{Z}$ is finitely presentable 
over $\mathbb{Z}[G]$, and so is $FP_\infty$ over $\mathbb{Z}[G]$, 
since that ring is coherent.
Hence $\mathbb{Z}$ is $FP_\infty$ over $\mathbb{Z}[H]$, 
i.e., $H$ is $FP_\infty$.
\end{proof}

Thus if either $G$ is coherent (as a group) or $\mathbb{Z}[G]$ is
coherent (as a ring) every finitely generated subgroup of $G$ is $FP_2$.
As the latter condition shall usually suffice for our purposes below,
we shall say that such a group is {\it almost coherent}. 
The connection between these notions has not been much studied.

The class of groups whose integral group ring is regular coherent contains 
them trivial group and is closed under generalised free products and 
HNN extensions with amalgamation over subgroups whose group rings are 
regular noetherian \cite[Theorem 19.1]{[Wd78]}.
If $[G:H]$ is finite and $G$ is torsion-free then 
$\mathbb{Z}[G]$ is regular coherent if and only if $\mathbb{Z}[H]$ is.
In particular, free groups and surface groups are coherent and their integral group rings
are regular coherent, while (torsion-free) virtually poly-$Z$ groups are coherent and their
integral group rings are (regular) noetherian.

\section{Ends and cohomology with free coefficients}

A finitely generated group $G$ has 0, 1, 2 or infinitely many ends.
It has 0 ends if and only if it is finite, 
in which case $H^0 (G;\mathbb{Z}[G])\cong\mathbb{Z}$ and 
$H^q (G;\mathbb{Z}[G])=0$ for $q>0$.                                     
Otherwise $H^0(G;\mathbb{Z}[G])=0$ and $H^1(G;\mathbb{Z}[G])$ 
is a free abelian group of rank $e(G)-1$, 
where $e(G)$ is the number of ends of $G$ \cite{[Sp49]}.
The group $G$ has more than one end if and only if it is a nontrivial 
generalised free product with amalgamation $G\cong A*_C B$ 
or an HNN extension $A*_C \phi$, where $C$ is a finite group.
In particular, it has two ends if and only if it is virtually $\mathbb{Z}$ 
if and only if it has a (maximal) finite normal subgroup $F$ such that 
$G/F\cong\mathbb{Z}$ or $D$, 
where $D=(Z/2Z)*(Z/2Z)$ is the {\it infinite dihedral
group\/} \cite{[St]} - see also \cite{[DD]}.

If $G$ is a group with a normal subgroup $N$, and $A$ is a left 
$\mathbb{Z}[G]$-module there is a {\it Lyndon-Hochschild-Serre spectral sequence} 
(LHSSS) for $G$ as an extension of $G/N$ by $N$ 
and with coefficients $A$:
\begin{equation*}
E_2 =H^p (G/N; H^q (N;A))\Rightarrow H^{p+q}(G;A),
\end{equation*}
the $r^{th}$ differential having bidegree $(r,1-r)$.
(See \cite[\S10.1]{[Mc]}.)
In several places below, when considering such spectral sequences
(e.g., in Theorems 2.12 and 8.1), 
we use without comment the fact that if
$M$ is a left $\mathbb{Z}[G]$-module and $M|_1$ is the underlying 
abelian group then $M\otimes\mathbb{Z}[G]$ (with the diagonal $G$-action) 
is canonically isomorphic to the induced module $M|_1\otimes\mathbb{Z}[G]$.
	                                        
\begin{theorem}
{\rm[Ro75]}\qua
If $G$ has a normal subgroup $N$ which is the union of an increasing sequence 
of subgroups $N_n $ such that $H^s (N_n ;\mathbb{Z}[G])=0$ for $s\leq r$ then
$H^s (G;\mathbb{Z}[G])=0$ for $s\leq r$.
\end{theorem}

\begin{proof}
Let $s\leq r$.
Let $f$ be an $s$-cocycle for $N$ with coefficients $\mathbb{Z}[G]$,
and let $f_n $ denote the restriction of $f$ to a cocycle on $N_n $.
Then there is an ${(s-1)}$-cochain $g_n $ on $N_n$ such that $\delta g_n =f_n$.
Since $\delta(g_{n+1} |_{N_n} -g_n )=0$ and $H^{s-1} (N_n ;\mathbb{Z}[G])=0$ there is
an $(s-2)$-cochain $h_n$ on $N_n$ with $\delta h_n =g_{n+1} |_{N_n} -g_n $.
Choose an extension $h_n '$ of $h_n $ to $N_{n+1} $ and let 
$\hat g_{n+1} =g_{n+1} -\delta h_n '$.
Then $\hat g_{n+1} |_{N_n} =g_n$ and $\delta\hat g_{n+1} =f_{n+1} $.
In this way we may extend $g_0$ to an $(s-1)$-cochain $g$ on $N$ 
such that $f=\delta g$ and so $H^s (N;\mathbb{Z}[G])=0$.
The LHSSS for $G$ as an extension of $G/N$ by $N$, 
with coefficients $\mathbb{Z}[G]$,
now gives $H^s (G;\mathbb{Z}[G])=0$ for $s\leq r$. 
\end{proof}

\begin{cor}
The hypotheses are satisfied if $N$ is the union 
of an increasing sequence of $FP_r$ subgroups $N_n $ such that 
$H^s (N_n ;\mathbb{Z}[N_n])=0$ for $s\leq r$.
In particular, if $N$ is the union of an increasing
sequence of finitely generated, one-ended subgroups then $G$ has one end.
\end{cor}

\begin{proof} 
We have $H^s (N_n;\mathbb{Z}[G])=
H^s (N_n ;\mathbb{Z}[N_n ])\otimes\mathbb{Z}[G/N_n]=0$,
for all $s\leq r$ and all $n$, since $N_n$ is $FP_r$.
\end{proof}

If the successive inclusions are finite this corollary may be sharpened further.

\begin{thm}
[Gildenhuys-Strebel] Let $G=\cup_{n\geq1}G_n$ be the union of an increasing 
sequence of $FP_r$ subgroups.
Suppose that $[G_{n+1}:G_n]<\infty$ and $H^s(G_n;\mathbb{Z}[G_n])=0$ 
for all $s<r$ and $n\geq1$.
If $G$ is not finitely generated then $H^s(G;F)=0$ for every free 
$\mathbb{Z}[G]$-module $F$ and all $s\leq{r}$. 
\qed
\end{thm}

The enunciation of this theorem in \cite{[GS81]} assumes also that 
$c.d.G_n=r$ for all $n\geq1$, 
and concludes that $c.d.G=r$ if and only if 
$G$ is finitely generated.
However the argument establishes the above assertion.

\begin{theorem}
Let $G$ be a finitely generated group with an infinite restrained normal
subgroup $N$ of infinite index. Then $e(G)=1$.
\end{theorem}

\begin{proof}
Since $N$ is infinite 
$H^1(G;\mathbb{Z}[G])\cong{H^0}(G/N;H^1(N;\mathbb{Z}[G]))$,
by the LHSSS.
If $N$ is finitely generated then
$H^1(N;\mathbb{Z}[G])\cong{H^1}(N;\mathbb{Z}[N])\otimes\mathbb{Z}[G/N]$,
with the diagonal $G/N$-action.
Since $G/N$ is infinite $H^1(G;\mathbb{Z}[G])=0$.
If $N$ is locally one-ended or locally virtually $\mathbb{Z}$ 
and not finitely generated then $H^1(N;\mathbb{Z}[G])=0$, 
by Theorem 1.15 and the Gildenhuys-Strebel Theorem,
respectively.
In all of these cases $e(G)=1$.

There remains the possibility that $N$ is locally finite.
If $e(G)>1$ then $G\cong A*_C B$ or $A*_C \phi$ with $C$ finite, 
by Stallings' characterization of such groups.
Suppose $G\cong A*_C B$. 
Since $N$ is infinite there is an $n\in{N\setminus{C}}$.
We may suppose that $n\in{A}$,
since elements of finite order in $A*_CB$ are conjugate to elements of
$A$ or $B$ \cite[Theorem 6.4.3]{[Ro]}.
Let $g\in{B\setminus{C}}$, and let $n'=gng^{-1}$.
Since $N$ is normal $nn'\in{N}$ also.
But $nn'=ngng^{-1}$ has infinite order in $G$, 
by the ``uniqueness of normal form" for such groups.
This contradicts the fact that $N$ is locally finite.
A similar argument shows that $G$ cannot be $A*_C \phi$.
Thus $G$ must have one end.
\end{proof}

In particular, a countable restrained group $N$ is either 
elementary amenable and $h(N)\leq1$
or is an increasing union of finitely generated, 
one-ended subgroups.
    
The second cohomology of a group with free coefficients ($H^2 (G;R[G])$, 
$R=\mathbb{Z}$ or a field) shall play an important role in our investigations.

\begin{thm}
[Farrell] Let $G$ be a finitely presentable group.
If $G$ has an element of infinite order and $R=\mathbb{Z}$ or is a 
field then $H^2 (G;R[G])$ is either 0 or $R$ or is not finitely generated.  
\qed
\end{thm}

Farrell also showed in \cite{[Fa74]} that if 
$H^2 (G;\mathbb{F}_2[G])\cong Z/2Z$ then 
every finitely generated subgroup of $G$ with one end has finite index in $G$.
Hence if $G$ is also torsion-free then subgroups of infinite index 
in $G$ are locally free.
Bowditch has since shown that such groups are virtually the
fundamental groups of aspherical closed surfaces 
(\cite{[Bo04]} - see \S8 below).

We would also like to know when $H^2(G;\mathbb{Z}[G])$ is 0
(for $G$ finitely presentable).
In particular, we expect this to be so if $G$ 
has an elementary amenable, normal subgroup $E$ such that
either $h(E)=1$ and $G/E$ has one end or $h(E)=2$ and $[G:E]=\infty$ 
or $h(E)\geq3$, 
or if $G$ is an ascending HNN extension over a finitely generated, 
one-ended base.
Our present arguments for these two cases require stronger finiteness 
hypotheses, and each use the following result of \cite{[BG85]}.

\begin{thm}
[Brown-Geoghegan] Let $G$ be an HNN extension $B*_\phi$ in which the base 
$B$ and associated subgroups $I$ and $\phi(I)$ are $FP_n$. 
If the homomorphism from $H^q(B;\mathbb{Z}[G])$ to $H^q(I;\mathbb{Z}[G])$
induced by restriction
is injective for some $q\leq n$ then the corresponding homomorphism 
in the Mayer-Vietoris sequence is injective,
so $H^q(G;\mathbb{Z}[G])$ is a quotient of $H^{q-1}(I;\mathbb{Z}[G])$.
\qed
\end{thm}

We begin with the case of ``large" elementary amenable normal subgroups. 

\begin{theorem}
Let $G$ be a finitely presentable group with a  
restrained normal subgroup $E$ of infinite index.
Suppose that either $E$ is abelian of rank $1$ and $G/E$ has one end 
or $E$ is torsion-free, elementary amenable and $h(E)>1$ or
$E$ is almost coherent, locally virtually indicable, 
and has a finitely generated, one-ended subgroup.
Then $H^s(G;\mathbb{Z}[G])=0$ for $s\leq2$.
\end{theorem}

\begin{proof}
If $E$ is abelian of positive rank and $G/E$ has one end then 
$G$ is 1-connected at $\infty$ by Theorem 1 of \cite{[Mi87]}, 
and so $H^s(G;\mathbb{Z}[G])=0$ for $s\leq2$, by \cite{[GM86]}.

Suppose next that $E$ is torsion-free, elementary amenable and $h(E)>1$.
If $E$ is virtually solvable it has a nontrivial characteristic 
abelian subgroup $A$.
If $A$ has rank 1 then $E/A$ is infinite, 
so $G/A$ has one end, by Theorem 1.16,
and then $H^s(G;\mathbb{Z}[G])=0$ for $s\leq2$, as before.
If $A\cong\mathbb{Z}^2$ then $H^2(A;\mathbb{Z}[G]))\cong\mathbb{Z}[G/A]$.
Otherwise, $A$ has $\mathbb{Z}^2$ as a subgroup of infinite index
and so $H^2(A;\mathbb{Z}[G])=0$.
If $E$ is not virtually solvable $H^s(E;\mathbb{Z}[G])=0$ for all $s$
\cite[Proposition 3]{[Kr93']}.
(The argument applies even if $E$ is not finitely generated.)
In all cases, an LHSSS argument gives $H^s(G;\mathbb{Z}[G])=0$ for $s\leq2$.

We may assume henceforth that $E$ is almost coherent and is an increasing 
union of finitely generated one-ended subgroups 
$E_n\subseteq E_{n+1}\dots\subseteq E=\cup E_n$.
Since $E$ is locally virtually indicable there are subgroups $F_n\leq E_n$ 
such that $[E_n:F_n]<\infty$ and which map onto $\mathbb{Z}$. 
Since $E$ is almost coherent these subgroups are $FP_2$.
Hence they are HNN extensions over $FP_2$ bases $H_n$,
by Theorem 1.13, and the extensions are ascending, since $E$ is restrained. 
Since $E_n$ has one end $H_n$ is infinite and so has one or two ends.

Suppose that $H_n$ has two ends, for all $n\geq1$.
Then $E_n$ is elementary amenable, $h(E_n)=2$ and $[E_{n+1}:E_n]<\infty$,
for all $n\geq1$.
Hence $E$ is elementary amenable and $h(E)=2$.
If $E$ is finitely generated it is $FP_2$ and so                                  
$H^s(G;\mathbb{Z}[G])=0$ for $s\leq2$, by an LHSSS argument.
This is also the case if $E$ is not finitely generated, for then 
$H^s(E;\mathbb{Z}[G])=0$ for $s\leq2$, by the Gildenhuys-Strebel Theorem,
and we may again apply an LHSSS argument.

Otherwise we may assume that $H_n$ has one end, for all $n\geq1$.
In this case $H^s(F_n;\mathbb{Z}[F_n])=0$ for $s\leq2$, 
by the Brown-Geoghegan Theorem.
Therefore $H^s(G;\mathbb{Z}[G])=0$ for $s\leq2$, by Theorem 1.15.
\end{proof}

The theorem applies if $E$ is almost coherent and elementary amenable,
since elementary amenable groups are restrained and locally virtually indicable.             
It also applies if $E=\sqrt G$ is large enough, 
since finitely generated nilpotent groups are virtually poly-$Z$.
Similar arguments show that if $h(\sqrt G)\geq r$ then 
$H^s (G;\mathbb{Z}[G])=0$ for $s<r$, 
and if also $[G:\sqrt G]=\infty$ then $H^r(G;\mathbb{Z}[G])=0$.                 

Are the hypotheses that $E$ be almost coherent 
and locally virtually indicable necessary?                     
Is it sufficient that $E$ be an increasing union 
of finitely generated, one-ended subgroups?

\begin{theorem}
Let $G=B*_\phi$ be an HNN extension with 
$FP_2$ base $B$ and associated subgroups $I$ and $\phi(I)=J$, and which has a 
restrained normal subgroup $N\leq \langle\langle B\rangle\rangle$.
Then $H^s(G;\mathbb{Z}[G])=0$ for $s\leq2$ if either
\begin{enumerate}
\item the HNN extension is ascending and $B=I\cong J$ has one end; or

\item $N$ is locally virtually $\mathbb{Z}$ and $G/N$ has one end; or

\item $N$ has a finitely generated subgroup with one end.
\end{enumerate}
\end{theorem}

\begin{proof}
The first assertion follows immediately from the Brown-Geogeghan Theorem.

Let $t$ be the stable letter, so that $tit^{-1}=\phi(i)$, for all $i\in I$.
Suppose that $N\cap J\not= N\cap B$, and let $b\in N\cap B\setminus{J}$.
Then $b^t=t^{-1}bt$ is in $N$, since $N$ is normal in $G$. 
Let $a$ be any element of $N\cap B$. 
Since $N$ has no noncyclic free subgroup there is a word $w\in F(2)$
such that $w(a,b^t)=1$ in $G$.
It follows from Britton's Lemma that $a$ must be in $I$, 
and so $N\cap B=N\cap I$. 
In particular, $N$ is the increasing union of copies of $N\cap B$.
                                                     
Hence $G/N$ is an HNN extension with base $B/N\cap B$ and associated subgroups
$I/N\cap I$ and $J/N\cap J$.
Therefore if $G/N$ has one end the latter groups are infinite,
and so $B$, $I$ and $J$ each have one end.
If $N$ is virtually $\mathbb{Z}$ then $H^s(G;\mathbb{Z}[G])=0$ for $s\leq2$, 
by an LHSSS argument.
If $N$ is locally virtually $\mathbb{Z}$ but is not finitely generated 
then it is the increasing union of a sequence of two-ended subgroups 
and $H^s(N;\mathbb{Z}[G])=0$ for $s\leq1$,
by the Gildenhuys-Strebel Theorem.
Since $H^2(B;\mathbb{Z}[G])\cong H^0(B;H^2(N\cap B;\mathbb{Z}[G]))$ 
and $H^2(I;\mathbb{Z}[G])\cong H^0(I;H^2(N\cap I;\mathbb{Z}[G]))$,
the restriction map from $H^2(B;\mathbb{Z}[G])$ to $H^2(I;\mathbb{Z}[G])$ 
is injective.
If $N$ has a finitely generated, one-ended subgroup $N_1$,
we may assume that $N_1\leq N\cap B$, and so $B$, $I$ and $J$ also have one end.
Moreover $H^s(N\cap B;\mathbb{Z}[G])=0$ for $s\leq 1$, by Theorem 1.15.
We again see that the restriction map from $H^2(B;\mathbb{Z}[G])$ to $H^2(I;\mathbb{Z}[G])$
is injective.
The result now follows in these cases from the Brown-Geoghegan Theorem.
\end{proof}

The final result of this section is Theorem 8.8 of \cite{[Bi]}.

\begin{thm}
[Bieri] Let $G$ be a nonabelian group with $c.d.G=n$.
Then $c.d.\zeta G\leq n-1$, and if $\zeta G$ has rank $n-1$ then $G'$ is free.
\qed
\end{thm}

\section{Poincar\'e duality groups}

A group $G$ is a $PD_n$-{\it group} if it is $FP$, 
$H^p(G;\mathbb{Z}[G])=Ext^p_{\mathbb{Z}[G]}(\mathbb{Z},\mathbb{Z}[G])$ 
is $0$ for $p\not=n$ and $H^n(G;\mathbb{Z}[G])$ is infinite cyclic.
The ``dualizing module" $H^n(G;\mathbb{Z}[G])$
is a right $\mathbb{Z}[G]$-module, 
with $G$-action determined by a homomorphism
$w=w_1(G):G\to\mathbb{Z}^\times$.
The group is {\it orientable} (or is a $PD^+_n$-{\it group}) if $w$ is trivial,
i.e., if $H^n(G;\mathbb{Z}[G])$ is isomorphic to the augmentation module $\mathbb{Z}$. 
(See \cite{[Bi]}.)

The only $PD_1$-group is $\mathbb{Z}$.
Eckmann, Linnell and M\"uller showed that every $PD_2 $-group is the
fundamental group of an aspherical closed surface.
(See Chapter VI of \cite{[DD]}.)
Bowditch has since found a much stronger result, 
which must be close to the
optimal characterization of such groups \cite{[Bo04]}.

\begin{thm}
[Bowditch] 
Let $G$ be an $FP_2$ group and $F$ a field. 
Then $G$ is virtually the fundamental group of an aspherical closed surface
if and only if $H^2(G;F[G])$ has a $1$-dimensional $G$-invariant subspace.
\qed
\end{thm}

In particular, this theorem applies if $H^2(G;\mathbb{Z}[G])$
is infinite cyclic,
for then its image in $H^2(G;\mathbb{F}_2[G])$
under reduction {\it mod} (2) is such a subspace.

The following result corresponds to the 
fact that an infinite covering space of
a PL $n$-manifold is homotopy equivalent to a complex of dimension $<n$
\cite{[St77]}.

\begin{thm}
[Strebel] 
Let $H$ be a subgroup of infinite index in a $PD_n $-group $G$. 
Then $c.d.H<n$.
\qed
\end{thm}

Let $S$ be a ring.
If $C$ is a left $S$-module and $R$ is a subring of $S$ let $C|_R$ be the left
$R$-module underlying $C$.
If $A$ is a left $R$-module the abelian group $Hom_R(S|_R,A)$ has a natural 
left $S$-module structure given by $((sf)(s')=f(s's)$ for all 
$f\in{Hom}_R (S|_R ,A)$ and $s,s'\in{S}$.
The groups $Hom_R (C|_R,A)$ and $Hom_S (C,Hom_R (S|_R ,A))$ are naturally 
isomorphic, 
for the maps $I$ and $J$ defined by $I(f)(c)(s)=f(sc)$ and 
$J(\theta)(c)=\theta(c)(1)$ for $f:C\to A$ and $\theta:C\to Hom_R (S,A)$ 
are mutually inverse isomorphisms.
When $K$ is a subgroup of $\pi$, $R=\mathbb{Z}[K]$ and $S=\mathbb{Z}[\pi]$ 
we may write $C|_K$ for $C|_R$, and
the module $Hom_{\mathbb{Z}[K]}(\mathbb{Z}[\pi]|_K,A)$ is said to 
be {\it coinduced from} $A$.
The above isomorphisms give rise to Shapiro's Lemma.
In our applications $\pi/K$ shall usually be infinite cyclic and $S$ is then 
a twisted Laurent extension of $R$.

If $G$ is a group and $A$ is a left $\mathbb{Z}[G]$-module let $A|_1$ be the
$\mathbb{Z}[G]$-module with the same underlying group and trivial $G$-action,
and let $A^G=Hom_\mathbb{Z}(\mathbb{Z}[G],A)$ be the module of functions 
$\alpha:G\to{A}$ with $G$-action given by $(g\alpha)(h)=g.\alpha(hg)$ 
for all $g,h\in{G}$.
Then ${A|_1}^G$ is coinduced from a module over the trivial group.

\begin{theorem}
Let $\pi$ be a $PD_n$-group with a
normal subgroup $K$ such that $\pi/K$ is a $PD_r$-group.
Then $K$ is a $PD_{n-r}$-group if and only if it is $FP_{[n/2]}$.
\end{theorem}

\begin{proof}
The condition is clearly necessary.
Assume that it holds, and let $C_*$ be a resolution of $\mathbb{Z}$ 
by free left $\mathbb{Z}[K]$-modules 
which are finitely generated in degrees $\leq\frac{n}2$.
After passing to a subgroup of index 2, if necessary,
we may assume that $G=\pi/K$ is orientable.
It is sufficient to show that the functors $H^s (K;-)$ 
from left $\mathbb{Z}[K]$-modules to abelian groups
commute with direct limit, for all $s\leq{n}$,
for then $K$ is $FP_{n-1}$ \cite{[Br75]}, 
and the result follows from \cite[Theorem 9.11]{[Bi]} 
(and an LHSSS corner argument to identify the dualizing module).
Since $K$ is $FP_{[n/2]}$ we may assume $s>n/2$.
If $A$ is a left $\mathbb{Z}[K]$-module 
and $W=Hom_{\mathbb{Z}[K]}(\mathbb{Z}[\pi],A)$
then $H^s(K;A)\cong{H^s}(\pi;W)\cong{H_{n-s}}(\pi;\overline { W})$, 
by Shapiro's Lemma and Poincar\'e duality.
Let $\sigma:G\to\pi$ be a (setwise) section of the projection.

Let $A_g$ be the left $\mathbb{Z}[K]$-module with the same underlying group 
as $A$ and $K$-action given by $k.a=\sigma(g)k\sigma(g)^{-1}a$ for all 
$a\in{A}$, $g\in{G}$ and $k\in{K}$.
The $\mathbb{Z}[K]$-epimorphisms $p_g:W\to{A_g}$ given by
$p_g(f)=f(\sigma(g))$ for all $f\in{W}$ and $g\in{G}$
determine an isomorphism $W\cong\Pi_{g\in{G}}A_g$.
Hence they induce $\mathbb{Z}$-linear isomorphisms 
$H_q(K;\overline {W})\cong\Pi_{g\in{G}}H_q(K;\overline {A_g})$ 
for $q\leq [n/2]$, since $C_*$ has finite $[n/2]$-skeleton.
The $\mathbb{Z}$-linear homomorphisms 
$t_{q,g}:\overline {A_g}\otimes_{\mathbb{Z}[K]}C_q\to 
\overline {A}\otimes_{\mathbb{Z}[K]}C_q$ given by
$t_{q,g}(a\otimes{c})=w(\sigma(g))a\otimes\sigma(g)c$ 
for all $a\in\overline {A}$ and 
$c\in{C_q}$ induce isomorphisms 
$H_q(K;\overline {A_g})\cong{H_q(K;\overline {A})}$ 
for all $q\geq0$ and $g\in{G}$.
Let $u_{q,g}=t_{q,g}(p_g\otimes{id}_{C_q})$.
Then $u_{q,g}(f\sigma(h)^{-1}\otimes\sigma(h){c})=u_{q,gh}(f\otimes{c})$
for all $g, h\in{G}$, $f\in\overline{W}$, $c\in{C_q}$ and $q\geq0$.
Hence these composites determine isomorphisms of left
$\mathbb{Z}[G]$-modules $H_q(K;\overline {W})\cong {A_q^G}$,
where $A_q=H_q(\overline{A}\otimes_{\mathbb{Z}[K]}C_*)=H_q(K;\overline{A})$
(with trivial $G$-action) for $q\leq[n/2]$.

Let $D(L)$ denote the conjugate of a left $\mathbb{Z}[G]$-module $L$
with respect to the canonical involution.
We shall apply the homology LHSSS
\[
E^2_{pq}=H_p(G;D(H_q(K;\overline{W}))\Rightarrow H_{p+q}(\pi;\overline{W}).
\]
Poincar\'e duality for $G$ and another application of Shapiro's Lemma
now give $H_p(G;D(A_q^G))\cong H^{r-p}(G;A_q^G)\cong H^{r-p}(1;A_q)$,
since $A_q^G$ is coinduced from a module over the trivial group.
If $s>[n/2]$ and $p+q=n-s$ then $q\leq[n/2]$ and so
$H_p(G;A_q^G)\cong A_q$ if $p=r$ and is 0 otherwise.
Thus the spectral sequence collapses to give 
$H_{n-s} (\pi;\overline{W})\cong H_{n-r-s}(K;\overline{A})$.
Since homology commutes with direct limits this proves the theorem.
\end{proof}

The finiteness condition cannot be relaxed further when $r=2$ and $n=4$,
for Kapovich has given an example of a pair $\nu<\pi$ with $\pi$ a $PD_4$-group,
$\pi/\nu$ a $PD_2$-group and $\nu$ finitely generated but not $FP_2$ 
\cite{[Ka98]}.

The most useful case of this theorem is when $G\cong\mathbb{Z}$.
The argument of the first paragraph of the theorem shows that if $K$ is
any normal subgroup such that $\pi/K\cong\mathbb{Z}$ then
$H^n(K;A)\cong H_0(\pi;\overline{W})=0$, and so $c.d.K<n$.
(This weak version of Strebel's Theorem suffices for some of the applications 
below.)

Let $R$ be a ring.
An $R$-chain complex has {\it finite $k$-skeleton} if it is chain homotopy
equivalent to a complex $P_*$ with $P_j$ a finitely generated free $R$-module 
for $j\leq{k}$.
If $R$ is a subring of $S$ and $C_*$ is an $S$-chain complex 
then $C_*$ is $R$-{\it finitely dominated} if $C_*|_R$
is chain homotopy equivalent to a finite projective $R$-chain complex.
The argument of Theorem 1.19 extends easily to the nonaspherical case
as follows. 
(See Chapter 2 for the definition of $PD_n$-space.)

\begin{thmm}
[1.19$^\prime$]
Let $M$ be a $PD_n$-space, $p:\pi_1(M)\to G$ be an epimorphism with $G$ 
a $PD_r$-group and $\nu=\mathrm{Ker}(p)$.
If $C_*(\widetilde{M})|_\nu$ has finite $[n/2]$-skeleton then
$C_*(\widetilde{M})$ is $\mathbb{Z}[\nu]$-finitely dominated
and $H^s (M_\nu;\mathbb{Z}[\nu])\cong H_{n-r-s}(M_\nu;\mathbb{Z}[\nu])$, 
for all $s$.
\qed
\end{thmm}

If $M$ is aspherical then $M_\nu=K(\nu,1)$ is a $PD_{n-r}$-space, 
by Theorem 1.19.
In Chapter 4 we shall show that this holds in general.

\begin{cor}
If either $r=n-1$ or $r=n-2$ and $\nu$ is infinite or $r=n-3$ and
$\nu$ has one end then $M$ is aspherical. 
\qed
\end{cor}

\section{Hilbert modules}

Let $\pi$ be a countable group and let $\ell^2(\pi)$ be the Hilbert space 
completion of $\mathbb{C}[\pi]$ with respect to the inner product given by 
$(\Sigma a_g g,\Sigma b_h h)=\Sigma a_g \overline{b_g} $.
Left and right multiplication by elements of $\pi$ determine left and right 
actions of $\mathbb{C}[\pi]$ as bounded operators on $\ell^2(\pi)$.
The (left) von Neumann algebra $\mathcal{ N}(\pi)$ is the algebra of bounded 
operators on 
$\ell^2(\pi)$ which are $\mathbb{C}[\pi]$-linear with respect to the left action. 
By the Tomita-Takesaki theorem this is also the bicommutant in $B(\ell^2(\pi))$
of the right action of $\mathbb{C}[\pi]$, i.e., the set of operators which
commute with every operator which is right $\mathbb{C}[\pi]$-linear.
(See \cite[pages 45-52]{[Su]}.) 
We may clearly use the canonical involution of $\mathbb{C}[\pi]$
to interchange the roles of left and right in these definitions.

If $e\in\pi$ is the unit element we may define the von Neumann trace on 
$\mathcal{ N}(\pi)$ by the inner product $tr(f)=(f(e),e)$. 
This extends to square matrices over $\mathcal{ N}(\pi)$ by taking the sum of 
the traces of the diagonal entries.
A {\it Hilbert $\mathcal{ N}(\pi)$-module} is a Hilbert space $M$ 
with a unitary left $\pi$-action which embeds isometrically 
and $\pi$-equivariantly into the completed tensor product
$H\widehat\otimes\ell^2(\pi)$ for some Hilbert space $H$. 
It is finitely generated if we may take $H\cong\mathbb{C}^n$ 
for some integer $n$.
(In this case we do not need to complete the ordinary tensor product over
$\mathbb{C}$.)
A {\it morphism} of Hilbert $\mathcal{ N}(\pi)$-modules is a $\pi$-equivariant
bounded linear operator $f:M\to N$.
It is a {\it weak isomorphism} if it is injective and has dense image.
A bounded $\pi$-linear operator on $\ell^2(\pi)^n=\mathbb{C}^n \otimes\ell^2(\pi)$
is represented by a matrix whose entries are in $\mathcal{ N}(\pi)$.
The {\it von Neumann dimension} of a finitely generated Hilbert 
$\mathcal{ N}(\pi)$-module $M$ is the real number
$dim_{\mathcal{ N}(\pi)}(M)=tr(P)\in [0,\infty)$, 
where $P$ is any projection operator on $H\otimes\ell^2(\pi)$ 
with image $\pi$-isometric to $M$. 
In particular, $dim_{\mathcal{ N}(\pi)}(M)=0$ if and only if $M=0$.
The notions of finitely generated Hilbert $\mathcal{ N}(\pi)$-module 
and finitely generated projective $\mathcal{ N}(\pi)$-module 
are essentially equivalent,
and arbitrary $\mathcal{ N}(\pi)$-modules have well-defined dimensions in
$[0,\infty]$ \cite{[Lu]}.

If $\pi$ is residually finite or satisfies the Strong Bass Conjecture 
and $P$ is a finitely generated
projective $\mathbb{Z}[\pi]$-module then
$\ell^2(\pi)\otimes{P}\cong\ell^2(\pi)^{\kappa(P)}$ \cite{[Ec96]}.

A sequence of bounded maps between Hilbert $\mathcal{ N}(\pi)$-modules 
\begin{equation*}
\begin{CD}
M@> j >> N@> p >> P
\end{CD}
\end{equation*}
is {\it weakly exact at} $N$ if $\mathrm{Ker}(p)$ is the closure of 
$\mathrm{Im}(j)$. 
If $0\to M\to N\to P\to 0$ is weakly exact then $j$ is injective, 
$\mathrm{Ker}(p)$ is the closure of $\mathrm{Im}(j)$ and $\mathrm{Im}(p)$ 
is dense in $P$,
and $dim_{\mathcal{ N}(\pi)}(N)=dim_{\mathcal{ N}(\pi)}(M)+dim_{\mathcal{ N}(\pi)}(P)$.
A finitely generated {\it Hilbert $\mathcal{ N}(\pi)$-complex} $C_*$ is a chain 
complex of finitely generated Hilbert $\mathcal{ N}(\pi)$-modules with bounded 
$\mathbb{C}[\pi]$-linear operators as differentials.
The {\it reduced} $L^2$-homology is defined to be 
$\bar H_p^{(2)}(C_*)=\mathrm{Ker}(d_p)/\overline{\mathrm{Im}(d_{p+1})}$.
The $p^{th}$ $L^2$-Betti number of $C_*$ is then 
$dim_{\mathcal{ N}(\pi)}\bar H_p^{(2)}(C_*)$. 
(As the images of the differentials need not be closed the {\it un}reduced 
$L^2$-homology modules $H_p^{(2)}(C_*)=\mathrm{Ker}(d_p)/\mathrm{Im}(d_{p+1})$
are not in general Hilbert modules.)
 
See \cite{[Lu]} for more on modules over von Neumann algebras
and $L^2$ invariants of complexes and manifolds.                                                                                          
    
[In this book $L^2$-Betti number arguments replace the localization
arguments used in \cite{[H2]}.
However we shall recall the definition of {\it safe extension} 
of a group ring used there.
An inclusion of rings $\mathbb{Z}[G]<S$ is a safe extension if it is flat, 
$S$ is weakly finite and $S\otimes_{\mathbb{Z}[G]}\mathbb{Z}=0$.
If $G$ has a nontrivial elementary amenable normal
subgroup whose finite subgroups have bounded order and which has no nontrivial
finite normal subgroup then $\mathbb{Z}[G]$ has a safe extension.]

%% file: m5-2.tex
\chapter{2-Complexes and $PD_3$-complexes}

This chapter begins with a review of the notation we use for (co)homology
with local coefficients and of the universal coefficient spectral sequence.
We then define the $L^2$-Betti numbers and present some useful vanishing
theorems of L\"uck and Gromov. 
These invariants are used in \S3, where they are used to estimate 
the Euler characteristics of finite $[\pi,m]$-complexes and to give 
a converse to the Cheeger-Gromov-Gottlieb Theorem on aspherical finite complexes.   
Some of the arguments and results here may be regarded as representing 
in microcosm the bulk of this book; the analogies and connections between 2-complexes and 4-manifolds are well known.                                                          
We then review Poincar\'e duality and $PD_n$-complexes.
In \S5-\S9 we shall summarize briefly what is known about the homotopy types 
of $PD_3 $-complexes.

\section{Notation}

Let $X$ be a connected cell complex and let $\widetilde X$ be its universal 
covering space. 
If $H$ is a normal subgroup of $G=\pi_1 (X)$ we may lift the cellular 
decomposition of $X$ to an equivariant cellular decomposition of the 
corresponding covering space $X_H$.
The cellular chain complex of $X_H$ with coefficients in a commutative 
ring $R$ is then a complex $C_*=C_*(X_H)$ of left $R[G/H]$-modules, 
with respect to the action of the covering group $G/H$. 
A choice of lifts of the $q$-cells of $X$ determines a free basis for $C_q$, 
for all $q$, and so $C_* $ is a complex of free modules.
If $X$ is a finite complex $G$ is finitely presentable and these modules 
are finitely generated.
If $X$ is finitely dominated, i.e., is a retract of a finite complex, 
then $G$ is again finitely presentable, by Lemma 1.12.
Moreover the chain complex of the universal cover is chain homotopy equivalent 
over $R[G]$ to a complex of finitely generated projective modules 
\cite{[Wl65]}.
The Betti numbers of $X$ with coefficients in a field $F$ shall be denoted by 
$\beta_i(X;F)=dim_F H_i(X;F)$ (or just $\beta_i(X)$, if $F=\mathbb{Q}$).

The $i^{th} $ {\it equivariant homology} module of $X$ with 
coefficients $R[G/H]$ is the left module $H_i (X;R[G/H])=H_i (C_*)$, 
which is clearly isomorphic to $H_i (X_H;R)$
as an $R$-module, with the action of the covering group determining its $R[G/H]$-module structure.
The $i^{th} $ {\it equivariant cohomology} module of $X$ 
with coefficients $R[G/H]$ is the right module $H^i (X;R[G/H])=H^i (C^* )$,
where  $C^* =Hom_{R[G/H]} (C_* ,R[G/H])$
is the associated cochain complex of right $R[G/H]$-modules.
More generally, if $A$ and $B$ are right and left $\mathbb{Z}[G/H]$-modules (respectively) we may define         
$H_j (X;A)=H_j (A\otimes_{\mathbb{Z}[G/H]} C_*)$ and $H^{n-j} (X;B)=H^{n-j} (Hom_{\mathbb{Z}[G/H]} (C_* ,B))$.
There is a {\it Universal Coefficient Spectral Sequence} (UCSS) 
relating equivariant homology and cohomology:
\begin{equation*}
E_2^{pq}=Ext_{R[G/H]}^q (H_p (X;R[G/H]),R[G/H])\Rightarrow H^{p+q}
(X;R[G/H]),
\end{equation*}
with $r^{th}$ differential $d_r $ of bidegree $(1-r,r)$.

If $J$ is a normal subgroup of $G$ which contains $H$ there is also a 
{\it Cartan-Leray} spectral sequence relating the homology of $X_H$ and $X_J$: 
\begin{equation*}
E^2_{pq}=Tor^{R[G/H]}_p (R[G/J],H_q(X;R[G/H]))\Rightarrow H_{p+q}(X;R[G/J]),
\end{equation*}
with $r^{th}$ differential $d^r$ of bidegree $(-r,r-1)$.
(See \cite{[Mc]} for more details on these spectral sequences.)

If $M$ is a cell complex let $c_M :M\to K(\pi_1 (M),1)$ denote the 
{\it classifying map} for the fundamental group and let $f_M :M\to P_2 (M)$ 
denote the second stage of the Postnikov tower for $M$.
(Thus $c_M =c_{P_2 (M)} f_M$.)
A map $f:X\to K(\pi_1 (M),1)$ lifts to a map from $X$ to
$P_2 (M)$ if and only if $f^* k_1 (M)=0$, where $k_1 (M)$ is the first $k$-invariant of 
$M$ in $H^3 (\pi_1 (M);\pi_2 (M))$.    
In particular, if $k_1 (M)=0$ then $c_{P_2 (M)} $ has a cross-section.
The {\it algebraic $2$-type} of $M$ is the triple $[\pi,\pi_2 (M),k_1 (M)]$.
Two such triples $[\pi,\Pi,\kappa]$ and $[\pi',\Pi',\kappa']$ 
(corresponding to $M$ and $M'$, respectively) are equivalent if there are 
isomorphisms $\alpha:\pi\to\pi'$ and $\beta:\Pi\to\Pi'$ such that 
$\beta(gm)=\alpha(g)\beta(m)$ for all $g\in\pi$ and $m\in\Pi$
and $\beta_* \kappa=\alpha^* \kappa'$ in $H^3 (\pi;\alpha^* \Pi' )$. 
Such an equivalence may be realized by a homotopy equivalence of 
$P_2(M)$ and $P_2(M')$.
(The reference \cite{[Ba]} gives a detailed treatment of 
Postnikov factorizations of nonsimple maps and spaces.)                                                       
Throughout this book {\it closed manifold} shall mean compact, 
connected TOP manifold without boundary.
Every closed manifold has the homotopy type of a finite 
Poincar\'e duality complex \cite{[KS]}.                                          
\section{$L^2$-Betti numbers}

Let $X$ be a finite complex with fundamental group $\pi$.
The $L^2$-{\it Betti numbers} of $X$ are defined by 
$\beta_i^{(2)} (X)=dim_{\mathcal{ N}(\pi)}(\bar H_2^{(2)}(\widetilde X))$, 
where the $L^2$-homology $\bar H_i^{(2)}(\widetilde X)=\bar H_i (C^{(2)}_*)$ 
is the reduced homology of the Hilbert $\mathcal{ N}(\pi)$-complex 
$C^{(2)}_*=\ell^2\otimes_{\mathbb{Z}[\pi]} C_*(\widetilde X)$ 
of square summable chains on $\widetilde X$. 
They are multiplicative in finite covers, 
and for $i=0$ or 1 depend only on $\pi$. 
(In particular, $\beta_0^{(2)}(\pi)=0$ if $\pi$ is infinite.)
The alternating sum of the $L^2$-Betti numbers is the Euler characteristic 
$\chi(X)$. (See \cite{[Lu]}.)

It may be shown that $\beta_i^{(2)}(X)=
dim_{\mathcal{ N}(\pi)}H_i(\mathcal{ N}(\pi)\otimes_{\mathbb{Z}[\pi]}C_*(\widetilde X))$,
and this formulation of the definition applies to arbitrary complexes 
\cite{[CG86],[Lu]}.
In particular, 
$\beta_i^{(2)}(\pi)=dim_{\mathcal{ N}(\pi)}H_i(\pi;\mathcal{ N}(\pi))$
is defined for all $\pi$.
If $X$ is finitely dominated then $\Sigma\beta_i^{(2)}(X)<\infty$,
and if also $\pi$ satisfies the Strong Bass Conjecture
then 
$\chi(X)=\Sigma(-1)^i\beta_i^{(2)}(X)$ \cite{[Ec96]}.
Moreover,  
$\beta_s^{(2)}(X)=\beta_s^{(2)} (\pi)$ for $s=0$ or 1,
and $\beta_2^{(2)}(X)\geq \beta_2^{(2)}(\pi)$.
(See Theorems 1.35 and 6.54 of \cite{[Lu]}.)
If $\pi\cong{A*_CB}$ the argument for \cite[Theorem 1.35.5]{[Lu]} extends
to give $\beta_1^{(2)}(\pi)\geq\frac1{|C|}-\frac1{|A|}-\frac1{|B|}$.
(Similarly for $A*_C\phi$.)
Thus if $\beta_1^{(2)}(\pi)=0$ then $e(\pi)$ is finite.
  
\begin{lemma}
Let $\pi=H*_\phi$ be a finitely presentable group 
which is an ascending HNN extension with finitely generated base $H$. 
Then $\beta_1^{(2)}(\pi)=0$.
\end{lemma}

\begin{proof}
Let $t$ be the stable letter and let $H_n$ be the subgroup generated by
$H$ and $t^n$, and suppose that $H$ is generated by $g$ elements.
Then $[\pi:H_n]=n$, so $\beta_1^{(2)} (H_n)=n\beta_1^{(2)}(\pi)$.
But each $H_n$ is also finitely presentable and generated by $g+1$ elements.
Hence $\beta_1^{(2)} (H_n)\leq g+1$, and so $\beta_1^{(2)} (\pi)=0$. 
\end{proof}

In particular, this lemma holds if $H$ is normal in $\pi$
and $\pi/H\cong\mathbb{Z}$.

\begin{theorem}
[L\"uck] 
Let $\pi$ be a group with a finitely generated infinite normal subgroup 
$\Delta$ such that $\pi/\Delta$ has an element of infinite order.
Then $\beta_1^{(2)}(\pi)=0$.
\end{theorem}

\begin{proof}
(Sketch)\qua
Let $\rho\leq \pi$ be a subgroup containing $\Delta$ such that 
$\rho/\Delta\cong\mathbb{Z}$.
The terms in the line $p+q=1$ of the homology LHSSS for $\rho$ 
as an extension of $\mathbb{Z}$ by $\Delta$ with coefficients 
$\mathcal{ N}(\rho)$ have dimension 0, 
by Lemma 2.1. 
Since $dim_{\mathcal{ N}(\rho)}M=
dim_{\mathcal{ N}(\pi)} (\mathcal{ N}(\pi)\otimes_{\mathcal{ N}(\rho)}M)$
for any $\mathcal{ N}(\rho)$-module $M$ 
the corresponding terms for the LHSSS for $\pi$ as an extension of $\pi/\Delta$
by $\Delta$ with coefficients $\mathcal{ N}(\pi)$ also have
dimension 0 and the theorem follows. 
\end{proof}

This is Theorem 7.2.6 of \cite{[Lu]}.
The hypothesis ``$\pi/\Delta$ has an element of infinite order" 
can be relaxed to ``$\pi/\Delta$ is infinite" \cite{[Ga00]}.
The next result also derives from \cite{[Lu]}.
(The case $s=1$ is extended further in \cite{[PT11]}.)

\begin{theorem}
Let $\pi$ be a group with an ascendant subgroup $N$ such that 
$\beta_i^{(2)}(N)=0$ for all $i\leq s$.
Then $\beta_i^{(2)}(\pi)=0$ for all $i\leq s$.
\end{theorem}

\begin{proof}
Let $N=N_0<N_1<\dots<N_\beth=\pi$ be an ascendant sequence.
Then we may show by transfinite induction on $\alpha$ that 
$\beta_i^{(2)}(N_\alpha)=0$ for all $i\leq s$ and $\alpha\leq\beth$,
using parts (2) and (3) of \cite[Theorem 7.2]{[Lu]} 
for the passages to successor ordinals and to limit ordinals,
respectively.
\end{proof}

\begin{cor}
[Gromov] Let $\pi$ be a group with an infinite 
amenable normal subgroup $A$. 
Then $\beta_i^{(2)}(\pi)=0$ for all $i$.
\end{cor}

\begin{proof}
If $A$ is an infinite amenable group $\beta_i^{(2)}(A)=0$ for all $i$ 
\cite{[CG86]}.                                                                   \end{proof}                     

Note that the normal closure of an amenable ascendant subgroup 
is amenable.

\section{2-Complexes and finitely presentable groups}

If a group $\pi$ has a finite presentation $P$ with $g$ generators and $r$ relators
then the {\it deficiency} of $P$ is $\mathrm{def}(P)=g-r$, and $\mathrm{def}(\pi)$ is the maximal 
deficiency of all finite presentations of $\pi$. 
Such a presentation determines a finite 2-complex $C(P)$ with 
one 0-cell, $g$ 1-cells and $r$ 2-cells and with $\pi_1(C(P))\cong\pi$.
Clearly $\mathrm{def}(P)=1-\chi(P)=\beta_1(C(P))-\beta_2(C(P))$ and so                                                  
$\mathrm{def}(\pi)\leq \beta_1(\pi)-\beta_2(\pi)$.
Conversely every finite 2-complex with one 0-cell arises in this way.
In general, any connected finite 2-complex $X$ is homotopy equivalent to one with a single 0-cell,
obtained by collapsing a maximal tree $T$ in the 1-skeleton $X^{[1]}$.

We shall say that $\pi$ has {\it geometric dimension at most} $2$, 
written $g.d.\pi\leq 2$,
if it is the fundamental group of a {\it finite} aspherical 2-complex. 

\begin{theorem}
Let $X$ be a connected finite $2$-complex with 
fundamental group $\pi$.
Then $\beta_2^{(2)}(X)\geq\beta_2^{(2)}(\pi)$, 
with equality if and only if $X$ is aspherical.
\end{theorem}

\begin{proof}
Since we may construct $K=K(\pi,1)$ by adjoining cells of dimension $\geq3$
to $X$ the natural homomorphism $\bar{H}_2(c_X)$ is an epimorphism, 
and so $\beta_2^{(2)}(X)\geq\beta_2^{(2)}(\pi)$.
Since $X$ is $2$-dimensional $\pi_2(X)=H_2(\widetilde X;\mathbb{Z})$ 
is a subgroup of $\bar H_2^{(2)}(\widetilde X)$, 
with trivial image in $\bar H_2^{(2)}(\widetilde{K})$.
If moreover $\beta_2^{(2)}(X)=\beta_2^{(2)}(\pi)$ then $\bar{H}_2(c_X)$
is an isomorphism \cite[Lemma 1.13]{[Lu]}, 
so $\pi_2(X)=0$ and $X$ is aspherical. 
\end{proof}
  
\begin{cor}
Let $\pi$ be a finitely presentable group.
Then\\ 
$\mathrm{def}(\pi)\leq 1+\beta_1^{(2)}(\pi)-\beta_2^{(2)}(\pi)$.
If equality holds then $g.d.\pi\leq 2$.
\end{cor}

\begin{proof}
This follows from the theorem and the $L^2$-Euler characteristic formula, 
applied to the 2-complex associated to an optimal presentation for $\pi$.
\end{proof}
                       
\begin{theorem}
Let $\pi$ be a finitely presentable group such that $\beta_1^{(2)}(\pi)=0$.
Then $\mathrm{def}(\pi)\leq1$, 
with equality if and only if $g.d.\pi\leq 2$ and 
$\beta_2(\pi)=\beta_1(\pi)-1$.
\end{theorem}
                      
\begin{proof}
The upper bound and the necessity of the conditions follow as in 
Corollary 2.4.1.
Conversely, if they hold and $X$ is a finite aspherical 2-complex 
with $\pi_1(X)\cong\pi$ then
$\chi(X)=1-\beta_1(\pi)+\beta_2(\pi)=0$. 
After collapsing a maximal tree in $X$ we may assume
it has a single 0-cell, 
and then the presentation read off the 1- and 2-cells has deficiency 1.
\end{proof}

This theorem applies if $\pi$ is finitely presentable and is 
an ascending HNN extension with finitely generated base $H$, or has an
infinite amenable normal subgroup.
In the latter case $\beta_i^{(2)}(\pi)=0$ for all $i$, by Theorem 2.3.
Thus if $X$ is a finite aspherical 2-complex with $\pi_1(X)\cong\pi$
then $\chi(X)=0$, and so the condition $\beta_2(\pi)=\beta_1(\pi)-1$ 
is redundant.

[Similarly, if $\mathbb{Z}[\pi]$ has a safe extension $\Psi$ and
$C_*$ is the equivariant cellular chain complex of the universal cover 
$\widetilde X$
then $\Psi\otimes_{\mathbb{Z}[\pi]}C_*$ is a complex of free left $\Psi$-modules
with bases corresponding to the cells of $X$.
Since $\Psi$ is a safe extension 
$H_i(X;\Psi)=\Psi\otimes_{\mathbb{Z}[\pi]} H_i(X;\mathbb{Z}[\pi])=0$ 
for all $i$,
and so again $\chi(X)=0$.]

\begin{cor}
Let $\pi$ be a finitely presentable group with an $FP_2$ normal subgroup $N$
such that $\pi/N\cong\mathbb{Z}$.
Then $\mathrm{def}(\pi)=1$ if and only if $N$ is free.
\end{cor}

\begin{proof} 
If $\mathrm{def}(\pi)=1$ then $g.d.\pi\leq2$, by Theorem 2.5,
and so $N$ is free \cite[Corollary 8.6]{[Bi]}.
The converse is clear.
\end{proof}
 
In fact it suffices to assume that $N$ is finitely generated 
(rather than $FP_2$) \cite{[Ko06]}. (See Corollary 4.3.1 below.) 

Let $G=F(2)\times F(2)$.
Then $g.d.G=2$ and $\mathrm{def}(G)\leq\beta_1(G)-\beta_2(G)=0$.
Hence $\langle u,v,x,y \mid ux=xu,\medspace uy=yu,\medspace vx=xv,\medspace
vy=yv\rangle$ is an optimal presentation, and $\mathrm{def}(G)=0$. 
The subgroup $N$ generated by $u$, $vx^{-1}$ and $y$
is normal in $G$ and $G/N\cong Z$, so $\beta_1^{(2)}(G)=0$, by Lemma 2.1.
However $N$ is not free,
since $u$ and $y$ generate a rank two abelian subgroup.
It follows from Corollary 2.5.1 that $N$ is not $FP_2$,
and so $F(2)\times F(2)$ is not almost coherent.

The next result is a version of the {\it Tits alternative\/} 
for coherent groups of cohomological dimension 2.
For each $m\in\mathbb{Z}$ let $Z*_m$ be the group with presentation 
$\langle a,t\mid tat^{-1}=a^m\rangle$.
(Thus $Z*_0\cong\mathbb{Z}$ and $Z*_{-1}\cong\mathbb{Z}\rtimes_{-1}\mathbb{Z}$.)                
\begin{theorem}
Let $\pi$ be a finitely generated group 
such that $c.d.\pi=2$.
Then $\pi\cong Z*_m$ for some $m\not=0$ if and only if it is almost coherent 
and restrained and $\pi/\pi'$ is infinite.
\end{theorem}

\begin{proof}
The conditions are easily seen to be necessary.
Conversely, if $\pi$ is almost coherent and $\pi/\pi'$ is infinite $\pi$ is an 
HNN extension with $FP_2$ base $H$, by Theorem 1.13.
The HNN extension must be ascending as $\pi$ has no noncyclic free subgroup.
Hence $H^2(\pi;\mathbb{Z}[\pi])$ is a quotient of 
$H^1(H;\mathbb{Z}[\pi])\cong H^1(H;\mathbb{Z}[H])\otimes\mathbb{Z}[\pi/H]$, 
by the Brown-Geoghegan Theorem.
Now $H^2(\pi;\mathbb{Z}[\pi])\not=0$, 
since $\pi$ is $FP_2$ and $c.d.\pi=2$, 
and so $H^1(H;\mathbb{Z}[H])\not=0$.
Since $H$ is restrained it must have two ends, 
so $H\cong\mathbb{Z}$ and $\pi\cong Z*_m$ for some $m\not=0$.
\end{proof}

\begin{cor}
Let $\pi$ be a finitely generated group. 
Then the following are equivalent:
\begin{enumerate}           
\item $\pi\cong Z*_m$ for some $m\in\mathbb{Z}$;

\item $\pi$ is torsion-free, elementary amenable, $FP_2$ and $h(\pi)\leq2$;

\item $\pi$ is elementary amenable and $c.d.\pi\leq2$;

\item $\pi$ is almost coherent, amenable and $c.d.\pi\leq2$;

\item $\pi$ is elementary amenable and $\mathrm{def}(\pi)=1$; 
and

\item $\pi$ is almost coherent, restrained and $\mathrm{def}(\pi)=1$.

\end{enumerate}   
\end{cor}                    

\begin{proof}
Condition (1) clearly implies the others.
Suppose (2) holds.    
We may assume that $h(\pi)=2$ and $h(\sqrt\pi )=1$ 
(for otherwise $\pi\cong\mathbb{Z}$, $\mathbb{Z}^2=Z*_1$ or 
$\mathbb{Z}\rtimes_{-1}\mathbb{Z}=Z*_{-1}$). 
Hence $h(\pi/\sqrt\pi )=1$, 
and so $\pi/\sqrt\pi$ is an extension of $\mathbb{Z}$ or $D$ by
a finite normal subgroup.
If $\pi/\sqrt\pi$ maps onto $D$ then $\pi\cong A*_CB$, 
where $[A:C]=[B:C]=2$ and $h(A)=h(B)=h(C)=1$, 
and so $\pi\cong{Z*_{-1}}$.
But then $h(\sqrt\pi)=2$.
Hence we may assume that $\pi$ maps onto $\mathbb{Z}$, and so
$\pi$ is an ascending HNN extension with finitely generated base $H$, 
by Theorem 1.13. 
Since $H$ is torsion-free, elementary amenable and $h(H)=1$ 
it must be infinite cyclic and so (2) implies (1).
If (3) holds $\pi$ is solvable, by Theorems 1.11, and 1.9, and so 
(1) follows from \cite{[Gi79]}.
If (4) holds then $\pi$ is restrained,
and $\chi(\pi)=0$ \cite{[Ec96]}, so $\pi/\pi'$ is infinite.
If $\mathrm{def}(\pi)=1$ then $\pi$ is an ascending HNN extension with 
finitely generated base,
so $\beta_1^{(2)}(\pi)=0$, by Lemma 2.1. 
Hence (4), (5) and (6) each imply (1), by Theorems 2.5 and 2.6. 
\end{proof}

If $\pi$ is $FP_2$ then $(3)\Rightarrow(2)$ 
(without invoking \cite{[Gi79]}).
Are these conditions equivalent to 
``$\pi$ is restrained and $c.d.\pi\leq2$" or
``$\pi$ is restrained and def$(\pi)=1$"? 
(Note that if $\mathrm{def}(\pi)>1$ then 
$\pi$ has noncyclic free subgroups \cite{[Ro77]}.)

Let $\mathcal{X}$ be the class of groups of finite graphs of groups, 
with all edge and vertex groups infinite cyclic.
A finitely generated,
noncyclic group $G$ is in $\mathcal{ X}$ if and only if 
$c.d.G=2$ and $G$ has an infinite cyclic subgroup $H$ 
which meets all its conjugates nontrivially \cite{[Kr90']}.
Moreover, $G$ is coherent, one ended, $\mathrm{def}(G)\geq1$ 
and $g.d.G=2$ \cite{[Kr90']},
while $\beta_1^{(2)}(G)=0$ \cite[Theorem 5.12]{[PT11]}.

\begin{theorem}
Let $\pi$ be a finitely generated group such that $c.d.\pi=2$.
If $\pi$ has a nontrivial normal subgroup $E$ which is either elementary 
amenable or almost coherent, 
locally virtually indicable and restrained then $\pi$ 
is in $\mathcal{ X}$, $\mathrm{def}(G)=1$ and either $E\cong\mathbb{Z}$ 
or $\pi'$ is abelian.
\end{theorem}

\begin{proof}
Since $c.d.E\leq c.d.\pi$, 
finitely generated subgroups of $E$ are metabelian,
by Theorems 1.11 and  2.6 and Corollary 2.6.1, 
and so all words in $E$ of the form $[[g,h],[g',h']]$ are trivial. 
Hence $E$ is metabelian also.
Therefore $A=\sqrt E$ is nontrivial, 
and as $A$ is characteristic in $E$ it is normal in $\pi$.
Since $A$ is the union of its finitely generated subgroups, 
which are torsion-free nilpotent groups of Hirsch length $\leq2$, it is abelian.
If $A\cong\mathbb{Z}$ then $[\pi:C_\pi(A)]\leq 2$.
Moreover $C_\pi(A)'$ is free, by Bieri's Theorem.
If $C_\pi(A)'$ is cyclic then $\pi\cong\mathbb{Z}^2$ or 
$\mathbb{Z}\!\rtimes_{-1}\!\mathbb{Z}$; 
if $C_\pi(A)'$ is nonabelian then $E=A\cong\mathbb{Z}$.
Otherwise $c.d.A=c.d.C_\pi(A)=2$ and so $C_\pi(A)=A$, by Bieri's Theorem.
If $A$ has rank 1 then $Aut(A)$ is abelian, so $\pi'\leq C_\pi(A)$ and $\pi$ 
is metabelian.
If $A\cong\mathbb{Z}^2$ then $\pi/A$ is isomorphic to a subgroup of $GL(2,\mathbb{Z})$,
and so is virtually free.                                          
As $A$ together with an element $t\in\pi$ of infinite order modulo $A$ would
generate a subgroup of cohomological dimension 3, which is impossible,
the quotient $\pi/A$ must be finite.
Hence $\pi\cong\mathbb{Z}^2$ or $\mathbb{Z}\!\rtimes_{-1}\!\mathbb{Z}$.
In all cases $\pi$ is in $\mathcal{ X}$ \cite[Theorem C]{[Kr90']}.
Since $\mathrm{def}(G)\geq1$ and $\beta_1^{(2)}(G)=0$, 
we see that $\mathrm{def}(G)=1$.
\end{proof}

If $c.d.\pi=2$, $\zeta\pi\not=1$ and $\pi$ is nonabelian then 
$\zeta\pi\cong\mathbb{Z}$ and $\pi'$ is free, by Bieri's Theorem.
On the evidence of his work on 1-relator groups 
Murasugi conjectured that if $G$ is a finitely presentable group 
other than $\mathbb{Z}^2$ and $\mathrm{def}(G)\geq1$
then $\zeta G\cong\mathbb{Z}$ or 1, and is trivial if $\mathrm{def}(G)>1$,
and he verified this for classical link groups \cite{[Mu65]}.
Theorems 2.3, 2.5 and 2.7 together imply that if $\zeta G$ is infinite 
then $\mathrm{def}(G)=1$ and $\zeta G\cong\mathbb{Z}$.

It remains an open question whether every finitely presentable group 
of cohomological dimension 2 has geometric dimension 2. 
The following partial answer to this question was first obtained by W.Beckmann 
under the additional assumptions that $\pi$ is $FF$ and $c.d.\pi\leq2$
(see \cite{[Dy87']}).

\begin{theorem}
Let $\pi$ be a finitely presentable group.
Then $g.d.\pi\leq2$ if and only if $c.d._{\mathbb{Q}}\pi\leq2$ and 
$\mathrm{def}(\pi)=\beta_1(\pi)-\beta_2(\pi)$.
\end{theorem}
                                                                      
\begin{proof}
The necessity of the conditions is clear.
Suppose that they hold and that $C(P)$ is the 2-complex corresponding to
a presentation for $\pi$ of maximal deficiency.
The cellular chain complex of $\widetilde{C(P)}$ gives an exact sequence
\[
0\to K=\pi_2(C(P))\to \mathbb{Z}[\pi]^r\to 
\mathbb{Z}[\pi]^g\to \mathbb{Z}[\pi]\to\mathbb{Z}\to 0.
\]
Extending coefficients to $\mathbb{Q}$ gives a similar exact sequence,
with kernel $\mathbb{Q}\otimes_{\mathbb{Z}}{K}$ on the left.
As $c.d._{\mathbb{Q}}\pi\leq2$ the image of $\mathbb{Q}[\pi]^r$ 
in $\mathbb{Q}[\pi]^g$ is projective, by Schanuel's Lemma.
Therefore the inclusion of $\mathbb{Q}\otimes_{\mathbb{Z}}{K}$ into 
$\mathbb{Q}[\pi]^r$ splits, 
and $\mathbb{Q}\otimes_{\mathbb{Z}}{K}$ is projective.
Moreover $dim_\mathbb{Q}(\mathbb{Q}\otimes_{\mathbb{Z}[\pi]}K)=0$,
and so $\mathbb{Q}\otimes_{\mathbb{Z}}{K}=0$, 
since the Weak Bass Conjecture holds for $\pi$ \cite{[Ec86]}.
Since $K$ is free as an abelian group it imbeds in 
$\mathbb{Q}\otimes_{\mathbb{Z}}{K}$, and so is also 0.
Hence $\widetilde{C(P)}$ is contractible, and so $C(P)$ is aspherical.
\end{proof}

The arguments of this section may easily be extended to other highly connected 
finite complexes.
A $[\pi,m]_f$-complex is a finite $m$-dimensional complex $X$ with 
$\pi_1(X)\cong\pi$ and with $(m-1)$-connected universal cover $\widetilde X$.
Such a $[\pi,m]_f$-complex $X$ is aspherical if and only if $\pi_m(X)=0$. 
In that case we shall say that $\pi$ has geometric dimension at most $m$, 
written $g.d.\pi\leq m$. 

\begin{thmm}
[2.4$^\prime$]
Let $X$ be a $[\pi,m]_f$-complex and suppose that $\beta_i^{(2)}(\pi)=0$ for 
$i<m$.
Then $(-1)^m\chi(X)\geq0$. 
If $\chi(X)=0$ then $X$ is aspherical. 
\qed
\end{thmm}

In general, the final implication of this theorem cannot be reversed.
For $S^1\vee S^1$ is an aspherical $[F(2),1]_f$-complex and $\beta_0^{(2)}(F(2))=0$,
but ${\chi(S^1\vee{S^1})}\not=0$. 

One of the applications of $L^2$-cohomology in \cite{[CG86]} 
was to show that if $X$ is a finite aspherical complex 
and $\pi_1(X)$ has an infinite amenable normal subgroup $A$ then $\chi(X)=0$.
(This generalised a theorem of Gottlieb, 
who assumed that $A$ was a central subgroup \cite{[Go65]}.)
We may similarly extend Theorem 2.5 to give a converse to the 
Cheeger-Gromov extension of Gottlieb's Theorem.

\begin{thmm}
[2.5$^\prime$]
Let $X$ be a $[\pi,m]_f$-complex and suppose that 
$\pi$ has an infinite amenable normal subgroup.
Then $X$ is aspherical if and only if $\chi(X)=0$.
\qed
\end{thmm}

\section{Poincar\'e duality}

The main reason for studying $PD$-complexes is that 
they represent the homotopy theory of manifolds. 
However they also arise in situations where the geometry does not 
immediately provide a corresponding manifold.
For instance, under suitable finiteness assumptions 
an infinite cyclic covering space of a closed 4-manifold 
with Euler characteristic 0 will be a $PD_3 $-complex,
but need not be homotopy equivalent to a closed 3-manifold. (See Chapter 11.)

A $PD_n $-{\it space} is a space homotopy equivalent to a cell complex 
which satisfies Poincar\'e duality of formal dimension $n$ 
with local coefficients. 
If $X$ is a $PD_n$-space with fundamental group $\pi$ then $C_*(\widetilde{X})$ 
is $\mathbb{Z}[\pi]$-finitely dominated, so $\pi$ is $FP_2$.
The $PD_n$-space $X$ is {\it finite} if 
$C_*(\widetilde{X})$ is $\mathbb{Z}[\pi]$-chain homotopy equivalent to a 
finite free $\mathbb{Z}[\pi]$-complex. 
It is a {\it $PD_n$-complex} if it is finitely dominated.
This is so if and only if $\pi$ is finitely presentable 
\cite{[Br72],[Br75]}.
Finite $PD_n$-complexes are homotopy equivalent to finite complexes.
(Note also that a cell complex $X$ is finitely dominated if and only if 
$X\times S^1$ is finite \cite[Proposition 3]{[Rn95]}.)
Although $PD_n$-complexes are most convenient for our purposes,
the broader notion of $PD_n$-space is occasionally useful.
All the $PD_n$-complexes that we consider shall be connected.

Let $P$ be a $PD_n$-complex.
We may assume that $P=P_o\cup D^n$, where $P_o$ is a complex of dimension
$\leq\max\{3,n-1\}$ \cite{[Wl67]}. 
Let $\pi=\pi_1(P)$, $w=w_1(P)$ and $\pi^+=\mathrm{Ker}(w)$,
and let $P^+=P_{\pi^+}$ be the associated covering space.
If $C_*=C_*(\widetilde{P})$ the Poincar\'e 
duality isomorphism may be described in terms of a 
chain homotopy equivalence $\overline{C^*}\cong{C_{n-*}}$, 
which induces isomorphisms from $H^j (\overline {C^*})$ 
to $H_{n-j} (C_* )$, given by cap product with a generator $[P]$ of
$H_n (P;\mathbb{Z}^w)=
H_n (\overline{\mathbb{Z}}\otimes_{\mathbb{Z}[\pi]} C_* )$.
From this point of view it is easy to see that Poincar\'e duality gives rise to 
($\mathbb{Z}$-linear) isomorphisms from $H^j (P;B)$ to $H_{n-j} (P;\bar B)$, 
where $B$ is any left $\mathbb{Z}[\pi]$-module of coefficients. 
(See \cite{[Wl67]} or \cite[Chapter II]{[Wl]} for further details.)
If $P$ is a Poincar\'e duality complex then the $L^2$-Betti numbers also
satisfy Poincar\'e duality. 
(This does not require that $P$ be finite or orientable!)

A group $G$ is a $PD_n$-group (as defined in Chapter 1)
if and only if $K(G,1)$ is a $PD_n$-space.
For every $n\geq4$ there are $PD_n$-groups 
which are not finitely presentable \cite{[Da98]}.

Dwyer, Stolz and Taylor have extended Strebel's Theorem to show that 
if $H$ is a subgroup of infinite index in $\pi$ then the corresponding 
covering space $P_H$ has homological dimension $<n$; hence if moreover 
$n\not=3$ then $P_H$ is homotopy equivalent to a complex of dimension $<n$ 
\cite{[DST96]}.

\section{$PD_3$-complexes}

In this section we shall summarize briefly what is known 
about $PD_n$-complexes of dimension at most 3.
It is easy to see that a connected $PD_1 $-complex 
must be homotopy equivalent to $S^1$.
The 2-dimensional case is already quite difficult, but has been settled by
Eckmann, Linnell and M\"uller, who showed that every $PD_2 $-complex is 
homotopy equivalent to a closed surface. (See \cite[Chapter VI]{[DD]}.
This result has been further improved by Bowditch's Theorem.)
There are $PD_3 $-complexes with finite fundamental group 
which are not homotopy equivalent to any closed 3-manifold. 
On the other hand, Turaev's Theorem below implies that every $PD_3$-complex 
with torsion-free fundamental group is homotopy equivalent to a closed 
3-manifold if every $PD_3$-group is a 3-manifold group.
The latter is so if the Hirsch-Plotkin radical of the group is nontrivial
(see \S7 below), but remains open in general.

The {\it fundamental triple} of a $PD_3$-complex $P$ is 
$(\pi_1(P),w_1(P),c_{P*}[P])$.
This is a complete homotopy invariant for such complexes.
(See also \S6 and \S9 below.)

\begin{thm}
[Hendriks] Two $PD_3 $-complexes are homotopy 
equivalent if and only if their fundamental triples are isomorphic.
\qed
\end{thm}

When $w_1(P)\not=0$ the class $[P]$ is only well-defined up to sign 
\cite{[Ta08]}.
(This issue has no major consequences for us.)    
Turaev has characterized the possible triples corresponding to 
a given finitely presentable group and orientation character,
and has used this result to deduce a basic splitting theorem 
\cite{[Tu90]}.

\begin{thm}
[Turaev] A $PD_3 $-complex is indecomposable with respect to connected sum 
if and only if its fundamental group is indecomposable with respect to 
free product.
\qed
\end{thm}

Wall asked whether every orientable $PD_3 $-complex whose fundamental group 
has infinitely many ends is a proper connected sum \cite{[Wl67]}. 
Since the fundamental group of a $PD_n$-space is $FP_2$ it is
the fundamental group of a finite graph of finitely generated groups 
in which each vertex group has at most one end and each edge group 
is finite \cite[Theorem VI.6.3]{[DD]}.
Crisp has given a substantial partial answer to Wall's question,
based on this observation \cite{[Cr00]}.

\begin{thm}
[Crisp] Let $P$ be an indecomposable orientable $PD_3$-complex.
If $\pi_1(P)$ is not virtually free then it has one end, 
and so $P$ is aspherical.
\qed
\end{thm}

The arguments of Turaev and Crisp for these theorems
extend to $PD_3$-spaces in a straightforward manner.
In particular, they imply that if $P$ is a $PD_3$-space 
then $\pi=\pi_1(P)$ is virtually torsion-free.
However, there is an indecomposable orientable $PD_3$-complex 
with $\pi\cong{S_3*_{Z/2Z}S_3}\cong{F(2)}\rtimes{S_3}$ and double cover 
homotopy equivalent to $L(3,1)\sharp{L(3,1)}$.
``Most" indecomposable $PD_3$-complexes with $\pi$ virtually free 
have double covers which are homotopy equivalent to connected sums of 
$\mathbb{S}^3$-manifolds \cite{[Hi12]}.

\section{The spherical cases}

Let $P$ be a $PD_3$-space with fundamental group $\pi$,
and let $w=w_1(P)$.
The Hurewicz Theorem, Poincar\'e duality and a choice of
orientation for $P$ together determine an isomorphism
$\pi_2(P)\cong\overline{H^1(\pi;\mathbb{Z}[\pi])}$.
In particular, $\pi_2(P)=0$ if and only if $\pi$ is finite or has one end.

The possible $PD_3$-complexes with $\pi$ finite are well understood.

\begin{theorem}{\rm[Wl67]}\qua
Let $X$ be a $PD_3$-complex with finite fundamental group $F$. Then
\begin{enumerate}                                              
\item $\widetilde X\simeq S^3$, $F$ has cohomological period dividing 4 
and $X$ is orientable;

\item the first nontrivial $k$-invariant $k(X)$ generates
$H^4(F;\mathbb{Z})\cong Z/|F|Z$.

\item the homotopy type of $X$ is determined by $F$ and 
the orbit of $k(X)$ under $Out(F)\times\{\pm1\}$.
\end{enumerate}          
\end{theorem} 
         
\begin{proof} 
Since the universal cover $\widetilde X$ is also a finite $PD_3$-complex it
is homotopy equivalent to $S^3$.
A standard Gysin sequence argument shows that $F$ has cohomological period
dividing 4.
Suppose that $X$ is nonorientable, and let $C$ be a cyclic subgroup of $F$
generated by an orientation reversing element.
Let $\widetilde{\mathbb{Z}}$ be the nontrivial infinite cyclic
$\mathbb{Z}[C]$-module.
Then $H^2(X_C;\widetilde{\mathbb{Z}})\cong H_1 (X_C;\mathbb{Z})\cong C$, 
by Poincar\'e duality.
But $H^2 (X_C;\widetilde{\mathbb{Z}})\cong{H^2 (C;\widetilde{\mathbb{Z}})}=0$,
since the classifying map from $X_C =\widetilde X/C$ to $K(C,1)$ 
is 3-connected.
Therefore $X$ must be orientable and $F$ must act trivially on $\pi_3(X)\cong
H_3(\widetilde X;\mathbb{Z})$.

The image of the orientation class of $X$ generates  
$H_3 (F;\mathbb{Z})\cong Z/|F|Z$.
The Bockstein $\beta:H^3(F;\mathbb{Q}/\mathbb{Z})\to{H^4(F;\mathbb{Z})}$
is an isomorphism, since $H^q(F;\mathbb{Q})=0$ for $q>0$,
and the bilinear pairing from
$H_3(F;\mathbb{Z})\times{H^4(F;\mathbb{Z})}$ to $\mathbb{Q}/\mathbb{Z}$ 
given by $(h,c)\mapsto\beta^{-1}(c)(h)$ is nonsingular.
Each generator $g$ of $H_3(F;\mathbb{Z})$ determines an unique
$k_g\in{H^4(F;\mathbb{Z})}$ such that $\beta^{-1}(k_g)(g)=\frac1{|F|}$ {\it mod}
$\mathbb{Z}$.
The element corresponding to $c_{X*}[X]$ is the first nontrivial 
$k$-invariant of $X$ \cite{[Th67]}.
Inner automorphisms of $F$ act trivially on $H^4 (F;\mathbb{Z})$,
while changing the orientation of $X$ corresponds to multiplication by $-1$.
Thus the orbit of $k(X)$ under $Out(F)\times\{\pm1\}$ is the significant
invariant.

We may construct the third stage of the Postnikov tower for $X$ 
by adjoining cells of dimension greater than 4 to $X$.
The natural inclusion $j:X\to P_3(X)$ is then 4-connected.
If $X_1 $ is another such $PD_3$-complex and $\theta:\pi_1 (X_1)\to F$
is an isomorphism which identifies the $k$-invariants then there is 
a 4-connected map $j_1 :X_1 \to P_3 (X)$ inducing $\theta$, 
which is homotopic to a map with image in the 
4-skeleton of $P_3 (X)$, and so there is a map $h:X_1 \to X$ 
such that $j_1 $ is homotopic to $jh$. 
The map $h$ induces isomorphisms on $\pi_i $ for $i\leq 3$, 
since $j$ and $j_1 $ are 4-connected, and so the lift 
$\tilde h:\widetilde X_1\simeq S^3 \to\widetilde X\simeq S^3$
is a homotopy equivalence, by the theorems of Hurewicz and Whitehead.
Thus $h$ is itself a homotopy equivalence.
\end{proof}

The list of finite groups with cohomological period dividing 4 is well known. 
Each such group $F$ and generator $k\in H^4(F;\mathbb{Z})$
is realized by some $PD_3^+$-complex \cite{[Sw60],[Wl67]}. 
(See also Chapter 11 below.)
In particular, there is an unique homotopy type of $PD_3$-complexes 
with fundamental group $S_3$,
but there is no 3-manifold with this fundamental group 
\cite{[Mi57]}.

The fundamental group of a $PD_3$-complex $P$ has two ends if and only if 
$\widetilde P\simeq S^2$, and then 
$P$ is homotopy equivalent to one of the four 
$\mathbb{S}^2\times\mathbb{E}^1$-manifolds 
$S^2\times S^1 $, $S^2\tilde\times S^1$, $RP^2 \times S^1$ or 
$RP^3\sharp RP^3 $.
The following simple lemma leads to an alternative characterization
of $RP^2\times{S^1}$.

\begin{lemma}
Let $X$ be a finite-dimensional complex with a connected regular covering 
space $\widehat{X}$ and covering group $C=Aut(\widehat{X}/X)$.
If $\widetilde{H}_q(\widehat{X};\mathbb{Z})=0$ for $q\not=m$
then $H_{s+m+1}(C;\mathbb{Z})\cong H_s(C;H_m(\widehat{X};\mathbb{Z}))$, 
for all $s>>0$.
\end{lemma}
                
\begin{proof}
The lemma follows by devissage applied to the homology of 
$C_*(\widehat{X})$, considered as a chain complex over $\mathbb{Z}[C]$.
(In fact $s\geq{dim(X)}-m$ suffices.)
\end{proof}

\begin{theorem}
Let $P$ be a $PD_3$-space whose fundamental group $\pi$ has a nontrivial 
finite normal subgroup $N$. 
Then either $P$ is homotopy equivalent to $RP^2 \times S^1 $ or 
$\pi$ is finite.
\end{theorem}

\begin{proof}
We may clearly assume that $\pi$ is infinite.
Then $H_q(\widetilde P;\mathbb{Z})=0$ for $q>2$, by Poincar\'e duality.
Let $\Pi=\pi_2 (P)$.
The augmentation sequence 
\begin{equation*}
0\to A(\pi)\to \mathbb{Z}[\pi]\to\mathbb{Z}\to0
\end{equation*}
gives rise to a short exact sequence 
\begin{equation*}
0\to Hom_{\mathbb{Z}[\pi]} (\mathbb{Z}[\pi],\mathbb{Z}[\pi])\to 
Hom_{\mathbb{Z}[\pi]} (A(\pi),\mathbb{Z}[\pi])\to H^1 (\pi;\mathbb{Z}[\pi])
\to0.
\end{equation*}
Let $f:A(\pi)\to \mathbb{Z}[\pi]$ be a homomorphism and $\zeta$ be a central element of $\pi$.
Then $f.\zeta(i)=f(i)\zeta=\zeta f(i)=f(\zeta i)=f(i\zeta)$ and so
$(f.\zeta-f)(i)=f(i(\zeta-1))=if(\zeta-1)$ for all $i\in A(\pi)$.
Hence $f.\zeta-f$ is the restriction of a homomorphism from $\mathbb{Z}[\pi]$ 
to $\mathbb{Z}[\pi]$.
Thus central elements of $\pi$ act trivially on $H^1 (\pi;\mathbb{Z}[\pi])$.

If $n\in N$ the centraliser $\gamma=C_\pi (\langle n\rangle )$ 
has finite index in $\pi$, 
and so the covering space $P_\gamma $ is again a $PD_3 $-complex 
with universal covering space $\widetilde{P}$.
Therefore $\Pi\cong\overline{ H^1 (\gamma;\mathbb{Z}[\gamma])}$ as a (left) $\mathbb{Z}[\gamma]$-module.
In particular, $\Pi$ is a free abelian group.
Since $n$ is central in $\gamma$ it acts trivially on 
${H^1(\gamma;\mathbb{Z}[\gamma])}$ and hence via $w(n)$ on $\Pi$. 
Suppose first that $w(n)=1$.
Then Lemma 2.10 (with $X=P_{\langle{n}\rangle}$, 
$\widehat{X}=\widetilde{P}$ and $m=2$) 
gives an exact sequence 
\begin{equation*}
0\to Z/o(n)Z\to\Pi\to\Pi\to0,
\end{equation*}
where $o(n)$ is the order of $n$ and
the right hand homomorphism is multiplication by $o(n)$,
since $n$ acts trivially on $\Pi$. 
As $\Pi$ is torsion-free we must have $n=1$.

Therefore if $n\in N$ is nontrivial it has order 2 and $w(n)=-1$.
In this case Lemma 2.10 gives an exact sequence
\begin{equation*}
0\to\Pi\to\Pi\to Z/2Z\to 0,
\end{equation*}
where the left hand homomorphism is multiplication by 2.
Since $\Pi$ is a free abelian group it must be infinite cyclic.
Hence $\widetilde P\simeq S^2 $ and $\widetilde P/(Z/2Z)\simeq{RP^2}$.
The theorem now follows,
since any self homotopy equivalence of $RP^2$ is homotopic to the identity
(compare \cite[Theorem 4.4]{[Wl67]}). 
\end{proof}

If $P$ is any $PD_3$-complex and $C_\pi(\langle{g}\rangle)$ is infinite,
for some $g\in\pi$, then $g^2=1$, $w(g)=-1$
and $C_\pi(\langle{g}\rangle)$ has two ends \cite{[Cr00]}.
In fact, $C_\pi(\langle{g}\rangle)\cong\langle{g}\rangle\times\mathbb{Z}$ \cite{[Hi17']}.

If $\pi_1(P)$ has a finitely generated infinite normal subgroup of infinite 
index then it has one end, and so $P$ is aspherical. 
We shall discuss this case next.

\section{$PD_3 $-groups}

As a consequence of the work of Turaev and Crisp the study of $PD_3$-complexes
reduces largely to the study of $PD_3$-groups. 
It is not yet known whether all such groups are 3-manifold groups,
or even whether they must be finitely presentable.
The fundamental groups of aspherical 3-manifolds which are 
Seifert fibred or are finitely covered by surface bundles 
may be characterized among all $PD_3 $-groups in simple group-theoretic terms.

\begin{theorem}
Let $G$ be a $PD_3$-group with a nontrivial $FP_2$
normal subgroup $N$ of infinite index. Then either 
\begin{enumerate}   
\item $N\cong Z$ and $G/N$ is virtually a $PD_2$-group; or

\item $N$ is a $PD_2$-group and $G/N$ has two ends.
\end{enumerate}   
\end{theorem}

\begin{proof}
Let $e$ be the number of ends of $N$.
If $N$ is free then $H^3(G;\mathbb{Z}[G])\cong H^2(G/N;H^1(N;\mathbb{Z}[G]))$.
Since $N$ is finitely generated and $G/N$ is $FP_2$
this is in turn isomorphic to $H^2(G/N;\mathbb{Z}[G/N])^{(e-1)}$.
Since $G$ is a $PD_3$-group we must have $e-1=1$ and so $N\cong\mathbb{Z}$.
We then have 
$H^2(G/N;\mathbb{Z}[G/N])\cong{H^3(G;\mathbb{Z}[G])}\cong\mathbb{Z}^{w_1(G)}$.
Hence $G/N$ is virtually a $PD_2$-group, by Bowditch's Theorem.

Otherwise $c.d.N=2$ and so $e=1$ or $\infty$.
The LHSSS gives an isomorphism
$H^2(G;\mathbb{Z}[G])\cong H^1(G/N;\mathbb{Z}[G/N])\otimes
H^1(N;\mathbb{Z}[N])\cong H^1(G/N;\mathbb{Z}[G/N])^{e-1}$.
Hence either $e=1$ or $H^1(G/N;\mathbb{Z}[G/N])=0$.
But in the latter case we have $H^3(G;\mathbb{Z}[G])\cong 
H^2(G/N;\mathbb{Z}[G/N])\otimes H^1(N;\mathbb{Z}[N])$
and so $H^3(G;\mathbb{Z}[G])$ is either 0 or infinite dimensional.
Therefore $e=1$, and so 
$H^3(G;\mathbb{Z}[G])\cong H^1(G/N;\mathbb{Z}[G/N])\otimes H^2(N;\mathbb{Z}[N])$.
Hence $G/N$ has two ends and $H^2(N;\mathbb{Z}[N])$ $\cong\mathbb{Z}^{w_1(G)|_N}$,
so $N$ is a $PD_2$-group.
\end{proof}

We shall strengthen this result in Theorem 2.17 below.
         
\begin{cor}
A $PD_3$-space $P$ is homotopy equivalent to the mapping torus of a
self homeomorphism of a closed surface if and only if there is an epimorphism
$\phi:\pi_1(P)\to\mathbb{Z}$ with finitely generated kernel.
\end{cor}

\begin{proof}
This follows from Theorems 1.19, 2.11 and 2.12.
\end{proof}

If $\pi_1(P)$ is infinite and is a nontrivial direct product
then $P$ is homotopy equivalent to the product of $S^1$ with a closed surface.

\begin{theorem}
Let $G$ be a $PD_3$-group.
If $S$ is an almost coherent, restrained, locally virtually indicable subgroup
then $S$ is virtually solvable.
If $S$ has infinite index in $G$ it is virtually abelian.
\end{theorem}

\begin{proof} 
Suppose first that $S$ has finite index in $G$, and so is again a $PD_3$-group.
Since $S$ is virtually indicable we may assume without loss of generality that 
$\beta_1(S)>0$. 
Then $S$ is an ascending HNN extension $H*_\phi$ with finitely generated base.
Since $G$ is almost coherent $H$ is finitely presentable, 
and since $H^3(S;\mathbb{Z}[S])\cong\mathbb{Z}^{w_1(S)}$ 
it follows from \cite[Lemma 3.4]{[BG85]} that $H$ is normal in $S$ 
and $S/H\cong\mathbb{Z}$. 
Hence $H$ is a $PD_2$-group, by Theorem 2.12.
Since $H$ has no noncyclic free subgroup it is virtually $\mathbb{Z}^2$ 
and so $S$ and $G$ are virtually poly-$Z$.
                    
If $[G:S]=\infty$ then $c.d.S\leq2$, by Strebel's Theorem.
Let $J$ be a finitely generated subgroup of $S$.
Then $J$ is $FP_2$ and virtually indicable, and hence is virtually solvable, 
by Theorem 2.6 and its Corollary.
Since $J$ contains a $PD_2$-group \cite[Corollary 1.4]{[KK05]},
it is virtually abelian.
Hence $S$ is virtually abelian also.
\end{proof}

As the fundamental groups of virtually Haken 3-manifolds are 
coherent and locally virtually indicable, this implies the
Tits alternative for such groups \cite{[EJ73]}.
A slight modification of the argument gives the following corollary.      

\begin{cor}
A $PD_3$-group $G$ is virtually poly-$Z$ if and only if it is coherent,
restrained and has a subgroup of finite index with infinite abelianization.
\qed
\end{cor}

If $\beta_1(G)\geq2$ the hypothesis of coherence is redundant,
for there is then an epimorphism $p:G\to\mathbb{Z}$ 
with finitely generated kernel \cite[Theorem D]{[BNS87]}, 
and the kernel is then $FP_2$ by Theorem 1.19.
    
The argument of Theorem 2.13 and its corollary extend to show by induction on 
$m$ that a $PD_m$-group is virtually poly-$Z$ if and only if it is restrained
and every finitely generated subgroup is $FP_{m-1}$ and virtually indicable.

\begin{theorem}
Let $G$ be a $PD_3 $-group. Then $G$ is the fundamental group
of an aspherical Seifert fibred $3$-manifold or a 
$\mathbb{S}ol^3$-manifold
if and only if $\sqrt G\not=1$. 
Moreover 
\begin{enumerate}   
\item $h(\sqrt G)=1$ if and only if $G$ is the group of an
$\mathbb{H}^2\times\mathbb{E}^1$- or $\widetilde{\mathbb{SL}}$-manifold;

\item $h(\sqrt G)=2$ if and only if $G$ is the group of a 
$\mathbb{S}ol^3 $-manifold; 

\item $h(\sqrt G)=3$ if and only if $G$ is the group of an $\mathbb{E}^3 $- or 
$\mathbb{N}il^3 $-manifold.
\end{enumerate}   
\end{theorem}

\begin{proof}
The necessity of the conditions is clear.
(See \cite{[Sc83']}, or \S2 and \S3 of Chapter 7 below.)
Certainly $h(\sqrt G)\leq c.d.\sqrt G\leq 3$. 
Moreover $c.d.\sqrt G=3$ if and only if $[G:\sqrt G]$ is finite, 
by Strebel's Theorem.
Hence $G$ is virtually nilpotent if and only if $h(\sqrt G)=3$.
If $h(\sqrt G)=2$ then $\sqrt G$ is locally abelian, and hence abelian.
Moreover $\sqrt G$ must be finitely generated, for otherwise $c.d.\sqrt G=3$.
Thus $\sqrt G\cong\mathbb{Z}^2$ and case (2) follows from Theorem 2.12.

Suppose now that $h(\sqrt G)=1$ and let $C=C_G (\sqrt G)$. 
Then $\sqrt G$ is torsion-free abelian of rank 1, 
so $Aut(\sqrt{G})$ is isomorphic to a subgroup of $\mathbb{Q}^\times$. 
If $G/C$ is infinite then $c.d.C\leq2$, by Strebel's Theorem.
Moreover, $Aut(\sqrt{G})$ is infinite, so $\sqrt{G}\not\cong\mathbb{Z}$.
Therefore $C$ is abelian, by \cite[Theorem 8.8]{[Bi]},
and hence $G$ is solvable.
But then $h(\sqrt G)>1$, which is contrary to our hypothesis.
Therefore $G/C$ is isomorphic to a finite subgroup of 
$\mathbb{Q}^\times \cong\mathbb{Z}^\infty \oplus (Z/2Z)$ 
and so has order at most 2.
In particular, if $A$ is an infinite cyclic subgroup of $\sqrt G$
then $A$ is normal in $G$, and so $G/A$ is virtually a $PD_2$-group, 
by Theorem 2.12.
If $G/A$ is a $PD_2$-group then $G$ is the fundamental group of an $S^1$-bundle
over a closed surface.
In general, a finite torsion-free extension of the 
fundamental group of a closed Seifert fibred 3-manifold is again
the fundamental group of a closed Seifert fibred 3-manifold, 
by \cite{[Sc83]} and \cite[\S63]{[Zi]}.
\end{proof}

The heart of this result is the deep theorem of Bowditch.
The weaker characterization of fundamental groups 
of $\mathbb{S}ol^3$-manifolds and aspherical Seifert fibred 3-manifolds
as $PD_3$-groups $G$ such that $\sqrt G\not=1$ and $G$ 
has a subgroup of finite index with infinite abelianization 
is much easier to prove \cite{[H2]}.
There is as yet no comparable characterization of the groups of 
$\mathbb{H}^3 $-manifolds,
although it may be conjectured that these are exactly the $PD_3 $-groups 
with no noncyclic abelian subgroups.
(It has been recently shown that every closed $\mathbb{H}^3$-manifold 
is finitely covered by a mapping torus \cite{[Ag13]}.)

$\mathbb{N}il^3 $- and $\widetilde{\mathbb{SL}}$-manifolds are orientable,
and so their groups are $PD_3^+$-groups.
This can also be seen algebraically, 
as every such group has a characteristic subgroup $H$ 
which is a nonsplit central extension of a $PD_2^+$-group $\beta$ 
by $\mathbb{Z}$.
An automorphism of such a group $H$ must be orientation preserving.

Theorem 2.14 implies that if a $PD_3$-group $G$ is not virtually poly-$Z$
then its maximal elementary amenable normal subgroup is $\mathbb{Z}$ or 1.
For this subgroup is virtually solvable, by Theorem 1.11,
and if it is nontrivial then so is $\sqrt G$.

\begin{lemma}
Let $G$ be a group such that $c.d.G=2$ and let $K$ be an ascendant
$FP_2$ subgroup of $G$. Then either $[G:K]$ is finite or $K$ is free.
\end{lemma}

\begin{proof} 
We may assume that $K$ is not free, and so $c.d.K=c.d.G=2$.
Suppose first that $K$ is normal in $G$.
Then $G/K$ is locally finite \cite[Corollary 8.6]{[Bi]},
and so $G$ is the increasing union of a (possibly finite) sequence
of $FP_2$ subgroups $K=U_0<U_1<\dots $ such that $[U_{i+1}:U_i]$ is
finite, for all $i\geq0$.
It follows from the Kurosh subgroup theorem that if $U\leq V$ are finitely
generated groups and $[V:U]$ is finite then $V$ has strictly 
fewer indecomposable factors than $U$ unless both groups are indecomposable.
(See \cite[Lemma 1.4]{[Sc76]}).
Hence if $K$ is a nontrivial free product then $[G:K]$ is finite.
Otherwise $K$ has one end, and so $H^s(U_i;\mathbb{Z}[U_i])=0$ 
for $s\leq1$ and $i\geq0$.
Since $K$ is $FP_2$, the successive indices are finite and  
$c.d.U_i=2=c.d.G$ for all $i\geq0$ the union is finitely generated, 
by the Gildenhuys-Strebel Theorem.
Hence the sequence terminates and $[G:K]$ is again finite.

If $K=K_0<K_1<\dots <K_\beth=G$ is an ascendant sequence
then $[K_{\alpha+1}:K_\alpha]$ is finite for all $\alpha$,
by the argument just given.
Let $\omega$ be the union of the finite ordinals in $\beth$.
Then $\cup_{\alpha<\omega}K_\alpha$ is finitely generated, 
by the Gildenhuys-Strebel Theorem, and so $\omega$ is finite.
Hence the chain is finite, and so $[G:K]<\infty$.
\end{proof}

\begin{theorem}
Let $G$ be a $PD_3$-group with an ascending sequence of subgroups
$K_0<K_1<\dots$ such that $K_n$ is normal in $K_{n+1}$ for all $n\geq0$.
If $K=K_0$ is one-ended and $FP_2$ then the sequence is finite and
either $[K_n:K]$ or $[G:K_n]$ is finite, for all $n\geq0$.
\end{theorem}

\begin{proof}
Suppose that $[K_1:K]$ and $[G:K_1]$ are both infinite.
Since $K$ has one end it is not free and so $c.d.K=c.d.K_1=2$,
by Strebel's Theorem.
Hence there is a free $\mathbb{Z}[K_1]$-module $W$ such that
$H^2(K_1;W)\not=0$ \cite[Proposition 5.1]{[Bi]}.
Since $K$ is $FP_2$ and has one end $H^q(K;W)=0$ for $q=0$ or 1 and
$H^2(K;W)$ is an induced $\mathbb{Z}[K_1/K]$-module.
Since $[K_1:K]$ is infinite $H^0(K_1/K;H^2(K;W))=0$ \cite[Lemma 8.1]{[Bi]}.
The LHSSS for $K_1$ as an extension of $K_1/K$ by $K$ now gives
$H^r(K_1;W)=0$ for $r\leq2$, which is a contradiction.
A similar argument applies to the other terms of the sequence.

Suppose that $[K_n:K]$ is finite for all $n\geq0$ and let
$\widehat{K}=\cup_{n\geq0}K_n$.
If $c.d.\widehat{K}=2$ then $[\widehat{K}:K]<\infty$, by Lemma 2.15.
Thus the sequence must be finite.
\end{proof}

\begin{cor}
Let $G$ be a $PD_3$-group with an $FP_2$ subgroup $H$ 
which has one end and is of infinite index in $G$.
Let $H_0=H$ and $H_{i+1}=N_G(H_i)$ for $i\geq0$.
Then $\widehat H=\cup H_i$ is $FP_2$ and has one end, 
and either $c.d.\widehat H=2$ and $N_G(\widehat H)=\widehat H$ or
$[G:\widehat H]<\infty$ and $G$ is virtually the group of a surface bundle.
\end{cor}

\begin{proof}
This follows immediately from Theorems 2.12 and 2.16.
\end{proof}

\begin{cor}
If $G$ has a subgroup $H$ which is a $PD_2$-group with $\chi(H)=0$ 
(respectively, $<0$) then either it has such a subgroup which is its own
normalizer in $G$ or it is virtually the group of a surface bundle.
\end{cor}

\begin{proof}
If $c.d.\widehat H=2$ then $[\widehat H:H]<\infty$, 
so $\widehat H$ is a $PD_2$-group,
and $\chi(H)=[\widehat H:H]\chi(\widehat H)$.
\end{proof}

When $\chi(H)<0$ the corollary follows easily from 
the finite divisibility of $\chi(H)$, 
but something like Theorem 2.16 seems necessary when $\chi(H)=0$.

\begin{theorem}
Let $G$ be a $PD_3$-group with a nontrivial $FP_2$
subgroup $H$ which is ascendant and of infinite index in $G$.
Then either $H\cong\mathbb{Z}$ and $H$ is normal in $G$ or 
$G$ is virtually poly-$Z$ or $H$ is a $PD_2$-group, 
${[G:N_G(H)]<\infty}$ and $N_G(H)/H$ has two ends.
\end{theorem}

\begin{proof} 
Let $H=H_0<H_1<\dots<H_\beth=G$ be an ascendant sequence
and let $\gamma={\min\{ \alpha<\beth\mid [H_\alpha:H]=\infty\}}$.
Let $\widehat{H}=\cup_{\alpha<\gamma}H_\alpha$.
Then $h.d.\widehat{H}\leq2$ and so $[G:\widehat{H}]=\infty$.
Hence $c.d.\widehat{H}\leq2$ also, by Strebel's Theorem,
and so either $H$ is free or $[\widehat{H}:H]<\infty$, by Lemma 2.15.

If $H$ is not free then $c.d.\widehat{H}=2$ and
$\widehat{H}$ is $FP_2$, normal and of infinite index in $H_\gamma$.
Therefore $[G:H_\gamma]<\infty$ and so $H_\gamma$ is a $PD_3$-group, 
by Theorem 2.16.
Hence $\widehat{H}$ is a $PD_2$-group and $H_\gamma/\widehat{H}$ has two ends,
by Theorem 2.12.
Since $[\widehat{H}:H]<\infty$ it follows easily that $H$ is a $PD_2$-group, 
${[G:N_G(H)]<\infty}$ and $N_G(H)/H$ has two ends.

If $H\cong{F(r)}$ for some $r>1$ then $\gamma$ and $[\widehat{H}:H]$ are finite,
since $[H_n:H]$ divides $\chi(H)=1-r$ for all $n<\gamma$.
A similar argument shows that $H_\gamma/\widehat{H}$ is not locally finite.
Let $K$ be a finitely generated subgroup of $H_\gamma$ which contains 
$\widehat{H}$ as a subgroup of infinite index.
Then $K/\widehat{H}$ is virtually free \cite[Theorem 8.4]{[Bi]}, 
and so $K$ is finitely presentable.
In particular, $\chi(K)=\chi(\widehat{H})\chi(K/\widehat{H})$.
Now $\chi(K)\leq0$ \cite[\S9]{[KK05]}.
Since $\chi(\widehat{H})<0$ this is only possible if $\chi(K/\widehat{H})\geq0$,
and so $K/\widehat{H}$ is virtually $\mathbb{Z}$.
Hence we may assume that $H_\gamma$ is the union of an increasing
sequence $N_0=H<N_1\leq\dots$ of finitely generated
subgroups with $N_i/H$ virtually $\mathbb{Z}$, for $i\geq1$. 
For each $i\geq1$ the group $N_i$ is $FP_2$, $c.d.N_i=2$,
$H^s(N_i;\mathbb{Z}[N_i])=0$ for $s\leq1$ and $[N_{i+1}:N_i]$ is finite.
Therefore $H_\gamma$ is finitely generated, 
by the Gildenhuys-Strebel Theorem.

In particular, $H_\gamma$ is virtually a semidirect product
$\widehat{H}\rtimes\mathbb{Z}$, and so it is $FP_2$ and $c.d.H_\gamma=2$.
Hence $H_\gamma$ is a $PD_2$-group, by the earlier argument.
But $PD_2$-groups do not have normal subgroups such as $\widehat{H}$.
Therefore if $H$ is free it is infinite cyclic: $H\cong\mathbb{Z}$.
Since $\sqrt{H_\alpha}$ is characteristic in $H_\alpha$ it is normal in
$H_{\alpha+1}$, for each $\alpha<\beth$.
Transfinite induction now shows that $H\leq\sqrt G$.
Therefore either $\sqrt G\cong\mathbb{Z}$, 
so $H\cong\mathbb{Z}$ and is normal in $G$, 
or $G$ is virtually poly-$Z$, by Theorem 2.14.
\end{proof}

If $H$ is a $PD_2$-group $N_G(H)$ is the fundamental group of a 3-manifold 
which is double covered by the mapping torus of a surface homeomorphism.
There are however $\mathbb{N}il^3$-manifolds whose groups have no normal 
$PD_2$-subgroup (although they always have subnormal copies of 
$\mathbb{Z}^2$).

The original version of this result assumed that $H$ is subnormal in $G$.
(See \cite{[BH91]} for a proof not using \cite{[Bo04]} or \cite{[KK05]}.)

\section{Subgroups of $PD_3$-groups and 3-manifold groups}

The central role played by incompressible surfaces in the geometric study 
of Haken 3-manifolds suggests strongly the importance of 
studying subgroups of infinite index in $PD_3$-groups.
Such subgroups have cohomological dimension $\leq2$, by Strebel's Theorem.

There are substantial constraints on 3-manifold groups and their subgroups.
Every finitely generated subgroup of a 3-manifold group
is the fundamental group of a compact 3-manifold (possibly with boundary),
by Scott's Core Theorem \cite{[Sc73]}, 
and thus is finitely presentable and is either a 3-manifold group 
or has finite geometric dimension 2 or is a free group. 
Aspherical closed 3-manifolds are Haken, 
hyperbolic or Seifert fibred, by the work of Perelman \cite{[B-P]}.
The groups of such 3-manifolds are residually finite \cite{[He87]},
and the centralizer of any element in the group is finitely generated 
\cite{[JS]}.
Solvable subgroups of such groups are virtually poly-$Z$ 
\cite{[EJ73]}.

In contrast, any group of finite geometric dimension 2 is the
fundamental group of a compact aspherical 4-manifold with boundary, 
obtained by attaching 1- and 2-handles to $D^4$.
On applying the reflection group trick of Davis \cite{[Da83]}
to the boundary we see that each such group embeds in a $PD_4 $-group. 
For instance, the product of two nonabelian $PD_2^+$-groups contains a copy of 
$F(2)\times F(2)$, and so is a $PD_4^+$-group which is not almost coherent.
No $PD_4$-group containing a Baumslag-Solitar group 
$\langle x,t\mid tx^pt^{-1}=x^q\rangle$ is residually finite, 
since this property is inherited by subgroups.
Thus the question of which groups of finite geometric dimension 2 are 
subgroups of $PD_3$-groups is critical.

Kapovich and Kleiner have given an algebraic Core Theorem,
showing that every one-ended $FP_2$ subgroup $H$ in a $PD_3$-group $G$ 
is the ``ambient group" of a $PD_3$-pair $(H,\mathcal{S})$ 
\cite{[KK05]}.
Using this the argument of \cite{[Kr90a]} may be adapted to show that
every strictly increasing sequence of centralizers in $G$ 
has length at most 4 \cite{[Hi06]}.
(The finiteness of such sequences and the fact that centralizers in $G$ 
are finitely generated or rank 1 abelian are due to Castel \cite{[Ca07]}.)
With the earlier work of Kropholler and Roller 
\cite{[KR88],[KR89],[Kr90],[Kr93]} 
it follows that if $G$ has a subgroup $H\cong\mathbb{Z}^2$ 
and $\sqrt G=1$ then it splits over a subgroup commensurate with $H$. 
It also follows easily from the algebraic Core Theorem that 
if a subgroup $H$ is an $\mathcal{ X}$-group then $H=\pi_1(N)$
for some Seifert fibred 3-manifold $N$ with $\partial{N}\not=\emptyset$.
In particular, no nontrivial Baumslag-Solitar relation holds in $G$ 
\cite{[Ca07]}.

The geometric conclusions of Theorem 2.14 and the coherence of 3-manifold 
groups suggest that Theorems 2.12 and 2.17 should hold under the weaker 
hypothesis that $N$ be finitely generated. 
(Compare Theorem 1.19.)
It is known that $F(2)\times F(2)$ is not a subgroup of any $PD_3 $-group 
\cite{[KR89]}.
This may be regarded as a weak coherence result.

Is there a characterization of  virtual $PD_3$-groups parallel to 
Bowditch's Theorem? (It may be relevant that
homology $n$-manifolds are manifolds for $n\leq2$.
There is no direct analogue in high dimensions.
For every $k\geq6$ there are $FP_k$ groups $G$ with 
$H^k(G;\mathbb{Z}[G])\cong\mathbb{Z}$ 
but which are not virtually torsion-free \cite{[FS93]}.)

\section{$\pi_2(P)$ as a $\mathbb{Z}[\pi]$-module}

Let $P$ be a $PD_3$-space with fundamental group $\pi$ and orientation
character $w$.
If $\pi$ is finite then $w=0$, $\pi_2(P)=0$ and $c_{P*}[P]\in{H_3(\pi;\mathbb{Z})}$ 
is essentially equivalent to the first nontrivial $k$-invariant of $P$, 
as outlined in Theorem 2.9.
Suppose that $\pi$ is infinite.
If $N$ is another $PD_3$-space and there is an isomorphism 
$\theta:\nu=\pi_1(N)\to\pi$ such that $w_1(N)=\theta^*w$
then $\pi_2(N)\cong\theta^*\pi_2(P)$ as $\mathbb{Z}[\nu]$-modules.
If moreover $k_1(N)=\theta^*k_1(P)$ (modulo automorphisms of the pair 
$(\nu,\pi_2(N))$) then $P_2(N)\simeq{P_2(P)}$.
Since we may construct these Postnikov 2-stages by adjoining cells of dimension
$\geq4$ it follows that there is a map $f:N\to P$ such that $\pi_1(f)=\theta$
and $\pi_2(f)$ is an isomorphism.
The homology of the universal covering spaces 
$\widetilde{N}$ and $\widetilde{P}$ is 0 above degree 2, 
and so $f$ is a homotopy equivalence, by the Whitehead Theorem.
Thus the homotopy type of $P$ is determined by the triple $(\pi,w,k_1(P))$.
One may ask how $c_{P*}[P]$ and $k_1(P)$ determine each other.

There is a facile answer:
in Turaev's realization theorem for homotopy triples
the element of $H_3(\pi;\mathbb{Z}^w)$ is used to construct a cell complex $X$ 
by attaching 2- and 3-cells to the 2-skeleton of $K(\pi,1)$.
If $C_*$ is the cellular chain complex of $\widetilde X$ then $k_1(X)$ is 
the class of
\begin{equation*}
\begin{CD}
0\to\pi_2(X)\to{C_2/\partial{C_3}}\to{C_1}\to{C_0}\to\mathbb{Z}\to0
\end{CD}
\end{equation*}
in $H^3(\pi;\pi_2(X))=Ext^3_{\mathbb{Z}[\pi]}(\mathbb{Z},\pi_2(X))$.
Conversely, a class $\kappa\in Ext^3_{\mathbb{Z}[\pi]}(\mathbb{Z},\Pi)$
corresponds to an extension 
\begin{equation*}
\begin{CD}
0\to\Pi\to{D_2}\to{D_1}\to{D_0}\to\mathbb{Z}\to0,
\end{CD}
\end{equation*}
with $D_1$ and $D_0$ finitely generated free $\mathbb{Z}[\pi]$-modules.
Let $\mathcal{D}_*$ be the complex ${D_2}\to{D_1}\to{D_0}$,
with augmentation $\varepsilon$ to $\mathbb{Z}$.
If $\kappa=k_1(P)$ for a $PD_3$-complex $P$ then 
$\mathbb{T}or_3^{\mathbb{Z}[\pi]}(\mathbb{Z}^w,\mathcal{D}_*)\cong
H_3(P_2(P);\mathbb{Z}^w)\cong\mathbb{Z}$
(where $\mathbb{T}or$ denotes hyperhomology), and the augmentation then
determines a class in $H_3(\pi;\mathbb{Z}^w)$ (up to sign).
Can these connections be made more explicit?
Is there a natural homomorphism from 
$H^3(\pi;\overline{H^1(\pi;\mathbb{Z}[\pi])})$ to $H_3(\pi;\mathbb{Z}^w)$?

If $P$ is an orientable 3-manifold which is the connected sum of 
a 3-manifold whose fundamental group is free of rank $r$ with $s\geq 1$ 
aspherical 3-manifolds then $\pi_2 (P)$ is a finitely generated free 
$\mathbb{Z}[\pi]$-module of rank $r+s-1$ \cite{[Sw73]}.
We shall give a direct homological argument that applies
for $PD_3 $-spaces with torsion-free fundamental group, and we shall also
compute $H^2 (P;\pi_2(P))$ for such spaces.
(This cohomology group arises in studying homotopy classes 
of self homotopy equivalences of $P$ \cite{[HL74]}.) 

\begin{theorem}
Let $P$ be a $PD_3 $-space with torsion-free fundamental group $\pi$
and orientation character $w=w_1(P)$.
Then
\begin{enumerate}                                                                                   
\item if $\pi$ is a nontrivial free group $\pi_2(P)$ is finitely generated 
and of projective dimension $1$ as a left $\mathbb{Z}[\pi]$-module
and $H^2 (P;\pi_2(P))$ is infinite cyclic;

\item if $\pi$ is not free $\pi_2(P)$ is a finitely generated free 
$\mathbb{Z}[\pi]$-module, $c.d.\pi=3$,
$H_3(c_P;\mathbb{Z}^w)$ is a monomorphism and $H^2 (P;\pi_2(P))=0$;

\item  $P$ is homotopy equivalent to a finite $PD_3$-complex if and only if 
$\pi$ is finitely presentable and $FF$.
\end{enumerate}   
\end{theorem}

\begin{proof}
Two applications of Poincar\'e duality give,
firstly,
$\pi_2(P)\cong\overline{H^1(\pi;\mathbb{Z}[\pi])}$
and then $H^2 (P;\pi_2(P))\cong{H_1(P;H^1(\pi;\mathbb{Z}[\pi]))=
H_1(\pi;H^1(\pi;\mathbb{Z}[\pi]))}$.
Since $\pi$ is $FP_2$ it is accessible, and so $\pi\cong\pi\mathcal{G}$,
where $\mathcal{G}$ is a finite graph of groups with all vertex groups
finite or one-ended and all edge groups finite \cite[Theorem VI.6.3]{[DD]}.
There is an associated Mayer-Vietoris presentation 
\[
0\to\oplus\mathbb{Z}[G_v\backslash\pi]\to\oplus\mathbb{Z}[G_e\backslash\pi]
\to{H^1(\pi;\mathbb{Z}[\pi])}\to0,
\]
where the sums involve only the {\it finite\/} vertex groups $G_v$ 
(and edge groups $G_e$) \cite[Theorem 2]{[Ch76]}.
If $\pi$ is free of rank $r>0$ we may assume there is one vertex,
with trivial vertex group, and $r$ edges.
The above presentation is then
\[
0\to\mathbb{Z}[\pi]\to\mathbb{Z}[\pi]^r\to{H^1(\pi;\mathbb{Z}[\pi])}\to0.
\]
On applying the functor $-\otimes_{\mathbb{Z}[\pi]}\mathbb{Z}$, 
the left hand homomorphism becomes the trivial homomorphism from
$\mathbb{Z}\to\mathbb{Z}^r$.
Hence
\[
H^2(P;\pi_2(P))\cong{H_1 (\pi;H^1 (\pi;\mathbb{Z}[\pi]))=
Tor_1^{\mathbb{Z}[\pi]}(H^1 (\pi;\mathbb{Z}[\pi]),\mathbb{Z})}
\cong\mathbb{Z},
\]
by the exact sequence of $Tor$.
Moreover $\pi_2(P)$ has projective dimension 1.
As $\pi$ is finitely presentable and projective 
$\mathbb{Z}[F(r)]$-modules are free \cite{[Ba64]}, 
$P$ is homotopy equivalent to a finite $PD_3$-complex \cite{[Wl65]}.
(In fact $P$ is homotopy equivalent to a connected sum of copies of
$S^2\times{S^1}$ and $S^2\tilde\times{S^1}$.)

If $\pi$ is torsion-free but not free then we may assume that
the vertex groups are finitely generated
and have one end, and the edge groups are trivial.
Hence $H^1(\pi;\mathbb{Z}[\pi])$ is a free right $\mathbb{Z}[\pi]$-module
with basis corresponding to the edges of $\mathcal{G}$,
and so $H^2 (P;\pi_2 (P))=0$.
We may assume that $P$ is 3-dimensional and
$C_*(\widetilde{P})$ is chain homotopy 
equivalent to a finitely generated projective $\mathbb{Z}[\pi]$-complex 
\begin{equation*}
0\to C_3 \to C_2 \to C_1 \to C_0 \to 0,
\end{equation*} 
where $C_i$ is free if $i\leq2$. 
Let $Z_2 $ be the module of 2-cycles.
Then the sequences 
\begin{gather*}
0\to Z_2 \to C_2 \to C_1 \to C_0 \to\mathbb{Z}\to 0\\
\text{and}\qquad\qquad 0\to C_3 \to Z_2 \to\pi_2 (P)\to 0
\qquad\qquad \phantom{\text{and}} 
\end{gather*}
are exact, since $H_3(\widetilde{P};\mathbb{Z})=0$.
Attaching 3-cells to $P$ along a basis for $\pi_2(P)$ gives
an aspherical 3-dimensional complex $K$ with fundamental group $\pi$.
The inclusion of $P$ into $K$ may be identified with $c_P$, and
clearly induces monomorphisms
$H_3(P;A)\to H_3(\pi;A)$ for any coefficient module $A$.
Hence $c.d.\pi=3$.

If $\pi$ is $FF$ there is a finite free resolution
\begin{equation*}
0\to D_3 \to D_2 \to D_1 \to D_0 \to\mathbb{Z}\to0.
\end{equation*} 
Therefore $Z_2$ is finitely generated and stably free, by Schanuel's Lemma.
Since $\pi_2(P)$ is free $Z_2 \cong\pi_2(P)\oplus C_3 $ and so $C_3$ is also
stably free. 
Hence if moreover $\pi$ is finitely presentable 
then $P$ is homotopy equivalent to a finite $PD_3$-complex.
The converse is clear, by the above construction of $K(\pi,1)\simeq{K}$. 
\end{proof}

We may remove the condition that $\pi$ be torsion-free.

\begin{cor}
If $P$ is a $PD_3$-space then $\pi_2(P)$ is finitely presentable 
as a $\mathbb{Z}[\pi]$-module.
Moreover, $H^2(P;\pi_2(P))$ is a finitely generated abelian group 
of rank $1$, if $\pi$ is infinite and virtually free, 
and is finite otherwise.
If $\pi$ is infinite but not torsion-free 
the projective dimension of $\pi_2(P)$ is infinite. 
\end{cor}

\begin{proof}
The first assertion follows from the theorem,
since $\pi$ is virtually torsion-free, by Crisp's Theorem.
The second follows easily from the Mayer-Vietoris
presentation for $H^1(\pi;\mathbb{Z}[\pi])$.
If $\pi$ is infinite and $\pi_2(P)$ has finite projective dimension
then so does $Z_2$, and so $c.d.\pi<\infty$, 
and then $\pi$ is torsion-free.
\end{proof}

Crisp uses an ingenious combinatorial argument based on the Mayer-Vietoris
presentation for $H^1(\pi;\mathbb{Z}[\pi])$
together with Lemma 2.10 to show that if $P$ is indecomposable, 
orientable and not aspherical the vertex groups must all be finite, 
and so $\pi$ is virtually free.
Elementary group theory then leads to the near-determination of the groups 
of such $PD_3$-complexes \cite{[Hi12]}.
(It is not yet clear what are the indecomposable {\it non}-orientable $PD_3$-complexes.)

%% file: m5-3.tex
\chapter{Homotopy invariants of $PD_4$-complexes}

The homotopy type of a $4$-manifold $M$ is largely determined 
(through Poincar\'e duality) by 
its algebraic 2-type and orientation character. 
In many cases the formally weaker invariants $\pi_1(M)$, 
$w_1(M)$ and $\chi(M)$ already suffice.
In \S1 we give criteria in such terms for a degree-1 map
between  $PD_4$-complexes to be a homotopy equivalence, 
and for a $PD_4$-complex to be aspherical.
We then show in \S2 that if the universal covering space of a $PD_4$-complex
is homotopy equivalent to a finite complex then it is either compact, contractible, 
or homotopy equivalent to $S^2$ or $S^3$.                             
In \S3 we obtain estimates for the minimal Euler characteristic
of $PD_4$-complexes with fundamental group of cohomological dimension at most 2
and determine the second homotopy groups of $PD_4$-complexes realizing 
the minimal value.
The class of such groups includes all surface groups and classical link groups, 
and the groups of many other (bounded) 3-manifolds.
The minima are realized by $s$-parallelizable PL 4-manifolds.
In \S4 we show that if $\chi(M)=0$ then $\pi_1(M)$ satisfies some stringent 
constraints, and in the final section we define
the reduced intersection pairing.

\section{Homotopy equivalence and asphericity}

Many of the results of this section depend on the following lemma, 
in conjunction with 
use of the Euler characteristic to compute the rank of the surgery kernel.
(Lemma 3.1 and Theorem 3.2 derive from \cite[ Lemmas 2.2 and 2.3]{[Wl]}.)

\begin{lemma}
Let $R$ be a ring and $C_* $ be a finite chain complex 
of projective $R$-modules. If $H_i(C_*)=0$ for $i<q$ and 
$H^{q+1} (Hom_R (C_* ,B))=0$ for any left $R$-module $B$ 
then $H_q (C_* )$ is projective. 
If moreover $H_i(C_*)=0 $ for $i>q$ then 
$H_q (C_* )\oplus\bigoplus_{i\equiv q+1~ (2)} C_i \cong 
\bigoplus_{i\equiv q~ (2)} C_i $.
\end{lemma}

\begin{proof}
We may assume without loss of generality that $q=0$ and $C_i=0$ for $i<0$.
We may factor $\partial_1 :C_1 \to C_0$ through $B=\mathrm{Im}\partial_1 $ 
as $\partial_1=j\beta$,
where $\beta$ is an epimorphism and $j$ is the natural inclusion of the 
submodule $B$.
Since $j\beta\partial_2=\partial_1\partial_2=0$ and $j$ is injective 
$\beta\partial_2=0$.
Hence $\beta $ is a 1-cocycle of the complex $Hom_R (C_* ,B)$.
Since $H^1 (Hom_R (C_* ,B))=0$ there is a homomorphism $\sigma:C_0\to B$ 
such that $\beta=\sigma\partial_1=\sigma j\beta$. 
Since $\beta$ is an epimorphism $\sigma j=id_B $ and so $B$ is a 
direct summand of $C_0$.
This proves the first assertion.

The second assertion follows by an induction on the length of the complex. 
\end{proof}

\begin{theorem} Let $M$ and $N$ be finite $PD_4$-complexes.
A map $f:M\to N$ is a homotopy equivalence if and only if $\pi_1(f)$ 
is an isomorphism,
$f^*w_1 (N)=w_1 (M)$, $f_* [M]=\pm [N]$ and $\chi(M)=\chi(N)$.
\end{theorem}

\begin{proof} The conditions are clearly necessary. Suppose that they hold. 
Up to homotopy type we may assume that $f$ is a cellular inclusion of finite 
cell complexes, 
and so $M$ is a subcomplex of $N$. 
We may also identify $\pi_1 (M)$ with $\pi=\pi_1 (N)$.
Let $C_* (M)$, $C_* (N)$ and $D_* $ be the cellular chain complexes of 
$\widetilde M$, $\widetilde N$ and $(\widetilde N,\widetilde M)$, 
respectively.
Then the sequence 
\begin{equation*}
0\to C_*(M)\to C_* (N)\to D_* \to 0
\end{equation*}
is a short exact sequence of finitely generated free 
$\mathbb{Z}[\pi]$-chain complexes.

By the projection formula $f_* (f^* a\cap[M])=a\cap f_* [M]=\pm a\cap [N]$ for
any cohomology class $a\in H^* (N;\mathbb{Z}[\pi])$.
Since $M$ and $N$ satisfy Poincar\'e duality it follows that $f$ induces
split surjections on homology and split injections on cohomology.
Hence $H_q (D_* )$ is the ``surgery kernel" in degree $q-1$, and the duality 
isomorphisms
induce isomorphisms from $H^r (Hom_{\mathbb{Z}[\pi]} (D_*,B))$ to 
$H_{6-r} (\overline{D_*}\otimes B)$,
where $B$ is any left $\mathbb{Z}[\pi]$-module.
Since $f$ induces isomorphisms on homology and cohomology in degrees $\leq 1$, 
with any coefficients, the hypotheses of Lemma 3.1 are satisfied for the 
$\mathbb{Z}[\pi]$-chain complex $D_* $, with $q=3$, 
and so $H_3 (D_*)=\mathrm{Ker}(\pi_2 (f))$ is projective. 
Moreover $H_3 (D_*)\oplus \bigoplus_{i~ odd} D_i \cong\bigoplus_{i~ even} D_i $.
Thus $H_3 (D_* )$ is a stably free $\mathbb{Z}[\pi]$-module of rank 
$\chi(E,M)=\chi(M)-\chi(E)=0$.
Hence $H_3 (D_* )=0$, since group rings are weakly finite,
and so $f$ is a homotopy equivalence. 
\end{proof}

If $M$ and $N$ are merely finitely dominated, rather than finite, then
$H_3 (D_* )$ is a finitely generated projective $\mathbb{Z}[\pi]$-module
such that $\mathbb{Z}\otimes_{\mathbb{Z}[\pi]}H_3 (D_*)=0$.
If the Wall finiteness obstructions satisfy $f_*\sigma(M)=\sigma(N)$ in 
$\tilde K_0(\mathbb{Z}[\pi])$ then $H_3 (D_* )$ is stably free, and the
theorem remains true.
The theorem holds as stated for arbitrary $PD_4$-spaces if $\pi$ satisfies 
the Weak Bass Conjecture. 
(Similar comments apply elsewhere in this section.)

We shall see that when $N$ is aspherical and $f=c_M$ we may drop
the hypotheses that $f^*w_1 (N)=w_1 (M)$ and $f$ has degree $\pm1$.

\begin{cor}
{\rm[Ha87]}
Let $N$ be orientable. Then a map $f:N\to N$ which induces automorphisms 
of $\pi_1(N)$ and $H_4(N;\mathbb{Z})$ is a homotopy equivalence.
\qed
\end{cor}

Any self-map of a geometric manifold of semisimple type
(e.g., an $\mathbb{H}^4$-, $\mathbb{H}^2(\mathbb{C})$- or 
$\mathbb{H}^2\times\mathbb{H}^2$-manifold) 
with nonzero degree is a homotopy equivalence \cite{[Re96]}.
      
If $X$ is a cell complex with fundamental group $\pi$ then
$\pi_2(X)\cong{H_2(X;\mathbb{Z}[\pi])}$, by the Hurewicz Theorem for
$\widetilde{X}$, and so there is an {\it evaluation} homomorphism
$ev:H^2(X;\mathbb{Z}[\pi])\to{Hom_{\mathbb{Z}[\pi]}(\pi_2(X),\mathbb{Z}[\pi])}$.
The latter module may be identified with 
$H^0(\pi;H^2(\widetilde X;\mathbb{Z}[\pi]))$,
the $\pi$-invariant subgroup of the cohomology of $\widetilde{X}$
with coefficients $\mathbb{Z}[\pi]$.

\begin{lemma}
Let $M$ be a $PD_4 $-space with fundamental group $\pi$
and let $\Pi=\pi_2(M)$.
Then $\Pi\cong\overline{H^2(M;\mathbb{Z}[\pi])}$ and
there is an exact sequence
\begin{equation*}
\begin{CD}
0\to H^2(\pi;\mathbb{Z}[\pi])\to{H^2(M;\mathbb{Z}[\pi])}@> ev >>
Hom_{\mathbb{Z}[\pi]}(\Pi,\mathbb{Z}[\pi])\to 
{H^3(\pi;\mathbb{Z}[\pi])}\to 0.
\end{CD}
\end{equation*}
\end{lemma}

\begin{proof} 
This follows from the Hurewicz Theorem, Poincar\'e duality and the UCSS,
since $H^3(M;\mathbb{Z}[\pi])\cong H_1(\widetilde M;\mathbb{Z})=0$.
\end{proof}

Exactness of much of this sequence can be derived without the UCSS.
When $\pi$ is finite the sequence reduces to the isomorphism 
$\pi_2(M)\cong\overline{Hom_{\mathbb{Z}[\pi]}(\pi_2(M),\mathbb{Z}[\pi])}$.

Let $ev^{(2)}:H^2_{(2)}(\widetilde M)\to 
{Hom_{\mathbb{Z}[\pi]}(\pi_2(M),\ell^2(\pi))}$ 
be the analogous evaluation defined on the {\it un}reduced $L^2$-cohomology
by $ev^{(2)}(f)(z)=\Sigma f(g^{-1}z)g$ 
for all square summable 2-cocycles $f$ and all 2-cycles $z$ representing 
elements of ${H_2(X;\mathbb{Z}[\pi])}\cong\pi_2(M)$.
Part of the next theorem is implicit in \cite{[Ec94]}.

\begin{theorem} 
Let $M$ be a $PD_4$-complex with fundamental group $\pi$.
Then
\begin{enumerate}

\item if $\beta^{(2)}_1(\pi)=0$ and either $M$ is finite 
or $\pi$ satisfies 
the Strong Bass Conjecture then $\chi(M)\geq0$;

\item $\mathrm{Ker}(ev^{(2)})$ is closed;

\item if $\beta^{(2)}_2(M)=\beta^{(2)}_2(\pi)$ then 
$H^2(c_M;\mathbb{Z}[\pi]):H^2 (\pi;\mathbb{Z}[\pi])\to H^2(M;\mathbb{Z}[\pi])$
is an isomorphism.
\end{enumerate}
\end{theorem} 

\begin{proof}
Since $M$ is a $PD_4$-complex $\beta_i^{(2)}(M)=\beta_{4-i}^{(2)}(M)$,
for all $i$.
If $M$ is finite or $\pi$ satisfies 
the Strong Bass Conjecture then $\chi(M)$ is
the alternating sum of the $L^2$-Betti numbers \cite{[Ec96]}.
Therefore if, moreover, $\beta_1^{(2)}(\pi)=0$ then
$\chi(M)=2\beta_0^{(2)}+\beta_2^{(2)}(M)\geq0$.
 
Let $z\in C_2(\widetilde M)$ be a 2-cycle and $f\in C_2^{(2)}(\widetilde M)$
a square-summable 2-cocycle.
As $||ev^{(2)}(f)(z)||_2\leq ||f||_2||z||_2$, 
the map $f\mapsto ev^{(2)}(f)(z)$ is continuous, for fixed $z$.
Hence if $f=lim f_n$ and $ev^{(2)}(f_n)=0$ for all $n$ then $ev^{(2)}(f)=0$.

The inclusion $\mathbb{Z}[\pi]<\ell^2(\pi)$ induces a homomorphism 
from the exact sequence of Lemma 3.3 to the corresponding sequence with
coefficients $\ell^2(\pi)$.
(See \cite[\S1.4]{[Ec94]}.
Note that we may identify $H^0(\pi;H^2(\widetilde{M};A))$ with 
$Hom_{\mathbb{Z}[\pi]}(\pi_2(M),A)$ 
for $A=\mathbb{Z}[\pi]$ or $\ell^2(\pi)$
since $\widetilde{M}$ is 1-connected.)
As $\mathrm{Ker}(ev^{(2)})$ is closed and
$ev^{(2)}(\delta g)(z)=ev^{(2)}(g)(\partial z)=0$ for any square summable
1-chain $g$, the homomorphism  $ev^{(2)}$ factors through 
the reduced $L^2$-cohomology $\bar H^2_{(2)}(\widetilde M)$.
If $\beta^{(2)}_2(M)=\beta^{(2)}_2(\pi)$
the classifying map $c_M:M\to K(\pi,1)$ induces weak isomorphisms on 
reduced $L^2$-cohomology 
$\bar H^i_{(2)} (\pi)\to\bar H ^i_{(2)} (\widetilde M)$ 
for $i\leq2$.
In particular, the image of $\bar{H}^2_{(2)}(\pi)$ is dense
in $\bar{H}^2_{(2)}(\widetilde{M})$.
Since $ev^{(2)}$ is trivial on $\bar{H}^2_{(2)}(\pi)$
and $\mathrm{Ker}(ev^{(2)})$ is closed it follows that $ev^{(2)}=0$.
Since the natural homomorphism from 
$Hom_{\mathbb{Z}[\pi]}(\pi_2(M),\mathbb{Z}[\pi])$ 
to $Hom_{\mathbb{Z}[\pi]}(\pi_2(M),\ell^2(\pi))$
is a monomorphism it follows that $ev=0$ also and 
so $H^2(c_M;\mathbb{Z}[\pi])$ is an isomorphism. 
\end{proof}

This gives a complete and natural criterion for asphericity
(which we state as a separate theorem to retain the enumeration of  
the original version of this book).

\begin{theorem} 
Let $M$ be a $PD_4$-complex with fundamental group $\pi$.
Then $M$ is aspherical if and only if $H^s(\pi;\mathbb{Z}[\pi])=0$ for $s\leq2$ 
and $\beta^{(2)}_2(M)=\beta^{(2)}_2(\pi)$.
\end{theorem}

\begin{proof} The conditions are clearly necessary.
If they hold then
$H^2(M;\mathbb{Z}[\pi])\cong{H}^2(\pi;\mathbb{Z}[\pi])=0$
and so $M$ is aspherical, by Poincar\'e duality. 
\end{proof}

Is it possible to replace the hypothesis
``$\beta^{(2)}_2(M)=\beta^{(2)}_2(\pi)$" 
by ``$\beta_2(M^+)=\beta_2(\mathrm{Ker}(w_1(M)))$", 
where $p_+:M^+\to M$ is the orientation cover?
It is easy to find examples to show that the homological 
conditions on $\pi$ cannot be relaxed further.

\begin{cor}
The $PD_4$-complex $M$ is finite and aspherical if and only if $\pi$ 
is a finitely presentable $PD_4$-group of type $FF$ and $\chi(M)=\chi(\pi)$.
\qed
\end{cor}

If $\beta_2(\pi)\not=0$ this follows from Theorem 3.2.
For we may assume $\pi$ and $M$ are orientable,
on replacing $\pi$ by
$K=\mathrm{Ker}(w_1(M))\cap\mathrm{Ker}(w_1(\pi))$ and $M$ by $M_K$.
As $H_2(c_M;\mathbb{Z})$ is onto it is an isomorphism, 
so $c_M$ has degree $\pm1$, by Poincar\'e duality.
Is $M$ always aspherical if $\pi$ is a $PD_4$-group and $\chi(M)=\chi(\pi)$?

\begin{cor}
If $\chi(M)=\beta_1^{(2)}(\pi)=0$ and 
$H^s(\pi;\mathbb{Z}[\pi])=0$ for $s\leq 2$ 
then $M$ is aspherical and $\pi$ is a $PD_4$-group.
\qed
\end{cor}

\begin{cor}
If $\pi\cong\mathbb{Z}^r$ then $\chi(M)\geq 0$, 
and is $0$ only if $r=1$, $2$ or $4$.
\end{cor}
                       
\begin{proof}
If $r>2$ then $H^s(\pi;\mathbb{Z}[\pi])=0$ for $s\leq 2$.
\end{proof}

Theorem 3.5 implies that if $\pi$ is a $PD_4$-group 
and $\chi(M)=\chi(\pi)$ then $c_{M*}[M]$ is nonzero.
If $\chi(M)>\chi(\pi)$ this need not be true.
Given any finitely presentable group $\pi$ there is a 
finite 2-complex $K$ with $\pi_1 (K)\cong\pi$.                                
The boundary of a regular neighbourhood $N$ 
of some embedding of $K$ in $\mathbb{R}^5 $ is a 
closed orientable 4-manifold $M$ with $\pi_1 (M)\cong\pi$.
As the inclusion of $M$ into $N$ is 2-connected and $K$ 
is a deformation retract of $N$
the classifying map $c_M$ factors through $c_K $ and so
induces the trivial homomorphism on homology in degrees $>2$.
However if $M$ and $\pi$ are orientable and $\beta_2(M)<2\beta_2 (\pi)$ then 
$c_M$ must have nonzero degree, for the image of 
$H^2(\pi;\mathbb{Q})$ in $H^2(M;\mathbb{Q})$ then cannot be self-orthogonal 
under cup-product.

\begin{theorem} 
Let $\pi$ be a $PD_4$-group of type $FF$. 
Then $\mathrm{def}(\pi)<1-\frac12\chi(\pi)$.
\end{theorem}
                      
\begin{proof} 
Suppose that $\pi$ has a presentation of deficiency $d\geq1-\frac12\chi(\pi)$, 
and let $X$ be the corresponding finite 2-complex. 
Then 
\[
\beta_2^{(2)}(\pi)-\beta_1^{(2)}(\pi)\leq 
\beta_2^{(2)}(X)-\beta_1^{(2)}(\pi)=\chi(X)=1-d.
\]
Since $\beta_2^{(2)}(\pi)-2\beta_1^{(2)}(\pi)=\chi(\pi)$ and
$\chi(\pi)\geq 2-2d$ it follows that $\beta_1^{(2)}(\pi)\leq{d-1}$.
Hence $\beta_2^{(2)}(X)=0$.  
Therefore $X$ is aspherical, by Theorem 2.4, and so $c.d.\pi\leq 2$.
But this contradicts the hypothesis that $\pi$ is a $PD_4$-group. 
\end{proof}

Note that if $\chi(\pi)$ is odd
the conclusion does not imply that $\mathrm{def}(\pi)\leq-\frac12\chi(\pi)$.
An old conjecture of H.Hopf asserts that if $M$ is an aspherical smooth $2k$-manifold then $(-1)^k\chi(M)\geq0$.
The first open case is when $k=2$.
If Hopf's conjecture is true then $\mathrm{def}(\pi_1(M))\leq0$.
Is $\mathrm{def}(\pi)\leq0$ for every $PD_4$-group $\pi$?
This bound is best possible for groups with $\chi=0$, 
since the presentation 
$\langle a,b\mid ba^2=a^3b^2,~b^2a=a^2b^3\rangle$
gives a Cappell-Shaneson 2-knot group $\mathbb{Z}^3\rtimes_A\mathbb{Z}$.

The hypothesis on orientation characters in Theorem 3.2 is often redundant.

\begin{theorem} 
Let $f:M\to N$ be a $2$-connected map between 
finite $PD_4 $-complexes with $\chi(M)=\chi(N)$. 
If $H^2 (N;\mathbb{F}_2)\not= 0$ then $f^* w_1 (N)=w_1 (M)$, 
and if moreover $N$ is orientable and $H^2(N;\mathbb{Q})\not=0$ then $f$ 
is a homotopy equivalence.
\end{theorem} 
                      
\begin{proof}
Since $f$ is 2-connected $H^2 (f;\mathbb{F}_2)$ is 
injective, and since $\chi(M)=\chi(N)$ it is an isomorphism.
Since $H^2 (N;\mathbb{F}_2)\not=0$, the nondegeneracy of Poincar\'e duality 
implies that $H^4 (f;\mathbb{F}_2)\not=0$, 
and so $f$ is a $\mathbb{F}_2$-(co)homology equivalence.
Since $w_1(M)$ is characterized by the Wu formula $x\cup w_1 (M)=Sq^1 x$ 
for all $x$ in $H^3 (M;\mathbb{F}_2)$, it follows that $f^*w_1 (N)=w_1 (M)$. 

If $H^2(N;\mathbb{Q})\not=0$ then $H^2 (N;\mathbb{Z})$ has positive rank and 
$H^2(N;\mathbb{F}_2)\not=0$, so $N$ orientable implies $M$ orientable. 
We may then repeat the above argument with integral coefficients,
to conclude that $f$ has degree $\pm 1$.
The result then follows from Theorem 3.2. 
\end{proof}

The argument breaks down if, for instance, $M=S^1 \tilde\times S^3 $ is the nonorientable
$S^3 $-bundle over $S^1 $, $N=S^1 \times S^3 $ and $f$ is the composite of the projection 
of $M$ onto $S^1 $ followed by the inclusion of a factor.

If $M$ and $N$ are closed 4-manifolds with isomorphic algebraic 2-types then 
there is a 3-connected map $f:M\to P_2 (N)$. 
The restriction of such a map to $M_o =M\setminus{D^4} $ 
is homotopic to a map $f_o :M_o \to N$ which induces isomorphisms on $\pi_i $ for $i\leq 2$.
In particular, $\chi(M)=\chi(N)$. 
Thus if $f_o $ extends to a map from $M$ to $N$ 
we may be able to apply Theorem 3.2.           
However we usually need more information on how the top cell is attached.
In fact the triple $(P_2(M),w_1(M),f_{M*}[M])$ is 
a complete invariant of the homotopy type \cite{[BB08]}.
(However which triples are thus realized is unknown.)
Can $f_{M*}[M]$ be replaced here by a more explicit ``primary" invariant,
such as the equivariant intersection pairing on 
$\pi_2 (M)$? (See also \cite{[Hi20]}.) 

The following criterion arises in studying the homotopy types
of circle bundles over 3-manifolds. (See Chapter 4.)
   
\begin{theorem} 
Let $E$ be a $PD_4 $-complex with fundamental group
$\pi $ and such that $H_4 (f_E ;\mathbb{Z}^{w_1 (E)} )$ is a monomorphism.
A $PD_4 $-complex $M$ is homotopy equivalent to $E$ if and 
only if there is an isomorphism $\theta $ from $\pi_1 (M)$ to $\pi $ 
such that $w_1 (M)=w_1 (E)\theta $, 
there is a lift $\hat c:M\to P_2 (E)$ of $\theta c_M $ such that
$\hat c_* [M]=\pm f_{E*} [E]$ and $\chi (M)=\chi (E)$.
\end{theorem}
                                                                    
\begin{proof} The conditions are clearly necessary. 
Conversely, suppose that they hold. 
We shall adapt to our situation the arguments of Hendriks in analyzing the 
obstructions to the existence of a degree 1 map between $PD_3 $-complexes realizing a 
given homomorphism of fundamental groups. 
For simplicity of notation we shall write $\widetilde{\mathbb{Z}}$ for 
$\mathbb{Z}^{w_1 (E)}$ 
and also for $\mathbb{Z}^{w_1 (M)} (=\theta^*\widetilde{\mathbb{Z}})$, 
and use $\theta $ to identify $\pi_1 (M)$ with $\pi $
and $K(\pi_1 (M),1)$ with $K(\pi ,1)$. 
We may suppose the sign of the fundamental class $[M]$ is so
chosen that $\hat c_* [M]=f_{E*} [E]$.

Let $E_o =E\setminus{D^4}$. 
Then $P_2 (E_o )=P_2 (E)$ and may be constructed as the union 
of $E_o $ with cells of dimension $\geq 4$. Let 
\begin{equation*}
h:\widetilde{\mathbb{Z}}\otimes_{\mathbb{Z}[\pi]} \pi_4 (P_2 (E_o ),E_o )
\to{H_4 (P_2 (E_o ),E_o ;\widetilde{\mathbb{Z}})}
\end{equation*}
be the $w_1 (E)$-twisted relative Hurewicz homomorphism, 
and let $\partial $ be the connecting 
homomorphism from $\pi_4 (P_2 (E_o ),E_o )$ to $\pi_3 (E_o )$ 
in the exact sequence of 
homotopy for the pair $(P_2 (E_o ),E_o )$. 
Then $h$ and $\partial $ are isomorphisms
since $f_{E_o } $ is 3-connected, and so the homomorphism
$\tau_E:H_4(P_2 (E);\widetilde{\mathbb{Z}})\to
{\widetilde{\mathbb{Z}}\otimes_{\mathbb{Z}[\pi]}\pi_3(E_o)}$
given by the composite of the inclusion 
\begin{equation*}
H_4 (P_2 (E);\widetilde{\mathbb{Z}})=H_4 (P_2 (E_o );\widetilde{\mathbb{Z}})
\to {H_4 (P_2 (E_o ),E_o ;\widetilde{\mathbb{Z}})}
\end{equation*} 
with $h^{-1}$ and $1\otimes_{\mathbb{Z}[\pi]}\partial $ is a monomorphism. 
Similarly $M_o =M\setminus{D^4}$ may be viewed as a subspace of $P_2 (M_o )$ 
and there is a monomorphism $\tau_M $ from 
$H_4 (P_2 (M);\widetilde{\mathbb{Z}})$ to 
$\widetilde{\mathbb{Z}}\otimes_{\mathbb{Z}[\pi]} \pi_3 (M_o )$. 
These monomorphisms are natural with respect to maps 
defined on the 3-skeleta (i.e., $E_o $ and $M_o $). 

The classes $\tau_E (f_{E*} [E])$ and $\tau_M (f_{M*} [M])$ 
are the images of the primary
obstructions to retracting $E$ onto $E_o $ and $M$ onto $M_o $, 
under the Poincar\'e duality isomorphisms from 
$H^4 (E,E_o ;\pi_3 (E_o ))$ to
$H_0(E\setminus{E_o};\widetilde{\mathbb{Z}}\otimes_{\mathbb{Z}[\pi]} 
\pi_3 (E_o ))=\widetilde{\mathbb{Z}}\otimes_{\mathbb{Z}[\pi]} \pi_3 (E_o )$
and $H^4 (M,M_o ;\pi_3 (M_o ))$ to 
$\widetilde{\mathbb{Z}}\otimes_{\mathbb{Z}[\pi]} \pi_3 (M_o )$, respectively. 
Since $M_o $ is homotopy equivalent to a cell complex of dimension $\leq 3$ 
the restriction of $\hat c$ to $M_o $ is homotopic to a map from $M_o $ to 
$E_o $. 
Let $\hat c_\sharp $ be the homomorphism from $\pi_3 (M_o )$ to $\pi_3 (E_o )$ 
induced by $\hat c|M_o $. 
Then $(1\otimes_{\mathbb{Z}[\pi]} \hat c_\sharp )\tau_M (f_{M*} [M])=
\tau_E (f_{E*} [E])$. 
It follows as in \cite{[Hn77]} that the obstruction to extending
$\hat c|M_o :M_o \to E_o $ to a map $d$ from $M$ to $E$ is trivial.
                                                       
Since $f_{E*} d_* [M]=\hat c_* [M]=f_{E*} [E]$ and $f_{E*} $ is a monomorphism 
in degree 4 the map $d$ has degree 1, 
and so is a homotopy equivalence, by Theorem 3.2. 
\end{proof}

If there is such a lift $\hat c $ then $c_M^* \theta^* k_1 (E)=0$ and 
$\theta_* c_{M*} [M]=c_{E*} [E]$. 

\section{Finitely dominated covering spaces}

In this section we shall show that if a $PD_4 $-complex $M$ has a
finitely dominated, infinite regular covering space then either $M$ is 
aspherical or its universal covering space is homotopy equivalent to $S^2$ 
or $S^3 $.
In Chapters 4 and 5 we shall see that such manifolds are close to being 
total spaces of fibre bundles.

\begin{theorem} 
Let $M$ be a $PD_4$-complex with fundamental group $\pi$,
and let $M_\nu$ be the covering space associated to $\nu=\mathrm{Ker}(p)$,
where ${p:\pi\to G}$ is an epimorphism.
Suppose that $M_\nu$ is finitely dominated.
Then 
\begin{enumerate}
\item $G$ has finitely many ends; 

\item if $M_\nu$ is acyclic then it is contractible and $M$ is aspherical;

\item if $G$ has one end and $\nu$ is infinite and $FP_3$ then $M$ is 
aspherical and $M_\nu$ is homotopy equivalent to an aspherical closed surface 
or to $S^1$;

\item if $G$ has one end and $\nu$ is finite 
then $M_\nu\simeq S^2$ or $RP^2$ or is acyclic;

\item $G$ has two ends if and only if $M_\nu$ is a $PD_3 $-complex.
\end{enumerate}
\end{theorem}

\begin{proof} We may clearly assume that $G$ is infinite.
As $\mathbb{Z}[G]$ has no nonzero left ideal (i.e., submodule) which 
is finitely generated as an abelian group
$Hom_{\mathbb{Z}[G]} (H_q (M_\nu;\mathbb{Z}),\mathbb{Z}[G])=0$ 
for all $q\geq 0$, 
and so the bottom row of the UCSS for the covering $p$ is 0.
From Poincar\'e duality and the UCSS we find that 
$H_4 (M_\nu;\mathbb{Z})=H^0(G;\mathbb{Z}[G])=0$ and
$H^1 (G;\mathbb{Z}[G])\cong\overline{H_3 (M_\nu;\mathbb{Z})}$. 
As this group is finitely generated, and as $G$ is infinite, 
$G$ has one or two ends. 
Similarly, $H^2(G;\mathbb{Z}[G])$ is finitely generated 
and so is $\mathbb{Z}$ or 0.

If $M_\nu$ is acyclic
$D_*=\mathbb{Z}[G]\otimes_{\mathbb{Z}[\pi]}C_*(\widetilde{M})$ 
is a resolution of the augmentation $\mathbb{Z}[G]$-module $\mathbb{Z}$ and
$H^q(D_*)\cong H_{4-q}(M_\nu;\mathbb{Z})$.
Hence $G$ is a $PD_4$-group, and so $H_s(\widetilde{M};\mathbb{Z})=
H_s(M_\nu;\mathbb{Z}[\nu])=H^{-s}(M_\nu;\mathbb{Z}[\nu])=0$ for $s>0$,
by Theorem 1.19$^\prime$.
Thus $M_\nu$ is contractible and so $M$ is aspherical.                                                                                                
Suppose that $G$ has one end. 
If $H^2(G;\mathbb{Z}[G])\cong\mathbb{Z}$ 
then $G$ is virtually a $PD_2$-group,
by Bowditch's Theorem, and so $M_\nu$ is a $PD_2$-complex \cite{[Go79]}.
In general, $C_* (\widetilde{M})|_\nu$ 
is chain homotopy equivalent to a finitely generated projective 
$\mathbb{Z}[\nu]$-chain complex $P_*$ and
$H_3 (M_\nu;\mathbb{Z})=H_4 (M_\nu;\mathbb{Z})=0$.
If $\nu$ is $FP_3 $ then the augmentation $\mathbb{Z}[\nu]$-module $\mathbb{Z}$
has a free resolution $F_* $ which is finitely generated in degrees $\leq 3$.
On applying Schanuel's Lemma to the exact sequences
\begin{gather*}
0\to Z_2 \to{P}_2 \to{P}_1 \to{P}_0 \to\mathbb{Z}\to 0\\
\andeq0\to\partial{F}_3 \to{F}_2 \to{F}_1 \to{F}_0 \to\mathbb{Z}\to 0\eqand
\end{gather*}
derived from these two chain complexes we find that $Z_2$
is finitely generated as a $\mathbb{Z}[\nu]$-module. 
Hence $\Pi=\pi_2 (M)=\pi_2 (M_\nu)$ is also finitely generated as
a $\mathbb{Z}[\nu]$-module and so $Hom_\pi (\Pi,\mathbb{Z}[\pi])=0$.
If moreover $\nu$ is infinite then $H^s (\pi;\mathbb{Z}[\pi])=0$ for $s\leq 2$,
so $\Pi=0$, by Lemma 3.3, and $M$ is aspherical.
If $H^2(G;\mathbb{Z}[G])=0$ a spectral sequence corner argument then 
shows that $H^3(G;\mathbb{Z}[G])\cong\mathbb{Z}$ and $M_\nu\simeq S^1$. 
(See the following theorem.)

If $\nu$ is finite but $M_\nu$ is not acyclic then the universal 
covering space $\widetilde M$ is also finitely dominated but not contractible, 
and $\Pi=H_2 (\widetilde M;\mathbb{Z})$ is a nontrivial finitely generated
abelian group, while $H_3(\widetilde M;\mathbb{Z})=H_4(\widetilde M;\mathbb{Z})=0$.
If $C$ is a finite cyclic subgroup of $\pi$ there are isomorphisms 
$H_{n+3}(C;\mathbb{Z})\cong H_n(C;\Pi)$, for all $n\geq4$, by Lemma 2.10.
Suppose that $C$ acts trivially on $\Pi$.
Then if $n$ is odd this isomorphism reduces to $0=\Pi/|C|\Pi$.
Since $\Pi$ is finitely generated, 
this implies that multiplication by $|C|$ is an isomorphism.
On the other hand, if $n$ is even we have $Z/|C|Z\cong\{ a\in\Pi\mid |C|a=0\}$.
Hence we must have $C=1$. 
Now since $\Pi$ is finitely generated any torsion subgroup 
of $Aut(\Pi)$ is finite.
(Let $T$ be the torsion subgroup of $\Pi$ and suppose 
that $\Pi/T$ has rank $r$. 
Then the natural homomorphism from $Aut(\Pi)$ to $Aut(\Pi/T)$ 
has finite kernel, 
and its image is isomorphic to a subgroup of $GL(r,\mathbb{Z})$, 
which is virtually torsion-free.)
Hence as $\pi$ is infinite it must have elements of infinite order.
Since $H^2 (\pi;\mathbb{Z}[\pi])\cong\overline{\Pi}$, by Lemma 3.3, 
it is a finitely generated abelian group.
Therefore it must be infinite cyclic \cite[Corollary 5.2]{[Fa74]}.
Hence $\widetilde M\simeq S^2 $ and $\nu$ has order at most 2, 
so $M_\nu\simeq S^2 $ or $RP^2 $.

Suppose now that $M_\nu$ is a $PD_3 $-complex.
After passing to a finite covering of $M$, if necessary,
we may assume that $M_\nu$ is orientable.
Then $H^1 (G;\mathbb{Z}[G])\cong\overline{ H_3 (M_\nu;\mathbb{Z})}$, 
and so $G$ has two ends. 
Conversely, if $G$ has two ends we may assume that $G\cong\mathbb{Z}$, 
after passing to a finite covering of $M$, if necessary.
Hence $M_\nu$ is a $PD_3$-complex \cite{[Go79]}.
\end{proof}

The hypotheses that $M$ be a $PD_4$-complex and $M_\nu$ be finitely dominated
can be relaxed to requiring that $M$ be a $PD_4$-space and
$C_*(\widetilde{M})$ be $\mathbb{Z}[\nu]$-finitely dominated, and
the appeal to \cite{[Go79]} can be avoided. (See Theorem 4.1.) 
It can be shown that the hypothesis in (3) that $\nu$ be $FP_3$ 
is redundant if $M$ is a {\it finite\/} $PD_4$-space.
(See \cite{[Hi13b]}.)

\begin{cor}
The covering space $M_\nu$ is homotopy equivalent to a closed surface 
if and only if it is finitely dominated 
and $H^2(G;\mathbb{Z}[G])\cong\mathbb{Z}$.
\qed
\end{cor}

In this case $M$ has a finite covering space which is homotopy equivalent 
to the total space of a surface bundle over an aspherical closed surface.
(See Chapter 5.)

\begin{cor}
The covering space $M_\nu$ is homotopy 
equivalent to $S^1 $ if and only if it is finitely dominated, $G$ has one end, 
$H^2 (G;\mathbb{Z}[G])=0$ and $\nu$ is 
a nontrivial finitely generated free group.
\end{cor}            

\begin{proof} If $M_\nu\simeq S^1 $ then it is finitely dominated and $M$ 
is aspherical,
and the conditions on $G$ follow from the LHSSS. 
The converse follows from part (3) of Theorem 3.9, 
since $\nu$ is infinite and $FP$. 
\end{proof}

In fact any finitely generated free normal subgroup $F$ of a $PD_n$-group 
$\pi$ must be infinite cyclic.
For $\pi/FC_\pi(F)$ embeds in $Out(F)$, so 
$v.c.d.\pi/FC_\pi(F)\leq v.c.d.Out(F(r))<\infty$.
If $F$ is nonabelian then $C_\pi(F)\cap F=1$ and so 
$\pi/F$ is an extension of $\pi/FC_\pi(F)$  by $C_\pi(F)$.
Hence $v.c.d.\pi/F<\infty$.
Since $F$ is finitely generated $\pi/F$ is $FP_\infty$.
Hence we may apply \cite[Theorem 9.11]{[Bi]}, 
and an LHSSS corner argument gives a contradiction.

In the simply connected case ``finitely dominated", ``homotopy equivalent to a
finite complex" and ``having finitely generated homology" are all equivalent.

\begin{cor}
If $H_* (\widetilde M;\mathbb{Z})$ is finitely generated then
either $M$ is aspherical or $\widetilde M$ is homotopy equivalent to $S^2 $ 
or $S^3 $ or $\pi_1 (M)$ is finite. 
\qed
\end{cor}

This was first stated (for $\pi_1 (M)$ satisfying a homological finiteness
condition) in \cite{[Ku78]}.                        
We shall examine the spherical cases more closely in Chapters 10 and 11.
(The arguments in these chapters may apply also to $PD_n$-complexes 
with universal covering space homotopy equivalent to $S^{n-1}$ or $S^{n-2}$. 
The analogues in higher codimensions appear to be less accessible.)

The following variation on the aspherical case shall be used in Theorem 4.8, 
but belongs naturally here.

\begin{theorem}
Let $\nu$ be a nontrivial $FP_3$ normal subgroup of infinite index 
in a $PD_4$-group $\pi$, and let $G=\pi/\nu$. Then either
\begin{enumerate}
\item $\nu$ is a $PD_3$-group and $G$ has two ends; 

\item $\nu$ is a $PD_2 $-group and $G$ is virtually a $PD_2$-group; or

\item $\nu\cong\mathbb{Z}$, $H^s(G;\mathbb{Z}[G])=0$ for $s\not=3$ and 
$H^3 (G;\mathbb{Z}[G])\cong\mathbb{Z}$.
\end{enumerate}
\end{theorem}
                 
\begin{proof} 
Since $c.d.\nu<4$, by Strebel's Theorem, 
$\nu$ is $FP$ and hence $G$ is $FP_\infty$.
The $E_2$ terms of the LHSSS with coefficients 
$\mathbb{Q}[\pi]$ can then be expressed as  
$E_2^{pq}=H^p (G;\mathbb{Q}[G])\otimes H^q (\nu;\mathbb{Q}[\nu])$. 
If $H^j(G;\mathbb{Q}[G])$ and $H^k(\nu;\mathbb{Q}[\nu])$ 
are the first nonzero such cohomology groups then $E_2^{jk}$
persists to $E_\infty$.
Hence $j+k=4$ and $H^j (G;\mathbb{Q}[G])$ and $H^{4-j} (\nu;\mathbb{Q}[\nu])$
each have dimension 1 over $\mathbb{Q}$,
since $\pi$ is a $PD_4$-group.
If $j=1$ then $G$ has two ends and so is virtually $\mathbb{Z}$, 
and then $\nu$ is a $PD_3$-group \cite[Theorem 9.11]{[Bi]}. 
If $j=2$ then $\nu$ and $G$ are virtually $PD_2$-groups, 
by Bowditch's Theorem.
Since $\nu$ is torsion-free it is then a $PD_2$-group.
The only remaining possibility is (3).
\end{proof}

In case (1) $\pi$ has a subgroup of index $\leq2$ which is a semidirect product 
$H\rtimes_\theta\mathbb{Z}$ with $\nu\leq H$ and $[H:\nu]<\infty$.
Is it sufficient that $\nu$ be $FP_2$ (as in Theorem 1.19)?
Must the quotient $\pi/\nu$ be virtually a $PD_3$-group in case (3)?

\begin{cor}
If $K$ is $FP_2$ and is ascendant in $\nu$ where $\nu$ is 
an $FP_3$ normal subgroup of infinite index in the $PD_4$-group $\pi$ 
then $K$ is a $PD_k$-group for some $k<4$.
\end{cor}
                     
\begin{proof} 
This follows from Theorem 3.10 together with Theorem 2.17. 
\end{proof}

What happens if we drop the hypothesis that the covering be regular? 
It follows easily from Theorem 2.18 that a $PD_3$-complex has a 
finitely dominated infinite covering space if and only if its fundamental 
group has one or two ends \cite{[Hi08]}. 
We might conjecture that if a $PD_4$-complex $M$ has a finitely dominated 
infinite covering space $\widehat M$ then either $M$ is aspherical or  
$\widetilde M$ is homotopy equivalent to $S^2$ or $S^3$ or $M$ 
has a finite covering space which is homotopy equivalent to the 
mapping torus of a self homotopy equivalence of a $PD_3$-complex.
(In particular, $\pi_1(M)$ has one or two ends.)
In \cite{[Hi08]} we extend the arguments of Theorem 3.9 to show that if 
$\pi_1(\widehat M)$ is $FP_3$ and ascendant in $\pi$ the only other 
possibility is that $\pi_1(\widehat M)$ has two ends, 
$h(\sqrt\pi)=1$ and $H^2(\pi;\mathbb{Z}[\pi])$ is not finitely generated.
This paper also considers in more detail $FP$ ascendant subgroups of 
$PD_4$-groups, corresponding to the aspherical case. 

\section{Minimizing the Euler characteristic}

It is well known that every finitely presentable group 
is the fundamental group of some closed orientable 4-manifold.
Such manifolds are far from unique, 
for the Euler characteristic may be made arbitrarily large 
by taking connected sums with simply connected manifolds.
Following Hausmann and Weinberger \cite{[HW85]}, 
we may define an invariant for any
finitely presentable group $\pi$ by
\[
q(\pi)=\min\{\chi(M)|M~is~a~PD_4~complex~with~\pi_1 (M)\cong\pi\}.
\] 
We may also define related invariants $q^X $ where the minimum 
is taken over the class
of $PD_4 $-complexes whose normal fibration has an $X$-reduction. 
There are the following basic estimates for $q^{SG}$,
which is defined in terms of $PD_4^+ $-complexes.

\begin{lemma}
Let $\pi$ be a finitely presentable group with a subgroup $H$ of finite index 
and let $F$ be a field. 
Then 

\begin{enumerate}
\item $1-\beta_1 (H;F)+\beta_2 (H;F)\leq [\pi:H](1-\mathrm{def}\pi)$;

\item $2-2\beta_1 (H;F)+\beta_2 (H;F)\leq [\pi:H]q^{SG} (\pi)$;

\item $q^{SG} (\pi)\leq 2(1-\mathrm{def}(\pi))$;

\item if all cup products of pairs of elements of $H^2(\pi;F)$ 
are trivial then 
$q^{SG}(\pi)\geq 2(1-\beta_1(\pi;F)+\beta_2(\pi;F))$,
and if moreover $H^4(\pi;\mathbb{F}_2)=0$ then 
$q(\pi)\geq 2(1-\beta_1(\pi;\mathbb{F}_2)+\beta_2(\pi;\mathbb{F}_2))$ also.

\end{enumerate}
\end{lemma}

\begin{proof}
Let $C$ be the 2-complex corresponding to a presentation for $\pi$ of maximal
deficiency and let $C_H$ be the covering space associated to the subgroup $H$.
Then $\chi(C)=1-\mathrm{def}\pi$ and $\chi(C_H)=[\pi:H]\chi(\pi)$.
Condition (1) follows since $\beta_1(H;F)=\beta_1(C_H;F)$ and
$\beta_2(H;F)\leq\beta_2(C_H;F)$.

Condition (2) follows similarly on considering the Euler characteristics of
a $PD_4^+$-complex $M$ with $\pi_1(M)\cong \pi$ and of the associated covering
space $M_H$.

The boundary of a regular neighbourhood of a PL embedding of $C$ in $R^5$ 
is a closed orientable 4-manifold realizing the upper bound in (3).

The image of $H^2(\pi;F)$ in $H^2(M;F)$ has dimension $\beta_2(\pi;F)$,
and is self-annihilating under cup-product if $H^4(\pi;F)=0$.
In that case $\beta_2 (M;F)\geq 2\beta_2(\pi;F)$, 
which implies the first part of (4).
The final observation follows since all $PD_n$-complexes
are orientable over $\mathbb{F}_2$.
\end{proof}

Condition (2) was used in \cite{[HW85]} to give examples 
of finitely presentable superperfect groups 
which are not fundamental groups of homology 4-spheres.
(See Chapter 14 below.)

If $\pi$ is a finitely presentable, orientable $PD_4 $-group we see 
immediately that $q^{SG} (\pi)\geq\chi(\pi)$.
Multiplicativity then implies that $q(\pi)=\chi(\pi)$ if $K(\pi,1)$ 
is a finite $PD_4 $-complex.

For groups of cohomological dimension at most 2 we can say more.

\begin{theorem}
Let $X$ be a $PD_4$-complex with fundamental group $\pi$ such that 
$c.d.\pi\leq2$, and let $C_*=C_*(X;\mathbb{Z}[\pi])$.
Then 
\begin{enumerate}

\item 
$C_*$ is $\mathbb{Z}[\pi]$-chain homotopy equivalent to 
$D_*\oplus{L[2]}\oplus D^{4-*}$, 
where $D_*$ is a projective resolution of $\mathbb{Z}$,
$L[2]$ is a finitely generated projective module $L$ concentrated 
in degree $2$ and $D^{4-*}$ is the conjugate dual of $D_*$, 
shifted to terminate in degree $2$;

\item $\pi_2(X)\cong{L}\oplus{\overline{H^2(\pi;\mathbb{Z}[\pi])}}$;

\item $\chi(X)\geq2\chi(\pi)$, with equality if and only if $L=0$;

\item $Hom_{\mathbb{Z}[\pi]}(\overline{H^2(\pi;\mathbb{Z}[\pi])},
\mathbb{Z}[\pi])=0$.

\end{enumerate}
\end{theorem}

\begin{proof}
The chain complex $C_*$ gives a resolution of the augmentation module
\[
0\to \mathrm{Im}(\partial_2^C)\to C_1\to C_0\to\mathbb{Z}\to0.
\]
Let $D_*$ be the corresponding chain complex with $D_0=C_0$, $D_1=C_1$ and 
$D_2=\mathrm{Im}(\partial_2^C)$.
Since $c.d.\pi\leq2$ and $D_0$ and $D_1$ are projective modules $D_2$ 
is projective, by Schanuel's Lemma.
Therefore the epimorphism from $C_2$ to $D_2$ splits,
and so $C_*$ is a direct sum $C_*\cong D_*\oplus(C/D)_*$.
Since $X$ is a $PD_4$-complex $C_*$ is chain homotopy equivalent to $C^{4-*}$.
The first two assertions follow easily.

On taking homology with simple coefficients $\mathbb{Q}$, we see that 
$\chi(X)=2\chi(\pi)+\mathrm{dim}_\mathbb{Q}\mathbb{Q}\otimes_\pi{L}$.
Hence $\chi(X)\geq2\chi(\pi)$.
Since $\pi$ satisfies the Weak Bass Conjecture \cite{[Ec86]} 
and $L$ is projective,
$L=0$ if and only if $\mathrm{dim}_\mathbb{Q}\mathbb{Q}\otimes_\pi{L}=0$.

Let $\delta:D_2\to D_1$ be the inclusion.
Then $\overline{H^2(\pi;\mathbb{Z}[\pi])}=\mathrm{Cok}(\delta^\dagger)$, 
where $\delta^\dagger$ is the conjugate transpose of $\delta$.
Hence 
$\overline{Hom_{\mathbb{Z}[\pi]}(\overline{H^2(\pi;\mathbb{Z}[\pi])},
\mathbb{Z}[\pi])}=\mathrm{Ker}(\delta^{\dagger\dagger)}$.
But $\delta^{\dagger\dagger}=\delta$, which is injective, and so
$Hom_{\mathbb{Z}[\pi]}(\overline{H^2(\pi;\mathbb{Z}[\pi])},\mathbb{Z}[\pi])=0$.
\end{proof}

The appeal to the Weak Bass Conjecture may be avoided if $X$ and $K(\pi,1)$
are homotopy equivalent to finite complexes.
For then $L$ is stably free, and so is 0 if and only if
$\mathbb{Z}\otimes_{\mathbb{Z}[\pi]}{L}=0$,
since group rings are weakly finite.

Similar arguments may be used to prove the following variation.

\begin{add}
Suppose that $c.d._R\pi\leq2$ for some ring $R$.
Then ${R\otimes\pi_2 (M)}\cong{P}\oplus\overline{H^2(\pi;R[\pi])}$, 
where $P$ is a projective $R[\pi]$-module, 
and $\chi(M)\geq2\chi(\pi;R)=$ ${2(1-\beta_1(\pi;R)+\beta_2(\pi;R))}$.
If $R$ is a subring of $\mathbb{Q}$ then $\chi(M)=2\chi(\pi;R)$
if and only if $\pi_2(M)\cong\overline{H^2(\pi;\mathbb{Z}[\pi])}$.
\qed
\end{add}

There are many natural examples of 4-manifolds with $\pi_1(M)=\pi$
having nontrivial torsion and such that
$c.d._{\mathbb{Q}}\pi\leq2$ and $\chi(M)=2\chi(\pi)$.
(See Chapters 10 and 11.) 
However all the known examples satisfy $v.c.d.\pi\leq2$.

\begin{cor}
If $H_2 (\pi;\mathbb{F}_2)\not=0$ the Hurewicz homomorphism from 
$\pi_2 (M)$ to $H_2 (M;\mathbb{F}_2)$ is nonzero.
\end{cor} 

\begin{proof} By the addendum to the theorem, 
$H_2 (M;\mathbb{F}_2)$ has dimension at least $2\beta_2 (\pi)$,
and so cannot be isomorphic to $H_2 (\pi;\mathbb{F}_2)$ unless both are 0. 
\end{proof}

\begin{cor}
If $\pi=\pi_1 (P)$ where $P$ is an aspherical finite $2$-complex then 
$q(\pi)=2\chi(P)$. 
The minimum is realized by an $s$-parallelizable PL $4$-manifold.
\end{cor}

\begin{proof} 
If we choose a PL embedding $j:P\to \mathbb{R}^5 $, the boundary of a regular neighbourhood $N$ 
of $j(P)$ is an $s$-parallelizable PL 4-manifold with fundamental group $\pi$ 
and with Euler characteristic $2\chi(P)$. 
\end{proof}

By Theorem 2.8 a finitely presentable group is the fundamental group 
of an aspherical finite 2-complex if and only if it has 
cohomological dimension $\leq 2$ and is efficient, 
i.e.\ has a presentation of deficiency 
$\beta_1 (\pi;\mathbb{Q})-\beta_2 (\pi;\mathbb{Q})$.
It is not known whether every finitely presentable group of cohomological dimension 2
is efficient.

In Chapter 5 we shall see that if $P$ is an aspherical closed surface and
$M$ is a closed 4-manifold with $\pi_1 (M)\cong \pi$ then $\chi(M)=q(\pi)$ 
if and only if $M$ is homotopy equivalent to the total space of an 
$S^2 $-bundle over $P$.
The homotopy types of such minimal 4-manifolds for $\pi$ may be 
distinguished by their Stiefel-Whitney classes.
Note that if $\pi$ is orientable then $S^2\times P$ is a minimal 
4-manifold for $\pi$ which is both $s$-parallelizable and also a 
projective algebraic complex surface.
Note also that the conjugation of the module structure in the theorem involves 
the orientation character of $M$ which may differ from that of the 
$PD_2$-group $\pi$. 

\begin{cor}
If $\pi$ is the group of an unsplittable $\mu$-component $1$-link 
then $q(\pi)=0$. 
\qed
\end{cor}
                           
If $\pi$ is the group of a $\mu$-component $n$-link with $n\geq 2$ then $H_2 (\pi;\mathbb{Q})=0$
and so $q(\pi)\geq2(1-\mu)$, with equality if and only if $\pi$ is the group of a 2-link.
(See Chapter 14.)

\begin{cor}
If $\pi$ is an extension of $\mathbb{Z}$ by a finitely
generated free normal subgroup then $q(\pi)=0$.
\qed 
\end{cor}        

In Chapter 4 we shall see that if $M$ is a closed 4-manifold with 
$\pi_1 (M)$ such an extension then $\chi(M)=q(\pi)$ 
if and only if $M$ is homotopy equivalent to a manifold 
which fibres over $S^1 $ with fibre a closed 3-manifold 
with free fundamental group, and then $\pi$ and $w_1(M)$ 
determine the homotopy type.

Finite generation of the normal subgroup is essential;
$F(2)$ is an extension of $\mathbb{Z}$ by $F(\infty)$, 
and $q(F(2))=2\chi(F(2))=-2$.

Let $\pi$ be the fundamental group of a closed orientable 3-manifold.
Then $\pi\cong F*\nu$ where $F$ is free of rank $r$ and 
$\nu$ has no infinite cyclic free factors.
Moreover $\nu=\pi_1(N)$ for some closed orientable 3-manifold $N$.
If $M_0$ is the closed 4-manifold obtained by surgery on $\{n\}\times S^1$ in 
$N\times S^1$
then $M=M_0\sharp(\sharp^r(S^1\times S^3)$ is a smooth $s$-parallelisable 4-manifold
with $\pi_1(M)\cong\pi$ and $\chi(M)=2(1-r)$.
Hence $q^{SG} (\pi)=2(1-r)$, by part (4) of Lemma 3.11.

The arguments of Theorem 3.12 give stronger results in this case also.

\begin{theorem} 
Let $\pi$ be a finitely presentable $PD_3$-group, 
and let $M$ be a $PD_4$-complex with fundamental group $\pi$ and
$w_1 (\pi)=w_1 (M)$.
Then $q(\pi)=2$, and there are finitely generated projective 
$\mathbb{Z}[\pi]$-modules $P$ and $P'$ such that 
$\pi_2 (M)\oplus{P}\cong A(\pi)\oplus{P'}$,
where $A(\pi)$ is the augmentation ideal of $\mathbb{Z}[\pi]$.
\end{theorem}

\begin{proof}
Let $N$ be a $PD_3$-complex with fundamental group $\pi$. 
We may suppose that $N=N_o \cup D^3$, where $N_o \cap D^3 =S^2 $.
Let $X=N_o \times S^1 \cup S^2 \times D^2 $.
Then $X$ is a $PD_4 $-complex, $\chi(X)=2$ and $\pi_1 (X)\cong \pi$.
Hence $q(\pi)\leq2$. 
On the other hand, $q(\pi)\geq2$ by part (4) of Lemma 3.11,
and so $q(\pi)=2$.

For any left $\mathbb{Z}[\pi]$-module $L$ let
$e^iL=Ext^i_{\mathbb{Z}[\pi]}(L,\mathbb{Z}[\pi])$, 
to simplify the notation. 
The cellular chain complex for the universal covering space of $M$ 
gives exact sequences
\begin{gather}
0\to C_4 \to C_3 \to Z_2 \to H_2 \to 0\\
\andeq0\to Z_2 \to C_2 \to C_1 \to C_0 \to\mathbb{Z}\to 0.\eqand
\end{gather}
Since $\pi$ is a $PD_3 $-group the augmentation module 
$\mathbb{Z}$ has a finite projective resolution of length 3. 
On comparing sequence 3.2 with such a resolution and applying
Schanuel's lemma we find that $Z_2 $ is a finitely generated 
projective $\mathbb{Z}[\pi]$-module.
Since $\pi$ has one end, the UCSS reduces to an exact sequence
\begin{equation}
0\to H^2 \to e^0 H_2 \to e^3\mathbb{Z}\to H^3 \to e^1 H_2 \to 0\
\end{equation}
and isomorphisms 
$H^4\cong e^2 H_2 $ and $e^3 H_2 =e^4 H_2 =0.$ 
Poincar\'e duality implies that 
$H^3=0$ and $H^4 \cong\overline{\mathbb{Z}}$.    
Hence sequence 3.3 reduces to
\begin{equation}
0\to H^2 \to e^0 H_2 \to e^3\mathbb{Z}\to 0
\end{equation}
and $e^1 H_2 =0$. 
On dualizing the sequence 3.1 and conjugating
we get an exact sequence of left modules
\begin{equation}
0\to\overline{e^0 H_2}\to\overline{e^0 Z_2}\to\overline{e^0 C_3} 
\to\overline{e^0 C_4}\to\overline{e^2 H_2}\cong\mathbb{Z}\to 0.
\end{equation}
Schanuel's lemma again implies that $\overline{e^0H_2}$ 
is a finitely generated projective module.
Now $\pi_2(M)\cong\overline{H^2}$, by Poincar\'e duality,
and $\overline{e^3\mathbb{Z}}\cong\mathbb{Z}$, since 
$\pi$ is a $PD_3 $-group and $w_1(M)=w_1(\pi)$.
Hence the final assertion follows from sequence 3.4 and Schanuel's Lemma. 
\end{proof}

The invariant $q^{SG}(\pi)$ has been determined for $\pi$ 
a 3-manifold group with no 2-torsion \cite{[SW21]}.
Does Theorem 3.13 extend to all free products of $PD_3 $-groups?

There has been some related work estimating the difference 
$\chi(M)-|\sigma(M)|$ where $M$ is a closed orientable 4-manifold $M$ 
with $\pi_1(M)\cong\pi$ and where $\sigma(M)$ is the signature of $M$.
In particular, this difference is always $\geq0$ if $\beta_1^{(2)}(\pi)=0$.
(See \cite{[JK93]} and \cite[Chapter 7.\S3]{[Lu]}.)
The minimum value of this difference ($p(\pi)=\min\{\chi(M)-|\sigma(M)|\}$)
is another numerical invariant of $\pi$, which is studied in \cite{[Ko94]}.

\section{Euler Characteristic 0}

In this section we shall consider the interaction of the fundamental group 
and Euler characteristic from another point of view.
We shall assume that $\chi(M)=0$ and show that if $\pi$ is an ascending HNN extension 
then it satisfies some very stringent conditions. 
The groups $Z*_m$ shall play an important role.
We shall approach our main result via several lemmas.

We begin with a simple observation relating Euler characteristic 
and fundamental group
which shall be invoked in several of the later chapters.
Recall that if $G$ is a group then $I(G)$ is the minimal normal subgroup 
such that $G/I(G)$ is free abelian.

\begin{lemma}
Let $M$ be a $PD_4 $-space with $\chi(M)\leq0$.
If $M$ is orientable then $H^1 (M;\mathbb{Z})\not=0$ and so 
$\pi=\pi_1 (M)$ maps onto $\mathbb{Z}$.
If $H^1(M;\mathbb{Z})=0$ then $\pi$ maps onto $D$.
\end{lemma}

\begin{proof} The covering space $M_W$ corresponding to 
$W=\mathrm{Ker}(w_1(M))$ is orientable and 
$\chi(M_W)=2-2\beta_1 (M_W)+\beta_2 (M_W)=[\pi:W]\chi(M)\leq0$.
Therefore $\beta_1 (W)=\beta_1 (M_W)>0$ and so $W/I(W)\cong\mathbb{Z}^r$, 
for some $r>0$.
Since $I(W)$ is characteristic in $W$ it is normal in $\pi$.
As $[\pi:W]\leq 2$ it follows easily that $\pi/I(W)$ maps onto 
$\mathbb{Z}$ or $D$. 
\end{proof}

Note that if $M=RP^4 \sharp RP^4 $, then $\chi(M)=0$ and $\pi_1 (M)\cong D$, 
but $\pi_1(M)$ does not map onto $\mathbb{Z}$.

\begin{lemma}
Let $M$ be a $PD_4^+$-complex such that 
$\chi(M)=0$ and $\pi=\pi_1(M)$ is an extension of $Z*_m$ 
by a finite normal subgroup $F$, for some $m\not=0$.
Then the abelian subgroups of $F$ are cyclic.
If $F\not=1$ then $\pi$ has a subgroup of finite index which is a central 
extension of $Z*_n$ by a nontrivial finite cyclic group, 
where $n$ is a power of $m$.
\end{lemma}

\begin{proof} 
Let $\widehat M$ be the infinite cyclic covering space corresponding to the 
subgroup $I(\pi)$.
Since $M$ is compact and $\Lambda=\mathbb{Z}[\mathbb{Z}]$ is noetherian 
the groups $H_i(\widehat M;\mathbb{Z})=H_i(M;\Lambda)$ are 
finitely generated as $\Lambda$-modules. 
Since $M$ is orientable, $\chi(M)=0$ and $H_1(M;\mathbb{Z})$ 
has rank 1 they are $\Lambda$-torsion modules, 
by the Wang sequence for the projection of $\widehat M$ onto $M$. 
Now $H_2(\widehat M;\mathbb{Z})\cong\overline{Ext^1_\Lambda(I(\pi)/I(\pi)',\Lambda)}$, by Poincar\'e duality.
There is an exact sequence 
\begin{equation*}
0\to T\to I(\pi)/I(\pi)'\to I(Z*_m)\cong \Lambda/(t-m)\to 0,
\end{equation*}
where $T$ is a finite $\Lambda$-module.
Therefore $Ext^1_\Lambda(I(\pi)/I(\pi)',\Lambda)\cong \Lambda/(t-m)$ and so $H_2(I(\pi);\mathbb{Z})$ is a quotient of $\Lambda/(mt-1)$, 
which is isomorphic to $\mathbb{Z}[\frac1m]$ as a group.
Now $I(\pi)/\mathrm{Ker}(f)\cong\mathbb{Z}[\frac1m]$ also, 
and $H_2 (\mathbb{Z}[\frac1m];\mathbb{Z})\cong 
\mathbb{Z}[\frac1m]\wedge\mathbb{Z}[\frac1m]=0$.
(See \cite[page 334]{[Ro]}.)
Hence $H_2(I(\pi);\mathbb{Z})$ is finite, by an LHSSS argument,
and so is cyclic, of order relatively prime to $m$.
                                               
Let $t$ in $\pi$ generate $\pi/I(\pi)\cong\mathbb{Z}$. 
Let $A$ be a maximal abelian subgroup of $F$ and let $C=C_\pi(A)$.
Then $q=[\pi:C]$ is finite, since $F$ is finite and normal in $\pi$.           
In particular, $t^q$ is in $C$ and $C$ maps onto $\mathbb{Z}$, 
with kernel $J$, say.
Since $J$ is an extension of $\mathbb{Z}[\frac1m]$ by a finite normal subgroup
its centre $\zeta J$ has finite index in $J$.
Therefore the subgroup $G$ generated by $\zeta J$ and $t^q$ 
has finite index in $\pi$,
and there is an epimorphism $f$ from $G$ onto $Z*_{m^q}$, with kernel $A$.
Moreover $I(G)=f^{-1}(I(Z*_{m^q}))$ is abelian, and is an extension of
$\mathbb{Z}[\frac1m]$ by the finite abelian group $A$.
Hence it is isomorphic to $A\oplus\mathbb{Z}[\frac1m]$.
(See \cite[page 106]{[Ro]}.)
Now $H_2(I(G);\mathbb{Z})$ is cyclic of order prime to $m$.
On the other hand $H_2(I(G);\mathbb{Z})\cong 
(A\wedge A)\oplus(A\otimes \mathbb{Z}[\frac1m])$
and so $A$ must be cyclic.

If $F\not=1$ then $A$ is cyclic, nontrivial, central in $G$ and $G/A\cong Z*_{m^q}$.
\end{proof}

\begin{lemma}
Let $M$ be a finite $PD_4$-complex with fundamental group $\pi$.
Suppose that $\pi$ has a nontrivial finite cyclic central subgroup $F$ 
with quotient $G=\pi/F$ such that $g.d.G=2$, $e(G)=1$ and $\mathrm{def}(G)=1$.
Then $\chi(M)\geq0$.
If $\chi(M)=0$ and $\Xi=\mathbb{F}_p[G]$ is a weakly finite ring
for some prime $p$ dividing $|F|$
then $\pi$ is virtually $\mathbb{Z}^2$.
\end{lemma}

\begin{proof}
Let $\widehat M$ be the covering space of $M$ with group $F$,
and let $c_q$ be the number of $q$-cells of $M$, for $q\geq0$.
Let $C_*=C_*(M;\Xi)=\mathbb{F}_p\otimes C_*(M)$ be the 
equivariant cellular chain complex of $\widehat M$ with coefficients
$\mathbb{F}_p$ and let $H_p=H_p(M;\Xi)=H_p(\widehat M;\mathbb{F}_p)$.
For any left $\Xi$-module $H$ let $e^qH=\overline{Ext^q_\Xi(H,\Xi)}$.

Since $\widehat M$ is connected and $F$ is cyclic 
$H_0\cong{H_1}\cong\mathbb{F}_p$ and 
since $G$ has one end Poincar\'e duality and the UCSS give $H_3=H_4=0$,
an exact sequence
\begin{equation*}
0\to{e^2\mathbb{F}_p}\to{H_2}\to{e^0H_2}\to{e^2H_1}\to{H_1}\to{e^1H_2}\to0
\end{equation*}
and an isomorphism $e^2H_2\cong\mathbb{F}_p$.
Since $g.d.G=2$ and $\mathrm{def}(G)=1$ the augmentation module has a resolution
\[
0\to\Xi^r\to\Xi^{r+1}\to\Xi\to\mathbb{F}_p\to0.
\]
The chain complex $C_*$ gives four exact sequences
\begin{gather*}
0\to Z_1\to C_1\to C_0\to\mathbb{F}_p\to0,\\
0\to Z_2\to C_2\to Z_1\to\mathbb{F}_p\to0,\\
0\to B_2\to Z_2\to H_2\to 0\\
\andeq0\to C_4\to C_3\to B_2\to 0.\eqand
\end{gather*}
Using Schanuel's Lemma several times we find that the cycle submodules $Z_1$ 
and $Z_2$ are stably free, of stable ranks $c_1-c_0$ and $c_2-c_1+c_0$,
respectively.
Dualizing the last two sequences gives two new sequences
\begin{gather*}
0\to e^0B_2\to e^0C_3\to e^0C_4\to e^1B_2\to 0\\
\andeq0\to e^0H_2\to e^0Z_2\to e^0 B_2\to e^1H_2\to 0,\eqand
\end{gather*}
and an isomorphism $e^1B_2\cong e^2H_2\cong\mathbb{F}_p$.
Further applications of Schanuel's Lemma show that $e^0B_2$
is stably free of rank $c_3-c_4$, and hence that $e^0H_2$ is stably free
of rank $c_2-c_1+c_0-(c_3-c_4)=\chi(M)$.
Since $\Xi$ maps onto the field $\mathbb{F}_p$ the rank must be non-negative,
and so $\chi(M)\geq0$.

If $\chi(M)=0$ and $\Xi=\mathbb{F}_p[G]$ is a weakly finite ring
then $e^0H_2=0$ and so $e^2\mathbb{F}_p=e^2H_1$ is a submodule of 
$\mathbb{F}_p\cong{H_1}$.
Moreover it cannot be 0, for otherwise the UCSS would give $H_2=0$ and then 
$H_1=0$, which is impossible.
Therefore $e^2\mathbb{F}_p\cong\mathbb{F}_p$.

Since $G$ is torsion-free and indicable it must be a 
$PD_2$-group \cite[Theorem V.12.2]{[DD]}.
Since $\mathrm{def}(G)=1$ it follows that $G\cong\mathbb{Z}^2$ 
or $\mathbb{Z}\rtimes_{-1}\mathbb{Z}$,
and hence that $\pi$ is also virtually $\mathbb{Z}^2$.
\end{proof}

The hypothesis ``$\Xi=\mathbb{F}_p[G]$ is weakly finite" is satisfied 
if $G$ is sofic \cite{[ES04]}.
In particular, this is so if $G'$ is finitely generated,
for then $G'$ is free, by Corollary 4.3.1 below,
and so $G$ is residually finite, and hence sofic.
(There are at present no known examples of groups which are not sofic!)


We may now give the main result of this section.          

\begin{theorem} 
Let $M$ be a finite $PD_4$-complex whose fundamental group $\pi$ is an 
ascending HNN extension with finitely generated base $B$.
Then $\chi(M)\geq 0$, and hence $q(\pi)\geq 0$.
If $\chi(M)=0$ and $B$ is $FP_2$ and finitely ended then either $\pi$ 
has two ends or $\pi\cong Z*_m$ or $Z*_m\rtimes(Z/2Z)$ for some $m\not=0$ 
or $\pm1$ or $\pi$ is virtually $\mathbb{Z}^2$ or $M$ is aspherical.
\end{theorem}

\begin{proof} 
The $L^2$ Euler characteristic formula gives $\chi(M)=\beta_2^{(2)}(M)\geq0$,
since $\beta_i^{(2)}(M)=\beta^{(2)}_i(\pi)=0$ for $i=0$ or 1,
by Lemma 2.1.

Let $\phi:B\to B$ be the monomorphism determining $\pi\cong B*_\phi$.
If $B$ is finite then $\phi$ is an automorphism and so $\pi$ has two ends.
If $B$ is $FP_2$ and  has one end then $H^s(\pi;\mathbb{Z}[\pi])=0$ for $s\leq2$, 
by the Brown-Geoghegan Theorem.
If moreover $\chi(M)=0$ then $M$ is aspherical, by Corollary 3.5.2.

If $B$ has two ends then it is an extension of $\mathbb{Z}$ or $D$ 
by a finite normal subgroup $F$. 
As $\phi$ must map $F$ isomorphically to itself, $F$ is normal in $\pi$,
and is the maximal finite normal subgroup of $\pi$.
Moreover $\pi/F\cong Z*_m$, for some $m\not=0$,
if $B/F\cong\mathbb{Z}$, and is a semidirect product
$Z*_m\rtimes(Z/2Z)$, with a presentation
$\langle a,t,u\mid tat^{-1}=a^m\! ,\medspace tut^{-1}=ua^r\! ,\medspace
u^2=1,\medspace uau=a^{-1}\rangle$,
for some $m\not=0$ and some $r\in\mathbb{Z}$, if $B/F\cong D$.
(On replacing $t$ by $a^{[r/2]}t$, if necessary,
we may assume that $r=0$ or 1.)

Suppose first that $M$ is orientable, and that $F\not=1$.
Then $\pi$ has a subgroup $\sigma$ of finite index which is a central extension
of $Z*_{m^q}$ by a finite cyclic group, for some $q\geq1$, by Lemma 3.15.
Let $p$ be a prime dividing $q$.
Since $Z*_{m^q}$ is a torsion-free solvable group the ring 
$\Xi=\mathbb{F}_p[Z*_{m^q}]$ 
has a skew field of fractions $L$, which as a right $\Xi$-module is the
direct limit of the system $\{\Xi_\theta\mid 0\not=\theta\in\Xi\}$, 
where each $\Xi_\theta=\Xi$,
the index set is ordered by right divisibility ($\theta\leq\phi\theta$) and 
the map from $\Xi_\theta$ to $\Xi_{\phi\theta}$ sends $\xi$ to $\phi\xi$ 
\cite{[KLM88]}.
In particular, $\Xi$ is a weakly finite ring and so 
$\pi$ is virtually $\mathbb{Z}^2$, by Lemma 3.16.

If $M$ is nonorientable then either $w_1(M)|_F$ is injective,
so $\pi\cong{Z*_m\rtimes(Z/2Z)}$, or $\pi$ is virtually $\mathbb{Z}^2$.
\end{proof}

Is $M$ still aspherical if $B$ is assumed only finitely generated and one ended? 

\begin{cor}
Let $M$ be a finite $PD_4$-complex such that $\chi(M)=0$ and
$\pi=\pi_1(M)$ is almost coherent and restrained.
Then either $\pi$ has two ends or 
$\pi\cong Z*_m$ or $Z\!*_m\rtimes(Z/2Z)$ for some $m\not=0$ or $\pm1$ 
or $\pi$ is virtually $\mathbb{Z}^2$ or $M$ is aspherical.
\end{cor}
                       
\begin{proof} 
Let $\pi^+=\mathrm{Ker}(w_1(M))$.
Then $\pi^+$ maps onto $\mathbb{Z}$, by Lemma 3.14, and so 
is an ascending HNN extension $\pi^+\cong B*_\phi$
with finitely generated base $B$.
Since $\pi$ is almost coherent $B$ is $FP_2$, 
and since $\pi$ has no nonabelian free subgroup $B$ has at most two ends. 
Hence Lemma 3.16 and Theorem 3.17 apply, so 
either $\pi$ has two ends or $M$ is aspherical or $\pi^+\cong Z*_m$ or 
$Z*_m\rtimes(Z/2Z)$ for some $m\not=0$ or $\pm1$.
In the latter case $\sqrt\pi$ is isomorphic to a subgroup of the additive 
rationals $\mathbb{Q}$, and $\sqrt\pi=C_\pi(\sqrt\pi)$.
Hence the image of $\pi$ in $Aut(\sqrt\pi)\leq\mathbb{Q}^\times$ is infinite.
Therefore $\pi$ maps onto $\mathbb{Z}$ and so is 
an ascending HNN extension $B*_\phi$,
and we may again use Theorem 3.17. 
\end{proof}

Does this corollary hold 
without the hypothesis that $\pi$ be almost coherent? 

There are nine groups which are virtually $\mathbb{Z}^2$ and 
are fundamental groups of $PD_4$-complexes with Euler 
characteristic 0. (See Chapter 11.)
Are any of the groups $Z*_m\rtimes(Z/2Z)$ with $|m|>1$
realized by $PD_4$-complexes with $\chi=0$?
If $\pi$ is restrained and $M$ is aspherical must $\pi$ be virtually poly-$Z$? 
(Aspherical 4-manifolds with virtually poly-$Z$ fundamental groups
are characterized in Chapter 8.)
                                                    
Let $G$ is a group with a presentation of deficiency $d$ 
and $w:G\to\{\pm1\}$ be a homomorphism, and let
$\langle x_i,\medspace 1\leq i\leq m\mid r_j,\medspace 1\leq j\leq n\rangle$
be a presentation for $G$ with $m-n=d$.
We may assume that $w(x_i)=+1$ for $i\leq m-1$.
Let $X=\natural^m(S^1\times D^3)$ if $w=1$ and
$X=(\natural^{m-1}(S^1\times D^3))\natural(S^1\tilde\times D^3)$ otherwise.
The relators $r_j$ may be represented by disjoint 
orientation preserving embeddings of $S^1$ in $\partial X$, 
and so we may attach 2-handles along product neighbourhoods,
to get a bounded 4-manifold $Y$ with $\pi_1(Y)=G$, $w_1(Y)=w$ and $\chi(Y)=1-d$.
Doubling $Y$ gives a closed 4-manifold
$M$ with $\chi(M)=2(1-d)$ and $(\pi_1(M),w_1(M))$ isomorphic to $(G,w)$.

Since the groups $Z*_m$ have deficiency 1 it follows that any homomorphism
$w:Z*_m\to \{\pm1\}$ may be realized as the orientation character of a 
closed 4-manifold $M$ with $\pi_1(M)\cong{Z*_m}$ and $\chi(M)=0$.
In the orientable case such manifolds are determined up to homeomorphism
by $\pi$ and $w_2$ \cite{[HKT09]}.

\section{The intersection pairing}

Let $X$ be a $PD_4$-complex with fundamental group $\pi$ and let $w=w_1(X)$.
In this section it shall be convenient to work with left modules.
Thus if $L$ is a left ${\mathbb{Z}[\pi]}$-module we shall
let $L^\dagger=\overline{Hom_{\mathbb{Z}[\pi]}(L,{\mathbb{Z}[\pi]})}$ 
be the conjugate dual module.
If $L$ is free, stably free or projective so is $L^\dagger$.

Let $H=\overline{H^2(X;\mathbb{Z}[\pi])}$ and $\Pi=\pi_2(X)$, and let
$D:H\to\Pi$ and ${{ev}:H\to\Pi^\dagger}$
be the Poincar\'e duality isomorphism and the evaluation homomorphism, 
respectively.
The {\it cohomology intersection pairing} 
$\lambda:H\times H\to\mathbb{Z}[\pi]$ is
defined by $\lambda(u,v)=ev(v)(D(u))$, for all $u,v\in{H}$. 
This pairing is $w$-hermitian: $\lambda(gu,hv)=g\lambda(u,v)\bar{h}$ and
$\lambda(v,u)=\overline{\lambda(u,v)}$ for all $u,v\in{H}$ and $g,h\in\pi$.
Since $\lambda(u,e)=0$ for all $u\in{H}$ and 
$e\in E=\overline{H^2(\pi;\mathbb{Z}[\pi])}$ 
the pairing $\lambda$ induces a pairing 
$\lambda_X:H/E\times{H/E}\to\mathbb{Z}[\pi]$, which we shall call the {\it
reduced} intersection pairing.
The adjoint homomorphism 
$\tilde\lambda_X:H/E\to(H/E)^\dagger$
is given by $\tilde\lambda_X([v])([u])=\lambda(u,v)=ev(v)(D(u))$, 
for all $u,v\in{H}$.
It is a monomorphism, and $\lambda_X$ is {\it nonsingular} if
$\tilde\lambda_X$ is an isomorphism.

\begin{lemma}
Let $X$ be a $PD_4$-complex with fundamental group $\pi$, and let
$E=\overline{H^2(\pi;\mathbb{Z}[\pi])}$.

\begin{enumerate}

\item If $\lambda_X$ is nonsingular then 
$\overline{H^3(\pi;\mathbb{Z}[\pi])}$ embeds as a submodule of $E^\dagger$;

\item if $\lambda_X$ is nonsingular and $H^2(c_X;\mathbb{Z}[\pi])$ 
splits then $E^\dagger\cong\overline{H^3(\pi;\mathbb{Z}[\pi])}$;

\item if $H^3(\pi;\mathbb{Z}[\pi])=0$ then $\lambda_X$ is nonsingular;

\item if $H^3(\pi;\mathbb{Z}[\pi])=0$ and $\Pi$ is a finitely generated
projective $\mathbb{Z}[\pi]$-module then $E=0$;

\item if $\pi$ is infinite and $H^1(\pi;\mathbb{Z}[\pi])$ and $\Pi$ are projective then $c.d.\pi=4$.
\end{enumerate}
\end{lemma}

\begin{proof}
Let $p:\Pi\to\Pi/D(E)$ and ${q:H\to{H/E}}$ be the canonical epimorphisms.
Poincar\'e duality induces an isomorphism $\gamma:H/E\cong\Pi/D(E)$.
It is straightforward to verify that
$p^\dagger(\gamma^\dagger)^{-1}\tilde\lambda_Xq=ev$.
If $\lambda_X$ is nonsingular then $\tilde\lambda_X$ is an isomorphism,
and so $\mathrm{Coker}(p^\dagger)=\mathrm{Coker}(ev)$.
The first assertion follows easily, since  
$\mathrm{Coker}(p^\dagger)\leq{E^\dagger}$.

If moreover $H^2(c_X;\mathbb{Z}[\pi])$ splits then so does $p$, and so
$E^\dagger\cong\mathrm{Coker}(p^\dagger)$.

If $H^3(\pi;\mathbb{Z}[\pi])=0$ then $ev$ is an epimorphism 
and so $p^\dagger$ is an epimorphism.
Since $p^\dagger$ is also a monomorphism it is an isomorphism.
Since $ev$ and $q$ are epimorphisms with the same kernel it folows that
$\tilde\lambda_X=\gamma^\dagger(p^\dagger)^{-1}$,
and so $\tilde\lambda_X$ is also an isomorphism.

If $\Pi$ is finitely generated and projective 
then so is $\Pi^\dagger$,
and $\Pi\cong\Pi^{\dagger\dagger}$.
If moreover $H^3(\pi;\mathbb{Z}[\pi])=0$ then
$\Pi\cong H\cong{E}\oplus\Pi^\dagger$.
Hence $E$ is also finitely generated and projective, 
and $E\cong E^{\dagger\dagger}=0$.

If $H^1(\pi;\mathbb{Z}[\pi])$ and $\Pi$ are projective then we may obtain 
a projective resolution of $\mathbb{Z}$ of length 4
from $C_*=C_*(\widetilde{X})$ by replacing $C_3$ and $C_4$ 
by $C_3\oplus\Pi$ and $C_4\oplus\overline{H^1(\pi;\mathbb{Z}[\pi])}$,
respectively, and modifying $\partial_3$ 
and $\partial_4$ appropriately.
Since $H_3(X;\mathbb{Z}[\pi])\cong\overline{H^1(\pi;\mathbb{Z}[\pi])}$ 
it is also projective.
It follows from the UCSS that $H^4(\pi;\mathbb{Z}[\pi])\not=0$. 
Hence $c.d.\pi=4$.
\end{proof}

In particular, the cohomology intersection pairing
is nonsingular if and only if $H^2(\pi;\mathbb{Z}[\pi])=
H^3(\pi;\mathbb{Z}[\pi])=0$.
If $X$ is a 4-manifold counting intersections of generic immersions
of $S^2$ in $\widetilde{X}$ gives an equivalent pairing on $\Pi$. 

We do not know whether the hypotheses in this lemma can be simplified.
For instance, is $H^2(\pi;\mathbb{Z}[\pi])^\dagger$ always 0?
Does ``$\Pi$ projective" imply that $H^3(\pi;\mathbb{Z}[\pi])=0$?
Projectivity of $\Pi^\dagger$ and $H^2(\pi;\mathbb{Z}[\pi])=0$ 
together do not imply this.
For if $\pi$ is a $PD_3^+$-group and $w=w_1(\pi)$ 
there are finitely generated projective $\mathbb{Z}[\pi]$-modules
$P$ and $P'$ such that $\Pi\oplus{P}\cong A(\pi)\oplus{P'}$,
where $A(\pi)$ is the augmentation ideal of $\mathbb{Z}[\pi]$,
by Theorem 3.13, and so $\Pi^\dagger$ is projective.
However $H^3(\pi;\mathbb{Z}[\pi])\cong\mathbb{Z}\not=0$.

The module $\Pi$ is finitely generated if and only if $\pi$ is of type $FP_3$.
As observed in the proof of Theorem 2.18,
if $\pi$ is a free product of infinite cyclic groups and groups with one end
and is not a free group then $H^1(\pi;\mathbb{Z}[\pi])$ is a free
$\mathbb{Z}[\pi]$-module.
An argument similar to that for part(5) of the lemma shows that $c.d.\pi\leq5$
if and only if $\pi$ is torsion-free and $p.d._{\mathbb{Z}[\pi]}\Pi\leq2$.

If $Y$ is a second $PD_4$-complex we write $\lambda_X\cong\lambda_Y$
if there is an isomorphism $\theta:\pi\cong\pi_1(Y)$ such that
$w_1(X)=w_1(Y)\theta$ and a $\mathbb{Z}[\pi]$-module isomorphism 
$\Theta:\pi_2(X)\cong\theta^*\pi_2(Y)$
inducing an isometry of cohomology intersection pairings.
If $f:X\to Y$ is a 2-connected degree-1 map the ``surgery kernel" 
$K_2(f)=\mathrm{Ker}(\pi_2(f))$ and ``surgery cokernel"
$K^2(f)=\overline{\mathrm{Cok}(H^2(f;\mathbb{Z}[\pi]))}$ 
are finitely generated and projective, 
and are stably free if $X$ and $Y$ are finite complexes \cite[Lemma 2.2]{[Wl]}.
(See also Theorem 3.2 above.)
Moreover cap product with $[X]$ induces an isomorphism from $K^2(f)$ 
to $K_2(f)$.
The pairing $\lambda_f=\lambda|_{K^2(f)\times K^2(f)}$ 
is nonsingular \cite[Theorem 5.2]{[Wl]}.

A comprehensive account of the pairing $\lambda_X$ and its role as a key invariant of the homotopy type of a $PD_4$-complex may be found in
\cite{[Hi20]}. 
The key results are briefly summarized in \S10.8 below.
The paper \cite{[BBH18]} extends the work of \cite{[Cr00],[Hi12],[Hi17']} 
to higher dimensions. 
In particular, it shows that if $X$ is a $PD_4$-complex such that $\pi_1(X)$ 
is torsion free and $\pi_2(X)=0$ then $X$ is homotopy equivalent to 
a connected sum of aspherical $PD_4$-complexes and copies of
$S^3\times{S^1}$ and $S^3\tilde\times{S^1}$.

%% file: m5-4.tex
\chapter{Mapping tori and circle bundles}

Stallings showed that if $M$ is a 3-manifold and $f:M\to S^1$ a map which 
induces an epimorphism $f_*:\pi_1(M)\to\mathbb{Z}$ with infinite kernel $K$
then $f$ is homotopic to a bundle projection if and only if 
$M$ is irreducible and $K$ is finitely generated.
Farrell gave an analogous characterization in dimensions $\geq6$,
with the hypotheses that the homotopy fibre of $f$ is finitely dominated 
and a torsion invariant $\tau(f)\in Wh(\pi_1(M))$ is 0.
The corresponding results in dimensions 4 and 5 are constrained by the
present limitations of geometric topology in these dimensions.
(In fact there are counter-examples to the most natural 4-dimensional 
analogue of Farrell's theorem \cite{[We87]}.)

Quinn showed that if the base $B$ and homotopy fibre $F$ of a fibration
$p:M\to{B}$ are finitely dominated then the total space $M$ is a Poincar\'e 
duality complex 
if and only if both the base and fibre are Poincar\'e duality complexes.
(The paper \cite{[Go79]} gives an elegant proof for the case when $M$, 
$B$ and $F$ are finite complexes. The general case follows
on taking products with copies of $S^1$ to reduce to the finite case
and using the K\"unneth theorem.)

We shall begin by giving a purely homological proof of a version of this result,
for the case when $M$ and $B$ are $PD$-spaces and $B=K(G,1)$ is aspherical.
The homotopy fibre $F$ is then the covering space associated to the kernel of
the induced epimorphism from $\pi_1(M)$ to $G$.
Our algebraic approach requires only that the equivariant chain complex of $F$ 
have finite $[n/2]$-skeleton.
In the next two sections we use the finiteness criterion of Ranicki
and the fact that Novikov rings associated to finitely generated groups
are weakly finite to sharpen this finiteness hypotheses when $B=S^1$, 
corresponding to infinite cyclic covers of $M$.
The main result of \S4.4 is a 4-dimensional homotopy fibration theorem 
with hypotheses similar to those of Stallings and a conclusion 
similar to that of Gottlieb and Quinn.
The next two sections consider products of 3-manifolds with $S^1$ 
and covers associated to ascendant subgroups.

We shall treat fibrations of $PD_4$-complexes over surfaces in Chapter 5,
by a different, more direct method.
In the final section of this chapter we consider instead bundles 
with fibre $S^1$.
We give conditions for a $PD_4$-complex to fibre over a $PD_3$-complex
with homotopy fibre $S^1$,
and show that these conditions are sufficient
if the fundamental group of the base is torsion-free but not free. 

\section{$PD_r$-covers of $PD_n$-spaces}

Let $M$ be a $PD_n$-space and $p:\pi=\pi_1(M)\to G$ an epimorphism with $G$ 
a $PD_r$-group, and let $M_\nu$ be the covering space corresponding to
$\nu=\mathrm{Ker}(p)$.
If $M$ is aspherical and $\nu$ is $FP_{[n/2]}$ then $\nu$ 
is a $PD_{n-r}$-group and $M_\nu=K(\nu,1)$ is a $PD_{n-r}$-space 
\cite[Theorem 9.11]{[Bi]}.
In general, there are isomorphisms
$H^q(M_\nu;\mathbb{Z}[\nu])\cong{H_{n-r-q}}(M_\nu;\mathbb{Z}[\nu])$,
by Theorem 1.19$^\prime$.
However in the nonaspherical case it is not clear that there are
such isomorphisms induced by cap product with a class in
$H_{n-r}(M_\nu;\mathbb{Z}[\nu])$.
If $M$ is a $PD_n$-complex and $\nu$ is finitely presentable
$M_\nu$ is finitely dominated, and we could apply the Gottlieb-Quinn Theorem
to conclude that $M_\nu$ is a $PD_{n-r}$-complex.
We shall give instead a purely homological argument which does not require
$\pi$ or $\nu$ to be finitely presentable, 
and so applies under weaker finiteness hypotheses.

A group $G$ is a {\it weak} $PD_r$-group if $H^r(G;\mathbb{Z}[G])$ 
is infinite cyclic if $q=r$ and is $0$ otherwise \cite{[Ba80]}.
If $r\leq2$ an $FP_2$ group is a weak $PD_r$-group 
if and only if it is virtually a $PD_r$-group. 
This is easy for $r\leq1$ and is due to Bowditch when $r=2$ \cite{[Bo04]}.
Barge has given a simple homological argument to show that 
if $G$ is a weak $PD_r$-group, $M$ is a $PD_n$-space and
$\eta_G\in H^r(M;\mathbb{Z}[G])$ is the image of a generator of
$H^r(G;\mathbb{Z}[G])$ then cap product with $[M_\nu]=\eta_G\cap[M]$ 
induces isomorphisms with simple coefficients \cite{[Ba80]}.
We shall extend his argument to the case of arbitrary local coefficients,
using coinduced modules to transfer arguments about subgroups 
and covering spaces to contexts where Poincar\'e duality applies. 

All tensor products $N\otimes{P}$ in the following theorem are taken over 
$\mathbb{Z}$.

\begin{theorem}
Let $M$ be a $PD_n$-space and $p:\pi=\pi_1(M)\to G$ an epimorphism with $G$ 
a weak $PD_r$-group, and let $\nu=\mathrm{Ker}(p)$.
If $C_*(\widetilde{M})$ is $\mathbb{Z}[\nu]$-finitely dominated then
$M_\nu$ is a $PD_{n-r}$-space.
\end{theorem}

\begin{proof}
Let $C_*$ be a finitely generated projective $\mathbb{Z}[\pi]$-chain complex
which is chain homotopy equivalent to $C_*(\widetilde{M})$.
Since $C_*(\widetilde{M})$ is $\mathbb{Z}[\nu]$-finitely dominated there is a 
finitely generated projective $\mathbb{Z}[\nu]$-chain complex $E_*$ and a pair 
of $\mathbb{Z}[\nu]$-linear chain homomorphisms $\theta:E_*\to{C_*|_\nu}$ and
$\phi:{C_*|_\nu}\to E_*$ such that $\theta\phi\sim{I}_{C_*}$ and
$\phi\theta\sim{I}_{E_*}$.
Let $C^q=Hom_{\mathbb{Z}[\pi]}(C_q,\mathbb{Z}[\pi])$ and 
$E^q=Hom_{\mathbb{Z}[\nu]}(E_q,\mathbb{Z}[\nu])$,
and let $\widehat{\mathbb{Z}[\pi]}=
Hom_{\mathbb{Z}[\nu]}(\mathbb{Z}[\pi]|_\nu,\mathbb{Z}[\nu])$
be the module coinduced from $\mathbb{Z}[\nu]$.
Then there are isomorphisms 
$\Psi:H^q(E^*)\cong{H^q}(C_*;\widehat{\mathbb{Z}[\pi]})$,
determined by $\theta$ and Shapiro's Lemma.

The complex $\mathbb{Z}[G]\otimes_{\mathbb{Z}[\pi]}C_*$ is an augmented 
complex of finitely generated projective $\mathbb{Z}[G]$-modules
with finitely generated integral homology.
Therefore $G$ is of type $FP_\infty$ \cite[Theorem 3.1]{[St96]}.
Hence the augmentation $\mathbb{Z}[G]$-module
$\mathbb{Z}$ has a resolution $A_*$ by finitely generated projective 
$\mathbb{Z}[G]$-modules.
Let $A^q= Hom_{\mathbb{Z}[G]}(A_q,\mathbb{Z}[G])$ 
and let $\eta\in{H}^r(A^*)=H^r(G;\mathbb{Z}[G])$ be a generator.
Let $\varepsilon_C:C_*\to A_*$ be a chain map corresponding to the
projection of $p$ onto $G$,
and let $\eta_G=\varepsilon_C^*\eta\in{H^r}(C_*;\mathbb{Z}[G])$.
The augmentation $A_*\to\mathbb{Z}$ determines a chain homotopy equivalence 
$p:C_*\otimes{A}_*\to{C_*}\otimes\mathbb{Z}=C_*$.
Let $\sigma:G\to\pi$ be a set-theoretic section.

We may define cup-products relating the cohomology of $M_\nu$ and $M$,
as follows.
Let $e:\widehat{\mathbb{Z}[\pi]}\otimes\mathbb{Z}[G]\to\mathbb{Z}[\pi]$
be the pairing given by $e(\alpha\otimes{g})=\sigma(g).\alpha(\sigma(g)^{-1})$
for all $\alpha:\mathbb{Z}[\pi]\to\mathbb{Z}[\nu]$ and $g\in{G}$.
Then $e$ is independent of the choice of section $\sigma$ and is
${\mathbb{Z}[\pi]}$-linear with respect to the diagonal left
$\pi$-action on $\widehat{\mathbb{Z}[\pi]}\otimes\mathbb{Z}[G]$.
Let $d:C_*\to C_*\otimes{C_*}$ be a $\pi$-equivariant diagonal,
with respect to the diagonal left $\pi$-action on $C_*\otimes{C_*}$,
and let $j=(1\otimes\varepsilon_C)d:C_*\to {C_*\otimes{A}_*}$.
Then $pj=Id_{C_*}$ and so $j$ is a chain homotopy equivalence.
We define the cup-product $[f]\cup\eta_G$ in
$H^{p+r}(C^*)=H^{p+r}(M;\mathbb{Z}[\pi])$ by 
$[f]\cup\eta_G=e_\#d^*(\Psi([f])\times\eta_G)
=e_\#j^*(\Psi([f])\times\eta)$ for all
$[f]\in{H^p}(E^*)=H^p(M_\nu;\mathbb{Z}[\nu])$. 

If $C$ is a left $\mathbb{Z}[\pi]$-module let 
$D=Hom_{\mathbb{Z}[\nu]}(C|_\nu,\mathbb{Z}[\pi])$
have the left $G$-action determined by 
$(g\lambda)(c)=\sigma(g)\lambda(\sigma(g)^{-1}c)$
for all $c\in{C}$ and $g\in{G}$.
If $C$ is free with basis $\{c_i|1\leq{i}\leq{n}\}$ there is 
an isomorphism of left $\mathbb{Z}[G]$-modules 
$\Theta:D\cong{\mathbb{Z}[\pi]^n|_1}^G$
given by $\Theta(\lambda)(g)=(\sigma(g).\lambda(\sigma(g)^{-1}c_1),
\dots ,\sigma(g).\lambda(\sigma(g)^{-1}c_n))$, 
for all $\lambda\in{D}$ and $g\in{G}$, 
and so $D$ is coinduced from a module over the trivial group.

Let $D^q=Hom_{\mathbb{Z}[\nu]}(C_q|_\nu,\mathbb{Z}[\pi])$ 
and let $\rho:E^*\otimes\mathbb{Z}[G]\to{D^*}$ be the 
$\mathbb{Z}$-linear cochain homomorphism defined by
$\rho(f\otimes{g})(c)=\sigma(g)f\phi(\sigma(g)^{-1}c)$ for all 
$c\in{C_q}$, $\lambda\in{D^q}$, $f\in{E^q}$, $g\in{G}$ and all $q$.
Then the $G$-action on $D^q$ and $\rho$ are independent of the choice 
of section $\sigma$,
and $\rho$ is $\mathbb{Z}[G]$-linear if $E^q\otimes\mathbb{Z}[G]$ 
has the left $G$-action given by $g(f\otimes{g'})=f\otimes{gg'}$
for all $g,g'\in{G}$ and $f\in{E^q}$.

If $\lambda\in{D^q}$ then $\lambda\theta_q(E_q)$ is a finitely generated 
$\mathbb{Z}[\nu]$-submodule of $\mathbb{Z}[\pi]$. 
Hence there is a family of homomorphisms $\{f_g\in{E^q}|g\in{F}\}$,
where $F$ is a finite subset of $G$,
such that 
$\lambda\theta_q(e)=\Sigma_{g\in{F}}f_g(e)\sigma(g)$ for all $e\in{E}_q$.
Let $\lambda_g(e)=
\sigma(g)^{-1}f_g(\phi\sigma(g)\theta(e))\sigma(g)$ for all $e\in{E}_q$ 
and $g\in{F}$.
Let $\Phi(\lambda)=\Sigma_{g\in{F}}{\lambda_g}\otimes{g}
\in{E}^q\otimes\mathbb{Z}[G]$.
Then $\Phi$ is a $\mathbb{Z}$-linear cochain homomorphism.
Moreover $[\rho\Phi(\lambda)]=[\lambda]$ for all $[\lambda]\in{H}^q(D^*)$
and $[\Phi\rho(f\otimes{g})]=[f\otimes{g}]$ for all 
$[f\otimes{g}]\in{H}^q(E^*\otimes\mathbb{Z}[G])$,
and so $\rho$ is a chain homotopy equivalence.
(It is not clear that $\Phi$ is $\mathbb{Z}[G]$-linear on the cochain level,
but we shall not need to know this).

We now compare the hypercohomology of $G$ with coefficients in the 
cochain complexes $E^*\otimes\mathbb{Z}[G]$ and $D^*$.
On one side we have ${\mathbb{H}^n(G;E^*\otimes\mathbb{Z}[G])}
={H}^n_{tot}(Hom_{\mathbb{Z}[G]}(A_*,E^*\otimes\mathbb{Z}[G]))$,
which may be identified with ${H}^n_{tot}(E^*\otimes{A}^*)$
since $A_q$ is finitely generated for all $q\geq0$.
This is in turn isomorphic to
$H^{n-r}(E^*)\otimes{H^r}(G;\mathbb{Z}[G])\cong{H^{n-r}}(E^*)$, 
since $G$ acts trivially on $E^*$ and is a weak $PD_r$-group.

On the other side we have 
${\mathbb{H}^n(G;D^*)}={H}^n_{tot}(Hom_{\mathbb{Z}[G]}(A_*,D^*))$.
The cochain homomorphism $\rho$ induces a morphism of double complexes 
from $E^*\otimes{A}^*$ to $Hom_{\mathbb{Z}[G]}(A_*,D^*)$ by
$\rho^{pq}(f\otimes\alpha)(a)=\rho(f\otimes\alpha(a))\in{D}^p$
for all $f\in{E}^p$, $\alpha\in{A}^q$ and $a\in{A}_q$ and all $p,q\geq0$.
Let $\hat\rho^p([f])=[\rho^{pr}(f\times\eta)]\in 
\mathbb{H}^{p+r}(G;{D^*})$ for all $[f]\in{H}^p(E^*)$.
Then $\hat\rho^p:H^p(E^*)\to\mathbb{H}^{p+r}(G;D^*)$ is an isomorphism, 
since $[f]\mapsto[f\times\eta]$
is an isomorphism and $\rho$ is a chain homotopy equivalence.
Since $C_p$ is a finitely generated projective $\mathbb{Z}[\pi]$-module
$D^p$ is a direct summand of a coinduced module.
Therefore $H^i(G;D^p)=0$ for all $i>0$, while 
$H^0(G;D^p)={Hom_{\mathbb{Z}[\pi]}(C_p,\mathbb{Z}[\pi])}$, for all $p\geq0$.
Hence $\mathbb{H}^n(G;{D^*})\cong{H^n(C^*)}$ for all $n$.

Let $f\in {E}^p$, $a\in A_r$ and $c\in{C}_p$, and suppose that
$\eta(a)=\Sigma{n_g}g$.
Since $\hat\rho^p([f])(a)(c)=\rho(f\otimes\eta(a))(c)=
\Sigma{n_g}\sigma(g)f\phi(\sigma(g)^{-1}c)=
([f]\cup\eta)(c,a)$
it follows that 
the homomorphisms from $H^p(E^*)$ to $H^{p+r}(C^*)$ given by cup-product 
with $\eta_G$ are isomorphisms for all $p$.

Let $[M]\in H_n(M;\mathbb{Z}^w)$ be a fundamental class for $M$,
and let $[M_\nu]=\eta_G\cap[M]\in{H}_{n-r}(M;\mathbb{Z}^w\otimes\mathbb{Z}[G])
=H_{n-r}(M_\nu;\mathbb{Z}^{w|_\nu})$.
Then cap product with $[M_\nu]$ induces isomorphisms
$H^p(M_\nu;\mathbb{Z}[\nu])\cong{H}_{n-r-p}(M_\nu;\mathbb{Z}[\nu])$ for all $p$,
since ${c\cap[M_\nu]}={(c\cup\eta_G)\cap[M]}$ 
in $H_{n-r-p}(M;\mathbb{Z}[\pi])=H_{n-r-p}(M_\nu;\mathbb{Z}[\nu])
=H_{n-r-p}(\widetilde{M};\mathbb{Z})$
for $c\in{H}^p(M_\nu;\mathbb{Z}[\nu])$.
Thus $M_\nu$ is a $PD_{n-r}$-space.
\end{proof}

Theorems 1.19$^\prime$ and 4.1 together give the following version of 
the Gottlieb-Quinn Theorem for covering spaces.

\begin{cor}
Let $M$ be a $PD_n$-space and $p:\pi=\pi_1(M)\to G$ an epimorphism with $G$ 
a $PD_r$-group, and let $\nu=\mathrm{Ker}(p)$.
Then $M_\nu$ is a $PD_{n-r}$-space if and only if 
$C_*(\widetilde{M})|_\nu$ has finite $[n/2]$-skeleton.
\end{cor}

\begin{proof}
The conditions are clearly necessary.
Conversely, if $M_\nu$ has finite $[n/2]$-skeleton 
then $C_*$ is $\mathbb{Z}[\nu]$-finitely dominated, by Theorem 1.19$^\prime$,
and hence is a $PD_{n-r}$-space, by Theorem 4.1.
\end{proof}

\begin{cor}
The space $M_\nu$ is a $PD_{n-r}$-complex if and only if it is homotopy 
equivalent to a complex with finite $[n/2]$-skeleton and 
$\nu$ is finitely presentable.
\qed
\end{cor}

\begin{cor}
If $\pi$ is a $PD_r$-group $\widetilde{M}$ is a $PD_{n-r}$-complex 
if and only if $H_q(\widetilde{M};\mathbb{Z})$ is finitely
generated for all $q\leq[n/2]$.
\qed
\end{cor}

Stark used \cite[Theorem 3.1]{[St96]} with the Gottlieb-Quinn Theorem 
to deduce that if $M$ is a $PD_n$-complex and $v.c.d.\pi/\nu<\infty$ 
then $\pi/\nu$ is of type $vFP$, and therefore is virtually a $PD$-group.
Is there a purely algebraic argument to show that if $M$ is a $PD_n$-space, 
$\nu$ is a normal subgroup of $\pi$ and $C_*(\widetilde{M})$ is 
$\mathbb{Z}[\nu]$-finitely dominated then $\pi/\nu$ must be a weak $PD$-group?

\section{Novikov rings and Ranicki's criterion}

The results of the above section apply in particular when $G=\mathbb{Z}$.
In this case however we may use an alternative finiteness criterion of
Ranicki to get a slightly stronger result, 
which can be shown to be best possible. 
The results of this section are based on
joint work with Kochloukova (in \cite{[HK07]}).

Let $\pi$ be a group, $\rho:\pi\to\mathbb{Z}$ an epimorphism with kernel $\nu$ 
and $t\in\pi$ an element such that $\rho(t)=1$.
Let $\alpha:\nu\to\nu$ be the automorphism determined by $\alpha(h)=tht^{-1}$ 
for all $h$ in $\nu$.
This automorphism extends to a ring automorphism (also denoted by $\alpha$)
of the group ring $R=\mathbb{Z}[\nu]$, and the ring $S=\mathbb{Z}[\pi]$ 
may then be viewed as a twisted Laurent extension, 
$\mathbb{Z}[\pi]=\mathbb{Z}[\nu]_\alpha [t,t^{-1} ]$.
The {\it Novikov ring} $\widehat{\mathbb{Z}[\pi]}_\rho$ associated to $\pi$ 
and $\rho$ is the ring of (twisted) Laurent series
$\Sigma_{j\geq a}\kappa_jt^{j}$, for some $a\in\mathbb{Z}$,
with coefficients $\kappa_j$ in $\mathbb{Z}[\nu]$.
Multiplication of such series is determined by conjugation in $\pi$:
if $g\in{\nu}$ then $tg=(tgt^{-1})t$.
If $\pi$ is finitely generated the Novikov rings 
$\widehat{\mathbb{Z}[\pi]}_\rho$ are weakly finite \cite{[Ko06]}.
Let $\widehat{S}_+=\widehat{\mathbb{Z}[\pi]}_\rho$ and 
$\widehat{S}_-=\widehat{\mathbb{Z}[\pi]}_{-\rho}$.

An $\alpha$-{\it twisted} endomorphism of an $R$-module $E$
is an additive function $h:E\to E$ such that 
$h(re)=\alpha(r)h(e)$ for all $e\in{E}$ and $r\in{R}$,
and $h$ is an $\alpha$-twisted automorphism if it is bijective. 
Such an endomorphism $h$ extends to $\alpha$-twisted endomorphisms 
of the modules
$S\otimes_RE$, $\widehat{E}_+=\widehat{S}_+\otimes_RE$ and
$\widehat{E}_-=\widehat{S}_-\otimes_RE$ by
$h(s\otimes{e})=tst^{-1}\otimes h(e)$ for all $e\in{E}$ and $s\in{S}$,
$\widehat{S}_+$ or $\widehat{S}_-$, respectively.
In particular, left multiplication by $t$ determines $\alpha$-twisted 
automorphisms of $S\otimes_RE$, $\widehat{E}_{+}$ and $\widehat{E}_{-}$
which commute with $h$.

If $E$ is finitely generated then $1-t^{-1}h$ is an automorphism 
of $\widehat{E}_-$, with inverse given by a geometric series:
$(1-t^{-1}h)^{-1}=\Sigma_{k\geq0}t^{-k}h^k$.
(If $E$ is not finitely generated this series may not give a function 
with values in $\widehat{E}_-$, and ${t-h}=t(1-t^{-1}h)$ may not be surjective).
Similarly, if $k$ is an $\alpha^{-1}$-twisted endomorphism of $E$
then $1-tk$ is an automorphism of $\widehat{E}_+$.

If $P_*$ is a chain complex with an endomorphism $\beta:P_*\to P_*$
let $P_*[1]$ be the suspension and $\mathcal{C}(\beta)_*$ be the mapping cone.
Thus $\mathcal{C}(\beta)_q=P_{q-1}\oplus P_q$, and 
$\partial_q(p,p')=(-\partial p,\beta(p)+\partial p')$, and 
there is a short exact sequence
\[
0\to P_*\to \mathcal{C}(\beta)_*\to P_*[1]\to0.
\]
The connecting homomorphisms in the associated long exact sequence of homology
are induced by $\beta$.
The {\it algebraic mapping torus} of an $\alpha$-twisted self chain homotopy
equivalence $h$ of an $R$-chain complex $E_*$ is the mapping cone
$\mathcal{C}(1-t^{-1}h)$ of the endomorphism $1-t^{-1}h$ of the $S$-chain 
complex $S\otimes_RE_*$.

\begin{lemma}
Let $E_*$ be a projective chain complex over $R$ which is finitely generated 
in degrees $\leq d$ and let $h:E_*\to E_*$ be an $\alpha$-twisted chain 
homotopy equivalence. 
Then $H_q({\widehat{S}_-\otimes_S\mathcal{C}(1-t^{-1}h)_*})=0$ for $q\leq d$.
\end{lemma}

\begin{proof}
There is a short exact sequence
\[
0\to S\otimes_RE_*\to \mathcal{C}(1-z^{-1}h)_*\to S\otimes_RE_*[1]\to0.
\]
Since $E_*$ is a complex of projective $R$-modules the sequence
\[
0\to\widehat{E}_{*-}\to{\widehat{S}_-\otimes_S\mathcal{C}(1-t^{-1}h)_*}
\to\widehat{E}_{*-}[1]\to0
\]
obtained by extending coefficients is exact.
Since $1-t^{-1}h$ induces isomorphisms on $\widehat{E}_{q-}$ for
$q\leq{d}$ it induces isomorphisms on homology in degrees $<d$ 
and an epimorphism on homology in degree $d$.
Therefore
$H_q({\widehat{S}_-\otimes_S\mathcal{C}(1-t^{-1}h)_*})=0$ for $q\leq d$,
by the long exact sequence of homology.
\end{proof}

The next theorem is our refinement of Ranicki's finiteness criterion 
\cite{[HK07]}.

\begin{theorem}
Let $C_*$ be a finitely generated projective $S$-chain complex.
Then $i^!C_*$ has finite $d$-skeleton
if and only if $H_q(\widehat{S}_\pm\otimes_SC_*)=0$ for $q\leq{d}$.
\end{theorem}

\begin{proof}
We may assume without loss of generality that $C_q$ is a 
finitely generated free
$S$-module for all $q\leq{d+1}$, with basis $X_i=\{c_{q,i}\}_{i\in I(q)}$.
We may also assume that $0 \notin \partial_i(X_i)$ for $i \leq d+1$,
where $\partial_i:C_i \rightarrow C_{i-1}$ is the differential of the complex.
Let $h_\pm$ be the $\alpha^{\pm1}$-twisted automorphisms of $i^!C_*$ 
induced by multiplication by $z^{\pm1}$ in $C_*$.
Let 
\[f_q(z^krc_{q,i})=(0,z^k\otimes{r}c_{q,i})\in
(S\otimes_RC_{q-1})\oplus(S\otimes_RC_q).
\]
Then $f_*$ defines $S$-chain homotopy equivalences from $C_*$ to each of
$\mathcal{C}(1-z^{-1}h_+)$ and $\mathcal{C}(1-zh_-)$.

Suppose first that $k_*:i^!C_*\to E_*$ and $g_*:E_*\to i^!C_*$ are chain 
homotopy equivalences, where $E_*$ is a projective $R$-chain complex
which is finitely generated in degrees $\leq{d}$.
Then $\theta_\pm=k_*h_\pm g_*$ are $\alpha^{\pm1}$-twisted self homotopy 
equivalences of $E_*$,
and $\mathcal{C}(1-z^{-1}h_+)$ and $\mathcal{C}(1-zh_-)$ 
are chain homotopy 
equivalent to $\mathcal{C}(1-z^{-1}\theta_+)$ and $\mathcal{C}(1-z\theta_-)$,
respectively.
Therefore 
$H_q(\widehat{S}_-\otimes_SC_*)=
{H_q}(\widehat{S}_-\otimes_S\mathcal{C}(1-z^{-1}\theta_+))=0$
and $H_q(\widehat{S}_+\otimes_SC_*)=
{H_q}(\widehat{S}_+\otimes_S\mathcal{C}(1-z\theta_-))=0$
for $q\leq{d}$, by Lemma 5, applied twice.

Conversely, suppose that $H_i(\widehat{S}_{\pm} \otimes_S C_*)= 0$ for all
$i \leq k$.
Adapting an idea from \cite{[BR88]},
we shall define inductively a support function $supp_{X}$ 
for $\lambda\in\cup_{i \leq d+1}C_i$ 
with values finite subsets of $\{ z^j \}_{j \in \mathbb{Z}}$ so that

\begin{enumerate}
\item $supp_X(0)=\emptyset$;

\item if $x\in{X}_0$ then  $supp_{X}(z^jx) = z^j$;

\item if $x\in{X}_i$ for $1\leq{i}\leq{d+1}$ then 
$supp_{X}(z^j x) = z^j .supp_{X}(\partial_i(x))$;

\item if $s = \sum_j r_j z^j \in S$,  where $r_j \in R$, 
$supp_X(sx) = \cup_{r_j \not= 0 } supp_X(z^j x)$;

\item if $0 \leq i \leq{d+1}$ and 
$\lambda=\sum_{s_x \in S, x \in X_i}s_x x$ then

$supp_{X}(\lambda)=\cup_{s_x\not=0,x\in{X_i}}supp_{X}(s_x x)$.
\end{enumerate}

Then $supp_X(\partial_i(\lambda)) \subseteq supp_X(\lambda)$
for all $\lambda\in{C_i}$ and all $ 1 \leq i \leq{d+1}$.
Since $X=\cup_{i\leq{d+1}} X_i$ is finite 
there is a positive integer $b$ such that 
\[
\cup_{x \in X_i, i \leq {d+1}} supp_{X}(x) 
\subseteq \{ z^j \}_{-b \leq j \leq b}.
\]
Define two subcomplexes $C^+$ and $C^-$ of $C$ which are 0 in degrees 
$ i \geq{d+2}$ as follows:

\begin{enumerate}
\item if $i \leq{d+1}$ an element $\lambda \in C_i $ is in $C^+$ 
if and only if $supp_{X}(\lambda) \subseteq \{z^j \}_{j \geq -b }$; and   

\item
if $i \leq{d+1}$ an element $\lambda \in C_i$ is in $C^-$ if and only if 
$supp_{X}(\lambda) \subseteq \{z^j \}_{j \leq b }$.

\end{enumerate}

Then $\cup_{i \leq d+1} X_i \subseteq (C^+)^{[d+1]} \cap (C^-)^{[d+1]} $
and so $(C^+)^{[d+1]} \cup (C^-)^{[d+1]} = C^{[d+1]},$ 
where the upper index $*$ denotes the $*$-skeleton. 
Moreover
$(C^+)^{[d+1]}$ is a complex of free finitely generated $R_{\alpha}[z]$-modules,
$(C^-)^{[d+1]}$ is a complex of free finitely generated 
$R_{\alpha}[z^{-1}]$-modules,
$(C^+)^{[d+1]} \cap (C^-)^{[d+1]}$
is a complex of free finitely generated $R$-modules and 
\[
C^{[d+1]} = S\otimes_{R_{\alpha}[z]} (C^+)^{[d+1]} = 
S\otimes_{R_{\alpha}[z^{-1}]} (C^-)^{[d+1]}.
\]
Furthermore there is a Mayer-Vietoris exact sequence
\[
0 \rightarrow (C^+)^{[d+1]} \cap (C^-)^{[d+1]}
\rightarrow (C^+)^{[d+1]} \oplus (C^-)^{[d+1]} \rightarrow C^{[d+1]} 
\rightarrow 0.
\]
Thus the $(d+1)$-skeletons of $C$, $C^+$ and $C^-$ satisfy
``algebraic transversality" in the sense of \cite{[Rn95]}.
To prove the theorem it suffices to show that $C^+$ and $C^-$ 
are each chain homotopy equivalent over $R$ to a complex of projective 
$R$-modules which is finitely generated in degrees $ \leq{d}$. 
As in \cite{[Rn95]}
there is an exact sequence of $R_{\alpha}[z^{-1}]$-module chain complexes
\[
0 \rightarrow (C^-)^{[d+1]}\rightarrow{C}^{[d+1]} 
\oplus R_{\alpha}[[z^{-1}]] 
\otimes_{R_{\alpha}[z^{-1}]} (C^-)^{[d+1]} \rightarrow
\widehat{S}_- \otimes_{S} C^{[d+1]} \rightarrow 0.
\]
Let $\tilde{i}$ denot the inclusion of $(C^-)^{[d+1]}$ into the central term.
Inclusions on each component define a chain homomorphism
\[
\tilde{j} :  (C^+)^{[d+1]} \cap (C^-)^{[d+1]}
\rightarrow (C^+)^{[d+1]} \oplus R_{\alpha}[[z^{-1}]] 
\otimes_{R_{\alpha}[z^{-1}]} (C^-)^{[d+1]}
\]
such that the mapping cones of $\tilde{i}$ and $\tilde{j}$ are chain equivalent $R$-module chain complexes.
The map induced by $\tilde{i}$ in homology is an epimorphism in degree $d$ 
and an isomorphism in degree $ < d$,
since $H_i (\widehat{S}_- \otimes_{S} C^{[d+1]}) = 0$ for $ i \leq{d}$.
In particular all homologies in degrees $ \leq{d}$ of the mapping cone 
of $\tilde{i}$ are 0. 
Hence all homologies of the mapping cone of $\tilde{j}$ are 0 in degrees 
$\leq d$. 
Then $(C^+)^{[d+1]}$ is homotopy equivalent over $R$ to a chain complex 
of projectives over $R$ whose $k$-skeleton is a summand of 
$(C^+)^{[d]} \cap (C^-)^{[d]}$. This completes the proof.
\end{proof}

The argument for the converse is entirely due to Kochloukova.

As an application we shall give a quick proof of Kochloukova's improvement of
Corollary 2.5.1.

\begin{cor} {\rm[Ko06]}
Let $\pi$ be a finitely presentable group with a finitely generated 
normal subgroup $N$ such that $\pi/N\cong\mathbb{Z}$.
Then $\mathrm{def}(\pi)=1$ if and only if $N$ is free.
\end{cor}

\begin{proof}
Let $X$ be the finite 2-complex
corresponding to an optimal presentation of $\pi$.
If $\mathrm{def}(G)=1$ then $\chi(X)=0$ and $X$ is aspherical, by Theorem 2.5.
Hence $C_*=C_*(\widetilde{X})$ is a finite free resolution of the augmentation
module $\mathbb{Z}$.
Let $A_\pm$ be the two Novikov rings corresponding to the two epimorphisms
$\pm{p}:\pi\to\mathbb{Z}$ with kernel $N$.
Then $H_j(A_\pm\otimes_{\mathbb{Z}[\pi]}{C}_*) = 0$ for $ j \leq 1$, 
by Theorem 4.3. 
But then $H_2(A_\pm\otimes_{\mathbb{Z}[\pi]} C_*)$ is stably free, by Lemma 3.1. 
Since $\chi(A_\pm\otimes_{\mathbb{Z}[\pi]} C_*)=\chi(C_*)=\chi(X)=0$ and
the rings $A_\pm$ are weakly finite \cite{[Ko06]} these modules are 0.
Thus $H_j(A_\pm\otimes_{\mathbb{Z}[\pi]}C_*)=0$ for all $j$,
and so $C_*|_\nu$ is chain homotopy equivalent to a finite projective 
$\mathbb{Z}[\nu]$-complex \cite[Theorem 2]{[Rn95]}.
In particular, $N$ is $FP_2$ and hence is free \cite[Corollary 8.6]{[Bi]}.

The converse is clear.
\end{proof}

\section{Infinite cyclic covers}

The {\it mapping torus} of a self homotopy equivalence $f:X\to X$ is the space
$M(f)=X\times [0,1]/\sim$, where $(x,0)\sim (f(x),1)$ for all $x\in X$.
The function $p([x,t])=e^{2\pi{i}t}$ defines a map $p:M(f)\to{S^1}$ with
homotopy fibre $X$, and the induced homomorphism $p_*:\pi_1(M(f))\to\mathbb{Z}$ 
is an epimorphism if $X$ is path-connected.
Conversely, let $E$ be a connected cell complex and let $f:E\to S^1 $ be a map 
which induces an epimorphism $f_*:\pi_1(E)\to\mathbb{Z}$, 
with kernel $\nu$.
Then $E_\nu =E\times_{S^1 } R=\{ (x,y)\in E\times R\mid f(x)=e^{2\pi iy}\}$,
and $E\simeq M(\phi)$, where 
$\phi:E_\nu \to E_\nu$ is the generator of the covering group given by
$\phi(x,y)=(x,y+1)$ for all $(x,y)$ in $E_\nu$.

\begin{theorem}
Let $M$ be a finite $PD_n$-space with fundamental group $\pi$ and let  
$p:\pi\to\mathbb{Z}$ be an epimorphism with kernel $\nu$.
Then $M_\nu$ is a $PD_{n-1}$-space if and only if $\chi(M)=0$ and 
$C_*(\widetilde{M_\nu})=C_*(\widetilde{M})|_\nu$ has finite $[(n-1)/2]$-skeleton.
\end{theorem}

\begin{proof}
If $M_\nu$ is a $PD_{n-1}$-space then $C_*(\widetilde{M_\nu})$ is 
$\mathbb{Z}[\nu]$-finitely dominated \cite{[Br72]}.
In particular, $H_*(M;\Lambda)=H_*(M_\nu;\mathbb{Z})$ is finitely generated.
The augmentation $\Lambda$-module $\mathbb{Z}$ has a short free resolution
$0\to\Lambda\to\Lambda\to\mathbb{Z}\to0$, 
and it follows easily from the exact sequence of homology 
for this coefficient sequence that $\chi(M)=0$ \cite{[Mi68]}.
Thus the conditions are necessary.

Suppose that they hold.
Let $A_\pm$ be the two Novikov rings corresponding to the two epimorphisms
$\pm{p}:\pi\to\mathbb{Z}$ with kernel $\nu$.
Then $H_j(A_\pm\otimes_{\mathbb{Z}[\pi]}{C}_*) = 0$ for $ j \leq [(n-1)/2]$, 
by Theorem 4.3. 
Hence $H_j(A_\pm\otimes_{\mathbb{Z}[\pi]} C_*) = 0$ for $j \geq n - [(n-1)/2]$,
by duality.
If $n$ is even there is one possible nonzero module, in degree $m=n/2$.
But then $H_m(A_\pm\otimes_{\mathbb{Z}[\pi]} C_*)$ 
is stably free, by the finiteness of $M$ and Lemma 3.1. 
Since $\chi(A_\pm\otimes_{\mathbb{Z}[\pi]} C_*)=\chi(C_*)=\chi(M)=0$ and
the rings $A_\pm$ are weakly finite \cite{[Ko06]} these modules are 0.
Thus $H_j(A_\pm\otimes_{\mathbb{Z}[\pi]}C_*)=0$ for all $j$,
and so $C_*|_\nu$ is chain homotopy equivalent to a finite projective 
$\mathbb{Z}[\nu]$-complex, by Theorem 4.4.
Thus the result follows from Theorem 4.1.
\end{proof}

When $n$ is odd $[n/2]=[(n-1)/2]$, so the finiteness condition on $M_\nu$ 
agrees with that of Corollary 4.1.1 (for $G=\mathbb{Z}$), 
but it is slightly weaker if $n$ is even.
Examples constructed by elementary surgery on simple $n$-knots
show that the $FP_{[(n-1)/2]}$ condition is best possible, 
even when $\pi\cong\mathbb{Z}$ and $\nu=1$. 

\begin{cor}
Under the same hypotheses on $M$ and $\pi$, 
if $n\not=4$ then $M_\nu$ is a $PD_{n-1}$-complex if and only if it is
homotopy equivalent to a complex with finite $[(n-1)/2]$-skeleton.
\end{cor}

\begin{proof}
If $n\leq3$ every $PD_{n-1}$-space is a $PD_{n-1}$-complex, 
while if $n\geq5$ then $[(n-1)/2]\geq2$ and so $\nu$ is finitely presentable.
\end{proof}

If $n\leq3$ we need only assume that $M$ is a $PD_n$-space and
$\nu$ is finitely {\it generated}.

It remains an open question whether every $PD_3$-space is finitely dominated.
The arguments of \cite{[Tu90]} and \cite{[Cr00]} on the factorization 
of $PD_3$-complexes into connected sums are essentially homological, 
and so every $PD_3$-space is a connected sum of aspherical $PD_3$-spaces and
a $PD_3$-complex with virtually free fundamental group.
Thus the question of whether every $PD_3$-space is finitely dominated reduces 
to whether every $PD_3$-group is finitely presentable.

\section{The case $n=4$}

If $M(f)$ is the mapping torus of a self homotopy equivalence of a
$PD_3$-space then $\chi(M)=0$ and $\pi_1(M)$ is an extension of $\mathbb{Z}$ 
by a finitely generated normal subgroup.
These conditions characterize such mapping tori, by Theorem 4.4.
We shall summarize various related results in the following theorem. 

\begin{theorem} 
Let $M$ be a finite $PD_4$-space whose fundamental group $\pi$ is 
an extension of $Z$ by a finitely generated normal subgroup $\nu$.
Then
\begin{enumerate}                                                                                     
\item $\chi(M)\geq0$, with equality if and only if $H_2(M_\nu;\mathbb{Q})$ is 
finitely generated;

\item $\chi(M)=0$ if and only if $M_\nu$ is a $PD_3$-space;

\item if $\chi(M)=0$ then $M$ is aspherical if and only if $\nu$ 
is a $PD_3$-group if and only if $\nu$ has one end;

\item if $M$ is aspherical then $\chi(M)=0$ if and only if $\nu$ is a 
$PD_3$-group if and only if $\nu$ is $FP_2$.
\end{enumerate}
\end{theorem}

\begin{proof} Since $C_*(\widetilde{M})$ is finitely dominated and 
$\mathbb{Q}\Lambda=\mathbb{Q}[t,t^{-1}]$ is noetherian
the homology groups $H_q(M_\nu;\mathbb{Q})$ 
are finitely generated as $\mathbb{Q}\Lambda$-modules.
Since $\nu$ is finitely generated they are finite dimensional as 
$\mathbb{Q}$-vector spaces if $q<2$, 
and hence also if $q>2$, by Poincar\'e duality. 
Now $H_2(M_\nu;\mathbb{Q})\cong \mathbb{Q}^r\oplus(\mathbb{Q}\Lambda)^s$ 
for some $r,s\geq0$, by the Structure Theorem for modules over a PID.
It follows easily from the Wang sequence for the covering projection from 
$M_\nu$ to $M$, that $\chi(M)=s\geq0$.

The space $M_\nu$ is a $PD_3$-space if and only if $\chi(M)=0$, 
by Theorem 4.4.

Since $M$ is aspherical if and only if $M_\nu$ is aspherical,
(3) follows from (2) and the facts that $PD_3$-groups have one end and
a $PD_3$-space is aspherical if and only if its fundamental group 
has one end. 

If $M$ is aspherical and $\chi(M)=0$ then $\nu$ is a $PD_3$-group.
If $\nu$ is a $PD_3$-group it is $FP_2$.
If $M$ is aspherical and $\nu$ is $FP_2$ then $\nu$ is a $PD_3$-group,
by Theorem 1.19 (or Theorem 4.4), and so $\chi(M)=0$.
\end{proof}

In particular, if $\chi(M)=0$ then $q(\pi)=0$. 
This observation and the bound $\chi(M)\geq0$ were given in Theorem 3.17.
(They also follow on counting bases for the cellular chain complex of 
$M_\nu$ and extending coefficients to $\mathbb{Q}(t)$.)

If $\chi(M)=0$ and $\nu$ is finitely presentable then 
$M_\nu$ is a $PD_3$-complex.
However $M_\nu$ need not be homotopy equivalent to a {\sl finite} complex. 
If $M$ is a {\it simple} $PD_4$-complex and a generator of 
$Aut(M_\nu/M)\cong\pi/\nu$ has finite order in the group of self homotopy 
equivalences of $M_\nu$ then $M$ is finitely covered by 
a simple $PD_4$-complex homotopy equivalent to $M_\nu\times S^1$.
In this case $M_\nu$ must be homotopy finite \cite{[Rn86]}.

If $\pi\cong\nu\rtimes\mathbb{Z}$ is a $PD_4$-group with 
$\nu$ finitely generated then $\chi(\pi)=0$ if and only if $\nu$ is $FP_2$, 
by Theorem 4.5.
However the latter conditions need not hold.
Let $F$ be the orientable surface of genus 2. 
Then $G=\pi_1(F)$ has a presentation
$\langle a_1,a_2,b_1,b_2\mid [a_1,b_1]=[a_2,b_2]\rangle$. 
The group $\pi=G\times G$ is a $PD_4$-group,
and the subgroup $\nu\leq\pi$ generated by the images of $(a_1,a_1)$ and
the six elements $(x,1)$ and $(1,x)$, for $x=a_2$, $b_1$ or $b_2$, 
is normal in $\pi$, with quotient $\pi/\nu\cong\mathbb{Z}$.
However $\chi(\pi)=4\not=0$ and so $\nu$ cannot be $FP_2$.

It can be shown that the finitely generated subgroup $N$ of $F(2)\times F(2)$ 
defined after Theorem 2.4 has one end.
However $H^2(F(2)\times F(2);\mathbb{Z}[F(2)\times F(2)])\not=0$.
(Note that $q(F(2)\times F(2))=2$, by Corollary 3.12.2.)

\begin{cor}
Let $M$ be a finite $PD_4$-space with $\chi(M)=0$ and 
whose fundamental group $\pi$ is an extension of $\mathbb{Z}$ 
by a normal subgroup $\nu$.
If $\pi$ has an infinite cyclic normal subgroup $C$ which 
is not contained in $\nu$ then the covering space $M_\nu$ 
with fundamental group $\nu$ is a $PD_3 $-complex.
\end{cor}

\begin{proof} We may assume without loss of generality that $M$ is 
orientable and that $C$ is central in $\pi$. 
Since $\pi/\nu$ is torsion-free $C\cap\nu=1$,
and so $C\nu\cong C\times\nu$ has finite index in $\pi$.  
Thus by passing to a finite cover we may assume that $\pi=C\times\nu$.
Hence $\nu$ is finitely presentable and so Theorem 4.5 applies.
\end{proof}

Since $\nu$ has one or two ends if it has an infinite cyclic normal subgroup,
Corollary 4.5.1 remains true if $C\leq\nu$ and $\nu$ is finitely presentable.
In this case $\nu$ is the fundamental group of a Seifert fibred 3-manifold,
by Theorem 2.14.

\begin{cor}
Let $M$ be a finite $PD_4$-space with $\chi(M)=0$ and 
whose fundamental group $\pi$ is an extension of $\mathbb{Z}$ by 
a finitely generated normal subgroup $\nu$.
If $\nu$ is finite then it has cohomological period dividing $4$.
If $\nu$ has one end then $M$ is aspherical and so $\pi$ is a $PD_4 $-group.
If $\nu$ has two ends then $\nu\cong\mathbb{Z}$, 
$\mathbb{Z}\oplus (Z/2Z)$ or $D$.
If moreover $\nu$ is finitely presentable the covering space $M_\nu$ with 
fundamental group $\nu$ is a $PD_3 $-complex.
\end{cor}
                     
\begin{proof} 
The final hypothesis is only needed if $\nu$ is one-ended, 
as finite groups and groups with two ends are finitely presentable. 
If $\nu$ is finite then $\widetilde M\simeq S^3$ and so the first assertion 
holds.
(See Chapter 11 for more details.)
If $\nu$ has one end we may use Theorem 4.5.
If $\nu$ has two ends and its maximal finite normal subgroup is nontrivial
then $\nu\cong\mathbb{Z}\oplus(Z/2Z)$, by Theorem 2.11 (
applied to the $PD_3$-complex $M_\nu$).
Otherwise $\nu\cong\mathbb{Z}$ or $D$.
\end{proof}

In Chapter 6 we shall strengthen this Corollary to obtain a fibration theorem
for 4-manifolds with torsion-free elementary amenable fundamental group.

\begin{cor}
Let $M$ be a finite $PD_4$-space with $\chi(M)=0$ and whose 
fundamental group $\pi$ is an extension of $\mathbb{Z}$ 
by a normal subgroup $\nu\cong F(r)$.
Then $M$ is homotopy equivalent to a closed PL $4$-manifold which fibres over 
the circle, with fibre $\sharp^r S^1\times S^2$ if $w_1 (M)|_\nu $ is trivial, 
and $\sharp^r S^1 \tilde\times S^2$ otherwise.
The bundle is determined by the homotopy type of $M$.
\end{cor} 

\begin{proof}
Since $M_\nu$ is a $PD_3$-complex with free fundamental group
it is homotopy equivalent to $N=\sharp^rS^1\times{S^2}$ if
$w_1 (M)|_\nu$ is trivial and to $\sharp^rS^1\tilde\times{S^2}$ otherwise.
Every self homotopy equivalence of a connected sum of $S^2$-bundles over $S^1$
is homotopic to a self-homeomorphism, and homotopy implies isotopy
for such manifolds \cite{[La]}.
Thus $M$ is homotopy equivalent to such a fibred 4-manifold,
and the bundle is determined by the homotopy type of $M$.
\end{proof}

The  homotopy types of such mapping tori are determined by $\pi$,
$w_1(M)$ and the orbit of $w_2(M)$ under the action of $Out(\pi)$.
It is easy to see that $Homeo(N)$ maps onto $Out(F(r)$,
and all such triples $(\pi,w_1,w_2)$ are realized \cite{[Hi20]}.

\begin{cor}
Let $M$ be a finite $PD_4$-space with $\chi(M)=0$ and whose fundamental group
$\pi$ is an extension of $\mathbb{Z}$ by a torsion-free normal subgroup 
$\nu$ which is the fundamental group of a closed 3-manifold $N$.
Then $M$ is homotopy equivalent to the mapping torus of a self 
homeomorphism of $N$.
\end{cor}
     
\begin{proof} 
There is a homotopy equivalence $f:N\to M_\nu$, by Turaev's Theorem. 
(See \S5 of Chapter 2.)
The indecomposable factors of $N$ are either Haken, hyperbolic 
or Seifert fibred $3$-manifolds, 
by the Geometrization Conjecture (see \cite{[B-P]}).
Let $t:M_\nu\to M_\nu$ be the generator of the covering transformations. 
Then there is a self homotopy equivalence $u:N\to N$ such that $fu\sim tf$. 
As each aspherical factor of $N$ has the property that self homotopy 
equivalences are homotopic to PL homeomorphisms 
(by \cite{[Hm]}, Mostow rigidity or \cite{[Sc83]}),
and a similar result holds for $\sharp^r(S^1\times S^2)$ 
(by \cite{[La]}), 
$u$ is homotopic to a homeomorphism \cite{[HL74]}, 
and so $M$ is homotopy equivalent to the mapping torus of this homeomorphism. 
\end{proof}

The hypothesis that $M$ be finite is redundant in each of the last two
corollaries, since $\tilde{K}_0(\mathbb{Z}[\pi])=0$. 
(See Theorem 6.3.)
All known $PD_3 $-complexes with torsion-free fundamental group are 
homotopy equivalent to 3-manifolds.

If the irreducible connected summands of the closed 3-manifold 
$N=\sharp_i N_i$ are $P^2$-irreducible and sufficiently large 
or have fundamental group $\mathbb{Z}$ then every self homotopy equivalence 
of $N$ is realized by an unique isotopy class 
of homeomorphisms \cite{[HL74]}. 
However if $N$ is not aspherical then it admits
nontrivial self-homeomorphisms (``rotations about 2-spheres") which induce the
identity on $\nu$, and so such bundles are not determined by the group alone.

Let $f:M\to E$ be a homotopy equivalence, where $E$ is a finite
$PD_4$-complex with $\chi(E)=0$ and fundamental group
$\pi\cong\nu\rtimes\mathbb{Z}$,
where $\nu$ is finitely presentable.
Then $w_1(M)=f^*w_1(E)$ and $c_{E*}f_*[M]=\pm{c_{E*}[E]}$ in 
$H_4(\pi;\mathbb{Z}^{w_1(E)})$.
Conversely, if $\chi(M)=0$ and there is an isomorphism
$\theta:\pi_1(M)\cong\pi$ such that $w_1(M)=\theta_i^*w$ and 
$\theta_{1*}c_{M*}[M]=c_{E*}[E]$ then $E_\nu$ and $M_\nu$ are $PD_3$-complexes, 
by Theorem 4.5.
A Wang sequence argument as in the next theorem shows that the fundamental
triples of $E_\nu$ and $M_\nu$ are isomorphic, 
and so they are homotopy equivalent, by Hendrik's Theorem.
What additional conditions are needed to determine the homotopy type of
such mapping tori?
Our next result is a partial step in this direction.

\begin{theorem} 
Let $E$ be a finite $PD_4 $-complex with $\chi(E)=0$ and whose fundamental 
group $\pi$ is an extension of $\mathbb{Z}$ by a finitely presentable 
normal subgroup $\nu$ which is not virtually free.
Let $\Pi=\overline{H^2(\pi;\mathbb{Z}[\pi])}$.
A $PD_4 $-complex $M$ is homotopy equivalent to $E$ if and only if $\chi(M)=0$,
there is an isomorphism $\theta$ from $\pi_1(M)$ to $\pi$ such that 
$w_1(M)=w_1(E)\theta $, ${\theta^*}^{-1}k_1(M)$ and $k_1(E)$ generate the same 
subgroup of $H^3(\pi;\Pi)$ under the action of 
${Out(\pi)\times{Aut_{\mathbb{Z}[\pi]}(\Pi)}}$,
and there is a lift $\hat c:M\to P_2 (E)$ of $\theta c_M $ 
such that $\hat c_* [M]=\pm f_{E*} [E]$ in $H_4(P_2(E);\mathbb{Z}^{w_1(E)})$.
\end{theorem}

\begin{proof} The conditions are clearly necessary.
Suppose that they hold.
The infinite cyclic covering spaces $N=E_\nu$ and $M_\nu$ are $PD_3$-complexes,
by Theorem 4.5, and $\pi_2(E)\cong\Pi$ and $\pi_2(M)\cong\theta^*\Pi$,
by Theorem 3.4.
The maps $c_N $ and $c_E $ induce a homomorphism between the Wang sequence
for the fibration of $E$ over $S^1$ and the corresponding Wang sequence for
$K(\pi,1)$.
Since $\nu$ is not virtually free 
$H_3(c_N;\mathbb{Z}^{w_1(E)})$ is a monomorphism.
Hence $H_4(c_E;\mathbb{Z}^{w_1(E)})$ and {\it a fortiori} 
$H_4(f_E;\mathbb{Z}^{w_1(E)})$ are monomorphisms, 
and so Theorem 3.8 applies.
\end{proof}

As observed in the first paragraph of \S9 of Chapter 2,
the conditions on $\theta$ and the $k$-invariants also imply that
$M_\nu\simeq{E_\nu}$.
 
The original version of this book gave an exposition of the extension of
Barge's argument to local coefficients for the case when $G\cong\mathbb{Z}$,
instead of the present Theorem 4.1, 
and used this together with an $L^2$-argument,
instead of the present Theorem 4.3, to establish the results corresponding 
to Theorem 4.5 for the case when $\nu$ was $FP_2$.
         
\section{Products}

If $M=N\times S^1$, where $N$ is a closed 3-manifold,
then $\chi(M)=0$, $\mathbb{Z}$ is a direct factor of $\pi_1(M)$,
$w_1(M)$ is trivial on this factor
and the $Pin^-$-condition $w_2=w_1^2$ holds.
These conditions almost characterize such products up to homotopy equivalence.
We need also a constraint on the other direct factor of the fundamental group.

\begin{theorem} 
Let $M$ be a finite $PD_4$-complex whose fundamental group $\pi$ 
has no $2$-torsion.
Then $M$ is homotopy equivalent to a product $N\times S^1$,
where $N$ is a closed $3$-manifold, if and only if $\chi (M)=0$, 
$w_2(M)=w_1(M)^2$ and there is an isomorphism 
$\theta:\pi\to\nu\times\mathbb{Z}$ such that 
$w_1(M)\theta^{-1}|_{\mathbb{Z}}=0$, 
where $\nu$ is a ($2$-torsion-free) $3$-manifold group.
\end{theorem}

\begin{proof} The conditions are clearly necessary,
since the $Pin^-$-condition holds for 3-manifolds.
 
If these conditions hold then the covering space $M_\nu$ with 
fundamental group $\nu$ is a $PD_3 $-complex, by Theorem 4.5 above.
Since $\nu$ is a 3-manifold group and has no 2-torsion it is
a free product of cyclic groups and groups of aspherical closed 3-manifolds.
Hence there is a homotopy equivalence $h:M_\nu\to N$,
where $N$ is a connected sum of lens spaces and aspherical closed 3-manifolds,
by Turaev's Theorem.
(See \S5 of Chapter 2.)
Let $\phi $ generate the covering group $Aut(M/M_\nu)\cong\mathbb{Z}$.
Then there is a self homotopy equivalence $\psi:N\to N$ such that
$\psi h\sim h\phi$, and $M$ is homotopy equivalent to the mapping torus 
$M(\psi)$.
We may assume that $\psi$ fixes a basepoint and induces the 
identity on $\pi_1(N)$, 
since $\pi_1(M)\cong\nu\times\mathbb{Z}$.                                  
Moreover $\psi$ preserves the local orientation, 
since $w_1(M)\theta^{-1}|_{\mathbb{Z}}=0$. 
Since $\nu$ has no element of order 2 there are no 
two-sided projective planes in $N$,
and so $\psi$ is homotopic to a rotation about a 2-sphere \cite{[Hn]}. 
Since $w_2(M)=w_1(M)^2$ the rotation is homotopic to the identity and so 
$M$ is homotopy equivalent to $N\times S^1$. 
\end{proof}

Let $\rho $ be an essential map from $S^1 $ to $SO(3)$,
and let $M=M(\tau)$, where $\tau:S^2\times{S^1}\to{S^2}\times{S^1}$ is 
the twist map, given by $\tau (x,y)=(\rho (y)(x),y)$ 
for all $(x,y)$ in $S^2\times{S^1}$.
Then $\pi_1(M)\cong\mathbb{Z}\times\mathbb{Z}$, $\chi(M)=0$, and $w_1(M)=0$,
but $w_2 (M)\not=w_1(M)^2=0$, so $M$ is not homotopy equivalent to a product. 
(Clearly however $M(\tau^2 )=S^2\times{S^1}\times{S^1}$.)

To what extent are the constraints on $\nu$ necessary?
There are orientable 4-manifolds which are homotopy equivalent to 
products $N\times S^1$ where $\nu=\pi_1 (N)$ is finite and is
{\it not} a 3-manifold group.
(See Chapter 11.)
Theorem 4.1 implies that $M$ is homotopy equivalent to a product of an aspherical 
$PD_3 $-complex with $S^1$ if and only if $\chi (M)=0$ and 
$\pi_1 (M)\cong\nu\times\mathbb{Z}$ where $\nu $ has one end.  

There are 4-manifolds which are simple homotopy equivalent to 
$S^1\times RP^3$ (and thus satisfy the hypotheses of our theorem) but which 
are not homeomorphic to mapping tori \cite{[We87]}.

\section{Ascendant subgroups}

In this brief section we shall give another characterization of aspherical 
$PD_4$-complexes with finite covering spaces which are homotopy equivalent to 
mapping tori.

\begin{theorem}
Let $M$ be a $PD_4$-complex.
Then $M$ is aspherical and has a finite cover which is homotopy equivalent 
to a mapping torus if and only if $\chi(M)=0$ and $\pi=\pi_1(M)$ 
has an ascendant $FP_3$ subgroup $G$ of infinite index
and such that $H^s(G;\mathbb{Z}[G])=0$ for $s\leq2$.
In that case $G$ is a $PD_3$-group,
$[\pi:N_\pi(G)]<\infty$ and $e(N_\pi(G)/G)=2$.
\end{theorem}

\begin{proof} 
The conditions are clearly necessary.
Suppose that they hold and that
 $G=G_0<G_1<\dots<G_\beth=\pi$ is an ascendant sequence.
Let $\gamma={\min\{ \alpha\mid [G_\alpha:G]=\infty\}}$.
Transfinite induction using the LHSSS with coefficients 
$\mathbb{Z}[\pi]$ and Theorem 1.15 shows that 
$H^s (\pi;\mathbb{Z}[\pi])=0$ for $s\leq2$.
If $\gamma$ is finite then $\beta_1^{(2)}(G_\gamma)=0$,
since it has a finitely generated normal subgroup 
of infinite index \cite{[Ga00]}.
Otherwise $\gamma$ is the first infinite ordinal,
and $[G_{j+1}:G_j]<\infty$ for all $j<\gamma$.
In this case $\beta_1^{(2)}(G_n)=\beta_1^{(2)}(G)/[G_n:G]$ and so
$\lim_{n\to\infty}\beta_1^{(2)}(G_n)=0$.
It then follows from \cite[Theorems 6.13 and 6.54(7)]{[Lu]} that
$\beta_1^{(2)}(G_\gamma)=0$.
In either case it then follows that $\beta_1^{(2)}(G_\alpha)=0$ 
for all $\gamma\leq\alpha\leq\beth$,
by Theorem 2.3 (which is part of \cite[Theorem 7.2]{[Lu]}).
Hence $M$ is aspherical, by Theorem 3.5.
                        
On the other hand $H^s(G_\gamma;W)=0$ for $s\leq3$ and any free 
$\mathbb{Z}[G_\gamma]$-module $W$, so $c.d.G_\gamma=4$.
Hence $[\pi:G_\gamma]<\infty$, by Strebel's Theorem.
Therefore $G_\gamma$ is a $PD_4$-group.
In particular, it is finitely generated and so $\gamma<\infty$.
If $\gamma=\beta+1$ then $[G_\beta:G]<\infty$.
It follows easily that $[\pi:N_\pi(G)]<\infty$.
Hence $G$ is a $PD_3$-group and $N_\pi(G)/G$ has two ends, 
by Theorem 3.10.
\end{proof}

The hypotheses on $G$ could be replaced by
``{\sl $G$ is a $PD_3$-group}", 
for then ${[\pi:G]=\infty}$, by Theorem 3.12.

Theorem 5.8 below gives an analogue for $PD_4$-complexes $M$ 
such that $\chi(M)=0$ and $\pi_1(M)$ has an ascendant subgroup 
of infinite index which is a $PD_2$-group.

\section{Circle bundles}

In this section we shall consider the ``dual" situation, of $PD_4$-complexes
which are homotopy equivalent to the total space of a $S^1$-bundle over 
a 3-dimensional base $N$.
Lemma 4.9 presents a number of conditions satisfied by such spaces.
(These conditions are not all independent.)
Bundles $c_N^*\xi$ induced from $S^1$-bundles over $K(\pi_1(N),1)$ 
are given equivalent characterizations in Lemma 4.10.
In Theorem 4.11 we shall show that the conditions of Lemmas 4.9 and 4.10 
characterize the homotopy types of such bundle spaces $E(c_N^*\xi)$, 
provided $\pi_1(N)$ is torsion-free but not free.
    
Since $BS^1\simeq K(\mathbb{Z},2)$ any $S^1 $-bundle over a connected base $B$ 
is induced from some bundle over $P_2(B)$.
For each epimorphism $\gamma :\mu\to\nu $ with cyclic kernel and such that 
the action of $\mu$ by conjugation on $\mathrm{Ker}(\gamma)$ 
factors through multiplication by $\pm1$ there is an $S^1 $-bundle 
$p(\gamma ):X(\gamma )\to Y(\gamma )$ whose fundamental group sequence 
realizes $\gamma$ and which is universal for such bundles; 
the total space $E(p(\gamma))$ is a $K(\mu ,1)$ space.
(See \cite[Proposition 11.4]{[Wl]}). 
                           
\begin{lemma} 
Let $p:E\to B$ be the projection of an $S^1 $-bundle $\xi $ 
over a connected cell complex $B$. Then
\begin{enumerate}
\item $\chi (E)=0$;

\item the natural map $p_* :\pi=\pi_1 (E)\to\nu=\pi_1 (B)$ is an 
epimorphism with cyclic kernel, 
and the action of $\nu $ on {\rm Ker}$(p_*)$ induced by conjugation 
in $\pi $ is given by 
$w=w_1 (\xi ):\pi_1 (B)\to\{\pm 1\}\leq {Aut(\mathrm{Ker}(p_*))}$; 

\item if $B$ is a $PD$-complex $w_1 (E)=p^*(w_1 (B)+w)$; 

\item if $B$ is a $PD_3 $-complex there are maps $\hat c:E\to P_2 (B)$ and 

\noindent $y:P_2 (B)\to Y(p_*)$ such that $c_{P_2 (B)} =c_{Y(p_*)} y$, 
$y\hat c=p(p_*)c_E $ and $(\hat c,c_E )_* [E]=\pm G(f_{B*} [B])$, 
where $G$ is the Gysin homomorphism from
$H_3 (P_2 (B);\mathbb{Z}^{w_1 (B)} )$ to 
$H_4 (P_2 (E);\mathbb{Z}^{w_1 (E)} )$;

\item If $B$ is a $PD_3 $-complex $c_{E*}[E]=\pm G(c_{B*}[B])$, 
where $G$ is the Gysin homomorphism from 
$H_3(\nu;\mathbb{Z}^{w_1(B)})$ to $H_4 (\pi ;\mathbb{Z}^{w_1(E)})$; 

\item $\mathrm{Ker}(p_*)$ acts trivially on $\pi_2 (E)$.
\end{enumerate}
\end{lemma}
          
\begin{proof}
Condition (1) follows from the multiplicativity of 
the Euler characteristic in a fibration.
If $\alpha $ is any loop in $B$ the total space of the induced bundle 
$\alpha^* \xi $ is the torus if $w(\alpha )=1$ and the Klein bottle 
if $w(\alpha )=-1$; hence $gzg^{-1} =z^{w(p_* (g))}$
for $g$ in $\pi_1 (E)$ and $z$ in $\mathrm{Ker}(p_*)$. 
Conditions (2) and (6) then follow from the exact homotopy sequence.
If the base $B$ is a $PD$-complex then so is $E$,
and we may use naturality and the Whitney sum formula (applied to the Spivak normal bundles)
to show that $w_1 (E) =p^* (w_1 (B) +w_1 (\xi ))$. 
(As $p^* :H^1 (B;\mathbb{F}_2)\to H^1 (E;\mathbb{F}_2)$ is 
a monomorphism this equation determines $w_1 (\xi )$.) 
        
Condition (4) implies (5), and follows from the observations in the
paragraph preceding the lemma.
(Note that the Gysin homomorphisms $G$ in (4) and (5) are well defined, 
since $H_1(\mathrm{Ker}(\gamma);\mathbb{Z}^{w_1(E)}$) 
is isomorphic to $\mathbb{Z}^{w_1(B)}$, by (3).) 
\end{proof}

Bundles with $\mathrm{Ker}(p_*)\cong\mathbb{Z}$ have the following equivalent characterizations.

\begin{lemma} 
Let $p:E\to B$ be the projection of an $S^1 $-bundle $\xi $ 
over a connected cell complex $B$. 
Then the following conditions are equivalent:
\begin{enumerate}
\item $\xi $ is induced from an $S^1 $-bundle over $K(\pi_1 (B),1)$ via $c_B $;

\item for each map $\beta :S^2 \to B$ the induced bundle $\beta^* \xi $ is trivial; 

\item the induced epimorphism $p_* :\pi_1 (E)\to\pi_1 (B)$ has 
infinite cyclic kernel.
\end{enumerate}
    
\noindent If these conditions hold then $c(\xi )=c_B^* \Xi $, 
where $c(\xi )$ is the characteristic class of $\xi $ 
in $H^2 (B;Z^w )$ and $\Xi $ is the class of the 
extension of fundamental groups 
in $H^2 (\pi_1 (B);\mathbb{Z}^w )=H^2 (K(\pi_1 (B),1);\mathbb{Z}^w )$, 
where $w=w_1 (\xi )$.
\end{lemma}

\begin{proof} Condition (1) implies condition (2) 
as for any such map $\beta $
the composite $c_B \beta $ is nullhomotopic. 
Conversely, as we may construct $K(\pi_1 (B),1)$ by adjoining cells 
of dimension $\geq 3$ to $B$ condition (2) implies 
that we may extend $\xi $ over the 3-cells, 
and as $S^1 $-bundles over $S^n $ are trivial 
for all $n>2$ we may then extend $\xi $ over the whole of $K(\pi_1 (B),1)$, 
so that (2) implies (1). 
The equivalence of (2) and (3) follows on observing that (3) holds
if and only if $\partial\beta=0$ for all such $\beta $, 
where $\partial $ is the
connecting map from $\pi_2 (B)$ to $\pi_1 (S^1 )$ in the exact sequence of homotopy
for $\xi $, and on comparing this with the corresponding sequence for $\beta^* \xi $.

As the natural map from the set of $S^1 $-bundles over $K(\pi ,1)$ 
with $w_1 =w$ (which are classified by $H^2 (K(\pi ,1);\mathbb{Z}^w )$) 
to the set of extensions of $\pi $ by $\mathbb{Z}$ with $\pi $ acting via $w$ 
(which are classified by $H^2 (\pi ;\mathbb{Z}^w )$) 
which sends a bundle to the extension of fundamental groups 
is an isomorphism we have $c(\xi )=c_B^* (\Xi )$. 
\end{proof}

If $N$ is a closed 3-manifold which has no summands of type $S^1 \times S^2 $ 
or $S^1\tilde\times S^2$ (i.e., if $\pi_1(N)$ has no infinite cyclic free 
factor) then every $S^1 $-bundle over $N$ with $w=0$ restricts to a trivial 
bundle over any map from $S^2 $ to $N$. 
For if $\xi $ is such a bundle, with characteristic class $c(\chi )$ in
$H^2 (N;\mathbb{Z})$, and $\beta :S^2 \to N$ is any map then
$\beta_* (c(\beta^* \xi )\cap [S^2 ])=\beta_* (\beta^* c(\xi )\cap [S^2 ])
=c(\xi )\cap\beta_* [S^2 ]=0$,
as the Hurewicz homomorphism is trivial for such $N$. 
Since $\beta_* $ is an isomorphism in degree 0 it follows that 
$c(\beta^* \xi )=0$ and so $\beta^* \xi $ is trivial. 
(A similar argument applies for bundles with $w\not= 0$, 
provided the induced 2-fold covering space $N^w $ has no summands 
of type $S^1 \times S^2 $ or $S^1 \tilde\times S^2 $.)

On the other hand, if $\eta $ is the Hopf fibration the bundle with total space 
$S^1 \times S^3 $, base $S^1 \times S^2 $ and projection $id_{S^1 }\times\eta$ 
has nontrivial pullback over any essential map from $S^2$ to $S^1 \times S^2$, 
and is not induced from any bundle over $K(\mathbb{Z},1)$. 
Moreover, $S^1 \times S^2 $ is a 2-fold covering space of $RP^3 \sharp RP^3 $, 
and so the above hypothesis on summands of $N$ is not stable
under passage to 2-fold coverings
(corresponding to a homomorphism $w$ from $\pi_1 (N)$ to $Z/2Z$).

\begin{theorem} 
Let $M$ be a $PD_4 $-complex and $N$ a $PD_3 $-complex whose fundamental group 
is torsion-free but not free.
Then $M$ is homotopy equivalent to the total space of an $S^1$-bundle over $N$ 
which satisfies the conditions of Lemma $4.10$ if and only if
\begin{enumerate}
\item $\chi (M)=0$;

\item there is an epimorphism $\gamma:\pi=\pi_1(M)\to\nu=\pi_1(N)$ 
with $\mathrm{Ker}(\gamma)\cong\mathbb{Z}$;

\item $w_1 (M)=(w_1 (N)+w)\gamma $, where 
$w:\nu\to{Aut(\mathrm{Ker}(\gamma))}$ is 
determined by the action of $\nu $ on $\mathrm{Ker}(\gamma)$ 
induced by conjugation in $\pi $;

\item  $k_1(M)=\gamma^* k_1(N)$ 
(and so $P_2(M)\simeq P_2(N)\times_{K(\nu,1)}K(\pi,1)$);

\item $f_{M*}[M]=\pm G(f_{N*}[N])$ in $H_4(P_2(M);\mathbb{Z}^{w_1(M)})$,
where $G$ is the Gysin homomorphism in degree $3$.
\end{enumerate}
\noindent If these conditions hold then $M$ has minimal Euler characteristic
for its fundamental group, i.e., $q(\pi)=0$.
\end{theorem}

{\bf Remark}\qua
The first three conditions and Poincar\'e duality imply that 
$\pi_2(M)\cong \gamma^*\pi_2(N)$, the $\mathbb{Z}[\pi]$-module with the same 
underlying group as $\pi_2(N)$ and with $\mathbb{Z}[\pi]$-action 
determined by the homomorphism $\gamma$.

\begin{proof} Since these conditions are homotopy invariant and hold 
if $M$ is the total space of such a bundle, they are necessary. 
Suppose conversely that they hold. 
As $\nu$ is torsion-free $N$ is the connected sum of a 3-manifold with
free fundamental group and some aspherical $PD_3 $-complexes 
\cite{[Tu90]}.
As $\nu$ is not free there is at least one aspherical summand.
Hence $c.d.\nu=3$ and $H_3(c_N;\mathbb{Z}^{w_1(N)})$ is a monomorphism.
         
Let $p(\gamma):K(\pi,1)\to K(\nu,1)$ be the $S^1 $-bundle corresponding to 
$\gamma$ and let $E=N\times_{K(\nu,1)}K(\pi,1)$ be the total space of the 
$S^1$-bundle over $N$ induced by the classifying map $c_N :N\to K(\nu,1)$.
The bundle map covering $c_N$ is the classifying map $c_E $.
Then $\pi_1(E)\cong\pi=\pi_1(M)$, $w_1(E)=(w_1(N)+w)\gamma=w_1(M)$, 
and $\chi(E)=0=\chi (M)$, 
by conditions (1) and (3). 
The maps $c_N $ and $c_E $ induce a homomorphism between the Gysin sequences 
of the $S^1 $-bundles. 
Since $N$ and $\nu$ have cohomological dimension 3
the Gysin homomorphisms in degree 3 are isomorphisms.
Hence $H_4(c_E;\mathbb{Z}^{w_1(E)})$ is a monomorphism, 
and so {\it a fortiori} $H_4(f_E;\mathbb{Z}^{w_1(E)})$ is also a monomorphism.

Since $\chi(M)=0$ and $\beta_1^{(2)}(\pi)=0$, by Theorem 2.3, part (3)
of Theorem 3.4 implies that 
$\pi_2 (M)\cong\overline{H^2(\pi;\mathbb{Z}[\pi])}$.  
It follows from conditions (2) and (3) and
the LHSSS that $\pi_2(M)\cong\pi_2(E)\cong\gamma^*\pi_2(N)$ as
$\mathbb{Z}[\pi]$-modules.  Conditions (4) and (5) then give us a map
$(\hat c,c_M )$ from $M$ to $P_2(E)=P_2(N)\times_{K(\nu,1)}K(\pi,1)$
such that $(\hat c,c_M )_* [M]=\pm f_{E*} [E]$.  Hence $M$ is homotopy
equivalent to $E$, by Theorem 3.8.

The final assertion now follows from part (1) of Theorem 3.4.
\end{proof}

As $\pi_2(N)$ is a projective $\mathbb{Z}[\nu]$-module, by Theorem 2.18, 
it is homologically trivial and so 
$H_q (\pi;\gamma^* \pi_2 (N)\otimes\mathbb{Z}^{w_1 (M)})=0$ if $q\geq 2$.
Hence it follows from the spectral sequence for $c_{P_2 (M)}$ that 
$H_4 (P_2 (M);\mathbb{Z}^{w_1 (M)})$ maps onto 
$H_4 (\pi;\mathbb{Z}^{w_1 (M)})$,
with kernel isomorphic to 
$H_0 (\pi;\Gamma(\pi_2 (M)))\otimes\mathbb{Z}^{w_1 (M)})$, 
where $\Gamma(\pi_2 (M))=H_4 (K(\pi_2 (M),2);\mathbb{Z})$ 
is Whitehead's universal quadratic construction on $\pi_2 (M)$.
(See \cite[Chapter I]{[Ba']}.)
This suggests that there may be another formulation of the theorem
in terms of conditions (1-3), together with some information on 
$k_1(M)$ and the intersection pairing on $\pi_2 (M)$.
If $N$ is aspherical conditions (4) and (5) are vacuous or
redundant.

Condition (4) is vacuous if $\nu$ is a free group,
for then $c.d.\pi\leq2$.
In this case the Hurewicz homomorphism from $\pi_3(N)$ to 
$H_3(N;\mathbb{Z}^{w_1(N)})$ is 0, 
and so $H_3(f_N;\mathbb{Z}^{w_1(N)})$ is a monomorphism.
The argument of the theorem would then extend if the Gysin map in degree 3 
for the bundle $P_2(E)\to P_2(N)$ were a monomorphism.
If $\nu=1$ then $M$ is orientable, $\pi\cong\mathbb{Z}$ and $\chi(M)=0$, 
so $M\simeq S^3\times S^1$.
In general, if the restriction on $\nu$ is removed it is not clear that there 
should be a degree 1 map from $M$ to such a bundle space $E$.

It would be of interest to have a theorem with hypotheses involving only $M$, 
without reference to a model $N$. 
There is such a result in the aspherical case.

\begin{theorem} 
A finite $PD_4 $-complex $M$ is homotopy equivalent to the total space 
of an $S^1 $-bundle over an aspherical $PD_3 $-complex if and only if 
$\chi (M)=0$ and $\pi=\pi_1(M)$ has an infinite cyclic normal subgroup $A$ 
such that $\pi/A$ has one end and finite cohomological dimension.
\end{theorem}

\begin{proof} The conditions are clearly necessary. 
Conversely, suppose that they hold. 
Since $\pi/A$ has one end $H^s (\pi/A;\mathbb{Z}[\pi/A])=0$ for $s\leq 1$,
and so an LHSSS calculation gives $H^t (\pi;\mathbb{Z}[\pi])=0$ for $t\leq 2$.
Moreover $\beta_1^{(2)}(\pi)=0$, by Theorem 2.3.
Hence $M$ is aspherical and $\pi$ is a $PD_4 $-group, 
by Corollary 3.5.2.
Since $A$ is $FP_\infty$ and $c.d.\pi/A<\infty$ the quotient $\pi/A$ is a 
$PD_3$-group \cite[Theorem 9.11]{[Bi]}. 
Therefore $M$ is homotopy equivalent to the total space of an 
$S^1 $-bundle over the $PD_3 $-complex $K(\pi/A,1)$. 
\end{proof}

Note that a finitely generated torsion-free group has one end if and only if it 
is indecomposable as a free product and is neither infinite cyclic nor trivial. 

In general, if $M$ is homotopy equivalent to the total space of an 
$S^1 $-bundle over some 3-manifold then $\chi (M)=0$ and $\pi_1 (M)$ has 
an infinite cyclic normal subgroup $A$ such that $\pi_1 (M)/A$ 
is virtually of finite cohomological dimension.
Do these conditions characterize such homotopy types?

%% file: m5-5.tex
\chapter{Surface bundles}

In this chapter we shall show that a closed 4-manifold $M$ is homotopy 
equivalent to the total space of a fibre bundle with base and fibre closed 
surfaces if and only if the obviously necessary conditions on the 
Euler characteristic and fundamental group hold.
When the base is $S^2$ we need also conditions on the 
characteristic classes of $M$, and when the base is $RP^2$ 
our results are incomplete.
We shall defer consideration of bundles over $RP^2$ with fibre 
$T$ or $Kb$ and $\partial\not=0$ to Chapter 11, 
and those with fibre $S^2$ or $RP^2$ to Chapter 12.

\section{Some general results}

If $B$, $E$ and $F$ are connected finite complexes and $p:E\to B$ is a Hurewicz 
fibration with fibre homotopy equivalent to $F$ then $\chi(E)=\chi(B)\chi(F)$ 
and the long exact sequence of homotopy gives an exact sequence 
\begin{equation*}
\pi_2 (B)\to\pi_1(F)\to\pi_1(E)\to\pi_1(B)\to1
\end{equation*}
in which the image of $\pi_2 (B)$ under the {\it connecting homomorphism} 
$\partial$ is in the centre of $\pi_1 (F)$. 
(See \cite[page 51]{[Go68]}.) 
These conditions are clearly homotopy invariant.                                

Hurewicz fibrations with base $B$ and fibre $X$ are classified by homotopy 
classes of maps from $B$ to the Milgram classifying space $BE(X)$, 
where $E(X)$ is the monoid of all self homotopy equivalences of $X$, 
with the compact-open topology \cite{[Mi67]}.
If $X$ has been given a base point the evaluation map from $E(X)$ to $X$ 
is a Hurewicz fibration with fibre the subspace (and submonoid) $E_0(X)$
of base point preserving self homotopy equivalences \cite{[Go68]}.

Let $T$ and $Kb$ denote the torus and Klein bottle, respectively.

\begin{lemma} 
Let $F$ be an aspherical closed surface and $B$ a closed smooth manifold.
There are natural bijections from the set of isomorphism classes of
smooth $F$-bundles over $B$ to the set of fibre homotopy equivalence classes of
Hurewicz fibrations with fibre $F$ over $B$ and to the set
$\coprod_{[\xi]}H^2(B;\zeta\pi_1(F)^\xi)$, where the union is over conjugacy
classes of homomorphisms $\xi:\pi_1(B)\to Out(\pi_1(F))$
and $\zeta\pi_1(F)^\xi$ is the $\mathbb{Z}[\pi_1(F)]$-module
determined by $\xi$.
\end{lemma}

\begin{proof} 
If $\zeta\pi_1(F)=1$ the identity components of $Diff(F)$ and $E(F)$ 
are contractible \cite{[EE69]}.
Now every automorphism of $\pi_1(F)$ is realizable by a diffeomorphism 
and homotopy implies isotopy for self diffeomorphisms of surfaces.
(See \cite[Chapter V]{[ZVC]}.)
Therefore $\pi_0(Diff(F))\cong\pi_0(E(F))\cong Out(\pi_1(F))$,
and the inclusion of $Diff(F)$ into $E(F)$ is a homotopy equivalence.
Hence $BDiff(F)\simeq BE(F)\simeq K(Out(\pi_1 (F),1)$,
so smooth $F$-bundles over $B$ and Hurewicz fibrations with fibre $F$ over $B$
are classified by the (unbased) homotopy set
\[
[B,K(Out(\pi_1 (F),1))]=Hom(\pi_1 (B),Out(\pi_1(F)))/\backsim,
\]
where $\xi\backsim\xi'$ if there is an  $\alpha\in Out(\pi_1(F))$
such that $\xi'(b)=\alpha\xi(b)\alpha^{-1}$ for all $b\in\pi_1(B)$.

If $\zeta\pi_1(F)\not=1$ then $F=T$ or $Kb$.
Left multiplication by $T$ on itself induces homotopy equivalences 
from $T$ to the identity components of $Diff(T)$ and $E(T)$.
(Similarly, the standard action of $S^1$ on $Kb$ induces homotopy
equivalences from $S^1$ to the identity components of $Diff(Kb)$ and $E(Kb)$.
See \cite[Theorem III.2]{[Go65]}.)
Let $\alpha:GL(2,\mathbb{Z})\to{Aut(T)}\leq{Diff(T)}$ 
be the standard linear action.
Then the natural maps from the semidirect product 
$T\rtimes_\alpha GL(2,\mathbb{Z})$ to $Diff(T)$ and to $E(T)$ are 
homotopy equivalences.                                
Therefore $BDiff(T)$ is a $K(\mathbb{Z}^2,2)$-fibration 
over $K(GL(2,\mathbb{Z}),1)$.
It follows that $T$-bundles over $B$ are classified by two invariants:
a conjugacy class of homomorphisms $\xi:\pi_1 (B)\to GL(2,\mathbb{Z})$ 
together with a cohomology class in $H^2 (B;(\mathbb{Z}^2)^\xi)$.
A similar argument applies if $F=Kb$.
\end{proof}

\begin{theorem} 
Let $M$ be a $PD_4$-complex and $B$ and $F$ aspherical closed surfaces.
Then $M$ is homotopy equivalent to the total space of an $F$-bundle over $B$ 
if and only if $\chi(M)=\chi(B)\chi(F)$ and $\pi=\pi_1 (M)$ is an extension of 
$\pi_1 (B)$ by $\pi_1 (F)$.
Moreover every extension of $\pi_1 (B)$ by $\pi_1 (F)$ is realized by some 
surface bundle, which is determined up to isomorphism by the extension.
\end{theorem}

\begin{proof} The conditions are clearly necessary. 
Suppose that they hold.
If $\zeta\pi_1 (F)=1$ each homomorphism $\xi:\pi_1(B)\to Out(\pi_1(F))$
corresponds to an unique equivalence class of extensions of $\pi_1(B)$ by 
$\pi_1(F)$ \cite[Proposition 11.4.21]{[Ro]}.
Hence there is an $F$-bundle $p:E\to B$ with $\pi_1(E)\cong\pi$ 
realizing the extension, and $p$ is unique up to bundle isomorphism.
If $F=T$ then every homomorphism $\xi:\pi_1(B)\to GL(2,\mathbb{Z})$
is realizable by an extension (for instance, 
the semidirect product $\mathbb{Z}^2\rtimes_\xi\pi_1(B)$) 
and the extensions realizing $\xi$ are classified 
up to equivalence by $H^2(\pi_1(B);(\mathbb{Z}^2)^\xi)$.
As $B$ is aspherical the natural map from bundles to group extensions 
is a bijection. 
Similar arguments apply if $F=Kb$.
In all cases the bundle space $E$ is aspherical, 
and so $\pi$ is an $FF$ $PD_4$-group. 
Such extensions satisfy the Weak Bass Conjecture \cite[Theorem 5.7]{[Co95]}.
Hence $M\simeq{E}$, by Corollary 3.5.1.
\end{proof}

Such extensions (with $\chi(F)<0$) were shown to be realizable by bundles 
in \cite{[Jo79]}.

\section{Bundles with base and fibre aspherical surfaces}

In many cases the group $\pi_1(M)$ determines the bundle up to diffeomorphism 
of its base.
Lemma 5.3 and Theorems 5.4 and 5.5 are based on \cite{[Jo94]}.
                                   
\begin{lemma} 
Let $G_1$ and $G_2$ be groups with no nontrivial 
abelian normal subgroup.
If $H$ is a normal subgroup of $G=G_1\times G_2$ which contains no nontrivial 
direct product then either $H\leq G_1\times\{1\}$ or $H\leq \{1\}\times G_2$.
\end{lemma}

\begin{proof}                       
Let $P_i$ be the projection of $H$ onto $G_i$, for $i=1,2$.
If $(h,h')\in H$, $g_1\in G_1$ and $g_2\in G_2$ then 
$([h,g_1],1)=[(h,h'),(g_1,1)]$ and $(1,[h',g_2])$ are in $H$.
Hence $[P_1,P_1]\times[P_2,P_2]\leq H$.
Therefore either $P_1$ or $P_2$ is abelian, and so is trivial,
since $P_i$ is normal in $G_i$, for $i=1,2$.
\end{proof}

\begin{theorem} 
Let $\pi$ be a group with a normal subgroup $K$ 
such that $K$ and $\pi/K$ are $PD_2$-groups with trivial centres.
\begin{enumerate}
\item If $C_\pi(K)=1$ and $K_1$ is a non-trivial finitely generated 
normal subgroup of $\pi$ then $C_\pi(K_1)=1$ also.

\item The index $[\pi:KC_\pi(K)]$ is finite if and only if $\pi$ is 
virtually a direct product of $PD_2$-groups.
\end{enumerate}
\end{theorem}                     

\begin{proof} 
(1)\qua Let $z\in C_\pi(K_1)$.
If $K_1\leq K$ then $[K:K_1]<\infty$ and $\zeta K_1=1$.
Let $M=[K:K_1]!$.
Then $f(k)=k^{-1}z^Mkz^{-M}$ is in $K_1$ for all $k$ in $K$.
Now $f(kk_1)=k_1^{-1}f(k)k_1$ and also $f(kk_1)=f(kk_1k^{-1}k)=f(k)$
(since $K_1$ is a normal subgroup centralized by $z$), for all $k$ in $K$ and
$k_1$ in $K_1$.
Hence $f(k)$ is central in $K_1$, and so $f(k)=1$ for all $k$ in $K$.
Thus $z^M$ centralizes $K$.
Since $\pi$ is torsion-free we must have $z=1$.
Otherwise the image of $K_1$ under the projection $p:\pi\to\pi/K$
is a nontrivial finitely generated normal subgroup of $\pi/K$, 
and so has trivial centralizer.
Hence $p(z)=1$.
Now $[K,K_1]\leq K\cap K_1$ and so $K\cap K_1\not=1$, for otherwise 
$K_1\leq C_\pi(K)$.
Since $z$ centralizes the nontrivial normal subgroup $K\cap K_1$
in $K$ we must again have $z=1$.

(2)\qua Since $K$ has trivial centre $KC_\pi(K)\cong K\times C_\pi(K)$ 
and so the condition is necessary.
Suppose that $f:G_1\times G_2\to \pi$ is an isomorphism onto a subgroup of 
finite index, where $G_1$ and $G_2$ are $PD_2$-groups. 
Let $H=K\cap f(G_1\times G_2)$.
Then $[K:H]<\infty$ and so $H$ is also a $PD_2$-group, 
and is normal in $f(G_1\times G_2)$.
We may assume that $H\leq f(G_1)$, by Lemma 5.3.                            
Then $f(G_1)/H$ is finite and is isomorphic to a subgroup of 
$f(G_1\times G_2)/K\leq \pi/K$, so $H=f(G_1)$.                       
Now $f(G_2)$ normalizes $K$ and centralizes $H$, and $[K:H]<\infty$.
Hence $f(G_2)$ has a subgroup of finite index which centralizes $K$, 
as in part (1).
Hence $[\pi:KC_\pi(K)]<\infty$.
\end{proof}

It follows immediately that if $\pi$ and $K$ are as in the theorem whether 
\begin{enumerate}
\item $C_\pi(K)\not=1$ and $[\pi:KC_\pi(K)]=\infty$; 

\item $[\pi:KC_\pi(K)]<\infty$; or 

\item $C_\pi(K)=1$
\end{enumerate} 
depends only on $\pi$ and not on the subgroup $K$.
In \cite{[Jo94]} these cases are labeled as types I, II and III, respectively.
(In terms of the action $\xi:\pi/K\to{Out}(K)$: 
if $\mathrm{Im}(\xi)$ is infinite and $\mathrm{Ker}(\xi)\not=1$ then 
$\pi$ is of type I, if $\mathrm{Im}(\xi)$ is finite then $\pi$ is of type II, 
and if $\xi$ is injective then $\pi$ is of type III.)

\begin{theorem} 
Let $\pi$ be a group with a normal subgroup $K$ such that 
$K$ and $\pi/K$ are virtually $PD_2$-groups 
with no non-trivial finite normal subgroup.
If $\sqrt\pi=1$ and $C_\pi(K)\not=1$
then $\pi$ has at most one other nontrivial finitely generated normal subgroup 
$K_1\not=K$ which contains no nontrivial direct product and is such that 
$\pi/K_1$ has no non-trivial finite normal subgroup.
In that case $K_1\cap K=1$ and  $[\pi:KC_\pi(K)]<\infty$.
\end{theorem}

\begin{proof}                           
Let $p:\pi\to\pi/K$ be the quotient epimorphism.
Then $p(C_\pi(K))$ is a nontrivial normal subgroup of $\pi/K$, 
since $K\cap C_\pi(K)=\zeta K=1$.
Suppose that $K_1<\pi$ is a nontrivial finitely generated normal subgroup 
which contains no nontrivial direct product and is such that $\pi/K_1$ 
has no non-trivial finite normal subgroup.
Let $\Sigma=K_1\cap(KC_\pi(K))$.
Since $\Sigma$ is normal in $KC_\pi(K)\cong{K}\times{C_\pi(K)}$ 
and $\Sigma\leq K_1$, 
either $\Sigma\leq{K}$ or $\Sigma\leq{C_\pi(K)}$, by Lemma 5.3.

If $\Sigma\leq K$ then $p(K_1)\cap p(C_\pi(K))=1$, and so $p(K_1)$ centralizes
the nontrivial normal subgroup $p(C_\pi(K))$ in $\pi/K$.
Therefore $K_1\leq K$ and so $[K:K_1]<\infty$.
Since $\pi/K_1$ has no non-trivial finite normal subgroup we find $K_1=K$.

If $\Sigma\leq{C_\pi(K)}$ then $K_1\cap{K}=1$.
Hence $[K,K_1]=1$, since each subgroup is normal in $\pi$,
and so $K_1\leq C_\pi(K)$.
Moreover $[\pi/K:p(K_1)]<\infty$ since $p(K_1)$ is a nontrivial finitely 
generated normal subgroup of $\pi/K$, and so
$K_1$ and $C_\pi(K)$ are $PD_2$-groups and 
$[\pi:KC_\pi(K)]=[\pi/K:p(C_\pi(K))]\leq [\pi/K:p(K_1)]<\infty$.

If $K_1\not=K$ and $K_2$ is another such subgroup of $\pi$  
then $K_2$ also has finite index in $C_\pi(K)$, by the same argument.
Since $\pi/K_1$ and $\pi/K_2$ have
no non-trivial finite normal subgroup it follows that $K_1=K_2$.
\end{proof}

\begin{cor}
{\rm[Jo93]}\qua
Let $\alpha$ and $\beta$ be automorphisms of $\pi$, and suppose that
${\alpha(K)\cap K=1}$.
Then $\beta(K)=K$ or $\alpha(K)$.
In particular, $Aut(K\times K)\cong Aut(K)^2\rtimes(Z/2Z)$.
\qed
\end{cor}

Groups of type I have an unique such normal subgroup $K$, 
while groups of type II have at most two such subgroups, by Theorem 5.5.                                                        
We shall obtain a somewhat weaker result for groups of type III as a
corollary of Theorem 5.6.

We shall use the following corollary in Chapter 9.

\begin{cor} 
Let $\pi$ be a $PD_4$-group such that $\sqrt\pi=1$. 
Then the following conditions are equivalent:
\begin{enumerate}
\item $\pi$ has a subgroup $\rho\cong\alpha\times\beta$ 
where $\alpha$ and $\beta$ are $PD_2$-groups;

\item $\pi$ has a normal subgroup $\sigma\cong{K}\times{L}$ of finite index 
where $K$ and $L$ are $PD_2$-groups and  $[\pi:N_\pi(K)]\leq2$;

\item $\pi$ has a subgroup $\tau$ such that $[\pi:\tau]\leq2$ and
$\tau\leq{G}\times{H}$ where $G$ and $H$ are virtually $PD_2$-groups.
\end{enumerate}
\end{cor}

\begin{proof} 
Suppose that (1) holds.
Then $[\pi:\rho]<\infty$, by Strebel's Theorem.
Let $N$ be the intersection of the conjugates of $\rho$ in $\pi$.
Then $N$ is normal in $\pi$ and $[\pi:N]<\infty$.
We shall identify $\alpha\cong\alpha\times\{1\}$ and
$\beta\cong\{1\}\times\beta$ with subgroups of $\pi$.
Let $K=\alpha\cap N $ and $L=\beta\cap N$.
Then $K$ and $L$ are $PD_2$-groups, $K\cap{L}=1$
and $\sigma=K.L\cong K\times L$ is normal in $N$ and has finite index in $\pi$.
Moreover $N/K$ and $N/L$ are isomorphic to subgroups of finite index in 
$\beta$ and $\alpha$, respectively, and so are also $PD_2$-groups.
If $\sqrt\pi=1$ all these groups have trivial centre, and so
any automorphism of $N$ must either fix $K$ and $L$
or interchange them, by Theorem 5.5.
Hence $\sigma$ is normal in $\pi$ and $[\pi:N_\pi(K)]\leq 2$.

If (2) holds then $N_\pi(K)=N_\pi(L)$.
Let $\tau=N_\pi(K)$ and let $p_G:\tau\to{G}=\tau/C_\pi(K)$ and 
$p_H:\tau\to{H}=\tau/C_\pi(L)$ be the natural epimorphisms.
Then $p_G|_K$, $p_H|_L$ and $(p_G,p_H)$ are injective and have images of 
finite index in $G$, $H$ and ${G}\times{H}$ respectively.
In particular, $G$ and $H$ are virtually $PD_2$-groups.

If (3) holds let $\alpha=\tau\cap(G\times\{1\})$ and
$\beta=\tau\cap(\{1\}\times{H})$.
Then $\alpha$ and $\beta$  have finite index in $G$ and $H$, respectively, 
and are torsion-free.
Hence they are $PD_2$-groups and clearly $\alpha\cap\beta=1$.
Therefore $\rho=\alpha.\beta\cong\alpha\times\beta$.
\end{proof}

It can be shown that these three conditions remain equivalent under the weaker 
hypothesis that $\pi$ be a $PD_4$-group which is not virtually abelian
(using Lemma 9.4 for the implication $(1)\Rightarrow(3)$).

\begin{theorem} 
Let $\pi$ be a group with normal subgroups $K$ 
and $K_1$ such that $K$, $K_1$ and $\pi/K$ are $PD_2$-groups, 
$\pi/K_1$ is torsion-free and $\chi(\pi/K)<0$. 
Then either $K_1=K$ or $K_1\cap K=1$ and $\pi\cong K\times K_1$ or 
$\chi(K_1)<\chi(\pi/K)$.
\end{theorem}

\begin{proof} 
Let $p:\pi\to\pi/K$ be the quotient epimorphism.
If $K_1\leq K$ then $K_1=K$, as in Theorem 5.5.
Otherwise $p(K_1)$ has finite index in $\pi/K$ and so
$p(K_1)$ is also a $PD_2$-group.
As the minimum number of generators of a $PD_2$-group $G$ is
$\beta_1(G;\mathbb{F}_2)$, we  have $\chi(K_1)\leq\chi(p(K_1))\leq\chi(\pi/K)$. 
We may assume that $\chi(K_1)\geq\chi(\pi/K)$.
Hence $\chi(K_1)=\chi(\pi/K)$ and so $p|_{K_1}$ is an epimorphism.
Therefore $K_1$ and $\pi/K$ have the same orientation type, 
by the nondegeneracy of Poincar\'e duality with coefficients $\mathbb{F}_2$
and the Wu relation $w_1\cup x=x^2$ for all $x\in H^1(G;\mathbb{F}_2)$ and
$PD_2$-groups $G$. 
Hence $K_1\cong\pi/K$.
Since $PD_2$-groups are hopfian $p|_{K_1}$ is an isomorphism.
Hence $[K,K_1]\leq K\cap K_1=1$ and so $\pi=K. K_1\cong K\times\pi/K$.
\end{proof}

\begin{cor}
{\rm[Jo99]}\qua
There are only finitely many such subgroups $K<\pi$.
\end{cor}

\begin{proof} 
We may assume that $\zeta{K}=1$ and $\pi$ is of type III.
There is an epimorphism $\rho:\pi\to{Z}/\chi(\pi)Z$ such that $\rho(K)=0$.
Then $\chi(\mathrm{Ker}(\rho))=\chi(\pi)^2$. 
Since $\pi$ is not virtually a product $K$ is the only normal
$PD_2$-subgroup of $\mathrm{Ker}(\rho)$ with quotient a $PD_2$-group and
such that $\chi(K)^2\leq\chi(\mathrm{Ker}(\rho))$.
The corollary follows since there are only finitely many such epimorphisms $\rho$.
\end{proof}

For each $n\geq1$ there are such groups $\pi$ 
with $\chi(\pi)=24n-8$ which have at least $2^n$ distinct subgroups $K$ 
such that $K$ and $\pi/K$ are orientable \cite{[Sa15]}. 
On the other hand, Cor. 5.6.1 leads to an upper bound of 
$\chi^{\frac\chi2}$, where $\chi=\chi(\pi)$.
Moreover, if $\chi(\pi)\geq16$ then at most $2^{\chi(\pi)}$ 
of these have $|\chi(\pi/K)|>\log_2\chi(\pi)$.

The next corollary follows by elementary arithmetic.

\begin{cor}
If $K_1\not=K$ and $\chi(K_1)=-1$ then $\pi\cong K\times K_1$.
\qed
\end{cor}

\begin{cor}
Let $M$ and $M'$ be the total spaces of 
bundles $\xi$ and $\xi'$ with the same base $B$ and fibre $F$,
where $B$ and $F$ are aspherical closed surfaces such that $\chi(B)<\chi(F)$.
Then $M'$ is diffeomorphic to $M$ via a 
fibre-preserving diffeomorphism if and only if $\pi_1(M')\cong\pi_1(M)$.
\qed
\end{cor}

Compare the statement of Melvin's Theorem on total spaces of $S^2$-bundles 
(Theorem 5.13 below.)

We can often recognise total spaces of aspherical surface bundles
under weaker hypotheses on the fundamental group.

\begin{theorem} 
Let $M$ be a $PD_4$-complex with fundamental group $\pi$. 
Then the following conditions are equivalent:
\begin{enumerate}
\item $M$ is homotopy equivalent to the total space of a bundle with base and 
fibre aspherical closed surfaces:

\item $\pi$ has an $FP_2$ normal subgroup $K$ such that $\pi/K$ is a 
$PD_2$-group and $\pi_2(M)=0$;

\item $\pi$ has a normal subgroup $N$ which is a $PD_2$-group, 
$\pi/N$ is torsion-free and $\pi_2(M)=0$. 
\end{enumerate}
\end{theorem}

\begin{proof} Clearly (1) implies (2) and (3).
Conversely they each imply that $\pi$ has one end and so $M$ is aspherical.
If $K$ is an $FP_3$ normal subgroup in $\pi$ and $\pi/K$ is a $PD_2$-group
then $K$ is a $PD_2$-group, by Theorem 1.19.
If $N$ is a normal subgroup which is a $PD_2$-group then 
$\pi/N$ is virtually a $PD_2$-group, by Theorem 3.10.
Since it is torsion-free it is a $PD_2$-group and so
the theorem follows from Theorem 5.2.
\end{proof}

If $\zeta N=1$ then $\pi/N$ is an extension of $C_\pi(N)$ 
by a subgroup of $Out(N)$.
Thus we may argue instead that $v.c.d.\pi/N<\infty$ and $\pi/N$ is 
$FP_\infty$, so $\pi/N$ is virtually a $PD_2$-group \cite[Theorem 9.11]{[Bi]}.

\begin{cor}
The $PD_4$-complex $M$ is homotopy equivalent to the total space 
of a $T$- or $Kb$-bundle over an aspherical closed surface, 
if and only if $\chi(M)=0$ and $\pi$ has a normal subgroup 
$A\cong\mathbb{Z}^2$ or $\mathbb{Z}\rtimes_{-1}\mathbb{Z}$ 
such that $\pi/A$ is torsion free.
\end{cor}

\begin{proof}
The conditions are clearly necessary.
If they hold then $M$ is aspherical, by Theorem 2.2 and Corollary 3.5.2,
and so this corollary follows from part (3) of Theorem 5.7
\end{proof}

Kapovich has given examples of aspherical closed 4-manifolds $M$ such that
$\pi_1(M)$ is an extension of a $PD_2$-group by a finitely generated normal 
subgroup which is not $FP_2$ \cite{[Ka13]}.

\begin{theorem} 
Let $M$ be a $PD_4$-complex whose fundamental group $\pi$ 
has an ascendant $FP_2$ subgroup $G$ of infinite index with one end
and such that $\chi(M)=0$.
Then $M$ is aspherical.
If moreover $c.d.G=2$ and $\chi(G)\not=0$ 
then $G$ is a $PD_2$-group and either $[\pi:N_\pi(G)]<\infty$ or
there is a subnormal chain $G<J<K\leq\pi$ such that 
$[\pi:K]<\infty$ and $K/J\cong{J/G}\cong\mathbb{Z}$.
\end{theorem}

\begin{proof} 
The argument of the first paragraph of the proof of Theorem 4.8 
applies equally well here to show that $M$ is aspherical.

Assume henceforth that $c.d.G=2$ and $\chi(G)<0$.
If $G<\tilde G<G_\gamma$ and $c.d.\tilde{G}=2$ then $[\tilde G:G]<\infty$, 
by Lemma 2.15.
Hence $\tilde{G}$ is $FP$ and $[\tilde G:G]\leq|\chi(G)|$,
since $\chi(G)=[\tilde{G}:G]\chi(\tilde{G})$.
We may assume that $\tilde{G}$ is maximal among all groups of cohomological
dimension 2 in an ascendant chain from $G$ to $\pi$.
Let $G=G_0<G_1<\dots<G_\beth=\pi$ be such an ascendant chain,
with $\tilde{G}=G_n$ for some finite ordinal $n$.
Then $[G_{n+1}:G]=\infty$ and $c.d.G_{n+1}\geq3$.

If $\tilde{G}$ is normal in $\pi$ then $\tilde{G}$ is a $PD_2$-group and
$\pi/\tilde{G}$ is virtually a $PD_2$-group, by Theorem 3.10.
Moreover $[\pi:N_\pi(G)]<\infty$,
since $\tilde{G}$ has only finitely many subgroups of index $[\tilde{G}:G]$.
Therefore $\pi$ has a normal subgroup $K\leq N_\pi(G)$ such that
$[\pi:K]<\infty$ and $K/G$ is a $PD_2^+$-group.

Otherwise, replacing $G_{n+1}$ by the union of the terms $G_\alpha$ 
which normalize $\tilde{G}$ and reindexing, if necessary, we may assume 
that $\tilde{G}$ is not normal in $G_{n+2}$.
Let $h$ be an element of $G_{n+2}$ such that 
$h\tilde{G}h^{-1}\not=\tilde{G}$, 
and let $H=\tilde{G}.h\tilde{G}h^{-1}$.
Then $\tilde{G}$ is normal in $H$ and $H$ is normal in $G_{n+1}$, 
so $[H:\tilde{G}]=\infty$ and $c.d.H=3$.
Moreover $H$ is $FP$ \cite[Proposition 8.3]{[Bi]},
and $H^s(H;\mathbb{Z}[H])=0$ for $s\leq2$, by an LHSSS argument.

If $c.d.G_{n+1}=3$ then $G_{n+1}/H$ is locally finite \cite[Theorem 8.2]{[Bi]}.
Hence it is finite, by the Gildenhuys-Strebel Theorem.
Therefore $G_{n+1}$ is $FP$ and 
$H^s(G_{n+1};\mathbb{Z}[G_{n+1}])=0$ for $s\leq2$.
Since $G_{n+1}$ is also ascendant in $\pi$ it is a $PD_3$-group,
$[\pi:N_\pi(G_{n+1})]<\infty$ and $N_\pi(G_{n+1})/G_{n+1}$ has two ends, 
by Theorem 4.8.
Hence $G_{n+1}/\tilde{G}$ has two ends also,  
and $\tilde{G}$ is a $PD_2$-group, by Theorem 2.12.
We may easily find subgroups $J\leq{G_{n+1}}$ and $K\leq N_\pi(G_{n+1})$ 
such that $G<J<K$, $J/G\cong{K/J}\cong\mathbb{Z}$ and $[\pi:K]<\infty$.

If $c.d.G_{n+1}=4$ then $[\pi:G_{n+1}]$ is again finite 
and $G_{n+1}$ is a $PD_4$-group.
Hence the result follows as for the case when $\tilde{G}$ is normal in $\pi$.
\end{proof} 

\begin{cor}
If $\chi(M)=0$,
$G$ is a $PD_2$-group, $\chi(G)\not=0$ and $G$ is normal in $\pi$ then 
$M$ has a finite covering space which is homotopy equivalent to the total space 
of a surface bundle over $T$.
\end{cor}

\begin{proof} Since $G$ is normal in $\pi$ and $M$ is aspherical
$M$ has a finite covering which is homotopy equivalent to a $K(G,1)$-bundle 
over an aspherical orientable surface, as in Theorem 5.7.
Since $\chi(M)=0$ the base must be $T$.
\end{proof}

If $\pi/G$ is virtually $\mathbb{Z}^2$ then it has a subgroup 
of index at most 6 which maps onto $\mathbb{Z}^2$ or 
$\mathbb{Z}\rtimes_{-1}\!\mathbb{Z}$.

Let $G$ be a $PD_2$-group such that $\zeta G=1$.
Let $\theta$ be an automorphism of $G$ whose class in $Out(G)$ has 
infinite order and let $\lambda:G\to\mathbb{Z}$ be an epimorphism. 
Let $\pi=(G\times\mathbb{Z})\rtimes_\phi\mathbb{Z}$ where
$\phi(g,n)=(\theta(g),\lambda(g)+n)$ for all $g\in G$ and $n\in\mathbb{Z}$.
Then $G$ is subnormal in $\pi$ but this group is not virtually the group of a 
surface bundle over a surface.

If $\pi$ has an ascendant subgroup $G$ which is a $PD_2$-group 
with $\chi(G)=0$ then $\sqrt G\cong\mathbb{Z}^2$ is ascendant in $\pi$ 
and hence contained in $\sqrt\pi$.
In this case $h(\sqrt\pi)\geq2$ and so either Theorem 8.1 or Theorem 9.2
applies, to show that $M$ has a finite covering space which is homotopy
equivalent to the total space of a $T$-bundle over an aspherical closed surface.

\section{Bundles with aspherical base and fibre $S^2$ or $RP^2$}

Let $E^+(S^2)$ denote the connected component of $id_{S^2} $ in $E(S^2)$,
i.e., the submonoid of degree 1 maps.
The connected component of $id_{S^2} $ in $E_0 (S^2) $ may be identified with 
the double loop space $\Omega^2 S^2 $. 
                                             
\begin{lemma} 
Let $X$ be a finite $2$-complex. 
Then there are natural bijections
$[X;BO(3)]\cong [X;BE(S^2)]\cong H^1(X;\mathbb{F}_2)\times H^2(X;\mathbb{F}_2)$.
\end{lemma}

\begin{proof} 
As a self homotopy equivalence of a sphere is homotopic to the identity if and 
only if it has degree $+1$ the inclusion of $O(3)$ into $E(S^2)$ is bijective 
on components.
Evaluation of a self map of $S^2$ at the basepoint determines 
fibrations of $SO(3)$ and $E^+(S^2)$ over $S^2$, 
with fibre $SO(2)$ and $\Omega^2 S^2$, respectively,
and the map of fibres induces an isomorphism on $\pi_1$.
On comparing the exact sequences of homotopy for these fibrations we see that 
the inclusion of $SO(3)$ in $E^+(S^2)$ also induces an isomorphism on $\pi_1 $.
Since the Stiefel-Whitney classes are defined for any spherical fibration
and $w_1$ and $w_2$ are nontrivial on suitable $S^2$-bundles over $S^1$ and 
$S^2$, respectively, the inclusion of $BO(3)$ into $BE(S^2)$ and the map
$(w_1,w_2):BE(S^2)\to K(Z/2Z,1)\times K(Z/2Z,2)$ induces isomorphisms
on $\pi_i$ for $i\leq2$.
The lemma follows easily.
\end{proof}

Thus there is a natural 1-1 correspondance between $S^2$-bundles and
spherical fibrations over such complexes, and any such bundle $\xi$ is 
determined up to isomorphism over $X$ by its
total Stiefel-Whitney class $w(\xi)=1+w_1(\xi)+w_2(\xi)$.
(From another point of view: if $w_1(\xi)=w_1(\xi')$ there is an 
isomorphism of the restrictions of $\xi$ and $\xi'$ 
over the 1-skeleton $X^{[1]}$.
The difference $w_2(\xi)-w_2(\xi')$ is the obstruction to extending
any such isomorphism over the 2-skeleton.)

\begin{theorem} 
Let $M$ be a $PD_4$-complex and 
$B$ an aspherical closed surface.
Then the following conditions are equivalent:
\begin{enumerate}
\item $\pi_1(M)\cong\pi_1(B)$ and $\chi(M)=2\chi(B)$;

\item $\pi_1 (M)\cong\pi_1 (B)$ and $\widetilde M\simeq S^2$;

\item $M$ is homotopy equivalent to the total space of an 
$S^2 $-bundle over $B$.
\end{enumerate}
\end{theorem}

\begin{proof} If (1) holds then 
$H_3 (\widetilde M;\mathbb{Z})=H_4 (\widetilde M;\mathbb{Z})=0$, 
as $\pi_1 (M)$ has one end, and 
$\pi_2(M)\cong\overline{H^2(\pi;\mathbb{Z}[\pi])}$, by Theorem 3.12. 
Since this is infinite cyclic, $\widetilde M$ is homotopy equivalent to $S^2$.
If (2) holds we may assume that there is a Hurewicz fibration $h:M\to B$ 
which induces an isomorphism of fundamental groups. 
As the homotopy fibre of $h$ is $\widetilde M$, Lemma 5.9 implies that 
$h$ is fibre homotopy equivalent to the projection of an $S^2 $-bundle over $B$.
Clearly (3) implies the other conditions. 
\end{proof}

We shall summarize some of the key properties of the Stiefel-Whitney classes
of such bundles in the following lemma.

\begin{lemma} 
Let $\xi$ be an $S^2 $-bundle over a closed surface $B$,
with total space $M$ and projection $p:M\to B$. Then
\begin{enumerate}              
\item $\xi$ is trivial if and only if $w(M)=p^*w(B)$;

\item $\pi_1(M)\cong \pi_1(B)$ acts on $\pi_2(M)$ by
multiplication by $w_1(\xi)$;

\item the intersection form on $H_2(M;\mathbb{F}_2)$ is even if and 
only if $w_2(\xi)=0$; 

\item if $q:B'\to B$ is a $2$-fold covering map with 
connected domain $B'$ then $w_2 (q^*\xi)=0$.
\end{enumerate}
\end{lemma}

\begin{proof} (1)\qua 
Applying the Whitney sum formula and naturality to the tangent bundle of 
the $B^3$-bundle associated to $\xi$ gives $w(M)=p^* w(B)\cup p^*w(\xi)$. 
Since $p$ is a 2-connected map the induced homomorphism $p^*$ is injective in 
degrees $\leq 2$ and so $w(M)=p^*w(B)$ if and only if $w(\xi)=1$.
By Lemma 5.9 this is so if and only if $\xi$ is trivial, 
since $B$ is 2-dimensional.

(2)\qua It is sufficient to consider the restriction of $\xi$ over loops
in $B$, where the result is clear.
                                                  
(3)\qua By Poincar\'e duality, the intersection form is even if and only if 
the Wu class $v_2(M)=w_2(M)+ w_1(M)^2$ is 0. 
Now 
\begin{equation*}
\begin{split}
v_2(M) &=p^*( w_1(B)+ w_1(\xi ))^2+
p^*(w_2(B)+w_1(B)\cup w_1(\xi )+w_2(\xi))\\
&=p^*(w_2(B)+w_1(B)\cup w_1(\xi )+w_2(\xi)+ w_1(B)^2+w_1(\xi )^2)\\
&=p^*(w_2(\xi )),\\
\end{split}
\end{equation*}
since $w_1(B)\cup\eta=\eta^2$ and $w_1(B)^2=w_2(B)$, 
by the Wu relations for $B$. 
Hence $v_2(M)=0$ if and only if $w_2(\xi )=0$, 
as $p^*$ is injective in degree 2.

\noindent (4) We have $q_* (w_2 (q^*\xi)\cap [B'])=
q_* ((q^*w_2 (\xi))\cap [B'])=w_2 (\xi)\cap q_* [B']$,
by the projection formula.
Since $q$ has degree 2 this is 0,
and since $q_* $ is an isomorphism in degree 0
we find $w_2 (q^*\xi)\cap [B']=0$.
Therefore $w_2 (q^*\xi)=0$, by Poincar\'e duality for $B'$.
\end{proof}

Melvin has determined criteria for the total spaces of $S^2 $-bundles over a 
compact surface to be diffeomorphic, in terms of their Stiefel-Whitney classes.
We shall give an alternative argument for the cases with aspherical base.

\begin{lemma} 
Let $B$ be a closed surface and $w$ be the 
Poincar\'e dual of $w_1(B)$. 
If $u_1$ and $u_2$ are elements of $H_1(B;\mathbb{F}_2)\setminus\{0,w\}$ 
such that 
$u_1.u_1=u_2.u_2$ then there is a diffeomorphism $f:B\to B$ which is 
a composite of Dehn twists about two-sided essential simple closed curves 
and such that $f_*(u_1)=u_2$.
\end{lemma}

\begin{proof} 
For simplicity of notation, we shall use the same symbol for a simple closed curve
$u$ on $B$ and its homology class in $H_1(B;\mathbb{F}_2)$. 
The curve $u$ is two-sided if and only if $u.u=0$.
In that case we shall let $c_u$ denote the
automorphism of $H_1(B;\mathbb{F}_2)$ induced by a Dehn twist about $u$.
Note also that $u.u=u.w$ and $c_v(u)=u+(u.v)v$ for all $u$ and two-sided $v$ in 
$H_1(B;\mathbb{F}_2)$.                                               

If $B$ is orientable it is well known that the group of isometries 
of the intersection form acts transitively on $H_1(B;\mathbb{F}_2)$, 
and is generated by the automorphisms $c_u$.
Thus the claim is true in this case.

If $w_1(B)^2\not=0$ then $B\cong RP^2\sharp T_g$, where $T_g$ is orientable.
If $u_1.u_1=u_2.u_2=0$ then $u_1$ and $u_2$ are represented by simple closed 
curves in $T_g$, and so are related by a diffeomorphism 
which is the identity on the $RP^2$ summand.
If $u_1.u_1=u_2.u_2=1$ let $v_i=u_i+w$. 
Then $v_i.v_i=0$ and this case follows from the earlier one.

Suppose finally that $w_1(B)\not=0$ but $w_1(B)^2=0$; equivalently, 
that $B\cong Kb\sharp T_g$, where $T_g$ is orientable.
Let $\{w,z\}$ be a basis for the homology of the $Kb$ summand.
In this case $w$ is represented by a 2-sided curve. 
If $u_1.u_1=u_2.u_2=0$ and $u_1.z=u_2.z=0$ then $u_1$ and $u_2$ are
represented by simple closed curves in $T_g$, and so are related by a diffeomorphism 
which is the identity on the $Kb$ summand. 
The claim then follows if $u.z=1$ for $u=u_1$ or $u_2$, since
we then have $c_w(u).c_w(u)=c_w(u).z=0$.
If $u.u\not=0$ and $u.z=0$ then $(u+z).(u+z)=0$ and $c_{u+z}(u)=z$.
If $u.u\not=0$, $u.z\not=0$ and $u\not= z$ then $c_{u+z+w}c_w(u)=z$.
Thus if $u_1.u_1=u_2.u_2=1$ both $u_1$ and $u_2$ are related to $z$.
Thus in all cases the claim is true.
\end{proof}

\begin{theorem} 
[Melvin] 
Let $\xi$ and $\xi'$ be two $S^2 $-bundles 
over an aspherical closed surface $B$. 
Then the following conditions are equivalent:
\begin{enumerate}
\item there is a diffeomorphism $f:B\to B$ such that $\xi=f^*\xi'$; 

\item the total spaces $E(\xi)$ and $E(\xi')$ are diffeomorphic; and

\item $w_1 (\xi)=w_1 (\xi')$ if $w_1 (\xi)=0$ or $w_1 (B)$, 
$w_1 (\xi)\cup w_1 (B)=w_1 (\xi')\cup w_1 (B)$ 
and $w_2 (\xi)=w_2 (\xi')$.
\end{enumerate}
\end{theorem} 
                      
\begin{proof} Clearly (1) implies (2).
A diffeomorphism $h:E\to E'$ induces an isomorphism on fundamental groups;
hence there is a diffeomorphism $f:B\to B$ such that $fp$ is homotopic to $p'h$.
Now $h^*w(E')=w(E)$ and $f^*w(B)=w(B)$.
Hence $p^*f^*w(\xi')=p^*w(\xi)$ and so $w(f^*\xi')=f^*w(\xi')=w(\xi)$.
Thus $f^*\xi'=\xi$, by Lemma 5.9, and so (2) implies (1).

If (1) holds then $f^*w(\xi')=w(\xi)$. 
Since $w_1(B)=v_1(B)$ is the characteristic element for the cup product pairing
from $H^1(B;\mathbb{F}_2)$ to $H^2(B;\mathbb{F}_2)$ and                                  
$H^2(f;\mathbb{F}_2)$ is the identity $f^*w_1(B)=w_1(B)$, 
$w_1 (\xi)\cup w_1 (B)=w_1 (\xi')\cup w_1 (B)$ and $w_2 (\xi)=w_2 (\xi')$.
Hence(1) implies (3).

If $w_1 (\xi)\cup w_1 (B)=w_1 (\xi')\cup w_1 (B)$ and $w_1 (\xi)$ and 
$w_1 (\xi')$ are neither 0 nor $w_1(B)$ then there is a diffeomorphism 
$f:B\to B$ such that $f^*w_1(\xi')=w_1(\xi)$,
by Lemma 5.12 (applied to the Poincar\'e dual homology classes).
Hence (3) implies (1).
\end{proof}

\begin{cor}
There are $4$ diffeomorphism classes of $S^2$-bundle spaces 
if $B$ is orientable and $\chi(B)\leq 0$, 
$6$ if $B=Kb$ and $8$ if $B$ is nonorientable and $\chi(B)<0$.
\qed
\end{cor}
         
See \cite{[Me84]} for a more geometric argument, 
which applies also to $S^2$-bundles over
surfaces with nonempty boundary.
The theorem holds also when $B=S^2$ or $RP^2$;
there are 2 such bundles over $S^2$ and 4 over $RP^2$. (See Chapter 12.)

\begin{theorem} 
Let $M$ be a $PD_4$-complex with fundamental group $\pi$. 
The following are equivalent:
\begin{enumerate}
\item $\pi\not=1$ and $\pi_2 (M)\cong\mathbb{Z}$.

\item $\widetilde{M}\simeq{S^2}$;

\item $M$ has a covering space of degree $\leq 2$ which is homotopy equivalent 
to the total space of an $S^2 $-bundle over an aspherical closed surface;

\item $\pi$ is virtually a $PD_2$-group and $\chi(M)=2\chi(\pi)$.
\end{enumerate}
\noindent If these conditions hold the kernel $K$ of the natural action of 
$\pi$ on $\pi_2(M)$ is a $PD_2$-group.
\end{theorem}

\begin{proof} 
Suppose that (1) holds.
If $\pi$ is finite and $\pi_2(M)\cong\mathbb{Z}$ then 
$\widetilde M\simeq CP^2$, and so admits no nontrivial free 
group actions, by the Lefshetz fixed point theorem.
Hence $\pi$ must be infinite.
Then $H_0(\widetilde M;\mathbb{Z})=\mathbb{Z}$, 
$H_1(\widetilde M;\mathbb{Z})=0$ and 
$H_2(\widetilde M;\mathbb{Z})=\pi_2 (M)$,
while $H_3 (\widetilde M;\mathbb{Z})\cong\overline{H^1 (\pi;\mathbb{Z}[\pi])}$ and $H_4 (\widetilde M;\mathbb{Z})=0$.                                                   
Now $Hom_{\mathbb{Z}[\pi]}(\pi_2(M),\mathbb{Z}[\pi])=0$, 
since $\pi$ is infinite and $\pi_2 (M)\cong\mathbb{Z}$. 
Therefore $H^2 (\pi;\mathbb{Z}[\pi])$ is infinite cyclic, by Lemma 3.3, 
and so $\pi$ is virtually a $PD_2$-group, by Bowditch's Theorem.
Hence $H_3 (\widetilde M;\mathbb{Z})=0$ and so $\widetilde M\simeq S^2 $.
If $C$ is a finite cyclic subgroup of $K$ then 
$H_{n+3} (C;\mathbb{Z})\cong  H_n (C;H_2 (\widetilde M;\mathbb{Z}))$ 
for all $n\geq 2$, by Lemma 2.10.
Therefore $C$ must be trivial, so $K$ is torsion-free.
Hence $K$ is a $PD_2$-group and (3) now follows from Theorem 5.10.
Clearly (3) implies (2) and (2) implies (1).
The equivalence of (3) and (4) follows from Theorem 5.10.
\end{proof}

A straightfoward Mayer-Vietoris argument may be used to show directly that 
if $H^2(\pi;\mathbb{Z}[\pi])\cong\mathbb{Z}$ then $\pi$ has one end.

\begin{lemma} 
Let $X$ be a finite $2$-complex. 
Then there are natural bijections
$[X;BSO(3)]\cong [X;BE(RP^2)]\cong  H^2(X;\mathbb{F}_2)$.
\end{lemma}

\begin{proof}
Let $(1,0,0)$ and $[1:0:0]$ be the base points for $S^2$ and $RP^2$
respectively.
A based self homotopy equivalence $f$ of $RP^2$ lifts to a based self homotopy
equivalence $f^+$ of $S^2$.
If $f$ is based homotopic to the identity then $deg(f^+)=1$.
Conversely, any based self homotopy equivalence is based homotopic to a map
which is the identity on $RP^1$; if moreover $deg(f^+)=1$ then this map is the
identity on the normal bundle and it quickly follows that $f$ 
is based homotopic to the identity.
Thus $E_0(RP^2)$ has two components. 
The diffeomorphism $g$ defined by $g([x:y:z])=[x:y:-z]$ is isotopic to the 
identity (rotate in the $(x,y)$-coordinates).
However $deg(g^+)=-1$.
It follows that $E(RP^2)$ is connected.
As every self homotopy equivalence of $RP^2$ is covered by a degree 1 self map 
of $S^2$,  there is a natural map from $E(RP^2)$ to $E^+(S^2)$.

We may use obstruction theory to show that $\pi_1 (E_0 (RP^2 ))$ has order 2. 
Hence $\pi_1 (E(RP^2))$ has order at most 4.
Suppose that there were a homotopy $f_t $ through self maps of $RP^2 $ with 
$f_0 =f_1 =id_{RP^2} $ and such that the loop $f_t (*)$ is essential, 
where $*$ is a basepoint. 
Let $F$ be the map from $RP^2 \times S^1 $ to $RP^2 $ determined by 
$F(p,t)=f_t (p)$, 
and let $\alpha$ and $\beta$ be the generators of $H^1 (RP^2;\mathbb{F}_2)$ 
and $H^1 (S^1;\mathbb{F}_2)$, respectively. 
Then $F^*\alpha=\alpha\otimes 1+1\otimes\beta$ and so
$(F^*\alpha)^3=\alpha^2\otimes\beta$ which is nonzero, 
contradicting $\alpha^3 =0$.
Thus there can be no such homotopy, and so the homomorphism from 
$\pi_1 (E(RP^2 ))$ to $\pi_1 (RP^2 )$ 
induced by the evaluation map must be trivial. 
It then follows from the exact sequence of homotopy for this evaluation map 
that the order of $\pi_1 (E(RP^2))$ is at most 2. 
The group $SO(3)\cong O(3)/(\pm I)$ acts isometrically on $RP^2$.
As the composite of the maps on $\pi_1 $ induced by the inclusions 
$SO(3)\subset E(RP^2)\subset E^+(S^2)$ is an isomorphism of groups of order 2 
the first map also induces an isomorphism.
It follows as in Lemma 5.9 that there are natural bijections
$[X;BSO(3)]\cong [X;BE(RP^2)]\cong  H^2(X;\mathbb{F}_2)$.
\end{proof}

Thus there is a natural 1-1 correspondance between $RP^2$-bundles and
orientable spherical fibrations over such complexes.
The $RP^2$-bundle corresponding to an orientable $S^2$-bundle is the quotient
by the fibrewise antipodal involution.
In particular, there are two $RP^2$-bundles over each closed aspherical surface.

\begin{theorem} 
Let $M$ be a $PD_4$-complex and 
$B$ an aspherical closed surface.
Then the following conditions are equivalent:
\begin{enumerate}
\item $\pi_1(M)\cong\pi_1(B)\times(Z/2Z)$ and $\chi(M)=\chi(B)$;

\item $\pi_1 (M)\cong\pi_1 (B)\times(Z/2Z)$ and $\widetilde M\simeq S^2$;

\item $M$ is homotopy equivalent to the total space of an 
$RP^2$-bundle over $B$.
\end{enumerate}
\end{theorem}

\begin{proof} 
Suppose that (1) holds, and let $w:\pi_1 (M)\to Z/2Z$ be the projection onto 
the $Z/2Z$ factor.
Then the covering space associated with the kernel of $w$ satisfies the 
hypotheses of Theorem 5.10 and so $\widetilde M\simeq S^2$. 

If (2) holds the homotopy fibre of the map $h$ from $M$ to $B$ inducing the 
projection of $\pi_1 (M)$ onto $\pi_1 (B)$ is homotopy equivalent to $RP^2 $.
The map $h$ is fibre homotopy equivalent 
to the projection of an $RP^2 $-bundle over $B$, by Lemma 5.15. 

If $E$ is the total space of an $RP^2 $-bundle over $B$, 
with projection $p$,
then $\chi(E)=\chi(B)$ and the long exact sequence of homotopy gives a short 
exact sequence $1\to Z/2Z\to\pi_1 (E)\to\pi_1 (B)\to 1$. 
Since the fibre has a product neighbourhood,
$j^* w_1 (E)=w_1 (RP^2)$, where $j:RP^2 \to E$ is the inclusion of the fibre 
over the basepoint of $B$, 
and so $w_1 (E)$ considered as a homomorphism from $\pi_1 (E)$ to $Z/2Z$
splits the injection $j_* $. 
Therefore $\pi_1 (E)\cong\pi_1 (B)\times(Z/2Z)$ and so (1) holds,
as these conditions are clearly invariant under homotopy.
\end{proof}

We may use the above results to refine some of the conclusions of Theorem 3.9
on $PD_4$-complexes with finitely dominated covering spaces.

\begin{theorem} 
Let $M$ be a $PD_4$-complex with fundamental group $\pi$,
and let $p:\pi\to{G}$ be an epimorphism with $FP_2$ kernel $\nu$.
Suppose that $H^2(G;\mathbb{Z}[G])$ is infinite cyclic.
Then the following conditions are equivalent:
\begin{enumerate}

\item $Hom_{\mathbb{Z}[\pi]}(\pi_2(M),\mathbb{Z}[\pi])=0$;

\item $C_*(\widetilde{M})|_\nu$ has finite $2$-skeleton;

\item the associated covering space $M_\nu$ is homotopy equivalent 
to a closed surface;

\item $M$ has a finite covering space which is homotopy equivalent 
to the total space of a surface bundle 
over an aspherical closed surface.
\end{enumerate}
\end{theorem}

\begin{proof} By Bowditch's Theorem $G$ is virtually a $PD_2$-group.
Hence $\pi$ has one end and $H^2(\pi;\mathbb{Z}[\pi])$ is infinite cyclic, 
if $\nu$ is finite, and is 0 otherwise, by an LHSSS argument.

If (1) holds $\pi_2(M)\cong\overline{H^2(\pi;\mathbb{Z}[\pi])}$,
by Lemma 3.3.
If (2) holds 
$\pi_2(M)\cong{H}_2(M_\nu;\mathbb{Z}[\nu])\cong{H}^0(M_\nu;\mathbb{Z}[\nu])$,
by Theorem 1.19$^\prime$.
In either case, if $\nu$ is finite $\pi_2(M)\cong\mathbb{Z}$,
while if $\nu$ is infinite $\pi_2(M)=0$ and $M$ is aspherical.
Condition (3) now follows from Theorems 5.10, 5.16 and 1.19,
and (4) follows easily.

If (4) holds then $\pi$ is infinite and
$\pi_2(M)=\pi_2(M_\nu)\cong\mathbb{Z}$ or is 0, and so (1) holds.
\end{proof}

The total spaces of such bundles with base an aspherical surface have minimal 
Euler characteristic for their fundamental groups (i.e., $\chi(M)=q(\pi)$),
by Theorem 3.12 and the remarks in the paragraph preceding it.

The $FP_2$ hypothesis is in general necessary,
as observed after Theorem 5.7. (See \cite{[Ka98]}.)
However it may be relaxed when $G$ is virtually $\mathbb{Z}^2$ and $\chi(M)=0$.

\begin{theorem} 
Let $M$ be a finite $PD_4$-complex with fundamental group $\pi$.
Then $M$ is homotopy equivalent to the total space of a surface 
bundle over $T$ or $Kb$ if and only if $\pi$ is an extension 
of $\mathbb{Z}^2$ or $\mathbb{Z}\rtimes_{-1}\!\mathbb{Z}$ (respectively) 
by a finitely generated normal subgroup 
$\nu$ and $\chi(M)=0$.
\end{theorem}

\begin{proof} The conditions are clearly necessary.
If they hold the covering space $M_\nu$ associated to the subgroup $\nu$
is homotopy equivalent to a closed surface,
by Corollaries 4.5.2 and 2.12.1.
The result then follows from Theorems 5.2, 5.10 and 5.16.
\end{proof}

In particular, if $\pi$ is the nontrivial extension of $\mathbb{Z}^2$ 
by $Z/2Z$ then $q(\pi)>0$.

\section{Bundles over $S^2 $}

Since $S^2 $ is the union of two discs along a circle, 
an $F$-bundle over $S^2 $ is determined by the homotopy class 
of the clutching function in $\pi_1(Diff(F))$.
(This group is isomorphic to $\zeta\pi_1(F)$ and hence to 
${H}^2(S^2;\zeta\pi_1(F))$.)
On the other hand, if $M$ is a $PD^4$-complex then
cellular approximation gives bijections 
$H^2(M;\mathbb{Z})=[M;CP^\infty]=[M;CP^2]$, and a map
$f:M\to{CP^2}$ factors through $CP^2\setminus{D^4}\sim{S^2}$ if 
and only if $\mathrm{deg}(f)=0$.
Thus if $u\in{H^2(M;\mathbb{Z})}$ and $i_2 $ generates $H^2 (S^2;\mathbb{Z})$
then $u=f^*i_2$ for some $f:M\to{S^2}$ if and only if $u^2=0$. 
The map is uniquely determined by $u$ \cite[Theorem 8.4.11]{[Sp]}.

\begin{theorem} 
Let $M$ be a $PD_4$-complex with
fundamental group $\pi$ and $F$ a closed surface.
Then $M$ is homotopy equivalent to the total space of an $F$-bundle over $S^2$ 
if and only if $\chi(M)=2\chi(F)$ and
\begin{enumerate}
\item (when $\chi(F)<0$) $\pi\cong\pi_1(F)$,
$w_1(M)=c_M^*w_1(F)$ and 
$w_2(M)=w_1(M)^2=(c_M^*w_1(F))^2$; or

\item (when $F=T$) $\pi\cong\mathbb{Z}^2$ and $w_1(M)=w_2(M)=0$, 
or $\pi\cong\mathbb{Z}\oplus(Z/nZ)$ for some $n>0$ and, 
if $n=1$ or $2$, $w_1(M)=0$; or

\item (when $F=Kb$) $\pi\cong\mathbb{Z}\rtimes_{-1}\!\mathbb{Z}$, 
$w_1(M)\not=0$ and $w_2(M)=w_1(M)^2=0$, or $\pi$ has a presentation 
$\langle x,y\mid yxy^{-1} =x^{-1},\medspace y^{2n}=1\rangle$ for some $n>0$, 
where $w_1(M)(x)=0$ and $w_1 (M)(y)=1$; or

\item (when $F=S^2$) $\pi=1$ and the index $\sigma(M)=0$; or

\item (when $F=RP^2$) $\pi=Z/2Z$, $w_1 (M)\not=0$ and there is a 
class $u$ of
infinite order in $H^2(M;\mathbb{Z})$ and such that $u^2=0$.
\end{enumerate}
\end{theorem}

\begin{proof}   
Let $p_E :E\to S^2$ be such a bundle. Then $\chi(E)=2\chi(F)$ and 
$\pi_1 (E)\cong \pi_1 (F)/\partial \pi_2 (S^2)$, where
Im$(\partial)\leq\zeta\pi_1(F)$ \cite{[Go68]}. 
The characteristic classes of $E$ restrict to the characteristic classes
of the fibre, as it has a product neighbourhood.
As the base is 1-connected $E$ is orientable if and only if the fibre is 
orientable.
Thus the conditions on $\chi$, $\pi$ and $w_1$ are all necessary.
We shall treat the other assertions case by case.
 
(1) \qua If $\chi(F)<0$ any $F$-bundle over $S^2 $ is trivial, 
by Lemma 5.1. 
Thus the conditions are necessary. 
Conversely, if they hold then $c_M$ is fibre homotopy equivalent to the 
projection of an $S^2 $-bundle $\xi$ with base $F$, by Theorem 5.10.
The conditions on the Stiefel-Whitney classes then imply that $w(\xi)=1$ 
and hence that the bundle is trivial, by Lemma 5.11.
Therefore $M$ is homotopy equivalent to $S^2 \times F$.
          
(2)\qua If $\partial=0$ there is a map $q:E\to T$ which induces an 
isomorphism of fundamental groups, and the map $(p_E ,q):E\to S^2 \times T$ 
is clearly a homotopy equivalence, so $w(E)=1$. 
Conversely, if $\chi(M)=0$, $\pi\cong\mathbb{Z}^2 $ and $w(M)=1$ 
then $M$ is homotopy equivalent to $S^2 \times T$, 
by Theorem 5.10 and Lemma 5.11.

If $\chi(M)=0$ and $\pi\cong\mathbb{Z}\oplus(Z/nZ)$ for some $n>0$ then 
the covering space $M_{Z/nZ} $ corresponding to the torsion subgroup $Z/nZ$
is homotopy equivalent to a lens space $L$, by Corollary 4.5.2. 
As observed in Chapter 4 the manifold $M$ is homotopy equivalent 
to the mapping torus of a generator of the group 
of covering transformations $Aut(M_{Z/nZ}/M)\cong\mathbb{Z}$.
Since the generator induces the identity on $\pi_1 (L)\cong Z/nZ$ 
it is homotopic to $id_L $, if $n>2$. 
This is also true if $n=1$ or 2 and $M$ is orientable.
(See \cite[\S29]{[Co]}.)
Therefore $M$ is homotopy equivalent to $L\times S^1 $, 
which fibres over $S^2 $ via the composition of the projection 
to $L$ with the Hopf fibration of $L$ over $S^2 $.
(Hence $w(M)=1$ in these cases also.)

(3)\qua As in part (2), 
if $\pi_1 (E)\cong\mathbb{Z}\rtimes_{-1}\!\mathbb{Z}=\pi_1 (Kb)$ then 
$E$ is homotopy equivalent to $S^2 \times Kb$ and so $w_1 (E)\not=0$,
while $w_2 (E)=0$.
Conversely, if $\chi(M)=0$, $\pi\cong\pi_1 (Kb)$, $M$ is nonorientable and
$w_1(M)^2=w_2(M)=0$ then $M$ is homotopy equivalent to $S^2\times Kb$. 
Suppose now that $\pi$ and $w_1$ satisfy
the second alternative (corresponding to bundles with $\partial\not=0$).            
Let $q:M^+ \to M$ be the orientation double cover.
Then $M^+ $ satisfies the hypotheses of part (3), and so 
there is a map $p^+:M^+\to S^2$ with homotopy fibre $T$.
Now $H^2(q;\mathbb{Z})$ is an epimorphism, since
$H^3(Z/2Z;\mathbb{Z})=H^2(Z/2Z;H^1(M^+;\mathbb{Z}))=0$.
Therefore $p^+=pq$ for some map $p:M\to{S^2}$.
Comparison of the exact sequences of homotopy for $p^+$ and $p$
shows that the homotopy fibre of $p$ must be $Kb$.
As in Theorem 5.2 above $p$ is fibre homotopy 
equivalent to a bundle projection.

(4)\qua There are just two $S^2 $-bundles over $S^2 $, with total spaces 
$S^2 \times S^2$ and $S^2 \tilde \times S^2 =CP^2 \sharp -CP^2 $, respectively. 
Thus the conditions are necessary. 
If $M$ satisfies these conditions then $H^2 (M;\mathbb{Z})\cong\mathbb{Z}^2$ 
and there is an element $u$ in $H^2 (M;\mathbb{Z})$ which generates 
an infinite cyclic direct summand and has square $u\cup u=0$. 
Thus $u=f^* i_2 $ for some map $f:M\to S^2 $.
Since $u$ generates a direct summand there is a homology class $z$ in 
$H_2 (M;\mathbb{Z})$ such that $u\cap z=1$, 
and therefore (by the Hurewicz theorem) there is a map $z:S^2 \to M$
such that $fz$ is homotopic to $id_{S^2 } $. 
The homotopy fibre of $f$ is 1-connected and has
$\pi_2\cong\mathbb{Z}$, by the long exact sequence of homotopy.
It then follows easily from the spectral sequence for $f$ 
that the homotopy fibre has the homology of $S^2 $. 
Therefore $f$ is fibre homotopy equivalent to the projection of
an $S^2 $-bundle over $S^2 $.

(5)\qua Since $\pi_1 (Diff(RP^2 ))=Z/2Z$ (see \cite[page 21]{[EE69]}) 
there are two $RP^2 $-bundles over $S^2 $. 
Again the conditions are clearly necessary.
If they hold we may assume that $u$ generates an infinite cyclic direct 
summand of $H^2(M;\mathbb{Z})$ and that $u=g^* i_2 $ for some map $g:M\to S^2 $.
Let $q:M^+ \to M$ be the orientation double cover and $g^+=gq$.
Since $H_2 (Z/2Z;\mathbb{Z})=0$ the second homology of $M$ is spherical.
Thus there is a map $z=qz^+:S^2\to M$ such that $gz=g^+z^+$ is 
homotopic to $id_{S^2} $.
Hence the homotopy fibre of $g^+ $ is $S^2 $, by case (5).
Since the homotopy fibre of $g$ has fundamental group $Z/2Z$ and is double 
covered by the homotopy fibre of $g^+$ it is homotopy equivalent to $RP^2 $. 
It follows as in Theorem 5.16 that $g$ is fibre homotopy equivalent to the 
projection of an $RP^2 $-bundle over $S^2 $. 
\end{proof}

Theorems 5.2, 5.10 and 5.16 may each be rephrased as giving criteria for maps 
from $M$ to $B$ to be fibre homotopy equivalent to fibre bundle projections.
With the hypotheses of Theorem 5.19 (and assuming also that $\partial=0$ if 
$\chi(M)=0$) we may conclude that a map $f:M\to S^2 $ is fibre homotopy 
equivalent to a fibre bundle projection if and only if $f^*i_2 $ 
generates an infinite cyclic direct summand of $H^2 (M;\mathbb{Z})$.

It follows from Theorem 5.10 that the conditions on the Stiefel-Whitney classes 
are independent of the other conditions when $\pi\cong\pi_1 (F)$. 
Note also that the nonorientable $S^3$- and $RP^3$-bundles over $S^1 $ are not 
$T$-bundles over $S^2 $, while if $M=CP^2\sharp CP^2$ then
$\pi=1$ and $\chi(M)=4$ but $\sigma(M)\not=0$. 
See Chapter 12 for further information on parts (4) and (5).

\section{Bundles over $RP^2$}

Since $RP^2 =Mb\cup D^2$ is the union of a M\"obius band $Mb$ and a disc $D^2$, 
a bundle $p:E\to RP^2$ with fibre $F$ is determined by a bundle over $Mb$ 
which restricts to a trivial bundle over $\partial Mb$,
i.e.\ by a conjugacy class of elements of order dividing 2 
in $\pi_0 (Homeo(F))$, 
together with the class of a gluing map over $\partial Mb=\partial D^2$ 
modulo those which extend across $D^2$ or $Mb$, 
i.e.\ an element of a quotient of $\pi_1 (Homeo(F))$. 
If $F$ is aspherical $\pi_0 (Homeo(F))\cong Out(\pi_1 (F))$, while
$\pi_1 (Homeo(F))\cong \zeta\pi_1 (F)$ \cite{[Go65]}.
                                   
We may summarize the key properties of the algebraic invariants of such 
bundles with $F$ an aspherical closed surface in the following lemma. 
Let $\widetilde{\mathbb{Z}}$ be the nontrivial infinite cyclic $Z/2Z$-module.
The groups $H^1 (Z/2Z;\widetilde{\mathbb{Z}})$, $H^1 (Z/2Z;\mathbb{F}_2)$
and $H^1 (RP^2 ;\widetilde{\mathbb{Z}})$ are canonically isomorphic to $Z/2Z$.

\begin{lemma} 
Let $p:E\to RP^2 $ be the projection of an $F$-bundle, 
where $F$ is an aspherical closed surface, and 
let $x$ be the generator of $H^1 (RP^2 ;\widetilde{\mathbb{Z}})$. Then
\begin{enumerate}
\item $\chi(E)=\chi(F)$;                               

\item $\partial(\pi_2(RP^2))\leq\zeta\pi_1(F)$ and there 
is an exact sequence of groups
\begin{equation*}
\begin{CD}
0\to\pi_2 (E)\to\mathbb{Z}@> \partial >> \pi_1 (F)\to\pi_1 (E)\to Z/2Z\to 1;
\end{CD}
\end{equation*}

\item if $\partial=0$ then $\pi_1 (E)$ acts nontrivially on 
$\pi_2 (E)\cong\mathbb{Z}$ and the covering space $E_F$ 
with fundamental group $\pi_1 (F)$ is homeomorphic to $S^2 \times F$, 
so 
$w_1 (E)|_{\pi_1 (F)} =w_1 (E_F)=w_1 (F)$ 
(as homomorphisms from $\pi_1 (F)$ to $Z/2Z$) and $w_2 (E_F)=w_1 (E_F )^2 $;

\item if $\partial\not=0$ then $\chi(F)=0$, $\pi_1 (E)$ has two ends, 
$\pi_2 (E)=0$ and $Z/2Z$ acts by inversion on $\partial(\mathbb{Z})$;

\item $p^* x^3=0\in H^3 (E;p^*\widetilde{\mathbb{Z}})$.
\end{enumerate}
\end{lemma}

\begin{proof} Condition (1) holds since the Euler characteristic 
is multiplicative in fibrations, while (2) is part of 
the long exact sequence of homotopy for $p$.
The image of $\partial$ is central by \cite{[Go68]}, 
and is therefore trivial unless $\chi(F)=0$.
Conditions (3) and (4) then follow as the homomorphisms in this 
sequence are compatible with the actions of the fundamental groups, 
and $E_F $ is the total space of an $F$-bundle over $S^2$, 
which is a trivial bundle if $\partial=0$, 
by Theorem 5.19.
Condition (5) holds since $H^3 (RP^2;\widetilde{\mathbb{Z}})=0$.
\end{proof}

Let $\pi$ be a group which is an extension of $Z/2Z$ by a normal subgroup $G$, 
and let $t\in\pi$ be an element which maps nontrivially to $\pi/G=Z/2Z$.
Then $u=t^2$ is in $G$ and conjugation by $t$ determines an automorphism 
$\alpha$ of $G$ such that $\alpha(u)=u$ and $\alpha^2$ 
is the inner automorphism given by conjugation by $u$.

Conversely, let $\alpha$ be an automorphism of $G$ whose square is inner, 
say $\alpha^2 (g)=ugu^{-1} $ for all $g\in G$. Let $v=\alpha(u)$. 
Then $\alpha^3 (g)=\alpha^2 (\alpha(g))=u\alpha(g)c^{-1} =\alpha(\alpha^2 (g))
=v\alpha(g)v^{-1} $ for all $g\in G$.
Therefore $vu^{-1} $ is central. 
In particular, if the centre of $G$ is trivial $\alpha$ fixes $u$, and we may 
define an extension
\begin{equation*}
\xi_\alpha : 1\to G\to \Pi_\alpha \to Z/2Z\to 1
\end{equation*}
in which $\Pi_\alpha $ has the presentation 
$\langle G,t_\alpha \mid t_\alpha gt_\alpha^{-1} =\alpha(g),\medspace t_\alpha^2 =u\rangle$.
If $\beta$ is another automorphism in the same outer automorphism class
then $\xi_\alpha$ and $\xi_\beta$ are equivalent extensions.
(Note that if $\beta=\alpha.c_h$, where $c_h $ is conjugation by $h$, then 
$\beta(\alpha(h)uh)=\alpha(h)uh$ and $\beta^2 (g)=\alpha(h)uh.g.(\alpha(h)uh)^{-1} $ 
for all $g\in G$.) 

\begin{lemma} 
If $\,\chi(F)<0$ or $\,\chi(F)=0$ and 
$\partial=0$ then an $F$-bundle over $RP^2$ is determined 
up to isomorphism by the corresponding extension of fundamental groups.
\end{lemma}
    
\begin{proof} 
If $\chi(F)<0$ such bundles and extensions are each determined by an element 
$\xi$ of order 2 in $Out(\pi_1 (F))$.
If $\chi(F)=0$ bundles with $\partial=0$ are the restrictions of bundles over
$RP^\infty=K(Z/2Z,1)$ (compare Lemma 4.10).
Such bundles are determined by an element 
$\xi$ of order 2 in $Out(\pi_1 (F))$ and a cohomology class 
in $H^2(Z/2Z;\zeta\pi_1(F)^\xi)$, by Lemma 5.1, and so correspond bijectively
to extensions also. 
\end{proof}

\begin{lemma} 
Let $M$ be a $PD_4$-complex with fundamental 
group $\pi$. 
A map $f:M\to RP^2 $ is fibre homotopy equivalent to the projection of a bundle over $RP^2 $ 
with fibre an aspherical closed surface if $\pi_1(f)$ is an epimorphism 
and either
\begin{enumerate}
\item $\chi(M)\leq 0$ and $\pi_2(f)$ is an isomorphism; or

\item $\chi(M)=0$, $\pi$ has two ends and $\pi_3(f)$ is an isomorphism.
\end{enumerate}
\end{lemma}
            
\begin{proof}
In each case $\pi$ is infinite, by Lemma 3.14.
In case (1) $H^2(\pi;\mathbb{Z}[\pi])$ is infinite cyclic(by Lemma 3.3) 
and so $\pi$ has one end, by Bowditch's Theorem.
Hence $\widetilde M\simeq S^2$.
Moreover the homotopy fibre of $f$ is aspherical, 
and its fundamental group is a surface group.
(See Chapter X for details.)
In case (2) $\widetilde M\simeq S^3 $, by Corollary 4.5.2.
Hence the lift $\tilde f:\widetilde M\to S^2 $ is fibre 
homotopy equivalent to the Hopf map, 
and so induces isomorphisms on all higher homotopy groups. 
Therefore the homotopy fibre of $f$ is aspherical. 
As $\pi_2 (M)=0$ the fundamental group of the homotopy fibre 
is a (torsion-free) infinite cyclic extension of $\pi$ 
and so must be either $\mathbb{Z}^2 $ or 
$\mathbb{Z}\rtimes_{-1}\!\mathbb{Z}$. 
Thus the homotopy fibre of $f$ is homotopy equivalent to $T$ or $Kb$.
In both cases the argument of Theorem 5.2 now shows 
that $f$ is fibre homotopy equivalent to a surface bundle projection. 
\end{proof}

\section{Bundles over $RP^2$ with $\partial=0$}

Let $F$ be a closed aspherical surface and 
$p:M\to{RP^2}$ be a bundle with fibre $F$,
and such that $\pi_2(M)\cong\mathbb{Z}$.
(This condition is automatic if $\chi(F)<0$.) 
Then $\pi=\pi_1(M)$ acts nontrivially on $\pi_2(M)$.
The covering space $M_\kappa$ associated to the kernel $\kappa$ 
of the action is an $F$-bundle over $S^2$,
and so $M_\kappa\cong{S^2}\times{F}$,
since all such bundles are trivial.
In particular, $v_2(M)\in{H^2(\pi;\mathbb{F}_2)}$, and $v_2(M)|_\kappa=0$.
The projection admits a section if and only if $\pi\cong\kappa\rtimes{Z/2Z}$. 

Our attempt (in the original version of this book)
to characterize more general surface bundles over $RP^2$ had an error
(in the claim that restriction from $H^2(RP^2;\mathbb{Z}^u)$ to 
$H^2(S^2;\mathbb{Z})$ is an isomorphism). 
We provide instead several partial results.
Further progress might follow from a better understanding of maps from
4-complexes to $RP^2$.
The reference \cite{[Si67]} cited in the former (flawed) theorem 
of this section remains potentially useful here.

The product $RP^2\times{F}$ is easily characterized.

\begin{theorem}
Let $M$ be a closed $4$-manifold with fundamental group $\pi$,
and let $F$ be an aspherical closed surface.
Then the following are equivalent.
\begin{enumerate}
\item$M\simeq{RP^2}\times{F}$;

\item$\pi\cong{Z/2Z}\times\pi_1(F)$, $\chi(M)=\chi(F)$ and $v_2(M)=0$;

\item$\pi\cong{Z/2Z}\times\pi_1(F)$, $\chi(M)=\chi(F)$ and $M\simeq{E}$,
where $E$ is the total space of an $F$-bundle over $RP^2$.
\end{enumerate}
\end{theorem}

\begin{proof}
Clearly $(1)\Rightarrow(2)$ and (3).
If (2) holds then $M$ is homotopy equivalent to the total space of an
$RP^2$-bundle over $F$, by Theorem 5.16.
This bundle must be trivial since $v_2(M)=0$.
If (3) holds then there are maps $q:M\to{F}$ and $p:M\to{RP^2}$
such that $\pi_1(p)$ and $\pi_1(q)$ are the projections of $\pi$ onto its
factors and $\pi_2(p)$ is surjective.
The map $(p,q):M\to{RP^2}\times{F}$ is then a homotopy equivalence.
\end{proof}

The implication $(3)\Rightarrow(1)$ fails if $F=RP^2$ or $S^2$.

The characterization of bundles with sections is based on a study of $S^2$-orbifold bundles. (See Chapter 10 below and \cite{[Hi13]}.)

\begin{theorem}
Let $F$ be an aspherical closed surface.
A closed orientable $4$-manifold $M$ is homotopy equivalent to 
the total space of an $F$-bundle over $RP^2$ with a section 
if and only if $\pi=\pi_1(M)$ has an element of order $2$,
$\pi_2(M)\cong\mathbb{Z}$ and 
$\kappa=\mathrm{Ker}(u)\cong\pi_1(F)$, 
where $u$ is the natural action of $\pi$ on $\pi_2(M)$.
\end{theorem}

\begin{proof}
The conditions are clearly necessary.
Suppose that they hold.
We may assume that $\pi$ is not a direct product $\kappa\times{Z/2Z}$.
Therefore $M$ is not homotopy equivalent to an $RP^2$-bundle space.
Hence it is homotopy equivalent to the total space $E$
of an $S^2$-orbifold bundle over a 2-orbifold $B$. 
(See Corollary 10.8.1 below.)
The involution $\zeta$ of $F$ corresponding to the orbifold
covering has non-empty fixed point set, since $\pi$ has torsion.
Let $M_{st}=S^2\times{F}/\sim$, where $(s,f)\sim(-s,\zeta(f))$.
Then $M_{st}$ is the total space of an $F$-bundle over $RP^2$,
and th e fixed points of $\zeta$ determine sections of this bundle.

The double cover of $E$ corresponding to $\kappa$ is an $S^2$-bundle over $F$.
Since $M$ is orientable and $\kappa$ acts trivially on $\pi_2(M)$,
$F$ must also be orientable and the covering involution 
of $F$ over $B$ must be orientation-reversing.
Since $\pi$ has torsion $\Sigma{B}$ is a non-empty union of 
reflector curves, and since $F$ is orientable these are ``untwisted".
Therefore $M\simeq{M_{st}}$ \cite[Corollary 4.8]{[Hi13]}.
\end{proof}

Orientability is used here mainly to ensure that $B$ has a reflector curve.

When $\pi$ is torsion-free $M$ is homotopy equivalent to 
the total space of an $S^2$-bundle over a surface $B$, 
with $\pi=\pi_1(B)$ acting nontrivially on the fibre.
Inspection of the geometric models for such bundle spaces 
given in Chapter 10 below
shows that if also $v_2(M)\not=0$ then the bundle space fibres over $RP^2$. 
Is the condition $v_2(M)\not=0$ necessary?

\section{Sections of surface bundles}

If a bundle $p:E\to{B}$ with base and fibre aspherical surfaces 
has a section then its fundamental group sequence splits.
The converse holds if the action $\xi$ can be realized by
a group of based self homeomorphisms of the fibre $F$.
(This is so if $F=T$ or $Kb$.)
The sequence splits if and only if the action factors through $Aut(\pi_1(F))$ 
and the class of the extension in $H^2(\pi_1(B);\zeta\pi_1(F))$ is 0.
This cohomology group is trivial if $\chi(F)<0$, 
and the class is easily computed if $\chi(F)=0$.
In particular, if  $B$ is orientable and $F=T$ 
then $p$ has a section if and only if 
$H_1(E;\mathbb{Z})\cong{H_1(B;\mathbb{Z})}\oplus{H_0(B;\pi_1(F))}$.
(The $T$-bundles over $T$ which are coset spaces of the nilpotent 
Lie groups $Nil^3\times\mathbb{R}$ and $Nil^4$ 
do not satisfy this criterion, 
and so do not have sections.)

If $p_*$ splits and $s$ and $s'$ are two sections determining 
the same lift $\widetilde\xi:\pi_1(B)\to{Aut}(\pi_1(F))$ 
then $s'(g)s(g)^{-1}$ is in $\zeta\pi_1(F)$, 
for all $g\in\pi_1(B)$.
Hence the sections are parametrized (up to conjugation 
by an element of $\pi_1(F)$) by $H^1(\pi_1(B);\zeta\pi_1(F))$.
In particular, if $\chi(F)<0$ and $p_*$ has a section then the section 
is unique up to conjugation by an element of $\pi_1(F)$.

It follows easily from Theorem 5.19 that nontrivial bundles over $S^2$ 
with aspherical fibre do not admit sections.

See also \cite{[Hi15]}.

%% file: m5-6.tex
\chapter{Simple homotopy type and surgery}

The problem of determining the high-dimensional manifolds within a given 
homotopy type has been successfully reduced to the determination of normal 
invariants and surgery obstructions.
This strategy applies also in dimension 4, 
provided that the fundamental group is in the class $SA$
generated from groups with subexponential growth by extensions 
and increasing unions \cite{[FT95]}.
(Essentially all the groups in this class that we shall discuss in this book 
are in fact virtually solvable.)
We may often avoid this hypothesis by using
5-dimensional surgery to construct $s$-cobordisms.

We begin by showing that the Whitehead group of the fundamental group is 
trivial for surface bundles over surfaces, most circle bundles 
over geometric 3-manifolds and for many mapping tori. 
In \S2 we define the modified surgery structure set, parametrizing 
$s$-cobordism classes of simply homotopy equivalences of closed 4-manifolds.
This notion allows partial extensions of surgery arguments to situations where 
the fundamental group is not elementary amenable.
Although many papers on surgery do not explicitly consider the 4-dimensional 
cases, their results may often be adapted to these cases. 
In \S3 we comment briefly on approaches to the $s$-cobordism theorem
and classification using stabilization by connected sum with copies 
of $S^2\times S^2$ or by cartesian product with $\mathbb{R}$.          

In \S4 we show that 4-manifolds $M$ such that $\pi=\pi_1 (M)$ is torsion-free
virtually poly-$Z$ and $\chi(M)=0$ are determined up to homeomorphism by their 
fundamental group (and Stiefel-Whitney classes, if $h(\pi)<4$). 
We also characterize 4-dimensional mapping tori with torsion-free, 
elementary amenable fundamental group and show that the structure sets 
for total spaces of $RP^2$-bundles over $T$ or $Kb$ are finite.
In \S5 we extend this finiteness to $RP^2$-bundle spaces over
closed hyperbolic surfaces and show that total spaces of bundles with 
fibre $S^2$ or an aspherical closed surface over aspherical bases
are determined up to $s$-cobordism by their homotopy type.
(We shall consider bundles with base or fibre geometric 3-manifolds 
in Chapter 13.)

\vfil\eject
\section{The Whitehead group}

In this section we shall rely heavily upon the work of Waldhausen 
in \cite{[Wd78]}.
The class of groups {\it Cl} is the smallest class of groups containing the 
trivial group and which is closed under generalised free products and HNN 
extensions with amalgamation over regular coherent subgroups and under 
filtering direct limit.
This class is also closed under taking subgroups 
\cite[Proposition 19.3]{[Wd78]}.
If $G$ is in {\it Cl} then so is $G\times\mathbb{Z}^n$,
and $Wh(G)=\tilde{K}(\mathbb{Z}[G])=0$ \cite[Theorem 19.4]{[Wd78]}.   
The argument for this theorem actually shows that if $G\cong A*_C B$ 
and $C$ is regular coherent then there are {\it Mayer-Vietoris sequences\/}:
\[
Wh(A)\oplus{Wh(B)}\!\to\!{Wh(G)}\!\to\!\tilde{K}_0(\mathbb{Z}[C])\!
\to\!\tilde{K}_0(\mathbb{Z}[A])\oplus\tilde{K}_0(\mathbb{Z}[B])\!\to\!
\tilde{K}_0(\mathbb{Z}[G])\!\to\!0
\]
and similarly if $G\cong A*_C $.
(See \cite[\S17.1.3 and \S17.2.3]{[Wd78]}.)

The class {\it Cl} contains all free groups and poly-$Z$ groups 
and the class $\mathcal{ X}$ of Chapter 2.
(In particular, all the groups $Z*_m$ are in {\it Cl}.)
Since every $PD_2 $-group is either poly-$Z$ or is the generalised free product
of two free groups with amalgamation over infinite cyclic subgroups it is 
regular coherent, and is in {\it Cl}.
Hence homotopy equivalences between $S^2$-bundles over aspherical surfaces 
are simple.
The following extension implies the corresponding result for quotients of 
such bundle spaces by free involutions.
                                    
\begin{theorem}
Let $\pi$ be a semidirect product $\rho\rtimes(Z/2Z)$ where $\rho$ 
is a surface group. Then $Wh(\pi)=0$.
\end{theorem}

\begin{proof}                                                    
Assume first that $\pi\cong\rho\times(Z/2Z)$.
Let $\Gamma=\mathbb{Z}[\rho]$. There is a cartesian square expressing 
$\Gamma[Z/2Z]=\mathbb{Z}[\rho\times(Z/2Z)]$
as the pullback of the reduction of coefficients map from $\Gamma$ to 
$\Gamma_2=\Gamma/2\Gamma =\mathbb{Z}/2\mathbb{Z}[\rho]$ over itself. 
(The two maps from $\Gamma[Z/2Z]$ to
$\Gamma$ send the generator of $Z/2Z$ to $+1$ and $-1$, respectively.)
The Mayer-Vietoris sequence for algebraic $K$-theory traps 
$K_1 (\Gamma[Z/2Z])$ between $K_2 (\Gamma_2)$ and $K_1 (\Gamma)^2 $ \cite[Theorem 6.4]{[Mi]}.
Now since $c.d.\rho=2$ the higher $K$-theory of $R[\rho]$ 
can be computed in terms
of the homology of $\rho$ with coefficients in the $K$-theory of $R$ 
(cf.\ the Corollary to Theorem 5 of the introduction of \cite{[Wd78]}). 
In particular, the map from $K_2 (\Gamma)$
to $K_2 (\Gamma_2 )$ is onto, while $K_1 (\Gamma)=K_1 (\mathbb{Z})\oplus(\rho/\rho')$ and
$K_1 (\Gamma_2 )=\rho/\rho'$. It now follows easily that 
$K_1 (\Gamma[Z/2Z])$ is generated by the
images of $K_1 (\mathbb{Z})=\{\pm 1\} $ and $\rho\times(Z/2Z)$, 
and so $Wh(\rho\times(Z/2Z))=0$. 

If $\pi=\rho\rtimes (Z/2Z)$ is not such a direct product it
is isomorphic to a discrete subgroup of $Isom(\mathbb{X})$ which acts properly 
discontinuously on $X$, where $\mathbb{X}=\mathbb{E}^2$ or $\mathbb{H}^2$. 
(See \cite{[Zi]}.)
The singularities of the corresponding 2-orbifold $\pi\backslash{X}$ 
are either cone points of order 2 or reflector curves; 
there are no corner points and no cone points of higher order.
Let $|\pi\backslash{X}|$ be the surface obtained by forgetting 
the orbifold structure of $\pi\backslash{X}$,
and let $m$ be the number of cone points.
Then $\chi(|\pi\backslash{X}|)-(m/2)=\chi_{orb} (\pi\backslash{X})\leq 0$, 
by the Riemann-Hurwitz formula \cite{[Sc83']}, so either
$\chi(|\pi\backslash{X}|)\leq0$ or $\chi(|\pi\backslash{X}|)=1$ and 
$m\geq 2$ or $|\pi\backslash{X}|\cong S^2$ and $m\geq 4$.

We may separate $\pi\backslash{X}$ along embedded circles 
(avoiding the singularities) into pieces which are either 
(i) discs with at least two cone points; 
(ii) annuli with one cone point;
(iii) annuli with one boundary a reflector curve; or
(iv) surfaces with nonempty boundary, other than $D^2$ or the annulus. 
In each case the inclusions of the separating circles induce monomorphisms on orbifold fundamental groups,
and so $\pi$ is a generalized free product with amalgamation 
over copies of $\mathbb{Z}$ of groups of the form 
(i) $*^m(Z/2Z)$ (with $m\geq 2$); 
(ii) $\mathbb{Z}*(Z/2Z)$; (iii) $\mathbb{Z}\oplus(Z/2Z)$; 
or (iv) $F(r)$ (with $r\geq2$), 
by the Van Kampen theorem for orbifolds \cite{[Sc83]}.
The Mayer-Vietoris sequences for algebraic $K$-theory now give $Wh(\pi)=0$.
\end{proof}

The argument for the direct product case is based on one for showing that 
$Wh(\mathbb{Z}\oplus(Z/2Z))=0$ from \cite{[Kw86]}.

Not all such orbifold groups arise in this way. For instance, the orbifold 
fundamental group of a torus with one cone point of order 2 has the
presentation  $\langle x,y\mid [x,y]^2=1\rangle$.
Hence it has torsion-free abelianization, 
and so cannot be a semidirect product as above.

The orbifold fundamental groups of flat 2-orbifolds 
are the 2-dimensional crystallographic groups.                                   Their finite subgroups are cyclic or dihedral, 
of order properly dividing 24, 
and have trivial Whitehead group. 
In fact $Wh(\pi)=0$ for $\pi$ any such 2-dimensional crystallographic group 
\cite{[Pe98]}. 
(If $\pi$ is the fundamental group of an orientable hyperbolic
2-orbifold with $k$ cone points of orders $\{n_1,\dots, n_k\}$
then $Wh(\pi)\cong\oplus_{i=1}^kWh(Z/n_iZ)$ \cite{[LS00]}.)

The argument for the next result is essentially due to F.T.Farrell.

\begin{theorem} 
If $\pi$ is an extension of $\pi_1 (B)$ by $\pi_1 (F)$
where $B$ and $F$ are aspherical closed surfaces then 
$Wh(\pi)=\tilde{K}_0(\mathbb{Z}[\pi])=0$.
\end{theorem}

\begin{proof} If $\chi(B)<0$ then $B$ admits 
a complete riemannian metric of constant negative curvature $-1$.
Moreover the only virtually poly-$Z$ subgroups of $\pi_1 (B)$ 
are 1 and $\mathbb{Z}$.
If $G$ is the preimage in $\pi$ of such a subgroup then 
$G$ is either $\pi_1 (F)$ or is the group of a Haken 3-manifold. 
It follows easily that for any $n\geq 0$ the group
$G\times\mathbb{Z}^n $ is in {\it Cl} and so $Wh(G\times\mathbb{Z}^n )=0$.
Therefore any such $G$ is $K$-flat and so the bundle is admissible, 
in the terminology of \cite{[FJ86]}. 
Hence $Wh(\pi)=\tilde{K}_0(\mathbb{Z}[\pi])=0$ by the main result of that paper.

If $\chi(B)=0$ then this argument does not work, 
although if moreover $\chi(F)=0$ then $\pi$ is poly-$Z$, 
so $Wh(\pi)=\tilde{K}_0(\mathbb{Z}[\pi])=0$ \cite[Theorem 2.13]{[FJ]}.
We shall sketch an argument of Farrell for the general case.
Lemma 1.4.2 and Theorem 2.1 of \cite{[FJ93]} together 
yield a spectral sequence                         
(with coefficients in a simplicial cosheaf) whose $E^2 $ term is
$H_i (X/\pi_1 (B);Wh'_j (p^{-1} (\pi_1 (B)^x)))$ and 
which converges to $Wh'_{i+j} (\pi)$.
Here $p:\pi\to\pi_1 (B)$ is the epimorphism of the extension 
and $X$ is a certain universal $\pi_1 (B)$-complex 
which is contractible and such that all the nontrivial isotropy subgroups
$\pi_1 (B)^x$ are infinite cyclic and the fixed point set of each infinite cyclic subgroup is a contractible (nonempty) subcomplex.
The Whitehead groups with negative indices are the lower $K$-theory of $\mathbb{Z}[G]$
(i.e., $Wh'_n (G)=K_n (\mathbb{Z}[G])$ for all $n\leq -1$), 
while $Wh'_0 (G)=\tilde{K}_0(\mathbb{Z}[G])$ and $Wh'_1 (G)=Wh(G)$.
Note that $Wh'_{-n} (G)$ is a direct summand of 
${Wh(G\times\mathbb{Z}^{n+1})}$.
If $i+j>1$ then $Wh'_{i+j}(\pi)$ agrees rationally with the higher Whitehead group $Wh_{i+j} (\pi)$.
Since the isotropy subgroups $\pi_1 (B)^x$ are infinite cyclic or trivial
$Wh(p^{-1} (\pi_1 (B)^x)\times\mathbb{Z}^n)=0$ for all $n\geq 0$, 
by the argument of the above paragraph, 
and so $Wh'_j (p^{-1} (\pi_1 (B)^x))=0$ if $j\leq 1$. 
Hence the spectral sequence gives $Wh(\pi)=\tilde{K}_0(\mathbb{Z}[\pi])=0$. 
\end{proof}

A closed 3-manifold is a {\it Haken manifold} if it is irreducible 
and contains an incompressible 2-sided surface.
Every aspherical closed 3-manifold $N$ is either Haken, hyperbolic
or Seifert-fibred, by the work of Perelman \cite{[B-P]}, and so 
either has an infinite solvable fundamental group or it has a $JSJ$
decomposition along a finite family of disjoint incompressible tori and Klein 
bottles so that the complementary components are Seifert fibred or hyperbolic.
Every closed 3-manifold with a metric of non-positive  curvature 
is {\it virtually fibred\/} (i.e., finitely covered by a mapping torus),
and so every aspherical closed 3-manifold is virtually Haken \cite{[Ag13],[PW18]}.

If an aspherical closed 3-manifold has a $JSJ$ decomposition
with at least one hyperbolic component then it has 
a metric of non-positive curvature \cite{[Lb95]}.
Otherwise it is a {\it graph manifold}: 
either it has solvable fundamental group or 
it has a $JSJ$ decomposition into Seifert fibred pieces. 
It is a {\it proper} graph manifold if the minimal such $JSJ$
decomposition is non-trivial.
A criterion for a proper graph manifold  to be virtually fibred
is given in \cite{[Ne97]}.

\begin{theorem}
Let $N$ be a connected sum of aspherical graph manifolds,
and let $\nu=\pi_1(N)$ and $\pi=\nu\rtimes_\theta\mathbb{Z}$,
where $\theta\in{Aut}(\nu)$.
Then $\nu\times\mathbb{Z}^n$ is regular coherent, and
$Wh(\pi\times\mathbb{Z}^n)=\tilde{K}_0(\mathbb{Z}[\pi\times\mathbb{Z}^n])=0$,
for all $n\geq0$.
\end{theorem}
                      
\begin{proof} 
The group $\nu$ is either polycyclic or is a generalized free product 
with amalgamation along poly-$Z$ subgroups 
(1, $\mathbb{Z}^2 $ or $\mathbb{Z}\rtimes_{-1}\!\mathbb{Z}$) 
of fundamental groups of Seifert fibred 3-manifolds (possibly with boundary).
The group rings of torsion-free polycyclic groups are regular noetherian, 
and hence regular coherent.
If $G$ is the fundamental group of a Seifert fibred 3-manifold then it
has a subgroup $G_o$ of finite index which is a central extension of 
the fundamental group of a surface $B$ (possibly with boundary) 
by $\mathbb{Z}$. 
We may assume that $G$ is not solvable and hence that $\chi(B)<0$.
If $\partial B$ is nonempty then $G_o\cong\mathbb{Z}\times\! F$ 
and so is an iterated generalized free product of copies of $\mathbb{Z}^2 $, 
with amalgamation along infinite cyclic subgroups.
Otherwise we may split $B$ along an essential curve and represent $G_o $
as the generalised free product of two such groups, 
with amalgamation along a copy of $\mathbb{Z}^2 $.
In both cases $G_o$ is regular coherent, and therefore so is $G$, 
since $[G:G_o]<\infty$ and $c.d.G<\infty$.
             
Since $\nu$ is the generalised free product with amalgamation 
of regular coherent groups,
with amalgamation along poly-$Z$ subgroups, 
it is also regular coherent.
Hence so is $\nu\times\mathbb{Z}^n$.
Let $N_i $ be an irreducible summand of $N$ and let $\nu_i =\pi_1 (N_i )$.
If $N_i $ is Haken then $\nu_i$ is in {\it Cl}
and so $Wh(\nu_i \times\mathbb{Z}^n)=0$, for all $n\geq0$.
Otherwise $N_i $ is a Seifert fibred 3-manifold which is not sufficiently large,
and the argument of [Pl80] extends easily to prove this.
Since $\tilde{K}_0(\mathbb{Z}[\sigma])$ is a direct summand of 
$Wh(\sigma\times\mathbb{Z})$, for any group $\sigma$,
we have $\tilde{K}_0(\mathbb{Z}[\nu_i\times\mathbb{Z}^n])=0$,
for all ${n\geq0}$.
The Mayer-Vietoris sequences for algebraic $K$-theory now give, 
firstly, $Wh(\nu\times\mathbb{Z}^n)=
\tilde K_0(\mathbb{Z}[\nu\times\mathbb{Z}^n])=0$, 
and then $Wh(\pi\times\mathbb{Z}^n)=
\tilde K_0 (\mathbb{Z}[\pi\times\mathbb{Z}^n])=0$ also. 
\end{proof}

All 3-manifold groups are coherent as {\it groups} \cite{[Hm]}.
If we knew that their group {\it rings} were regular coherent then 
we could use \cite{[Wd78]} instead of \cite{[FJ86]} 
to give a purely algebraic proof of Theorem 6.2, 
for as surface groups are free products of free groups with amalgamation 
over an infinite cyclic subgroup,
an extension of one surface group by another is a free product of groups 
with $Wh=0$,
amalgamated over the group of a surface bundle over $S^1$.
Similarly, we could deduce from \cite{[Wd78]} 
and the work of Perelman \cite{[B-P]}
that $Wh(\nu\rtimes_\theta\mathbb{Z})=0$
for any torsion-free 3-manifold group $\nu=\pi_1(N)$ 
where $N$ is a closed 3-manifold.

\begin{theorem}
Let $N$ be a closed 3-manifold such that $\nu=\pi_1(N)$ is torsion-free,
and let $\mu$ be a group with an infinite cyclic normal subgroup
$A$ such that $\mu/A\cong\nu$. 
Then $Wh(\mu )=Wh(\nu )=0$.
\end{theorem}

\begin{proof} 
Let $N=\sharp_{1\leq i\leq n}{N_i}$ be the factorization of $N$ into
irreducibles, and let $\nu\cong*_{1\leq i\leq n}\nu_i $, 
where $\nu_i=\pi_1(N_i)$.
The irreducible factors are either Haken, hyperbolic or Seifert fibred,
by the work of Perelman \cite{[B-P]}.
Let $\mu_i $ be the preimage of $\nu_i $ in $\mu $, for $1\leq i\leq n$. 
Then $\mu $ is the generalized free product of the $\mu_i $'s, 
amalgamated over infinite cyclic subgroups. For all $1\leq i\leq n$ we have 
$Wh(\mu_i )=0$, by \cite[Lemma 1.1]{[St84]} if $K(\nu_i ,1)$ is Haken, 
by the main result of \cite{[FJ86]} if it is hyperbolic, 
by an easy extension of the argument of \cite{[Pl80]} if it is Seifert fibred 
but not Haken and by \cite[Theorem 19.5]{[Wd78]} 
if $\nu_i $ is infinite cyclic. 
The Mayer-Vietoris sequences for algebraic $K$-theory now give 
$Wh(\mu )=Wh(\nu )=0$ also. 
\end{proof}

Theorem 6.4 may be used to strengthen Theorem 4.11 to give criteria 
for a closed 4-manifold $M$ to be {\it simple} homotopy equivalent 
to the total space of an $S^1 $-bundle, 
if $\pi_1 (M)$ is torsion-free.

\section{The $s$-cobordism structure set}

The TOP structure set for a closed 4-manifold $M$ 
with fundamental group $\pi$ and orientation 
character $w:\pi\to\{\pm1\}$ is
\[
S_{TOP} (M)=\{ f:N\to M\mid N~a~TOP~4\mathrm{-}manifold,
~f~a~simple~h.e.\}/\!\sim,
\]
where $f_1\sim{f_2} $ if $f_1=f_2h$ for some homeomomorphism $h:N_1\to{N_2}$.
If $\pi$ is ``good" (e.g., if it is in $SA$)
then $L^s_5(\pi,w)$ acts on the structure set $S_{TOP}(M)$,
and the orbits of the action $\omega$ correspond to the normal invariants 
$\eta(f)$ of simple homotopy equivalences \cite{[FQ],[FT95]}.
The surgery sequence                                                          
\[
[SM;G/TOP]\buildrel\sigma_5\over\longrightarrow {L_5^s(\pi,w)}
\buildrel\omega\over\longrightarrow {S_{TOP}(M)}
\buildrel\eta\over\longrightarrow [M;G/TOP]
\buildrel\sigma_4\over\longrightarrow {L_4^s(\pi,w)}
\]
may then be identified with the algebraic surgery sequence of \cite{[Rn]}.
The additions on the homotopy sets $[X,G/TOP]$ derive from an
$H$-space structure on $G/TOP$. 
(In low dimensions this is unambiguous,
as $G/TOP$ has Postnikov 5-stage $K(Z/2Z,2)\times{K(\mathbb{Z},4)}$,
which has an unique $H$-space structure.) 
We shall not need to specify the addition on $S_{TOP}(M)$. 

As it is not yet known whether 5-dimensional $s$-cobordisms over other 
fundamental groups are products,
we shall redefine the structure set by setting
\[
S^s_{TOP} (M)=\{ f:N\to M\mid N~a~TOP~4\mathrm{-}manifold,
~f~a~simple~h.e.\}/\!\approx,
\]
where $f_1\approx f_2 $ if there is a map $F:W\to M$ with domain $W$ an 
$s$-cobordism with $\partial W=N_1 \cup N_2 $ and $F|_{N_i}=f_i $ for $i=1,2$.
If the $s$-cobordism theorem holds over $\pi$ 
this is the usual TOP structure set for $M$.
We shall usually write $L_n (\pi,w)$ for $L_n^s (\pi,w)$ if $Wh(\pi)=0$ and 
$L_n (\pi)$ if moreover $w$ is trivial. 
When the orientation character is nontrivial and otherwise clear from 
the context we shall write $L_n (\pi,-)$. 
We shall say that a closed 4-manifold is {\it $s$-rigid}
if it is determined up to $s$-cobordism by its homotopy type.                    
The homotopy set $[M;G/TOP]$ may be identified with the set of normal maps 
$(f,b)$, where $f:N\to M$ is a degree 1 map and 
$b$ is a stable framing of $T_N\oplus f^*\xi$, 
for some TOP $R^n$-bundle $\xi$ over $M$.
If $f:N\to M$ is a homotopy equivalence, with homotopy inverse $h$, 
let $\xi=h^*\nu_N$ and $b$ be the framing determined by a homotopy 
from $hf$ to $id_N$.
Let $\hat{f}\in[M,G/TOP]$ be the homotopy class corresponding to $(f,b)$. 
Let $k_2$ generate $H^2(G/TOP;\mathbb{F}_2)\cong Z/2Z$
and $l_4$ generate $H^4(G/TOP;\mathbb{Z})\cong\mathbb{Z}$,
with image $[l_4]$ in $H^4(G/TOP;\mathbb{F}_2)$.
The function from $[M;G/TOP]$ to $H^2(M;\mathbb{F}_2)\oplus H^4(M;\mathbb{Z})$
which sends $\hat f$ to $(\hat f^*(k_2),\hat f^*(l_4))$ is an isomorphism. 
Let $KS(M)\in H^4(M;\mathbb{F}_2)$ be the Kirby-Siebenmann obstruction 
to lifting the TOP normal fibration of $M$ to a vector bundle.
If $\hat{f}$ is a normal map then 
\[KS(M)-(f^*)^{-1}KS(N)=\hat{f}^*(k_2^2+[l_4]),\]
and $\hat{f}$ factors through $G/PL$ if and only if this difference is 0
\cite{[KT02]}.
If $M$ is orientable then $\hat{f}^*(l_4)([M])=(\sigma(M)-\sigma(N))/8$,
where $\sigma(M)$ is the signature of the intersection pairing on
$H_2(M;\mathbb{Z})$, and so
\[
(KS(M)-(f^*)^{-1}KS(N)-\hat{f}^*(k_2)^2)([M])\equiv(\sigma(M)-\sigma(N))/8
\quad {mod}~(2).
\]
The {\sl Kervaire-Arf invariant} of a normal map $\hat g:N^{2q}\to G/TOP$
is the image of the surgery obstruction in $L_{2q}(Z/2Z,-)=Z/2Z$
under the homomorphism induced by the orientation character,
$c(\hat g)=L_{2q}(w_1(N))(\sigma_{2q}(\hat g))$.
The argument of Theorem 13.B.5 of [Wl] may be adapted to show that
there are universal classes $K_{4i+2}$
in $H^{4i+2}(G/TOP;\mathbb{F}_2)$ (for $i\geq0$) such that 
\[c(\hat g)=(w(M)\cup \hat g^*((1+Sq^2+Sq^2Sq^2)\Sigma K_{4i+2}))\cap[M].
\]
Moreover $K_2=k_2$, since $c$ induces the isomorphism $\pi_2(G/TOP)=Z/2Z$.
In the 4-dimensional case this expression simplifies to
\[
c(\hat g)=(w_2(M)\cup\hat g^*(k_2)+\hat g^*(Sq^2k_2))([M])=
(w_1(M)^2\cup\hat g^*(k_2))([M]).
\]
The {\sl codimension-2 Kervaire invariant} of a 4-dimensional
normal map $\hat g$ is $kerv(\hat g)=\hat g^*(k_2)$.
Its value on a 2-dimensional homology class represented by an immersion 
$y:Y\to M$ is the Kervaire-Arf invariant of the normal map induced over 
the surface $Y$.

The structure set may overestimate the number of homeomorphism types within 
the homotopy type of $M$, if $M$ has self homotopy equivalences which are not 
homotopic to homeomorphisms.
Such ``exotic" self homotopy equivalences may often
be constructed as follows. 
Given $\alpha:S^2 \to M$, let $\beta:S^4\to M$ be the composition 
$\alpha\eta{S\eta}$, where $\eta$ is the Hopf map, 
and let $s:M\to M\vee S^4 $ be the pinch map obtained 
by shrinking the boundary of a 4-disc in $M$.
Then the composite $f_\alpha =(id_M\vee\beta)s$ is a self homotopy equivalence 
of $M$.
         
\begin{lemma} 
{\rm[No64]}\qua 
Let $M$ be a closed $4$-manifold and let 
$\alpha:S^2\to M$ be a map such that $\alpha_*[S^2]\not=0$ in 
$H_2(M;\mathbb{F}_2)$ and $\alpha^*w_2(M)=0$.
Then $kerv(\widehat{f_\alpha})\not=0$ and so 
$f_\alpha$ is not normally cobordant to a homeomorphism.
\end{lemma}

\begin{proof} 
Since $\alpha_*[S^2]\not=0$ there is a $u\in H_2(M;\mathbb{F}_2)$ 
such that $\alpha_* [S^2].u=1$.
This class may be realized as $u=g_*[Y]$ where $Y$ is a closed surface
and $g:Y\to M$ is transverse to $f_\alpha$.
Then $g^*kerv(\widehat{f_\alpha})[Y]$ is the Kervaire-Arf invariant of the 
normal map induced over $Y$ and is nontrivial.
(See \cite[Theorem 5.1]{[CH90]} for details.) 
\end{proof}

The family of surgery obstruction maps may be identified with 
a natural transformation from $\mathbb{L}_0 $-homology to $L$-theory. 
(In the nonorientable case we must use $w$-twisted $\mathbb{L}_0 $-homology.)
In dimension 4 the cobordism invariance of surgery obstructions
(as in \cite[\S 13B]{[Wl]}) leads to the following formula.
 
\begin{theorem} 
{\rm[Da05]}\qua
There are homomorphisms $I_0:\!H_0(\pi;\mathbb{Z}^w)\!\to\! L_4(\pi,w)$ and
$\kappa_2:H_2(\pi;\mathbb{F}_2)\to L_4(\pi,w)$ such that for any
$\hat f:M\to G/TOP$ the surgery obstruction is 
$\sigma_4(\hat f)=I_0(c_{M*}(\hat f^*(l_4)\cap[M]))+
\kappa_2(c_{M*}(kerv(\hat f)\cap[M]))$.
\qed
\end{theorem}

In the orientable case the signature homomorphism from $L_4(\pi)$ 
to $\mathbb{Z}$ is a left inverse for $I_0:\mathbb{Z}\to L_4(\pi)$, 
but in general $I_0$ is not injective.
This formula can be made somewhat more explicit as follows.

\begin{thmm}
[6.6$^\prime$]
{\rm[Da05]}\qua
If $\hat f=(f,b)$ where $f:N\to M$ is a degree $1$ map
then the surgery obstructions are given by

$
\sigma_4(\hat f)=I_0((\sigma(N)-\sigma(M))/8)+
\kappa_2(c_{M*}(kerv(\hat f)\cap[M]))$,\qquad{if} $w=1$,\quad {and}

$\sigma_4(\hat f)=I_0(KS(N)-KS(M)+kerv(\hat f)^2)
+\kappa_2(c_{M*}(kerv(\hat f)\cap[M]))$,
if $w\not=1$.

(In the latter case we identify $H^4(M;\mathbb{Z})$, 
$H^4(N;\mathbb{Z})$ and $H^4(M;\mathbb{F}_2)$ with \nl
$H_0(\pi;\mathbb{Z}^w)=Z/2Z$.)
\qed
\end{thmm}

The homomorphism $\sigma_4$ is trivial on the image of $\eta$,
but in general we do not know whether a 4-dimensional normal map 
with trivial surgery obstruction must be normally cobordant to 
a simple homotopy equivalence.
(See however \cite{[Kh07]} and \cite{[Ym07]}.) 
In our applications we shall always have a simple homotopy equivalence 
in hand.

A more serious problem is that it is not clear how to define the action 
$\omega$ in general.
We shall be able to circumvent this problem by {\it ad hoc} arguments 
in some cases.
(There is always an action on the homological structure set,
defined in terms of $\mathbb{Z}[\pi]$-homology equivalences \cite{[FQ]}.)

If we fix an isomorphism $i_{\mathbb{Z}}:\mathbb{Z}\to L_5(\mathbb{Z})$ 
we may define a function $I_\pi:\pi\to L^s_5(\pi)$ 
for any group $\pi$ by $I_\pi(g)=g_*(i_{\mathbb{Z}}(1))$,
where $g_*:\mathbb{Z}=L_5(\mathbb{Z})\to{L^s_5(\pi)}$ 
is induced by the homomorphism sending 1 in $\mathbb{Z}$ to $g$ in $\pi$. 
Then $I_{\mathbb{Z}}=i_{\mathbb{Z}}$ and $I_\pi$ is natural 
in the sense that if $f:\pi\to H$ is a homomorphism then 
$L_5(f)I_\pi=I_Hf$.
As abelianization and projection to the summands of $\mathbb{Z}^2$ 
induce an isomorphism from $L_5(\mathbb{Z}*\mathbb{Z})$ to 
$L_5(\mathbb{Z})^2$ \cite{[Ca73]},
it follows easily from naturality that $I_\pi$ 
is a homomorphism (and so factors through $\pi/\pi'$) \cite{[We83]}.
We shall extend this to the nonorientable case by defining
$I_\pi^+:\mathrm{Ker}(w)\to L_5^s(\pi;w)$ as the composite of 
$I_{{\mathrm {Ker}}(w)}$
with the homomorphism induced by inclusion.
   
\begin{theorem} 
Let $M$ be a closed $4$-manifold with 
fundamental group $\pi$ and let $w=w_1(M)$. 
Given any $\gamma\in{\mathrm {Ker}}(w)$ there is a normal cobordism from
$id_M$ to itself with surgery obstruction $I_\pi^+(\gamma)\in L_5^s(\pi,w)$.
\end{theorem}

\begin{proof}
We may assume that $\gamma$ is represented 
by a simple closed curve with a product neighbourhood $U\cong S^1 \times D^3$.
Let $P$ be the $E_8 $ manifold \cite{[FQ]} and delete the interior of a submanifold
homeomorphic to $D^3 \times [0,1]$ to obtain $P_o$.
There is a normal map $p:P_o\to D^3\times [0,1]$ ({\it rel} boundary).
The surgery obstruction for $p\times id_{S^1}$ 
in $L_5(\mathbb{Z})\cong L_4(1)$ is given by 
a codimension-1 signature \cite[\S12B]{[Wl]},
and generates $L_5(\mathbb{Z})$.
Let $Y=(M\setminus{int}U)\times [0,1]\cup P_o \times S^1 $, 
where we identify
$(\partial U)\times [0,1]=S^1 \times S^2 \times [0,1]$ with 
$S^2 \times [0,1]\times S^1 $ in $\partial P_o\times S^1 $.
Matching together $id|_{(M\setminus{int}U)\times[0,1]} $ and 
$p\times id_{S^1}$ gives a normal cobordism $Q$ from $id_M$ to itself. 
The theorem now follows by the additivity of surgery obstructions
and naturality of the homomorphisms $I_\pi^+$.
\end{proof}

In particular, if $\pi$ is in $SA$ then the image of $I_\pi^+$
acts trivially on $S_{TOP}(M)$.

\begin{cor}
Let $\lambda_* :L_5^s (\pi)\to L_5 (\mathbb{Z})^d =\mathbb{Z}^d$ 
be the homomorphism induced by a basis $\{\lambda_1 ,\dots,\lambda_d\}$ 
for $Hom(\pi,\mathbb{Z})$.
If $M$ is orientable, $f:M_1\to M$ is a simple homotopy equivalence
and $\theta\in L_5 (\mathbb{Z})^d$ there is a normal cobordism 
from $f$ to itself whose surgery obstruction in $L_5 (\pi)$ 
has image $\theta$ under $\lambda_* $.
\end{cor}
                      
\begin{proof} 
If $\{\gamma_1 ,\dots,\gamma_d\}\in\pi$ represents a ``dual basis" 
for $H_1(\pi;\mathbb{Z})$ modulo torsion (so that 
$\lambda_i (\gamma_j)=\delta_{ij}$ for $1\leq i,j\leq d$),
then $\{\lambda_*(I_\pi(\gamma_1)),\dots,\lambda_*(I_\pi(\gamma_d))\}$ 
is a basis for $L_5(\mathbb{Z})^d$.
\end{proof}

If $\pi$ is free or is a $PD_2^+$-group the homomorphism 
$\lambda_* $ is an isomorphism \cite{[Ca73]}. 
In most of the other cases of interest to us the following corollary
applies.

\begin{cor}
If $M$ is orientable and ${\mathrm {Ker}}(\lambda_*)$  is finite then 
$S^s_{TOP} (M)$ is finite.
In particular, this is so if ${\mathrm {Coker}}(\sigma_5)$ is finite.
\end{cor}

\begin{proof} The signature difference maps 
$[M;G/TOP]=H^4 (M;\mathbb{Z})\oplus H^2 (M;\mathbb{F}_2)$ 
onto $L_4 (1)=\mathbb{Z}$ and so there are only 
finitely many normal cobordism classes of simple homotopy 
equivalences $f:M_1 \to M$. 
Moreover, $\mathrm{Ker}(\lambda_*)$ is finite 
if $\sigma_5$ has finite cokernel, 
since $[SM;G/TOP]\cong\mathbb{Z}^d\oplus (Z/2Z)^d$. 
Suppose that $F:N\to M\times I$ is a normal cobordism between two 
simple homotopy equivalences $F_- =F|\partial_- N$ and $F_+ =F|\partial_+ N$.
By Theorem 6.7 there is another normal cobordism $F':N'\to M\times I$
from $F_+ $ to itself with $\lambda_* (\sigma_5(F'))=\lambda_* (-\sigma_5(F))$.
The union of these two normal cobordisms along $\partial_+ N=\partial_- N'$
is a normal cobordism from $F_-$ to $F_+$ with surgery obstruction in 
$\mathrm{Ker}(\lambda_*)$. 
If this obstruction is 0 we may obtain an $s$-cobordism $W$ by 
5-dimensional surgery (rel $\partial$).
\end{proof}

The surgery obstruction groups for a semidirect product 
$\pi\cong G\rtimes_\theta\mathbb{Z}$,
may be related to those of the (finitely presentable) normal subgroup 
$G$ by means of \cite[Theorem 12.6]{[Wl]}.
If $Wh(\pi)=Wh(G)=0$ this theorem asserts that there is an exact sequence
\begin{equation*}
\dots L_m(G,w|_G)\buildrel1-w(t)\theta_* \over\longrightarrow  
L_m(G,w|_G)\to L_m(\pi,w)\to L_{m-1}(G,w|_G)\dots,
\end{equation*}
where $t$ generates $\pi$ modulo $G$ and $\theta_*=L_m(\theta,w|_G)$.
The following result is based on \cite[Theorem 15.B.1]{[Wl]}.

\begin{theorem} 
Let $M$ be a $4$-manifold which is homotopy equivalent 
to a mapping torus $M(\theta)$,
where $\theta$ is a self-homeomorphism of an aspherical closed 
$3$-manifold $N$.
If $Wh(\pi_1(M))=Wh(\pi_1(M)\times\mathbb{Z})=0$ 
then $M$ is $s$-cobordant to $M(\theta)$ 
and $\widetilde M$ is homeomorphic to $\mathbb{R}^4$.
\end{theorem}

\begin{proof}
The surgery obstruction homomorphisms $\sigma_i^N$ 
are isomorphisms for all large $i$ \cite{[Ro11]}.
Comparison of the Mayer-Vietoris sequences for $\mathbb{L}_0$-homology 
and $L$-theory (as in \cite[Proposition 2.6]{[St84]}) 
shows that $\sigma_i^M$ and $\sigma_i^{M\times S^1}$ 
are also isomorphisms for all large $i$,
and so $S_{TOP}(M(\Theta)\times S^1)$ has just one element.
If $h:M\to{M(\Theta)}$ is a homotopy equivalence then
$h\times{id}$ is homotopic to a homeomorphism
$M\times{S^1}\cong{M(\Theta)}\times{S^1}$,
and so $M\times\mathbb{R}\cong{M(\Theta)}\times\mathbb{R}$.
This product contains $s$-cobordisms bounded by disjoint copies of
$M$ and $M(\Theta)$.

The final assertion follows from \cite[Corollary 7.3B]{[FQ]} 
since $M$ is aspherical and $\pi$ is 1-connected at $\infty$ \cite{[Ho77]}.
\end{proof} 

It remains an open question whether aspherical closed manifolds 
with isomorphic fundamental groups must be homeomorphic. 
This has been verified in higher dimensions in many cases, 
in particular under geometric assumptions \cite{[FJ]}, 
and under assumptions on the combinatorial structure 
of the group \cite{[Ca73],[St84],[NS85]}. 
We shall see that many aspherical 4-manifolds are determined 
up to $s$-cobordism by their groups.

There are more general ``Mayer-Vietoris" sequences 
which lead to calculations of the surgery obstruction groups 
for certain generalized free products and HNN extensions 
in terms of those of their building blocks \cite{[Ca73],[St87]}.

A subgroup $H$ of a group $G$ is {\it square-root closed\/} in $G$
if $g^2\in{H}$ implies $g\in{H}$, for $g\in{G}$.
A group $\pi$ is {\it square-root closed accessible\/} if it can be obtained
from the trivial group by iterated HNN extensions 
with associated subgroups square-root closed in the base group
and amalgamated products over square-root closed subgroups.
In particular, finitely generated free groups and poly-$Z$ groups are
square-root closed accessible. 
A geometric argument implies that cuspidal subgroups of 
the fundamental group $\Gamma$ of a complete hyperbolic manifold
of finite volume are maximal parabolic subgroups, 
and hence are square root closed in $\Gamma$.
If $S$ is a closed surface with $\chi(S)<0$ it may be decomposed as
the union of two subsurfaces with connected boundary and hyperbolic interior.
Therefore all $PD_2$-groups are square-root closed accessible.
      
\begin{lemma} 
Let $\pi$ be either the group of a 
finite graph of groups, all of whose vertex groups are infinite cyclic,
or a square root closed accessible group of cohomological dimension $2$.  
Then $I_\pi^+$ is an epimorphism.
If $M$ is a closed $4$-manifold with fundamental group $\pi$ 
the surgery obstruction maps $\sigma_4 (M)$ and $\sigma_5(M)$ 
are epimorphisms.
\end{lemma}

\begin{proof} 
Since $\pi$ is in {\it Cl} we have $Wh(\pi)=0$ and a comparison of 
Mayer-Vietoris sequences shows that the assembly map from 
$H_* (\pi;\mathbb{L}_0^w )$ to $L_* (\pi,w)$ is an isomorphism \cite{[Ca73],[St87]}. 
Since $c.d.\pi\leq 2$ and $H_1({\mathrm {Ker}}(w);\mathbb{Z})$ maps onto
$H_1(\pi;Z^w)$
the component of this map in degree 1 may be identified with $I_\pi^+$.
In general, the surgery obstruction maps factor through the assembly map.
Since $c.d.\pi\leq 2$ the homomorphism 
$c_{M*}:H_*(M;\mathcal{D})\to H_*(\pi;\mathcal{D})$
is onto for any local coefficient module $\mathcal{D}$, 
and so the lemma follows.
\end{proof}


The class of groups considered in this lemma includes free groups, 
$PD_2 $-groups and the groups $Z*_m$.
Note however that if $\pi$ is a $PD_2$-group $w$ 
need not be the canonical orientation character.

\section{Stabilization and $h$-cobordism}

It has long been known that many results of high dimensional differential 
topology hold for smooth 4-manifolds after stabilizing by connected sum 
with copies of $S^2\times S^2$ \cite{[CS71],[FQ80],[La79],[Qu83]}.
In particular, if $M$ and $N$ are $h$-cobordant closed smooth 4-manifolds
then $M\sharp(\sharp^kS^2\times S^2)$ is diffeomorphic to 
$N\sharp(\sharp^kS^2\times S^2)$ for some $k\geq0$.
In the spin case $w_2(M)=0$ this is an elementary consequence of the 
existence of a well-indexed handle decomposition of the $h$-cobordism 
\cite{[Wl64]}.
In \cite[Chapter VII]{[FQ]} it is shown that 5-dimensional TOP cobordisms
have handle decompositions relative to a component of their boundaries,
and so a similar result holds for $h$-cobordant closed TOP 4-manifolds.
Moreover, if $M$ is a TOP 4-manifold then $KS(M)=0$ if and only if
$M\sharp (\sharp^k S^2\times S^2)$ is smoothable for some $k\geq0$ 
\cite{[LS71]}.

These results suggest the following definition.
Two 4-manifolds $M_1 $ and $M_2 $ are {\it stably homeomorphic} 
if $M_1\sharp (\sharp^k S^2\times S^2)$ 
and $M_2\sharp (\sharp^l S^2\times S^2)$ are homeomorphic,
for some $k$, $l\geq0$.
(Thus $h$-cobordant closed 4-manifolds are stably homeomorphic.)
Clearly $\pi_1(M)$, $w_1(M)$, the orbit of $c_{M*}[M]$
in $H_4(\pi_1(M);\mathbb{Z}^{w_1(M)})$ under the action of $Out(\pi_1(M))$, 
and the parity of $\chi(M)$ are invariant under stabilization.
If $M$ is orientable $\sigma(M)$ is also invariant.

Kreck has shown that (in any dimension) classification up to
stable homeomorphism (or diffeomorphism) can be reduced to bordism theory.
There are three cases:
If $w_2(\widetilde M)\not=0$ and $w_2(\widetilde N)\not=0$
then $M$ and $N$ are stably homeomorphic if and only if for some choices of
orientations and identification of the fundamental groups
the invariants listed above agree (in an obvious manner).
If $w_2(M)=w_2(N)=0$ then $M$ and $N$ are stably homeomorphic if and
only if for some choices of orientations, Spin structures and identification 
of the fundamental group they represent the same element in
$\Omega_4^{SpinTOP}(K(\pi,1))$.
The most complicated case is when $M$ and $N$ are not Spin, but
the universal covers are Spin. 
(See \cite{[Kr99],[Te]} for expositions of Kreck's ideas, 
and see \cite{[Po13]} for an application to 4-manifolds
determined by Tietze-equivalent presentations.)

We shall not pursue this notion of stabilization further
(with one minor exception, in Chapter 14),
for it is somewhat at odds with the tenor of this book.
The manifolds studied here usually have minimal Euler characteristic,
and often are aspherical. Each of these properties disappears after
stabilization.
We may however also stabilize by cartesian product with 
the real line $\mathbb{R}$,
and there is then the following simple but satisfying result.
             
\begin{lemma} 
Closed $4$-manifolds $M$ and $N$ are $h$-cobordant 
if and only if $M\times\mathbb{R}$ and $N\times\mathbb{R}$ are homeomorphic.
\end{lemma}
                      
\begin{proof}
If $W$ is an $h$-cobordism from $M$ to $N$ 
(with fundamental group $\pi=\pi_1(W)$) 
then $W\times S^1$ is an $h$-cobordism from $M\times S^1$ to 
$N\times S^1$.
The torsion is 0 in $Wh(\pi\times\mathbb{Z})$ \cite[Theorem 23.2]{[Co]}, 
and so there is a homeomorphism from $M\times S^1$ 
to $N\times S^1$ which carries $\pi_1(M)$ to $\pi_1(N)$.
Hence $M\times\mathbb{R}\cong N\times\mathbb{R}$.
Conversely, if $M\times\mathbb{R}\cong N\times\mathbb{R}$ 
then $M\times\mathbb{R}$ contains a copy of 
$N$ disjoint from $M\times\{0\}$, 
and the region $W$ between $M\times\{0\}$ and $N$ is an $h$-cobordism.
\end{proof}

\section{Manifolds with $\pi_1$ elementary amenable and $\chi=0$}

In this section we shall show that closed manifolds satisfying the 
hypotheses of Theorem 3.17 and with torsion-free fundamental group 
are determined up to homeomorphism by their homotopy type. 
As a consequence, closed 4-manifolds with torsion-free elementary 
amenable fundamental group and Euler characteristic 0 are homeomorphic to 
mapping tori.
We also estimate the structure sets for $RP^2$-bundles over $T$ or $Kb$.
In the remaining cases involving torsion
computation of the surgery obstructions is much more difficult.
We shall comment briefly on these cases in Chapters 10 and 11.

\begin{theorem} 
Let $M$ be a closed $4$-manifold with
$\chi(M)=0$ and whose fundamental group $\pi$ is torsion-free, coherent, 
locally virtually indicable and restrained.
Then $M$ is determined up to homeomorphism by its homotopy type.
If, moreover, $h(\pi)=4$ then every automorphism of $\pi$ is realized by a 
self homeomorphism of $M$.
\end{theorem}
                       
\begin{proof} 
By Theorem 3.17 either $\pi\cong\mathbb{Z}$ or $Z*_m$ for some $m\not=0$, or
$M$ is aspherical, $\pi$ is virtually poly-$Z$ and $h(\pi)=4$.
Hence $Wh(\pi)=0$, in all cases.
If $\pi\cong\mathbb{Z}$ or $Z*_m$ then the surgery obstruction homomorphisms 
are epimorphisms, by Lemma 6.9. 
We may calculate $L_4(\pi,w)$ by means of \cite[Theorem 12.6]{[Wl]}, 
or more generally \cite[\S3]{[St87]}, 
and we find that if $\pi\cong\mathbb{Z}$ or $Z*_{2n}$ then 
$\sigma_4(M)$ is in fact an isomorphism. 
If $\pi\cong Z*_{2n+1}$ then there are two normal cobordism classes of 
homotopy equivalences $h:X\to M$. 
Let $\xi$ generate the image of $H^2(\pi;\mathbb{F}_2)\cong Z/2Z$ in 
$H^2(M;\mathbb{F}_2)\cong (Z/2Z)^2$,
and let $j:S^2\to M$ represent the unique nontrivial spherical class in 
$H_2(M;\mathbb{F}_2)$.
Then $\xi^2=0$, since $c.d.\pi=2$, and $\xi\cap j_*[S^2]=0$, 
since $c_Mj$ is nullhomotopic.
It follows that $j_*[S^2]$ is Poincar\'e dual to $\xi$, and so 
$v_2(M)\cap j_*[S^2]=\xi^2\cap[M]=0$.
Hence $j^*w_2(M)=j^*v_2(M)+(j^*w_1(M))^2=0$ and so $f_j$ has nontrivial 
normal invariant, by Lemma 6.5.
Therefore each of these two normal cobordism classes contains a self 
homotopy equivalence of $M$.

If $M$ is aspherical, $\pi$ is virtually poly-$Z$ and $h(\pi)=4$ 
then $S_{TOP}(M)$ has just one element \cite[Theorem 2.16]{[FJ]}.
The theorem now follows.
\end{proof}

\begin{cor}
Let $M$ be a closed $4$-manifold with $\chi(M)=0$ and 
fundamental group $\pi\cong\mathbb{Z}$, $\mathbb{Z}^2$ or 
$\mathbb{Z}\rtimes_{-1}\!\mathbb{Z}$. 
Then $M$ is determined up to homeomorphism by $\pi$ and $w(M)$.
\end{cor}
                       
\begin{proof} If $\pi\cong\mathbb{Z}$ then $M$ is homotopy equivalent to 
$S^1\times{S^3}$ or $S^1\tilde\times{S^3}$, by Corollary 4.5.3,
while if $\pi\cong\mathbb{Z}^2$ or $\mathbb{Z}\rtimes_{-1}\!\mathbb{Z}$ 
it is homotopy equivalent to the total space of an $S^2$-bundle 
over $T$ or $Kb$,
by Theorem 5.10.
\end{proof}

Closed orientable 4-manifolds $M$ with $\chi(M)=0$ and $\pi\cong{Z*_m}$
are also determined up to homeomorphism by $\pi$ and $w(M)$  
\cite{[Hi09],[HKT09]}.

We may now give an analogue of the Farrell and Stallings fibration theorems for
4-manifolds with torsion-free elementary amenable fundamental group.
                                                        
\begin{theorem} 
Let $M$ be a closed $4$-manifold whose
fundamental group $\pi$ is torsion-free and elementary amenable. 
A map $f:M\to S^1$ is homotopic to a fibre bundle projection if and only if 
$\chi(M)=0$ and $f$ induces an epimorphism from $\pi$ to $\mathbb{Z}$ with  
finitely generated kernel.
\end{theorem}

\begin{proof} The conditions are clearly necessary. 
Suppose that they hold.
Let $\nu=\mathrm{Ker}(\pi_1 (f))$, 
let $M_\nu$ be the infinite cyclic covering space of $M$
with fundamental group $\nu$ and 
let $t:M_\nu\to M_\nu$ be a generator of the group of covering transformations.
By Corollary 4.5.2 either $\nu=1$ (so $M_\nu\simeq S^3 $)
or $\nu\cong\mathbb{Z}$ (so $M_\nu\simeq S^2 \times S^1 $ 
or $S^2 \tilde\times S^1 $) or $M$ is aspherical.
In the latter case $\pi$ is a torsion-free virtually poly-$Z$ group, 
by Theorems 1.11 and 9.23 of \cite{[Bi]}. 
Thus in all cases there is a homotopy equivalence $f_\nu$ from $M_\nu$
to a closed 3-manifold $N$. 
Moreover the self homotopy equivalence $f_\nu tf_\nu^{-1} $ 
of $N$ is homotopic to a homeomorphism, $g$ say, 
and so $f$ is fibre homotopy equivalent to
the canonical projection of the mapping torus $M(g)$ onto $S^1 $.
It now follows from Theorem 6.11 that any homotopy equivalence 
from $M$ to $M(g)$ is homotopic to a homeomorphism.                                                    
\end{proof}

The structure sets of the $RP^2$-bundles over $T$ or $Kb$ are also finite. 

\begin{theorem} 
Let $M$ be the total space of an $RP^2$-bundle
over $T$ or $Kb$. Then $S_{TOP} (M)$ has order at most $32$.
\end{theorem} 

\begin{proof}                    
As $M$ is nonorientable $H^4(M;\mathbb{Z})=Z/2Z$ and as $\beta_1(M;\mathbb{F}_2)=3$
and $\chi(M)=0$ we have $H^2(M;\mathbb{F}_2)\cong(Z/2Z)^4 $. 
Hence $[M;G/TOP]$ has order 32.
Let $w=w_1(M)$.
It follows from the Shaneson-Wall splitting theorem 
\cite[Theorem 12.6]{[Wl]} that 
$L_4(\pi,w)\cong L_4 (Z/2Z,-)\oplus L_2 (Z/2Z,-)\cong(Z/2Z)^2$,
detected by the Kervaire-Arf invariant and the codimension-2 Kervaire invariant.
Similarly $L_5(\pi,w)\cong L_4 (Z/2Z,-)^2$ and the projections
to the factors are Kervaire-Arf invariants of normal maps induced over
codimension-1 submanifolds.
(In applying the splitting theorem, 
note that $Wh(\mathbb{Z}\oplus(Z/2Z))=Wh(\pi)=0$, 
by Theorem 6.1 above.)
Hence $S_{TOP}(M)$ has order at most 128.

The Kervaire-Arf homomorphism $c$ is onto,
since $c(\hat g)=(w^2 \cup\hat g^*(k_2))\cap [M]$,
$w^2\not=0$ and every element of 
$H^2 (M;\mathbb{F}_2)$ is equal to $\hat g^*(k_2)$ for some 
normal map $\hat g:M\to G/TOP$.
Similarly there is a normal map $f_2 :X_2 \to RP^2 $
with $\sigma_2(f_2 )\not= 0$ in $L_2 (Z/2Z,-)$. 
If $M=RP^2 \times B$, where $B=T$ or $Kb$ is the base of the bundle, 
then $f_2 \times id_B :X_2 \times B\to RP^2 \times B$ 
is a normal map with surgery obstruction 
$(0,\sigma_2(f_2 ))\in L_4 (Z/2Z,-)\oplus L_2 (Z/2Z,-)$.
We may assume that $f_2$ is a homeomorphism over a disc $\Delta\subset RP^2 $. 
As the nontrivial bundles may be obtained from the product bundles by cutting 
$M$ along $RP^2 \times \partial\Delta $ and regluing via the twist map of 
$RP^2 \times S^1 $, the normal maps for the product bundles 
may be compatibly modified to give normal maps with nonzero obstructions
in the other cases. 
Hence $\sigma_4$ is onto and so $S_{TOP} (M)$ has order at most 32. 
\end{proof}

In each case $H_2(M;\mathbb{F}_2)\cong H_2(\pi;\mathbb{F}_2)$,
so the argument of Lemma 6.5 does not apply. 
However we can improve our estimate in the abelian case.                                         

\begin{theorem} 
Let $M$ be the total space of an $RP^2$-bundle over $T$.
Then $S_{TOP}(M)$ has order 8.
\end{theorem} 

\begin{proof} 
Since $\pi$ is abelian the surgery sequence may be identified 
with the algebraic surgery sequence of \cite{[Rn]}, 
which is an exact sequence of abelian groups.
Thus it shall suffice to show that $L_5(\pi,w)$ acts trivially 
on the class of $id_M$ in $S_{TOP}(M)$.

Let $\lambda_1 ,\lambda_2:\pi\to\mathbb{Z}$ be epimorphisms generating 
$Hom(\pi,\mathbb{Z})$ and let $t_1, t_2\in\pi$ represent a dual basis for
$\pi/(torsion)$ (i.e., $\lambda_i (t_j)=\delta_{ij}$ for $i=1,2$).
Let $u$ be the element of order 2 in $\pi$ and let 
$k_i:\mathbb{Z}\oplus(Z/2Z)\to\pi$ 
be the monomorphism defined by $k_i(a,b)=at_i+bu$, for $i=1,2$.
Define splitting homomorphisms $p_1,p_2$ by 
$p_i (g)=k_i^{-1}(g-\lambda_i (g)t_i)$ for all $g\in\pi$.
Then $p_ik_i=id_{\mathbb{Z}\oplus(Z/2Z)}$ and $p_i k_{3-i} $ 
factors through $Z/2Z$, 
for $i=1,2$.
The orientation character $w=w_1 (M)$ maps the torsion subgroup of $\pi$ 
onto $Z/2Z$, by Theorem 5.13, and $t_1$ and $t_2$ are in $\mathrm{Ker}(w)$.
Therefore $p_i$ and $k_i$ are compatible with $w$, for $i=1,2$.
As $L_5 (Z/2Z,-)=0$ it follows that $L_5 (k_1)$ and $L_5 (k_2)$ are inclusions 
of complementary summands of $L_5 (\pi,w)\cong(Z/2Z)^2 $,
split by the projections $L_5 (p_1 )$ and $L_5 (p_2 )$.

Let $\gamma_i$ be a simple closed curve in $T$ which represents $t_i \in\pi$.
Then $\gamma_i$ has a product neighbourhood $N_i\cong S^1 \times [-1,1]$ 
whose preimage $U_i\subset M$ is homeomorphic to $RP^2\times S^1\times [-1,1]$.
As in Theorem 6.13 there is a normal map $f_4:X_4 \to RP^2 \times [-1,1]^2 $ ({\it rel} boundary) 
with $\sigma_4(f_4 )\not= 0$ in $L_4 (Z/2Z,-)$.
Let $Y_i=(M\setminus{int}U_i)\times [-1,1]\cup X_4 \times S^1 $, where we 
identify $(\partial U_i)\times [-1,1]=RP^2 \times S^1\times S^0\times [-1,1]$ 
with $RP^2 \times [-1,1]\times S^0 \times S^1 $ in $\partial X_4\times S^1 $.
If we match together $id_{(M\setminus{int}U_i)\times [-1,1]} $ and 
$f_4 \times id_{S^1} $ we obtain a normal cobordism $Q_i$ from $id_M$ to itself. 
The image of $\sigma_5(Q_i)$ in 
$L_4 (\mathrm{Ker}(\lambda_i),w)\cong L_4 (Z/2Z,-)$ 
under the splitting homomorphism is $\sigma_4(f_4)$.
On the other hand its image in $L_4 (\mathrm{Ker}(\lambda_{3-i}),w)$ is 0, 
and so it generates the image of $L_5 (k_{3-i} )$. 
Thus $L_5 (\pi, w)$ is generated by $\sigma_5 (Q_1)$ and 
$\sigma_5(Q_2)$, and so acts trivially on $id_M$.
\end{proof}

Does $L_5 (\pi,w)$ act trivially on each class in $S_{TOP} (M)$
when $M$ is an $RP^2 $-bundle over $Kb$?
If so, then $S_{TOP} (M)$ has order 8 in each case. 
Are these manifolds determined up to homeomorphism by their homotopy type?

\section{Bundles over aspherical surfaces}

The fundamental groups of total spaces of bundles over 
hyperbolic surfaces all contain nonabelian free subgroups.
Nevertheless, such bundle spaces are determined up to
$s$-cobordism by their homotopy type, except when the fibre is $RP^2$,
in which case we can only show that the structure sets are finite.

\begin{theorem} 
Let $M$ be a closed $4$-manifold which is homotopy equivalent to 
the total space $E$ of an 
$F$-bundle over $B$ where $B$ and $F$ are aspherical closed surfaces.
Then $M$ is $s$-cobordant to $E$ and $\widetilde M$ 
is homeomorphic to $\mathbb{R}^4$.
\end{theorem} 
                      
\begin{proof} 
If $\chi(B)=0$ then $\pi\times\mathbb{Z}$ is an extension of a poly-$Z$ group
(of Hirsch length 3) by $\pi_1(F)$.
Otherwise, $\pi_1(B)\cong{F*_{\mathbb{Z}}F'}$, 
where the amalgamated subgroup $\mathbb{Z}$ 
is square-root closed in each of the free groups $F$ and $F'$.
(See the final paragraph on page 120.)
In all cases $\pi\times\mathbb{Z}$ is a square root closed 
generalised free product with amalgamation of groups in {\it Cl}.
Comparison of the Mayer-Vietoris sequences 
for $\mathbb{L}_0 $-homology and $L$-theory 
(as in \cite[Proposition 2.6]{[St84]})
shows that $S_{TOP}(E\times S^1)$ has just one element.
(Note that even when $\chi(B)=0$ the groups arising in intermediate stages of 
the argument all have trivial Whitehead groups.)
Hence $M\times S^1\cong E\times S^1$, and so $M$ is $s$-cobordant to $E$
by Lemma 6.10 and Theorem 6.2.

The final assertion follows from \cite[Corollary 7.3B]{[FQ]}, 
since $M$ is aspherical and $\pi$ is 1-connected at $\infty$ \cite{[Ho77]}.
\end{proof}

Davis has constructed aspherical 4-manifolds whose universal covering space
is not 1-connected at $\infty$ \cite{[Da83]}.

\begin{theorem} 
Let $M$ be a closed $4$-manifold which is homotopy equivalent to 
the total space $E$ of an $S^2 $-bundle over an aspherical closed surface $B$. 
Then $M$ is $s$-cobordant to $E$, and
$\widetilde M$ is homeomorphic to $S^2 \times\mathbb{R}^2$.
\end{theorem}

\begin{proof} 
Let $\pi=\pi_1 (E)\cong\pi_1 (B)$.
Then $Wh(\pi)=0$, and $H_* (\pi;\mathbb{L}_0^w )\cong L_* (\pi,w)$, 
as in Lemma 6.9.
Hence $L_4(\pi,w)\cong\mathbb{Z}\oplus (Z/2Z)$ if $w=0$ 
and $(Z/2Z)^2$ otherwise. 
The surgery obstruction map $\sigma_4 (E)$ is onto, by Lemma 6.9.
Hence there are two normal cobordism classes of maps
$h:X\to E$ with $\sigma_4(h)=0$. 
The kernel of the natural homomorphism 
from $H_2 (E;\mathbb{F}_2)\cong (Z/2Z)^2 $ to $H_2 (\pi;\mathbb{F}_2)\cong Z/2Z$
is generated by $j_* [S^2 ]$, where $j:S^2 \to E$ is the inclusion of a fibre.
As $j_* [S^2 ]\not=-0$, while $w_2 (E)(j_*[S^2])=j^*w_2 (E)=0$
the normal invariant of $f_j$ is nontrivial, by Lemma 6.5. 
Hence each of these two normal cobordism classes contains 
a self homotopy equivalence of $E$.

Let $f:M\to E $ be a homotopy equivalence (necessarily simple).
Then there is a normal cobordism $F:V\to E\times [0,1]$ 
from $f$ to some self homotopy equivalence of $E$. 
As $I_\pi^+$ is an isomorphism, by Lemma 6.9,
there is an $s$-cobordism $W$ from $M$ to $E$,
as in Corollary 6.7.2.

The universal covering space $\widetilde W$ is a proper $s$-cobordism from 
$\widetilde M$ to $\widetilde E\cong S^2 \times\mathbb{R}^2$. 
Since the end of $\widetilde E$ is tame and has fundamental group 
$\mathbb{Z}$ we may apply \cite[Corollary 7.3B]{[FQ]} to conclude that 
$\widetilde W$ is homeomorphic to a product. 
Hence $\widetilde M$ is homeomorphic to $S^2 \times\mathbb{R}^2$. 
\end{proof}

Let $\rho$ be a $PD_2$-group.
As $\pi=\rho\times(Z/2Z)$ is square-root closed accessible from $Z/2Z$,
the Mayer-Vietoris sequences of \cite{[Ca73]} imply that $L_4(\pi,w)\cong
L_4(Z/2Z,-)\oplus L_2(Z/2Z,-)$ and that $L_5(\pi,w)\cong L_4(Z/2Z,-)^\beta$,
where $w=pr_2:\pi\to Z/2Z$ and $\beta=\beta_1(\rho;\mathbb{F}_2)$.
Since these $L$-groups are finite the structure sets of total spaces 
of $RP^2$-bundles over aspherical surfaces are also finite.
(Moreover the arguments of Theorems 6.13 and 6.14 can be extended to show that
$\sigma_4$ is an epimorphism and that most of $L_5(\pi,w)$ acts trivially
on $id_E$, where $E$ is such a bundle space.)

%% file: m5-7.tex
\part{4-dimensional Geometries}

\chapter{Geometries and decompositions}

Every closed connected surface is geometric, i.e., 
is a quotient of one of the 
three model 2-dimensional geometries $\mathbb{E}^2$, 
$\mathbb{H}^2$ or $\mathbb{S}^2$ by a free and properly 
discontinuous action of a discrete group of isometries. 
Every closed irreducible 3-manifold admits a finite decomposition 
into geometric pieces.
(That this should be so was the Geometrization Conjecture of Thurston,
which was proven in 2003 by Perelman, 
through an analysis of the Ricci flow, introduced by Hamilton.)
In \S1 we shall recall Thurston's definition of geometry, 
and shall describe briefly the 19 4-dimensional geometries.
Our concern in the middle third of this book is not 
to show how this list arises
(as this is properly a question of differential geometry; see 
\cite{[Is55],[Fi],[Pa96]} and \cite{[Wl85],[Wl86]}), 
but rather to describe the geometries sufficiently well 
that we may subsequently characterize geometric manifolds 
up to homotopy equivalence or homeomorphism. 
In \S2 and \S3 we relate the notions of ``geometry of solvable Lie type" 
and ``infrasolvmanifold".
The limitations of geometry in higher dimensions are illustrated in \S4, 
where it is shown that a closed 4-manifold which admits a finite decomposition 
into geometric pieces is (essentially) either geometric or aspherical.           The geometric viewpoint is nevertheless of considerable interest 
in connection with complex surfaces
\cite{[Ue90],[Ue91],[Wl85],[Wl86]}.                                              With the exception of the geometries $\mathbb{S}^2\times\mathbb{H}^2$, 
$\mathbb{H}^2\times\mathbb{H}^2$,
$\mathbb{H}^2\times\mathbb{E}^2$ and 
$\widetilde{\mathbb{SL}}\times\mathbb{E}^1$ 
no closed geometric manifold has a proper geometric decomposition.                                                         
A number of the geometries support natural Seifert fibrations or compatible 
complex structures.
In \S4 we characterize the groups of aspherical 4-manifolds which are
orbifold bundles over flat or hyperbolic 2-orbifolds.
We outline what we need about Seifert fibrations and complex surfaces
in \S5 and \S6.

Subsequent chapters shall consider in turn geometries 
whose models are contractible (Chapters 8 and 9), 
geometries with models diffeomorphic to $S^2\times\mathbb{R}^2$ (Chapter 10),
the geometry $\mathbb{S}^3\times\mathbb{E}^1$ (Chapter 11) 
and the geometries with compact models (Chapter 12).
In Chapter 13 we shall consider geometric structures and decompositions of bundle spaces.
In the final chapter of the book we shall consider knot manifolds which admit geometries.

\vfil\eject
\section{Geometries}

An $n${\it -dimensional geometry} $\mathbb{X}$ in the sense of Thurston is 
represented by a pair $(X,G_X)$ where $X$ is a complete 1-connected 
$n$-dimensional Riemannian manifold and $G_X$ is a Lie group which acts 
effectively, transitively and isometrically on $X$ and which has discrete 
subgroups $\Gamma$ which act freely on $X$ so that $\Gamma\backslash X$ 
has finite volume. 
(Such subgroups are called {\it lattices}.)
Since the stabilizer of a point in $X$ is isomorphic to a closed subgroup of 
$O(n)$ it is compact, and so $\Gamma\backslash X$ is compact if and only if 
$\Gamma\backslash G_X$ is compact.
Two such pairs $(X,G)$ and $(X',G')$ define the same geometry if there is a
diffeomorphism $f:X\to X'$ which conjugates the action of $G$ onto that of $G'$.
(Thus the metric is  only an adjunct to the definition.)
We shall assume that $G$ is maximal among Lie groups acting thus on $X$, 
and write $Isom(\mathbb{X})=G$, 
and $Isom_o(\mathbb{X})$ for the component of the identity.
A closed manifold $M$ is an $\mathbb{X}${\it -manifold} 
if it is a quotient $\Gamma\backslash X$ for some lattice in $G_X$.
Under an equivalent formulation, $M$ is an $\mathbb{X}$-manifold
if it is a quotient $\Gamma\backslash X$ for some discrete group $\Gamma$ 
of isometries acting freely on a 1-connected homogeneous space $X=G/K$, 
where $G$ is a connected Lie group and $K$ is a compact subgroup of $G$ such 
that the intersection of the conjugates of $K$ is trivial, 
and $X$ has a $G$-invariant metric. 
The manifold {\it admits a geometry of type} $\mathbb{X}$ if it 
is homeomorphic to such a quotient.
If $G$ is solvable we shall say that the geometry is of {\it solvable Lie type}.
If $\mathbb{X}=(X,G_X)$ and $\mathbb{Y}=(Y,G_Y)$ are two geometries then 
$X\times Y$ supports a geometry in a natural way; 
$G_{X\times Y}=G_X\times G_Y$ if $X$ and $Y$ are irreducible 
and not isomorphic,
but otherwise $G_{X\times Y}$ may be larger.

The geometries of dimension 1 or 2 are the familiar geometries of constant 
curvature: $\mathbb{E}^1$, $\mathbb{E}^2$, $\mathbb{H}^2$ and $\mathbb{S}^2$. 
Thurston showed that there are eight maximal 3-dimensional geometries 
($\mathbb{E}^3$, $\mathbb{N}il^3$, $\mathbb{H}^2\times\mathbb{E}^1$,
$\widetilde{\mathbb{SL}}$, $\mathbb{S}ol^3$, $\mathbb{H}^3$, 
$\mathbb{S}^2\times\mathbb{E}^1$ and $\mathbb{S}^3$.)
Manifolds with one of the first five of these geometries are
aspherical Seifert fibred 3-manifolds or $\mathbb{S} ol^3$-manifolds. 
These are determined 
by their fundamental groups,
which are the $PD_3 $-groups with nontrivial Hirsch-Plotkin radical.
A closed 3-manifold $M$ is hyperbolic if and only if it is aspherical and
$\pi_1(M)$ has no rank 2 abelian subgroup, while it is an
$\mathbb{S}^3$-manifold if and only if $\pi_1(M)$ is finite,
by the work of Perelman \cite{[B-P]}.
(See \S11.4 below for more on $\mathbb{S}^3$-manifolds and their groups.)
There are just four $\mathbb{S}^2\times\mathbb{E}^1 $-manifolds. 
For a detailed and lucid account of the 3-dimensional geometries see 
see \cite{[Sc83']}.

There are 19 maximal 4-dimensional geometries; 
one of these ($\mathbb{S}ol^4_{m,n}$) 
is in fact a family of closely related geometries, 
and one ($\mathbb{F}^4$) is not realizable by any closed manifold \cite{[Fi]}. 
We shall see that the geometry is determined by the fundamental group 
(cf. \cite{[Wl86],[Ko92]}).
In addition to the geometries of constant curvature and products of lower 
dimensional geometries there are seven ``new" 4-dimensional geometries. 
Two of these have models the irreducible Riemannian symmetric spaces 
$CP^2=U(3)/U(2)$ and $H^2(C)=SU(2,1)/S(U(2)\times U(1))$. 
The model for the geometry $\mathbb{F}^4$ is 
the tangent bundle of the hyperbolic plane, 
which we may identify with $\mathbb{R}^2\times H^2$.
Its isometry group is the semidirect product 
$\mathbb{R}^2\times_\alpha SL^{\pm}(2,\mathbb{R})$, 
where $SL^{\pm}(2,\mathbb{R})=\{A\in{GL(2,\mathbb{R})}\mid\det{A}=\pm1\}$,
and $\alpha$ is the natural action of $SL^{\pm}(2,\mathbb{R})$ 
on $\mathbb{R}^2$.
The identity component acts on $\mathbb{R}^2\times{H}^2$ as follows:
if $u\in\mathbb{R}^2$ and 
$A=\left(\smallmatrix a&b\\
c&d\endsmallmatrix\right)\in SL(2,\mathbb{R})$
then $u(w,z)=(u+w,z)$ and 
$A(w,z)=(Aw,\frac{az+b}{cz+d})$
for all $(w,z)\in\mathbb{R}^2\times{H}^2$. 
The matrix $D=\left(\smallmatrix 1&0\\
0&-1\endsmallmatrix\right)$ 
acts via $D(w,z)=(Dw,-\bar{z})$.
All $\mathbb{H}^2(\mathbb{C})$- and $\mathbb{F}^4$-manifolds are orientable.
The other four new geometries are of solvable Lie type, 
and shall be described in \S2 and \S3.

In most cases the model $X$ is homeomorphic to $\mathbb{R}^4$, 
and the corresponding geometric manifolds are aspherical.
Six of these geometries ($\mathbb{E}^4$, $\mathbb{N}il^4$, 
$\mathbb{N}il^3\times\mathbb{E}^1$, $\mathbb{S} ol^4_{m,n}$, 
$\mathbb{S}ol^4_0$ and $\mathbb{S}ol^4_1$) are of solvable Lie type; 
in Chapter 8 we shall show manifolds admitting such geometries
have Euler characteristic 0 and fundamental group
a torsion-free virtually poly-$Z$ group of Hirsch length 4.
Such manifolds are determined up to homeomorphism by their fundamental groups,
and every such group arises in this way.
In Chapter 9 we shall consider closed 4-manifolds admitting one of the other geometries 
of aspherical type ($\mathbb{H}^3\times\mathbb{E}^1$, 
$\widetilde{\mathbb{SL}}\times\mathbb{E}^1$, 
$\mathbb{H}^2\times\mathbb{E}^2$, $\mathbb{H}^2\times\mathbb{H}^2$, 
$\mathbb{H}^4$, $\mathbb{H}^2(\mathbb{C})$ and $\mathbb{F}^4$).
These may be characterised up to $s$-cobordism by their fundamental group 
and Euler characteristic.
However it is unknown to what extent surgery arguments apply in these cases,
and we do not yet have good characterizations of the possible 
fundamental groups.                         
Although no closed 4-manifold admits the geometry $\mathbb{F}^4$, 
there are such manifolds with
proper geometric decompositions involving this geometry; 
we shall give examples in Chapter 13.

Three of the remaining geometries ($\mathbb{S}^2\times\mathbb{E}^2$, 
$\mathbb{S}^2\times\mathbb{H}^2$ and $\mathbb{S}^3\times\mathbb{E}^1$) 
have models homeomorphic to $S^2\times\mathbb{R}^2$ or $S^3\times\mathbb{R}$.
The final three ($\mathbb{S}^4$, $\mathbb{CP}^2$ and $\mathbb{S}^2\times\mathbb{S}^2$) 
have compact models, and there are only eleven such manifolds.
We shall discuss these nonaspherical geometries in Chapters 10, 11 and 12.

\section{Infranilmanifolds}

The notions of ``geometry of solvable Lie type" and ``infrasolvmanifold" 
are closely related.
We shall describe briefly the latter class of manifolds, from a 
rather utilitarian point of view. 
As we are only interested in closed manifolds, 
we shall frame our definitions accordingly.
We consider the easier case of infranilmanifolds in this section,
and the other infrasolvmanifolds in \S3.

A {\it flat $n$-manifold} is a quotient of $\mathbb{R}^n$ 
by a discrete torsion-free subgroup 
of $E(n)=Isom(\mathbb{E}^n)=\mathbb{R}^n\rtimes_\alpha O(n)$ 
(where $\alpha$ is the natural action of $O(n)$ on $\mathbb{R}^n$).
A group $\pi$ is a {\it flat $n$-manifold group} if it is torsion-free and
has a normal subgroup of finite index which is isomorphic to $\mathbb{Z}^n$.
(These are necessary and sufficient conditions for $\pi$ to be the fundamental 
group of a closed flat $n$-manifold.)
The {\it translation subgroup} $T(\pi)=\pi\cap\mathbb{R}^n$
is the maximal abelian normal subgroup of $\pi$. 
Conjugation in $\pi$ induces a faithful action of $\pi/T(\pi)$ on $T(\pi)$. 
On choosing an isomorphism $T(\pi)\cong\mathbb{Z}^n$ 
we may identify $\pi/T(\pi)$ with a subgroup of $GL(n,\mathbb{Z})$;
this subgroup is called the {\it holonomy group} of $\pi$, 
and is well defined up to conjugacy in $GL(n,\mathbb{Z})$. 
We say that $\pi$ is {\it orientable} if the holonomy group lies in 
$SL(n,\mathbb{Z})$; equivalently, 
$\pi$ is orientable if the flat $n$-manifold $\mathbb{R}^n/\pi$ 
is orientable or if $\pi\leq{E(n)^+}=\mathbb{R}^n\rtimes_\alpha{SO(n)}$.
If two discrete torsion-free cocompact subgroups of $E(n)$ are isomorphic
then they are conjugate in the larger group 
$Aff(\mathbb{R}^n)=\mathbb{R}^n\rtimes_\alpha GL(n,\mathbb{R})$,
and the corresponding flat $n$-manifolds are ``affinely" diffeomorphic.
There are only finitely many isomorphism classes of such flat $n$-manifold 
groups for each $n$.

A nilmanifold is a coset space of a 1-connected nilpotent Lie group by a discrete subgroup.
More generally, an {\it infranilmanifold} is a quotient $\pi\backslash N$ 
where $N$ is a 1-connected nilpotent Lie group and $\pi$ is a discrete 
torsion-free subgroup of the semidirect product 
$Aff(N)=N\rtimes_\alpha Aut(N)$ such 
that the translation subgroup $T(\pi)=\pi\cap N$ is a lattice in $N$ 
and $\pi/\pi\cap N$ is finite.      
Thus infranilmanifolds are finitely covered by nilmanifolds.
It follows from the work of Mal'cev that $T(\pi)=\sqrt\pi$ and that
the Lie group $N$ is determined by $\sqrt\pi$ \cite{[Ma49]}.
Two infranilmanifolds are diffeomorphic if and only if 
their fundamental groups  are isomorphic. 
The isomorphism may then be induced by an affine diffeomorphism.
The infranilmanifolds derived from the abelian Lie groups $\mathbb{R}^n$ 
are just the flat manifolds.
It is not hard to see that there are just three 
4-dimensional (real) nilpotent Lie algebras.
(Compare the analogous argument of Theorem 1.4.)
Hence there are three 1-connected 4-dimensional nilpotent Lie groups,
$\mathbb{R}^4$, $Nil^3\times\mathbb{R}$ and $Nil^4$. 

The group $Nil^3$ is the subgroup of $SL(3,\mathbb{R})$ 
consisting of upper triangular matrices
$[r,s,t]=
\begin{pmatrix}
1 & r & t\\
0 & 1 & s\\
0 & 0 & 1
\end{pmatrix}.$                
It has abelianization $\mathbb{R}^2$ and centre 
$\zeta Nil^3={Nil^3}'\cong\mathbb{R}$.
The elements $[1,0,0]$, $[0,1,0]$ and $[0,0,1/q]$ generate 
a discrete cocompact subgroup of $Nil^3$ isomorphic to $\Gamma_q$, 
and these are essentially the only such subgroups.
(Since they act orientably on $\mathbb{R}^3$ they are $PD_3^+$-groups.)
The coset space $N_q=Nil^3/\Gamma_q$ is the total space of the $S^1$-bundle 
over $S^1\times S^1$ with Euler number $q$, 
and the action of $\zeta Nil^3$ on $Nil^3$ induces a free action
of $S^1=\zeta Nil/\zeta\Gamma_q$ on $N_q$.
The group $Nil^4$ is the semidirect product 
$\mathbb{R}^3\rtimes_\theta\mathbb{R}$, 
where $\theta(t)=[t,t,t^2/2]$.                                             
It has abelianization $\mathbb{R}^2$ and central series 
$\zeta{Nil^4}\cong\mathbb{R}< \zeta_2{Nil^4}={Nil^4}'\cong\mathbb{R}^2$.

These Lie groups have natural left invariant metrics,
and the isometry groups are generated by left translations 
and the stabilizer of the identity.
For $\mathbb{N}il^3$ this stabilizer is $O(2)$,
and $Isom(\mathbb{N}il^3)$ is an extension of $E(2)$ by $\mathbb{R}$.
Hence $Isom(\mathbb{N}il^3\times\mathbb{E}^1)=Isom(\mathbb{N}il^3)\times{E(1)}$.
For $\mathbb{N}il^4$ the stabilizer is $(Z/2Z)^2$,
and is generated by two involutions, 
which send $((x,y,z),t)$ to $(-(x,y,z),t)$ and $((-x,y,z),-t)$, respectively.
(See \cite{[Sc83'],[Wl86]}.)

\section{Infrasolvmanifolds}

The situation for (infra)solvmanifolds is more complicated.
An {\it infrasolvmanifold} is a quotient $M=\Gamma\backslash S$ 
where $S$ is a 1-connected solvable Lie group and 
$\Gamma$ is a closed torsion-free subgroup of the semidirect product 
$Aff(S)=S\rtimes_\alpha Aut(S)$ such that $\Gamma_o$ 
(the component of the identity of $\Gamma$) is contained 
in the {\it nilradical} of $S$ 
(the maximal connected nilpotent normal subgroup of $S$),
$\Gamma/\Gamma\cap S$ has compact closure in $Aut(S)$ and $M$ is compact.
The pair $(S,\Gamma)$ is called a presentation for $M$, 
and is discrete if $\Gamma$ is a discrete subgroup of $Aff(S)$,
in which case $\pi_1(M)=\Gamma$.
Every infrasolvmanifold has a presentation such that 
$\Gamma/\Gamma\cap S$ is finite \cite{[FJ97]}, but
$\Gamma$ need not be discrete, and $S$ is not determined by $\pi$. 
(For example, $\mathbb{Z}^3$ is a lattice in both $\mathbb{R}^3$ and
$\widetilde{E(2)^+}=\mathbb{C}\rtimes_{\tilde\alpha}\mathbb{R}$,
where $\tilde\alpha(t)(z)=e^{2\pi it}z$ for all $t\in\mathbb{R}$ and
$z\in\mathbb{C}$.)

Since $S$ and $\Gamma_o$ are each contractible, 
$X=\Gamma_o\backslash S$ is contractible also.                            
It can be shown that $\pi=\Gamma/\Gamma_o$ acts freely on $X$, 
and so is a $PD_m$ group, 
where $m$ is the dimension of $M=\pi\backslash X$. 
(See \cite[Chapter III.3]{[Au73]} for the solvmanifold case.)
Since $\pi$ is also virtually solvable it is thus virtually poly-$Z$ 
of Hirsch length $m$ \cite[Theorem 9.23]{[Bi]}, 
and $\chi(M)=\chi(\pi)=0$.

Working in the context of real algebraic groups, 
Baues has shown in \cite{[Ba04]} that
\begin{enumerate}
\item every infrasolvmanifold has a discrete presentation with
$\Gamma/\Gamma\cap S$ finite; 

\item infrasolvmanifolds with isomorphic fundamental groups are diffeomorphic.
\end{enumerate}
He has also given a new construction which realizes each torsion-free 
virtually poly-$Z$ group as the fundamental group of an infrasolvmanifold,
a result originally due to Auslander and Johnson \cite{[AJ76]}.
Farrell and Jones had shown earlier that (2) holds
in all dimensions {\it except} perhaps 4.
However there is not always an affine diffeomorphism \cite{[FJ97]}.
(Theorem 8.9 below gives an {\it ad hoc\/} argument,
using the Mostow orbifold bundle associated to a presentation of 
an infrasolvmanifold and standard 3-manifold theory,
which covers most of the 4-dimensional cases.)
Other notions of infrasolvmanifold are related in \cite{[KY13]}.
 
An important special case includes most infrasolvmanifolds 
of dimension $\leq4$ (and all infranilmanifolds).
Let $T_n^+(\mathbb{R})$ be the subgroup of $GL(n,\mathbb{R})$ consisting of
upper triangular matrices with positive diagonal entries.
A Lie group $S$ is {\it triangular} if is isomorphic to a closed subgroup of
$T_n^+(\mathbb{R})$ for some $n$.
The eigenvalues of the image of each element of $S$ 
under the adjoint representation are then all real,
and so $S$ is {\it of type} {\sl R} in the terminology of \cite{[Go71]}.
(It can be shown that a Lie group is triangular if and only if it is 
1-connected and solvable of type {\sl R}.)
Two infrasolvmanifolds with discrete presentations $(S_i,\Gamma_i)$ 
where each $S_i$ is triangular (for $i=1$, 2) are affinely diffeomorphic 
if and only if their fundamental groups are isomorphic 
\cite[Theorem 3.1]{[Le95]}.
The {\it translation subgroup} $S\cap\Gamma$ of a discrete pair with $S$ 
triangular can be characterised intrinsically as the subgroup of $\Gamma$ 
consisting of the elements $g\in\Gamma$ such that all the eigenvalues of 
the automorphisms of the abelian sections of the lower central series for 
$\sqrt\Gamma$ induced by conjugation by $g$ are positive \cite{[De97]}.

Let $S$ be a connected solvable Lie group of dimension $m$, and
let $N$ be its nilradical.
If $\pi$ is a lattice in $S$ then it is torsion-free and virtually poly-$Z$ 
of Hirsch length $m$ and $\pi\cap N=\sqrt\pi$ is a lattice in $N$.     
If $S$ is 1-connected then $S/N$ is isomorphic to some 
vector group $\mathbb{R}^n$,
and $\pi/\sqrt\pi\cong\mathbb{Z}^n$.
A complete characterization of such lattices is not known, but
a torsion-free virtually poly-$Z$ group $\pi$ is a lattice in a
connected solvable Lie group $S$ if and only if $\pi/\sqrt\pi$ is abelian.
(See \cite[Sections 4.29-31]{[Rg]}.)                               

The 4-dimensional solvable Lie geometries other than the infranil
geometries are $\mathbb{S}ol^4_{m,n}$, $\mathbb{S} ol^4_0$ and 
$\mathbb{S} ol^4_1$, and the
model spaces are solvable Lie groups with left invariant metrics.
The following descriptions are based on \cite{[Wl86]}. 
The Lie group is the identity component of the isometry group 
for the geometries $\mathbb{S}ol^4_{m,n}$ and $\mathbb{S}ol^4_1$;
the identity component of $Isom(\mathbb{S}ol^4_0)$ is isomorphic 
to the semidirect product
$(\mathbb{C}\oplus\mathbb{R})\rtimes_\gamma\mathbb{C}^\times$, 
where $\gamma(z)(u,x)=(zu,|z|^{-2}x)$ for all $(u,x)$ 
in $\mathbb{C}\oplus\mathbb{R}$ and $z$ in $\mathbb{C}^\times$, 
and thus $\mathbb{S}ol_0^4$ admits an additional 
isometric action of $S^1$, 
by rotations about an axis in $\mathbb{C}\oplus\mathbb{R}\cong\mathbb{R}^3$, 
the radical of $Sol^4_0$. 

$Sol^4_{m,n}=\mathbb{R}^3\rtimes_{\theta_{m,n}}\mathbb{R}$, where 
$m$ and $n$ are integers such that the polynomial $f_{m,n}=X^3-mX^2+nX-1$ 
has distinct roots $e^a$, $e^b$ and $e^c$ (with $a<b<c$ real) and 
$\theta_{m,n}(t)$ is the diagonal matrix $diag[e^{at},e^{bt},e^{ct}]$.
Since $\theta_{m,n}(t)=\theta_{n,m}(-t)$ we may assume that $m\leq n$;
the condition on the roots then  holds if and only if $2\sqrt n\leq m\leq n$.
The metric given by $ds^2=e^{-2at}dx^2+e^{-2bt}dy^2+e^{-2ct}dz^2+dt^2$
(in the obvious global coordinates) is left invariant, and the 
automorphism of $Sol^4_{m,n}$ which sends $(x,y,z,t)$ to $(px,qy,rz,t)$
is an isometry if and only if $p^2=q^2=r^2=1$.
Let $G$ be the subgroup of $GL(4,\mathbb{R})$ of bordered matrices
$\left(\smallmatrix D&\xi\\
0&1\endsmallmatrix\right)$,
where $D=diag[\pm e^{at},\pm e^{bt},\pm e^{ct}]$ and $\xi\in\mathbb{R}^3$.
Then $Sol^4_{m,n}$ is the subgroup of $G$ with positive diagonal entries,
and $G=Isom(\mathbb{S}ol^4_{m,n})$ if $m\not= n$.
If $m=n$ then $b=0$ and $\mathbb{S}ol^4_{m,m}=\mathbb{S}ol^3\times\mathbb{E}^1$,
which admits the additional isometry sending $(x,y,z,t)$ to $(z,y,x,-t)$,
and $G$ has index 2 in $Isom(\mathbb{S}ol^3\times \mathbb{E}^1)$.
The stabilizer of the identity in the full isometry group is $(Z/2Z)^3$ 
for $Sol^4_{m,n}$ if $m\not=n$ and $D_8\times(Z/2Z)$ 
for $Sol^3\times\mathbb{R}$.
In all cases $Isom(\mathbb{S}ol^4_{m,n})\leq Aff(Sol^4_{m,n})$.

In general $\mathbb{S}ol^4_{m,n}=\mathbb{S}ol^4_{m',n'}$ if and only if 
$(a,b,c)=\lambda(a',b',c')$ for some $\lambda\not=0$.
Must $\lambda$ be rational?
(This is a case of the ``problem of the four exponentials" of transcendental 
number theory.)
If $m\not= n$ then $F_{m,n}=\mathbb{Q}[X]/(f_{m,n})$ is a totally real cubic
number field, generated over $\mathbb{Q}$ by the image of $X$.
The images of $X$ under embeddings of $F_{m,n}$ in $\mathbb{R}$
are the roots $e^a$, $e^b$ and $e^c$, and so it represents a unit of norm 1.
The group of such units is free abelian of rank 2.
Therefore if $\lambda=r/s\in\mathbb{Q}^\times$ this unit is an $r^{th}$ power 
in $F_{m,n}$ (and its $r^{th}$ root satisfies another such cubic).
It can be shown that $|r|\leq\log_2(m)$, and so
(modulo the problem of the four exponentials) there is a canonical
``minimal" pair $(m,n)$ representing each such geometry.

$Sol^4_0=\mathbb{R}^3\rtimes_\xi\mathbb{R}$, where $\xi(t)$                    
is the diagonal matrix $diag[e^t,e^t,e^{-2t}]$.
Note that if $\xi(t)$ preserves a lattice in $\mathbb{R}^3$ 
then its characteristic 
polynomial has integral coefficients and constant term $-1$. 
Since it has $e^t$ as a repeated root we must have $\xi(t)=I$.
Therefore $Sol_0^4$ does not admit any lattices.  
The metric given by the expression
$ds^2=e^{-2t}(dx^2+dy^2)+e^{4t}dz^2+dt^2$
is left invariant, and $O(2)\times O(1)$ acts via rotations and reflections
in the $(x,y)$-coordinates and reflection in the $z$-coordinate,
to give the stabilizer of the identity.
These actions are automorphisms of $Sol^4_0$,
so $Isom(\mathbb{S}ol^4_0)=Sol^4_0\rtimes(O(2)\times O(1))\leq Aff(Sol^4_0)$.
The identity component of $Isom(\mathbb{S}ol^4_0)$ is not triangular.

$Sol^4_1$ is the group of real matrices
$\{ 
\begin{pmatrix}
1& y& z\\
0& t & x\\
0&0&1
\end{pmatrix} : t>0,~x,~y,~z\in\mathbb{R} \}$.
The metric given by 
$ds^2=t^{-2}((1+x^2)(dt^2+dy^2)+t^2(dx^2+dz^2)-2tx(dtdx+dydz))$
is left invariant, and the stabilizer of the identity is $D_8$,
generated by the isometries which send $(t,x,y,z)$ to $(t,-x, y, -z)$ 
and to $t^{-1}(1,-y,-x,xy-tz)$.
These are automorphisms.
(The latter one is the restriction of the involution 
$\Omega$ of $GL(3,\mathbb{R})$ 
which sends $A$ to $J(A^{tr})^{-1}J$, where $J$ reverses the order of the 
standard basis of $\mathbb{R}^3$.)
Thus $Isom(\mathbb{S}ol^4_1)\cong{Sol_1^4}\rtimes{D_8}\leq Aff(Sol^4_1)$.
The orientation-preserving subgroup is isomorphic to the subgroup 
$\mathfrak{G}$ of $GL(3,\mathbb{R})$ generated by $Sol^4_1$ and the 
diagonal matrices $diag[-1,1,1]$ and $diag[1,1,-1]$. 
(Note that these diagonal matrices act by {\it conjugation\/} on $Sol_1^4$.)

Closed $\mathbb{S}ol^4_{m,n}$- or $\mathbb{S}ol^4_1$-manifolds  
are clearly infrasolvmanifolds.
The $\mathbb{S}ol^4_0$ case is more complicated.
Let $\tilde\gamma(z)(u,x)=(e^zu,e^{-2\mathrm{Re}(z)}x)$ 
for all $(u,x)$ in $\mathbb{C}\oplus\mathbb{R}$ and $z$ in $\mathbb{C}$.
Then 
$\widetilde I=(\mathbb{C}\oplus\mathbb{R})\rtimes_{\tilde\gamma}\mathbb{C}$
is the universal covering group of $Isom(\mathbb{S}ol^4_0)$.
If $M$ is a closed $\mathbb{S}ol^4_0$-manifold its fundamental group $\pi$
is a semidirect product $\mathbb{Z}^3\rtimes_\theta\mathbb{Z}$,
where $\theta(1)\in GL(3,\mathbb{Z})$ has two complex conjugate eigenvalues
$\lambda\not=\bar\lambda$ with $|\lambda|\not=0$ or 1 and one real eigenvalue 
$\rho$ such that $|\rho|=|\lambda|^{-2}$.
(See Chapter 8.)
If $M$ is orientable (i.e., $\rho>0$) then $\pi$ is a lattice in 
$S_\pi=(\mathbb{C}\oplus\mathbb{R})\rtimes_{\tilde\theta}\mathbb{R}
<\widetilde I$,
where $\tilde\theta(r)=\tilde\gamma(r\mathrm{log}(\lambda))$.
In general, $\pi$ is a lattice in $Aff(S_{\pi^+})$.
The action of $\widetilde I$ on $Sol^4_0$ determines 
a diffeomorphism $S_{\pi^+}/\pi\cong M$, 
and so $M$ is an infrasolvmanifold with a discrete presentation. 

\section{Orbifold bundles}

An $n$-dimensional orbifold $B$ has an open covering by subspaces of the form
$D^n/G$, where $G$ is a finite subgroup of $O(n)$.
The orbifold $B$ is {\it good} if $B=\Gamma\backslash{M}$, 
where $\Gamma$ is a discrete group acting properly 
discontinuously on a 1-connected manifold $M$; 
we then write $\pi^{orb}(B)=\Gamma$.
It is {\it aspherical} if the universal cover $M$ is contractible.
A good 2-orbifold $B$ is a quotient of $\mathbb{R}^2$, $S^2$ or $H^2$ 
by an isometric action of $\pi^{orb}(B)$, 
and so is flat, spherical or hyperbolic.
Moreover, if $B$ is compact it is a quotient of a closed surface by 
the action of a finite group. 
An orbifold is {\it bad} if it is not good.

Let $F$ be a closed manifold.
An {\it orbifold bundle with general fibre $F$ over $B$} is a map $f:M\to B$ 
which is locally equivalent to a projection 
${G\backslash(F\times D^n)}\to G\backslash {D^n}$, 
where $G$ acts freely on $F$ and effectively and orthogonally on $D^n$.
If the base $B$ has a finite regular covering $\hat B$ which is a manifold,
then $p$ induces a fibre bundle projection $\hat p:\hat M\to \hat B$
with fibre $F$, and the action of the covering group maps fibres to fibres.
Conversely, if $p_1:M_1\to B_1$ is a fibre bundle projection with
fibre $F_1$ and $G$ is a finite group which acts freely on $M_1$ 
and maps fibres to fibres then passing to orbit spaces gives an 
orbifold bundle $p:M=G\backslash M_1\to B=H\backslash B_1$ with general fibre
$F=K\backslash F_1$, where $H$ is the induced group of homeomorphisms of $B_1$ 
and $K$ is the kernel of the epimorphism from $G$ to $H$.
We shall also say that $f:M\to{B}$ is an $F$-{\it orbifold bundle}
and $M$ is an $F$-{\it orbifold bundle space}.

\begin{theorem}
{\rm[Cb]}
Let $M$ be an infrasolvmanifold.
Then there is an orbifold bundle $p:M\to B$ with general fibre an
infranilmanifold and base a flat orbifold.
\end{theorem}

\begin{proof}
Let $(S,\Gamma)$ be a presentation for $M$ and let $R$ be the nilradical of $S$.
Then $A=S/R$ is a 1-connected abelian Lie group,
and so $A\cong\mathbb{R}^d$ for some $d\geq0$.
Since $R$ is characteristic in $S$ there is a natural projection 
$q:Aff(S)\to Aff(A)$.
Let $\Gamma_S=\Gamma\cap S$ and $\Gamma_R=\Gamma\cap R$.
Then the action of $\Gamma_S$ on $S$ induces an action of the discrete group
$q(\Gamma_S)=R\Gamma_S/R$ on $A$.
The Mostow fibration for $M_1=\Gamma_S\backslash S$ is the quotient map to 
$B_1=q(\Gamma_S)\backslash A$,
which is a bundle projection with fibre $F_1=\Gamma_R\backslash R$.
Now $\Gamma_o$ is normal in $R$ \cite[Theorem 2.3, Corollary 3]{[Rg]},
and $\Gamma_R/\Gamma_o$ is a lattice in the nilpotent Lie group $R/\Gamma_o$.
Therefore $F_1$ is a nilmanifold, while $B_1$ is a torus.

The finite group $\Gamma/\Gamma_S$ acts on $M_1$, respecting the Mostow 
fibration.
Let $\overline\Gamma=q(\Gamma)$, $K=\Gamma\cap\mathrm{Ker}(q)$
and $B=\overline\Gamma\backslash A$.
Then the induced map $p:M\to B$ is an orbifold bundle projection
with general fibre the infranilmanifold
$F=K\backslash R=(K/\Gamma_o)\backslash(R/\Gamma_o)$,
and base a flat orbifold.
\end{proof}

We shall call $p:M\to B$ the {\it Mostow orbifold bundle} 
corresponding to $(S,\Gamma)$.
In Theorem 8.9 we shall use this construction to show that 
most 4-dimensional infrasolvmanifolds are determined up 
to diffeomorphism by their fundamental groups.
(Our argument fails for two virtually abelian fundamental groups.)

The {\it Gluck reconstruction} of an $S^2$-orbifold bundle $p:M\to{B}$ 
is the orbifold bundle $p^\tau:M^\tau\to{B}$
obtained by removing a product neighbourhood $S^2\times{D^2}$ 
of a general fibre and reattaching it via the nontrivial 
twist $\tau$ of $S^2\times{S^1}$.
It can be shown that $p$ is determined up to Gluck reconstruction 
by $\pi=\pi^{orb}(B)$ and the action $u:\pi\to{Aut(\pi_2(M))}$, 
and that if $B$ has a reflector curve then $p$ 
and $p^\tau$ are isomorphic as orbifold bundles over $B$.
(See \cite{[Hi13]}.)

\section{Geometric decompositions}

An $n$-manifold $M$ admits a {\it geometric decomposition} 
if it has a finite family $\mathcal{S}$ of disjoint 2-sided hypersurfaces $S$ 
such that each component of ${M\setminus\cup{S}}$ 
is geometric of finite volume,
i.e., is homeomorphic to $\Gamma\backslash X$,
for some geometry $\mathbb{X}$ and lattice $\Gamma$.
We shall call the hypersurfaces $S$ {\it cusps} and the components of 
$M\setminus\cup S$ {\it pieces} of $M$. 
The decomposition is {\it proper} if $\mathcal{S}$ is nonempty.

\begin{theorem}
If a closed $4$-manifold $M$ admits a geometric decomposition 
then either 
\begin{enumerate}
\item $M$ is geometric; or

\item $M$ is the total space of an orbifold bundle 
with general fibre $S^2$ over a hyperbolic $2$-orbifold;
or

\item the components of $M\setminus\cup S$ all have geometry 
$\mathbb{H}^2\times\mathbb{H}^2$; or

\item the components of $M\setminus\cup S$ have geometry $\mathbb{H}^4$, 
$\mathbb{H}^3\times\mathbb{E}^1$, $\mathbb{H}^2\times\mathbb{E}^2$ or
$\widetilde{\mathbb{SL}}\times\mathbb{E}^1$; or

\item the components of $M\setminus\cup S$ have geometry $\mathbb{H}^2(\mathbb{C})$ 
or $\mathbb{F}^4$.
\end{enumerate}   
In cases $(3)$, $(4)$ or $(5)$ $\chi(M)\geq0$ and in cases $(4)$ or $(5)$ 
$M$ is aspherical.
\end{theorem}                            

\begin{proof} The proof consists in considering the possible ends 
(cusps) of complete geometric 4-manifolds of finite volume. 
The hypersurfaces bounding a component of $M\setminus\cup S$ correspond 
to the ends of its interior. 
If the geometry is of solvable or compact type then there are no ends, 
since every lattice is then cocompact \cite{[Rg]}.
Thus we may concentrate on the eight geometries $\mathbb{S}^2\times\mathbb{H}^2$, 
$\mathbb{H}^2\times\mathbb{E}^2$, $\mathbb{H}^2\times\mathbb{H}^2$, 
$\widetilde{\mathbb{SL}}\times\mathbb{E}^1$, $\mathbb{H}^3\times\mathbb{E}^1$, 
$\mathbb{H}^4$, $\mathbb{H}^2(\mathbb{C})$ and $\mathbb{F}^4$.
The ends of a geometry of constant negative curvature 
$\mathbb{H}^n$ are flat \cite{[Eb80]};
since any lattice in a Lie group must meet the radical in a lattice it follows 
easily that the ends are also flat in the mixed euclidean cases 
$\mathbb{H}^3\times\mathbb{E}^1$, $\mathbb{H}^2\times\mathbb{E}^2$ and 
$\widetilde{\mathbb{SL}}\times\mathbb{E}^1$.                     
Similarly, the ends of $\mathbb{S}^2\times\mathbb{H}^2$-manifolds are 
$\mathbb{S}^2\times\mathbb{E}^1$-manifolds. 
Since the elements of $PSL(2,\mathbb{C})$ corresponding to the cusps 
of finite area hyperbolic surfaces are parabolic, the ends of 
$\mathbb{F}^4$-manifolds are $\mathbb{N}il^3$-manifolds. 
The ends of $\mathbb{H}^2(\mathbb{C})$-manifolds are also
$\mathbb{N}il^3$-manifolds \cite{[Go]},
while the ends of $\mathbb{H}^2\times\mathbb{H}^2$-manifolds 
are $\mathbb{S}ol^3$-manifolds in the irreducible cases \cite{[Sh63]},
$\mathbb{H}^2\times\mathbb{E}^1$-manifolds if the pieces are virtually products
of a closed surface with a punctured surface, 
and non-trivial graph manifolds otherwise.
Clearly if two pieces are contiguous their common cusps must be homeomorphic.
If the piece is not virtually a product of two punctured surfaces then
the inclusion of a cusp into the closure of the piece induces 
a monomorphism on fundamental group.

If $M$ is a closed 4-manifold with a geometric decomposition of type (2)
the inclusions of the cusps into the closures of the pieces induce isomorphisms on $\pi_2$, 
and a Mayer-Vietoris argument in the universal covering space $\widetilde M$ 
shows that $\widetilde M$ is homotopy equivalent to $S^2$. 
The natural foliation of $S^2\times H^2$ by 2-spheres induces a codimension-2 foliation 
on each piece, with leaves $S^2$ or $RP^2$. 
The cusps bounding the closure of a piece are $\mathbb{S}^2\times\mathbb{E}^1$-manifolds,
and hence also have codimension-1 foliations, with leaves $S^2$ or $RP^2$. 
Together these foliations give a foliation of the closure of the piece,                                              
so that each cusp is a union of leaves.
The homeomorphisms identifying cusps of contiguous pieces are isotopic
to isometries of the corresponding $\mathbb{S}^2\times\mathbb{E}^1$-manifolds. 
As the foliations of the cusps are preserved by isometries $M$ 
admits a foliation with leaves $S^2$ or $RP^2$. 
If all the leaves are homeomorphic then the projection 
to the leaf space is a submersion and so $M$ is the total 
space of an $S^2$- or $RP^2$-bundle over a hyperbolic surface.
In Chapter 10 we shall show that $S^2$- and $RP^2$-bundles over 
aspherical surfaces are geometric.
Otherwise, $M$ is the total space of an $S^2$-orbifold bundle 
over a hyperbolic 2-orbifold.

If at least one piece has an aspherical geometry other than 
$\mathbb{H}^2\times\mathbb{H}^2$ 
then all do and $M$ is aspherical.         
Since all the pieces of type $\mathbb{H}^4$, $\mathbb{H}^2(\mathbb{C})$ or $\mathbb{H}^2\times\mathbb{H}^2$ 
have strictly positive Euler characteristic while those of type              
$\mathbb{H}^3\times\mathbb{E}^1$, $\mathbb{H}^2\times\mathbb{E}^2$, 
$\widetilde{\mathbb{SL}}\times\mathbb{E}^1$ 
or $\mathbb{F}^4$ have Euler characteristic 0 we must have $\chi(M)\geq 0$. 
\end{proof}

In Theorem 10.2 we shall show every pair $(B,u)$ with $B$ an
aspherical 2-orbifold and $u$ an action of $\pi=\pi^{orb}(B)$ 
on $\mathbb{Z}$ with torsion-free kernel is realized by 
an $S^2$-orbifold bundle $p:M\to{B}$, 
with geometric total space.
It can be shown that every $S^2$-orbifold bundle over $B$ 
is isomomorphic to $p$ or to its Gluck reconstruction $p^\tau$,
and that $M^\tau$ is also geometric if and only if either $B$ 
has a reflector curve, 
in which case $p^\tau\cong{p}$ as orbifold bundles,
or $\pi$ is not generated by involutions.
However, if $B$ is the 2-sphere with $2k\geq4$ 
cone points of order 2 then $M^\tau$ is not geometric.
(See \cite{[Hi11]} and \cite{[Hi13]}.)

If in case (3) one piece is finitely covered by the product 
of a closed surface and a punctured surface then all are,
and $M$ is aspherical.
If one piece is finitely covered by the product of two punctured surfaces 
then its boundary is connected, and so there must be just two pieces.
Such manifolds are never aspherical, for as the product of two free groups
has cohomological dimension 2 and the cusp is a nontrivial graph manifold
the homomorphisms of fundamental groups induced by the inclusions of the cusp
into either piece have nontrivial kernel.
The simplest such example is the double of $T_o\times T_o$, 
where $T_o=T\setminus{int}D^2$ is the once-punctured torus.

Is there an essentially unique minimal decomposition?
Since hyperbolic surfaces are connected sums of tori, 
and a punctured torus admits a complete hyperbolic geometry of finite area,
we cannot expect that there is an unique decomposition, even in dimension 2.
Any $PD_n$-group satisfying {\it max-c} (the maximal condition on centralizers) 
has an essentially unique minimal finite splitting along virtually poly-$Z$ 
subgroups of Hirsch length $n-1$ \cite[Theorem A2]{[Kr90]}. 
(The {\it max-c} condition is unnecessary \cite{[SS]}.)
A compact non-positively curved $n$-manifold ($n\geq3$) with convex boundary 
is either flat or has a canonical decomposition along totally
geodesic closed flat hypersurfaces into pieces which are Seifert fibred or
codimension-1 atoroidal \cite{[LS00]}.
Which 4-manifolds with geometric decompositions admit such metrics?
(Closed $\widetilde{\mathbb{SL}}\times\mathbb{E}^1$-manifolds do not \cite{[Eb82]}.)

If an aspherical closed 4-manifold $M$ has 
a nontrivial geometric decomposition then the subgroups of $\pi_1(M)$ 
corresponding to the cusps contain copies of $\mathbb{Z}^2$,
and if $M$ has no pieces which are reducible
$\mathbb{H}^2\times\mathbb{H}^2$-manifolds the cuspidal subgroups 
are polycyclic of Hirsch length 3.
Closed $\mathbb{H}^4$- or $\mathbb{H}^2(\mathbb{C})$-manifolds admit no proper 
geometric decompositions,
since their fundamental groups have no noncyclic abelian subgroups \cite{[Pr43]}.
Similarly, closed $\mathbb{H}^3\times\mathbb{E}^1$-manifolds admit no proper 
decompositions, since they are finitely covered by cartesian products of 
$\mathbb{H}^3$-manifolds with $S^1$. 
Thus closed 4-manifolds with a proper geometric decomposition involving pieces of types other than 
$\mathbb{S}^2\times\mathbb{H}^2$, $\mathbb{H}^2\times\mathbb{E}^2$,
$\mathbb{H}^2\times\mathbb{H}^2$ or $\widetilde{\mathbb{SL}}\times\mathbb{E}^1$ are never geometric.
    
Many $\mathbb{S}^2\times\mathbb{H}^2$-, $\mathbb{H}^2\times\mathbb{H}^2$-,
$\mathbb{H}^2\times\mathbb{E}^2$- and $\widetilde{\mathbb{SL}}\times\mathbb{E}^1$-manifolds
admit proper geometric decompositions.  
On the other hand, a manifold with a geometric decomposition into pieces of 
type $\mathbb{H}^2\times\mathbb{E}^2$ need not be geometric.
For instance, let $G=\langle u,v,x,y\mid [u,v]=[x,y]\rangle$ 
be the fundamental group of $T\sharp T$, 
the closed orientable surface of genus 2, 
and let $\theta:G\to SL(2,\mathbb{Z})$ be the epimorphism determined by 
$\theta(u)=\left(\smallmatrix 
0& -1\\
1& 0
\endsmallmatrix\right)$,
$\theta(x)=\left(\smallmatrix 
0& 1\\
-1& 1
\endsmallmatrix\right)$
and $\theta(v)=\theta(y)=1$.
Then the semidirect product $\pi=\mathbb{Z}^2\rtimes_\theta G$ 
is the fundamental group of a torus bundle over $T\sharp T$ 
which has a geometric decomposition into two 
pieces of type $\mathbb{H}^2\times\mathbb{E}^2$, but is not geometric, 
since $\pi$ has no subgroup of finite index with centre $\mathbb{Z}^2$.

It is easily seen that each $\mathbb{S}^2\times\mathbb{E}^1$-manifold may be 
realized as the end of a complete $\mathbb{S}^2\times\mathbb{H}^2$-manifold 
with finite volume and a single end. 
However, if the manifold is orientable the ends must be orientable, 
and if it is complex analytic then they must be $S^2\times S^1$.
Every flat 3-manifold is a cusp of some (complete, finite volume)
$\mathbb{H}^4$-manifold \cite{[Ni98]}.
However if such a manifold has only one cusp the cusp cannot have
holonomy $Z/3Z$ or $Z/6Z$ \cite{[LR00]}.
The fundamental group of a cusp of an 
$\widetilde{\mathbb{SL}}\times\mathbb{E}^1$-manifold must 
have a chain of abelian normal subgroups 
$\mathbb{Z}<\mathbb{Z}^2<\mathbb{Z}^3$;
thus if the cusp is orientable the group is $\mathbb{Z}^3$ 
or $\mathbb{Z}^2\rtimes_{-I}\mathbb{Z}$.
Every $\mathbb{N}il^3$-manifold is a cusp of some 
$\mathbb{H}^2(\mathbb{C})$-manifold \cite{[McR09]}.
The ends of complex analytic $\mathbb{H}^2\times\mathbb{H}^2$-manifolds 
with irreducible fundamental group are orientable 
$\mathbb{S}ol^3$-manifolds which are mapping tori \cite{[Sh63]}.
Every orientable $\mathbb{S}ol^3$-manifold is a cusp of some
$\mathbb{H}^2\times\mathbb{H}^2$-manifold \cite{[McR08]}.

Let $\mathbb{X}=\mathbb{H}^4$, $\mathbb{H}^2(\mathbb{C})$ or
$\mathbb{H}^2\times\mathbb{H}^2$.
A nontrivial element of $Isom(\mathbb{X}^n)$ which acts freely on $X$
must have a fixed point on $\partial_\infty{X}$, 
by the Brouwer fixed point theorem.
It is parabolic if it has an unique fixed point.
Clearly if $g^k$ is parabolic for some $k>1$ then so is $g$,
and they fix the same point of $\partial_\infty{X}$.
A point $v\in\partial_\infty{X}$ represents a cusp of $\Gamma\backslash{X}$
if $\Gamma_v=\{g\in\Gamma\mid g(v)=v\}$ is a nontrivial subgroup
whose non-identity elements are parabolic.
It follows immediately that $\Gamma_v$ is square-root closed in $\Gamma$.
However cuspidal subgroups of pieces of type $\mathbb{H}^2\times\mathbb{E}^2$,
$\mathbb{F}^4$ or $\mathbb{H}^3\times\mathbb{E}^1$
need not be square-root closed in the groups of such pieces.
For instance, the natural Seifert fibration of 
a piece of type $\mathbb{H}^2\times\mathbb{E}^2$
may have base with cone point singularities of even order, 
or reflector curves.

\section{Realization of virtual bundle groups}

Every extension of one $PD_2$-group by another
may be realized by some surface bundle, by Theorem 5.2.
The study of Seifert fibred 4-manifolds and singular fibrations of complex
surfaces lead naturally to consideration of the larger class of torsion-free
groups which are virtually such extensions.
Johnson has asked whether such {\it virtual bundle groups\/}
may be realized by aspherical 4-manifolds.

\begin{theorem}
Let $\pi$ be a torsion-free group with normal 
subgroups $K<G<\pi$ such that $K$ and $G/K$ are $PD_2$-groups and 
$[\pi:G]<\infty$.
Then $\pi$ is the fundamental group of an aspherical closed smooth $4$-manifold
which is the total space of an orbifold bundle with general fibre an aspherical
closed surface over a $2$-dimensional orbifold.
\end{theorem}

\begin{proof} Let $p:\pi\to\pi/K$ be the quotient homomorphism.
Since $\pi$ is torsion-free the preimage in $\pi$ of any finite subgroup
of $\pi/K$ is a $PD_2$-group.
As the finite subgroups of $\pi/K$ have order at most $[\pi:G]$, 
we may assume that $\pi/K$ has no nontrivial finite normal subgroup,
and so is the orbifold fundamental group of some 2-dimensional orbifold $B$,
by the solution to the Nielsen realization problem for surfaces \cite{[Ke83]}.
Let $F$ be the aspherical closed surface with $\pi_1(F)\cong K$.
If $\pi/K$ is torsion-free then $B$ is a closed aspherical surface,
and the result follows from Theorem 5.2.
In general, $B$ is the union of a punctured
surface $B_o$ with finitely many cone discs and regular neighborhoods
of reflector curves (possibly containing corner points).
The latter may be further decomposed as the union of squares with a reflector
curve along one side and with at most one corner point, with
two such squares meeting along sides adjacent to the reflector curve.
These suborbifolds $U_i$ (i.e., cone discs and squares) are quotients of $D^2$
by finite subgroups of $O(2)$.
Since $B$ is finitely covered (as an orbifold) by the aspherical surface
with fundamental group $G/K$ these finite groups embed in 
$\pi_1^{\mathrm{orb}}(B)\cong\pi/K$, 
by the Van Kampen Theorem for orbifolds.

The action of $\pi/K$ on $K$ determines an action of $\pi_1(B_o)$ on $K$
and hence an $F$-bundle over $B_o$.
Let $H_i$ be the preimage in $\pi$ of $\pi_1^{orb}(U_i)$.
Then $H_i$ is torsion-free and $[H_i:K]<\infty$,
so $H_i$ acts freely and cocompactly on $\mathbb{X}^2$, 
where $\mathbb{X}^2=\mathbb{R}^2$ if $\chi(K)=0$ and 
$\mathbb{X}^2=\mathbb{H}^2$ otherwise,
and $F$ is a finite covering space of $H_i\backslash\mathbb{X}^2$.
The obvious action of $H_i$ on $\mathbb{X}^2\times D^2$
determines a bundle with general fibre $F$ over the orbifold $U_i$.
Since self homeomorphisms of $F$ are determined up to isotopy by the induced
element of $Out(K)$, bundles over adjacent suborbifolds
have isomorphic restrictions along common edges.
Hence these pieces may be assembled to give a bundle with general fibre $F$ 
over the orbifold $B$, whose total space is an aspherical closed smooth 
4-manifold with fundamental group $\pi$.
\end{proof}

We can improve upon Theorem 5.7 as follows.

\begin{cor}
Let $M$ be a $PD_4$-complex with fundamental group $\pi$. 
Then the following are equivalent.
\begin{enumerate}
\item $M$ is homotopy equivalent to the total space of an orbifold bundle with
general fibre an aspherical surface over an $\mathbb{E}^2$- or 
$\mathbb{H}^2$-orbifold;

\item $\pi$ has an $FP_2$ normal subgroup $K$ such that $\pi/K$ is 
virtually a $PD_2$-group and $\pi_2(M)=0$;

\item $\pi$ has a normal subgroup $N$ which is a $PD_2$-group 
and $\pi_2 (M)=0$.
\end{enumerate}
\end{cor} 

\begin{proof} 
Condition (1) clearly implies (2) and (3).
Conversely, if they hold the argument of Theorem 5.7 shows that 
$K$ is a $PD_2$-group and $N$ is virtually a $PD_2$-group. 
In each case (1) now follows from Theorem 7.3.
\end{proof}

It follows easily from the argument of part (1) of Theorem 5.4
that if $\pi$ is a group with a normal subgroup $K$ such that
$K$ and $\pi/K$ are $PD_2$-groups with $\zeta K=\zeta(\pi/K)=1$,
$\rho$ is a subgroup of finite index in $\pi$
and $L=K\cap\rho$ then $C_\rho(L)=1$ if and only if $C_\pi(K)=1$.
Since $\rho$ is virtually a product of $PD_2$-groups with trivial centres 
if and only if $\pi$ is, Johnson's trichotomy extends to groups
commensurate with extensions of one centreless $PD_2$-group by another.

Theorem 7.3 settles the realization question for groups of type I.
(For suppose $\pi$ has a subgroup $\sigma$ of finite index 
with a normal subgroup $\nu$ such that $\nu$ and $\sigma/\nu$
are $PD_2$-groups with $\zeta\nu=\zeta(\sigma/\nu)=1$. 
Let $G=\cap h\sigma h^{-1}$ and $K=\nu\cap G$.
Then $[\pi:G]<\infty$, $G$ is normal in $\pi$, and $K$ and $G/K$
are $PD_2$-groups.
If $G$ is of type I then $K$ is characteristic in $G$, by Theorem 5.5,
and so is normal in $\pi$.)
Groups of type II need not have such normal $PD_2$-subgroups --
although this is almost true.
In Theorem 9.9 we show that such groups are also realizable.
It is not known whether every type III extension of centreless $PD_2$-groups 
has a characteristic $PD_2$-subgroup.

If $\pi$ is an extension of $\mathbb{Z}^2$ by a normal $PD_2$-subgroup $K$
with $\zeta K=1$ then $C_\pi(K)=\sqrt\pi$, and $[\pi:KC_\pi(K)]<\infty$ 
if and only if $\pi$ is virtually ${K\times\mathbb{Z}^2}$, 
so Johnson's trichotomy extends to such groups.
The three types may be characterized by (I) $\sqrt\pi\cong\mathbb{Z}$, 
(II) $\sqrt\pi\cong\mathbb{Z}^2$, and (III) $\sqrt\pi=1$.
If $\beta_1(\pi)=2$ the subgroup $K$ is canonical, 
but in general $\pi$ may admit such subgroups of arbitrarily large genus
\cite{[Bu07]}.
Every such group is realized by an iterated mapping torus construction.
As these properties are shared by commensurate torsion-free groups, 
the trichotomy extends further to torsion-free groups 
which are virtually such extensions, 
but it is not known whether every group of this larger class 
is realized by some aspherical closed 4-manifold.

The Johnson trichotomy is inappropriate if $\zeta K\not=1$,
as there are then nontrivial extensions with trivial action ($\theta=1$).
Moreover $Out(K)$ is virtually free and so $\theta$ is never injective.
However all such groups $\pi$ may be realized by aspherical 4-manifolds, 
for either $\sqrt\pi\cong\mathbb{Z}^2$ and Theorem 7.3 applies,
or $\pi$ is virtually poly-$Z$ and is
the fundamental group of an infrasolvmanifold. (See Chapter 8.)

Aspherical orbifold bundles (with 2-dimensional base and fibre)
are determined up to fibre-preserving diffeomorphism 
by their fundamental groups, 
subject to conditions on $\chi(F)$ and $\chi^{orb}(B)$ 
analogous to those of \S2 of Chapter 5 \cite{[Vo77]}.

\section{Seifert fibrations}

A 4-manifold $S$ is {\it Seifert fibred} if it is the total space of an 
orbifold bundle with general fibre a torus or Klein bottle over a 2-orbifold.
(In \cite{[Zn85],[Ue90],[Ue91]} it is required that the general fibre be a torus.
This is always so if the manifold is orientable.
In \cite{[Vo77]} ``Seifert fibration" means 
``orbifold bundle over a 2-dimensional base", in our terms.)
It is easily seen that $\chi(S)=0$.
(In fact $S$ is finitely covered by the total space of a torus bundle 
over a surface.
This is clear if the base orbifold is good, 
and follows from the result of Ue quoted below if the base is bad.)

Let $p:S\to{B}$ be a Seifert fibration with closed aspherical base,
and let $j:F\to{S}$ be the inclusion of the fibre over the basepoint of $B$.
Let $H=j_*(\pi_1(F))$ and $A=\sqrt{H}\cong\mathbb{Z}^2$.
Then $j_*:\pi_1(F)\to\pi=\pi_1(S)$ is injective,
$A$ is a normal subgroup of $\pi$ and $\pi/A$ is virtually a surface group.
If moreover $B$ is hyperbolic $H$ is the unique maximal
solvable normal subgroup of $\pi$, and $\sqrt\pi=A$.
Let $\alpha:\pi/A\to{Aut(A)}\cong{GL(2,\mathbb{Z})}$ be the
homomorphism induced by conjugation in $\pi$,
$\mathcal{A}=\mathbb{Q}\otimes_\mathbb{Z}\sqrt\pi$ the corresponding
$\mathbb{Q}[\pi/A]$-module and $e^\mathbb{Q}(p)\in{H^2(\pi/A;\mathcal{A})}$
the class corresponding to $\pi$ as an extension of $\pi/A$ by $A$.
We shall call $\alpha$ and $e^\mathbb{Q}(p)$ the {\it action} and 
the (rational) {\it Euler class} of the Seifert fibration, respectively.
(When $A=\sqrt\pi$ we shall write $e^\mathbb{Q}(\pi)$ for
$e^\mathbb{Q}(p)$).
Let $\hat\pi$ be a normal subgroup of finite index in $\pi$ which contains $A$
and such that $\hat\pi/A$ is a $PD_2^+$-group.
Then $H^2(\pi/A;\mathcal{A})\cong{H^0(\pi/\hat\pi;H^2(\hat\pi/A;\mathcal{A}))}
\cong{H^0(\pi/\hat\pi;\mathcal{A}/\hat{I}\mathcal{A})}$,
where $\hat{I}$ is the augmentation ideal of $\mathbb{Q}[\hat\pi/A]$.
It follows that restriction to subgroups of finite index which contain $A$
is injective, and so whether $e^\mathbb{Q}(p)$ is 0 or not is invariant
under passage to such subgroups.
If $\alpha(\hat\pi)=1$ (so $\alpha$ has finite image) then
$H^2(\pi/A;\mathcal{A})\leq\mathcal{A}\cong\mathbb{Q}^2$.
(Note that if the general fibre is the Klein bottle the action is
diagonalizable, with image of order $\leq4$,
and $H^2(\pi/A;\mathcal{A})\cong\mathbb{Q}$ or 0.
The action and the rational Euler class may also be defined when the base is 
not aspherical, but we shall not need to do this.)

If $\mathbb{X}$ is one of the geometries $\mathbb{N}il^4$, 
$\mathbb{N}il^3\times\mathbb{E}^1$, $\mathbb{S}ol^3\times\mathbb{E}^1$,
$\mathbb{S}^2\times\mathbb{E}^2$, $\mathbb{H}^2\times\mathbb{E}^2$,
${\widetilde{\mathbb{SL}}\times\mathbb{E}^1}$ or
$\mathbb{F}^4$ its model space $X$ has 
a canonical foliation with leaves diffeomorphic to $R^2$ 
and which is preserved by isometries.
(For the Lie groups $Nil^4$, $Nil^3\times\mathbb{R}$ and 
$Sol^3\times\mathbb{R}$ we may take 
the foliations by cosets of the normal subgroups $\zeta_2Nil^4$, 
$\zeta{Nil^3}\times\mathbb{R}$ and ${Sol^3}'$.)
These foliations induce Seifert fibrations on quotients by lattices.
All $\mathbb{S}^3\times\mathbb{E}^1$-manifolds 
are also Seifert fibred.
Case-by-case inspection of the 74 flat 4-manifold groups shows that 
all but three have rank 2 free abelian normal subgroups,
and the representations given in \cite{[B-Z]} may be used to show 
that the corresponding manifolds are Seifert fibred.
The exceptions are the semidirect products $G_6\rtimes_\theta\mathbb{Z}$ 
where $\theta=j$, $cej$ and $abcej$. 
(See \S3 of Chapter 8 for definitions of these automorphisms.)
Closed 4-manifolds with one of the other geometries are not Seifert fibred.
(Among these, only $\mathbb{S}ol^4_{m,n}$ (with $m\not=n$),
$\mathbb{S}ol^4_0$, $\mathbb{S}ol^4_1$ and
$\mathbb{H}^3\times\mathbb{E}^1$ have closed quotients 
$M=\Gamma\backslash X$ with $\chi(M)=0$,
and for these the lattices $\Gamma$
do not have $\mathbb{Z}^2$ as a normal subgroup.)

The relationship between Seifert fibrations and geometries for orientable
4-manifolds is as follows \cite{[Ue90],[Ue91]}:

\begin{thm}
{\rm[Ue]} 
Let $S$ be a closed orientable $4$-manifold 
which is Seifert fibred over the $2$-orbifold $B$. Then
\begin{enumerate}   
\item If $B$ is spherical or bad $S$ has geometry $\mathbb{S}^3\times\mathbb{E}^1$ or 
$\mathbb{S}^2\times\mathbb{E}^2$;

\item If  $B$ is flat then  $S$ has geometry $\mathbb{E}^4$,
$\mathbb{N}il^4$, $\mathbb{N}il^3\times\mathbb{E}^1$ or $\mathbb{S}ol^3\times\mathbb{E}^1$;

\item If $B$ is hyperbolic then $S$ is geometric if and only if the action 
$\alpha$ has finite image.
The geometry is then $\mathbb{H}^2\times\mathbb{E}^2$ if 
$e^\mathbb{Q}(\pi_1(S))=0$ and ${\widetilde{\mathbb{SL}}\times\mathbb{E}^1}$ 
otherwise.

\item If $B$ is hyperbolic then $S$ has a complex structure
if and only if $B$ is orientable and $S$ is geometric. 
\end{enumerate}   
Conversely, excepting only two flat $4$-manifolds, any orientable $4$-manifold 
admitting one of these geometries is Seifert fibred.  

If the base is aspherical $S$ is determined up to diffeomorphism by $\pi_1(S)$; 
if moreover the base is hyperbolic or $S$ is geometric of type
$\mathbb{N}il^4$ or $\mathbb{S}ol^3\times\mathbb{E}^1$ 
there is a fibre-preserving diffeomorphism. 
\qed
\end{thm}
 
We have corrected a minor oversight in \cite{[Ue90]}; 
there are in fact {\it two} orientable flat 4-manifolds 
which are not Seifert fibred.
If the base is bad or spherical then $S$ may admit 
many inequivalent Seifert fibrations.
(See also \S10 of Chapter 8 and \S2 of Chapter 9 for further discussion of
the flat base and hyperbolic base cases, respectively.)

In general, 4-manifolds which are Seifert fibred over aspherical bases are
determined up to diffeomorphism by their fundamental groups.
This was first shown by Zieschang for the cases with base a hyperbolic orbifold
with no reflector curves and general fibre a torus \cite{[Zi69]},
and the full result is due to Vogt \cite{[Vo77]}.
Kemp has shown that a nonorientable aspherical Seifert fibred 4-manifold 
is geometric if and only if its orientable double covering space is geometric
[Ke]. (See also Theorems 9.4 and 9.5).
Closed 4-manifolds which fibre over $S^1$ with fibre a small Seifert fibred 
3-manifold are also determined by their fundamental groups \cite{[Oh90]}.
This class includes many nonorientable Seifert fibred 4-manifolds over 
bad, spherical or flat bases, but not all. 
    
The homotopy type of a $\mathbb{S}^2\times\mathbb{E}^2$-manifold is determined 
up to finite ambiguity by its fundamental group 
(which is virtually $\mathbb{Z}^2$), 
Euler characteristic (which is 0) and Stiefel-Whitney classes.
There are just nine possible fundamental groups. 
Six of these have infinite abelianization, and the above invariants
determine the homotopy type in these cases. 
(See Chapter 10.)
The homotopy type of a $\mathbb{S}^3\times\mathbb{E}^1$-manifold is determined 
by its fundamental group (which has two ends), Euler characteristic 
(which is 0), orientation character $w_1$ and first $k$-invariant in 
$H^4(\pi;\pi_3)$. 
(See Chapter 11.)

Let $S$ be a Seifert fibred 4-manifold with base an flat orbifold,
and let $\pi=\pi_1(S)$.
Then $\chi(S)=0$ and $\pi$ is solvable of Hirsch length 4,
and so $S$ is homeomorphic to an infrasolvmanifold, 
by Theorem 6.11 and \cite{[AJ76]}.
Every such group $\pi$ is the fundamental group of some 
Seifert fibred geometric 4-manifold, 
and so $S$ is in fact diffeomorphic to an infrasolvmanifold \cite{[Vo77]}.
(See Chapter 8.\S9 and Theorem 8.10 below.)
The general fibre must be a torus if the geometry is 
$\mathbb{N}il^4$ or $\mathbb{S}ol^3\times\mathbb{E}^1$,
since $Out(\pi_1(Kb))$ is finite.

As $\mathbb{H}^2\times\mathbb{E}^2$- and 
$\widetilde{\mathbb{SL}}\times\mathbb{E}^1$-manifolds are aspherical, 
they are determined up to homotopy equivalence by their fundamental groups.  
(See Chapter 9.)
Theorem 7.3 specializes to give the following characterization of the 
fundamental groups of Seifert fibred 4-manifolds over hyperbolic bases.    

\begin{theorem}
A group $\pi$ is the fundamental group of a closed $4$-manifold 
which is Seifert fibred over a hyperbolic $2$-orbifold 
if and only if it is torsion-free, 
$\sqrt\pi\cong\mathbb{Z}^2$, 
$\pi/\sqrt\pi$ is virtually a $PD_2$-group and the
maximal finite normal subgroup of $\pi/\sqrt\pi$ has order at most $2$.
\qed
\end{theorem}

If $\sqrt\pi$ is central ($\zeta\pi\cong\mathbb{Z}^2$) 
the corresponding Seifert fibred manifold $M(\pi)$ 
admits an effective torus action with finite isotropy subgroups. 

\section{Complex surfaces and related structures}

In this section we shall summarize what we need from 
\cite{[BHPV],[Ue90],[Ue91],[Wl86]}
and \cite{[GS]}, and we refer to these sources for more details.
    
A {\it complex surface} shall mean a compact connected 
complex analytic manifold $S$ of complex dimension 2.                          
It is {\it K\"ahler\/} (and thus diffeomorphic to 
a projective algebraic surface) if and only if $\beta_1(S)$ is even.
Since the K\"ahler condition is local, all finite covering spaces 
of such a surface must also have $\beta_1$ even.
If $S$ has a complex submanifold $L\cong{CP^1}$ 
with self-intersection $-1$ then $L$ may be blown down: 
there is a complex surface $S_1$ and a holomorphic map $p:S\to S_1$ 
such that $p(L)$ is a point and $p$ restricts to
a biholomorphic isomorphism from $S\setminus{L}$ to $S_1\setminus{p(L)}$.
In particular, $S$ is diffeomorphic to $S_1\sharp\overline{CP^2}$.
If there is no such embedded projective line $L$ the surface is {\it minimal}.
Every surface has a minimal model, 
and the model is unique if it is neither rational nor ruled.
    
For many of the 4-dimensional geometries $(X,G)$ the identity component $G_o$ 
of the isometry group preserves a natural complex structure on $X$,
and so if $\pi$ is a discrete subgroup of $G_o$ which acts freely on $X$ 
the quotient $\pi\backslash X$ is a complex surface.
This is clear for the geometries $\mathbb{CP}^2$, 
$\mathbb{S}^2\times\mathbb{S}^2$, 
$\mathbb{S}^2\times\mathbb{E}^2$,
$\mathbb{S}^2\times\mathbb{H}^2$, $\mathbb{H}^2\times\mathbb{E}^2$, 
$\mathbb{H}^2\times\mathbb{H}^2$ and $\mathbb{H}^2(\mathbb{C})$. 
(The corresponding model spaces may be identified with $CP^2$, 
$CP^1\times CP^1$, $CP^1\times\mathbb{C}$, $CP^1\times H^2$,
$H^2\times\mathbb{C}$, 
$H^2\times H^2$ and the unit ball in $\mathbb{C}^2$, respectively, 
where $H^2$ is identified with the upper half plane.)
It is also true for $\mathbb{N}il^3\times\mathbb{E}^1$, 
$\mathbb{S}ol^4_0$, $\mathbb{S}ol^4_1$, 
$\widetilde{\mathbb{SL}}\times\mathbb{E}^1$ and $\mathbb{F}^4$.
In addition, the subgroups $\mathbb{R}^4\rtimes{U(2)}$ of $E(4)$ and
$U(2)\times\mathbb{R}$ of $Isom(\mathbb{S}^3\times\mathbb{E}^1)$ 
act biholomorphically on $\mathbb{C}^2$ and $\mathbb{C}^2\setminus\{0\}$,
respectively, 
and so some $\mathbb{E}^4$- and $\mathbb{S}^3\times\mathbb{E}^1$-manifolds 
have complex structures.
No other geometry admits a compatible complex structure.
Since none of the model spaces contain an embedded $S^2$ 
with self-intersection $-1$,
any complex surface which admits a compatible geometry must be minimal.

Complex surfaces may be coarsely classified by their Kodaira dimension $\kappa$,
which may be $-\infty$, 0, 1 or 2. 
Within this classification, minimal surfaces may be 
further classified into a number of families. 
We have indicated in parentheses where the geometric complex surfaces 
appear in this classification.
(The dashes signify families which include nongeometric surfaces.)

$\kappa=-\infty$: Hopf surfaces ($\mathbb{S}^3\times\mathbb{E}^1$, --);
Inoue surfaces ($\mathbb{S}ol^4_0$, $\mathbb{S}ol^4_1$); 
(other) surfaces of class VII with $\beta_2>0$ (--);
rational surfaces ($\mathbb{CP}^2$, $\mathbb{S}^2\times\mathbb{S}^2$); 
ruled surfaces ($\mathbb{S}^2\times\mathbb{E}^2$, 
$\mathbb{S}^2\times\mathbb{H}^2$, --).
        
$\kappa=0$: complex tori ($\mathbb{E}^4$); 
hyperelliptic surfaces ($\mathbb{E}^4$); 
Enriques surfaces (--); 
K3 surfaces (--); Kodaira surfaces ($\mathbb{N}il^3\times\mathbb{E}^1$). 

$\kappa=1$: minimal properly elliptic surfaces ($\widetilde{\mathbb{SL}}\times\mathbb{E}^1$, 
$\mathbb{H}^2\times\mathbb{E}^2$).

$\kappa=2$: minimal (algebraic) surfaces of general type 
($\mathbb{H}^2\times\mathbb{H}^2$, $\mathbb{H}^2(\mathbb{C})$, --).

A {\it Hopf surface} is a complex surface whose universal covering space is 
homeomorphic to $S^3\times\mathbb{R}\cong\mathbb{C}^2\setminus\{0\}$.
Some Hopf surfaces admit no compatible geometry, and there are
$\mathbb{S}^3\times\mathbb{E}^1 $-manifolds that admit no complex structure.
The {\it Inoue} surfaces are exactly the complex surfaces admitting one of the
geometries $\mathbb{S}ol^4_0$ or $\mathbb{S}ol^4_1$. 
Surfaces of {\it class VII} have $\kappa=-\infty$ and $\beta_1=1$,
and are not yet fully understood.
(A theorem of Bogomolov asserts that every minimal complex surface 
of class VII with $\beta_2(S)=0$ is either a Hopf surface or an Inoue surface. 
See \cite{[Tl94]} for a complete proof.)

A {\it rational surface} is a complex surface birationally equivalent to $CP^2$.
Minimal rational surfaces are diffeomorphic to $CP^2$ or to $CP^1\times CP^1$.
A {\it ruled surface} is a complex surface which is holomorphically fibred over
a smooth complex curve (closed orientable 2-manifold) of genus $g>0$ 
with fibre $CP^1$.
Rational and ruled surfaces may be characterized as the complex surfaces 
$S$ with $\kappa(S)=-\infty$ and $\beta_1(S)$ even.
Not all ruled surfaces admit geometries compatible with 
their complex structures.
    
A {\it complex torus} is a quotient of $\mathbb{C}^2$ by a lattice, 
and a {\it hyperelliptic surface} is one properly covered by a complex torus. 
If $S$ is a complex surface which is homeomorphic to a flat 4-manifold 
then $S$ is a complex torus or is hyperelliptic, 
since it is finitely covered by a complex torus.  
Since $S$ is orientable and $\beta_1(S)$ is even $\pi=\pi_1 (S)$ 
must be one of the eight flat 4-manifold groups of orientable type 
and with $\pi\cong\mathbb{Z}^4$ or $I(\pi)\cong\mathbb{Z}^2$.
In each case the holonomy group is cyclic, and so is conjugate (in $GL^+(4,\mathbb{R})$)
to a subgroup of $GL(2,\mathbb{C})$. 
(See Chapter 8.)
Thus all of these groups may be realized by complex surfaces.
A {\it Kodaira surface} is a surface with $\beta_1$ odd
and which has a finite cover which fibres holomorphically 
over an elliptic curve with fibres of genus 1.
                    
An {\it elliptic surface} $S$ is a complex surface 
which admits a holomorphic map $p$ to a complex curve 
such that the generic fibres of $p$ are diffeomorphic to the torus $T$.
If the elliptic surface $S$ has no singular fibres it is Seifert fibred, 
and it then has a geometric structure if and only if the base is a good orbifold.
An orientable Seifert fibred 4-manifold over a hyperbolic base 
has a geometric structure if and only if it is an elliptic surface 
without singular fibres \cite{[Ue90]}.
The elliptic surfaces $S$ with $\kappa(S)=-\infty$ and $\beta_1(S)$ odd are
the geometric Hopf surfaces. The elliptic surfaces $S$ with $\kappa(S)=-\infty$ 
and $\beta_1(S)$ even are the cartesian products of elliptic curves with $CP^1$.

All rational, ruled and hyperelliptic surfaces are projective algebraic surfaces, 
as are all surfaces with $\kappa=2$. 
Complex tori and surfaces with geometry $\mathbb{H}^2\times\mathbb{E}^2$ are 
diffeomorphic to projective algebraic surfaces.
Hopf, Inoue and Kodaira surfaces and surfaces with geometry 
$\widetilde{\mathbb{SL}}\times\mathbb{E}^1$ all have $\beta_1$ odd,
and so are not K\"ahler, let alone projective algebraic.

An {\it almost complex structure} on a smooth $2n$-manifold $M$ is a reduction 
of the structure group of its tangent bundle to
$GL(n,\mathbb{C})<GL^+(2n,\mathbb{R})$.
Such a structure determines an orientation on $M$.
If $M$ is a closed oriented 4-manifold and $c\in H^2(M;\mathbb{Z})$
then there is an almost complex structure on $M$ with first Chern class $c$
and inducing the given orientation if and only if 
$c\equiv w_2(M)$ {\it mod} $(2)$ and $c^2\cap[M]=3\sigma(M)+2\chi(M)$,
by a theorem of Wu. 
(See \cite[Chapter 1. Appendix]{[GS]} 
for a recent account.)

A {\it symplectic structure} on a closed smooth manifold $M$ is a closed 
nondegenerate 2-form $\omega$.
Nondegenerate means that for all $x\in M$ and all $u\in T_xM$ there is 
a $v\in T_xM$ such that $\omega(u,v)\not=0$.
Manifolds admitting symplectic structures are even-dimensional and orientable.
A condition equivalent to nondegeneracy is that the $n$-fold wedge 
$\omega^{\wedge n}$ is nowhere 0, where $2n$ is the dimension of $M$.
The $n^{th}$ cup power of the corresponding cohomology class 
$[\omega]$ is then a nonzero element of $H^{2n}(M;\mathbb{R})$.
Any two of a riemannian metric, a symplectic structure and an almost complex
structure together determine a third, if the given two are compatible.
In dimension 4, this is essentially equivalent to the fact that
$SO(4)\cap Sp(4)=SO(4)\cap GL(2,\mathbb{C})=Sp(4)\cap GL(2,\mathbb{C})=U(2)$,
as subgroups of $GL(4,\mathbb{R})$.
(See \cite{[GS]} for a discussion of relations between these structures.)
In particular, K\"ahler surfaces have natural symplectic structures, 
and symplectic 4-manifolds admit compatible almost complex tangential 
structures.
However orientable $\mathbb{S}ol^3\times\mathbb{E}^1$-manifolds which fibre 
over $T$ are symplectic \cite{[Ge92]} but have no complex structure 
(by the classification of surfaces) and Hopf surfaces are complex manifolds 
with no symplectic structure (since $\beta_2=0$).

%% file: m5-8.tex
\chapter{Solvable Lie geometries}

The main result of this chapter is the characterization of 
4-dimensional infrasolvmanifolds up to homeomorphism, given in \S1.
All such manifolds are either mapping tori of self homeomorphisms of 
3-dimensional infrasolvmanifolds
or are unions of two twisted $I$-bundles over such 3-manifolds.
In the rest of the chapter we consider each of the possible 
4-dimensional geometries of solvable Lie type.

In \S2 we determine the automorphism groups of the flat 3-manifold groups,
while in \S3 and \S4 we determine {\it ab initio} the 74 flat 4-manifold groups.
There have been several independent computations of these groups; 
the consensus reported on page 126 of \cite{[Wo]} is that there are 
27 orientable groups and 48 nonorientable groups. 
However the tables of 4-dimensional crystallographic groups in 
\cite{[B-Z]} list only 74 torsion-free groups,
of which 27 are orientable and 47 are nonorientable.
As these computer-generated tables give little insight into 
how these groups arise, 
and as the earlier computations were never published in detail, 
we shall give a direct and elementary computation, motivated by Lemma 3.14.
Our conclusions as to the numbers of groups with abelianization of given rank, 
isomorphism type of holonomy group and orientation type agree with those of \cite{[B-Z]} and \cite{[LRT13]}.
(We refer to \cite{[LRT13]} for details of some gaps relating to the cases
with $\beta_1=0$ in earlier versions of this book.)

There are infinitely many examples for each of the other geometries.
In \S5 we show how these geometries may be distinguished, in terms
of the group theoretic properties of their lattices.
In \S6, \S7 and \S8 we consider mapping tori of self homeomorphisms of 
$\mathbb{E}^3$-, $\mathbb{N}il^3$- and $\mathbb{S}ol^3$-manifolds, respectively.
In \S9 we show directly that ``most" groups allowed by Theorem 8.1
are realized geometrically and outline classifications for them,
while in \S10 we show that ``most" 4-dimensional infrasolvmanifolds
are determined up to diffeomorphism by their fundamental groups.

\section{The characterization}

In this section we show that 4-dimensional infrasolvmanifolds may be 
characterized up to homeomorphism in terms of the fundamental group and Euler 
characteristic.

\begin{theorem} 
Let $M$ be a closed $4$-manifold with fundamental group $\pi$
and such that $\chi(M)=0$. The following conditions are equivalent:
\begin{enumerate}
\item $\pi$ is torsion-free and virtually poly-$Z$ and $h(\pi)=4$;

\item $h(\sqrt\pi)\geq3$;

\item $\pi$ has an elementary amenable normal subgroup $\rho$ with 
$h(\rho)\geq 3$, and $H^2 (\pi;\mathbb{Z}[\pi])=0$; and

\item $\pi$ is restrained, every finitely generated subgroup of $\pi$ is $FP_3$ 
and $\pi$ maps onto a virtually poly-$Z$ group $Q$ with $h(Q)\geq3$.
\end{enumerate}
Moreover, if these conditions hold then $M$ is aspherical, 
and is determined up to homeomorphism by $\pi$,
and every automorphism of $\pi$ may be realized by a self homeomorphism of $M$.
\end{theorem}
                       
\begin{proof} 
If (1) holds then $h(\sqrt\pi)\geq 3$, by Theorem 1.6, and so (2) holds.
There is an epimorphism $\lambda:\pi\to\mathbb{Z}$ or $D$, by Lemma 3.14.
Then $E=\mathrm{Ker}(\lambda|_{\sqrt\pi})$ is locally nilpotent, 
normal, $[\pi:E]=\infty$ and $h(E)\geq{h(\sqrt\pi)-1}$.
Hence (2) implies (3), by Theorem 1.17.
If (3) holds then $\pi$ has one end, by Lemma 1.15,
and $\beta_1^{(2)}(\pi)=0$, by Corollary 2.3.1.
Hence $M$ is aspherical, by Corollary 3.5.2.
Hence $\pi$ is a $PD_4 $-group and $3\leq h(\rho)\leq c.d.\rho\leq 4$.
In particular, $\rho$ is virtually solvable, by Theorem 1.11.
If $c.d.\rho=4$ then $[\pi:\rho]$ is finite, by Strebel's Theorem, 
and so $\pi$ is virtually solvable also.
If $c.d.\rho=3$ then $c.d.\rho=h(\rho)$ and so $\rho$ is a duality group 
and is $FP$ \cite{[Kr86]}.
Therefore $H^q (\rho;\mathbb{Q}[\pi])\cong 
H^q (\rho;\mathbb{Q}[\rho])\otimes \mathbb{Q}[\pi/\rho]$ 
and is 0 unless $q=3$.
It then follows from the LHSSS for $\pi$ as an extension of $\pi/\rho$ by 
$\rho$ (with coefficients $\mathbb{Q}[\pi]$) that 
$H^4 (\pi;\mathbb{Q}[\pi])\cong H^1 (\pi/\rho;\mathbb{Q}[\pi/\rho])\otimes 
H^3 (\rho;\mathbb{Q}[\rho])$.
Therefore $H^1 (\pi/\rho;\mathbb{Q}[\pi/\rho])\cong\mathbb{Q}$, 
so $\pi/\rho$ has two ends and we again find that $\pi$ is virtually solvable.
In all cases $\pi$ is torsion-free and virtually poly-$Z$
\cite[Theorem 9.23]{[Bi]}, and $h(\pi)=4$. 

If (4) holds then $\pi$ is an ascending HNN extension $\pi\cong B*_\phi$ 
with base $FP_3$ and so $M$ is aspherical, by Theorem 3.16. 
As in Theorem 2.13 we may deduce from \cite{[BG85]} 
that $B$ must be a $PD_3$-group and $\phi$ an isomorphism,
and hence $B$ and $\pi$ are virtually poly-$Z$. 
Conversely (1) clearly implies (4).

The final assertions follow from Theorem 2.16 of \cite{[FJ]}.
\end{proof}

Does the hypothesis $h(\rho)\geq 3$ in (3)
imply $H^2 (\pi;\mathbb{Z}[\pi])=0$?
The examples $F\times S^1\times S^1$ where $F=S^2$ or is a closed hyperbolic 
surface show that the condition that $h(\rho)>2$ is necessary. 
(See also \S1 of Chapter 9.)

\begin{cor}
The $4$-manifold $M$ is homeomorphic to an infrasolvmanifold if and only if 
the equivalent conditions of Theorem $8.1$ hold.
\end{cor} 

\begin{proof} 
If $M$ is homeomorphic to an infrasolvmanifold then $\chi(M)=0$, 
$\pi$ is torsion-free and virtually poly-$Z$ and $h(\pi)=4$. (See Chapter 7.)
Conversely, if these conditions hold then $\pi$ is the fundamental 
group of an infrasolvmanifold, by \cite{[AJ76]}. 
\end{proof}

It is easy to see that all such groups are realizable by closed smooth
4-manifolds with Euler characteristic 0.

\begin{theorem} 
If $\pi$ is torsion-free and virtually poly-$Z$ of Hirsch length $4$ 
then it is the fundamental group of a closed smooth $4$-manifold $M$
which is either a mapping torus of a self homeomorphism of a closed 
$3$-dimensional infrasolvmanifold or is the union of two twisted $I$-bundles 
over such a $3$-manifold.
Moreover, the $4$-manifold $M$ is determined up to homeomorphism by the group.
\end{theorem}

\begin{proof} The Eilenberg-Mac Lane space $K(\pi,1)$ is a $PD_4 $-complex
with Euler characteristic 0.
By Lemma 3.14, either there is an epimorphism $\phi:\pi\to\mathbb{Z}$, 
in which case $\pi$ is a semidirect product 
$G\rtimes_\theta\mathbb{Z}$ where $G=\mathrm{Ker}(\phi)$,
or $\pi\cong G_1 *_G G_2 $ where $[G_1:G]=[G_2:G]=2$. 
The subgroups $G$, $G_1$ and $G_2$ are torsion-free and virtually poly-$Z$.
Since in each case $\pi/G$ has Hirsch length 1 these subgroups have Hirsch length 3
and so are fundamental groups of closed 3-dimensional infrasolvmanifolds. 
The existence of such a manifold now follows by standard 3-manifold topology,
while its uniqueness up to homeomorphism was proven in Theorem 6.11.
\end{proof}

The first part of this theorem may be stated and proven 
in purely algebraic terms, 
since torsion-free virtually poly-$Z$ groups are Poincar\'e duality groups.
(See \cite[Chapter III]{[Bi]}.) 
If $\pi$ is such a group then either it is virtually nilpotent 
or $\sqrt\pi\cong\mathbb{Z}^3$ or $\Gamma_q$ for some $q$, 
by Theorems 1.5 and 1.6.
In the following sections we shall consider how such 
groups may be realized geometrically.
The geometry is largely determined by $\sqrt\pi$.
We shall consider first the virtually abelian cases.

\section{Flat 3-manifold groups and their automorphisms}

The flat $n$-manifold groups for $n\leq 2$ are $\mathbb{Z}$, 
$\mathbb{Z}^2$ and $K=\mathbb{Z}\rtimes_{-1}\! \mathbb{Z}$, 
the Klein bottle group.
There are six orientable and four nonorientable flat 3-manifold groups. 
The first of the orientable flat 3-manifold groups $G_1$ - $G_6$ 
s $G_1=\mathbb{Z}^3$.
The next four have $I(G_i)\cong\mathbb{Z}^2$ and are 
semidirect products $\mathbb{Z}^2\rtimes_T\mathbb{Z}$ where 
$T=-I$, 
$\left(\begin{smallmatrix}
0&-1\\
1&-1
\end{smallmatrix}\right),$
$\left(\begin{smallmatrix}
0&-1\\
1&0
\end{smallmatrix}\right)$
or
$\left(\begin{smallmatrix}
0&-1\\
1&1
\end{smallmatrix}\right)$,
respectively, is an element of finite order in $SL(2,\mathbb{Z})$. 
These groups all have cyclic holonomy groups, of orders 2, 3, 4 and 6, respectively.
The group $G_6$ is the group of the {\it Hantzsche-Wendt} flat 3-manifold, 
and has a presentation
\begin{equation*} 
\langle x,y\mid xy^2x^{-1}=y^{-2},\medspace yx^2y^{-1}=x^{-2}\rangle.
\end{equation*}  
Its maximal abelian normal subgroup is generated by $x^2,y^2$ and $(xy)^2$
and its holonomy group is the diagonal subgroup of $SL(3,\mathbb{Z})$, 
which is isomorphic to $(Z/2Z)^2$.
(This group is the generalized free product of two copies of $K$,
amalgamated over their maximal abelian subgroups, and so maps onto $D$.)

The nonorientable flat 3-manifold groups $B_1$ - $B_4$ 
are semidirect products $K\rtimes_\theta\mathbb{Z}$, 
corresponding to the classes in $Out(K)\cong (Z/2Z)^2$.
In terms of the presentation $\langle x,y\mid xyx^{-1} =y^{-1}\rangle$ 
for $K$ these classes 
are represented by the automorphisms $\theta$ which fix $y$ and send $x$ to 
$x,xy, x^{-1}$ and $x^{-1}y$, respectively.
The groups $B_1$ and $B_2$ are also semidirect products 
$\mathbb{Z}^2 \rtimes_T\mathbb{Z}$, where
$T=\left(\begin{smallmatrix}
1&0\\
0&-1
\end{smallmatrix}\right)$ or
$\left(\begin{smallmatrix}
0&1\\
1&0
\end{smallmatrix}\right)$
has determinant $-1$ and $T^2=I$. They have holonomy groups of order 2,
while the holonomy groups of $B_3$ and $B_4 $ are isomorphic to $(Z/2Z)^2$.

All the flat 3-manifold groups either map onto $\mathbb{Z}$ or map onto $D$.
The methods of this chapter may be easily adapted to find all such groups. 
Assuming these are all known we may use Sylow theory and some calculation
to show that there are no others. We sketch here such an argument.
Suppose that $\pi$ is a flat 3-manifold group with finite abelianization.
Then $0=\chi(\pi)=1+\beta_2 (\pi)-\beta_3 (\pi)$, 
so $\beta_3 (\pi)\not=0$ and $\pi$ must be orientable.
Hence the holonomy group $F=\pi/T(\pi)$ is a subgroup of $SL(3,\mathbb{Z})$.
Let $f$ be a nontrivial element of $F$.
Then $f$ has order 2, 3, 4 or 6, and has a $+1$-eigenspace of rank 1, 
since it is orientation preserving.
This eigenspace is invariant under the action of the normalizer 
$N_F (\langle f\rangle)$,
and the induced action of $N_F (\langle f\rangle)$ 
on the quotient space is faithful. 
Thus $N_F (\langle f\rangle)$ is isomorphic to a subgroup of 
$GL(2,\mathbb{Z})$ and so is cyclic or dihedral of order dividing 24. 
This estimate applies to the Sylow subgroups of $F$,
since $p$-groups have nontrivial centres,
and so the order of $F$ divides 24.
If $F$ has a nontrivial cyclic normal subgroup then $\pi$ has a normal subgroup 
isomorphic to $\mathbb{Z}^2$ and hence maps onto $\mathbb{Z}$ or $D$.
Otherwise $F$ has a nontrivial Sylow 3-subgroup $C$ which is not normal in $F$.
The number of Sylow 3-subgroups is congruent to 1 $mod~(3)$ 
and divides the order of $F$.
The action of $F$ by conjugation on the set of such subgroups is transitive.
It must also be faithful.
(For otherwise $\cap_{g\in F}gN_F(C)g^{-1}\not=1$.
As $N_F(C)$ is cyclic or dihedral it would follow that
$F$ must have a nontrivial cyclic normal subgroup, contrary to hypothesis.)
Hence $F$ must be $A_4 $ or $S_4$, and so contains $V\cong(Z/2Z)^2$ as a normal subgroup.
Suppose that $G$ is a flat 3-manifold group with holonomy $A_4$.
It is easily seen that $G_6$ is the only flat 3-manifold group 
with holonomy $(Z/2Z)^2$, and so
we may assume that the images in $SL(3,\mathbb{Z})$ 
of the elements of order 2 are diagonal matrices.
It then follows easily that the images of the elements of order 3 are (signed)
permutation matrices. 
(Solve the linear equations $wu=vw$ and $wv=uvw$ in $SL(3,\mathbb{Z})$, 
where $u=diag[1,-1,-1]$ and $v=diag[-1,-1,1]$.)
Hence $G$ has a presentation of the form 

\centerline{$\langle\mathbb{Z}^3,u,v,w\mid ux=xu, yuy=u, zuz=u, 
xvx=v, yvy=v, zv=vz,
wx=zw,$}

\centerline{$ wy=xw, wz=yw, wu=vw, u^2=x,w^3=x^ay^bz^c,(uw)^3=x^py^qz^r\rangle.$}

\noindent It may be checked that no such group is torsion-free.
Therefore neither $A_4$ nor $S_4$ can be
the holonomy group of a flat 3-manifold. 

We shall now determine the (outer) automorphism groups of each of the flat 3-manifold groups.
Clearly $Out(G_1)=Aut(G_1)=GL(3,\mathbb{Z})$.                   
If $2\leq i\leq 5$ let $t\in G_i$ represent a generator of the quotient $G_i/I(G_i)\cong\mathbb{Z}$.
The automorphisms of $G_i$ must preserve the characteristic subgroup 
$I(G_i)$ and so may be identified with triples 
$(v,A,\epsilon)\in\mathbb{Z}^2\times GL(2,\mathbb{Z})\times\{\pm1\}$ 
such that $ATA^{-1} =T^\epsilon$ and which act via $A$ 
on $I(G_i)=\mathbb{Z}^2$ and send $t$ to $t^\epsilon v$.
Such an automorphism is orientation preserving if and only if $\epsilon=det(A)$.
The multiplication is given by $(v,A,\epsilon)(w,B,\eta)=(\Xi v+Aw,AB,\epsilon\eta)$,
where $\Xi=I$ if $\eta=1$ and $\Xi=-T^\epsilon $ if $\eta=-1$.
The inner automorphisms are generated by $(0,T,1)$ and 
$((T-I)\mathbb{Z}^2,I,1)$.

In particular, 
$Aut(G_2)\cong(\mathbb{Z}^2\rtimes_\alpha GL(2,\mathbb{Z}))\times\{\pm1\}$,
where $\alpha$ is the natural action of $GL(2,\mathbb{Z})$ on $\mathbb{Z}^2$,
for $\Xi$ is always $I$ if $T=-I$.
The involution $(0,I,-1)$ is central in $Aut(G_2 )$, 
and is orientation reversing.
Hence $Out(G_2)$ is isomorphic to $((Z/2Z)^2\rtimes_{P\alpha} PGL(2,\mathbb{Z}))\times (Z/2Z)$, where $P\alpha$
is the induced action of $PGL(2,\mathbb{Z})$ on $(Z/2Z)^2$.

If $n=3$, 4 or 5 the normal subgroup $I(G_i)$ may be viewed as a module 
over the ring $R=\mathbb{Z}[t]/(\phi(t))$, where $\phi(t)=t^2+t+1$, 
$t^2+1$ or $t^2-t+1$, respectively.
As these rings are principal ideal domains and $I(G_i)$ is torsion-free of 
rank 2 as an abelian group, in each case it is free of rank 1 as an $R$-module.
Thus matrices $A$ such that $AT=TA$ correspond to units of $R$.
Hence automorphisms of $G_i$ which induce the identity on $G_i/I(G_i)$ have the form
$(v,\pm T^m,1)$, for some $m\in\mathbb{Z}$ and $v\in\mathbb{Z}^2$.
There is also an involution 
$(0,\left(\begin{smallmatrix}
0&1\\
1&0
\end{smallmatrix}\right),-1)$
which sends $t$ to $t^{-1}$. In all cases $\epsilon=det(A)$.
It follows that $Out(G_3)\cong S_3\times(Z/2Z)$, $Out(G_4)\cong (Z/2Z)^2$ and $Out(G_5)=Z/2Z $.
All these automorphisms are orientation preserving.

The subgroup $A$ of $G_6$ generated by $\{ x^2,y^2,(xy)^2\} $ is
the maximal abelian normal subgroup of $G_6$, and $G_6/A\cong (Z/2Z)^2$.
Let $a$, $b$, $c$, $d$, $e$, $f$, $i$ and $j$ be the automorphisms 
of $G_6$ which send $x$ to $x^{-1},x,x,x,y^2x,(xy)^2x,y,xy$ and 
$y$ to $y,y^{-1},(xy)^2y,x^2y,y,(xy)^2y,x,x$, respectively.
The natural homomorphism from $Aut(G_6)$ to 
$Aut(G_6/A)\cong GL(2,\mathbb{F}_2)$ is onto, 
as the images of $i$ and $j$ generate $GL(2,\mathbb{F}_2)$,
and its kernel $E$ is generated by $\{ a,b,c,d,e,f\}$.
(For an automorphism which induces the identity on $G_6/A$ must send $x$ to $x^{2p}y^{2q}(xy)^{2r}x$,
and $y$ to $x^{2s}y^{2t}(xy)^{2u}y$. 
The images of $x^2$, $y^2$ and $(xy)^2$ are then
$x^{4p+2}$, $y^{4t+2}$ and $(xy)^{4(r-u)+2}$, which generate $A$ if and only if
$p=0$ or $-1$, $t=0$ or $-1$ and $r=u-1$ or $u$.
Composing such an automorphism appropriately with $a$, $b$ and $c$ we may acheive $p=t=0$ and $r=u$.
Then by composing with powers of $d$, $e$ and $f$ we may obtain the identity automorphism.)
The inner automorphisms are generated by $bcd$ (conjugation by $x$) and 
$acef$ (conjugation by $y$). 
Then $Out(G_6)$ has a presentation
\begin{gather*} 
\langle a,b,c,e,i,j\mid a^2=b^2=c^2=e^2=i^2=j^6=1,
~a,b,c,e \text{ commute, } iai=b,\\
ici=ae,\medspace jaj^{-1}=c,\medspace jbj^{-1}=abc,\medspace jcj^{-1}=be,
\medspace
j^3=abce,\medspace (ji)^2=bc\rangle.
\end{gather*} 
The generators $a,b,c,$ and $j$ represent orientation reversing 
automorphisms.
(Note that $jej^{-1}=bc$ follows from the other relations. 
See \cite{[Zn90]} for an alternative description.)
       
The group $B_1=\mathbb{Z}\times\! K$ has a presentation 
\begin{equation*} 
\langle t,x,y\mid tx=xt,\medspace ty=yt,\medspace xyx^{-1}=y^{-1}\rangle.
\end{equation*} 
An automorphism of $B_1$ must preserve the centre $\zeta B_1$ 
(which has basis $t,x^2$) and $I(B_1)$ (which is generated by $y$).
Thus the automorphisms of $B_1$ may be identified with triples 
$(A,m,\epsilon)\in \Upsilon_2\times\mathbb{Z}\times\{\pm1\}$, 
where $\Upsilon_2 $ is the subgroup 
of $GL(2,\mathbb{Z})$ consisting of matrices congruent $mod~(2)$ to upper triangular matrices.
Such an automorphism sends $t$ to $t^ax^b$, $x$ to $t^cx^dy^m$ and $y$ to $y^\epsilon$,
and induces multiplication by $A$ on $B_1 /I(B_1)\cong\mathbb{Z}^2$.
Composition of automorphisms is given by
$(A,m,\epsilon)(B,n,\eta)=(AB,m+\epsilon n,\epsilon\eta)$.
The inner automorphisms are generated by $(I,1,-1)$ and $(I,2,1)$, and so
$Out(B_1)\cong \Upsilon_2\times (Z/2Z)$.
 
The group $B_2$ has a presentation  
\begin{equation*} 
\langle t,x,y\mid txt^{-1}=xy,\medspace ty=yt,\medspace xyx^{-1}=y^{-1}\rangle.
\end{equation*} 
Automorphisms of $B_2$ may be identified with triples 
$(A,(m,n),\epsilon)$, where $A\in\Upsilon_2$, 
$m, n\in\mathbb{Z}$, $\epsilon=\pm1$
and $m=(A_{11}-\epsilon)/2$.
Such an automorphism sends $t$ to $t^ax^by^m$, $x$ to $t^cx^dy^n$ and $y$ to 
$y^\epsilon$,
and induces multiplication by $A$ on $B_2 /I(B_2)\cong\mathbb{Z}^2$.
The automorphisms which induce the identity on $B_2 /I(B_2)$ are all
inner, and so $Out(B_2)\cong\Upsilon_2$.

The group $B_3 $ has a presentation
\begin{equation*} 
\langle t,x,y\mid txt^{-1}=x^{-1},\medspace ty=yt,
\medspace xyx^{-1}=y^{-1}\rangle.
\end{equation*} 
An automorphism of $B_3 $ must preserve $I(B_3)\cong K$ 
(which is generated by $x,y$) and $I(I(B_3))$ (which is generated by $y$).
It follows easily that $Out(B_3)\cong (Z/2Z)^3$, and is generated 
by the classes of the automorphisms 
which fix $y$ and send $t$ to $t^{-1}, 
t, tx^2$ and $x$ to $x,xy,x$,
respectively.                                           

A similar argument using the presentation
\begin{equation*} 
\langle t,x,y\mid txt^{-1}=x^{-1}y,\medspace ty=yt,
\medspace xyx^{-1}=y^{-1}\rangle
\end{equation*} 
for $B_4$ shows that $Out(B_4)\cong (Z/2Z)^3$, and is generated by 
the classes of the automorphisms
which fix $y$ and send $t$ to $t^{-1}y^{-1}, 
t, tx^2$ and $x$ to $x,x^{-1},x$,
respectively.                                           

\section{Flat 4-manifold groups with infinite abelianization}

We shall organize our determination of the flat 4-manifold groups 
$\pi$ in terms of $I(\pi)$.
Let $\pi$ be a flat 4-manifold group, $\beta=\beta_1 (\pi)$ and $h=h(I(\pi))$.
Then $\pi/I(\pi)\cong\mathbb{Z}^\beta$ and $h+\beta=4$. 
If $I(\pi)$ is abelian then $C_\pi (I(\pi))$ is a nilpotent normal subgroup of $\pi$
and so is a subgroup of the Hirsch-Plotkin radical $\sqrt\pi$,
which is here the maximal abelian normal subgroup $T(\pi)$.
Hence $C_\pi (I(\pi))=T(\pi)$ and the holonomy group is isomorphic to $\pi/C_\pi(I(\pi))$.

\medskip
$h=0$\qua In this case $I(\pi)=1$, so $\pi\cong\mathbb{Z}^4$ and is orientable. 

\medskip
$h=1$\qua In this case $I(\pi)\cong\mathbb{Z}$ and $\pi$ is nonabelian, 
so $\pi/C_\pi (I(\pi))=Z/2Z$.
Hence $\pi$ has a presentation of the form 
\begin{equation*} 
\langle t,x,y,z\mid txt^{-1}=xz^a,\medspace tyt^{-1}=yz^b,
\medspace tzt^{-1}=z^{-1},~x,y,z~commute\rangle,
\end{equation*} 
for some integers $a$, $b$.
On replacing $x$ by $xy$ or interchanging $x$ and $y$ if necessary we may assume that $a$ is even. 
On then replacing $x$ by $xz^{a/2}$ and $y$ by $yz^{[b/2]}$ we may assume that $a=0$
and $b=0$ or 1. 
Thus $\pi$ is a semidirect product $\mathbb{Z}^3\rtimes_T\mathbb{Z}$, 
where the normal subgroup $\mathbb{Z}^3$ is generated 
by the images of $x$, $y$ and $z$, and the action of $t$ is determined
by a matrix $T=
\left(\begin{smallmatrix}
I_2&0\\
(0,b)&-1
\end{smallmatrix}\right)$ in $GL(3,\mathbb{Z})$.
Hence $\pi\cong\mathbb{Z}\times B_1 =\mathbb{Z}^2\times K$ 
or $\mathbb{Z}\times B_2$.
Both of these groups are nonorientable.

\medskip  
$h=2$\qua If $I(\pi)\cong\mathbb{Z}^2$ and $\pi/C_\pi (I(\pi))$ 
is cyclic then we may again assume that $\pi$ is 
a semidirect product $\mathbb{Z}^3\rtimes_T\mathbb{Z}$, 
where $T=\left(\begin{smallmatrix}
1&0\\
\mu& U
\end{smallmatrix}\right)$, 
with $\mu=\left(\begin{smallmatrix}
a\\
b
\end{smallmatrix}\right)$ 
and $U\in GL(2,\mathbb{Z})$ is of order 2, 3, 4 or 6 
and does not have 1 as an eigenvalue.
Thus $U=-I_2$, 
$\left(\begin{smallmatrix}
0&-1\\
1&-1
\end{smallmatrix}\right)$, 
$\left(\begin{smallmatrix}
0&-1\\
1&0
\end{smallmatrix}\right)$
or $\left(\begin{smallmatrix}
0&-1\\
1&1
\end{smallmatrix}\right)$.                                                                   
Conjugating $T$ by 
$\left(\begin{smallmatrix}
1&0\\
\nu& I_2
\end{smallmatrix}\right)$ replaces $\mu$ by $\mu+(I_2 -U)\nu$.
In each case the choice $a=b=0$ leads to a group 
of the form $\pi\cong\mathbb{Z}\times G$,
where $G$ is an orientable flat 3-manifold group with $\beta_1 (G)=1$.
For each of the first three of these matrices there is one other possible group.
However if $U=\left(\begin{smallmatrix}
0&-1\\
1&1
\end{smallmatrix}\right)$ then $I_2 -U$ is invertible and so
$\mathbb{Z}\times G_5$ is the only possibility. 
All seven of these groups are orientable.

If $I(\pi)\cong\mathbb{Z}^2$ and $\pi/C_\pi (I(\pi))$ is not cyclic then $\pi/C_\pi(I(\pi))\cong(Z/2Z)^2$.
There are two conjugacy classes of embeddings of $(Z/2Z)^2$ in $GL(2,\mathbb{Z})$.
One has image the subgroup of diagonal matrices. 
The corresponding groups $\pi$ have presentations of the form
\begin{gather*} 
\langle t,u,x,y\mid tx=xt,\medspace tyt^{-1}=y^{-1},\medspace
uxu^{-1}=x^{-1},\medspace uyu^{-1}=y^{-1},\medspace xy=yx,\\
tut^{-1}u^{-1}=x^my^n\rangle,
\end{gather*} 
for some integers $m$, $n$. 
On replacing $t$ by $tx^{-[m/2]}y^{[n/2]}$ if necessary we may assume that $0\leq m,n\leq 1$. 
On then replacing $t$ by $tu$ and interchanging $x$ and $y$ if necessary we may assume that $m\leq n$.
The only infinite cyclic subgroups of $I(\pi)$ which are normal in $\pi$ are the subgroups 
$\langle x\rangle$ and $\langle y\rangle$. 
On comparing the quotients of these groups $\pi$ by such subgroups
we see that the three possibilities are distinct.
The other embedding of $(Z/2Z)^2$ in $GL(2,\mathbb{Z})$ has image generated by
$-I$ and $\left(\begin{smallmatrix}
0&1\\
1&0
\end{smallmatrix}\right)$.
The corresponding groups $\pi$ have presentations of the form
\begin{gather*} 
\langle t,u,x,y\mid txt^{-1}=y,\medspace tyt^{-1}=x,\medspace
uxu^{-1}=x^{-1},\medspace uyu^{-1}=y^{-1},\medspace xy=yx,\\
tut^{-1}u^{-1}=x^my^n\rangle,
\end{gather*} 
for some integers $m$, $n$.                                                      On replacing $t$ by $tx^{[(m-n)/2]}$ and $u$ by $ux^{-m}$, 
if necessary, we may assume that $m=0$ and $n=0$ or 1. 
Thus there two such groups.
All five of these groups are nonorientable.

Otherwise, $I(\pi)\cong K$, $I(I(\pi))\cong\mathbb{Z}$ 
and $G=\pi/I(I(\pi))$ is a flat 3-manifold group 
with $\beta_1(G)=2$, but with $I(G)=I(\pi)/I(I(\pi))$ not contained in $G'$ 
(since it acts nontrivially on $I(I(\pi))$).
Therefore $G\cong B_1=\mathbb{Z}\times K$, 
and so has a presentation
\begin{equation*}
\langle t,x,y\mid tx=xt,\medspace ty=yt,\medspace xyx^{-1}=y^{-1}\rangle.
\end{equation*}                                
If $w:G\to Aut(\mathbb{Z})$ is a homomorphism which restricts nontrivially 
to $I(G)$ then we may assume (up to isomorphism of $G$) that $w(x)=1$ 
and $w(y)=-1$.
Groups $\pi$ which are extensions of $\mathbb{Z}\times K$ by $\mathbb{Z}$
corresponding to the action with $w(t)=w$ ($=\pm 1$) 
have presentations of the form
\begin{gather*} 
\langle t,x,y,z\mid txt^{-1}=xz^a,\medspace tyt^{-1}=yz^b,
\medspace tzt^{-1}=z^w,\medspace xyx^{-1}=y^{-1} z^c,\medspace xz=zx,\\
yzy^{-1}=z^{-1}\rangle,
\end{gather*} 
for some integers $a,b$ and $c$.
Any group with such a presentation is easily seen to be an extension 
of $\mathbb{Z}\times K$ by a cyclic normal subgroup. 
However conjugating the fourth relation leads to the equation
\begin{equation*} 
txt^{-1}tyt^{-1}(txt^{-1})^{-1}=txyx^{-1}t^{-1}=ty^{-1}z^ct^{-1}=
tyt^{-1}(tzt^{-1})^c
\end{equation*} 
which simplifies to $xz^a yz^b z^{-a}x^{-1}=(yz^b)^{-1}z^{wc}$ and 
hence to $z^{c-2a}=z^{wc}$.
Hence this cyclic normal subgroup is finite unless $2a=(1-w)c$.

Suppose first that $w=1$. Then $z^{2a}=1$ and so we must have $a=0$.
On replacing $t$ by $tz^{[b/2]}$ and $x$ by $xz^{[c/2]}$, if necessary, 
we may assume that $0\leq b,c\leq 1$.
If $b=0$ then $\pi\cong\mathbb{Z}\times B_3$ or $\mathbb{Z}\times B_4$.
Otherwise, after further replacing $x$ by $txz$,
if necessary, we may assume that $b=1$ and $c=0$.
The three possibilities may be distinguished by their abelianizations,
and so there are three such groups. 
In each case the subgroup generated by
$\{ t,x^2,y^2,z\}$ is maximal abelian, 
and the holonomy group is isomorphic to $(Z/2Z)^2$.
            
If instead $w=-1$ then $z^{2(c-a)}=1$ and so we must have $c=a$.
On replacing $x$ by $xz^{[a/2]}$ and $y$ by $yz^{[b/2]}$,  
if necessary, 
we may assume that $0\leq{a,b}\leq 1$.
If $b=1$ then after replacing $x$ by $txy$, if necessary, 
we may assume that $a=0$.
If $a=b=0$ then $\pi/\pi'\cong\mathbb{Z}^2\oplus (Z/2Z)^2$.
The two other possibilities each have abelianization 
$\mathbb{Z}^2\oplus (Z/2Z)$,
but one has centre of rank 2 and the other has centre of rank 1.
Thus there are three such groups.
The subgroup generated by $\{ ty,x^2,y^2,z\}$ is maximal abelian, 
and the holonomy group is isomorphic to $(Z/2Z)^2$.
All of these groups $\pi$ with $I(\pi)\cong K$ are nonorientable.

\medskip
$h=3$\qua In this case $\pi$ is uniquely a semidirect product 
$\pi\cong I(\pi)\rtimes_\theta\mathbb{Z}$, 
where $I(\pi)$ is a flat 3-manifold group and $\theta$ 
is an automorphism of $I(\pi)$ such that the induced automorphism 
of $I(\pi)/I(I(\pi))$ has no eigenvalue 1, 
and whose image in $Out(I(\pi))$ has finite order. 
(The conjugacy class of the image of $\theta$ in $Out(I(\pi))$ 
is determined up to inversion by $\pi$.)
         
Since $T(I(\pi))$ is the maximal abelian normal subgroup of $I(\pi)$ 
it is normal in $\pi$.
It follows easily that $T(\pi)\cap I(\pi)=T(I(\pi))$.
Hence the holonomy group of $I(\pi)$ is isomorphic to a normal subgroup of
the holonomy subgroup of $\pi$, with quotient cyclic of order dividing
the order of $\theta$ in $Out(I(\pi))$. 
(The order of the quotient can be strictly smaller.)

If $I(\pi)\cong\mathbb{Z}^3$ then $Out(I(\pi))\cong GL(3,\mathbb{Z})$. 
If $T\in GL(3,\mathbb{Z})$ has finite order $n$ 
and $\beta_1 (\mathbb{Z}^3\rtimes_T\mathbb{Z})=1$ 
then either $T=-I$ or $n=4$ or 6 and the characteristic polynomial of $T$ 
is $(t+1)\phi(t)$ with $\phi(t)=t^2+1$, $t^2+t+1$ or $t^2-t+1$.
In the latter cases $T$ is conjugate to a matrix of the form 
$\left(\begin{smallmatrix}
-1&\mu\\
0& A
\end{smallmatrix}\right)$, where 
$A=\left(\begin{smallmatrix}
0&-1\\
1&0
\end{smallmatrix}\right)$,
$\left(\begin{smallmatrix}
0&-1\\
1&-1
\end{smallmatrix}\right)$ or 
$\left(\begin{smallmatrix}
0&-1\\
1&1
\end{smallmatrix}\right)$, respectively.
The row vector $\mu=(m_1, m_2)$ is well defined $mod~\mathbb{Z}^2(A+I)$.
Thus there are seven such conjugacy classes.
All but one pair (corresponding to 
$\left(\begin{smallmatrix}
0&-1\\
1&1
\end{smallmatrix}\right)$ 
and $\mu\notin\mathbb{Z}^2(A+I)$) are self-inverse, 
and so there are six such groups. 
The holonomy group is cyclic, of order equal to the order of $T$.
As such matrices all have determinant $-1$ all of these groups are nonorientable.     

If $I(\pi)\cong G_i$ for $2\leq i\leq 5$ the automorphism
$\theta=(v,A,\epsilon)$ must have $\epsilon=-1$, for otherwise $\beta_1 (\pi)=2$.
We have $Out(G_2)\cong ((Z/2Z)^2\rtimes PGL(2,\mathbb{Z}))\times (Z/2Z)$.
The five conjugacy classes of finite order in $PGL(2,\mathbb{Z})$ are
represented by the matrices $I$, 
$\left(\begin{smallmatrix}
0&-1\\
1&0
\end{smallmatrix}\right)$, 
$\left(\begin{smallmatrix}
0&1\\
1&0
\end{smallmatrix}\right)$, 
$\left(\begin{smallmatrix}
1&0\\
0&-1
\end{smallmatrix}\right)$
and $\left(\begin{smallmatrix}
0&1\\
-1&1
\end{smallmatrix}\right)$. 
The numbers of conjugacy classes in $Out(G_2)$ with $\epsilon=-1$ corresponding to these
matrices are two, two, two, three and one, respectively.
All of these conjugacy classes are self-inverse.
Of these, only the two conjugacy classes corresponding to 
$\left(\begin{smallmatrix}
0&1\\
1&0
\end{smallmatrix}\right)$  
and the three conjugacy classes corresponding to 
$\left(\begin{smallmatrix}
1&0\\
0&-1
\end{smallmatrix}\right)$ 
give rise to orientable groups.                     
The holonomy groups are all isomorphic to $(Z/2Z)^2$, except when
$A=\left(\begin{smallmatrix}
0&-1\\
1&0
\end{smallmatrix}\right)$ or 
$\left(\begin{smallmatrix}
0&1\\
-1&1
\end{smallmatrix}\right)$, 
when they are isomorphic to $Z/4Z$ or $Z/6Z\oplus Z/2Z$, respectively.
There are five orientable groups and five nonorientable groups.

As $Out(G_3)\cong S_3\times(Z/2Z)$, 
$Out(G_4)\cong (Z/2Z)^2$ and $Out(G_5)=Z/2Z $,
there are three, two and one conjugacy classes corresponding to 
automorphisms with $\epsilon=-1$, respectively, 
and all these conjugacy classes are closed under inversion.          
The holonomy groups are dihedral of order 6, 8 and 12, respectively.
The six such groups are all orientable.

The centre of $Out(G_6)$ is generated by the image of $ab$, and the image of
$ce$ in the quotient $Out(G_6)/\langle ab\rangle $ generates a central $Z/2Z$ direct factor.
The quotient $Out(G_6)/\langle ab,ce\rangle $ is isomorphic to the semidirect product of
a normal subgroup $(Z/2Z)^2$ (generated by the images of $a$ and $c$) with
$S_3$ (generated by the images of $ia$ and $j$), and has five conjugacy classes,
represented by $1,a,i,j$ and $ci$.
Hence $Out(G_6)/\langle ab\rangle $ has ten conjugacy classes, 
represented by $1,a,ce,e=iacei,i,cei,j,cej,ci$ and $cice=ei$.
Thus $Out(G_6)$ itself has between 10 and 20 conjugacy classes.
In fact $Out(G_6)$ has 14 conjugacy classes, of which those represented by
$1,ab,bce,e,i,cej$, $abcej$ and $ei$ are orientation preserving, 
and those represented by
$a,ce,cei,j,abj$ and $ci$ are orientation reversing.
All of these classes are self inverse, except for $j$ and $abj$, 
which are mutually inverse ($j^{-1}=ai(abj)ia$). 
The holonomy groups corresponding to the classes $1,ab,bce$ and $e$ 
are isomorphic to $(Z/2Z)^2$,
those corresponding to $a$ and $ce$ are isomorphic to $(Z/2Z)^3$,
those corresponding to $i,ei,cei$ and $ci$ are dihedral of order 8,
those corresponding to $cej$ and $abcej$ are isomorphic to $A_4$ 
and the one corresponding to $j$ is isomorphic to 
$(Z/2Z)^2\rtimes{Z/6Z}\cong{A_4\times{Z/2Z}}$.
There are eight orientable groups and five nonorientable groups.
        
All the remaining cases give rise to nonorientable groups.

\noindent $I(\pi)\cong\mathbb{Z}\times K$.  
If a matrix $A$ in $\Upsilon_2$ has finite order then as its trace is even
the order must be 1, 2 or 4. If moreover $A$ does not have 1 as an eigenvalue then
either $A=-I$ or $A$ has order 4 and is conjugate (in $\Upsilon_2 $) to
$\left(\begin{smallmatrix}
-1&1\\
-2&1
\end{smallmatrix}\right)$. 
Each of the four corresponding conjugacy classes in $\Upsilon_2\times{Z/2Z}$
is self inverse, and so there are four such groups.  
The holonomy groups are isomorphic to $Z/nZ\oplus{Z/2Z}$, where $n=2$ or 4 is the order of $A$.

\noindent $I(\pi)\cong B_2$. As $Out(B_2)\cong \Upsilon_2$ there are two relevant
conjugacy classes and hence two such groups.                                     
The holonomy groups are again isomorphic to ${Z/nZ\oplus{Z/2Z}}$, 
where $n=2$ or 4 is the order of $A$.

\noindent $I(\pi)\cong B_3$ or $B_4 $. 
In each case $Out(H)\cong(Z/2Z)^3$, 
and there are four outer automorphism classes determining 
semidirect products with $\beta=1$.
(Note that here conjugacy classes are singletons and are self-inverse.) 
The holonomy groups are all isomorphic to $(Z/2Z)^3$.

\section{Flat 4-manifold groups with finite abelianization}

There remains the case when $\pi/\pi'$ is finite (equivalently, $h=4$).
By Lemma 3.14 if $\pi$ is such a flat 4-manifold group it is nonorientable and 
is isomorphic to a generalized free product $J*_\phi \tilde J$,
where $\phi$ is an isomorphism from $G<J$ to $\tilde G<\tilde J$
and $[J:G]=[\tilde J:\tilde G]=2$.
The groups $G$, $J$ and $\tilde J$ are then flat 3-manifold groups.
If $\lambda$ and $\tilde\lambda$ are automorphisms of $G$ and $\tilde G$ 
which extend to $J$ and $\tilde J$, respectively, then $J*_\phi\tilde J$ and
$J*_{\tilde\lambda\phi\lambda}\tilde J$ are isomorphic, 
and so we shall say that 
$\phi$ and $\tilde\lambda\phi\lambda$ are equivalent.
The major difficulty is that some such groups 
split as a generalised free product in several essentially distinct ways.

It follows from the Mayer-Vietoris sequence 
for $\pi\cong J*_\phi \tilde J$ 
that $H_1 (G;\mathbb{Q})$ maps onto 
$H_1 (J;\mathbb{Q})\oplus H_1 (\tilde J;\mathbb{Q})$, and hence that 
$\beta_1 (J)+\beta_1 (\tilde J)\leq \beta_1 (G)$.
Since $G_3 $, $G_4 $, $B_3 $ and $B_4 $ are only subgroups of other
flat 3-manifold groups via maps inducing isomorphisms on $H_1 (-;\mathbb{Q})$
and $G_5$ and $G_6$ are not index 2 subgroups of any flat 3-manifold group,
we may assume that $G\cong\mathbb{Z}^3$, $G_2$, $B_1$ or $B_2$.
If $G$ is orientable at least one of $J$ and $\tilde{J}$ is orientable,
for otherwise $J*_\phi\tilde{J}$ is orientable.
If $j$ and $\tilde j$ are the automorphisms of $T(J)$ and $T(\tilde J)$ 
determined by conjugation in $J$ and $\tilde J$, respectively, 
then $\pi$ is a flat 4-manifold 
group if and only if $\Phi=jT(\phi)^{-1}\tilde jT(\phi)$ has finite order. 
In particular, $|tr(\Phi)|\leq3$.
At this point detailed computation seems unavoidable.                 
(We note in passing that any generalised free product $J*_G\tilde J$ with 
$G\cong G_3$, $G_4$, $B_3$ or $B_4$, $J$ and $\tilde J$ torsion-free and 
$[J:G]=[\tilde J:G]=2$ is a flat 4-manifold group, since $Out(G)$ is then finite. 
However all such groups have infinite abelianization.)
    
Suppose first that $G\cong\mathbb{Z}^3$, with basis $\{ x,y,z\}$.
Then $J$ and $\tilde J$ must have holonomy of order $\leq 2$, and 
$\beta_1(J)+\beta_1(\tilde J)\leq 3$. 
Hence we may assume that $J\cong G_2$ and $\tilde J\cong G_2$, $B_1$ or $B_2$.
In each case we have $G=T(J)$ and $\tilde G=T(\tilde J)$.
We may assume that $J$ and $\tilde J$ are generated by $G$ and elements $s$ and $t$, respectively,
such that $s^2=x$ and $t^2\in \tilde G$,
and that the action of $s$ on $G$ has matrix 
$j=\left(\begin{smallmatrix}
1&0\\
0& -I
\end{smallmatrix}\right)$ 
with respect to the basis $\{ x,y,z\}$.
Fix an isomorphism $\phi:G\to\tilde G$ and let 
$T=T(\phi)^{-1}\tilde jT(\phi)=
\left(\begin{smallmatrix}
a&\delta\\
\gamma& D
\end{smallmatrix}\right)$ 
be the matrix corresponding to the action of $t$ on $\tilde G$.  
(Here $\gamma$ is a $2\times 1$ column vector, $\delta$ is a $1\times 2$ row vector
and $D$ is a $2\times 2$ matrix, possibly singular.)
Then $T^2=I$ and so the trace of $T$ is odd.            
Since $j\equiv I~mod~(2)$ the trace of $\Phi=jT$ is also odd, and so
$\Phi$ cannot have order 3 or 6. Therefore $\Phi^4=I$.
If $\Phi=I$ then $\pi/\pi'$ is infinite.
If $\Phi$ has order 2 then $jT=Tj$ and so $\gamma=0$, $\delta=0$ and $D^2=I_2$. 
Moreover we must have $a=-1$ for otherwise $\pi/\pi'$ is infinite.
After conjugating $T$ by a matrix commuting with $j$ if necessary we may assume
that $D=I_2$ or 
$\left(\begin{smallmatrix}
1&0\\
0&-1
\end{smallmatrix}\right)$.
(Since $\tilde J$ must be torsion-free we cannot have 
$D=\left(\begin{smallmatrix}
0&1\\
1&0
\end{smallmatrix}\right)$.)
These two matrices correspond to the generalized free products 
$G_2 *_\phi B_1$ and
$G_2 *_\phi G_2$, with presentations 
\begin{gather*} 
\langle s,t,z\mid st^2s^{-1} =t^{-2},\medspace szs^{-1}=z^{-1},
\medspace ts^2t^{-1}=s^{-2},\medspace tz=zt\rangle\\
\text{and}\quad\langle s,t,z\mid st^2s^{-1} =t^{-2},\medspace szs^{-1}=z^{-1},
\medspace ts^2t^{-1}=s^{-2},\medspace tzt^{-1}=z^{-1}\rangle,
\end{gather*} 
respectively. 
These groups each have holonomy group isomorphic to $(Z/2Z)^2$.
If $\Phi$ has order 4 then we must have $(jT)^2=(jT)^{-2} =(Tj)^2$ and so 
$(jT)^2$ commutes with $j$. 
After conjugating $T$ by a matrix commuting with $j$, 
if necessary, 
we may assume that $T$ is the elementary matrix which interchanges 
the first and third rows.
The corresponding group $G_2 *_\phi{B_2}$ has a presentation
\begin{equation*} 
\langle s,t,z\mid st^2s^{-1}=t^{-2},\medspace 
{szs^{-1}=z^{-1}},\medspace {ts^2t^{-1}=z},\medspace{tzt^{-1}=s^2}\rangle.
\end{equation*}
(The final relation is redundant, and so $G_2 *_\phi{B_2}$ has deficiency 0.) 
Its holonomy group is isomorphic to the dihedral group of order 8.

If $G\cong G_2 $ then $\beta_1(J)+\beta_1(\tilde J)\leq 1$, so we may assume that $J\cong G_6$.
The other factor $\tilde J$ must then be one of $G_2 $, 
$G_4 $, $G_6 $, $B_3 $ or $B_4$, and then
every amalgamation has finite abelianization. 
If $\tilde{J}\cong{G_2}$ there are two index-2 embeddings of $G_2$ 
in $\tilde{J}$ up to composition with an automorphism of $\tilde{J}$.
(One of these was overlooked in earlier versions of this book.)
In all other cases the image of $G_2$ in $\tilde{J}$ is canonical,
and the matrices for $j$ and $\tilde j$ have the form
$\left(\begin{smallmatrix}
\pm1&0\\
0& N
\end{smallmatrix}\right)$ where $N^4 =I\in GL(2,\mathbb{Z})$,
and 
$T(\phi)=\left(\begin{smallmatrix}
\epsilon&0\\
0& M
\end{smallmatrix}\right)$ for some $M\in GL(2,\mathbb{Z})$.
Calculation shows that $\Phi$ has finite order if and only if $M$ is in the 
dihedral subgroup $D_8$ of $GL(2,\mathbb{Z})$ generated by the diagonal 
matrices and 
$\left(\begin{smallmatrix}
0&1\\
1&0
\end{smallmatrix}\right)$.
(In other words, either $M$ is diagonal or both diagonal elements of $M$ are 0.)
Now the subgroup of $Aut(G_2)$ consisting of automorphisms which extend to
$G_6$ is $(\mathbb{Z}^2\rtimes_\alpha D_8) \times\{\pm1\}$.
Hence any two such isomorphisms $\phi$ from $G$ to $\tilde G$ are equivalent,
and so there is an unique such flat 4-manifold group $G_6*_\phi\tilde J$ 
for each of these choices of $\tilde J=G_4 $, $G_6 $, $B_3 $ or $B_4$.
The corresponding presentations are
\begin{gather*} 
\langle u,x,y\mid xux^{-1}=u^{-1},\medspace y^2=u^2,
\medspace yx^2y^{-1}=x^{-2},\medspace u(xy)^2=(xy)^2 u\rangle,\\
\langle u,x,y\mid xux^{-1}=u^{-1},
\medspace xy^2x^{-1}=y^{-2},
\medspace y^2=u^2(xy)^2,
\medspace yx^2y^{-1}=x^{-2},\\
 u(xy)^2=(xy)^2 u\rangle,\\
\langle u,x,y\mid yx^2y^{-1}=x^{-2},\medspace uy^2u^{-1}=(xy)^2,
\medspace u(xy)^2u^{-1}=y^{-2},\medspace x=u^2\rangle,\\
\langle u,x,y\mid xy^2x^{-1}=y^{-2},\medspace 
yx^2y^{-1}=ux^2u^{-1}=x^{-2},\medspace y^2=u^2,\medspace
yxy=uxu\rangle,\\
\langle t,x,y\mid xy^2x^{-1}=y^{-2},\medspace yx^2y^{-1}=x^{-2},
\medspace x^2=t^2,\medspace y^2=(t^{-1}x)^2,\medspace
t(xy)^2=(xy)^2 t\rangle,\\
\text{and}\quad
\langle t,x,y\mid xy^2x^{-1}=y^{-2},\medspace 
yx^2y^{-1}=x^{-2},\medspace x^2=t^2(xy)^2,\medspace y^2=(t^{-1}x)^2,\\
t(xy)^2=(xy)^2t\rangle,
\end{gather*}
respectively.
The corresponding holonomy groups are isomorphic to $(Z/2Z)^3$, 
$(Z/2Z)^3$, $D_8$, $(Z/2Z)^2$,
$(Z/2Z)^3$ and $(Z/2Z)^3$, respectively. 

If $G\cong B_1$ or $B_2$ then $J$ and $\tilde J$ are nonorientable and 
$\beta_1(J)+\beta_1(\tilde J)\leq 2$.
Hence $J$ and $\tilde J$ are $B_3 $ or $B_4 $.
There are two essentially different embeddings of $B_1$ 
as an index 2 subgroup in each of $B_3$ or $B_4$. 
(The image of one contains $I(B_i)$ while the other does not.) 
The group $B_2$ is not a subgroup of $B_4$. 
However, it embeds in  $B_3$ as the index 2 subgroup generated by $\{t,x^2,xy\}$.
Contrary to our claim in earlier versions of this book that 
no flat 4-manifold groups with finite abelianization
are such amalgamations, there are in fact three, all with $G=B_1$.
However they are each isomorphic to one of the groups already given: $B_3*_\phi{B_3}\cong{B_3*_\phi{B_4}}\cong{G_6*_\phi{B_4}}$ 
and $B_4*_\phi{B_4}\cong{G_6*_\phi{B_3}}$.
See \cite{[LRT13]}.
                                                                       
There remain nine generalized free products $J*_G \tilde J$ which are 
flat 4-manifold groups with $\beta=0$ and with $G$ orientable.
The groups $G_2 *_\phi B_1$, $G_2 *_\phi G_2$ and $G_6*_\phi G_6$ 
are all easily seen to be semidirect products of $G_6$ 
with an infinite cyclic normal subgroup, on which $G_6$ acts nontrivially.
It follows easily that these three groups are in fact isomorphic, 
and so there is just one flat 4-manifold group with finite abelianization 
and holonomy isomorphic to $(Z/2Z)^2$.

The groups $G_2*_\phi{B_2}$ and $G_6*_\phi{G_4}$ are in fact isomorphic;
the function sending $s$ to $y$, $t$ to $yu^{-1}$ and $z$ to $uy^2u^{-1}$
determines an equivalence between the above presentations.
Thus there is just one flat 4-manifold group with finite abelianization and
holonomy isomorphic to $D_8$.

The first amalgamation  $G_6*_\phi G_2$ and $G_6*_\phi B_4$ are also isomorphic,
via the function sending $u$ to $xy^{-1}$, $x$ to $xt^{-1}$ and $y$ to $yt$.
Similarly, the second amalgamation  $G_6*_\phi G_2$ and $G_6*_\phi B_3$ are isomorphic;
via the function sending $u$ to $tyx$, $x$ to $tx^{-1}$ and $y$ to $ty$.
(These isomorphisms and the one in the paragraph above 
were found by Derek Holt, 
using the program described in \cite{[HR92]}.)
The translation subgroups of $G_6*_\phi B_3$ and $G_6*_\phi B_4$
are generated by the images of $U=(ty)^2$, $X=x^2$, $Y=y^2$ and
$Z=(xy)^2$ (with respect to the above presentations).
In each case the images of $t$, $x$ and $y$ act diagonally, via the
matrices $diag[-1,1,-1,1]$, $diag[1,1,-1,-1]$ and $diag[-1,-1,1,-1]$, 
respectively.
However the maximal orientable subgroups
have abelianization $\mathbb{Z}\oplus (Z/2)^3$ and
$\mathbb{Z}\oplus(Z/4Z)\oplus(Z/2Z)$, respectively, and so
$G_6*_\phi B_3$ is not isomorphic to $G_6*_\phi B_4$.
Thus there are two flat 4-manifold groups with finite abelianization
and holonomy isomorphic to $(Z/2Z)^3$.

In summary, there are 27 orientable flat 4-manifold groups 
(all with $\beta>0$), 43 nonorientable flat 4-manifold groups with $\beta>0$ 
and 4 (nonorientable) flat 4-manifold groups with $\beta=0$.
All orientable flat 4-manifolds are Spin,
excepting for those with $\pi\cong{G_6}\rtimes_\theta\mathbb{Z}$,
where $\theta=bce, e$ or $ei$ \cite{[PS10]}.

\section{Distinguishing between the geometries}

If $\Gamma$ is a lattice in a 1-connected solvable Lie group $G$ 
with nilradical $N$ then $\Gamma\cap{N}$ and $\Gamma\cap{N}'$ are
lattices in $N$ and $N'$, respectively, and $h(\Gamma)=dim(G)$.
If $G$ is also a linear algebraic group then $\Gamma$ is Zariski dense in $G$.
In particular, $\Gamma\cap{N}$ and $N$ have the same nilpotency class,
and $\zeta\Gamma=\Gamma\cap\zeta{N}$.
(See \cite[Chapter 2]{[Rg]}.)
These observations imply that the geometry of a closed 4-manifold $M$ 
with a geometry of solvable Lie type is largely determined by 
the structure of $\sqrt\pi$.
(See also \cite[Proposition 10.4]{[Wl86]}.)
As each covering space has the same geometry 
it shall suffice to show that the geometries 
on suitable finite covering spaces (corresponding to subgroups 
of finite index in $\pi$) can be recognized.

If $M$ is an infranilmanifold then $[\pi:\sqrt\pi]<\infty$.
If it is flat then $\sqrt\pi\cong\mathbb{Z}^4$, while if it has the geometry
$\mathbb{N}il^3\times\mathbb{E}^1$ or $\mathbb{N}il^4$ then 
$\sqrt\pi$ is nilpotent of class 2 or 3 respectively.
(These cases may also be distinguished by the rank of $\zeta\sqrt\pi$.) 
All such groups have been classified, and may be realized geometrically.
(See \cite{[De]} for explicit representations of the 
$\mathbb{N}il^3\times\mathbb{E}^1$- and $\mathbb{N}il^4$-groups as
lattices in $Aff(Nil^3\times\mathbb{R})$ and $Aff(Nil^4)$, respectively.)
If $M$ is a $\mathbb{S}ol^4_{m,n}$- or $\mathbb{S}ol^4_0$-manifold then 
$\sqrt\pi\cong Z^3$.
Hence $h(\pi/\sqrt\pi)=1$ and so $\pi$ has a normal subgroup of finite 
index which is a semidirect product $\sqrt\pi\rtimes_\theta\mathbb{Z}$.
It is easy to give geometric realizations of such subgroups.

\begin{theorem} 
Let $\pi$ be a torsion-free group with $\sqrt\pi\cong\mathbb{Z}^3$ 
and such that $\pi/\sqrt\pi\cong\mathbb{Z}$.
Then $\pi$ is the fundamental group of a $\mathbb{S}ol^4_{m,n}$- or
$\mathbb{S}ol^4_0$-manifold.
\end{theorem}

\begin{proof}
Let $t\in\pi$ represent a generator of $\pi/\sqrt\pi$,
and let $\theta$ be the automorphism of $\sqrt\pi\cong\mathbb{Z}^3$
determined by conjugation by $t$.
Then $\pi\cong\sqrt\pi\rtimes_\theta\mathbb{Z}$.
If the eigenvalues of $\theta$ were roots of unity of order dividing $k$
then the subgroup generated by $\sqrt\pi$ and $t^k$ would be nilpotent,
and of finite index in $\pi$.
Therefore we may assume that the eigenvalues 
$\kappa,\lambda,\mu$ of $\theta$ are distinct and that 
neither $\kappa$ nor $\lambda$ is a root of unity.

Suppose first that the eigenvalues are all real.
Then the eigenvalues of $\theta^2$ are all positive,
and $\theta^2$ has characteristic polynomial $X^3 -mX^2 +nX-1$, 
where $m=trace(\theta^2)$ and $n=trace(\theta^{-2} )$.                                            
Since $\sqrt\pi\cong\mathbb{Z}^3$ there is a monomorphism $f:\sqrt\pi\to\mathbb{R}^3$ 
such that $f\theta=\Psi{f}$, 
where $\Psi=diag[\kappa,\lambda,\mu]\in{GL}(3,\mathbb{R})$.
Let $F(g)=\left(\smallmatrix I_3&f(n)\\
0&1\endsmallmatrix\right)$, for $g\in\sqrt\pi$.
If $\nu=\sqrt\pi$ we extend $F$ to $\pi$ by setting
$F(t)=\left(\smallmatrix\Psi&0\\
0&1\endsmallmatrix\right)$. 
In this case $F$ defines a discrete cocompact embedding of 
$\pi$ in $Isom(\mathbb{S}ol^4_{m,n})$.
(See \S3 of Chapter 7.
If one of the eigenvalues is $\pm1$ then $m=n$ and
the geometry is $\mathbb{S}ol^3\times\mathbb{E}^1$.)

If the eigenvalues are not all real we may assume that 
$\lambda=\bar\kappa$ and $\mu\not=\pm1$.
Let $R_\phi\in{SO(2)}$ be rotation of $\mathbb{R}^2$ through the angle
$\phi=Arg(\kappa)$.
There is a monomorphism $f:\sqrt\pi\to\mathbb{R}^3$ such that
$f\theta=\Psi{f}$ where $\Psi=\left(\smallmatrix |\kappa|R_\phi &0\\
0&\mu\endsmallmatrix\right)$. 
Let $F(n)=\left(\smallmatrix I_3&f(n)\\
0&1\endsmallmatrix\right)$, for $n\in\sqrt\pi$,
and let $F(t)=\left(\smallmatrix\Psi&0\\
0&1\endsmallmatrix\right)$.
Then $F$ defines a discrete cocompact embedding of 
$\pi$ in $Isom(\mathbb{S}ol^4_0)$.
\end{proof}

If $M$ is a $\mathbb{S}ol^4_{m,n}$-manifold the eigenvalues of $\theta$ are 
distinct and real.
The geometry is $\mathbb{S}ol^3\times\mathbb{E}^1 (=\mathbb{S}ol^4_{m,m}$ 
for any $m\geq4$) if and only if $\theta$ has 1 as a simple eigenvalue.
If $M$ is a $\mathbb{S}ol^4_0$-manifold two of the eigenvalues are complex 
conjugates, and none are roots of unity.

The groups of $\mathbb{E}^4$-, $\mathbb{N}il^3\times\mathbb{E}^1$- and 
$\mathbb{N}il^4$-manifolds also have finite index subgroups 
$\sigma\cong\mathbb{Z}^3\rtimes_\theta\mathbb{Z}$.
We may assume that all the eigenvalues of $\theta$ are 1, 
so $N=\theta-I$ is nilpotent.
If the geometry is $\mathbb{E}^4$ then $N=0$; if it is
$\mathbb{N}il^3\times\mathbb{E}^1$ then $N\not=0$ but $N^2=0$,
while if it is $\mathbb{N}il^4$ then $N^2\not=0$ but $N^3=0$.
(Conversely, it is easy to see that such semidirect products may be realized by
lattices in the corresponding Lie groups.)

Finally, if $M$ is a $\mathbb{S}ol_1^4$-manifold then $\sqrt\pi\cong\Gamma_q$ 
for some $q\geq1$ (and so is nonabelian, of Hirsch length 3).
Every group $\pi\cong\Gamma_q\rtimes_\theta\mathbb{Z}$ 
may be realized geometrically.
(See Theorem 8.7 below.)

If $h(\sqrt\pi)=3$ then $\pi$ is an extension of $\mathbb{Z}$ 
or $D$ by a normal subgroup $\nu$ which contains $\sqrt\pi$ 
as a subgroup of finite index. 
Hence either $M$ is the mapping torus of a self homeomorphism of a flat 
3-manifold or a $\mathbb{N}il^3 $-manifold,
or it is the union of two twisted $I$-bundles over such 3-manifolds
and is doubly covered by such a mapping torus.
(Compare Theorem 8.2.)

We shall consider further the question of realizing geometrically such 
torsion-free virtually poly-$Z$ groups $\pi$ (with $h(\pi)=4$ and
$h(\sqrt\pi)=3$) in \S9.

\section{Mapping tori of self homeomorphisms of $\mathbb{E}^3 $-manifolds}

It follows from the above that a 4-dimensional infrasolvmanifold $M$ admits 
one of the product geometries of type $\mathbb{E}^4$, 
$\mathbb{N}il^3\times\mathbb{E}^1$ or $\mathbb{S}ol^3\times\mathbb{E}^1$ 
if and only if $\pi_1 (M)$ has a subgroup of finite index of the form 
$\nu\times\mathbb{Z}$, where $\nu$ is abelian, 
nilpotent of class 2 or solvable but not virtually nilpotent, respectively.      
In the next two sections we shall examine when $M$ is the mapping torus of 
a self homeomorphism of a 3-dimensional infrasolvmanifold. 
(Note that if $M$ is orientable then it must be a mapping torus, 
by Lemma 3.14 and Theorem 6.11.)

\begin{theorem} 
Let $\nu$ be the fundamental group of a flat $3$-manifold,
and let $\theta$ be an automorphism of $\nu$. Then 
\begin{enumerate}
\item $\sqrt\nu$ is the maximal abelian subgroup of $\nu$ and 
$\nu/\sqrt\nu$ embeds in $Aut(\sqrt\nu)$;

\item $Out(\nu)$ is finite if and only if $[\nu:\sqrt\nu ]>2$;

\item the restriction homomorphism from $Out(\nu)$ 
to $Aut(\sqrt\nu)$ has finite kernel;

\item if $[\nu:\sqrt\nu ]=2$ then $(\theta|_{\sqrt\nu})^2 $ 
has $1$ as an eigenvalue;

\item if $[\nu:\sqrt\nu ]=2$ and $\theta|_{\sqrt\nu} $ 
has infinite order but all of its eigenvalues 
are roots of unity then $((\theta|_{\sqrt\nu})^2 -I)^2=0$;

\item if $\theta$ is orientation-preserving and 
$((\theta|_{\sqrt\nu})^2 -I)^3=0$ but $(\theta|_{\sqrt\nu})^2\not=I$
then $(\theta|_{\sqrt\nu}-I)^3=0$, 
and so $\beta_1(\nu\rtimes_\theta\mathbb{Z})\geq2$.

\end{enumerate}
\end{theorem}

\begin{proof} It follows immediately from Theorem 1.5 
that $\sqrt\nu\cong\mathbb{Z}^3$ 
and is thus the maximal abelian subgroup of $\nu$. 
The kernel of the homomorphism from $\nu$ to $Aut(\sqrt\nu)$ 
determined by conjugation is the centralizer $C=C_\nu (\sqrt\nu)$.
As $\sqrt\nu$ is central in $C$ and $[C:\sqrt\nu]$ is finite,
$C$ has finite commutator subgroup, by Schur's Theorem 
(Proposition 10.1.4 of \cite{[Ro]}).
Since $C$ is torsion-free it must be abelian and so $C=\sqrt\nu$.                                            
Hence $H=\nu/\sqrt\nu$ embeds in $Aut(\sqrt\nu)\cong GL(3,\mathbb{Z})$. 
(This is just the holonomy representation.)

If $H$ has order 2 then $\theta$ induces the identity on $H$; 
if $H$ has order greater than 2 then some power of $\theta$ 
induces the identity on $H$,
since $\sqrt\nu$ is a characteristic subgroup of finite index. 
The matrix $\theta|_{\sqrt\nu} $ then commutes with each element 
of the image of $H$ in $GL(3,\mathbb{Z})$, 
and assertions (2)--(5) follow from simple calculations, 
on considering the possibilities for $\pi$ 
and $H$ listed in \S3 above. 
The final assertion follows on considering the Jordan normal form of 
$\theta|_{\sqrt\nu}$.
\end{proof}

\begin{cor}
The mapping torus $M(\phi)=N\times_\phi S^1$ of a self homeomorphism 
$\phi$ of a flat $3$-manifold $N$ is flat if and only if the outer automorphism $[\phi_* ]$ 
induced by $\phi$ has finite order. 
\qed
\end{cor}

If $N$ is flat and $[\phi_* ]$ has infinite order then 
$M(\phi)$ may admit one of 
the other product geometries $\mathbb{S}ol^3 \times\mathbb{E}^1 $ or $\mathbb{N}il^3 \times\mathbb{E}^1 $;
otherwise it must be a $\mathbb{S}ol^4_{m,n}$-, $\mathbb{S}ol^4_0$- or $\mathbb{N}il^4 $-manifold. 
(The latter can only happen if $N=\mathbb{R}^3/\mathbb{Z}^3 $, 
by part ({\sl v}) of the theorem.)

\begin{theorem} 
Let $M$ be a closed $4$-manifold with a geometry of solvable Lie type and
fundamental group $\pi$.
If $\sqrt\pi\cong\mathbb{Z}^3$ and $\pi/\sqrt\pi$ 
is an extension of $D$ by a finite normal subgroup 
then $M$ is a $\mathbb{S}ol^3\times\mathbb{E}^1 $-manifold.
\end{theorem}

\begin{proof} Let $p:\pi\to D$ be an epimorphism with kernel $K$ 
containing $\sqrt\pi$ as a subgroup of finite index,
and let $t$ and $u$ be elements of $\pi$ whose images 
under $p$ generate $D$ and such that 
$p(t)$ generates an infinite cyclic subgroup of index 2 in $D$.
Then there is an $N>0$ such that the image of $s=t^N $ in $\pi/\sqrt\pi$
generates a normal subgroup. In particular, the subgroup generated by $s$ and $\sqrt\pi$ 
is normal in $\pi$ and $usu^{-1}$ and $s^{-1}$ have the same image in $\pi/\sqrt\pi$.
Let $\theta$ be the matrix of the action of $s$ on $\sqrt\pi$, with respect to some basis $\sqrt\pi\cong\mathbb{Z}^3 $.
Then $\theta$ is conjugate to its inverse, since $usu^{-1}$ and $s^{-1}$ agree modulo $\sqrt\pi$.
Hence one of the eigenvalues of $\theta$ is $\pm 1$.
Since $\pi$ is not virtually nilpotent the eigenvalues of $\theta$ must be distinct,
and so the geometry must be of type $\mathbb{S}ol^3 \times\mathbb{E}^1 $. 
\end{proof}

\begin{cor}
If $M$ admits one of the geometries 
$\mathbb{S}ol^4_0$ or $\mathbb{S}ol^4_{m,n}$ 
with $m\not= n$ then it is the mapping torus 
of a self homeomorphism of $\mathbb{R}^3/\mathbb{Z}^3$,
and so $\pi\cong\mathbb{Z}^3\rtimes_\theta\mathbb{Z}$ 
for some $\theta$ in $GL(3,\mathbb{Z})$
and is a metabelian poly-$Z$ group.
\end{cor} 

\begin{proof} This follows immediately from Theorems 8.3 and 8.4. 
\end{proof}

We may use the idea of Theorem 8.2 to give examples of $\mathbb{E}^4$-, 
$\mathbb{N}il^4$-, $\mathbb{N}il^3\times\mathbb{E}^1$- and 
$\mathbb{S}ol^3\times\mathbb{E}^1$-manifolds which are not mapping tori.
For instance, the groups with presentations
\begin{gather*} 
\langle u,v,x,y,z\mid xy=yx,\medspace xz=zx,\medspace yz=zy,\medspace
uxu^{-1} =x^{-1},\medspace u^2=y,\medspace uzu^{-1} =z^{-1},\\
v^2=z,\medspace vxv^{-1} =x^{-1},\medspace vyv^{-1} =y^{-1} \rangle,\\
\langle u,v,x,y,z\mid xy=yx,\medspace xz=zx,\medspace yz=zy,\medspace
u^2= x,\medspace uyu^{-1}= y^{-1},\medspace uzu^{-1}=z^{-1},\\
v^2=x,\medspace vyv^{-1}= v^{-4} y^{-1},\medspace 
vzv^{-1}= z^{-1}\rangle\\
\text{and}\qquad\langle u,v,x,y,z\mid xy=yx,\medspace xz=zx,\medspace 
yz=zy,\medspace u^2=x,\medspace v^2=y,\qquad\phantom{\text{and}}\\
uyu^{-1}= x^4 y^{-1},\medspace vxv^{-1}= x^{-1} y^2,\medspace uzu^{-1}= vzv^{-1}= z^{-1}\rangle
\end{gather*} 
are each generalised free products of two copies of 
$\mathbb{Z}^2\rtimes_{-I}\mathbb{Z}$ amalgamated 
over their maximal abelian subgroups.
The Hirsch-Plotkin radicals of these groups are isomorphic 
to $\mathbb{Z}^4 $ (generated by $\{ (uv)^2 ,x,y,z\}$), 
$\Gamma_2 \times\mathbb{Z}$ (generated by $\{ uv,x,y,z\}$) 
and $\mathbb{Z}^3$ (generated by $\{x,y,z\}$), respectively. 
The group with presentation 
\begin{gather*} 
\langle u,v,x,y,z\mid xy=yx,\medspace xz=zx,\medspace yz=zy,
\medspace u^2=x,\medspace uz=zu,\medspace uyu^{-1} =xy^{-1},\\
vxv^{-1}=x^{-1},\medspace v^2=y,\medspace vzv^{-1}=v^2 z^{-1}\rangle
\end{gather*} 
is a generalised free product of copies of 
$(\mathbb{Z}\rtimes_{-1}\! \mathbb{Z})\!\times\!\mathbb{Z}$ 
(generated by $\{ u,y,z\} $) and 
$\mathbb{Z}^2\rtimes_{-I}\!\mathbb{Z}$ (generated by $\{ v,x,z\}$) 
amalgamated over their maximal abelian subgroups.
Its Hirsch-Plotkin radical is the subgroup of index 4 generated by 
$\{(uv)^2,x,y,z\}$,
and is nilpotent of class 3. 
The manifolds corresponding to these groups admit 
the geometries $\mathbb{E}^4 $, 
$\mathbb{N}il^3\times\mathbb{E}^1 $, 
$\mathbb{S}ol^3 \times\mathbb{E}^1 $ and $\mathbb{N}il^4 $, respectively. 
However they cannot be mapping tori, as these groups each have finite abelianization.

\section{Mapping tori of self homeomorphisms of $\mathbb{N}il^3$-manifolds}

Let $\phi$ be an automorphism of $\Gamma_q$, 
sending $x$ to $x^ay^bz^m$ and $y$ to $x^cy^dz^n$ for some $a,\dots,n$ in
$\mathbb{Z}$.
The induced automorphism of $\Gamma_q/I(\Gamma_q)\cong\mathbb{Z}^2$ 
has matrix
$A=\left(\smallmatrix a&c\\
b&d\endsmallmatrix\right)\in{GL}(2,\mathbb{Z})$ and
$\phi(z)=z^{det(A)}$.
(In particular, the $PD_3$-group $\Gamma_q$ is orientable,
and $\phi$ is orientation preserving,
as observed in \S2 of Chapter 7.
See also \S3 of Chapter 18 below.)
Every pair $(A,\mu)$ in the set $GL(2,\mathbb{Z})\times\mathbb{Z}^2$ 
determines an automorphism (with $\mu=(m,n)$).
However $Aut(\Gamma_q)$ is not a semidirect product,
as 
\begin{equation*}
(A,\mu)(B,\nu)=(AB,\mu B+det(A)\nu+q\omega(A,B)),
\end{equation*}                                                           
where $\omega(A,B)$ is biquadratic in the entries of $A$ and $B$.
The natural map $p:Aut(\Gamma_q)\to{Aut(\Gamma_q/\zeta\Gamma_q)}=
GL(2,\mathbb{Z})$ sends $(A,\mu)$ to $A$ and is an epimorphism, with 
$\mathrm{Ker}(p)\cong\mathbb{Z}^2$.
The inner automorphisms are represented by $q\mathrm{Ker}(p)$,
and $Out(\Gamma_q)\cong(Z/qZ)^2\rtimes{GL(2,\mathbb{Z})}$. 
(Let $[A,\mu]$ be the image of $(A,\mu)$ in $Out(\Gamma_q)$. 
Then $[A,\mu][B,\nu]=[AB,\mu B+det(A)\nu]$.)
In particular, $Out(\Gamma_1)=GL(2,\mathbb{Z})$.

\begin{theorem} 
Let $\nu$ be the fundamental group of a $\mathbb{N}il^3$-manifold $N$.
Then
\begin{enumerate}
\item $\nu/\sqrt\nu$ embeds in $Aut(\sqrt\nu/\zeta\sqrt\nu)\cong GL(2,\mathbb{Z})$;

\item $\bar\nu=\nu/\zeta\sqrt\nu$ is a 2-dimensional crystallographic group;

\item the images of elements of $\bar\nu$ of finite order under the holonomy 

\noindent representation in $Aut(\sqrt{\bar\nu})\cong GL(2,\mathbb{Z})$ 
have determinant $1$;
 
\item $Out(\bar\nu)$ is infinite if and only if $\bar\nu\cong\mathbb{Z}^2$ 
or $\mathbb{Z}^2\rtimes _{-I} (Z/2Z)$;

\item the kernel of the natural homomorphism from $Out(\nu)$ to 
$Out(\bar\nu)$ is finite.

\item $\nu$ is orientable and every automorphism of $\nu$ is
orientation preserving.
\end{enumerate}
\end{theorem}

\begin{proof} 
Let $h:\nu\to Aut(\sqrt\nu/\zeta\sqrt\nu)$ be the homomorphism determined by 
conjugation, and let $C=\mathrm{Ker}(h)$.
Then $\sqrt\nu/\zeta\sqrt\nu$ is central in $C/\zeta\sqrt\nu$
and $[C/\zeta\sqrt\nu:\sqrt\nu/\zeta\sqrt\nu]$ is finite,
so $C/\zeta\sqrt\nu$ has finite commutator subgroup, by Schur's Theorem 
(Proposition 10.1.4 of \cite{[Ro]}.)
Since $C$ is torsion-free it follows easily that $C$ is nilpotent and 
hence that $C=\sqrt\nu$. This proves (1) and (2). 
In particular, $h$ factors through the holonomy representation for $\bar\nu$,
and $gzg^{-1}=z^{d(g)}$ for all $g\in\nu$ and $z\in\zeta\sqrt\nu$,
where $d(g)=det(h(g))$.
If $g\in\nu$ is such that $g\not=1$ and $g^k\in\zeta\sqrt\nu$ for some $k>0$
then $g^k\not=1$ and so $g$ must commute with elements of $\zeta\sqrt\nu$,
i.e., the determinant of the image of $g$ is 1.
Condition (4) follows as in Theorem 8.4, 
on considering the possible finite subgroups of $GL(2,\mathbb{Z})$. 
(See Theorem 1.3.)

If $\zeta\nu\not=1$ then $\zeta\nu=\zeta\sqrt\nu\cong\mathbb{Z}$ 
and so the kernel of the natural homomorphism from $Aut(\nu)$ 
to $Aut(\bar\nu)$ is isomorphic to $Hom(\nu/\nu',\mathbb{Z})$.
If $\nu/\nu'$ is finite this kernel is trivial.
If $\bar\nu\cong\mathbb{Z}^2 $ then $\nu=\sqrt\nu\cong \Gamma_q $, 
for some $q\geq 1$,
and the kernel is isomorphic to $(Z/qZ)^2 $. 
Otherwise $\bar\nu\cong\mathbb{Z}\rtimes_{-1}\!\mathbb{Z}$, 
$\mathbb{Z}\times\! D$ or $D\rtimes_\tau\mathbb{Z}$ 
(where $\tau$ is the automorphism of $D=(Z/2Z)*(Z/2Z)$ 
which interchanges the factors).
But then $H^2 (\bar\nu;\mathbb{Z})$ is finite, 
and so any central extension of
such a group by $\mathbb{Z}$ is virtually abelian, 
and thus not a $\mathbb{N}il^3$-manifold group. 

If $\zeta\nu=1$ then $\nu/\sqrt\nu<GL(2,\mathbb{Z})$ 
has an element of order 2 with determinant $-1$.
No such element can be conjugate to 
$\left(\begin{smallmatrix} 
0&1\\
1&0
\end{smallmatrix}\right),$
for otherwise $\nu$ would not be torsion-free.
Hence the image of $\nu/\sqrt\nu$ in $GL(2,\mathbb{Z})$ is conjugate to a 
subgroup of the group of diagonal matrices
$\left(\begin{smallmatrix} 
\epsilon&0\\
0&\epsilon'
\end{smallmatrix}\right),$
with $|\epsilon|=|\epsilon'|=1$.
If $\nu/\sqrt\nu$ is generated by 
$\left(\begin{smallmatrix} 
1&0\\
0&-1
\end{smallmatrix}\right)$
then $\nu/\zeta\sqrt\nu\cong\mathbb{Z}\rtimes_{-1}\!\mathbb{Z}$ 
and $\nu\cong\mathbb{Z}^2\rtimes_\theta\mathbb{Z}$,
where $\theta=
\left(\begin{smallmatrix}
-1& r\\0&-1
\end{smallmatrix}\right)$
for some nonzero integer $r$,
and $N$ is a circle bundle over the Klein bottle. 
If $\nu/\sqrt\nu\cong(Z/2Z)^2$ then $\nu$ has a presentation
\begin{equation*} 
\langle t,u,z\mid u^2=z,tzt^{-1}=z^{-1},\medspace ut^2u^{-1}=t^{-2}z^s\rangle,
\end{equation*} 
and $N$ is a Seifert bundle over the orbifold $P(22)$.
It may be verified in each case that the kernel of the natural homomorphism 
from $Out(\nu)$ to $Out(\bar\nu)$ is finite.
Therefore (5) holds.

Since $\sqrt\nu\cong\Gamma_q$ is a $PD_3^+$-group, $[\nu:\sqrt\nu]<\infty$
and every automorphism of $\Gamma_q$ is orientation preserving
$\nu$ must also be orientable.
Since $\sqrt\nu$ is characteristic in $\nu$ and the image of 
$H_3(\sqrt\nu;\mathbb{Z})$ in $H_3(\nu;\mathbb{Z})$ has index
$[\nu:\sqrt\nu]$ it follows easily that any automorphism of $\nu$ must be
orientation preserving.
\end{proof}

In fact every $\mathbb{N}il^3 $-manifold is a Seifert bundle 
over a 2-dimensional euclidean orbifold \cite{[Sc83']}. 
The base orbifold must be one of the seven such with no reflector curves,
by (3).

\begin{theorem} 
Let $\theta=(A,\mu)$ be an automorphism of $\Gamma_q$ and
$\pi=\Gamma_q\rtimes_\theta\mathbb{Z}$. Then
\begin{enumerate}
\item If $A$ has finite order $h(\sqrt\pi)=4$ and $\pi$ is a lattice in
$Isom(\mathbb{N}il^3\times\mathbb{E}^1)$;

\item if $A$ has infinite order and equal eigenvalues 
$h(\sqrt\pi)=4$ and $\pi$ is a lattice in $Isom(\mathbb{N}il^4)$;

\item Otherwise $\sqrt\pi=\Gamma_q$ and $\pi$ is a 
lattice in $Isom(\mathbb{S}ol^4_1)$.
\end{enumerate}
\end{theorem}

\begin{proof}
Let $t\in\pi$ represent a generator for $\pi/\Gamma_q\cong\mathbb{Z}$.
The image of $\theta$ in $Out(\Gamma_q)$ has finite order 
if and only if $A$ has finite order.
If $A^k=1$ for some $k\geq1$ the subgroup generated by $\Gamma_q$ 
and $t^{kq}$ is isomorphic to $\Gamma_q\times\mathbb{Z}$.
If $A$ has infinite order and equal eigenvalues then $A^2$ is conjugate to
$\left(\smallmatrix 1& n\\
0&1\endsmallmatrix\right)$, for some $n\not=0$, and the subgroup generated by
$\Gamma_q$ and $t^2$ is nilpotent of class 3.
In each of these cases $\pi$ is virtually nilpotent, 
and may be embedded as a lattice in $Isom(\mathbb{N}il^3\times\mathbb{E}^1)$ 
or $Isom(\mathbb{N}il^4)$ \cite{[De]}.

Otherwise the eigenvalues $\alpha,\beta$ of $A$ are distinct and not $\pm1$.
Let $e,f\in\mathbb{R}^2$ be the corresponding eigenvectors.
Let $(1,0)={x_1e+x_2f}$, $(0,1)={y_1e+y_2f}$, $\mu={z_1e+z_2f}$
and $h=\frac1q(x_2y_1-x_1y_2)$.
\quad
Let $F(x)=
\left(\begin{matrix}
1& x_2&0\\
0&1& x_1\\
0&0&1
\end{matrix}\right)$,\\
$F(y)=\left(\begin{matrix}
1& y_2&0\\
0&1& y_1\\
0&0&1
\end{matrix}\right)$,
$F(z)= 
\left(\begin{matrix}
 1& 0 & h\\
0&1& 0\\
0&0&1
\end{matrix}\right)$ and
$F(t)=\left(\begin{matrix}
\alpha\beta& z_2&0\\
0&\alpha& z_1\\
0&0&1
\end{matrix}\right)$.

Then $F$ defines an embedding of $\pi$ as a lattice in 
$Isom(\mathbb{S}ol^4_1)$.
\end{proof}

\begin{theorem} 
The mapping torus $M(\phi)=N\times_\phi S^1 $ of a self homeomorphism 
$\phi$ of a $\mathbb{N}il^3$-manifold $N$ is orientable, and is a 
$\mathbb{N}il^3\times\mathbb{E}^1 $-manifold if and only if the outer automorphism 
$[\phi_* ]$ induced by $\phi$ has finite order.
\end{theorem}

\begin{proof} 
Since $N$ is orientable and $\phi$ is orientation preserving (by part (6) of
Theorem 8.6) $M(\phi)$ must be orientable.

The subgroup $\zeta\sqrt\nu$ is characteristic in $\nu$ 
and hence normal in $\pi$, 
and $\nu/\zeta\sqrt\nu$ is virtually $\mathbb{Z}^2 $.
If $M(\phi)$ is a $\mathbb{N}il^3 \times\mathbb{E}^1 $-manifold then 
$\pi/\zeta\sqrt\nu$ is also virtually abelian.
It follows easily that that the image of $\phi_* $ in 
$Aut(\nu/\zeta\sqrt\nu)$ has finite order. 
Hence $[\phi_* ]$ has finite order also, by Theorem 8.6.
Conversely, if $[\phi_* ]$ has finite order in $Out(\nu)$ 
then $\pi$ has a subgroup of finite index 
which is isomorphic to $\nu\times\mathbb{Z}$, 
and so $M(\phi)$ has the product geometry, by the discussion above. 
\end{proof}

Theorem 4.2 of \cite{[KLR83]} (which extends Bieberbach's theorem 
to the virtually nilpotent case) may be used to show directly 
that every outer automorphism class of finite order of the
fundamental group of an $\mathbb{E}^3 $- or 
$\mathbb{N}il^3 $-manifold is realizable by an isometry 
of an affinely equivalent manifold. 

\begin{theorem} 
Let $M$ be a closed $\mathbb{N}il^3\times\mathbb{E}^1$-, 
$\mathbb{N}il^4$- or $\mathbb{S}ol_1^4 $-manifold.
Then $M$ is the mapping torus of a self homeomorphism of a 
$\mathbb{N}il^3 $-manifold if and only if it is orientable.
\end{theorem}
                       
\begin{proof} If $M$ is such a mapping torus then it is orientable, 
by Theorem 8.8.
Conversely, if $M$ is orientable then $\pi=\pi_1 (M)$ has infinite abelianization, by Lemma 3.14. 
Let $p:\pi\to\mathbb{Z}$ be an epimorphism with kernel $K$, and let 
$t$ be an element of $\pi$ such that $p(t)=1$.
If $K$ is virtually nilpotent of class 2 we are done, by Theorem 6.12.
(Note that this must be the case if $M$ is a $\mathbb{S}ol_1^4 $-manifold.)
If $K$ is virtually abelian then $K\cong\mathbb{Z}^3$ or $G_2$, 
by part (5) of Theorem 8.4.
The action of $t$ on $\sqrt{K}$ by conjugation must be orientation
preserving, since $M$ is orientable. 
Since $\det(t)=1$ and $\pi$ is virtually nilpotent but not virtually abelian,
at least one eigenvalue must be $+1$.
It follows easily that $\beta_1(\pi)\geq2$.
Hence there is another epimorphism with kernel nilpotent of class 2,
and so the theorem is proven. 
\end{proof}

\begin{cor}
Let $M$ be a closed $\mathbb{S}ol^4_1$-manifold with fundamental group $\pi$.
Then $\beta_1(M)\leq1$ and $M$ is orientable if and only if $\beta_1(M)=1$.
\end{cor}                                    

\begin{proof} 
The first assertion is clear if $\pi$ is a semidirect product 
$\Gamma_q\rtimes_\theta\mathbb{Z}$,
and then follows in general. 
Hence if  $p\colon\pi\to\mathbb{Z}$ is an epimorphism $\mathrm{Ker}(p)$ 
must be virtually nilpotent of class 2 and the result follows from the theorem.
\end{proof}

If $M$ is a $\mathbb{N}il^3\times\mathbb{E}^1$- or 
$\mathbb{N}il^4$-manifold then 
$\beta_1(\pi)\leq3$ or 2, respectively,
with equality if and only if $\pi$ is nilpotent.
In the latter case $M$ is orientable, and is a mapping torus,
both of a self homeomorphism of $\mathbb{R}^3/\mathbb{Z}^3$ 
and also of a self 
homeomorphism of a $\mathbb{N}il^3 $-manifold.
We have already seen that $\mathbb{N}il^3\times\mathbb{E}^1 $- and 
$\mathbb{N}il^4 $-manifolds need not be mapping tori at all. 
We shall round out this discussion with examples illustrating the remaining
combinations of mapping torus structure and orientation compatible with 
Lemma 3.14 and Theorems 8.8 and 8.9.
As the groups have abelianization of rank 1 the corresponding manifolds 
are mapping tori in an essentially unique way.
The groups with presentations                                                  
\begin{gather*} 
\langle t,x,y,z\mid xz=zx,\medspace yz=zy,\medspace txt^{-1}=x^{-1}, 
\medspace tyt^{-1}=y^{-1},\medspace tzt^{-1}= yz^{-1} \rangle\\
\text{and}\quad 
\langle t,x,y,z\mid xyx^{-1}y^{-1}=z,\medspace xz=zx,\medspace yz=zy,
\medspace txt^{-1}=x^{-1},\medspace tyt^{-1}=y^{-1}\rangle
\end{gather*} 
are each virtually nilpotent of class 2.
The corresponding $\mathbb{N}il^3 \times\mathbb{E}^1 $-manifolds are 
mapping tori of self homeomorphisms of $\mathbb{R}^3/\mathbb{Z}^3 $ 
and a $\mathbb{N}il^3 $-manifold, 
respectively.
The groups with presentations 
\begin{gather*} 
\langle t,x,y,z\mid xz=zx,\medspace yz=zy,\medspace txt^{-1}=x^{-1},
\medspace tyt^{-1}=xy^{-1},\medspace tzt^{-1}= yz^{-1} \rangle\\
\text{and}\quad\langle t,x,y,z\mid xyx^{-1}y^{-1}=z,
\medspace xz=zx,\medspace yz=zy,\medspace txt^{-1} =x^{-1},
\medspace tyt^{-1} =xy^{-1}\rangle
\end{gather*} 
are each virtually nilpotent of class 3.
The corresponding $\mathbb{N}il^4 $-manifolds are mapping tori 
of self homeomorphisms of $\mathbb{R}^3/\mathbb{Z}^3 $ 
and of a $\mathbb{N}il^3 $-manifold, 
respectively.
The group with presentation
\begin{gather*} 
\langle t,u,x,y,z\mid xyx^{-1}y^{-1}=z^2,\medspace xz=zx,\medspace 
yz=zy,\medspace txt^{-1} =x^2y,\medspace tyt^{-1}=xy,\\
tz=zt,\medspace u^4=z,\medspace uxu^{-1}=y^{-1},\medspace uyu^{-1} =x,
\medspace utu^{-1}=t^{-1}\rangle
\end{gather*} 
has Hirsch-Plotkin radical isomorphic to $\Gamma_2 $ (generated by 
$\{ x,y,z\}$), and has finite abelianization.
The corresponding $\mathbb{S}ol_1^4 $-manifold is nonorientable 
and is not a mapping torus.

\section{Mapping tori of self homeomorphisms of $\mathbb{S}ol^3 $-manifolds}

The arguments in this section are again analogous to those of \S6. 

\begin{theorem} 
Let $\sigma$ be the fundamental 
group of a $\mathbb{S}ol^3$-manifold.
Then
\begin{enumerate}
\item $\sqrt\sigma\cong\mathbb{Z}^2 $ and 
$\sigma/\sqrt\sigma\cong\mathbb{Z}$ or $D$;

\item $Out(\sigma)$ is finite.
\end{enumerate} 
\end{theorem}

\begin{proof} The argument of Theorem 1.6 implies that $h(\sqrt\sigma)>1$.
Since $\sigma$ is not virtually nilpotent $h(\sqrt\sigma)<3$.
Hence $\sqrt\sigma\cong\mathbb{Z}^2 $, by Theorem 1.5. 
Let $\tilde F$ be the preimage in $\sigma$ 
of the maximal finite normal subgroup of $\sigma/\sqrt\nu$, 
let $t$ be an element of $\sigma$ whose image generates the maximal
abelian subgroup of $\sigma/\tilde F$ and let 
$\tau$ be the automorphism of $\tilde F$ determined by conjugation by $t$. 
Let $\sigma_1 $ be the subgroup of $\sigma$ 
generated by $\tilde F$ and $t$.                         
Then $\sigma_1 \cong\tilde F\rtimes_\tau\mathbb{Z}$,
$[\sigma:\sigma_1 ]\leq 2$, $\tilde F$ is torsion-free and $h(\tilde F)=2$. 
If $\tilde F\not=\sqrt\sigma$ then 
$\tilde F\cong\mathbb{Z}\rtimes_{-1}\!\mathbb{Z}$.
But extensions of $\mathbb{Z}$ by $\mathbb{Z}\rtimes_{-1}\!\mathbb{Z}$ 
are virtually abelian, 
since $Out(\mathbb{Z}\rtimes_{-1}\!\mathbb{Z})$ is finite.
Hence $\tilde F=\sqrt\sigma$ and so 
$\sigma/\sqrt\sigma\cong\mathbb{Z}$ or $D$.

Every automorphism of $\sigma$ induces automorphisms of 
$\sqrt\sigma$ and of $\sigma/\sqrt\sigma$.
Let $Out^+ (\sigma)$ be the subgroup of $Out(\sigma)$ 
represented by automorphisms
which induce the identity on $\sigma/\sqrt\sigma$.
The restriction of any such automorphism to $\sqrt\sigma$ commutes with $\tau$.
We may view $\sqrt\sigma$ as a module over the ring $R=\mathbb{Z}[X]/(\lambda(X))$, where
$\lambda(X)=X^2-tr(\tau)X+det(\tau)$ is the characteristic polynomial of $\tau$.
The polynomial $\lambda$ is irreducible and has real roots 
which are not roots of unity,
for otherwise $\sqrt\sigma\rtimes_\tau\mathbb{Z}$ would be virtually nilpotent.
Therefore $R$ is a domain and its field of fractions $\mathbb{Q}[X]/(\lambda(X))$ 
is a real quadratic number field.
The $R$-module $\sqrt\sigma$ is clearly finitely generated, 
$R$-torsion-free and of rank 1.
Hence the endomorphism ring $End_R (\sqrt\sigma)$ is a subring of $\tilde R$,
the integral closure of $R$. 
Since $\tilde R$ is the ring of integers in $\mathbb{Q}[X]/(\lambda(X))$
the group of units $\tilde R^\times $ is isomorphic to 
$\{ \pm 1\}\times\mathbb{Z}$.
Since $\tau$ determines a unit of infinite order in $R^\times $ 
the index $[\tilde R^\times :\langle\tau\rangle ]$ is finite.

Suppose now that $\sigma/\sqrt\sigma\cong\mathbb{Z}$.
If $f$ is an automorphism which induces the identity on $\sqrt\sigma$ and on $\sigma/\sqrt\sigma$
then $f(t)=tw$ for some $w$ in $\sqrt\sigma$. 
If $w$ is in the image of $\tau-1$ then $f$ is an inner automorphism.
Now $\sqrt\sigma/(\tau-1)\sqrt\sigma$ is finite, of order $det(\tau-1)$. 
Since $\tau$ is the image of an inner automorphism of $\sigma$ 
it follows that $Out^+(\sigma)$ is an extension of a subgroup of 
$\tilde R^\times /\langle\tau\rangle$ 
by $\sqrt\sigma/(\tau-1)\sqrt\sigma$.
Hence $Out(\sigma)$ has order dividing 
$2[\tilde R^\times :\langle\tau\rangle]det(\tau-1)$. 

If $\sigma/\sqrt\sigma\cong D$ then $\sigma$ has 
a characteristic subgroup $\sigma_1 $
such that $[\sigma:\sigma_1 ]=2$, 
$\sqrt\sigma<\sigma_1 $ and $\sigma_1 /\sqrt\sigma\cong\mathbb{Z}=\sqrt D$.
Every automorphism of $\sigma$ restricts to an automorphism of $\sigma_1 $.
It is easily verified that the restriction from $Aut(\sigma)$ 
to $Aut(\sigma_1 )$ is a monomorphism. 
Since $Out(\sigma_1 )$ is finite it follows that $Out(\sigma)$ is also finite. 
\end{proof}

\begin{cor}
The mapping torus of a self homeomorphism of a 
$\mathbb{S}ol^3 $-manifold is a $\mathbb{S}ol^3 \times\mathbb{E}^1 $-manifold. 
\qed
\end{cor}

The group with presentation 
\begin{equation*} 
\langle{x,y,t}\mid{xy=yx,}\medspace{txt^{-1}=xy,}\medspace{tyt^{-1}=x}\rangle
\end{equation*} 
is the fundamental group of a nonorientable $\mathbb{S}ol^3 $-manifold $\Sigma$.
The nonorientable $\mathbb{S}ol^3 \times\mathbb{E}^1 $-manifold 
$\Sigma\times S^1$ is the mapping torus of $id_\Sigma$
and is also the mapping torus of a self homeomorphism 
of $\mathbb{R}^3 /\mathbb{Z}^3 $.

The groups with presentations                                                 
 \begin{gather*} 
\langle t,x,y,z\mid xy=yx,\medspace zxz^{-1} =x^{-1},\medspace zyz^{-1} =y^{-1},
\medspace txt^{-1}=xy,\medspace tyt^{-1}=x,\\ 
tzt^{-1}= z^{-1}\rangle,\\
\langle t,x,y,z\mid xy=yx,\medspace zxz^{-1} =x^2 y,\medspace
zyz^{-1} =xy,\medspace tx=xt,\medspace tyt^{-1}=x^{-1} y^{-1},\\
tzt^{-1}= z^{-1}\rangle,
\\
\langle t,x,y,z\mid xy=yx,\medspace xz=zx,\medspace yz=zy,\medspace
txt^{-1} =x^2 y,\medspace tyt^{-1} =xy,\medspace tzt^{-1}=z^{-1}\rangle
\end{gather*}
\begin{gather*}
\text{and}\qquad\langle t,u,x,y\mid xy=yx,\medspace txt^{-1} =x^2 y,
\medspace tyt^{-1} =xy,\medspace uxu^{-1}=y^{-1},\qquad\phantom{\text{and}}\\  
uyu^{-1} =x,\medspace utu^{-1} =t^{-1}\rangle
\end{gather*} 
have Hirsch-Plotkin radical $\mathbb{Z}^3 $ and abelianization of rank 1.
The corresponding $\mathbb{S}ol^3\times\mathbb{E}^1 $-manifolds are mapping 
tori in an essentially unique way.
The first two are orientable, and are mapping tori of self homeomorphisms of 
the orientable flat 3-manifold with holonomy of order 2 and of an 
orientable $\mathbb{S}ol^3 $-manifold, respectively.
The latter two are nonorientable, and are mapping tori of orientation reversing
self homeomorphisms of $\mathbb{R}^3 /\mathbb{Z}^3 $ 
and of the same orientable $\mathbb{S}ol^3 $-manifold, 
respectively. 

\section{Realization and classification}

Let $\pi$ be a torsion-free virtually poly-$Z$ group of Hirsch length 4.
If $\pi$ is virtually abelian then it is the fundamental group of a flat
4-manifold, by the work of Bieberbach, and such groups are listed in
\S2-\S4 above.

If $\pi$ is virtually nilpotent but not virtually abelian then $\sqrt\pi$ 
is nilpotent of class 2 or 3.
In the first case it has a characteristic chain
$\sqrt\pi'\cong\mathbb{Z}<C=\zeta\sqrt\pi\cong\mathbb{Z}^2$.
Let $\theta:\pi\to Aut(C)\cong GL(2,\mathbb{Z})$ be the homomorphism induced 
by conjugation in $\pi$.
Then $\mathrm{Im}(\theta)$ is finite and triangular, 
and so is 1, $Z/2Z$ or $(Z/2Z)^2$.
Let $K=C_\pi(C)=\mathrm{Ker}(\theta)$.
Then $K$ is torsion-free and $\zeta K=C$, 
so $K/C$ is a flat 2-orbifold group. 
Moreover as $K/\sqrt K$ acts trivially on $\sqrt\pi'$ 
it must act orientably on $\sqrt K/C$,
and so $K/\sqrt K$ is cyclic of order 1, 2, 3, 4 or 6.
As $\sqrt\pi$ is the preimage of $\sqrt K$ in $\pi$ we see that 
$[\pi:\sqrt\pi]\leq 24$.
(In fact $\pi/\sqrt\pi\cong F$ or $F\oplus(Z/2Z)$, 
where $F$ is a finite subgroup of $GL(2,\mathbb{Z})$,
excepting only direct sums of the dihedral groups of order 6, 8 or 12 with
$(Z/2Z)$ \cite{[De]}.)
Otherwise (if $\sqrt\pi'\nleq\zeta\sqrt\pi$) it has a subgroup 
of index $\leq2$ which is a semidirect product 
$\mathbb{Z}^3\rtimes_\theta\mathbb{Z}$, 
by part (5) of Theorem 8.4.
Since $(\theta^2-I)$ is nilpotent it follows that 
$\pi/\sqrt\pi=1$, $Z/2Z$ or $(Z/2Z)^2$.
All these possibilities occur.

Such virtually nilpotent groups are fundamental groups of 
$\mathbb{N}il^3\times\mathbb{E}^1$- and $\mathbb{N}il^4$-manifolds
(respectively), and are classified in \cite{[De]}.
Dekimpe observes that $\pi$ has a characteristic subgroup 
$A\cong\mathbb{Z}$ such that $Q=\pi/A$ is a 
$\mathbb{N}il^3$- or $\mathbb{E}^3$-orbifold group and 
classifies the torsion-free extensions of such $Q$ by $\mathbb{Z}$.
There are 61 families of $\mathbb{N}il^3\times\mathbb{E}^1$-groups and 
7 families of $\mathbb{N}il^4$-groups.
He also gives a faithful affine representation for each such group.

We shall sketch an alternative approach for the geometry $\mathbb{N}il^4$, 
which applies also to $\mathbb{S}ol^4_{m,n}$, $\mathbb{S}ol^4_0$ 
and $\mathbb{S}ol^4_1$. 
Each such group $\pi$ has a characteristic subgroup $\nu$ 
of Hirsch length 3, and such that $\pi/\nu\cong\mathbb{Z}$ or $D$.
The preimage in $\pi$ of $\sqrt{\pi/\nu}$ is characteristic, 
and is a semidirect product $\nu\rtimes_\theta\mathbb{Z}$.
Hence it is determined up to isomorphism by the union of the conjugacy 
classes of $\theta$ and $\theta^{-1}$ in $Out(\nu)$, by Lemma 1.1.
All such semidirect products may be realized as lattices 
and have faithful affine representations.

If the geometry is $\mathbb{N}il^4$ then 
$\nu=C_{\sqrt\pi}(\zeta_2\sqrt\pi)\cong\mathbb{Z}^3$, 
by Theorem 1.5 and part (5) of Theorem 8.4.
Moreover $\nu$ has a basis $x,y,z$ such that
$\langle z\rangle=\zeta\sqrt\pi$ and $\langle y,z\rangle=\zeta_2\sqrt\pi$.
As these subgroups are characteristic the matrix of $\theta$
with respect to such a basis is $\pm (I+N)$, 
where $N$ is strictly lower triangular and $n_{21}n_{32}\not=0$.
(See \S5 above.)
The conjugacy class of $\theta$ is determined by 
$(det(\theta),|n_{21}|,|n_{32}|,[n_{31}~mod~(n_{32})])$.
(Thus $\theta$ is conjugate to $\theta^{-1}$ if and only if
$n_{32}$ divides $2n_{31}$.)
The classification is more complicated if $\pi/\nu\cong D$.

If the geometry is $\mathbb{S}ol^4_{m,n}$ for some $m\not=n$ then
$\pi\cong\mathbb{Z}^3\rtimes_\theta\mathbb{Z}$, 
where the eigenvalues of $\theta$ are distinct and real, 
and not $\pm1$, by Corollary 8.5.1.
The translation subgroup $\pi\cap Sol^4_{m,n}$ is 
$\mathbb{Z}^3\rtimes_A\mathbb{Z}$, 
where $A=\theta$ or $\theta^2$ is the least nontrivial power of $\theta$ 
with all eigenvalues positive, and has index $\leq2$ in $\pi$.
Conversely, every such group is a lattice in $Isom(\mathbb{S}ol^4_{m,n})$,
by Theorem 8.3.
The conjugacy class of $\theta$ is determined by its characteristic polynomial
$\Delta_\theta(t)$ and the ideal class of $\nu\cong\mathbb{Z}^3$, 
considered as a rank 1 module over the order 
$\Lambda/(\Delta_\theta(t))$, by Theorem 1.4.
(No such $\theta$ is conjugate to its inverse, as neither 1 nor -1 is an
eigenvalue.)

A similar argument applies for $\mathbb{S}ol^4_0$,
where we again have $\pi\cong\mathbb{Z}^3\rtimes_\theta\mathbb{Z}$.
Although $Sol_0^4$ has no lattice subgroups, any semidirect product
$\mathbb{Z}^3\rtimes_\theta\mathbb{Z}$ where $\theta$ 
has a pair of complex conjugate roots 
which are not roots of unity is a lattice in $Isom(\mathbb{S}ol^4_0)$,
by Theorem 8.3.
Such groups are again classified by the characteristic polynomial 
and an ideal class.

If the geometry is $\mathbb{S}ol^4_1$ then $\sqrt\pi\cong \Gamma_q$ 
for some $q\geq1$,
and either $\nu=\sqrt\pi$ or $\nu/\sqrt\pi=Z/2Z$ and 
$\nu/\zeta\sqrt\pi\cong\mathbb{Z}^2\rtimes_{-I}(Z/2Z)$.
(In the latter case $\nu$ is uniquely determined by $q$.)
Moreover $\pi$ is orientable if and only if $\beta_1(\pi)=1$.
In particular, $\mathrm{Ker}(w_1(\pi))\cong\nu\rtimes_\theta\mathbb{Z}$, 
for some $\theta\in Aut(\nu)$.
Let $A=\theta|_{\sqrt\pi}$ and let $\overline{A}$ be its image in 
$Aut(\sqrt\pi/\zeta\sqrt\pi)\cong GL(2,\mathbb{Z})$.
If $\nu=\sqrt\pi$ the translation subgroup $\pi\cap Sol^4_1$ is 
$T=\Gamma_q\rtimes_B\mathbb{Z}$, 
where $B=A$ or $A^2$ is the least nontrivial power of $A$ 
such that both eigenvalues of $\overline A$ are positive.
If $\nu\not=\sqrt\pi$ the conjugacy class of $\overline A$ 
is only well-defined up to sign.
If moreover $\pi/\nu\cong D$ then $\overline{A}$ is conjugate to its inverse, 
and so $det(\overline{A})=1$, since $\overline{A}$ has infinite order.
We can then choose $\theta$ and hence $A$ so that 
$T=\sqrt\pi\rtimes_A\mathbb{Z}$.
The quotient $\pi/T$ is a subgroup of $D_8$,
since $Isom(\mathbb{S}ol^4_1)\cong{Sol_1^4}\rtimes{D_8}$.

Conversely, any torsion-free group with a subgroup of index $\leq2$
which is such a semidirect product $\nu\rtimes_\theta\mathbb{Z}$ 
(with $[\nu:\Gamma_q]\leq2$ and $\nu$ as above) 
and which is not virtually nilpotent
is a lattice in $Isom(\mathbb{S}ol^4_1)$,
by an argument extending that of Theorem 8.7.
(See \cite{[Hi07]}.)
The conjugacy class of $\theta$ is determined up to a finite ambiguity 
by the characteristic polynomial of $\overline A$.
The $\mathbb{S}ol_1^4$-lattices are classified in \cite{[LT15],[VT19]}.
In particular, it is shown there that each of the possible 
isomorphism classes of subgroups of $D_8$
is realized as $\pi/T$ for some $\pi$.

In the remaining case $\mathbb{S}ol^3\times\mathbb{E}^1$
the subgroup $\nu$ is one of the four flat 3-manifold groups
$\mathbb{Z}^3$, $\mathbb{Z}^2\rtimes_{-I}\mathbb{Z}$, $B_1$ or $B_2$, 
and $\theta|_{\sqrt\nu}$ has distinct real eigenvalues, one being $\pm1$.
The index of the translation subgroup $\pi\cap (Sol^3\times R)$ 
in $\pi$ divides 8. 
(Note that $Isom(\mathbb{S}ol^3\times\mathbb{E}^1)$ has 16 components.)
Conversely any torsion-free group with a subgroup of index $\leq2$
which is such a semidirect product $\nu\rtimes_\theta\mathbb{Z}$ 
is a lattice in $Isom(Sol^3\times \mathbb{E}^1)$, 
by an argument extending that of Theorem 8.3.
(See \cite{[Hi07]}.)
The semidirect products $N\rtimes_\theta\mathbb{Z}$
(with $N=\nu$ or a $\mathbb{S}ol^3$-group)
may be classified in terms of conjugacy classes in $Aut(N)$.
(See also \cite{[Cb]}.)

Every $\mathbb{N}il^4$- or $\mathbb{S}ol^3\times\mathbb{E}^1$-manifold
has an essentially unique Seifert fibration.
The general fibre is a torus, 
and so $\pi$ has an unique normal subgroup $A\cong\mathbb{Z}^2$ 
such that $\pi/A$ is a flat 2-orbifold group.
It is easy to see that $A\leq\pi'$.
The action has infinite image in $Aut(A)\cong{GL(2,\mathbb{Z})}$,
for these geometries. 
Since $GL(2,\mathbb{Z})$ is virtually free 
the base orbifold $B$ must itself fibre over $S^1$ or 
the reflector interval $\mathbb{I}$.
Moreover, $B=T$ if and only if $\beta_1(\pi)=2$,
while $B=Kb,\mathbb{A}$ or $\mathbb{M}b$ if and only if $\beta_1(\pi)=1$.
In the $\mathbb{N}il^4$ case, if $B=T,\mathbb{A}$ or $\mathbb{M}b$
then $M$ must be orientable, by Theorem 8.7.
Inspection of the tables of \cite[Chapter 7]{[De]} shows that
there are $\mathbb{N}il^4$-manifolds with $B=T,Kb,\mathbb{A},\mathbb{M}b$,
$S(2,2,2,2)$, $P(2,2)$ or $\mathbb{D}(2,2)$.
In \cite{[Hi18]} it is shown that the same seven bases may be
realized by $\mathbb{S}ol^3\times\mathbb{E}^1$-manifolds.
Moreover, if $M$ is a $\mathbb{S}ol^3\times\mathbb{E}^1$-manifold 
with $\pi/\pi'$ finite then $B=S(2,2,2,2)$, $P(2,2)$ or $\mathbb{D}(2,2)$,
and such manifolds 
may be classified in terms of certain matrices in $GL(2,\mathbb{Z})$.

We mention briefly another aspect of these groups.
All $\mathbb{E}^4$-, $\mathbb{N}il^3\times\mathbb{E}^1$-, $\mathbb{N}il^4$-,
$\mathbb{S}ol^4_0$-, $\mathbb{S}ol^4_1$- and 
$\mathbb{S}ol^3\times\mathbb{E}^1$-lattices are arithmetic.
However, if $m\not=n$ then no $\mathbb{S}ol^4_{m,n}$-lattice is arithmetic 
\cite[Theorem B]{[Ts20]}.

\section{Diffeomorphism}

Geometric 4-manifolds of solvable Lie type are infrasolvmanifolds
(see \S3 of Chapter 7), and infrasolvmanifolds are the total spaces 
of orbifold bundles with infranilmanifold fibre and flat base, by Theorem 7.2.
Baues showed that infrasolvmanifolds are determined up to
diffeomorphism by their fundamental groups \cite{[Ba04]}.
In dimensions $\leq3$ this follows from standard results of low dimensional
topology. We shall show that related arguments also cover most
4-dimensional orbifold bundle spaces.
The following theorem extends the main result of \cite{[Cb]}
(in which it was assumed that $\pi$ is not virtually nilpotent).

\begin{theorem}
Let $M$ and $M'$ be $4$-manifolds which are total spaces of orbifold bundles 
$p:M\to B$ and $p':M'\to B'$ with fibres infranilmanifolds $F$ and $F'$
(respectively) and bases flat orbifolds,
and suppose that $\pi_1(M)\cong\pi_1(M')\cong\pi$.
If $\pi$ is virtually abelian and $\beta_1(\pi)=1$ assume that $\pi$ is
orientable.
Then $M$ and $M'$ are diffeomorphic.
\end{theorem}

\begin{proof}
We may assume that $d=dim(B)\leq d'=dim(B')$.
Suppose first that $\pi$ is not virtually abelian
or virtually nilpotent of class 2.
Then all subgroups of finite index in $\pi$ have
$\beta_1\leq2$, and so $1\leq d\leq d'\leq2$.
Moreover $\pi$ has a characteristic nilpotent 
subgroup $\tilde\nu$ such that $h(\pi/\tilde\nu)=1$, by Theorems 1.5 and 1.6.
Let $\nu$ be the preimage in $\pi$ of the maximal finite normal
subgroup of $\pi/\tilde\nu$.
Then $\nu$ is a characteristic virtually nilpotent subgroup 
(with $\sqrt\nu=\tilde\nu$) and $\pi/\nu\cong\mathbb{Z}$ or $D$.
If $d=1$ then $\pi_1(F)=\nu$ and $p:M\to B$ induces this isomorphism.
If $d=2$ the image of $\nu$ in $\pi_1^{orb}(B)$ is normal.
Hence there is an orbifold map $q$ from $B$ to the circle 
$S^1$ or the reflector interval $\mathbb{I}$ such that
$qp$ is an orbifold bundle projection.
A similar analysis applies to $M'$.
In either case, $M$ and $M'$ are canonically mapping tori or  
unions of two twisted $I$-bundles, and the theorem follows via standard
3-manifold theory.

If $\pi$ is virtually nilpotent it is realized 
by an infranilmanifold $M_0$ \cite{[De]}.
Hence we may assume that $M'=M_0$, $d'=4$, $h(\sqrt\pi)=4$ and 
${\sqrt\pi}'\cong\mathbb{Z}$ or 1.
If $d=0$ or 4 then $M$ is also an infranilmanifold and the result is clear.
If there is an orbifold bundle projection from $B$ to $S^1$ or $\mathbb{I}$
then $M$ is a mapping torus or a union of twisted $I$-bundles,
and $\pi$ is a semidirect product $\kappa\rtimes\mathbb{Z}$ 
or a generalized free product with amalgamation $G*_JH$, 
where $[G:J]=[H:J]=2$.
The model $M_0$ then has a corresponding structure as
a mapping torus or a union of twisted $I$-bundles,
and we may argue as before.

If  $\beta_1(\pi)+d>4$ then $\pi_1^{orb}(B)$ maps onto $\mathbb{Z}$,
and so $B$ is an orbifold bundle over $S^1$.
Hence the above argument applies.
If there is no such orbifold bundle projection then $d\not=1$.
Thus we may assume that $d=2$ or 3 and that $\beta_1(\pi)\leq{4-d}$.
(If moreover $\beta_1(\pi)=4-d$ and there is no such projection
then $\pi'\cap\pi_1(F)=1$ and so $\pi$ is virtually abelian.)
If $d=2$ then $M$ is Seifert fibred.
Since $M'$ is an infranilmanifold (and $\pi$ cannot be one of the
three exceptional flat 4-manifold groups $G_6\rtimes_\theta\mathbb{Z}$ 
with $\theta=j$, $cej$ or $abcej$) it is also Seifert fibred, 
and so $M$ and $M'$ are diffeomorphic, by \cite{[Vo77]}.

If $d=3$ then $\pi_1(F)\cong\mathbb{Z}$.
The group $\pi$ has a normal subgroup $K$ such that 
$\pi/K\cong\mathbb{Z}$ or $D$, by Lemma 3.14.
If $\pi_1(F)<K$ then $\pi_1^{orb}(B)$ maps onto $\mathbb{Z}$ 
or $D$ and we may argue as before.
Otherwise $\pi_1(F)\cap K=1$, 
since $\mathbb{Z}$ and $D$ have no nontrivial
finite normal subgroups, and so $\pi$ is virtually abelian.
If $\beta_1(\pi)=1$ then $\pi_1(F)\cap\pi'=1$ 
(since $\pi/K$ does not map onto $\mathbb{Z}$) 
and so $\pi_1(F)$ is central in $\pi$.
It follows that $p$ is the orbit map of an $S^1$-action on $M$.
Once again, the model $M_0$ has an $S^1$-action inducing the same
orbifold fundamental group sequence.
Orientable 4-manifolds with $S^1$-action are determined up to diffeomorphism 
by the orbifold data and an Euler class corresponding to the central extension 
of $\pi_1^{orb}(B)$ by $\mathbb{Z}$ \cite{[Fi78]}.
Thus $M$ and $M'$ are diffeomorphic.
It is not difficult to determine the maximal infinite cyclic normal subgroups 
of the flat 4-manifold groups $\pi$ with $\beta_1(\pi)=0$,
and to verify that in each case the quotient maps onto $D$.
\end{proof}

It is highly probable that the arguments of Fintushel can be extended 
to all 4-manifolds which admit smooth $S^1$-actions,
and the theorem is surely true without any restrictions on $\pi$.
(Note that the algebraic argument of the final sentence of Theorem 8.11
does not work for nine of the 30 nonorientable flat 4-manifold groups $\pi$ 
with $\beta_1(\pi)=1$.)

If $\pi$ is orientable then it is realized geometrically and determines the 
total space of such an orbifold bundle up to diffeomorphism.
Hence orientable smooth 4-manifolds admitting such orbifold fibrations are 
diffeomorphic to geometric 4-manifolds of solvable Lie type. 
Is this also so in the nonorientable case?

%% file: m5-9.tex
\chapter{The other aspherical geometries}

The aspherical geometries of nonsolvable type which are realizable
by closed 4-manifolds are the ``mixed" geometries
$\mathbb{H}^2\times\mathbb{E}^2$, $\widetilde{\mathbb{SL}}\times\mathbb{E}^1$, 
$\mathbb{H}^3\times\mathbb{E}^1$ and the ``semisimple" geometries 
$\mathbb{H}^2\times\mathbb{H}^2$, $\mathbb{H}^4$ and $\mathbb{H}^2(\mathbb{C})$.
(We shall consider the geometry $\mathbb{F}^4$ briefly in Chapter 13.)
Closed $\mathbb{H}^2\times\mathbb{E}^2$- or 
$\widetilde{\mathbb{SL}}\times\mathbb{E}^1$-manifolds are Seifert fibred,
have Euler characteristic 0 and their fundamental groups have Hirsch-Plotkin 
radical $\mathbb{Z}^2$. 
In \S1 and \S2 we examine to what extent these properties characterize such 
manifolds and their fundamental groups.
Closed $\mathbb{H}^3\times\mathbb{E}^1$-manifolds also have 
Euler characteristic 0, 
but we have only a conjectural characterization 
of their fundamental groups (\S3).
In \S4 we determine the mapping tori of self homeomorphisms of geometric 
3-manifolds which admit one of these mixed geometries. 
(We return to this topic in Chapter 13.)
In \S5 we consider the three semisimple geometries.
All closed 4-manifolds with product geometries other than 
$\mathbb{H}^2\times\mathbb{H}^2$ are finitely covered by cartesian products.
We characterize the fundamental groups of
$\mathbb{H}^2\times\mathbb{H}^2$-manifolds with this property;
there are also ``irreducible"
$\mathbb{H}^2\times\mathbb{H}^2$-manifolds which are not virtually products.
Relatively little is known about manifolds admitting one 
of the two hyperbolic geometries.

Although it is not yet known whether the disk embedding theorem 
holds over lattices for such geometries, 
we can show that the fundamental group and Euler characteristic 
determine the manifold up to $s$-cobordism (\S6).
Moreover an aspherical closed 4-manifold which is finitely covered 
by a geometric manifold is homotopy equivalent to a geometric manifold.

\section{Aspherical Seifert fibred 4-manifolds}

In Chapter 8 we saw that if $M$ is a closed 4-manifold with fundamental
group $\pi$ such that $\chi(M)=0$ and $h(\sqrt\pi)\geq3$ then $M$ is 
homeomorphic to an infrasolvmanifold.
Here we shall show that if $\chi(M)=0$, $h(\sqrt\pi)=2$ and
$[\pi:\sqrt\pi]=\infty$ then $M$ is homotopy equivalent to a 4-manifold which
is Seifert fibred over a hyperbolic 2-orbifold.
(We shall consider the case when $\chi(M)=0$, $h(\sqrt\pi)=2$ and
$[\pi:\sqrt\pi]<\infty$ in Chapter 10.)

\begin{theorem} 
Let $M$ be a $PD_4$-complex with fundamental group $\pi$. 
If $\chi(M)=0$, $\pi$ has an elementary amenable normal subgroup $\rho$ 
with $h(\rho)=2$ and  $H^2(\pi;\mathbb{Z}[\pi])=0$ 
then $M$ is aspherical and $\rho$ is virtually abelian.
\end{theorem}

\begin{proof} 
Since $\pi$ has one end, by Corollary 1.15.1,
$\beta_1^{(2)}(\pi)=0$, by Theorem 2.3, 
and $H^2(\pi;\mathbb{Z}[\pi])=0$,
$M$ is aspherical, by Corollary 3.5.2.
In particular, $\rho$ is torsion-free and $[\pi:\rho]=\infty$.

Since $\rho$ is torsion-free elementary amenable and $h(\rho)=2$ 
it is virtually solvable, by Theorem 1.11. 
Therefore $A=\sqrt\rho$ is nontrivial, and as it is characteristic in
$\rho$ it is normal in $\pi$.
Since $A$ is torsion-free and $h(A)\leq 2$ it is abelian, by Theorem 1.5.

If $h(A)=1$ then $A$ is isomorphic to a subgroup of $\mathbb{Q}$ 
and the homomorphism from $B=\rho/A$ to $Aut(A)$ induced by 
conjugation in $\rho$ is injective. 
Since $Aut(A)$ is isomorphic to a subgroup of $\mathbb{Q}^\times $
and $h(B)=1$ either $B\cong\mathbb{Z}$ or $B\cong\mathbb{Z}\oplus(Z/2Z)$.
We must in fact have $B\cong\mathbb{Z}$, since $\rho$ is torsion-free. 
Moreover $A$ is not finitely generated and the centre of $\rho$ is trivial.
The quotient group $\pi/A$ has one end as the image of $\rho$ is an infinite
cyclic normal subgroup of infinite index.

As $A$ is a characteristic subgroup every automorphism of $\rho$ restricts 
to an automorphism of $A$. 
This restriction from $Aut(\rho)$ to $Aut(A)$ is an epimorphism,      
with kernel isomorphic to $A$, and so $Aut(\rho)$ is solvable.
Let $C=C_\pi (\rho)$ be the centralizer of $\rho$ in $\pi$.
Then $C$ is nontrivial, for otherwise $\pi$ would be isomorphic to a subgroup
of $Aut(\rho)$ and hence would be virtually poly-$Z$.
But then $A$ would be finitely generated, $\rho$ would be virtually abelian and $h(A)=2$.
Moreover $C\cap\rho=\zeta\rho=1$, so $C\rho\cong C\times\rho$ and
$c.d.C+c.d.\rho=c.d.C\rho\leq c.d.\pi=4$. 
The quotient group $\pi/C\rho$ is
isomorphic to a subgroup of $Out(\rho)$.

If $c.d.C\rho\leq 3$ then as $C$ is nontrivial and $h(\rho)=2$ we must have
$c.d.C=1$ and $c.d.\rho=h(\rho)=2$. 
Therefore $C$ is free and $\rho$ is of type
$FP$ \cite{[Kr86]}. 
By Theorem 1.13 $\rho$ is an ascending HNN group with base a
finitely generated subgroup of $A$ and so has a presentation 
$\langle a,t\mid tat^{-1} =a^n\rangle$ for some nonzero integer $n$.
We may assume $|n|>1$, as $\rho$ is not virtually abelian.
The subgroup of $Aut(\rho)$ represented by $(n-1)A$ consists of inner automorphisms.
Since $n>1$ the quotient $A/(n-1)A\cong Z/(n-1)Z$ is finite,
and as $Aut(A)\cong \mathbb{Z}[1/n]^\times $ it follows that $Out(\rho)$ is virtually abelian.
Therefore $\pi$ has a subgroup $\sigma$ of finite index which contains $C\rho$
and such that $\sigma/C\rho$ is a finitely generated free abelian group,
and in particular $c.d.\sigma/C\rho$ is finite.
As $\sigma$ is a $PD_4 $-group it follows from \cite[Theorem 9.11]{[Bi]} 
that $C\rho$ is a $PD_3 $-group and hence that $\rho$ is a $PD_2 $-group.
We reach the same conclusion if $c.d.C\rho=4$, for then $[\pi:C\rho]$ is finite,
by Strebel's Theorem, and so $C\rho$ is a $PD_4 $-group.
As a solvable $PD_2 $-group is virtually $\mathbb{Z}^2$ our original assumption
must have been wrong. 

Therefore $h(A)=2$.
As every finitely generated subgroup of $\rho$ is either 
isomorphic to $\mathbb{Z}\rtimes_{-1}\mathbb{Z}$ 
or is abelian $[\rho:A]\leq 2$. 
\end{proof}

If $h(\rho)=2$, $\rho$ is torsion-free and $[\pi:\rho]=\infty$ then
$H^2(\pi;\mathbb{Z}[\pi])=0$, by Theorem 1.17.
Can the latter hypothesis in the above theorem
be replaced by ``$[\pi:\rho]=\infty$"?
Some such hypothesis is needed,
for if $M=S^2\times{T}$ then $\chi(M)=0$ and $\pi\cong\mathbb{Z}^2$.

\begin{theorem}
Let $M$ be a $PD_4$-complex with fundamental group $\pi$.
If $h(\sqrt\pi)=2$, $[\pi:\sqrt\pi]=\infty$ and $\chi(M)=0$ then
$M$ is aspherical and $\sqrt\pi\cong \mathbb{Z}^2$.
\end{theorem}

\begin{proof} As $H^s(\pi;\mathbb{Z}[\pi])=0$ for $s\leq2$, 
by Theorem 1.17, $M$ is aspherical, by Theorem 9.1.
We may assume henceforth that $\sqrt\pi$ is a torsion-free abelian group 
of rank 2 which is not finitely generated. 

Suppose first that $[\pi:C]=\infty$, where $C=C_\pi (\sqrt\pi)$. 
Then $c.d.C\leq 3$, by Strebel's Theorem.
Since $\sqrt\pi$ is not finitely generated $c.d.\sqrt\pi=h(\sqrt\pi)+1=3$ 
\cite[Theorem 7.14]{[Bi]}.
Hence $C=\sqrt\pi$ \cite[Theorem 8.8]{[Bi]}, 
so the homomorphism from $\pi/\sqrt\pi$ to
$Aut(\sqrt\pi)$ determined by conjugation in $\pi$ is a monomorphism.
Since $\sqrt\pi$ is torsion-free abelian of rank 2, 
$Aut(\sqrt\pi)$ is isomorphic to a subgroup of $GL(2,\mathbb{Q})$ 
and therefore any torsion subgroup of $Aut(\sqrt\pi)$ is finite,
by Corollary 1.3.1.
Thus if $\pi'\sqrt\pi/\sqrt\pi$ is a torsion group $\pi'\sqrt\pi$ 
is elementary amenable and so $\pi$ is itself elementary amenable, 
contradicting our assumption. 
Hence we may suppose that there is an element $g$ in $\pi'$ 
which has infinite order modulo $\sqrt\pi$.
The subgroup $\langle \sqrt\pi,g\rangle$ generated by $\sqrt\pi$ 
and $g$ is an extension of $\mathbb{Z}$ by $\sqrt\pi$ and 
has infinite index in $\pi$, for otherwise $\pi$ would be virtually solvable.
Hence $c.d.\langle \sqrt\pi,g\rangle=3=h(\langle \sqrt\pi,g\rangle)$, 
by Strebel's Theorem.
By \cite[Theorem 7.15]{[Bi]}, 
$L=H_2 (\sqrt\pi;\mathbb{Z})$ is the underlying abelian group 
of a subring $\mathbb{Z}[m^{-1} ]$ of $\mathbb{Q}$,
and the action of $g$ on $L$ is multiplication by a rational number $a/b$, 
where $a$ and $b$ are relatively prime and $ab$ and $m$ 
have the same prime divisors.
But $g$ acts on $\sqrt\pi$ as an element of $GL(2,\mathbb{Q})'\leq SL(2,\mathbb{Q})$.
Since $L=\sqrt\pi\wedge \sqrt\pi$ \cite[11.14.16]{[Ro]}, 
$g$ acts on $L$ via $det(g)=1$.
Therefore $m=1$ and so $L$ must be finitely generated.
But then $\sqrt\pi$ must also be finitely generated, again contradicting our assumption.

Thus we may assume that $C$ has finite index in $\pi$.
Let $A <\sqrt\pi$ be a subgroup of $\sqrt\pi$ which is free abelian of rank 2.
Then $A_1 $ is central in $C$ and $C/A$ is finitely presentable.
Since $[\pi:C]$ is finite $A $ has only finitely many distinct conjugates 
in $\pi$,
and they are all subgroups of $\zeta C$. Let $N$ be their product. 
Then $N$ is a finitely generated torsion-free abelian normal subgroup of $\pi$
and $2\leq h(N)\leq h(\sqrt C)\leq h(\sqrt\pi)=2$.
An LHSSS argument gives $H^2(\pi/N;\mathbb{Z}[\pi/N])\cong\mathbb{Z}$,
and so $\pi/N$ is virtually a $PD_2$-group, by Bowditch's Theorem.
Since $\sqrt\pi/N$ is a torsion group it must be finite, 
and so $\sqrt\pi\cong\mathbb{Z}^2$.
\end{proof}

\begin{cor}
Let $M$ be a closed $4$-manifold with fundamental group $\pi$. 
Then $M$ is homotopy equivalent to one
which is Seifert fibred with general fibre $T$ or $Kb$ over a hyperbolic 
2-orbifold if and only if $h(\sqrt\pi)=2$, 
$[\pi:\sqrt\pi]=\infty$ and $\chi(M)=0$.
\end{cor}

\begin{proof} This follows from the theorem together with 
Theorem 7.3.
\end{proof}

\section{The Seifert geometries: $\mathbb{H}^2\times\mathbb{E}^2$ and 
$\widetilde{\mathbb{SL}}\times\mathbb{E}^1$}

A manifold with geometry $\mathbb{H}^2\times\mathbb{E}^2$ or 
$\widetilde{\mathbb{SL}}\times\mathbb{E}^1$ 
is Seifert fibred with base a hyperbolic orbifold. 
However not all such Seifert fibred 4-manifolds are geometric.
We shall show that geometric Seifert fibred 4-manifolds 
may be characterized in terms of their fundamental groups.
With \cite{[Vo77]}, Theorems 9.5 and 9.6 imply the main result of \cite{[Ke]},
which is Corollary 9.6.2 below.

\begin{theorem} 
Let $M$ be a closed $\mathbb{H}^3\times\mathbb{E}^1$-,
$\widetilde{\mathbb{SL}}\times\mathbb{E}^1$- or $\mathbb{H}^2\times\mathbb{E}^2$-manifold. 
Then $M$ has a finite covering space which is diffeomorphic to a product
$N\times S^1$.
\end{theorem}

\begin{proof} 
If $M$ is an $\mathbb{H}^3\times\mathbb{E}^1$-manifold then $\pi=\pi_1(M)$ 
is a discrete cocompact subgroup of $G=Isom(\mathbb{H}^3\times\mathbb{E}^1)$.
The radical of this group is $Rad(G)\cong\mathbb{R}$, 
and $G_o/Rad(G)\cong PSL(2,\mathbb{C})$,
where $G_o$ is the component of the identity in $G$.
Since $PSL(2,\mathbb{C})$ has no compact factor,
$A=\pi\cap Rad(G)$ is a lattice subgroup \cite[Proposition 8.27]{[Rg]}.
Since $\mathbb{R}/A$ is compact the image of $\pi/A$ 
in $Isom(\mathbb{H}^3)$ is again a discrete cocompact subgroup.
Hence $\sqrt\pi=A\cong\mathbb{Z}$.

On passing to a 2-fold covering space, if necessary, we may assume
that $\pi\leq Isom(\mathbb{H}^3)\times\mathbb{R}$ 
and (hence) $\zeta\pi=\sqrt\pi$.
Projection to the second factor maps $\sqrt\pi$ 
monomorphically to $\mathbb{R}$.
Hence on passing to a further finite covering space, if necessary,
we may assume that $\pi\cong\nu\times\mathbb{Z}$, 
where $\nu=\pi/\sqrt\pi\cong\pi_1(N)$
for some closed orientable $\mathbb{H}^3$-manifold $N$.
(Note that we do {\it not} claim that $\pi=\nu\times\mathbb{Z}$ 
as a subgroup of $PSL(2,\mathbb{C})\times\mathbb{R}$.)
The foliation of $H^3\times\mathbb{R}$ by lines is preserved by $\pi$,
and so induces an $S^1$-bundle structure on $M$, with base $N$.
As such bundles (with aspherical base) are determined by their
fundamental groups, $M$ is diffeomorphic to $N\times S^1$.

Similar arguments apply if the geometry is $\mathbb{X}^4=
\mathbb{H}^2\times\mathbb{E}^2$ or 
$\widetilde{\mathbb{SL}}\times\mathbb{E}^1$.
If $G=Isom(\mathbb{X}^4)$ then $Rad(G)\cong\mathbb{R}^2$, 
and $PSL(2,\mathbb{R})$ is a cocompact subgroup of $G_o/Rad(G)$. 
The intersection $A=\pi\cap{Rad(G)}$ is again a lattice subgroup,
and $\pi/A$ has a subgroup of finite index which 
is a discrete cocompact subgroup of $PSL(2,\mathbb{R})$.
Hence $\sqrt\pi=A\cong\mathbb{Z}^2$.
Moreover, on passing to a finite covering space we may assume that
$\zeta\pi=\sqrt\pi$ and $\pi/\sqrt\pi$ is a $PD_2$-group.
If $\mathbb{X}^4=\mathbb{H}^2\times\mathbb{E}^2$ then 
projection to the second factor maps $\sqrt\pi$ monomorphically 
and $\pi$ preserves the foliation of $H^2\times\mathbb{R}^2$ by planes. 
If $\mathbb{X}^4=\widetilde{\mathbb{SL}}\times\mathbb{E}^1$
then $\sqrt\pi\cap{Isom}(\widetilde{\mathbb{SL}})$ must be nontrivial,
since ${Isom}(\widetilde{\mathbb{SL}})$ has no subgroups which are 
$PD_2$-groups. (See \cite[page 466]{[Sc83']}.)
Hence $\pi$ is virtually a product $\nu\times\mathbb{Z}$ 
with $\nu=\pi_1(N)$ for some 
closed orientable $\widetilde{\mathbb{SL}}$-manifold $N$.
In each case, $M$ is virtually a product.
\end{proof}

There may not be such a covering which is geometrically
a cartesian product.
Let $\nu$ be a discrete cocompact subgroup of $Isom(\mathbb{X})$ 
where $\mathbb{X}=\mathbb{H}^3$ or $\widetilde{\mathbb{SL}}$
which admits an epimorphism $\alpha:\nu\to\mathbb{Z}$.
Define a homomorphism 
$\theta:\nu\times\mathbb{Z}\to Isom(\mathbb{X}\times\mathbb{E}^1)$
by $\theta(g,n)(x,r)=(g(x),r+n+\alpha(g)\sqrt2)$ for all $g\in\nu$, 
$n\in\mathbb{Z}$, $x\in X$ and $r\in\mathbb{R}$.
Then $\theta$ is a monomorphism onto a discrete subgroup
which acts freely and cocompactly on $X\times\mathbb{R}$, 
but its image in $E(1)$ has rank 2.

\begin{lemma}
Let $\pi$ be a finitely generated group with normal subgroups $A\leq N$
such that $A$ is free abelian of rank $r$, $[\pi:N]<\infty$ and $N\cong
A\times{N/A}$. Then there is a homomorphism $f:\pi\to{E(r)}$ with image
a discrete cocompact subgroup and such that $f|_A$ is injective.
\end{lemma}

\begin{proof} 
Let $G=\pi/N$ and $M=N^{ab}\cong{A}\oplus(N/AN')$.
Then $M$ is a finitely generated $\mathbb{Z}[G]$-module and the
image of $A$ in $M$ is a $\mathbb{Z}[G]$-submodule. 
Extending coefficients to the rationals $\mathbb{Q}$ gives a natural inclusion 
$\mathbb{Q}A\leq\mathbb{Q}M$, since $A$ is a direct summand of $M$ 
(as an abelian group), 
and $\mathbb{Q}A$ is a $\mathbb{Q}[G]$-submodule of $\mathbb{Q}M$. 
Since $G$ is finite $\mathbb{Q}[G]$ is semisimple,
and so $\mathbb{Q}A$ is a $\mathbb{Q}[G]$-direct summand of $\mathbb{Q}M$.
Let $K$ be the kernel of the homomorphism from $M$ to $\mathbb{Q}A$
determined by a splitting homomorphism from $\mathbb{Q}M$ to $\mathbb{Q}A$,
and let $\tilde{K}$ be the preimage of $K$ in $\pi$.
Then $K$ is a $\mathbb{Z}[G]$-submodule of $M$ and $M/K\cong\mathbb{Z}^r$,
since it is finitely generated and torsion-free of rank $r$.
Moreover $\tilde K$ is a normal subgroup of $\pi$ and $A\cap\tilde{K}=1$.
Hence $H=\pi/\tilde{K}$ is an extension of $G$ by $M/K$ and
$A$ maps injectively onto a subgroup of finite index in $H$.
Let $T$ be the maximal finite normal subgroup of $H$. 
Then $H/T$ is isomorphic to a discrete cocompact subgroup of $E(r)$,
and the projection of $\pi$ onto $H/T$ is clearly injective on $A$.
\end{proof}

\begin{theorem} 
Let $M$ be a closed $4$-manifold with fundamental group $\pi$. Then
the following are equivalent:
\begin{enumerate}
\item 
$M$ is homotopy equivalent to a $\mathbb{H}^2\times\mathbb{E}^2$-manifold;

\item $\pi$ has a finitely generated infinite subgroup $\rho$ 
such that $[\pi:N_\pi(\rho)]<\infty$, $\sqrt\rho=1$, 
$\zeta C_\pi(\rho)\cong\mathbb{Z}^2$ and $\chi(M)=0$;

\item $\sqrt\pi\cong\mathbb{Z}^2$, $[\pi:\sqrt\pi]=\infty$, 
$[\pi:C_\pi (\sqrt\pi)]<\infty$, $e^\mathbb{Q}(\pi)=0$ and $\chi(M)=0$.
\end{enumerate}
\end{theorem}
                    
\begin{proof} 
If $M$ is a $\mathbb{H}^2\times\mathbb{E}^2$-manifold
it is finitely covered by $B\times{T}$, 
where $B$ is a closed hyperbolic surface.
Thus (1) implies (2),
on taking $\rho=\pi_1(B)$.

If (2) holds $M$ is aspherical and so $\pi$ is a $PD_4$-group,
by Theorem 9.1. 
Let $C=C_\pi(\rho)$. 
Then $C$ is also normal in $\nu=N_\pi(\rho)$, and $C\cap\rho=1$, 
since $\sqrt\rho=1$.
Hence $\rho\times C\cong \rho.C\leq \pi$.
Now $\rho$ is nontrivial.
If $\rho$ were free then an argument using the LHSSS for 
$H^*(\nu;\mathbb{Q}[\nu])$ would imply that $\rho$ has two ends, 
and hence that $\sqrt\rho=\rho\cong\mathbb{Z}$.
Hence $c.d.\rho\geq 2$.
Since moreover $\mathbb{Z}^2\leq C$ we must have $c.d.\rho=c.d.C=2$ and 
$[\pi:\rho.C]<\infty$.
It follows easily that $\sqrt\pi\cong\mathbb{Z}^2$ 
and $[\pi:C_\pi (\sqrt\pi)]<\infty$. 
Moreover $\pi$ has a normal subgroup $K$ of finite index which contains
$\sqrt\pi$ and is such that $K\cong\sqrt\pi\times{K/\sqrt\pi}$.
In particular, $e^\mathbb{Q}(K)=0$ and so $e^\mathbb{Q}(\pi)=0$.
Thus (2) implies (3).

If (3) holds $M$ is homotopy equivalent to a manifold which
is Seifert fibred over a hyperbolic orbifold, by Corollary 9.2.1.
Since $e^\mathbb{Q}(\pi)=0$ this manifold has a finite regular covering 
which is a product $B\times{T}$, with $\pi_1(T)=\sqrt\pi$. 
Let $H$ be the maximal solvable normal subgroup of $\pi$.
Since $\pi/\sqrt\pi$ has no infinite solvable normal subgroup $H/\sqrt\pi$ is 
finite, and since $\pi$ is torsion-free the preimage 
of any finite subgroup of $\pi/\sqrt\pi$ is $\sqrt\pi$ 
or $\mathbb{Z}\rtimes_{-1}\mathbb{Z}$.
Then $[H:\sqrt\pi]\leq2$,
$\pi_1(B)$ embeds in $\pi/H$ as a subgroup of finite index
and $\pi/H$ has no nontrivial finite normal subgroup.
Therefore there is a homomorphism $h:\pi\to Isom(\mathbb{H}^2)$ with kernel $H$ 
and image a discrete cocompact subgroup,
by the solution to the Nielsen realization problem for surfaces 
\cite{[Ke83]}.
By the lemma there is also a homomorphism $f:\pi\to{E(2)}$ 
which maps $\sqrt\pi$ to a lattice.
The homomorphism $(h,f):\pi\to Isom(\mathbb{H}^2\times\mathbb{E}^2)$
is injective, since $\pi$ is torsion-free,
and its image is discrete and cocompact.
Therefore it is a lattice, and so (3) implies (1).
\end{proof}

A similar argument may be used to characterize 
${\widetilde{\mathbb{SL}}\times\mathbb{E}^1}$-manifolds.

\begin{theorem} 
Let $M$ be a closed $4$-manifold with fundamental group $\pi$. Then
the following are equivalent:
\begin{enumerate}
\item 
$M$ is homotopy equivalent to a 
$\widetilde{\mathbb{SL}}\times\mathbb{E}^1$-manifold;

\item $\sqrt\pi\cong\mathbb{Z}^2$, ${[\pi:\sqrt\pi]}=\infty$, 
${[\pi:C_\pi (\sqrt\pi)]}<\infty$, 
$e^\mathbb{Q}(\pi)\not=0$ and $\chi(M)=0$.
\end{enumerate}
\end{theorem}

\begin{proof} 
(Sketch)
These conditions are clearly necessary.
If they hold $M$ is aspherical and $\pi$ has a normal subgroup $K$ of finite 
index which is a central extension of a $PD_2^+$-group $G$ by $\sqrt\pi$.
Let $e(K)\in H^2(G;\mathbb{Z}^2)\cong\mathbb{Z}^2$ 
be the class of this extension.
There is an epimorphism $\lambda:\mathbb{Z}^2\to\mathbb{Z}$ 
such that $\lambda_\sharp(e(K))=0$ in $H^2(G;\mathbb{Z})$,
and so $K/\mathrm{Ker}(\lambda)\cong G\times\mathbb{Z}$.
Hence $K\cong\nu\times\mathbb{Z}$, 
where $\nu$ is a $\widetilde{\mathbb{SL}}$-manifold group.
Let $A<\sqrt\pi$ be an infinite cyclic normal subgroup of $\pi$ which
maps onto $K/\nu$, and let $H$ be the preimage in $\pi$ of the maximal finite
normal subgroup of $\pi/A$.
Then $[H:A]\leq2$, $\nu$ embeds in $\pi/H$ as a subgroup of finite index
and $\pi/H$ has no nontrivial finite normal subgroup.
Hence $\pi/H$ is a $\widetilde{\mathbb{SL}}$-orbifold group 
\cite[Satz 2.1]{[ZZ82]}.
(This is another application of \cite{[Ke83]}).
A homomorphism $f:\pi\to{E(1)}$ which is injective on $H$ and with image a
lattice may be constructed and the sufficiency of these conditions 
may then be established as in Theorem 9.5.
\end{proof}


\begin{cor} 
A group $\pi$ is the fundamental group of a closed 
$\mathbb{H}^2\times\mathbb{E}^2$- or
$\widetilde{\mathbb{SL}}\times\mathbb{E}^1$-manifold if and only if 
it is a $PD_4$-group, $\sqrt\pi\cong\mathbb{Z}^2$
and the action $\alpha$ has finite image in $GL(2,\mathbb{Z})$.
The geometry is $\mathbb{H}^2\times\mathbb{E}^2$ if and only if 
$e^\mathbb{Q}(\pi)=0$.
\qed
\end{cor}

\begin{cor} 
{\rm[Ke]}\qua
An aspherical Seifert fibred 4-manifold is geometric if and only if
it is finitely covered by a geometric 4-manifold.
\qed
\end{cor}

A closed $4$-manifold $M$ is an $\mathbb{H}^2\times\mathbb{E}^2$-manifold
if and only if it is both Seifert fibred and also the total space of an
orbifold bundle over a flat $2$-orbifold and with general fibre
a hyperbolic surface, for the two projections determine a direct product
splitting of a subgroup of finite index, and so $e^\mathbb{Q}(\pi)=0$.

Similarly, $M$ is a product $T\times{B}$ with $\chi(B)<0$ if and only if 
it fibres both as a torus bundle and as a bundle with hyperbolic fibre.

\section{$\mathbb{H}^3\times\mathbb{E}^1$-manifolds}

An argument related to that of Theorem 9.5 
(using the Virtual Fibration Theorem \cite{[Ag13]},
and using Mostow rigidity instead of \cite{[Ke83]}) 
shows that a 4-manifold $M$ is homotopy equivalent 
to an $\mathbb{H}^3\times\mathbb{E}^1 $-manifold 
if and only if $\chi(M)=0$ and $\pi=\pi_1(M)$ has a normal subgroup of 
finite index which is isomorphic to $\rho\times\mathbb{Z}$,
where $\rho$ has an infinite $FP_2$ normal subgroup $\nu$ of infinite index,
but has no noncyclic abelian subgroup.
Moreover, $\nu$ is then a $PD_2$-group, 
and every torsion free group $\pi$ with such subgroups is 
the fundamental group of an $\mathbb{H}^3\times\mathbb{E}^1$-manifold.

If every $PD_3 $-group is the fundamental group of a closed 
3-manifold we could replace the condition on $\nu$ by the simpler 
condition that $\rho$ have one end.
For then $M$ would be aspherical and hence $\rho$ would be a $PD_3 $-group.
An aspherical 3-manifold whose fundamental group
has no noncyclic abelian subgroup is atoroidal, and hence hyperbolic 
\cite{[B-P]}. 

The foliation of $H^3\times\mathbb{R}$ by copies of $H^3$ induces a
codimension 1 foliation of any closed $\mathbb{H}^3\times\mathbb{E}^1$-manifold.
If all the leaves are compact then the manifold is either a mapping torus 
or the union of two twisted $I$-bundles. 
Is this always the case?

\begin{theorem}  
Let $M$ be a closed $\mathbb{H}^3\times\mathbb{E}^1$-manifold.
If $\zeta\pi\cong\mathbb{Z}$ then $M$ is homotopy equivalent to a
mapping torus of a self homeomorphism of an $\mathbb{H}^3$-manifold; 
otherwise $M$ is homotopy equivalent to the union
of two twisted $I$-bundles over $\mathbb{H}^3$-manifold bases.
\end{theorem}

\begin{proof} 
There is a homomorphism $\lambda:\pi\to{E(1)}$ with image a discrete cocompact 
subgroup and with $\lambda(\sqrt\pi)\not=1$, by Lemma 9.4.
Let $K=\mathrm{Ker}(\lambda)$.
Then $K\cap\sqrt\pi=1$, so $K$ is isomorphic to a subgroup of finite index in
$\pi/\sqrt\pi$.
Therefore $K\cong \pi_1(N)$ for some closed $\mathbb{H}^3$-manifold,
since it is torsion-free.
If $\zeta\pi=\mathbb{Z}$ then 
$\mathrm{Im}(\lambda)\cong\mathbb{Z}$ (since $\zeta D=1$);
if $\zeta\pi=1$ then $\mathrm{Im}(\lambda)\cong D$.
The theorem now follows easily.
\end{proof}

\section{Mapping tori}

In this section we shall use 3-manifold theory to characterize 
mapping tori with one of the geometries $\mathbb{H}^3\times\mathbb{E}^1$, 
$\widetilde{\mathbb{SL}}\times\mathbb{E}^1$ or $\mathbb{H}^2\times\mathbb{E}^2$. 
\begin{theorem} 
Let $\phi$ be a self homeomorphism of a 
closed $3$-manifold $N$ which admits the geometry $\mathbb{H}^2\times\mathbb{E}^1$ 
or $\widetilde{\mathbb{SL}}$. 
Then the mapping torus $M(\phi)=N\times_\phi S^1 $ admits the corresponding 
product geometry if and only if the outer automorphism $[\phi_* ]$ 
induced by $\phi$ has finite order. 
The mapping torus of a self homeomorphism $\phi$ of an $\mathbb{H}^3$-manifold $N$ 
admits the geometry $\mathbb{H}^3\times\mathbb{E}^1$.
\end{theorem}

\begin{proof} Let $\nu=\pi_1 (N)$ and let $t$ be an element of 
$\pi=\pi_1 (M(\phi))$ which projects to a generator of $\pi_1 (S^1)$. 
If $M(\phi)$ has geometry $\widetilde{\mathbb{SL}}\times\mathbb{E}^1 $ then 
after passing to the 2-fold covering space $M(\phi^2)$, 
if necessary, we may assume that $\pi$ is a discrete cocompact subgroup 
of $Isom(\widetilde{\mathbb{SL}})\times\mathbb{R}$. 
As in Theorem 9.3 the intersection of $\pi$ with the centre of this group is 
a lattice subgroup $L\cong\mathbb{Z}^2$.
Since the centre of $\nu$ is $\mathbb{Z}$ the image of $L$ 
in $\pi/\nu$ is nontrivial,
and so $\pi$ has a subgroup $\sigma$ of finite index which is 
isomorphic to $\nu\times\mathbb{Z}$.
In particular, 
conjugation by $t^{[\pi:\sigma]}$ induces an inner automorphism of $\nu$.

If $M(\phi)$ has geometry $\mathbb{H}^2\times\mathbb{E}^2$ 
a similar argument implies that
$\pi$ has a subgroup $\sigma$ of finite index which is isomorphic to 
$\rho\times\mathbb{Z}^2$, 
where $\rho$ is a discrete cocompact subgroup of $PSL(2,\mathbb{R})$,
and is a subgroup of $\nu$. It again follows that             
$t^{[\pi:\sigma]} $ induces an inner automorphism of $\nu$.

Conversely, suppose that $N$ has a geometry of type 
$\mathbb{H}^2\times\mathbb{E}^1$
or $\widetilde{\mathbb{SL}}$ and that $[\phi_* ]$ has finite order in $Out(\nu)$.
Then $\phi$ is homotopic to a self homeomorphism of (perhaps larger) 
finite order \cite{[Zn80]} and is therefore isotopic to 
such a self homeomorphism \cite{[Sc85],[BO91]},
which may be assumed to preserve the geometric structure \cite{[MS86]}.
Thus we may assume that $\phi$ is an isometry.
The self homeomorphism of $N\times R$ sending $(n,r)$ to $(\phi(n),r+1)$
is then an isometry for the product geometry and the mapping torus 
has the product geometry.

If $N$ is hyperbolic then $\phi$ is homotopic to an isometry of finite order,
by Mostow rigidity \cite{[Ms68]}, 
and is therefore isotopic to such an isometry \cite{[GMT03]},
so the mapping torus again has the product geometry. 
\end{proof}

A closed 4-manifold $M$ which admits an effective $T$-action with hyperbolic 
base orbifold is homotopy equivalent to such a mapping torus.
For then $\zeta\pi=\sqrt\pi$ and the LHSSS for homology gives an exact 
sequence 
\begin{equation*}
H_2 (\pi/\zeta\pi;\mathbb{Q})\to H_1 (\zeta\pi;\mathbb{Q})
\to H_1 (\pi;\mathbb{Q}).
\end{equation*}
As $\pi/\zeta\pi$ is virtually a $PD_2$-group
$H_2 (\pi/\zeta\pi;\mathbb{Q})\cong\mathbb{Q}$ or 0, so 
$\zeta\pi/\zeta\pi\cap\pi'$ has rank at least 1.
Hence $\pi\cong \nu\rtimes_\theta\mathbb{Z}$ where $\zeta\nu\cong\mathbb{Z}$, 
$\nu/\zeta\nu$ is virtually a $PD_2$-group
and $[\theta]$ has finite order in $Out(\nu)$.
If moreover $M$ is orientable then it is geometric 
(\cite{[Ue90],[Ue91]} -- see also \S7 of Chapter 7).
Note also that if $M$ is a
$\widetilde{\mathbb{SL}}\times\mathbb{E}^1$-manifold then $\zeta\pi=\sqrt\pi$
if and only if $\pi\leq Isom_o(\widetilde{\mathbb{SL}}\times\mathbb{E}^1)$.
        
Let $F$ be a closed hyperbolic surface and $\alpha:F\to F$ a pseudo-Anasov 
homeomorphism.
Let $\Theta(f,z)=(\alpha(f),\bar z)$ for all $(f,z)$ in $N=F\times S^1$. 
Then $N$ is an $\mathbb{H}^2\times\mathbb{E}^1$-manifold.
The mapping torus of $\Theta$ is homeomorphic to
an $\mathbb{H}^3\times\mathbb{E}^1$-manifold which is
not a mapping torus of any self-homeomorphism of an $\mathbb{H}^3$-manifold.
In this case $[\Theta_*]$ has infinite order.
However if $N$ is a $\widetilde{\mathbb{SL}}$-manifold and $[\phi_* ]$ has infinite order
then $M(\phi)$ admits no geometric structure, 
for then $\sqrt\pi\cong\mathbb{Z}$ but
is not a direct factor of any subgroup of finite index.     

If $\zeta\nu\cong\mathbb{Z}$ and $\zeta(\nu/\zeta\nu)=1$ 
then $Hom(\nu/\nu',\zeta\nu)$ embeds in $Out(\nu)$, 
and thus $\nu$ has outer automorphisms of infinite order, 
in most cases \cite{[CR77]}.

Let $N$ be an aspherical closed $\mathbb{X}^3$-manifold where 
$\mathbb{X}^3=\mathbb{H}^3$, $\widetilde{\mathbb{SL}}$ or $\mathbb{H}^2\times\mathbb{E}^1$, 
and suppose that $\beta_1 (N)>0$ but $N$ is not a mapping torus.        
Choose an epimorphism $\lambda:\pi_1 (N)\to\mathbb{Z}$ 
and let $\widehat N$ be the 2-fold covering space associated to 
the subgroup $\lambda^{-1} (2\mathbb{Z})$.
If $c:\widehat N\to \widehat N$ is the covering involution then 
$\mu(n,z)=(c(n),\bar z)$ defines a free
involution on $N\times S^1$, and the orbit space $M$ is an $\mathbb{X}^3\times\mathbb{E}^1$-manifold with
$\beta_1 (M)>0$ which is not a mapping torus. 

\section{The semisimple geometries: $\mathbb{H}^2\times\mathbb{H}^2$,
$\mathbb{H}^4$ and $\mathbb{H}^2(\mathbb{C})$}

In this section we shall consider the remaining three geometries realizable by
closed 4-manifolds. (Not much is known about $\mathbb{H}^4$ or 
$\mathbb{H}^2(\mathbb{C})$.)

Let $P=PSL(2,\mathbb{R})$ be the group of orientation preserving isometries 
of $\mathbb{H}^2$. 
Then $Isom(\mathbb{H}^2\times\mathbb{H}^2)$ contains $P\times P$ as a normal 
subgroup of index 8.
If $M$ is a closed $\mathbb{H}^2\times\mathbb{H}^2 $-manifold then 
$\sigma(M)=0$ and $\chi(M)>0$, 
and $M$ is a complex surface if (and only if) $\pi_1(M)$ is a subgroup of 
$P\times P$.
It is {\it reducible} if it has a finite cover isometric to a product 
of closed surfaces. 
The fundamental groups of such manifolds may be characterized as follows.

\begin{theorem} 
A group $\pi$ is the fundamental group of
a reducible $\mathbb{H}^2\times\mathbb{H}^2$-manifold if and only 
if it is torsion-free, $\sqrt\pi=1$ and $\pi$ has a subgroup of finite index 
which is isomorphic to a product of $PD_2$-groups.
\end{theorem}

\begin{proof} 
The conditions are clearly necessary.
Suppose that they hold.
Then $\pi$ is a $PD_4$-group and has a normal subgroup of finite index 
which is a direct product $K.L\cong{K}\times{L}$, 
where $K$ and $L$ are $PD_2$-groups and $\nu=N_\pi(K)=N_\pi(L)$
has index at most 2 in $\pi$, by Corollary 5.5.2.
After enlarging $K$ and $L$, if necessary,
we may assume that $L=C_\pi(K)$ and $K=C_\pi(L)$.
Hence $\nu/K$ and $\nu/L$ have no nontrivial finite normal subgroup.
(For if $K_1$ is normal in $\nu$ and contains $K$ as a subgroup of finite index
then $K_1\cap L$ is finite, hence trivial, and so $K_1\leq C_\pi(L)$.)
The action of $\nu/L$ by conjugation on $K$ has finite image in $Out(K)$,
and so $\nu/L$ embeds as a discrete cocompact subgroup of $Isom(\mathbb{H}^2)$, 
by the Nielsen conjecture \cite{[Ke83]}. 
Together with a similar embedding for $\nu/K$ we obtain a homomorphism 
from $\nu$ to a discrete cocompact subgroup of $Isom(\mathbb{H}^2\times\mathbb{H}^2)$.

If $[\pi:\nu]=2$ let $t$ be an element of $\pi\setminus\nu$, and let 
$j:\nu/K\to Isom(\mathbb{H}^2)$ be an embedding onto a discrete cocompact 
subgroup $S$. 
Then $tKt^{-1}=L$ and conjugation by $t$ induces an isomorphism 
$f:\nu/K\to\nu/L$. 
The homomorphisms $j$ and $j\circ f^{-1}$ determine an embedding 
$J:\nu\to Isom(\mathbb{H}^2\times\mathbb{H}^2)$ onto 
a discrete cocompact subgroup of finite index in $S\times S$. 
Now $t^2\in\nu$ and $J(t^2)=(s,s)$, where $s=j(t^2K)$.
We may extend $J$ to an embedding of $\pi$ in $Isom(\mathbb{H}^2\times\mathbb{H}^2)$ 
by defining $J(t)$ to be the isometry sending $(x,y)$ to $(y,s.x)$.                       
Thus (in either case) $\pi$ acts isometrically and properly discontinuously 
on $H^2\times H^2$.
Since $\pi$ is torsion-free the action is free, and so $\pi=\pi_1(M)$,
where $M=\pi\backslash (H^2\times H^2)$. 
\end{proof}


\begin{cor}
Let $M$ be a $\mathbb{H}^2\times\mathbb{H}^2$-manifold.
Then $M$ is reducible if and only if it has a $2$-fold covering space
which is homotopy equivalent to the total space of an orbifold bundle over 
a hyperbolic $2$-orbifold.
\end{cor}

\begin{proof}
That reducible manifolds have such coverings was proven in the theorem.
Conversely, an irreducible lattice in $P\times P$ cannot have any nontrivial 
normal subgroups of infinite index, by Theorem IX.6.14 of \cite{[Ma]}. 
Hence an $\mathbb{H}^2\times\mathbb{H}^2$-manifold which is finitely covered by 
the total space of a surface bundle is virtually a cartesian product.
\end{proof}

Is the 2-fold covering space itself such a bundle space over a 2-orbifold?
In general we cannot assume that $M$ is itself fibred over a 2-orbifold.
Let $G$ be a $PD_2$-group with $\zeta{G}=1$ 
and let $\lambda:G\to\mathbb{Z}$ be an epimorphism.
Choose $x\in\lambda^{-1}(1)$.
Then $y=x^2$ is in $K=\lambda^{-1}(2\mathbb{Z})$, 
but is not a square in $K$,
and so 
\[\pi=\langle K\times{K},t\mid t(k,l)t^{-1}=(xlx^{-1},k)\,~\mathrm{for~all}\,
~(k,l)\in{K\times{K}}, t^4=(y,y))\rangle\]
is torsion-free.
A cocompact free action of $G$ on $H^2$ determines 
a cocompact free action of 
$\pi$ on $H^2\times{H}^2$ by $(k,l).(h_1,h_2)=(k.h_1,l.h_2)$ and
$t(h_1,h_2)=(x.h_2,h_1)$, for all $(k,l)\in{K\times{K}}$ and
$(h_1,h_2)\in{H^2\times{H^2}}$.
The group $\pi$ has no normal subgroup which is a $PD_2$-group.
(Note also that if $K$ is orientable $\pi\backslash(H^2\times{H^2})$
is a compact complex surface.)

We may use Theorem 9.9 to give several characterizations of the homotopy types 
of such manifolds.

\begin{theorem} 
Let $M$ be a closed $4$-manifold with fundamental 
group $\pi$.
Then the following are equivalent:
\begin{enumerate}
\item $M$ is homotopy equivalent to a reducible 
$\mathbb{H}^2\times\mathbb{H}^2$-manifold;

\item $\pi$ has an ascendant subgroup $G$ which is $FP_2$, has one end and 
such that $C_\pi(G)$ is not a free group, $\pi_2(M)=0$ and $\chi(M)\not=0$;

\item $\pi$ has a subgroup $\rho$ of finite index which is isomorphic 
to a product of two $PD_2$-groups 
and $\chi(M)[\pi:\rho]=\chi(\rho)\not=0$.

\item $\pi$ is virtually a $PD_4 $-group, $\sqrt\pi=1$ and $\pi$ has a 
torsion-free subgroup of 
finite index which is isomorphic to a nontrivial product $\sigma\times\tau$ 
where 
$\chi(M)[\pi:\sigma\times\tau]=(2-\beta_1(\sigma))(2-\beta_1 (\tau))$.
\end{enumerate}              
\end{theorem}

\begin{proof} 
As $\mathbb{H}^2\times\mathbb{H}^2$-manifolds are aspherical (1) implies (2), 
by Theorem 9.8.
                              
Suppose now that (2) holds. Then $\pi$ has one end,
by transfinite induction, as in Theorem 4.8.
Hence $M$ is aspherical and $\pi$ is a $PD_4 $-group, since $\pi_2 (M)=0$. 
Since $\chi(M)\not=0$ we must have $\sqrt\pi=1$. 
(For otherwise $\beta_i^{(2)}(\pi)=0$ for all $i$, by Theorem 2.3,
and so $\chi(M)=0$.)
In particular, every ascendant subgroup of $\pi$ has trivial centre.
Therefore $G\cap C_\pi(G)=\zeta G=1$ and so 
$G\times C_\pi(G)\cong \rho=G.C_\pi (G)\leq \pi$.
Hence $c.d.C_\pi(G)\leq 2$. 
Since $C_\pi(G)$ is not free $c.d.G\times C_\pi(G)=4$ and so 
$\rho$ has finite index in $\pi$. 
(In particular, $[C_\pi(C_\pi(G)):G]$ is finite.)
Hence $\rho$ is a $PD_4$-group and $G$ and $C_\pi(G)$ are $PD_2$-groups,
so $\pi$ is virtually a product.
Thus (2) implies (1), by Theorem 9.9. 

It is clear that (1) implies (3).
If (3) holds then on applying Theorems 2.2 and 3.5 to the 
finite covering
space associated to $\rho$ we see that $M$ is aspherical, so $\pi$ is
a $PD_4$-group and (4) holds.
Similarly, $M$ is asperical if (4) holds.
In particular, $\pi$ is a $PD_4$-group and so is torsion-free. 
Since $\sqrt\pi=1$ neither $\sigma$ nor $\tau$ can be infinite cyclic, 
and so they are each $PD_2 $-groups. Therefore $\pi$ is the fundamental group
of a reducible $\mathbb{H}^2\times\mathbb{H}^2$-manifold, by
Theorem 9.9, and $M\simeq\pi\backslash {H}^2\times {H}^2$,
by asphericity.
\end{proof}

The asphericity of $M$ could be ensured by assuming that $\pi$ be $PD_4$ and 
$\chi(M)=\chi(\pi)$, instead of assuming that $\pi_2(M)=0$.

For $\mathbb{H}^2\times\mathbb{H}^2$-manifolds we can give more 
precise criteria for reducibility.

\begin{theorem} 
Let $M$ be a closed $\mathbb{H}^2\times\mathbb{H}^2$-manifold 
with fundamental group $\pi$. Then the following are equivalent:
\begin{enumerate}
\item $\pi$ has a subgroup of finite index which is a nontrivial direct product;

\item $\mathbb{Z}^2<\pi$;

\item $\pi$ has a nontrivial element with nonabelian centralizer;

\item $\pi\cap(\{1\}\times P)\not=1$;

\item $\pi\cap(P\times \{1\})\not=1$; 

\item $M$ is reducible.
\end{enumerate}
\end{theorem}

\begin{proof} Since $\pi$ is torsion-free each of the above conditions
is invariant under passage to subgroups of finite index, and so we may assume 
without loss of generality that $\pi\leq P\times P$.
Suppose that $\sigma$ is a subgroup of finite index in $\pi$ which is a nontrivial direct product. 
Since $\chi(\sigma)\not=0$ neither factor can be infinite cyclic, and so the factors must be $PD_2$-groups. 
In particular, $\mathbb{Z}^2<\sigma$ and the centraliser of any element of either direct factor is nonabelian.
Thus (1) implies (2) and (3).

Suppose that $(a,b)$ and $(a',b')$ generate a subgroup of $\pi$ 
isomorphic to $\mathbb{Z}^2$.
Since centralizers of elements of infinite order in $P$ are cyclic
the subgroup of $P$ generated by $\{a,a'\}$ is infinite cyclic or is finite.
We may assume without loss of generality that $a'=1$, 
and so (2) implies (4).
Similarly, (2) implies (5).

Let $g=(g_1,g_2)\in P\times P$ be nontrivial.
Since centralizers of elements of infinite order in $P$ 
are infinite cyclic and
$C_{P\times P}(\langle g\rangle)=
C_P(\langle g_1\rangle)\times C_P(\langle g_2\rangle)$ 
it follows that 
if $C_\pi(\langle g\rangle)$ is nonabelian 
then either $g_1$ or $g_2$ has finite order. 
Thus (3) implies (4) and (5).

Let $K_1=\pi\cap(\{1\}\times P)$ and $K_2=\pi\cap(P\times\{1\})$.
Then $K_i$ is normal in $\pi$, and there are exact sequences
\begin{equation*}
1\to K_i\to\pi\to L_i\to 1,
\end{equation*}
where $L_i=pr_i (\pi)$ is the image of $\pi$ under 
projection to the $i^{th}$ factor of $P\times P$, for $i=1$ and $2$.
Moreover $K_i$ is normalised by $L_{3-i}$, for $i=1$ and 2.
Suppose that $K_1\not=1$.
Then $K_1$ is nonabelian, since it is normal in $\pi$ and $\chi(\pi)\not=0$.
If $L_2$ were not discrete then elements of $L_2$ sufficiently close to the
identity would centralize $K_1$. 
As centralizers of nonidentity elements of $P$ are abelian, this would imply that $K_1$ is abelian.
Hence $L_2$ is discrete.
Now $L_2\backslash H^2$ is a quotient of $\pi\backslash H\times H$ and so is compact.
Therefore $L_2$ is virtually a $PD_2$-group.
Now $c.d.K_2+v.c.d.L_2\geq c.d.\pi=4$, so $c.d.K_2\geq 2$.
In particular, $ K_2\not=1$ and so a similar argument now shows that $c.d.K_1\geq 2$.
Hence $c.d.K_1\times K_2\geq 4$. Since $K_1\times K_2\cong K_1.K_2 \leq \pi$
it follows that $\pi$ is virtually a product, and $M$ is finitely covered
by $(K_1\backslash H^2)\times (K_2\backslash H^2)$.  
Thus (4) and (5) are equivalent, and imply (6). 
Clearly (6) implies (1). 
\end{proof}

The idea used in showing that (4) implies (5) and (6) 
derives from one used in the proof of \cite[Theorem 6.3]{[Wl85]}.
Orientable reducible $\mathbb{H}^2\times\mathbb{H}^2$-manifolds 
with isomorphic fundamental group are diffeomorphic \cite{[Ca00]}.

If $\Gamma$ is a discrete cocompact subgroup of $P\times P$
such that $M=\Gamma\backslash H^2\times H^2$ is irreducible 
then $\Gamma\cap P\times\{1\}=\Gamma\cap\{1\}\times P=1$, by the theorem. 
Hence the natural foliations of $H^2\times H^2$ descend to give 
a pair of transverse foliations of $M$ by copies of $H^2$. 
Conversely, if $M$ is a closed Riemannian 4-manifold with a codimension 2 
metric foliation by totally geodesic surfaces then $M$ has a finite cover 
which either admits the geometry $\mathbb{H}^2\times\mathbb{E}^2$ 
or $\mathbb{H}^2\times\mathbb{H}^2 $, or is the total space of an 
$S^2$- or $T$-bundle over a closed surface,
or is the mapping torus of a self homeomorphism of 
$\mathbb{R}^3/\mathbb{Z}^3$, 
$S^2\times S^1$ or a lens space \cite{[Ca90]}.

An irreducible $\mathbb{H}^2\times\mathbb{H}^2 $-lattice 
is an arithmetic subgroup of $Isom(\mathbb{H}^2\times\mathbb{H}^2)$, 
and has no nontrivial normal subgroups of infinite index 
\cite[Theorems IX.6.5 and 14]{[Ma]}.
Such irreducible lattices are rigid, and so the argument of  
\cite[Theorem 8.1]{[Wa72]} implies that there are only finitely many 
irreducible $\mathbb{H}^2\times\mathbb{H}^2$-manifolds 
with given Euler characteristic.
What values of $\chi$ are realized by such manifolds?
If $M$ is a closed orientable $\mathbb{H}^2\times\mathbb{H}^2$-manifold 
then $\sigma(M)=0$, so $\chi(M)$ is even, and $\chi(M)>0$
\cite{[Wl86]}.
There are examples (``fake quadrics") with $\beta_1=0$ and $\chi(M)=4$ 
\cite{[Dz13]}.

Irreducible arithmetic $\mathbb{H}^2\times\mathbb{H}^2$-lattices are
commensurable with lattices constructed as follows.
Let $F$ be a totally real number field, with ring of integers $O_F$.
Let $H$ be a skew field which is a quaternion algebra over $F$
such that $H\otimes_\sigma\mathbb{R}\cong M_2(\mathbb{R})$ 
for exactly two embeddings $\sigma$ of $F$ in $\mathbb{R}$.
If $A$ is an order in $H$ (a subring which is also 
a finitely generated $O_F$-submodule and such that $F.A=H$) 
then the quotient of the group of units $A^\times$ by $\pm1$ 
embeds as a cocompact irreducible 
$\mathbb{H}^2\times\mathbb{H}^2$-lattice $\Gamma(A)$. 
(It is difficult to pin down a reference for the claim
that this construction realizes all commensurability classes,
but it appears to be ``well-known to the experts".
See \cite{[Sh63],[Bo81]}.)

Much less is known about closed $\mathbb{H}^4$- or $\mathbb{H}^2(\mathbb{C})$-manifolds.
There are only finitely many such
manifolds with a given Euler characteristic. 
(See \cite[Theorem 8.1]{[Wa72]}.)
If $M$ is a closed orientable $\mathbb{H}^4$-manifold 
then $\sigma(M)=0$, so $\chi(M)$ is even,
and $\chi(M)>0$ \cite{[Ko92]}.
The examples of \cite{[CM05]} and \cite{[Da85]} have $\beta_1>0$,
and so covers of these realize all positive multiples of 16 and of 26.
No closed $\mathbb{H}^4$-manifold admits a complex structure.
If $M$ is a closed $\mathbb{H}^2(\mathbb{C})$-manifold 
it is orientable and $\chi(M)=3\sigma(M)>0$ \cite{[Wl86]}. 
The isometry group of $\mathbb{H}^2(\mathbb{C})$ has two components; 
the identity component is $SU(2,1)$ 
and acts via holomorphic isomorphisms on the unit ball 
\begin{equation*}
\{(w,z)\in C^2 :|w|^2+|z|^2<1\}.
\end{equation*}
There are $\mathbb{H}^2(\mathbb{C})$-manifolds with $\beta_1=2$ and $\chi=3$
\cite{[CS10]},
and so all positive multiples of 3 are realized.
Since $H^4$ and $H^2(\mathbb{C})$ are rank 1 symmetric spaces 
the fundamental groups can contain no noncyclic abelian subgroups 
\cite{[Pr43]}.
In each case there are cocompact lattices which are not arithmetic.
At present there are not even conjectural intrinsic characterizations of 
such groups.
(See also \cite{[Rt]} for the geometries $\mathbb{H}^n$ and 
\cite{[Go]} for the geometries $\mathbb{H}^n(\mathbb{C})$.)

Each of the geometries $\mathbb{H}^2\times\mathbb{H}^2$, 
$\mathbb{H}^4$ and $\mathbb{H}^2(\mathbb{C})$ 
admits cocompact lattices which are not coherent.
(See \S4 of Chapter 4 above, \cite{[BM94]} and \cite{[Ka13]}, respectively.)
Is this true of every such lattice for one of these geometries?
(Lattices for the other geometries are coherent.)

\section{Miscellany}

A homotopy equivalence between two closed $\mathbb{H}^n$- or 
$\mathbb{H}^n(\mathbb{C})$-manifolds of dimension $\geq3$ 
is homotopic to an isometry, 
by {\it Mostow rigidity} \cite{[Ms68]}.
Farrell and Jones have established ``topological" analogues of Mostow rigidity, 
for manifolds with a metric of nonpositive sectional curvature 
and dimension $\geq5$.
By taking cartesian products with $S^1$, 
we can use their work in dimension 4 also.

\begin{theorem} 
Let $\mathbb{X}^4$ be a geometry of aspherical type.
A closed 4-manifold $M$ with fundamental group $\pi$
is $s$-cobordant to an $\mathbb{X}^4$-manifold if and only if 
$\pi$ is isomorphic to a cocompact lattice in
$Isom(\mathbb{X}^4)$ and $\chi(M)=\chi(\pi)$.
\end{theorem}

\begin{proof} 
The conditions are clearly necessary.
If they hold $c_M:M\to\pi\backslash X$ is a homotopy equivalence, 
by Theorem 3.5.
If $\mathbb{X}^4$ is of solvable type $c_M$ is homotopic to a homeomorphism,
by Theorem 8.1.
In most of the remaining cases (excepting only
$\widetilde{\mathbb{SL}}\times\mathbb{E}^1$ -- see \cite{[Eb82]})
the geometry has nonpositive sectional curvatures,
so $Wh(\pi)=Wh(\pi\times\mathbb{Z})=0$ and $M\times S^1$ 
is homeomorphic to $(\pi\backslash X)\times S^1$ \cite{[FJ93']}.
Hence $M$ and $\pi\backslash X$ are $s$-cobordant, by Lemma 6.10.
The case $\mathbb{X}^4=\widetilde{\mathbb{SL}}\times\mathbb{E}^1$ 
follows from \cite{[NS85]} if 
$\pi\leq Isom_o(\widetilde{\mathbb{SL}}\times\mathbb{E}^1)$,
so that $\pi\backslash(\widetilde{SL}\times{R})$ 
admits an effective $T$-action,
and from \cite{[HR11]} in general.
\end{proof}

If an aspherical closed 4-manifold $M$ has a geometric decomposition
then $\pi=\pi_1(M)$ is built from the fundamental groups of the pieces by
amalgamation along torsion-free virtually poly-$Z$ subgroups.
As the Whitehead groups of the geometric pieces are trivial
(by the argument of \cite{[FJ86]}) and the amalgamated subgroups are 
regular noetherian it follows from the $K$-theoretic
Mayer-Vietoris sequence of Waldhausen that $Wh(\pi)=0$.
Is $M$ $s$-rigid?
This is so if all the pieces are $\mathbb{H}^4$- or
$\mathbb{H}^2(\mathbb{C})$-manifolds or irreducible
$\mathbb{H}^2\times\mathbb{H}^2$-manifolds,
for then the inclusions of the cuspidal subgroups into 
the fundamental groups of the pieces are square-root closed.
(See \cite{[BJS17]}.)

For the semisimple geometries we may avoid the appeal to $L^2$-methods 
to establish asphericity as follows.
Since $\chi(M)>0$ and $\pi$ is infinite and residually finite there is a 
subgroup 
$\sigma$ of finite index such that the associated covering spaces $M_\sigma $ 
and $\sigma\backslash X$ are orientable and $\chi(M_\sigma)=\chi(\sigma)>2$. 
In particular, $H^2 (M_\sigma;\mathbb{Z})$ has elements of infinite order.
Since the classifying map $c_{M_\sigma} :M_\sigma\to\sigma\backslash X$ is 
2-connected it induces an isomorphism on $H^2$ and hence is a degree-1 map, 
by Poincar\'e duality. 
Therefore it is a homotopy equivalence, by Theorem 3.2.

\begin{theorem} 
An aspherical closed 4-manifold $M$ which is finitely covered by
a geometric manifold is homotopy equivalent to a geometric 4-manifold.
\end{theorem}

\begin{proof} 
The result is clear for infrasolvmanifolds, and follows from 
Theorems 9.5 and 9.6 if the geometry is
$\mathbb{H}^2\times\mathbb{E}^2$ or
$\widetilde{\mathbb{SL}}\times\mathbb{E}^1$,
and from Theorem 9.9 if $M$ is finitely covered by 
a reducible $\mathbb{H}^2\times\mathbb{H}^2$-manifold.
It holds for the other closed $\mathbb{H}^2\times\mathbb{H}^2$-manifolds 
and for the geometries $\mathbb{H}^4$ and 
$\mathbb{H}^2(\mathbb{C})$ by Mostow rigidity. 
If the geometry is $\mathbb{H}^3\times\mathbb{E}^1$ 
then $\sqrt\pi\cong\mathbb{Z}$ and 
$\pi/\sqrt\pi$ is virtually the group of a $\mathbb{H}^3$-manifold.
Hence $\pi/\sqrt\pi$ acts isometrically and properly discontinuously
on $\mathbb{H}^3$, by Mostow rigidity.
Moreover as the hypotheses of Lemma 9.4 are satisfied, by Theorem 9.3,
there is a homomorphism
$\lambda:\pi\to D<Isom(\mathbb{E}^1)$ which maps $\sqrt\pi$ injectively.
Together these actions determine a discrete and cocompact action of $\pi$
by isometries on $H^3\times\mathbb{R}$.
Since $\pi$ is torsion-free this action is free, and so
$M$ is homotopy equivalent to an $\mathbb{H}^3\times\mathbb{E}^1$-manifold. 
\end{proof}

The result holds also for $\mathbb{S}^4$ and $\mathbb{CP}^2$,
but is not yet clear for 
$\mathbb{S}^2\times\mathbb{E}^2$ or $\mathbb{S}^2\times\mathbb{H}^2$.
It fails for $\mathbb{S}^3\times\mathbb{E}^1$
or $\mathbb{S}^2\times\mathbb{S}^2$. 
In particular, there is a closed nonorientable 4-manifold which is doubly 
covered by $S^2\times S^2$ but is not homotopy equivalent 
to an $\mathbb{S}^2\times\mathbb{S}^2$-manifold. (See Chapters 11 and 12.)

If $\pi$ is the fundamental group of an aspherical closed geometric 4-manifold 
then $\beta_s^{(2)}(\pi)=0$ for $s=0$ or 1,
and so $\beta_2^{(2)}(\pi)=\chi(\pi)$, by Theorem 1.35 of \cite{[Lu]}.
Hence $\mathrm{def}(\pi)\leq{\min\{0, -\chi(\pi)\}}$, 
by Theorems 2.4 and 2.5 (and since $c.d.\pi>2$).
If $\pi$ is orientable this gives
$\mathrm{def}(\pi)\leq2\beta_1(\pi)-\beta_2(\pi)-2$.
When $\beta_1(\pi)\leq1$ this is an improvement on the estimate 
$\mathrm{def}(\pi)\leq\beta_1(\pi)-\beta_2(\pi)$ 
derived from the ordinary homology of
a 2-complex with fundamental group $\pi$.
(In particular, the fundamental groups of  
``fake complex projective planes"
-- compact complex surfaces with the rational homology of $\mathbb{CP}^2$,
but with geometry $\mathbb{H}^2(\mathbb{C})$ --
have deficiency $\leq-3$.)

%% file: m5-10.tex
\chapter{Manifolds covered by $S^2 \times\mathbb{R}^2$}

If the universal covering space of a closed 4-manifold with infinite 
fundamental group is homotopy equivalent to a finite complex 
then it is either contractible or homotopy equivalent to $S^2$ or $S^3$,
by Theorem 3.9.
The aspherical cases have been considered in Chapters 8 and 9.
In this chapter and the next we shall consider the spherical cases.
We show first that if $\widetilde M\simeq S^2 $ then $M$ has a finite covering
space which is $s$-cobordant to a product $S^2\times B$, 
where $B$ is aspherical, 
and $\pi$ is the group of a $\mathbb{S}^2\times{\mathbb E}^2$- 
or $\mathbb{S}^2\times{\mathbb H}^2$-manifold. 
In \S2 and \S3 we show that there are at most two homotopy types of
such manifolds for each such group $\pi$ and action $u:\pi\to{Z/2Z}$. 
In \S4 we show that all $S^2$- and $RP^2$-bundles over aspherical closed 
surfaces are geometric.
We determine the nine possible elementary amenable groups
(corresponding to the geometry $\mathbb{S}^2\times{\mathbb E}^2$) in \S5.
Six have infinite abelianization, 
and in \S6 we use Stiefel-Whitney classes to
distinguish the related homotopy types.  
After some remarks on the homeomorphism classification, 
we show finally that every 4-manifold with fundamental group a
$PD_2$-group admits a 2-connected degree-1 map 
to the total space of an $S^2$-bundle.
For brevity, we shall let $\mathbb{X}^2$ denote both $\mathbb{E}^2$ and $\mathbb{H}^2$.

\section{Fundamental groups}

The determination of the closed 4-manifolds with universal covering space 
homotopy equivalent to $S^2$ rests on Bowditch's Theorem, via Theorem 5.14.

\begin{theorem} 
Let $M$ be a closed $4$-manifold with fundamental 
group $\pi$. Then the following conditions are equivalent:
\begin{enumerate}
\item $\pi$ is virtually a $PD_2 $-group and $\chi(M)=2\chi(\pi)$;

\item  $\pi\not=1$ and $\pi_2 (M)\cong\mathbb{Z}$;

\item $M$ has a covering space of degree dividing $4$ which is
$s$-cobordant to $S^2\times B$, 
where $B$ is an aspherical closed orientable surface;

\item $M$ is virtually $s$-cobordant to an 
$\mathbb{S}^2\times \mathbb{X}^2$-manifold.

\end{enumerate}
If these conditions hold then $\widetilde M$ 
is homeomorphic to $S^2\times\mathbb{R}^2 $,
and the kernel of the action
$u:\pi\to{Aut(\pi_2(M))}=Z/2Z$ is a $PD_2$-group.
\end{theorem}

\begin{proof} If (1) holds then $\pi_2(M)\cong\mathbb{Z}$,
by Theorem 5.10, and so (2) holds.
If (2) holds then the covering space associated to
the kernel of the natural action of $\pi$ on $\pi_2(M)$ is
homotopy equivalent to the total space of an $S^2$-bundle $\xi$
over an aspherical closed surface with $w_1(\xi)=0$, by 
Lemma 5.11 and Theorem 5.14.
On passing to a 2-fold covering space,
if necessary, we may assume that $w_2 (\xi)=w_1 (M)=0$ also. 
Hence $\xi$ is trivial and so the corresponding covering space of $M$ 
is $s$-cobordant to a product $S^2\times B$ with $B$ orientable.
It is clear that $(3)\Rightarrow(4)$ and $(4)\Rightarrow(1)$.

The final observations follow from Theorems 6.16 and 5.14, respectively.
\end{proof}

This follows also from \cite{[Fa74]} instead of \cite{[Bo04]}, 
if we know also that $\chi(M)\leq0$.

\begin{theorem} 
Let $u:\pi\to{Z/2Z}$ be a homomorphism such that 
$\kappa=\mathrm{Ker}(u)$ is a $PD_2$-group.
Then the pair $(\pi,u)$ is realized by a closed 
$\mathbb{S}^2\times\mathbb{E}^2 $-manifold, 
if $\pi$ is virtually $\mathbb{Z}^2$, 
and by a closed $\mathbb{S}^2\times\mathbb{H}^2$-manifold otherwise.
\end{theorem}

\begin{proof} 
Let $X$ be $\mathbb{R}^2$, if $\pi$ is virtually $\mathbb{Z}^2$,
or the hyperbolic plane, otherwise. 
If $\pi$ is torsion-free then it is itself a surface group.
If $\pi$ has a nontrivial finite normal subgroup 
then it is a direct product $(Z/2Z)\times\kappa$.
In either case $\pi$ is the fundamental group of 
a corresponding product of surfaces.
Otherwise $\pi$ is a semidirect product $\kappa\rtimes(Z/2Z)$ 
and is a plane motion group,
by a theorem of Nielsen \cite{[Zi]}.
Thus there is a monomorphism $f:\pi\to{Isom(\mathbb{X}^2)}$ 
with image a discrete subgroup which acts cocompactly on $X$,
with quotient $B=\pi\backslash{X}$ an $\mathbb{X}^2$-orbifold.
The homomorphism
\[(u,f):\pi\to\{\pm I\}\times Isom(\mathbb{X}^2)\leq 
Isom(\mathbb{S}^2\times\mathbb{X}^2)\] 
is then a monomorphism onto a discrete subgroup which acts freely and 
cocompactly on $S^2\times{X}$.
In all cases the pair $(\pi,u)$ may be realised geometrically. 
\end{proof}

The manifold $M$ constructed in Theorem 10.2 is a cartesian product with $S^2$ 
if $u$ is trivial and fibres over $RP^2$ otherwise.
If $\pi\not\cong(Z/2Z)\times\kappa$ projection to $X$ induces
an orbifold bundle projection from $M$ to $B$ with general fibre $S^2$.
 
\section{The first $k$-invariant}

The main result of this section is that if $\pi=\pi_1(M)$ is not a product 
then $k_1(M)=\beta^u(U^2)$, where 
$U\in{H^1(\pi;\mathbb{F}_2)}=Hom(\pi,Z/2Z)$ corresponds 
to the action $u:\pi\to Aut(\pi_2(M))$ and $\beta^u$ is a
``twisted Bockstein" described below.

We shall first show that the orientation character 
$w$ and the action $u$ of $\pi$ on $\pi_2$ determine each other.

\begin{lemma} 
Let $M$ be a $PD_4$-complex with fundamental group $\pi\not=1$ 
and such that $\pi_2 (M)\cong\mathbb{Z}$.
Then $H^2(\pi;\mathbb{Z}[\pi])\cong\mathbb{Z}$ and $w_1(M)=u+v$, 
where $u:\pi\to Aut(\pi_2(M))=Z/2Z$ and 
$v:\pi\to Aut(H^2(\pi;\mathbb{Z}[\pi]))=Z/2Z$ are the natural actions.
\end{lemma} 

\begin{proof} Since $\pi\not=1$ it is infinite, by Theorem 10.1.
Thus $Hom_{\mathbb{Z}[\pi]}(\pi_2(M),\mathbb{Z}[\pi])$ $=0$ and so 
Poincar\'e duality determines an isomorphism
$D:\overline{H^2(\pi;\mathbb{Z}[\pi])}\cong\pi_2(M)$, by Lemma 3.3.
Let $w=w_1(M)$. Then $gD(c)=(-1)^{w(g)}cg^{-1}=(-1)^{v(g)+w(g)}c$ for all
$c\in H^2(\pi;\mathbb{Z}[\pi])$ and $g\in\pi$, and so $u=v+w$.
\end{proof}

Note that $u$ and $w_1(M)$ are constrained by the further conditions that 
$\kappa=\mathrm{Ker}(u)$ is torsion-free and $\mathrm{Ker}(w_1(M))$ 
has infinite abelianization if $\chi(M)\leq 0$.
If $\pi< Isom(\mathbb{X}^2)$ is a plane motion group then $v(g)$ detects
whether $g\in\pi$ preserves the orientation of $X^2$.
If $\pi\cong(Z/2Z)\times\kappa$ then $v|_{Z/2Z}=0$ and $v|_\kappa=w_1(\kappa)$.
If $\pi$ is torsion-free then $M$ is homotopy equivalent to the total space of
an $S^2$-bundle $\xi$ over an aspherical closed surface $B$,
and the equation $u=w_1(M)+v$ follows from Lemma 5.11.
        
Let $\beta^u$ be the Bockstein operator associated with the 
coefficient sequence 
\[
0\to\mathbb{Z}^u\to\mathbb{Z}^u\to\mathbb{F}_2\to 0, 
\]
and let $\overline{\beta^u}$ be the composition 
with reduction {\it mod} (2).  
In general $\beta^u$ is NOT the Bockstein operator for the untwisted sequence 
$0\to\mathbb{Z}\to\mathbb{Z}\to\mathbb{F}_2\to 0$, 
and $\overline{\beta^u}$ is not $Sq^1$,
as can be seen already for cohomology of the group $Z/2Z$ acting nontrivially 
on $\mathbb{Z}$, as $\beta^u(H^1(Z/2Z:\mathbb{F}_2))=0$ if $u$ is nontrivial.

\begin{lemma} 
Let $M$ be a $PD_4$-complex with fundamental group $\pi$
and such that $\pi_2 (M)\cong\mathbb{Z}$.
If $\pi$ has nontrivial torsion
$H^s(M;\mathbb{F}_2)\cong H^s(\pi;\mathbb{F}_2)$ for $s\leq2$.
The Bockstein operator 
$\beta^u:H^2(\pi;\mathbb{F}_2)\to{H^3(\pi;\mathbb{Z}^u)}$ 
is onto, and reduction {\it mod} $(2)$ from $H^3(\pi;\mathbb{Z}^u)$ to 
$H^3(\pi;\mathbb{F}_2)$ is a monomorphism.
The restriction of $k_1(M)$ to each subgroup of order $2$ is nontrivial.
Its image in  $H^3(M;\mathbb{Z}^u)$ is $0$.
\end{lemma} 
           
\begin{proof}
These assertions hold vacuously if $\pi$ is torsion-free,
so we may assume that $\pi$ has an element of order 2.
Then $M$ has a covering space $\widehat{M}$ homotopy equivalent to $RP^2$,
and so the mod-2 Hurewicz homomorphism from $\pi_2(M)$ to
$H_2(M;\mathbb{F}_2)$ is trivial, 
since it factors through $H_2(\widehat{M};\mathbb{F}_2)$. 
Since we may construct $K(\pi,1)$ from $M$ by adjoining cells to kill the
higher homotopy of $M$ the first assertion follows easily.

The group $H^3(\pi;\mathbb{Z}^u)$ has exponent dividing 2, 
since the composition of restriction to $H^3(\kappa;\mathbb{Z})=0$ 
with the corestriction back to $H^3(\pi;\mathbb{Z}^u)$ 
is multiplication by the index $[\pi:\kappa]$. 
Consideration of the long exact sequence
associated to the coefficient sequence shows that $\beta^u$ is onto.
If $f:Z/2Z\to\pi$ is a monomorphism then $f^*k_1(M)$ is 
the first $k$-invariant of $\widetilde M/f(Z/2Z)\simeq RP^2$, 
which generates $H^3(Z/2Z;\pi_2(M))=Z/2Z$. 
The final assertion is clear.
\end{proof}

\begin{lemma}
Let $\alpha=*^kZ/2Z=
\langle{x_i, 1\leq{i}\leq{k}}\mid{x_i^2=1~\forall~i}\rangle$ 
and let $u(x_i)=-1$ for all $i$.
Then restriction from $\alpha$ to $\phi=\mathrm{Ker}(u)$ 
induces an epimorphism from $H^1(\alpha;\mathbb{Z}^u)$ to $H^1(\phi;\mathbb{Z})$.
\end{lemma}

\begin{proof}
Let $x=x_1$ and $y_i=x_1x_i$ for all $i>1$.
Then $\phi=\mathrm{Ker}(u)$ is free with basis $\{y_2,\dots,y_k\}$
and so $\alpha\cong{F(k-1)\rtimes{Z/2Z}}$.
If $k=2$ then $\alpha$ is the infinite dihedral group $D$ and
the lemma follows by direct calculation with resolutions.
In general, the subgroup $D_i$ generated by $x$ and $y_i$ is an
infinite dihedral group, and is a retract of $\alpha$.
The retraction is compatible with $u$,
and so restriction maps $H^1(\alpha;\mathbb{Z}^u)$ onto $H^1(D_i;\mathbb{Z}^u)$.
Hence restriction maps $H^1(\alpha;\mathbb{Z}^u)$ onto each summand
$H^1(\langle{y_i}\rangle;\mathbb{Z})$ of $H^1(\phi;\mathbb{Z})$, 
and the result follows.
\end{proof}

In particular, if $k$ is even then $z=\Pi{x_i}$ generates 
a free factor of $\phi$, and restriction 
maps $H^1(\alpha;\mathbb{Z}^u)$ onto $H^1(\langle{z}\rangle;\mathbb{Z})$.

\begin{theorem}
Let $B$ be an aspherical $2$-orbifold with non-empty singular locus, 
and let $u:\pi=\pi_1^{orb}(B)\to{Z/2Z}$ be an
epimorphism with torsion-free kernel $\kappa$. 
Suppose that $B$ has $r$ reflector curves and $k$ cone points.
Then $H^2(\pi;\mathbb{Z}^u)\cong(Z/2Z)^r$ if $k>0$ and
$H^2(\pi;\mathbb{Z}^u)\cong\mathbb{Z}\oplus(Z/2Z)^{r-1}$ if $k=0$.
In all cases $\beta^u(U^2)$ is the unique element of $H^3(\pi;\mathbb{Z}^u)$    
which restricts non-trivially to each subgroup of order $2$.
\end{theorem}

\begin{proof}
If $g$ has order 2 in $\pi$ then $u(g)=-1$, 
and so $\beta^u(U^2)$ restricts non-trivially 
to the subgroup generated by $g$. 

Suppose first that $B$ has no reflector curves.
Then $B$ is the connected sum of a closed surface $G$
with $S(2_k)$, the sphere with $k$ cone points of order 2.
If $B=S(2_k)$ then $k\geq4$, since $B$ is aspherical.
Hence $\pi\cong\mu*_{\mathbb{Z}}\nu$, 
where $\mu=*^{k-2}Z/2Z$ and $\nu=Z/2Z*Z/2Z$
are generated by cone point involutions.
Otherwise $\pi\cong\mu*_{\mathbb{Z}}\nu$,
where $\mu=*^kZ/2Z$ and $\nu=\pi_1(G\setminus{D^2})$ 
is a non-trivial free group.
Every non-trivial element of finite order in such a generalized 
free product must be conjugate to one of the involutions.
In each case a generator of the amalgamating subgroup is identified with 
the product of the involutions which generate the factors of $\mu$
and which is in $\phi=\mathrm{Ker}(u|_\mu)$.

Restriction from $\mu$ to $\mathbb{Z}$ induces an epimorphism from
$H^1(\mu;\mathbb{Z}^u)$ to $H^1(\mathbb{Z};\mathbb{Z})$, by Lemma 10.5,
and so 
\[H^2(\pi;\mathbb{Z}^u)\cong{H^2(\mu;\mathbb{Z}^u)}\oplus
{H^2(\nu;\mathbb{Z}^u)}=0,\]
by the Mayer-Vietoris sequence with coefficients $\mathbb{Z}^u$.
Similarly,
\[H^2(\pi;\mathbb{F}_2)\cong
{H^2(\mu;\mathbb{F}_2)}\oplus{H^2(\nu;\mathbb{F}_2)},\]
by the Mayer-Vietoris sequence with coefficients $\mathbb{F}_2$.
Let $e_i\in{H^2(\pi;\mathbb{F}_2)}$
$=Hom(H_2)(\pi);\mathbb{F}_2),\mathbb{F}_2)$ 
correspond to restriction to the $i^{th}$ cone point.
Then $\{e_1,\dots,e_{2g+2}\}$ forms a basis for 
$H^2(\pi;\mathbb{F}_2)\cong\mathbb{F}_2^{2g+2}$,
and $\Sigma{e_i}$ is clearly the only element with nonzero
restriction to all the cone point involutions.
Since $H^2(\pi;\mathbb{Z}^u)=0$ the $u$-twisted Bockstein maps
$H^2(\pi;\mathbb{F}_2)$ isomorphically onto $H^3(\pi;\mathbb{Z}^u)$,
and so $\beta^u(U^2)$ is the unique such element 
of $H^3(\pi;\mathbb{Z}^u)$.

Suppose now that $r>0$.
Then $B=B_o\cup{r}\mathbb{J}$, where $B_o$ is a connected 2-orbifold 
with $r$ boundary components and $k$ cone points,
and $\mathbb{J}=S^1\times[[0,1)$ is a product neighbourhood 
of a reflector circle.
Hence $\pi=\pi\mathcal{G}$,
where $\mathcal{G}$ is a graph of groups with underlying graph
a tree having one vertex of valency $r$ with group 
$\nu=\pi_1^{orb}(B_o)$, $r$ terminal vertices,
with groups $\gamma_i\cong\pi_1^{orb}(\mathbb{J})=\mathbb{Z}\oplus{Z/2Z}$,
and $r$ edge groups $\omega_i\cong\mathbb{Z}$.
If $k>0$ then restriction maps $H^1(\nu;\mathbb{Z}^u)$ 
onto $\oplus{H^1(\omega_i;\mathbb{Z})}$, and then
$H^2(\pi;\mathbb{Z}^u)\cong\oplus{H^2(\gamma_i;\mathbb{Z}^u)}\cong{Z/2Z}^r$.
However if $k=0$ then 
$H^2(\pi;\mathbb{Z}^u)\cong\mathbb{Z}\oplus(Z/2Z)^{r-1}$.
The Mayer-Vietoris sequence now gives 
\[H^2(\pi;\mathbb{F}_2)\cong{H^2(\nu;\mathbb{F}_2)}\oplus
(H^2(\mathbb{Z}\oplus{Z/2Z};\mathbb{F}_2))^r\cong\mathbb{F}_2^{2r+k}.\]
The generator of the second summand of 
$H^2(\mathbb{Z}\oplus{Z/2Z};\mathbb{F}_2)$
is in the image of reduction modulo $(2)$ from
$H^2(\mathbb{Z}\oplus{Z/2Z};\mathbb{Z}^u)$, 
and so is in the kernel of $\beta^u$.
Therefore the image of $\beta^u$ has a basis corresponding to the 
cone points and reflector curves,
and we again find that $\beta^u(U^2)$ is the unique element 
of $H^3(\pi;\mathbb{Z}^u)$ with non-trivial restriction 
to each subgroup of order 2.
\end{proof}

If $\pi$ is torsion-free it is a $PD_2$-group and so $k_1(M)$ and
$\beta^u(U^2)$ are both 0.
If $\pi\cong(Z/2Z)\times\kappa$ then the torsion subgroup is unique,
and $u$ and $w$ each split the inclusion of this subgroup.
In this case $H^3(\pi;\mathbb{Z}^u)\cong(Z/2Z)^2$, 
and there are two classes with nontrivial restriction to the torsion subgroup.
Can it be seen {\it a priori\/} (i.e., without appealing to Theorem 5.16)
that the $k$-invariant must be standard?
 
\section{Homotopy type}
 
Let $M$ be the 4-manifold realizing $(\pi,u)$ in Theorem 10.2.
If $u$ is trivial then $P_2(M)\simeq{CP^\infty}\times{K(\pi,1)}$.
Otherwise, we may construct a model for $P_2(M)$ as follows.
Let $\zeta$ be the involution of $K(\pi,1)$
inducing the action of $\pi/\kappa$ on $\kappa$,
and let $\sigma$ be the involution of $CP^\infty$ given by
$\sigma(\vec{z})=[-\overline{z_1}:\overline{z_0}:\dots:\overline{z_n}]$
for all $\vec{z}=[z_0:z_1:\dots:z_n]$ in $CP^\infty=\varinjlim{CP^n}$.
Then $\sigma$ extends the antipodal map of $S^2=CP^1$.
Let
\[P={CP^\infty\times{S^\infty}\times{K(\kappa,1)}/
(z,s,k)\sim(\sigma(z),-s,\zeta(k))}.\]
The diagonal map of $S^2$ into $CP^1\times{S^2}$ gives rise to a natural
inclusion of $M$ into $P$.
This map is 3-connected, and so may be identified with $f_M$.

Every $PD_4$-complex $X$ with $\pi_1(X)\cong(Z/2Z)\times\kappa$ 
and $\pi_2(X)\cong\mathbb{Z}$ is homotopy equivalent to the total space 
of an $RP^2$-bundle, 
and there are two such bundle spaces for each pair $(\pi,u)$, 
distinguished by $w_2$. 
(See \S5.3.)
As this case is well-understood, 
we shall assume in this section that $\pi\not\cong(Z/2Z)\times\kappa$.
Hence $M$ is an $S^2$-orbifold bundle space,
and $k_1(M)=\beta^u(U^2)$, by Theorem 10.6.

\begin{theorem} 
Let $M$ be an $S^2$-orbifold bundle space with $\pi=\pi_1(M)\not=1$. 
Then $(\pi,u)$ is realized by at most two homotopy types 
of $PD_4$-complexes $X$.
\end{theorem}

\begin{proof} 
Let $p:\widetilde P\simeq K(Z,2)\to{P}$ be the universal covering of $P=P_2(M)$.
The action of $\pi$ on $\pi_2(M)$ also determines $w_1(M)$, by Lemma 10.3.
As $f_M :M\to P$ is 3-connected we may define a class $w$ in 
$H^1 (P;\mathbb{Z}/2\mathbb{Z})$ by $f^*_M w=w_1 (M)$. 
Let $S^{PD}_4 (P)$ be the set of ``polarized" $PD_4 $-complexes $(X,f)$, 
where $f:X\to P$ is 3-connected and $w_1 (X)=f^* w$, 
modulo homotopy equivalence over $P$.
(Note that as $\pi$ is one-ended the universal cover of $X$ 
is homotopy equivalent to $S^2 $).
Let $[X]$ be the fundamental class of $X$ in $H_4 (X;\mathbb{Z}^w)$. 
It follows as in \cite[Lemma 1.3]{[HK88]} 
that given two such polarized complexes 
$(X,f)$ and $(Y,g)$ there is a map $h:X\to Y$ with $gh=f$ if and only if
$f_* [X]=g_* [Y]$ in $H_4 (P;\mathbb{Z}^w )$. 
Since $\widetilde X\simeq\widetilde Y\simeq S^2 $ and $f$ and $g$ are
3-connected such a map $h$ must be a homotopy equivalence.

From the Cartan-Leray homology spectral sequence for the classifying map 
$c_P:P\to{K(\pi,1)}$ we see that there is an exact sequence
\begin{equation*}
0\to H_2(\pi;H_2(\widetilde P)\otimes\mathbb{Z}^w)/\mathrm{Im}(d^2_{5,0})
\to{H_4(P;\mathbb{Z}^w)/J}\to{H_4(\pi;\mathbb{Z}^w)},
\end{equation*}
where $J=H_0(\pi;H_4(\widetilde P;\mathbb{Z})\otimes\mathbb{Z}^w)/
\mathrm{I}m(d^2_{3,2}+d^4_{5,0})$ 
is the image of ${H_4(\widetilde P;\mathbb{Z})\otimes\mathbb{Z}^w}$ 
in $H_4(P;\mathbb{Z}^w)$.
On comparing this spectral sequence with that for $c_X$ we see that 
$H_3(f;\mathbb{Z}^w)$ is an isomorphism and that $f$ 
induces an isomorphism from
$H_4(X;\mathbb{Z}^w)$ to $H_4 (P;\mathbb{Z}^w)/J$.
Hence 
\[
J\cong\mathrm{Coker}(H_4(f;\mathbb{Z}^w))=H_4(P,X;\mathbb{Z}^w)\cong 
H_0(\pi;H_4(\widetilde P,\widetilde X;\mathbb{Z})\otimes\mathbb{Z}^w),
\]
by the exact sequence of homology with coefficients $\mathbb{Z}^w$ 
for the pair $(P,X)$.     
Since $H_4(\widetilde P,\widetilde X;\mathbb{Z})\cong\mathbb{Z}$ 
as a $\pi$-module 
this cokernel is $\mathbb{Z}$ if $w=0$ and $Z/2Z$ otherwise.
(In particular, $d^2_{32}=d^4_{50}=0$, 
since $E^2_{04}$ is $\mathbb{Z}$ or $Z/2Z$ and $H_i(\pi;\mathbb{Z})$ 
is finite for $i>2$.)
In other words, $p$ and $f_M$ induce an isomorphism
\[
H_0(\pi;H_4(\widetilde{P};\mathbb{Z}^w))\oplus H_4(M;\mathbb{Z}^w)=
H_0(\pi;\mathbb{Z}^w)\oplus H_4(M;\mathbb{Z}^w)\cong{H_4(P;\mathbb{Z}^w)}.
\]
Let $\mu=f_{M*}[M]\in H_4(P;\mathbb{Z}^w)$,
and let $G$ be the group of (based) self homotopy equivalences 
of $P$ which induce the identity on $\pi$ and $\pi_2(P)$. 
Then $G\cong{H^2(\pi;\mathbb{Z}^u)}$ \cite{[Ru92]}.
We shall show that there are at most 2 orbits of fundamental classes 
of such polarized complexes (up to sign) under the action of $G$.

This is clear if $w\not=0$, so we may assume that $w=0$.
Suppose first that $u=0$, so $\pi$ is a $PD_2^+$-group and $k_1(M)=0$.
Let $S=K(\pi,1)$ and $C=CP^\infty=K(\mathbb{Z},2)$.
Then $P\simeq S\times{C}$, so $[P,P]\cong [P,S]\times [P,C]$ and $G\cong [S,C]$.
The group structure on $[S,C]$ is determined by the loop-space multiplication 
$m:C\times{C}\to{C}\simeq\Omega{K}(\mathbb{Z},3)$.
This is characterized by the property $m^*z=z\otimes1+1\otimes{z}$,
where $z$ is a generator of $H^2(CP^\infty;\mathbb{Z})$.
The action of $G$ on $P$ is given by 
$\tilde g(s,c)=(s,m(g(s),c))$ for all $g\in G$ and
$(s,c)\in S\times{C}$.

Let $\sigma$ and $\gamma$ be fundamental classes for $S$
and $CP^1$, respectively.
The inclusion of $CP^1$ into $C$ induces a bijection $[S,CP^1]=[S,C]$,
and the degree of a representative map of surfaces
determines an isomorphism $d:[S,C]=H^2(\pi;\mathbb{Z})$.
Let $j:S\times{CP^1}\to{S\times{C}}$ be the natural inclusion.
Then $\omega=j_*(\sigma\otimes\gamma)$ is the image of the fundamental class 
of $S\times{C}$ in $H_4(P;\mathbb{Z}^w)$ and
$\mu\equiv\omega$ modulo $H_0(\pi;H_4(\widetilde{P};\mathbb{Z})^w)$.
Since $\tilde{g}_*\sigma=\sigma+d(g)\gamma$ and $\tilde{g}_*\gamma=\gamma$,
it follows that $\tilde{g}_*\omega=\omega+d(g)m_*[\gamma\otimes\gamma]$.

Since $m^*(z^2)=z^2\otimes1+2z\otimes{z}+1\otimes{z^2}$ and the restriction of 
$z^2$ to $CP^1$ is trivial it follows that $m^*(z^2)([\gamma\otimes\gamma])=2$,
and so $m_*[\gamma\otimes\gamma]=2[CP^2]$,
where $[CP^2]$ is the canonical generator of $H_4(CP^\infty;\mathbb{Z})$.
Hence there are two $G$-orbits of elements in $H_4(P;\mathbb{Z}^w)$ 
whose images agree with $\mu$ modulo $H_0(\pi;H_4(\widetilde{P};\mathbb{Z}))$.

In general let $M_\kappa$ and $P_\kappa$ denote the covering spaces 
corresponding to the subgroup $\kappa$,
and let $G_\kappa$ be the group of self homotopy equivalences of $P_\kappa$.
Lifting self homotopy equivalences defines a homomorphism 
from $G$ to $G_\kappa$, 
which may be identified with the restriction from 
$H^2(\pi;\mathbb{Z}^u)$ to $H^2(\kappa;\mathbb{Z})\cong\mathbb{Z}$,
and which has image of index $\leq2$ (see \cite{[Ts80]}).
Let $q:P_\kappa\to P$ and $q_M:M_\kappa\to{M}$ be the projections.
Then $q_*f_{M_\kappa*}[M_\kappa]=2\mu$ modulo 
$H_0(\pi;H_4(\widetilde{P};\mathbb{Z}))$.
It follows easily that if $g\in G$ and $d(g|_\kappa)=d$ then
$\tilde g_*(\mu)=\mu+d[CP^2]$.
Thus there are again at most two $G$-orbits of elements 
in $H_4(P;\mathbb{Z}^w)$ whose images agree with $\mu$ modulo
$H_0(\pi;H_4(\widetilde{P};\mathbb{Z})^w)$.
This proves the theorem. 
\end{proof}

It is clear from the above argument that the polarized complexes 
are detected by the image of the {\it mod}-(2) fundamental class in
$H_4(P;\mathbb{F}_2)\cong\mathbb{F}_2^2$, 
which is generated by the images of $[M]$ and $[CP^2]$.

\begin{theorem}
Let $M$ be a $S^2$-orbifold bundle space with $\pi=\pi_1(M)\not=1$.
Then the images of $[M]$ and $[M^\tau]$ in $H_4(P_2(M);\mathbb{F}_2)$ 
are distinct.
\end{theorem}

\begin{proof} 
We may assume that $M$ is the geometric 4-manifold realizing
$(\pi,u)$, as constructed in Theorem 10.2.
Let $\kappa=\mathrm{Ker}(u)$.
The Postnikov map $f_M$ (given by $f_M([s,x])=[s,s,x]$ 
for $(s,x)\in{S^2}\times{X}$) embeds $M$ as a submanifold of
${CP^1\times{S^2}\times{K(\kappa,1)}/\sim}$ in $P=P_2(M)$.
The projection of ${CP^\infty\times{S^\infty}\times{K(\kappa,1)}}$
onto its first two factors induces a map 
$g:P\to{Q}=CP^\infty\times{S^\infty}/(z,s)\sim(\sigma(z),-s)$
which is in fact a bundle projection with fibre $K(\kappa,1)$.
Since $gf_M$ factors through $S^2$ the image of $[M]$
in $H_4(Q;\mathbb{F}_2)$ is trivial.
 
Let $v:S^2\times{D^2}\to{V}\subset{M}$ be a fibre-preserving
homeomorphism onto a regular neighbourhood of a general fibre.
Since $V$ is 1-connected $f_M|_V$ factors through 
${CP^\infty\times{S^\infty}\times{K(\kappa,1)}}$.
Let $f_1$ and $f_2$ be the composites of a fixed lift of 
$f_Mv\tau:S^2\times{S^1}\to{P}$ with the projections to $CP^\infty$ 
and $S^\infty$, respectively.
Let $F_1$ be the extension of $f_1$ given by
\[F_1([z_0:z_1],d)=[dz_0:z_1:(1-|d|)z_0]\]
for all $[z_0:z_1]\in{S^2}=CP^1$ and $d\in{D^2}$.
Since $f_2$ maps $S^2\times{S^1}$ to $S^2$ it is nullhomotopic in $S^3$,
and so extends to a map $F_2:S^2\times{D^2}\to{S^3}$.
Then the map $F:M^\tau\to{P}$ given by $f$ 
on $M\setminus{N}$ and $F(s,d)=[F_1(s),F_2(s),d]$ 
for all $(s,d)\in{S^2\times{D^2}}$ is 3-connected,
and so $F=f_{M^\tau}$.

Now $F_1$ maps the open subset $U=\mathbb{C}\times{intD^2}$ 
with $z_0\not=0$ bijectively onto its image in $CP^2$,
and maps $V$ onto $CP^2$.
Let $\Delta$ be the image of $CP^1$ under the diagonal embedding in
$CP^1\times{CP^1}\subset{CP^2\times{S^3}}$.
Then $(F_1,F_2)$ carries $[V,\partial{V}]$ to the image of $[CP^2,CP^1]$ in
$H_4(CP^2\times{S^3},\Delta;\mathbb{F}_2)$.
The image of $[V,\partial{V}]$ generates $H_4(M,M\setminus{U};\mathbb{F}_2)$.
A diagram chase now shows that $[M^\tau]$ and $[CP^2]$ have the same
image in $H_4(Q;\mathbb{F}_2)$, and so $[M^\tau]\not=[M]$
in $H_4(P_2(M);\mathbb{F}_2)$.
\end{proof}

It remains to consider the action of $Aut(P)$.
Since $M$ is geometric $Aut(\pi)$ acts isometrically.
The antipodal map on the fibres defines a self-homeomorphism which
induces $-1$ on $\pi_2(M)$.
These automorphisms clearly fix $H_4(P;\mathbb{F}_2)$.
Thus it is enough to consider the action of 
$G=H^2(\pi;\mathbb{Z}^u)$ on $H^2(\pi;\mathbb{Z}^u)$.

\begin{cor}
Every $4$-manifold realizing $(\pi,u)$ is homotopy equivalent 
to $M$ or $M^\tau$.
If $B=X/\pi$ has no reflector curves then $M^\tau\not\simeq{M}$.
\end{cor}

\begin{proof}
The first assertion holds since the image of the fundamental class
in $H_4(P_2(M);\mathbb{F}_2)$ must generate {\it mod} $[CP^2]$,
and so be $[M]$ or $[M]+[CP^2]$.

If $B$ is nonsingular then Gluck reconstruction changes the
self-intersection of a section,
and hence changes the Wu class $v_2(M)$. 
If $B$ has cone points but no reflector curves then $H^2(\pi;\mathbb{Z}^u)=0$,
by Theorem 10.6, and so $M^\tau\not\simeq{M}$, by Theorem 10.8.
\end{proof} 

If the base $B$ has a reflector curve which is ``untwisted" for $u$,
then $p$ and $p^\tau$ are isomorphic as orbifold bundles over $B$, 
and so $M^\tau\cong{M}$. (See \cite{[Hi13]}.)

\section{Bundle spaces are geometric}

All $\mathbb{S}^2\times\mathbb{X}^2$-manifolds are total spaces of orbifold
bundles over $\mathbb{X}^2$-orbifolds.
We shall determine the $S^2$- and $RP^2$-bundle spaces among them in terms of
their fundamental groups, 
and then show that all such bundle spaces are geometric.

\begin{lemma} 
Let $J=(A,\theta)\in O(3)\times Isom(\mathbb{X}^2)$ be 
an isometry of order $2$ which is fixed point free. Then $A=-I$. 
If moreover $J$ is orientation reversing then $\theta=id_X$ 
or has a single fixed point.
\end{lemma} 

\begin{proof} Since any involution of $\mathbb{R}^2$ (such as $\theta$) 
must fix a point, a line or be the identity,
$A\in O(3)$ must be a fixed point free involution, and so $A=-I$. 
If $J$ is orientation reversing 
then $\theta$ is orientation preserving, and so must fix a point or 
be the identity.
\end{proof}

\begin{theorem}                                  
Let $M$ be a closed $\mathbb{S}^2\times \mathbb{X}^2$-manifold with fundamental group $\pi$. Then
\begin{enumerate}
\item $M$ is the total space of an orbifold bundle with base an
$\mathbb{X}^2$-orbifold and general fibre $S^2$ or $RP^2$;

\item $M$ is the total space of an $S^2$-bundle over a closed 
aspherical surface if and only if $\pi$ is torsion-free; 

\item $M$ is the total space of an $RP^2$-bundle over a closed aspherical 
surface if and only if $\pi\cong (Z/2Z)\times K$, where $K$ is torsion-free. 
\end{enumerate}
\end{theorem}

\begin{proof}
(1)\qua 
The group $\pi$ acts freely and cocompactly on ${S^2\times\mathbb{R}^2}$,
and is a discrete subgroup of 
$Isom({\mathbb{S}^2\times\,\mathbb{X}^2})=O(3)\times Isom(\mathbb{X}^2)$. 
In particular, $N=\pi\cap(O(3)\times\{ 1\})$ is finite and acts freely
on $S^2$, so has order $\leq2$.
Let $p_1$ and $p_2$ be the projections of 
$Isom(\mathbb{S}^2\times\mathbb{X}^2)$ onto $O(3)$ and $Isom(\mathbb{X}^2)$,
respectively.
Then $p_2(\pi)$ is a discrete subgroup of $Isom(\mathbb{X}^2)$ 
which acts cocompactly on $\mathbb{R}^2$, 
and so has no nontrivial finite normal subgroup.
Hence $N$ is the maximal finite normal subgroup of $\pi$.
Projection of $S^2\times\mathbb{R}^2$ onto $\mathbb{R}^2$ induces 
an orbifold bundle projection of $M$ onto $p_2(\pi)\backslash\mathbb{R}^2$ 
with general fibre $N\backslash S^2$.
If $N\not=1$ then $N\cong Z/2Z$ and $\pi\cong (Z/2Z)\times\kappa$, where
$\kappa=\mathrm{Ker}(u)$ is a $PD_2$-group, by Theorem 5.14.

(2)\qua  The condition is clearly necessary.
(See Theorem 5.10).
The kernel of the projection of $\pi$ onto its image 
in $Isom(\mathbb{X}^2)$ is the subgroup $N$.
Therefore if $\pi$ is torsion-free it is isomorphic to its image
in $Isom(\mathbb{X}^2)$, which acts freely on $\mathbb{R}^2$.
The projection $\rho:S^2\times\mathbb{R}^2\to\mathbb{R}^2$ induces a map 
$r:M\to \pi\backslash\mathbb{R}^2$,
and we have a commutative diagram:
\begin{equation*}
\begin{CD}
S^2\times\mathbb{R}^2           @>\rho>>       \mathbb{R}^2 \\
@VVfV                                   @VV{\bar f}V \\
{M=\pi\backslash (S^2\times\mathbb{R}^2)}  @>r>>   \pi\backslash\mathbb{R}^2
\end{CD}
\end{equation*}
where $f$ and $\bar f$ are covering projections.
It is easily seen that $r$ is an $S^2$-bundle projection.

(3)\qua  The condition is necessary, by Theorem 5.16.
Suppose that it holds.
Then $K$ acts freely and properly discontinuously on $\mathbb{R}^2$, 
with compact quotient. 
Let $g$ generate the torsion subgroup of $\pi$. 
Then $p_1(g)=-I$, by Lemma 10.9.
Since $p_2(g)^2=\mathrm{id}_{\mathbb{R}^2}$ the fixed point set
$F=\{ x\in\mathbb{R}^2\mid p_2(g)(x)=x\}$ is nonempty, and is either a point, 
a line, or the whole of $\mathbb{R}^2$. 
Since $p_2(g)$ commutes with the action of $K$ on $\mathbb{R}^2$ 
we have $KF=F$, and so $K$ acts freely and properly discontinuously on $F$. 
But $K$ is neither trivial nor infinite cyclic, 
and so we must have $F=\mathbb{R}^2$. 
Hence $p_2(g)=\mathrm{id}_{\mathbb{R}^2}$. 
The result now follows, as $K\backslash(S^2\times\mathbb{R}^2)$ 
is the total space of an $S^2$-bundle over $K\backslash\mathbb{R}^2$, 
by part (1), and $g$ acts as the antipodal involution on the fibres.
\end{proof}

If the $\mathbb{S}^2\times\mathbb{X}^2$-manifold $M$ is the total space of an 
$S^2$-bundle $\xi$ then $w_1(\xi )=u$ and is detected by the determinant: 
$\det (p_1(g))=(-1)^{w_1(\xi )(g)}$ for all $g\in \pi$.

The total space of an $RP^2$-bundle over $B$ is the quotient 
of its orientation double cover (which is an $S^2$-bundle over $B$) 
by the fibrewise antipodal involution, 
and so there is a bijective correspondance between 
orientable $S^2$-bundles over $B$ 
and $RP^2$-bundles over $B$.

Let $\mathbf{i}= 
\left(\begin{smallmatrix}
1\\
0
\end{smallmatrix}\right))$ and
 $\mathbf{j}=
\left(\begin{smallmatrix}
0\\
1
\end{smallmatrix}\right))$, 
and let $(A,\beta,C)\in O(3)\times E(2)=O(3)\times(\mathbb{R}^2\rtimes O(2))$ 
be the $\mathbb{S}^2\times\mathbb{E}^2$-isometry which sends 
$(v,x)\in S^2 \times \mathbb{R}^2 $ to $(Av,Cx+\beta)$.

\begin{theorem}                                   
Let $M$ be the total space of an $S^2$- or $RP^2$-bundle over $T$ or $Kb$.
Then $M$ admits the geometry $\mathbb{S}^2\times\mathbb{E}^2$.
\end{theorem}

\begin{proof}
Let $R_i \in O(3)$ be the reflection of $\mathbb{R}^3 $ 
which changes the sign of the $i^{th}$ coordinate, for $i=1,2,3$. 
If $A$ and $B$ are products of such reflections then the subgroups of 
$Isom(\mathbb{S}^2\times\mathbb{E}^2)$ generated by 
$\alpha=(A,\mathbf{i},I)$ and $\beta=(B,\mathbf{j},I)$ are discrete, 
isomorphic to $\mathbb{Z}^2$ and act freely and cocompactly on 
$S^2 \times\mathbb{R}^2 $. Taking 
\begin{enumerate}
\item $A=B=I$; 
\item $A=R_1 R_2 ,B=R_1 R_3 $;
\item $A=R_1,B=I$; and 
\item $A=R_1 ,B=R_1 R_2 $
\end{enumerate}
gives four $S^2 $-bundles $\eta_i$ over the torus. 
If instead we use the isometries
$\alpha=(A,\frac12\mathbf{i},
\left(\begin{smallmatrix}
1&0\\
0&-1
\end{smallmatrix}\right))$
and $\beta=(B,\mathbf{j},I)$
we obtain discrete subgroups isomorphic to 
$\mathbb{Z}\rtimes_{-1}\!\mathbb{Z}$ which act 
freely and cocompactly.
Taking 
\begin{enumerate}
\item $A=R_1,B=I$; 
\item $A=R_1, B=R_2 R_3 $; 
\item $A=I,B=R_1$; 
\item $A=R_1 R_2 ,B=R_1 $; 
\item $A=B=I$; and 
\item  $A=I,B=R_1R_2$
\end{enumerate}
gives six $S^2 $-bundles $\xi_i$ over the Klein bottle.

To see that these are genuinely distinct, 
we check first the fundamental groups, 
then the orientation character of the total space; 
consecutive pairs of generators determine 
bundles with the same orientation character, 
and we distinguish these by means of the second Stiefel-Whitney classes, 
by computing the self-intersections of sections of the bundle.
(See Lemma 5.11.(3).)
We shall use the stereographic projection of 
$S^2\subset\mathbb{R}^3=\mathbb{C}\times\mathbb{R}$ onto 
$CP^1=\mathbb{C}\cup \{ \infty \}$, 
to identify the reflections $R_i: S^2\to S^2$ with the antiholomorphic involutions:
\[z {\buildrel R_1\over \longmapsto } -\overline{z},\quad
{z} {\buildrel R_2\over \longmapsto }\overline{z},\quad
{z} {\buildrel R_3\over \longmapsto}{\overline{z}}^{-1}.\]
Let $\mathfrak{T}=\{ (s,t)\in\mathbb{R}^2\mid 0\leq s,t\leq 1\}$
be the fundamental domain for the standard action of 
$\mathbb{Z}^2$ on $\mathbb{R}^2$.
Sections of $\xi_i$ correspond to 
functions $\sigma:\mathfrak{T}\to{S^2}$ such that $\sigma(1,t)=A(\sigma(0,t))$ 
and $\sigma(s,1)=B(\sigma(s,0))$ for all $(s,t)\in\partial\mathfrak{T}$.
In particular, points of $S^2$ fixed by both $A$ and $B$ determine sections.

As the orientable cases ($\eta_1$, $\eta_2$, $\xi_1$ and $\xi_2$)
have been treated in \cite{[Ue90]} 
we may concentrate on the nonorientable cases.
In the case of $\eta_3$ fixed points of $A$ determine sections, since $B=I$,
and the sections corresponding to distinct fixed points are disjoint.
Since the fixed-point set of $A$ on $S^2$ is a circle 
all such sections are isotopic.
Therefore $\sigma\cdot\sigma=0$, so $v_2(M)=0$ and hence $w_2(\eta_3)=0$.

We may define a 1-parameter family of sections for $\eta_4$ by
\begin{equation*}
\sigma_\lambda(s,t)=
\lambda(2t^2-1+i(2s-1)(2t-1)).
\end{equation*}
Now $\sigma_0$ and $\sigma_1$ intersect transversely in a single point,
corresponding to $s=\frac12$ and $t=\frac1{\sqrt2}$.
Hence $\sigma \cdot \sigma =1$, so $v_2(M)\not=0$ and $w_2(\eta_4 )\not=0$.

The remaining cases correspond to $S^2$-bundles over $Kb$ with
nonorientable total space.                            
We now take 
$\mathfrak{K}=
\{(s,t)\in\mathbb{R}^2\mid 0\leq{s}\leq\frac12, |t|\leq\frac12 \}$
as the fundamental domain for the action of 
$\mathbb{Z}\rtimes_{-1}\! \mathbb{Z}$ on $\mathbb{R}^2$,
and seek functions $\sigma:\mathfrak{K}\to{S^2}$ such that
$\sigma(\frac12,t)=A(\sigma(0,-t))$ and $\sigma(s,\frac12 )=B(\sigma(s,-\frac12))$
for all $(s,t)\in\partial\mathfrak{K}$.

The cases of $\xi_3$ and $\xi_5$ are similar to that of $\eta_3$: 
there are obvious one-parameter families of disjoint sections,
and so $w_2(\xi_3)=w_2(\xi_5)=0$. However $w_1(\xi_3)\not= w_1(\xi_5)$.
(In fact $\xi_5$ is the product bundle).
               
The functions $\sigma_{\lambda}(s,t)=\lambda (4s-1+it)$
define a 1-parameter family of sections for $\xi_4$ such that 
$\sigma_0$ and $\sigma_1$ 
intersect transversely in one point, so that $\sigma \cdot \sigma =1$.
Hence $v_2(M)\not=0$ and so $w_2(\xi_4)\not=0$. 

For $\xi_6$ the functions 
$\sigma_{\lambda}(s,t)=\lambda(4s-1)t+i(1-\lambda)(4t^2-1)$
define a 1-parameter family of sections such that $\sigma_0$ 
and $\sigma_1(s,t)$ intersect 
transversely in one point, so that $\sigma \cdot \sigma=1$. 
Hence $v_2(M)\not=0$ and so $w_2(\xi_6)\not=0$. 

Thus these bundles are all distinct, and so all $S^2$-bundles over $T$ or $Kb$ 
are geometric of type $\mathbb{S}^2\times\mathbb{E}^2$.

Adjoining the fixed point free involution $(-I,0,I)$ 
to any one of the above ten sets of generators 
for the $S^2 $-bundle groups amounts to dividing out 
the $S^2 $ fibres by the antipodal map 
and so we obtain the corresponding $RP^2 $-bundles. 
(Note that there are just four such $RP^2$-bundles -- 
but each has several distinct double covers which are $S^2 $-bundles). 
\end{proof}

\begin{theorem}                                                    
Let $M$ be the total space of an $S^2$- or $RP^2$-bundle over a closed hyperbolic surface.
Then $M$ admits the geometry $\mathbb{S}^2\times\mathbb{H}^2$.
\end{theorem}

\begin{proof} 
Let $T_g$ be the closed orientable surface of genus $g$, 
and let ${\mathfrak T}^g\subset \mathbb {H}^2$ be 
a $2g$-gon representing the fundamental domain of $T_g$.
The map $\Omega : {\mathfrak T}^g\to \mathfrak T $ that collapses $2g-4$ 
sides of ${\mathfrak T}^g$ 
to a single vertex in the rectangle ${\mathfrak T}$ induces a degree-1 map 
$\widehat\Omega$ from $T_g$ to $T$ 
that collapses $g-1$ handles on $T^g$ to a single point on $T$.                                                               
(We may assume the induced epimorphism from 
\begin{equation*}
\pi_1(T_g)=\langle a_1,b_1,\dots ,a_g,b_g\mid \Pi_{i=1}^g[a_i,b_i]=1\rangle
\end{equation*}
 to $Z^2$ kills the generators $a_j,b_j$ for $j>1$).
Hence given an $S^2$-bundle $\xi$ over $T$ with total space
$M_{\xi }=\Gamma_\xi\backslash(S^2\times E^2)$, where 
\begin{equation*}
\Gamma_\xi=\{ (\xi(h),h)\mid h\in \pi_1(T)\} \leq Isom(S^2\times E^2)
\end{equation*}
and $\xi:\mathbb{Z}^2\to O(3)$ is as in Theorem 10.11,
the pullback $\widehat\Omega^*(\xi )$ is an $S^2$-bundle over $T_g$, 
with total space 
$M_{\xi\Omega }=\Gamma_{\xi\Omega}\backslash (S^2\times \mathbb{H}^2)$, 
where
$\Gamma_{\xi\Omega}=\{ (\xi{\widehat\Omega}_*(h),h)\mid h\in \Pi_1(T^g)\} 
\leq Isom(S^2\times \mathbb{H}^2)$.
As $\widehat\Omega $ is a degree-1 map,
it induces monomorphisms in cohomology,
so $w(\xi)$ is nontrivial if and only if  
$w(\widehat\Omega^* (\xi ))=\widehat\Omega^* w(\xi )$
is nontrivial.
Hence all $S^2$-bundles over $T^g$ for $g\geq 2$ are geometric of type $S^2\times \mathbb{H}^2$.

Suppose now that $B$ is the closed surface 
$\sharp^3{RP^2}=T\sharp{RP^2}=Kb\sharp{RP^2}$. 
Then there is a map $\widehat\Omega:T\sharp{RP^2}\to RP^2$
that collapses the torus summand to a single point. 
This map $\widehat\Omega$ is again a degree-1 map,
and so induces monomorphisms in cohomology. 
In particular $\widehat\Omega^*$ preserves the orientation character, that is 
$w_1(\widehat\Omega^*(\xi))=\widehat\Omega^*w_1(RP^2)= w_1(B)$, and is an isomorphism on $H^2$.
We may pull back the four $S^2$-bundles $\xi$ over $RP^2$ along 
$\widehat\Omega$ to obtain the four bundles over $B$ with first Stiefel-Whitney 
class $w_1(\widehat\Omega^* \xi)$ either 0 or $w_1(B)$.
   
Similarly there is a map $\widehat\Upsilon: Kb\sharp{RP^2}\to RP^2$
that collapses the Klein bottle summand to a single point. 
This map $\widehat\Upsilon$ has degree 1 {\it mod} (2) so that 
${\widehat\Upsilon}^*w_1(RP^2)$ has nonzero square since $w_1(RP^2)^2\not=0$. 
Note that in this case $\widehat\Upsilon^*w_1(RP^2)\not=w_1(B)$.
Hence we may pull back the two $S^2$-bundles $\xi$ over $RP^2$ with 
$w_1(\xi )=w_1(RP^2)$ to obtain a further two bundles over $B$ with 
$w_1(\widehat\Upsilon^* (\xi ))^2=\widehat\Upsilon^*w_1(\xi )^2\not=0$, 
as $\widehat\Upsilon$ is a ring monomorphism.

There is again a map $\widehat\Theta: Kb\sharp{RP^2}\to Kb$
that collapses the projective plane summand to a single point. 
Once again $\widehat\Theta$ is of degree 1 {\it mod} (2),
so that we may pull back the two $S^2$-bundles $\xi$ 
over $Kb$ with $w_1(\xi)=w_1(Kb)$ along $\widehat\Theta$ 
to obtain the remaining two $S^2$-bundles over $B$. 
These two bundles $\widehat\Theta^*(\xi )$ have 
$w_1(\widehat\Theta^*(\xi ))\not=0$ but $w_1(\widehat\Theta^*(\xi ))^2=0$;
as $w_1(Kb)\not=0$ but $w_1(Kb)^2=0$ and $\widehat\Theta^*$ is a monomorphism.

Similar arguments apply to bundles over $\sharp^n RP^2$ where $n>3$. 

Thus all $S^2$-bundles over all closed aspherical surfaces are geometric. 
Furthermore since the antipodal involution of a geometric $S^2$-bundle 
is induced by an isometry $(-I,id_{\mathbb{H}^2})\in O(3)\times Isom(\mathbb{H}^2)$
we have that all $RP^2$-bundles over closed aspherical surfaces are geometric. 
\end{proof}

An alternative route to Theorems 10.11 and 10.12 would be to first show that 
orientable 4-manifolds which are total spaces of $S^2$-bundles are geometric, 
then deduce that $RP^2$-bundles are geometric (as above); 
and finally observe that every $S^2$-bundle space double covers 
an $RP^2$-bundle space.

The other $\mathbb{S}^2\times \mathbb{X}^2$-manifolds are 
$S^2$-orbifold bundle spaces. 
It can be shown that there are at most two such orbifold 
bundles with given base orbifold $B$ and action $u$,  
differing by (at most) Gluck reconstruction.
If $B$ has a reflector curve they are isomorphic as orbifold bundles.
Otherwise $B$ has an even number of cone points and the bundles are distinct.
The bundle space is geometric except when $B$ is orientable and $\pi$ 
is generated by involutions,
in which case the action is unique and there 
is one non-geometric orbifold bundle.
(See \cite{[Hi13]}.)

If $\chi(F)<0$ or $\chi(F)=0$ and $\partial=0$ then every $F$-bundle 
over $RP^2$ is geometric, by
Lemma 5.21 and the remark following Theorem 10.2.

It is not generally true that the projection of $S^2\times X$ onto 
$S^2$ induces an orbifold bundle projection from $M$ to an 
$\mathbb{S}^2$-orbifold.
For instance, if $\rho$ and $\rho'$ are rotations of $S^2$ about a common axis
which generate a rank 2 abelian subgroup of $SO(3)$
then $(\rho,(1,0))$ and $(\rho',(0,1))$ generate a discrete subgroup of
$SO(3)\times\mathbb{R}^2$ which acts freely, cocompactly and isometrically
on $S^2\times\mathbb{R}^2$.
The orbit space is homeomorphic to $S^2\times T$.
(It is an orientable $S^2$-bundle over the torus,
with disjoint sections, detemined by the ends of the axis of the rotations).
Thus it is Seifert fibred over $S^2$, but the fibration is not
canonically associated to the metric structure, for $\langle\rho,\rho'\rangle$ 
does not act properly discontinuously on $S^2$.

\section{Fundamental groups of $\mathbb{S}^2\times\mathbb{E}^2$-manifolds}

We shall show first that if $M$ is a closed 4-manifold 
any two of the conditions ``$\chi(M)=0$", 
``$\pi_1(M)$ is virtually $\mathbb{Z}^2$" and 
``$\pi_2(M)\cong\mathbb{Z}$" imply the third, 
and then determine the possible fundamental groups.

\begin{theorem} 
Let $M$ be a closed $4$-manifold with fundamental group $\pi$.
Then the following conditions are equivalent:
\begin{enumerate}
\item $\pi$ is virtually $\mathbb{Z}^2 $ and $\chi(M)=0$;

\item $\pi$ has an ascendant infinite restrained subgroup and 
$\pi_2 (M)\cong\mathbb{Z}$;

\item $\chi(M)=0$ and $\pi_2 (M)\cong\mathbb{Z}$; and

\item $M$ has a covering space of degree
dividing $4$ which is homeomorphic to $S^2\times T$.

\item $M$ is virtually homeomorphic to an $\mathbb{S}^2\times\mathbb{E}^2$-manifold.
\end{enumerate}
\end{theorem} 

\begin{proof} 
If $\pi$ is virtually a $PD_2$-group and either $\chi(\pi)=0$ or $\pi$ has an
ascendant infinite restrained subgroup then $\pi$ is virtually $\mathbb{Z}^2$.
Hence the equivalence of these conditions follows from Theorem 10.1, with the
exception of the assertions regarding homeomorphisms, which then follow from
Theorem 6.11.
\end{proof}

We shall assume henceforth that the conditions of Theorem 10.13 hold,
and shall show next that there are nine possible groups. 
Seven of them are 2-dimensional crystallographic groups, and we shall give 
also the name of the corresponding $\mathbb{E}^2$-orbifold, 
following \cite[Appendix A]{[Mo]}.
(The restriction on finite subgroups eliminates the remaining ten
$\mathbb{E}^2$-orbifold groups from consideration).
                                             
\begin{theorem} 
Let $M$ be a closed $4$-manifold 
such that $\pi=\pi_1(M)$ is virtually $\mathbb{Z}^2$ and $\chi(M)=0$.   
Let $F$ be the maximal finite normal subgroup of $\pi$. 
If $\pi$ is torsion-free then either 
\begin{enumerate}
\item $\pi=\sqrt\pi\cong\mathbb{Z}^2$ (the torus); or 

\item $\pi\cong\mathbb{Z}\rtimes_{-1}\mathbb{Z}$ (the Klein bottle).
\end{enumerate}
If $F=1$ but $\pi$ has nontrivial torsion and $[\pi:\sqrt\pi]=2$ 
then either 
\begin{enumerate}
\addtocounter{enumi}{2}
\item $\pi\cong D\times\mathbb{Z}\cong
(\mathbb{Z}\oplus (Z/2Z))*_{\mathbb{Z}} (\mathbb{Z}\oplus(Z/2Z))$, 
with the presentation 

$\langle s,x,y\mid x^2=y^2=1,\medspace sx=xs,\medspace sy=ys\rangle$  
(the silvered annulus $\mathbb{A}$); or 

\item $\pi\cong D\rtimes_\tau\mathbb{Z}\cong
\mathbb{Z}*_{\mathbb{Z}} (\mathbb{Z}\oplus(Z/2Z))$, 
with the presentation

$\langle t,x\mid x^2=1,\medspace t^2x=xt^2\rangle$ 
(the silvered M\"obius band $\mathbb{M}b$); or

\item $\pi\cong (\mathbb{Z}^2)\rtimes_{-I}(Z/2Z)\cong D*_{\mathbb{Z}}D$, 
with the presentations

$\langle s,t,x\mid x^2=1,\medspace xsx=s^{-1}\! ,\medspace xtx=t^{-1}\! ,
\medspace st=ts\rangle$ and (setting $y=xt$) 

$\langle s,x,y\mid x^2=y^2=1,\medspace xsx=ysy=s^{-1}\rangle$ 
(the pillowcase $S(2222)$).
\end{enumerate}
If $F=1$ and $[\pi:\sqrt\pi]=4$ then either
\begin{enumerate}
\addtocounter{enumi}{5}
\item $\pi\cong D*_{\mathbb{Z}}(\mathbb{Z}\oplus(Z/2Z))$, 
with the presentations

$\langle s,t,x\mid x^2=1,\medspace xsx=s^{-1}\! ,\medspace xtx=t^{-1}\! ,
\medspace tst^{-1}=s^{-1}\rangle$ and 

(setting $y=xt$) 
$\langle s,x,y\mid x^2=y^2=1,\medspace xsx=s^{-1}\! ,\medspace
ys=sy\rangle$ ($\mathbb{D}(22)$); or

\item $\pi\cong\mathbb{Z}*_{\mathbb{Z}}D$, with the presentations 

$\langle r,s,x\mid x^2=1,\medspace xrx=r^{-1}\! ,\medspace xsx=rs^{-1}\! ,
\medspace srs^{-1}=r^{-1}\rangle$ and 

(setting $t=xs$) 
$\langle t,x\mid x^2 =1,\medspace xt^2 x=t^{-2}\rangle$ ($P(22)$).
\end{enumerate}
If $F$ is nontrivial then either
\begin{enumerate}
\addtocounter{enumi}{7}

\item $\pi\cong\mathbb{Z}^2\oplus(Z/2Z)$; or 

\item $\pi\cong(\mathbb{Z}\rtimes_{-1}\mathbb{Z})\times(Z/2Z)$.
\end{enumerate}
\end{theorem}

\begin{proof}                                                  
Let $u:\pi\to{Aut(\pi_2 (M))}=\mathbb{Z}^\times$ be the natural homomorphism.
Since $\kappa=\mathrm{Ker}(u)$ is torsion-free it is either 
$\mathbb{Z}^2 $ or $\mathbb{Z}\rtimes_{-1}\mathbb{Z}$; 
since it has index at most 2 it follows that $[\pi:\sqrt\pi]$ divides 4
and $F$ has order at most 2.
If $F=1$ then $\sqrt\pi\cong\mathbb{Z}^2$ and $\pi/\sqrt\pi$ 
acts effectively on $\sqrt\pi$, 
so $\pi$ is a 2-dimensional crystallographic group.
If $F\not=1$ then it is central in $\pi$ and $u$ maps $F$ isomorphically to 
$\mathbb{Z}^\times$, so $\pi\cong (Z/2Z)\times\kappa$. 
\end{proof}

Each of these groups may be realised geometrically, by Theorem 10.2.
It is easy to see that any $\mathbb{S}^2\times{E}^2$-manifold whose fundamental 
group has infinite abelianization is a mapping torus, 
and hence is determined up to diffeomorphism by its homotopy type.
(See Theorems 10.11 and 10.15).
We shall show next that there are four affine diffeomorphism classes of
$\mathbb{S}^2\times\mathbb{E}^2$-manifolds whose fundamental groups have finite 
abelianization.

Let $\Omega$ be a discrete subgroup of 
$Isom(\mathbb{S}^2\times\mathbb{E}^2 )=O(3)\times E(2)$ 
which acts freely and cocompactly on $S^2 \times\mathbb{R}^2$.
If $\Omega\cong D*_{\mathbb{Z}}D$ or $D*_{\mathbb{Z}}(\mathbb{Z}\oplus(Z/2Z))$
it is generated by elements of order 2, 
and so $p_1(\Omega)=\{ \pm I\}$, 
by Lemma 10.9.
Since $p_2(\Omega)<E(2)$ is a 2-dimensional crystallographic group 
it is determined up to conjugacy in 
$Aff(2)=\mathbb{R}^2\rtimes GL(2,\mathbb{R})$ by its isomorphism type,
$\Omega$ is determined up to conjugacy in $O(3)\times Aff(2)$ and the 
corresponding geometric 4-manifold is determined up to affine diffeomorphism.

Although $\mathbb{Z}*_{\mathbb{Z}}D$ is not generated by involutions, 
a similar argument applies.
The isometries $T=(\tau,\frac12\mathbf{j},
\left(\begin{smallmatrix}
-1&0\\
0&1
\end{smallmatrix}\right) )$
and $X=(-I,\frac12(\mathbf{i}+\mathbf{j}),-I)$
generate a discrete subgroup of 
$Isom(\mathbb{S}^2\times\mathbb{E}^2)$
isomorphic to $\mathbb{Z}*_{\mathbb{Z}}D$ and which acts freely 
and cocompactly on $S^2\times\mathbb{R}^2$,
provided $\tau^2=I$. 
Since $x^2=(xt^2)^2=1$ this condition is necessary, by Lemma 10.6.
Conjugation by the reflection across the principal diagonal of $\mathbb{R}^2$
induces an automorphism which fixes $X$ and carries $T$ to $XT$.
Thus we may assume that $T$ is orientation preserving,
i.e., that $\mathrm{det}(\tau)=-1$.
(The isometries $T^2$ and $XT$ then generate $\kappa=\mathrm{Ker}(u)$).
Thus there are two affine diffeomorphism classes of such manifolds,
corresponding to the choices $\tau=-I$ or $R_3$.

None of these manifolds fibre over $S^1$, 
since in each case $\pi/\pi'$ is finite.
However if $\Omega$ is a $\mathbb{S}^2\times\mathbb{E}^2$-lattice such that
$p_1(\Omega)\leq\{\pm I\}$ then 
$\Omega\backslash(S^2\times\mathbb{R}^2)$ fibres over $RP^2$,
since the map sending  $(v,x)\in S^2 \times\mathbb{R}^2$ 
to $[\pm v]\in RP^2$ is compatible with the action of $\Omega$.
If $p_1(\Omega)=\{\pm I\}$ the fibre is $\omega\backslash\mathbb{R}^2$, 
where $\omega=\Omega\cap(\{ 1\}\times E(2))$; 
otherwise it has two components.
Thus three of these four manifolds fibre over $RP^2$
(excepting perhaps only the case
$\Omega\cong\mathbb{Z}*_{\mathbb{Z}}D$ and $R_3\in p_1(\Omega)$).
                             
\section{Homotopy types of $\mathbb{S}^2\times\mathbb{E}^2$-manifolds}

Our next result shows that if $M$ satisfies the conditions of Theorem 10.13 
and its fundamental group has infinite abelianization then it is determined 
up to homotopy by $\pi_1(M)$ and its Stiefel-Whitney classes.

\begin{theorem} 
Let $M$ be a closed $4$-manifold 
such that $\pi=\pi_1(M)$ is virtually $\mathbb{Z}^2$ and $\chi(M)=0$.   
If $\pi/\pi'$ is infinite then $M$ is 
homotopy equivalent to an $\mathbb{S}^2\times\mathbb{E}^2$-manifold which 
fibres over $S^1 $.
\end{theorem}

\begin{proof} The infinite cyclic covering space of $M$ determined by an 
epimorphism $\lambda:\pi\to\mathbb{Z}$ is a $PD_3 $-complex, by Theorem 4.5,
and therefore is homotopy equivalent to 
\begin{enumerate}
\item $S^2 \times S^1 $ 
(if $\mathrm{Ker}(\lambda)\cong\mathbb{Z}$ is torsion-free and 
$w_1 (M)|_{\mathrm{Ker}(\lambda)}=0$), 
\item $S^2\tilde\times S^1 $ 
(if $\mathrm{Ker}(\lambda)\cong\mathbb{Z}$ and 
$w_1 (M)|_{\mathrm{Ker}(\lambda)}\not= 0$), 
\item $RP^2\times S^1 $ 
(if $\mathrm{Ker}(\lambda)\cong\mathbb{Z}\oplus(Z/2Z)$) 
or 
\item $RP^3\sharp RP^3 $ (if $\mathrm{Ker}(\lambda)\cong D$).
\end{enumerate}
Therefore $M$ is homotopy equivalent to the mapping torus $M(\phi)$ 
of a self homotopy equivalence of one of these spaces.

The group of free homotopy classes of self homotopy equivalences 
$E(S^2 \times S^1)$ is generated by the reflections in each factor 
and the twist map, and has order 8.
The groups $E(S^2\tilde\times{S^1})$ and $E(RP^2\times{S^1})$
are each generated by the reflection in the 
second factor and a twist map, and have order 4.
(See \cite{[KR90]} for the case of $S^2\tilde\times{S^1}$.) 
Two of the corresponding mapping tori of self-homeomorphisms of 
$S^2\tilde\times{S^1}$ also arise from self homeomorphisms 
of $S^2\times{S^1}$. 
The other two have nonintegral $w_1$. 
As all these mapping tori are also $S^2$- or $RP^2$-bundles over the torus 
or Klein bottle,
they are geometric by Theorem 10.11.

The group $E(RP^3 \sharp RP^3)$ is generated by the reflection 
interchanging the summands and the fixed point free involution 
(see \cite[page 251]{[Ba']}), 
and has order 4.                                         
Let $\alpha=(-I,0,
\left(\begin{smallmatrix}
-1&0\\
0&1
\end{smallmatrix}\right))$,
$\beta=(I,\mathbf{i},I)$
$\gamma=(I,\mathbf{j},I)$
and 
$\delta=(-I,\mathbf{j},I)$.
Then the subgroups generated by $\{\alpha,\beta,\gamma\}$, $\{\alpha,\beta,\delta\}$, 
$\{\alpha,\beta\gamma\}$ and $\{\alpha,\beta\delta\}$, respectively, give the four $RP^3\sharp RP^3 $-bundles. 
(Note that these may be distinguished by their groups and orientation characters).
\end{proof}

A $T$-bundle over $RP^2$ with $\partial=0$ which does not also fibre 
over $S^1$ has fundamental group $D*_{\mathbb{Z}}D$, 
while the group of a $Kb$-bundle over $RP^2$ which does not also fibre over 
$S^1$ is $D*_{\mathbb{Z}}(\mathbb{Z}\oplus(Z/2Z))$ 
or $\mathbb{Z}*_{\mathbb{Z}}D$.

When $\pi$ is torsion-free every homomorphism from $\pi$ 
to $\mathbb{Z}^\times$ is the orientation character of some $M$ 
with fundamental group $\pi$. 
However if $\pi\cong D\times\mathbb{Z}$ or $D\rtimes_\tau\mathbb{Z}$ 
the orientation character must be trivial on all elements of order 2, 
while it is determined up to composition 
with an automorphism of $\pi$ if $F\not=1$.

\begin{theorem} 
Let $M$ be a closed $4$-manifold such that $\chi(M)=0$ and
$\pi=\pi_1(M)$ is an extension of $\mathbb{Z}$ by a finitely generated 
infinite normal subgroup $\nu$ with a nontrivial finite normal subgroup $F$. 
Then $M$ is homotopy equivalent to the mapping torus of a self homeomorphism 
of $RP^2\times S^1$.
\end{theorem}
                     
\begin{proof} 
The covering space $M_\nu$ corresponding
to the subgroup $\nu$ is a $PD_3$-space, by Theorem 4.5.
Therefore $M_\nu\simeq{RP^2\times S^1}$, by Theorem 2.11.
Since every self-homotopy equivalence of $RP^2\times S^1$
is homotopic to a homeomorphism $M$ is homotopy equivalent to a mapping torus.
\end{proof}

The possible orientation characters for the groups with finite abelianization 
are restricted by Lemma 3.14, which implies that $\mathrm{Ker}(w_1)$
must have infinite abelianization.
For $D*_{\mathbb{Z}}D$ we must have $w_1 (x)=w_1 (y)=-1$ and $w_1 (s)=1$.
For $D*_{\mathbb{Z}}(\mathbb{Z}\oplus (Z/2Z))$ we must have $w_1 (s)=1$ 
and $w_1 (x)=-1$; 
since the subgroup generated by the commutator subgroup and $y$ is isomorphic 
to $D\times\mathbb{Z}$ we must also have $w_1 (y)=1$.
Thus the orientation characters are uniquely determined for these groups.
For $\mathbb{Z}*_{\mathbb{Z}}D$ we must have $w_1 (x)=-1$, 
but $w_1 (t)$ may be either $-1$ or 1. 
As there is an automorphism $\phi$ of $\mathbb{Z}*_{\mathbb{Z}}D$ 
determined by $\phi(t)=xt$ and $\phi(x)=x$ 
we may assume that $w_1 (t)=1$ in this case.

In all cases, to each choice of orientation character there corresponds 
an unique action of $\pi$ on $\pi_2 (M)$, by Lemma 10.9.
However the homomorphism from $\pi$ to $Z/2Z$ determining the action 
may differ from $w_1 (M)$.
(Note also that elements of order 2 must act nontrivially,
by Theorem 10.1).

In the first version of this book we used elementary arguments 
involving cochain computations of cup product in low degrees  
and Poincar\'e duality to compute the cohomology rings
$H^*(M;\mathbb{F}_2)$, $k_1(M)=\beta^u(U^2)$ and $v_2(M)=U^2$,
for the cases with $\pi/\pi'$ finite.
The calculation of $k_1(M)$ has been subsumed into the
(new) Theorem 10.6 above,
while $v_2(M)$ is determined by $\pi$ whenever $\pi$ 
has torsion but is not a direct product.
(See \cite{[Hi13]}.)

Gluck reconstruction of the $\mathbb{S}^2\times\mathbb{E}^2$-manifolds 
with group $D*_{\mathbb{Z}}D$ or $\mathbb{Z}*_{\mathbb{Z}}D$ 
changes the homotopy type, by Corollary 10.8.1.
The two geometric manifolds with $\pi\cong\mathbb{Z}*_{\mathbb{Z}}D$ 
are Gluck reconstructions of each other, 
but there is a nongeometric 4-manifold with $\pi\cong{D*_{\mathbb{Z}}D}$ 
and $\chi=0$.
The $\mathbb{S}^2\times\mathbb{E}^2$-manifold
with $\pi\cong{D*_{\mathbb{Z}}(\mathbb{Z}\oplus(Z/2Z))}$
is isomorphic to its Gluck reconstruction,
and thus every closed manifold with 
$\pi\cong{D*_{\mathbb{Z}}(\mathbb{Z}\oplus(Z/2Z))}$ 
and $\chi=0$ is homotopy equivalent to this manifold.
In summary, there are 22 affine diffeomorphism classes of closed
$\mathbb{S}^2\times\mathbb{E}^2$-manifolds 
and 23 homotopy types of closed 4-manifolds covered by 
$S^2\times\mathbb{R}^2$ and with Euler characteristic 0.
(See \cite{[Hi13]}.)

\section{Some remarks on the homeomorphism types}

In Chapter 6 we showed that if $\pi$ is $\mathbb{Z}^2$ 
or $\mathbb{Z}\rtimes_{-1}\!\mathbb{Z}$ 
then $M$ must be homeomorphic to the total space of an $S^2$-bundle 
over the torus or Klein bottle, and we were able to estimate the size 
of the structure sets when $\pi\cong\kappa\times{Z/2Z}$.
We may apply the Shaneson-Wall exact sequence and results 
on $L_*(D,w)$ from \cite{[CD04]} to obtain $L_*(\pi,w)$ 
for $\pi=D\times\mathbb{Z}$ or $D\rtimes_\tau\mathbb{Z}$.
It follows also that $L_1(\pi,w)$ is not finitely generated
if $\pi=D*_{\mathbb{Z}}D$ or $D*_{\mathbb{Z}}(\mathbb{Z}\oplus(Z/2Z))$ 
and $w$ is as in \S5 above.
For these groups retract onto $D$ compatibly with $w|_D$.
Although $\mathbb{Z}*_{\mathbb{Z}}D$ does not map onto $D$,
L\"uck has shown that $L_1(\mathbb{Z}*_{\mathbb{Z}}D,w)$ 
is not finitely generated (private communication).
(The groups $L_*(\pi)\otimes\mathbb{Z}[\frac12]$ have been computed for all
cocompact planar groups $\pi$, with $w$ trivial \cite{[LS00]}.)

In particular, if $\pi\cong D\times\mathbb{Z}$ 
and $M$ is orientable then $[SM;G/TOP]$ has rank 1, 
while $L_1(D\times\mathbb{Z})$ has rank 3, and so $S_{TOP}(M)$ is infinite.
On the other hand, if $\pi\cong D\times\mathbb{Z}$ and $M$ is non-orientable
then $L_1(D\times\mathbb{Z},w)=0$ and so $S_{TOP}(M)$ has order at most 32.

If $M$ is geometric then $Aut(\pi)$ acts isometrically on $M$.
The natural map from $\pi_0(E(M))$ to $\pi_0(E(P_3(M)))$
has kernel of order at most 2, 
and hence the subgroup which induces the identity on all homotopy groups
is finitely generated, by Corollary 2.9 of \cite{[Ru92]}.
In particular, if $\pi=D*_{\mathbb{Z}}D$, 
${\mathbb{Z}*_{\mathbb{Z}}D}$ or
$D*_{\mathbb{Z}}(\mathbb{Z}\oplus(Z/2Z))$ this subgroup is finite.
Thus $\pi_0(E(M))$ acts on $S_{TOP}(M)$ through a finite group,
and so there are infinitely many manifolds 
in the homotopy type of each such geometric manifold.
   
\section{Minimal models}

Let $X$ be a $PD_4$-complex with fundamental group $\pi$.
A $PD_4$-complex $Z$ is a {\it model} for $X$ 
if there is a 2-connected degree-1 map $f:X\to Z$.
It is {\it strongly minimal} if $\lambda_X=0$.
A strongly minimal $PD_4$-complex $Z$ is minimal with respect 
to the partial order given by $X\geq Y$ 
if $Y$ is a model for $X$.

We shall show that every $PD_4$-complex with fundamental group 
a $PD_2$-group $\pi$ has a strongly minimal model which is the 
total space of an $S^2$-bundle over the surface $F\simeq{K}(\pi,1)$.
(More generally, a $PD_4$-complex $X$ has a strongly minimal model if and only
if $\lambda_X$ is nonsingular and $\mathrm{Cok}(H^2(c_X;\mathbb{Z}[\pi]))$
is a finitely generated projective $\mathbb{Z}[\pi]$-module 
\cite{[Hi20]}.)

The group $Z/2Z$ acts on $CP^\infty$ via complex conjugation,
and so a homomorphism $u:\pi\to Z/2Z$ determines a
product action of $\pi$ on $\widetilde F\times CP^\infty$.
Let $L=L_\pi(\mathbb{Z}^u,2)=(\widetilde F\times CP^\infty)/\pi$ be the 
quotient space.
Projection on the first factor induces a map $q_u=c_L:L\to F$.
In all cases the fixed point set of the action of $u$ on $CP^\infty$
is connected and contains $RP^\infty$. 
Thus $q_u$ has cross-sections $\sigma$, and any two are isotopic.
Let $j:CP^\infty\to L$ be the inclusion of the fibre over the
basepoint of $F$.
If $u$ is trivial $L_\pi(\mathbb{Z}^u,2)\cong F\times CP^\infty$.

The (co)homology of $L$ with coefficients in a $\mathbb{Z}[\pi]$-module 
is split by the homomorphisms induced by $q_u$ and $\sigma$. 
In particular, $H^2(L;\mathbb{Z}^u)\cong 
H^2(F;\mathbb{Z}^u)\oplus H^2(CP^\infty;\mathbb{Z})$,
with the projections to the summands induced by $\sigma$ and $j$.
Let $\omega_F$ be a generator of $H^2(F;\mathbb{Z}^{w_1(F)})\cong\mathbb{Z}$,
and let $\phi$ be a generator of $H^2(F;\mathbb{Z}^u)$.
If $u=w_1(F)$ choose $\phi=\omega_F$; 
otherwise $H^2(F;\mathbb{Z}^u)$ has order 2.
Let $[c]_2$ denote the reduction {\it mod} (2) of a cohomology class $c$
(with coefficients $\mathbb{Z}^{w_1(F)}$ or $\mathbb{Z}^u$).
Then $[\phi]_2=[\omega_F]_2$ in $H^2(F;\mathbb{Z}/2\mathbb{Z})$.
Let $\iota_u\in H^2(L;\mathbb{Z}^u)$ generate the complementary 
$\mathbb{Z}$ summand.
Then $(\iota_u)^2$ and $\iota_u\cup q_u^*\phi$ generate $H^4(L;\mathbb{Z})$.

The space $L=L_\pi(\mathbb{Z}^u,2)$ is a {\it generalized Eilenberg-Mac Lane 
complex of type $(\mathbb{Z}^u,2)$ over $K(\pi,1)$}, 
with {\it characteristic element} $\iota_u$.
Homotopy classes of maps from spaces $X$ into $L$ compatible with 
a fixed homomorphism
$\theta:\pi_1(X)\to\pi$ correspond bijectively to elements of
$H^2(X;\mathbb{Z}^{u\theta})$, via the correspondance 
$f\leftrightarrow{f}^*\iota_u$.
(See Chapter III.\S6 of \cite{[Ba']}).
In particular, $E_\pi(L)$ is the subset of $H^2(L;\mathbb{Z}^u)$ 
consisting of elements of the form $\pm(\iota_u+k\phi)$, 
for $k\in\mathbb{Z}$.
(Such classes restrict to generators of 
$H^2(\widetilde L;\mathbb{Z})\cong\mathbb{Z}$).
As a group $E_\pi(L)\cong{H^2}(\pi;\mathbb{Z}^u)\rtimes\{\pm1\}$.

Let $p:E\to F$ be an $S^2$-bundle over $F$.
Then $\widetilde{E}\cong\widetilde{F}\times{S^2}$ and
$p$ may be identified with the classifying map $c_E$.
If $E_u$ is the image of $\widetilde F\times CP^1$ in $L$ 
and $p_u=q_u|_{E_u}$ then $p_u$ is an $S^2$-bundle over $F$ with $w_1(p_u)=u$, 
and $w_2(p_u)=v_2(E_u)=0$, 
since cross-sections determined by distinct
fixed points are isotopic and disjoint.
(From the dual point of view,
the 4-skeleton of $L$ is $E_u\cup_{CP^1} j(CP^2)=E_u\cup_\eta D^4$,
where $\eta\in\pi_3(E_u)\cong\pi_3(S^2)$ is the Hopf map). 

\begin{theorem} 
Let $E$ be the total space of an $S^2$-bundle over 
an aspherical closed surface $F$, 
and let $X$ be a $PD_4$-complex with $\pi_1(X)\cong\pi=\pi_1(F)$.
Then there is a $2$-connected degree-$1$ map $h:X\to E$ such that $c_E=c_Xh$
if and only if $(c_X^*)^{-1}w_1(X)=(c_E^*)^{-1}w_1(E)$ and 
$\xi\cup{c_X^*H^2(\pi;\mathbb{F}_2)}\not=0$
for some $\xi\in H^2(X;\mathbb{F}_2)$
such that $\xi^2=0$ if $v_2(E)=0$ and $\xi^2\not=0$ if $v_2(E)\not=0$.
\end{theorem}

\begin{proof}
Compatibility of the orientation characters is clearly necessary in order that
the degree be defined as an integer; we assume this henceforth.
Since $c_X$ is 2-connected there is an $\alpha\in H_2(X;\mathbb{Z}^{c_X^*w_1(F)})$
such that $c_{X*}\alpha=[F]$, 
and since $\mathbb{Z}^u\otimes\mathbb{Z}^{c_X^*w_1(F)}\cong\mathbb{Z}^{w_1(X)}$
there is an $x\in H^2(X;\mathbb{Z}^u)$ such that $x\cap[X]=\alpha$, 
by Poincar\'e duality.
Hence $(x\cup c_X^*\omega_F)\cap[X]=1$.
Clearly either $[x]_2^2=0$ or $[x]_2^2=[x]_2\cup c_X^*[\omega_F]_2$.

The map $f:E\to L=L_\pi(\mathbb{Z}^u,2)$ corresponding to a class 
$f^*\iota_u\in H^2(E;\mathbb{Z}^u)$ which restricts to a generator for 
$H^2(S^2;\mathbb{Z})$ induces isomorphisms on $\pi_1$ and $\pi_2$,
and so $f=f_E$.
(We may vary this map by composition with a self homotopy equivalence
of $L$, replacing $f^*\iota_u$ by $f^*\iota_u+kf^*\phi$).
Note also that $f_E^*\iota_u\cup c_E^*\omega_F$ generates
$H^4(E;\mathbb{Z}^{w_1(X)})\cong\mathbb{Z}$.

The action of $\pi$ on $\pi_3(E)\cong\pi_3(S^2)=\Gamma_W(\mathbb{Z})$
is given by $\Gamma_W(u)$, and so is trivial.
Therefore the third stage of the Postnikov tower for $E$ is a simple
$K(\mathbb{Z},3)$-fibration over $L$, 
determined by a map $\kappa:L\to K(\mathbb{Z},4)$
corresponding to a class in $H^4(L;\mathbb{Z})$.
If $L(m)$ is the space induced by
$\kappa_m=\iota_u^2+m\iota_u\cup\phi$ then 
$\widetilde{L(m)}$ is induced from $\widetilde L\simeq CP^\infty$ by
the canonical generator of $H^4(CP^\infty)$, and so
$H_3(\widetilde{L(m)};\mathbb{Z})=H_4(\widetilde{L(m)};\mathbb{Z})=0$,
by a spectral sequence argument. 
Hence $\Gamma_W(\mathbb{Z})\cong\pi_3(L(m))=\mathbb{Z}$.

The map $f_E$ factors through a map $g_E:E\to L(m)$
if and only if $f_E^*\kappa_m=0$.
We then have $\pi_3(g_E)=\Gamma_W(f_E)$, which is an isomorphism.
Thus $g_E$ is 4-connected, and so is the third stage of the Postnikov 
tower for $E$.
If $v_2(E)=0$ then $f_E^*\iota_u^2=2kf_E^*(\iota_u\cup\phi)$ for some
$k\in\mathbb{Z}$, and so $f_E$ factors through $L(-2k)$;
otherwise $f_E$ factors through $L(-2k-1)$,
and thus these spaces provide models for the third stages $P_3(E)$
of such $S^2$-bundle spaces.
The self homotopy equivalence of $L$ corresponding to the class
$\pm(\iota_u+k\phi)$ in  $H^2(L;\mathbb{Z}^u)$ carries 
$\kappa_m=\iota_u^2+m\iota_u\cup\phi$ 
to $\kappa_{m\pm2k}$,
and thus $c_{L(m)}$ is fibre homotopy equivalent to $c_{L(0)}$ if $m$ is even 
and to $c_{L(1)}$ otherwise.

Since $P_3(E)$ may also be obtained from $E$ 
by adjoining cells of dimension $\geq5$, 
maps from a complex $X$ of dimension at most 4 to $E$ 
compatible with $\theta:\pi_1(X)\to\pi$ correspond to maps 
from $X$ to $P_3(E)$ compatible with $\theta$ and thus to elements 
$y\in H^2(X;\mathbb{Z}^{u\theta})$ 
such that $[y]_2^2=0$ if $v_2(E)=0$ 
and $[y]_2^2=[y]_2\cup c_X^*[\omega_F]_2$ otherwise.
(In the next paragraph we omit $\theta=\pi_1(c_X)$ from the notation.)

If $g:X\to E$ is a 2-connected degree-1 map then
$\xi=g^*f_E^*[\iota_u]_2$ satisfies the conditions of the theorem,
since $c_X\sim c_Eg$, which factors through $P_3(E)$.
Conversely, let $\xi$ be such a class.
Reduction {\it mod} (2) maps 
$H^2(X;\mathbb{Z}^u)$ onto $H^2(X;\mathbb{F}_2)$,
since $H^3(X;\mathbb{Z}^u)\cong
{H_1(X;\mathbb{Z}^{w_1(F)})}\cong{H^1(\pi;\mathbb{Z})}$
is torsion free.
Therefore there is an 
$x\in{H^2(X;\mathbb{Z}^u)}$ such that $[x]_2=\xi$.
Since $\omega_F$ generates a direct summand of $H^2(X;\mathbb{Z}^{w_1(F)})$,
and $(x\cup\omega_F)[X]$ is odd, 
we may choose $x$ so that $(x\cup\omega_F)[X]=1$.
Then $x=h^*f_E^*\iota_u$ for some $h:X\to E$ such that $c_Eh=c_X$ and
$(f_E^*\iota_u\cup c_E^*\omega_F)h_*[X]=(x\cup c_X^*\omega_F)[X]=1$.
Thus $\pi_1(h)$ is an isomorphism and $h$ is a degree-1 map, 
and so $h$ is 2-connected \cite[Lemma 2.2]{[Wl]}.
\end{proof}

We shall summarize related work on the homotopy types of
$PD_4$-complexes.

\begin{thm} 
{\rm[Hi20]}
There is a strongly minimal model for $X$ if 
$H^3(\pi;\mathbb{Z}[\pi])=0$ and
$Hom_{\mathbb{Z}[\pi]}(\pi_2(X),\mathbb{Z}[\pi])$ is a finitely generated 
projective $\mathbb{Z}[\pi]$-module. 
If $Z$ is strongly minimal $\pi_2(Z)\cong\overline{H^2(\pi;\mathbb{Z}[\pi])}$
and $\lambda_Z=0$, and if $v_2(X)=0$ and $Z$ is a model for $X$ then $v_2(Z)=0$.
\qed
\end{thm}

These conditions hold if $c.d.\pi\leq2$, 
and then $\chi(Z)=q(\pi)$, by Theorem 3.12.
The strongly minimal $PD_4$-complexes with $\pi$ free,
$F(r)\rtimes\mathbb{Z}$ or a $PD_2$-group are given by Theorems 14.9,
4.5 and 5.10, respectively.
(See also Theorem 10.17.)
If $v_2(\widetilde{X})\not=0$ the minimal model may not be unique.
For example, if $C$ is a compact complex curve of genus $\geq1$
the ruled surface $C\times\mathbb{CP}^1$ is strongly minimal,
but the blowup $(C\times\mathbb{CP}^1)\sharp\overline{\mathbb{CP}}^2$ also has
the nontrivial bundle space as a strongly minimal model.
(Many of the other minimal complex surfaces in the Enriques-Kodaira 
classification are aspherical, and hence strongly minimal in our sense.
However 1-connected complex surfaces are never strongly minimal,
since the unique minimal 1-connected $PD_4$-complex is $S^4$, which
has no almost complex structure,  
by the theorem of Wu cited on page 149 above.)

\begin{thm}
{\rm[Hi20]}
Let $\pi$ be a finitely presentable group with $c.d.\pi\leq2$.
Two $PD_4$-complexes $X$ and $Y$ with fundamental group $\pi$,
$w_1(X)=w_1(Y)=w$ and $\pi_2(X)\cong\pi_2(Y)$ are homotopy equivalent 
if and only if $X$ and $Y$ have a common strongly minimal model $Z$ and
$\lambda_X\cong\lambda_Y$.
Moreover $\lambda_X$ is nonsingular and every nonsingular $w$-hermitean pairing 
on a finitely generated projective $\mathbb{Z}[\pi]$-module is the 
reduced intersection pairing of some such $PD_4$-complex.
\qed
\end{thm}

In particular, a $Spin$ 4-manifold with fundamental group a $PD_2$-group $\pi$ 
has a well-defined strongly minimal model and so two such $Spin$ 4-manifolds 
$X$ and $Y$ are homotopy equivalent if and only if $\lambda_X\cong\lambda_Y$.

%% file: m5-11.tex
\chapter{Manifolds covered by $S^3 \times\mathbb{R}$}

In this chapter we shall show that a closed 4-manifold $M$
is covered by $S^3\times\mathbb{R}$ if and only if $\pi=\pi_1(M)$ 
has two ends and $\chi(M)=0$.
Its homotopy type is then determined by $\pi$ and the first 
nonzero $k$-invariant $k(M)$.
The maximal finite normal subgroup of $\pi$ is either
the group of a $\mathbb{S}^3 $-manifold or one of the groups 
$P_{48r}''$ with $r$ odd or
$Q(8a,b,c)\times Z/dZ$ with $a,b,c$ and $d$ odd.
(There are examples of the latter type, 
and no such $M$ is homotopy equivalent 
to a $\mathbb{S}^3\times\mathbb{E}^1 $-manifold.)
The possibilities for $\pi$ are not yet known
even when $F$ is a $\mathbb{S}^3$-manifold group and $\pi/F\cong\mathbb{Z}$.
Solving this problem may involve first determining which 
$k$-invariants are realizable when $F$ is cyclic.
 
Manifolds which fibre over $RP^2$ with fibre $T$ or $Kb$ and $\partial\not=0$ 
have universal cover $S^3\times\mathbb{R}$. 
In \S6 we determine the possible fundamental groups, 
and show that an orientable 4-manifold $M$
with such a group and with $\chi(M)=0$ must be homotopy equivalent to
a $\mathbb{S}^3\times\mathbb{E}^1$-manifold which fibres over $RP^2$.

As groups with two ends are virtually solvable, surgery techniques may 
be used to study manifolds covered by $S^3\times\mathbb{R}$.
However computing $Wh(\pi)$ and $L_*(\pi; w_1)$ is a major task.
Simple estimates suggest that there are usually infinitely many 
nonhomeomorphic manifolds within a given homotopy type.

\section{Invariants for the homotopy type}

The determination of the closed 4-manifolds with universal covering space 
homotopy equivalent to $S^3$ is based on the structure of groups with two ends.

\begin{theorem}
Let $M$ be a closed $4$-manifold with 
fundamental group $\pi$.
Then $\widetilde M\simeq S^3 $ if and only if $\pi$ has two ends and 
$\chi(M)=0$. If so
\begin{enumerate}                                             
\item $M$ is finitely covered by $S^3 \times S^1 $ and so 
$\widetilde M\cong S^3 \times\mathbb{R}\cong\mathbb{R}^4\setminus\{0\} $;

\item the maximal finite normal subgroup $F$ of $\pi$ has 
cohomological period dividing $4$, 
acts trivially on $\pi_3 (M)\cong\mathbb{Z}$ and the corresponding covering 
space $M_F$ has the homotopy type of an orientable finite $PD_3 $-complex;

\item if $v:\pi\to{Aut}(H^1(\pi;\mathbb{Z}[\pi]))$ is the natural action and 
$w=w_1(M)$ then the action $u:\pi\to{Aut}(\pi_3(M))$ is given by $u=v+w$;

\item the homotopy type of $M$ is determined by $\pi$, $w_1(M)$ and the 
orbit of the first nontrivial $k$-invariant 
$k(M)\in H^4(\pi;\mathbb{Z}^u)$ under $Out(\pi)\times\{\pm1\}$;

\item the restriction of $k(M)$ to $H^4(F;\mathbb{Z})$ is a generator;

\item if $\pi/F\cong\mathbb{Z}$ then 
$H^4 (\pi;\pi_3 (M))\cong H^4 (F;\mathbb{Z})\cong Z/|F|Z$.
\end{enumerate}          
\end{theorem}

\begin{proof} 
If $\widetilde M\simeq S^3 $ then $H^1 (\pi;\mathbb{Z}[\pi])$ 
is infinite cyclic and so $\pi$ has two ends.
Hence $\pi$ is virtually $\mathbb{Z}$. 
The covering space $M_A$ corresponding to an infinite cyclic subgroup 
$A$ is homotopy equivalent to the mapping torus of a self 
homotopy equivalence of $S^3 \simeq \widetilde M$, and so $\chi(M_A)=0$.
As $[\pi:A]<\infty$ it follows that $\chi(M)=0$ also.

Suppose conversely that $\chi(M)=0$ and $\pi$ is virtually $\mathbb{Z}$.
Then $H_3 (\widetilde M;\mathbb{Z})\cong\mathbb{Z}$ and 
$H_4 (\widetilde M;\mathbb{Z})=0$.
Let $M_{\mathbb{Z}}$ be an orientable finite covering space 
with fundamental group $\mathbb{Z}$.
Then $\chi(M_{\mathbb{Z}})=0$ and so $H_2 (M_{\mathbb{Z}};\mathbb{Z})=0$. 
The homology groups of $\widetilde M=\widetilde M_{\mathbb{Z}}$ 
may be regarded as modules over 
$\mathbb{Z}[t,t^{-1} ]\cong\mathbb{Z}[\mathbb{Z}]$. 
Multiplication by $t-1$ maps $H_2 (\widetilde M;\mathbb{Z})$ onto itself,
by the Wang sequence for the projection of $\widetilde M$ 
onto $M_{\mathbb{Z}}$.
Therefore 
$Hom_{\mathbb{Z}[\mathbb{Z}]}(H_2(\widetilde M;\mathbb{Z}),
\mathbb{Z}[\mathbb{Z}])=0$
and so $\pi_2(M)=\pi_2(M_{\mathbb{Z}})=0$, by Lemma 3.3.
Therefore the map from $S^3 $ to $\widetilde M$ representing a generator of 
$\pi_3 (M)$ is a homotopy equivalence. 
Since $M_{\mathbb{Z}}$ is orientable the generator of the group of covering 
translations $Aut(\widetilde M/M_{\mathbb{Z}})\cong\mathbb{Z}$ 
is homotopic to the identity,
and so $M_{\mathbb{Z}}\simeq \widetilde M\times S^1\simeq S^3 \times S^1 $.
Therefore $M_{\mathbb{Z}}\cong S^3 \times S^1 $, 
by surgery over $\mathbb{Z}$.
Hence $\widetilde M\cong S^3 \times\mathbb{R}$.

Let $F$ be the maximal finite normal subgroup of $\pi$.
Since $F$ acts freely on $\widetilde M\simeq S^3 $ it has cohomological
period dividing 4 and $M_F =\widetilde M/F$ is a $PD_3$-complex.
In particular, $M_F$ is orientable and $F$ acts trivially on $\pi_3(M)$. 
The image of the finiteness obstruction for $M_F $ under the 
``geometrically significant injection" of $K_0(\mathbb{Z}[F])$ 
into $Wh(F\times\mathbb{Z})$ of \cite{[Rn86]} is the obstruction 
to $M_F \times S^1 $ being a simple $PD$-complex.
If $f:M_F\to M_F$ is a self homotopy equivalence which induces 
the identity on $\pi_1(M_F)\cong F$ and on $\pi_3(M_F)\cong\mathbb{Z}$ 
then $f$ is homotopic to the identity, by obstruction theory \cite{[Pl82]}.
Therefore $\pi_0(E(M_F))$ is finite and so $M$ has a finite cover which is 
homotopy equivalent to $M_F\times S^1$. 
Since manifolds are simple $PD_n$-complexes $M_F$ must be finite. 

The third assertion follows from the Hurewicz Theorem
and Poincar\'e duality, as in Lemma 10.3.
The first nonzero $k$-invariant lies in $H^4(\pi;\mathbb{Z}^u)$,
since $\pi_2(M)=0$ and $\pi_3 (M)\cong\mathbb{Z}^u$, 
and it restricts to the $k$-invariant for $M_F$ in $H^4(F;\mathbb{Z})$.
Thus (4) and (5) follow as in Theorem 2.9.
The final assertion follows from the LHSSS (or Wang sequence) 
for $\pi$ as an extension of $\mathbb{Z}$ by $F$,
since $\pi/F$ acts trivially on $H^4(F;\mathbb{Z}^u)$.
\end{proof}
 
The list of finite groups with cohomological period dividing 4 
is well known (see \cite{[Mi57],[AM]}).
There are the generalized quaternionic groups $Q(2^na,b,c)$ 
(with $n\geq3$ and $a,b,c$ odd),
the extended binary tetrahedral groups $T_k^*$,
the extended binary octahedral groups $O_k^*$,
the groups $P''_{48r}$ (with $r$ odd $>1$),
the binary icosahedral group $I^*$,
the dihedral groups $A(m,e)$ (with $m$ odd $>1$),
and the direct products of
any one of these with a cyclic group $Z/dZ$ of relatively prime order.
(In particular, a $p$-group with periodic cohomology
is cyclic if $p$ is odd and cyclic or quaternionic if $p=2$.)
We shall give presentations for these groups in \S2.

Each such group $F$ is the fundamental group of some $PD_3$-complex 
\cite{[Sw60]}. 
Such {\it Swan complexes} for $F$ are orientable, 
and are determined up to homotopy equivalence by their $k$-invariants, 
which are generators of $H^4 (F;\mathbb{Z})\cong Z/|F|Z$,
by Theorem 2.9.
Thus they are parametrized up to homotopy by the quotient of 
$(Z/|F|Z)^\times $ under the action of $Out(F)\times\{\pm1\}$.
The set of finiteness obstructions for all such complexes forms a coset of the 
``Swan subgroup" of $\tilde K_0 (\mathbb{Z}[F])$ and there is a finite complex 
of this type if and only if the coset contains 0.
(This condition fails if $F$ has a subgroup isomorphic to $Q(16,3,1)$ 
and hence if $F\cong O^*_k \times(Z/dZ)$ for some $k>1$ 
or $P''_{48r}$ with $3|r$ \cite[Corollary 3.16]{[DM85]}.)
If $X$ is a Swan complex for $F$ then $X\times S^1$ is a finite 
$PD_4^+$-complex with $\pi_1(X\times S^1)\cong F\times\mathbb{Z}$ 
and $\chi(X\times S^1)=0$.

\begin{theorem}
Let $M$ be a closed $4$-manifold 
such that $\pi=\pi_1(M)$ has two ends and with $\chi(M)=0$.
Then the group of unbased homotopy classes of self homotopy equivalences 
of $M$ is finite.
\end{theorem}
                     
\begin{proof}
We may assume that $M$ has a finite cell structure with a single 4-cell.
Suppose that $f:M\to M$ is a self homotopy equivalence which fixes a base 
point and induces the identity on $\pi$ and on $\pi_3(M)\cong\mathbb{Z}$.
Then there are no obstructions to constructing a homotopy from $f$ to 
$id_{\widetilde M}$ on the 3-skeleton $M_0=M\backslash int D^4$, 
and since $\pi_4(M)=\pi_4(S^3)=Z/2Z$ there are just two possibilities for $f$.
It is easily seen that $Out(\pi)$ is finite.
Since every self map is homotopic to one which fixes a basepoint the group 
of unbased homotopy classes of self homotopy equivalences of $M$ is finite.
\end{proof}
 
If $\pi$ is a semidirect product $F\rtimes_\theta\mathbb{Z}$ 
then $Aut(\pi)$ is finite and the group of {\it based} homotopy classes 
of based self homotopy equivalences is also finite. 

\section{The action of $\pi/F$ on $F$}

Let $F$ be a finite group with cohomological period dividing 4.
Automorphisms of $F$ act on $H_*(F;\mathbb{Z})$ and $H^*(F;\mathbb{Z})$
through $Out(F)$,
since inner automorphisms induce the identity on (co)homology.
Let $J_+(F)$ be the kernel of the action on $H_3(F;\mathbb{Z})$,
and let $J(F)$ be the subgroup of $Out(F)$ which acts by $\pm1$.

An outer automorphism class induces a well defined action on 
$H^4(S;\mathbb{Z})$ for each Sylow subgroup $S$ of $F$, 
since all $p$-Sylow subgroups are conjugate in $F$ and the inclusion
of such a subgroup induces an isomorphism from the $p$-torsion of 
$H^4(F;\mathbb{Z})\cong Z/|F|Z$ to $H^4(S;\mathbb{Z})\cong Z/|S|Z$, 
by Shapiro's Lemma.
Therefore an outer automorphism class of $F$ induces multiplication by $r$ 
on $H^4 (F;\mathbb{Z})$ if and only if it does so for each 
Sylow subgroup of $F$, by the Chinese Remainder Theorem. 

The map sending a self homotopy equivalence $h$ of a Swan complex $X_F$ 
for $F$ to the induced outer automorphism class determines 
a homomorphism from the group of (unbased) homotopy classes 
of self homotopy equivalences $E(X_F)$ to $Out(F)$.
The image of this homomorphism is $J(F)$, 
and it is a monomorphism if $|F|>2$ \cite[Corollary 1.3]{[Pl82]}. 
(Note that \cite{[Pl82]} works with {\it based} homotopies.)
If $F=1$ or $Z/2Z$ the orientation reversing involution of $X_F$ 
($\simeq S^3$ or $RP^3$, respectively) induces the identity on $F$.

\begin{lemma}          
Let $M$ be a closed $4$-manifold with universal
cover $S^3\times\mathbb{R}$, 
and let $F$ be the maximal finite normal subgroup of $\pi=\pi_1(M)$.
The quotient $\pi/F$ acts on $\pi_3 (M)$ and 
$H^4 (F;\mathbb{Z})$ through multiplication by $\pm 1$. 
It acts trivially if the order of $F$ is divisible by $4$ or by 
any prime congruent to $3$ {\it mod} $(4)$.
\end{lemma}

\begin{proof} 
The group $\pi/F$ must act through $\pm 1$ on the 
infinite cyclic groups $\pi_3 (M)$ and $H_3 (M_F;\mathbb{Z})$. 
By the universal coefficient theorem $H^4 (F;\mathbb{Z})$ 
is isomorphic to $H_3 (F;\mathbb{Z})$, which is the cokernel of the Hurewicz 
homomorphism from $\pi_3 (M)$ to $H_3 (M_F;\mathbb{Z})$.
This implies the first assertion.

To prove the second assertion we may pass to the Sylow subgroups of $F$, 
by Shapiro's Lemma. 
Since the $p$-Sylow subgroups of $F$ also have cohomological period 4
they are cyclic if $p$ is an odd prime and are cyclic or quaternionic 
($Q(2^n)$) if $p=2$. 
In all cases an automorphism induces multiplication 
by a square on the third homology \cite{[Sw60]}. 
But $-1$ is not a square modulo 4 nor modulo any prime $p=4n+3$. 
\end{proof}

Thus the groups $\pi\cong F\rtimes\mathbb{Z}$ realized by such 4-manifolds
correspond to outer automorphisms in $J(F)$ or $J_+(F)$.
We shall next determine these subgroups of $Out(F)$
for $F$ a group of cohomological period dividing 4.
If $m$ is an integer let $l(m)$ be the number of odd prime divisors of $m$.
\[
Z/dZ=\langle x\mid x^d=1\rangle.
\]
$Out(Z/dZ)=(Z/dZ)^\times $. 
Hence $J(Z/dZ)=\{ s\in (Z/dZ)^\times \mid s^2 =\pm 1\}$.
$J_+ (Z/dZ)=(Z/2Z)^{l(d)}$ if $d\not\equiv 0$ {\it mod\/} $(4)$, 
$(Z/2Z)^{l(d)+1}$ if $d\equiv 4$ {\it mod\/} $(8)$,
and $(Z/2Z)^{l(d)+2}$ if $d\equiv 0$ {\it mod\/} $(8)$.
\[
Q(8)=\langle x,y\mid x^2 =y^2 =(xy)^2\rangle.
\]
An automorphism of $Q=Q(8)$ induces the identity on $Q/Q'$ if and only if it 
is inner, and every automorphism of $Q/Q'$ lifts to one of $Q$. 
In fact $Aut(Q)$ is the semidirect product of 
$Out(Q)\cong Aut(Q/Q')\cong SL(2,\mathbb{F}_2)$ with the normal subgroup 
$Inn(Q)=Q/Q'\cong (Z/2Z)^2$. 
Moreover $J(Q)=Out(Q)$, generated by the images of the automorphisms 
$\sigma$ and $\tau$, where $\sigma$ sends 
$x$ and $y$ to $y$ and $xy$, respectively, and $\tau$ interchanges $x$ and $y$.
\[Q(8k)=\langle x,y\mid x^{2k}=y^2,\medspace
yxy^{-1}=x^{-1}\rangle,\quad\mathrm
{where}~k>1.
\]
(The relations imply that $x^{4k}=y^4=1$.)
All automorphisms of $Q(8k)$ are of the form $[i,s]$,
where $(s,2k)=1$, $[i,s](x)=x^s$ and $[i,s](y)=x^iy$, 
and $Aut(Q(8k))$ is the semidirect product of $(Z/4kZ)^\times $ 
with the normal subgroup $\langle[1,1]\rangle\cong Z/4kZ$. 
$Out(Q(8k))=(Z/2Z)\oplus((Z/4kZ)^\times/(\pm 1))$, 
generated by the images of the $[0,s]$ and [1,1]. 
The automorphism $[i,s]$ induces multiplication by $s^2$ on
$H^4(Q(2^n);\!\mathbb{Z})$ \cite{[Sw60]}.
Hence $J(Q(8k))=(Z/2Z)^{l(k)+1}$ if $k$ is odd and $(Z/2Z)^{l(k)+2}$ 
if $k$ is even.
\[
T_k^* =\langle Q(8),z\mid z^{3^k} =1,\medspace zxz^{-1} =y,\medspace
zyz^{-1} =xy\rangle,
\quad\mathrm
{where}~ k\geq 1.
\]
Setting $t=zx^2$ gives the balanced presentation
$\langle{t,x}\mid{t^2x=xtxt},\medspace{t^{3^k}=x^2}\rangle$.
Let $\rho $ be the automorphism sending $x$, $y$ and $z$ to $y^{-1}$, 
$x^{-1} $ and $z^2 $ respectively. 
Let $\xi$, $\eta$ and $\zeta$ be the inner automorphisms of $T^*_k$
determined by conjugation by $x$, $y$ and $z$, respectively.
An induction on $k$ shows that the image of 2 generates $(Z/3^kZ)^\times$,
and so $\rho$ has order $2.3^{k-1}$.
Then $Aut(T_k^*)$ has the presentation
\[
\langle
\rho,\xi,\eta,\zeta\mid\rho^{2.3^{k-1}}\! =\eta^2=\zeta^3=(\eta\zeta)^3=1,\medspace
\rho\zeta\rho^{-1}=\zeta^2,\medspace 
\rho\eta\rho^{-1}=\zeta^{-1}\eta\zeta=\xi\rangle.
\]
$Out(T_k^*)\cong (Z/3^kZ)^\times$.
The 3-Sylow subgroup generated by $z$ is preserved by $\rho$,
and it follows that $J(T_k^*)=Z/2Z$ 
(generated by the image of $\rho^{3^{k-1}}$).
\[
O_k^* =\langle T_k^*,w\mid w^2 =x^2,\medspace wxw^{-1} =yx,\medspace
wzw^{-1} =z^{-1}\rangle,
\quad\mathrm
{where}~ 
k\geq 1.
\]
(The relations imply that $wyw^{-1} =y^{-1}$.)
Let $\omega,\xi,\eta$ and $\zeta$ be the inner automorphisms 
of $O_k^*$ determined by
conjugation by $w,x,y$ and $z$, respectively.
As we may extend $\rho$ to an automorphism of $O^*_k$ 
via $\rho(w)=w^{-1} z^2$, 
the restriction from $Aut(O_k^* )$ to $Aut(T_k^*)$ is onto.
An automorphism in the kernel sends $w$ to $wv$ for some $v\in T_k^*$,
which must be central in $T^*_k$.
Hence the kernel is generated by $\alpha$ and $\zeta^3$, 
where $\alpha(w)=w^{-1}=wx^2$.
Now $\alpha=\rho^{3^{k-1}}\omega\zeta^{-2}$, 
and so the image of $\rho$ generates $Out(O_k^*)$.
The 2-Sylow subgroup is generated by $u=xw$ and $x$, 
and is isomorphic to $Q(16)$.
As $\alpha(u)=u^3$ and $\alpha(x)=x$ this subgroup is preserved by $\alpha$,
and $H^4(\alpha|_{\langle u,x\rangle};\mathbb{Z})$ is multiplication by 9.
Hence $H_4(\rho^{3^{k-1}};\mathbb{Z})$ is multiplication by 9 on the 2-primary torsion.
As the order of $Out(O_k^*)$ divides $2.3^{k-1}$ it follows that $J(O_k^* )=1$.
\[
Z/aZ\rtimes_{-1}O^*_k=\langle{O^*_k,u}\mid{u^a=1},~wuw^{-1}=u^{-1},~
x,y,z\leftrightharpoons{u}\rangle
\]
Here $(a,6)=1$.
This is the group $P''_{48.3^{k-1}a}$ of \cite{[Mi57]}.
An automorphism of this group induces automorphisms of 
the normal subgroup $Z/aZ$ and of the quotient $O^*_k$.
It is easily seen that
$Aut(P''_{48.3^{k-1}a})\cong{Aut(Z/aZ)\times{Aut}(O^*_k)}$.
Since conjugation by $w$ induces the inversion on $Z/aZ$,
$Out(P''_{48.3^{k-1}a})\cong{Aut(Z/aZ)\rtimes_{-1}{Out}(O^*_k)}$.
Since $J(O^*_k)=1$, we have $J(P''_{48.3^{k-1}a})\cong{J_+(Z/aZ)}$.
\[
I^* =\langle{y,z}\mid(yz)^2 =y^3 =z^5\rangle.
\]
The map sending the generators $y,z$ to
$\left(\smallmatrix 2&2\\
1&4\endsmallmatrix\right)$ and
$\left(\smallmatrix 4&1\\
0&4\endsmallmatrix\right)$, respectively,
induces an isomorphism from $I^*$ to $SL(2,\mathbb{F}_5)$.
Conjugation in $GL(2,\mathbb{F}_5)$ induces a monomorphism from 
$PGL(2,\mathbb{F}_5)$ to $Aut(I^*)$. 
The natural map from $Aut(I^*)$ to $Aut(I^*/\zeta I^*)$
is injective, since $I^*$ is perfect.
Now $I^*/\zeta I^*\cong PSL(2,\mathbb{F}_5)\cong A_5$. 
The alternating group $A_5$ is generated by 3-cycles, 
and has ten 3-Sylow subgroups, each of order 3.
It has five subgroups isomorphic to $A_4$ generated by pairs of such 3-Sylow 
subgroups.
The intersection of any two of them has order 3, and is invariant under
any automorphism of $A_5$ which leaves invariant each of these subgroups.
It is not hard to see that such an automorphism must fix the 3-cycles.
Thus $Aut(A_5)$ embeds in the group $S_5$ of permutations of these subgroups.
Since $|PGL(2,\mathbb{F}_5)|=|S_5|=120$ it follows that 
$Aut(I^*)\cong S_5$ and $Out(I^* )=Z/2Z$.
The outer automorphism class is represented by the matrix
$\omega=\left(\smallmatrix 2&0\\
0&1\endsmallmatrix\right)$ in $GL(2,\mathbb{F}_5)$.

\begin{lemma} {\rm[Pl83]}\qua $J(I^*)=1$.
\end{lemma}

\begin{proof} 
The element $\gamma=z^4=\left(\smallmatrix 1&1\\
0&1\endsmallmatrix\right)$
generates a 5-Sylow subgroup of $I^*$.
It is easily seen that $\omega\gamma\omega^{-1}=\gamma^2$,
and so $\omega$ induces multiplication by 2 on 
$H^2(Z/5Z;\mathbb{Z})\cong H_1(Z/5Z;\mathbb{Z})=Z/5Z$.
Since $H^4(Z/5Z;\mathbb{Z})\cong Z/5Z$ is generated by 
the square of a generator for $H^2(Z/5Z;\mathbb{Z})$
we see that $H^4(\omega;\mathbb{Z})$ is multiplication by $4=-1$ on 5-torsion.
Hence $J(I^*)=1$.
\end{proof}
 
In fact $H^4(\omega;\mathbb{Z})$ is multiplication by 49 
\cite{[Pl83]}.

\smallskip
\noindent $A(m,e)=\langle x,y\mid x^m =y^{2^e} =1,\medspace 
yxy^{-1} =x^{-1}\rangle$, 
where $e\geq1$ and $m>1$ is odd. 

\smallskip
If $m=2n+1$ then $\langle x,y\mid x^m =y^{2^e},\medspace 
yx^ny^{-1} =x^{n+1}\rangle$ is a balanced presentation.
All automorphisms of $A(m,e)$ are of the form $[s,t,u]$,
where $(s,m)=(t,2)=1$, $[s,t,u](x)=x^s$ and $[s,t,u](y)=x^uy^t$.
$Out(A(m,e))$ is generated by the images of $[s,1,0]$ and $[1,t,0]$ 
and is isomorphic to $(Z/2^e )^\times \oplus ((Z/mZ)^\times/(\pm 1)) $. 
Hence $J(A(m,1))=\{ s\in (Z/mZ)^\times \mid s^2 =\pm 1\}/(\pm 1)$,
$J(A(m,2))=(Z/2Z)^{l(m)}$, and $J(A(m,e))=(Z/2Z)^{l(m)+1}$ if $e>2$.

\smallskip                                        
\noindent $Q(2^n a,b,c)=\langle Q(2^n ),u\mid u^{abc} =1,\medspace
xu^{ab} =u^{ab} x,\medspace xu^cx^{-1} =u^{-c},\medspace
yu^{ac} =u^{ac} y,\medspace$ 
$yu^b y^{-1} =u^{-b}\rangle$, where $a$, $b$ and $c$ are odd and relatively prime, 
and either $n=3$ and at most one of $a$, $b$ and $c$ is 1 or $n>3$ and $bc>1$. 

\smallskip
An automorphism of $G=Q(2^n a,b,c)$ must induce the identity on $G/G'$. 
If it induces the identity on the characteristic subgroup 
$\langle u\rangle\cong Z/abcZ$ and on $G/\langle u\rangle\cong Q(2^n)$ 
it is inner, and so $Out(Q(2^n a,b,c))$ is a subquotient of 
$Out(Q(2^n ))\times(Z/abcZ)^\times $.
In particular, $Out(Q(8a,b,c))\cong (Z/abcZ)^\times$,
and $J(Q(8a,b,c))\cong (Z/2Z)^{l(abc)}$.
(We need only consider $n=3$, by \S5 below.)
                                
\smallskip
\noindent 
As $Aut(G\times H)=Aut(G)\times Aut(H)$ and 
$Out(G\times H)=Out(G)\times Out(H)$
if $G$ and $H$ are finite groups of relatively prime order, 
we have $J_+ (G\times Z/dZ)=J_+ (G)\times J_+ (Z/dZ)$.
In particular, if $G$ is not cyclic or dihedral  
$J(G\times Z/dZ)=J_+ (G\times Z/dZ)=J(G)\times J_+ (Z/dZ)$. 
In all cases except when $F$ is cyclic or $Q(8)\times Z/dZ$ the group $J(F)$ 
has exponent 2 and hence $\pi $ has a subgroup of index at most 4 which is 
isomorphic to $F\times\mathbb{Z}$. 

\section{Extensions of $D$}

We shall now assume that $\pi/F\cong D$, 
and so $\pi\cong G*_F H$, where $[G:F]=[H:F]=2$.
Let $u,v\in D$ be a pair of involutions which generate $D$ and
let $s=uv$. Then $s^{-n} us^n =us^{2n} $, 
and any involution in $D$ is conjugate to $u$ or to $v=us$.
Hence any pair of involutions $\{ u',v'\}$ which generates $D$ 
is conjugate to the pair 
$\{ u,v\}$, up to change of order.

\begin{theorem}
Let $M$ be a closed $4$-manifold 
with $\chi(M)=0$, and such that there is an epimorphism 
$p:\pi=\pi_1(M)\to D$ with finite kernel $F$.
Let $\hat u$ and $\hat v$ be a pair of elements of $\pi$ whose images 
$u=p(\hat u)$ and $v=p(\hat v)$ in $D$ are involutions which together generate $D$. Then
\begin{enumerate}          
\item $M$ is nonorientable and $\hat u,\hat v$ each represent orientation reversing loops;

\item the subgroups $G$ and $H$ generated by $F$ and $\hat u$ and by $F$ and 
$\hat v$, respectively, 
each have cohomological period dividing $4$, and the unordered pair $\{ G,H\}$ 
of groups is determined up to isomorphisms by $\pi $ alone;

\item conversely, $\pi $ is determined up to isomorphism by the unordered pair 
$\{ G,H\}$ of groups with index $2$ subgroups isomorphic to $F$ as the free 
product with amalgamation $\pi =G*_F H$;

\item $\pi$ acts trivially on $\pi_3 (M)$;

\item the restrictions of $k(M)$ generate the groups 
$H^4 (G;\mathbb{Z})$ and $H^4 (H;\mathbb{Z})$,
and $H^4(\pi;\mathbb{Z})\cong\{(\zeta,\xi)\in(Z/|G|Z)\oplus(Z/|H|Z)\mid
\zeta\equiv\xi~\mathrm{mod}~(|F|)\}\cong (Z/2|F|Z)\oplus( Z/2Z)$.
\end{enumerate}          
\end{theorem}

\begin{proof}                           
Let $\hat s=\hat u\hat v$.
Suppose that $\hat u$ is orientation preserving. 
Then the subgroup $\sigma$ generated by $\hat u$ and $\hat s^2$ 
is orientation preserving,
so the corresponding covering space $M_\sigma$ is orientable.
As $\sigma$ has finite index in $\pi$ and $\sigma/\sigma'$ 
is finite this contradicts Lemma 3.14.
Similarly, $\hat v$ must be orientation reversing.

By assumption, $\hat u^2 $ and $\hat v^2 $ are in $F$, and $[G:F]=[H:F]=2$.
If $F$ is not isomorphic to $Q\times Z/dZ$ then $J(F)$ is abelian 
and so the (normal) subgroup generated by $F$ and $\hat s^2$ 
is isomorphic to $F\times\mathbb{Z}$.
In any case the subgroup generated by $F$ 
and $\hat s^k$ is normal, and is isomorphic to $F\times\mathbb{Z}$ if 
$k$ is a nonzero multiple of 12.
The uniqueness up to isomorphisms of the pair $\{ G,H\}$ follows 
from the uniqueness up to conjugation and order of the pair of generating involutions for $D$. 
Since $G$ and $H$ act freely on $\widetilde M$ they also have cohomological period dividing 4. 
On examining the list above we see that $F$ must be cyclic or the product of $Q(8k)$, $T(v)$ 
or $A(m,e)$ with a cyclic group of relatively prime order, 
as it is the kernel of a map from $G$ to $Z/2Z$. 
It is easily verified that in all such cases every automorphism 
of $F$ is the restriction of automorphisms of $G$ and $H$. 
Hence $\pi $ is determined up to isomorphism as the amalgamated free product 
$G*_F H$ by the unordered pair $\{ G,H\} $ of groups with index 2 subgroups 
isomorphic to $F$ (i.e., it is unnecessary to specify the 
identifications of $F$ with these subgroups). 

The third assertion follows because each of the spaces $M_G =\widetilde M/G$ 
and $M_H =\widetilde M/H$ are $PD_3 $-complexes 
with finite fundamental group and 
therefore are orientable, and $\pi$ is generated by $G$ and $H$.

The final assertion follows from a Mayer-Vietoris argument, 
as for parts (5) and (6) of Theorem 11.1. 
\end{proof}
 
Must the spaces $M_G $ and $M_H $ be homotopy equivalent to finite complexes?

In particular, if $\pi\cong D$ the $k$-invariant is unique,
and so any closed 4-manifold $M$ with $\pi_1 (M)\cong D$ and $\chi(M)=0$ 
is homotopy equivalent to $RP^4 \sharp RP^4 $.

\section{$\mathbb{S}^3 \times\mathbb{E}^1 $-manifolds}

With the exception of $A(m,1)$ (with $m>1$),
$O_k^* $ (with $k>1$), $P''_{48r}$ (with $r>1$)
and $Q(2^n a,b,c)$ (with either $n=3$ and at most one of $a$, $b$ and $c$ 
being 1 or $n>3$ and $bc>1$) 
and their products with cyclic groups,
all of the groups listed in \S2 have fixed point free representations 
in $SO(4)$ and so act freely on $S^3 $.
(Cyclic groups, the binary dihedral groups $D^*_{4m}=A(m,2)$, with $m$ odd, 
and $D^*_{8k}=Q(8k,1,1)$, with $k\geq1$ and
the three binary polyhedral groups $T_1^*$, 
$O^*_1$ and $I^*$ are subgroups of $S^3$.) 
We shall call such groups $\mathbb{S}^3$-{\it groups}.
A $k$-invariant in $H^4(F;\mathbb{Z})$ is {\it linear}
if it is realized by an $\mathbb{S}^3$-manifold $S^3/F$,
while it is {\it almost linear\/} if all covering spaces corresponding to 
subgroups isomorphic to $A(m,e)\times Z/dZ$ or $Q(8k)\times Z/dZ$ are 
homotopy equivalent to $\mathbb{S}^3$-manifolds \cite{[HM86]}.

Let $N$ be a $\mathbb{S}^3$-manifold with $\pi_1(N)=F$.
Then the projection of $Isom(N)$ onto its group 
of path components splits, and the inclusion of $Isom(N)$ into $Diff(N)$
induces an isomorphism on path components.
Moreover if $|F|>2$ isometries which induce the identity outer
automorphism are isotopic to the identity,
and so $\pi_0(Isom(N))$ maps injectively to $Out(F)$.
The group $\pi_0(Isom(N))$ has order 2 or 4, except when $F=Q(8)\times(Z/dZ)$,
in which case it has order 6 (if $d=1$) or 12 (if $d>1$).
(See \cite{[Mc02]}.)

\begin{theorem}
Let $M$ be a closed $4$-manifold with $\chi(M)=0$ and 
$\pi=\pi_1(M)\cong{F\rtimes_\theta\mathbb{Z}}$, where $F$ is finite.
Then $M$ is homeomorphic to a $\mathbb{S}^3\times\mathbb{E}^1$-manifold 
if and only if $M$ is the mapping torus 
of a self homeomorphism of a $\mathbb{S}^3$-manifold 
with fundamental group $F$, 
and such mapping tori are determined up to homeomorphism 
by their homotopy type.
\end{theorem}

\begin{proof} 
Let $p_1$ and $p_2$ be the projections of 
$Isom(\mathbb{S}^3\times\mathbb{E}^1)=
O(4)\times E(1)$ onto $O(4)$ and $E(1)$ respectively.
If $\pi$ is a discrete subgroup of $Isom(\mathbb{S}^3\times\mathbb{E}^1)$ 
which acts freely on $S^3\times\mathbb{R}$ then 
$p_1$ maps $F$ monomorphically and $p_1(F)$ acts freely on $S^3$,
since every isometry of $\mathbb{R}$ of finite order 
has nonempty fixed point set.
Moreover, $p_2(\pi)<E(1)$ acts discretely and cocompactly on $\mathbb{R}$, 
and so has no nontrivial finite normal subgroup.
Hence $F=\pi\cap(O(4)\times\{ 1\})$.
If $t\in\pi$ maps to a generator of $\pi/F\cong\mathbb{Z}$ 
then conjugation by $t$ induces an isometry $\theta$ of $S^3/F$, 
and $M\cong M(\theta)$.
Conversely, any self homeomorphism $h$ of a $\mathbb{S}^3$-manifold
is isotopic to an isometry of finite order,
and so $M(h)$ is homeomorphic to a
$\mathbb{S}^3\times\mathbb{E}^1$-manifold.
The final assertion follows from Theorem 3 of \cite{[Oh90]}.
\end{proof}

Every Swan complex for $Z/dZ$ is homotopy equivalent to a lens space $L(d,s)$.
(This follows from Theorem 2.9.)
All lens spaces have isometries which induce inversion on the group.
If $s^2\equiv\pm1$ {\it mod\/} $(d)$ there are also isometries of $L(d,s)$ 
which induce multiplication by $\pm{s}$.
No other nontrivial automorphism of $Z/dZ$ is realized by a
simple self homotopy equivalence of $L(d,s)$.
(See \S30 of \cite{[Co]}.)
However $(Z/dZ)\rtimes_s\mathbb{Z}$ may also be realized by mapping tori 
of self homotopy equivalences of other lens spaces.
If $d>2$ a $PD_4$-complex with this group and Euler characteristic 0 
is orientable if and only if $s^2\equiv 1$ {\it mod\/} $(d)$.

If $F$ is a noncyclic $\mathbb{S}^3$-group there is an unique orbit of
linear $k$-invariants under the action of $Out(F)\times\{\pm1\}$,
and so for each $\theta\in Aut(F)$ at most one homeomorphism class of 
$\mathbb{S}^3\times\mathbb{E}^1$-manifolds has fundamental group 
$\pi=F\rtimes_\theta\mathbb{Z}$.
If $F=Q(2^k)$ or $T_k^*$ for some $k>1$ then 
$S^3/F$ is the unique finite Swan complex for $F$ \cite{[Th80]}.
In general, there may be other finite Swan complexes. 
(In particular, there are exotic finite Swan complexes for $T_1^*$.)

Suppose now that $G$ and $H$ are $\mathbb{S}^3$-groups with index 2 subgroups 
isomorphic to $F$.
If $F$, $G$ and $H$ are each noncyclic then the corresponding 
$\mathbb{S}^3$-manifolds are uniquely determined, 
and we may construct a nonorientable $\mathbb{S}^3\times\mathbb{E}^1$-manifold 
with fundamental group $\pi=G*_FH$ as follows. 
Let $u$ and $v:S^3 /F\to S^3 /F$ be the covering involutions with 
quotient spaces $S^3/G$ and $S^3/H$, respectively, and let $\phi =uv$. 
(Note that $u$ and $v$ are isometries of $S^3 /F$.) 
Then $U([x,t])=[u(x),1-t]$ defines a fixed point free involution on the 
mapping torus $M(\phi )$ and the quotient space has fundamental group $\pi$. 
A similar construction works if $F$ is cyclic and $G\cong H$ or if $G$ is
cyclic.

\section{Realization of the invariants}

Let $F$ be a finite group with cohomological period dividing 4, and
let $X_F$ denote a finite Swan complex for $F$.
If $\theta$ is an automorphism of $F$ which induces $\pm1$ on
$H_3(F;\mathbb{Z})$ there is a self homotopy equivalence $h$ of $X_F$ which 
induces $[\theta]\in J(F)$.
The mapping torus $M(h)$ is a finite $PD_4 $-complex with $\pi_1(M)\cong
F\rtimes_\theta\mathbb{Z}$ and $\chi(M(h))=0$.
Conversely, every $PD_4$-complex $M$ with $\chi(M)=0$ and such that
$\pi_1(M)$ is an extension of $\mathbb{Z}$ by a finite normal subgroup $F$
is homotopy equivalent to such a mapping torus.
Moreover, if $\pi\cong F\times\mathbb{Z}$ and $|F|>2$ 
then $h$ is homotopic to the 
identity and so $M(h)$ is homotopy equivalent to $X_F\times S^1$.

The question of interest here is which such groups $\pi$
(and which $k$-invariants in $H^4(F;\mathbb{Z})$) 
may be realized by closed 4-manifolds.
The Spivak normal fibration of a $PD_3$-complex $X$ is trivial \cite{[Hb18]}. 
Hence every 3-dimensional Swan complex $X_F$ has a TOP reduction, 
i.e., there are normal maps $(f,b):N^3 \to X_F$.
Such a map has a ``proper surgery" obstruction $\lambda^p (f,b)$ in $L_3^p (F)$,
which is 0 if and only if $(f,b)\times id_{S^1} $ is normally cobordant to a 
simple homotopy equivalence. 
In particular, a surgery semicharacteristic must be 0.
Hence all subgroups of $F$ of order $2p$ (with $p$ prime) are cyclic, 
and $Q(2^n a,b,c)$ (with $n>3$ and $b$ or $c>1$) cannot occur \cite{[HM86]}.
As the $2p$ condition excludes groups with subgroups isomorphic to $A(m,1)$ 
with $m>1$ and as there are no finite Swan complexes for $O^*_k$ with $k>1$
or $P''_{48r}$ with $3|r$,
the cases remaining to be decided are when 
$F\cong Q(8a,b,c)\times Z/dZ$, 
where $a,b$ and $c$ are odd and at most one of them is 1,
and the groups $P''_{48a}$ with $(a,6)=1$. 
The main result of \cite{[HM86]} is that in such a case 
$F\times\mathbb{Z}$ acts freely and properly with almost linear $k$-invariant 
if and only if some arithmetical conditions depending on subgroups 
of $F$ of the form $Q(8a,b,1)$ hold.
(The constructive part of the argument may be
extended to the 4-dimensional case by reference to \cite{[FQ]}.)
  
The following more direct argument for the existence of a free proper action 
of $F\times\mathbb{Z}$ on $S^3\times\mathbb{R}$ was outlined in \cite{[KS88]}, 
for the cases when $F$ acts freely on an homology 3-sphere $\Sigma$.
Let $\Sigma$ and its universal covering space $\widetilde\Sigma$
have equivariant cellular decompositions lifted from
a cellular decomposition of $\Sigma/F$,
and let $\Pi=\pi_1 (\Sigma/F)$.
Then $C_* (\Sigma)=\mathbb{Z}[F]\otimes_\Pi C_* (\widetilde\Sigma)$ is 
a finitely generated free $\mathbb{Z}[F]$-complex,
and may be realized by a finite Swan complex $X$.
The chain map (over the epimorphism $:\Pi\to F$) from 
$C_* (\widetilde\Sigma)$ to $C_* (\widetilde X)$ may be realized by a
map $h:\Sigma/F\to X$, since these spaces are 3-dimensional.
As $h\times id_{S^1} $ is a simple 
$\mathbb{Z}[F\times\mathbb{Z}]$-homology equivalence 
it has surgery obstruction 0 in $L^s_4 (F\times\mathbb{Z})$,
and so is normally cobordant to a simple homotopy equivalence.
Is there a simple, explicit example of a free action of some
$Q(8a,b,1)$ (with $a,b>1$) on an homology 3-sphere?

Although $Q(24,13,1)$ cannot act freely on any homology 3-sphere 
\cite{[DM85]},
there is a closed orientable 4-manifold $M$ with fundamental group
$Q(24,13,1)\times\mathbb{Z}$, by the argument of \cite{[HM86]}.
The infinite cyclic cover $M_F$ is finitely dominated; 
is the Farrell obstruction to fibration in $Wh(\pi_1(M))$ nonzero?
No such 4-manifold can fibre over $S^1$, 
since $Q(24,13,1)$ is not a 3-manifold group.

If $F=T_k^*$ (with $k>1$), $Q(2^nm)$ or $A(m,2)$ (with $m$ odd)
$F\times\mathbb{Z}$ can only act freely and properly on 
$\mathbb{R}^4 \backslash\{ 0\}$ with the linear $k$-invariant.
This follows from Corollary C of \cite{[HM86']} for $A(m,2)$.
Since $Q(2^nm)$ contains $A(m,2)$ as a normal subgroup this implies that
the restriction of the $k$-invariant for a $Q(2^nm)$ action
to the cyclic subgroup of order $m$ must also be linear.
The nonlinear $k$-invariants for $T_k^*$ (with $k>1$) and
$Q(2^n)$ have nonzero finiteness obstruction.
As the $k$-invariants of free linear representations of $Q(2^nm)$ 
are given by elements in $H^4(Q(2^nm);\mathbb{Z})$ whose restrictions 
to $Z/mZ$ are squares and whose restrictions to $Q(2^n)$ are squares 
times the basic generator (see \cite[page 120]{[Wl78]}),
only the linear $k$-invariant is realizable in this case also.
However in general it is not known which $k$-invariants are realizable. 
Every group of the form $Q(8a,b,c)\times Z/dZ\times\mathbb{Z}$ 
admits an almost linear $k$-invariant, but there may be other actions.
(See \cite{[HM86],[HM86']} for more on this issue.)

When $F=Q(8)$, $T_k^*$, $O_1^*$, $I^*$ or $A(p^i,e)$ 
(for some odd prime $p$ and $e\geq2$) 
each element of $J(F)$ is realized by an isometry of $S^3/F$ 
\cite{[Mc02]}.
The study of more general groups may largely reduce to the case 
when $F$ is cyclic.

\begin{theorem}
Let $M$ be a $PD_4$-complex with $\chi(M)=0$ and
$\pi_1(M)\cong{F}\rtimes_\theta\mathbb{Z}$, 
where $F$ is a noncyclic $\mathbb{S}^3$-group,
and with $k(M)$ linear.
If the covering space associated to $C\rtimes_{\theta|C}\mathbb{Z}$
is homotopy equivalent to a $\mathbb{S}^3\times\mathbb{E}^1$-manifold 
for all characteristic cyclic subgroups $C<F$ then so is $M$.
\end{theorem}

\begin{proof}
Let $N=S^3/F$ be the $\mathbb{S}^3$-manifold with $k(N)=k(M)$,
and let $N_C=S^3/C$ for all subgroups $C<F$.
We may suppose that $F=G\times{Z/dZ}$ where $G$ is an $\mathbb{S}^3$-group
with no nontrivial cyclic direct factor.
If $C=Z/dZ$ then $N_C\cong{L(d,1)}$, and so $\theta|_C=\pm1$.
Hence $\theta$ is in $J(G)\times\{\pm1\}$.
If $G\cong{Q(8)}$, $T_1^*$, $O_1^*$, $I^*$ or $A(p^i,2)$ then 
$\pi_0(Isom(N))=J(G)$ (if $d=1$) or $J(G)\times\{\pm1\}$ (if $d>1$) 
\cite{[Mc02]}.
If $G\cong{T_k^*}$ with $k>1$ then $N_{\zeta{F}}\cong{L(3^{k-1}2d,1)}$,
so $\theta|_{\zeta{F}}=\pm1$,
and the nontrivial isometry of $N$ induces the involution of $\zeta{F}$.
In the remaining cases $F$ has a characteristic cyclic subgroup 
$C\cong{Z/4kdZ}$ of index 2, 
and $F$ acts on $C$ through $\sigma(c)=c^s$ for $c\in{C}$,
where $s\equiv-1$ {\it mod\/} $(4k)$ and $s\equiv1$ {\it mod\/} $(d)$.
Hence $N_C\cong{L(4kd,s)}$.
Restriction induces an epimorphism from $J(N)$ to
$J(C)/\langle\sigma\rangle$ with kernel of order 2,
and which maps $\pi_0(Isom(N))$ onto
$\pi_0(Isom(N_C))/\langle\sigma\rangle$.
In all cases $\theta$ is realized by an isometry of $N$.
Since $E(N)\cong{Out(F)}$ the result now follows from Theorem 11.6.
\end{proof}

When $F$ is cyclic the natural question is whether $k(M)$ and $k(L(d,s))$ agree,
if $M$ is a 4-manifold and $\pi_1(M)\cong(Z/dZ)\rtimes_s\mathbb{Z}$. 
If so, then every 4-manifold $M$ with $\chi(M)=0$,
$\pi_1(M)\cong{F}\rtimes_\theta\mathbb{Z}$ for some $\mathbb{S}^3$-group $F$
and $k(M)$ linear is homotopy equivalent to a 
$\mathbb{S}^3\times\mathbb{E}^1$-manifold.
Davis and Weinberger have settled the case when $M$ is {\it not\/} orientable.
Since $-1$ is then a square {\it mod} $(d)$, 
elementary considerations show that it suffices to assume that 
$d$ is prime and $d\equiv1$ {\it mod} $(4)$,
and thus to prove the following lemma.

\begin{lem}
{\rm[DW07]}\qua
Let $p$ be a prime such that $p\equiv1$ {\it mod} $(4)$.
Let $M(g)$ be the mapping torus of an orientation-reversing 
self-homotopy equivalence $g$ of $L(p,q)$, where $(p,q)=1$.
Suppose that $M(g)$ is homotopy equivalent to a $4$-manifold.
If $p\equiv1$ {\it mod} $(8)$ then $q$ is a quadratic residue
{\it mod} $(p)$, while if $p\equiv5$ {\it mod} $(8)$ then $q$ 
is a quadratic nonresidue {\it mod} $(p)$.
\qed
\end{lem}

Since $g$ is a self homotopy equivalence of a closed orientable
3-manifold there is a degree-1 normal map $F:W\to{L(p,q)}\times[0,1]$
with $F|_{\partial{W}}=g\amalg{id_{L(p,q)}}$.
(See \cite[Theorem 2]{[JK03]}.)
Gluing the ends gives a degree-1 normal map $G:Y\to{M(g)}$,
with domain a closed manifold.
The key topological observation is that if $M(g)$ 
is also homotopy equivalent to a closed manifold then 
$\sigma_4^h(G)\in{L_4^h(\pi,w)}$ is in the image 
of the assembly map from $H_4(M(g);\mathbb{L}\langle1\rangle_w)$.
(See \cite[Proposition 18.3]{[Rn]}.)
The rest of the argument involves delicate calculations.

If $M$ is a closed non-orientable 4-manifold with 
$\pi_1(M)\cong(Z/dZ)\rtimes_s\mathbb{Z}$ and $\chi(M)=0$ then $M\simeq{M(g)}$, 
where $g$ is a self homotopy equivalence of some lens space $L(d,q)$,
and $s^2\equiv-1$ {\it mod} $(d)$.
Hence if a prime $p$ divides $d$ then $p\equiv1$ {\it mod} $(4)$.
The Lemma and the Chinese Remainder Theorem imply
that $q\equiv{sz^2}$ {\it mod} $(d)$, for some $(z,d)=1$,
and so $M$ is homotopy equivalent to the mapping torus 
of an isometry of $L(d,s)$. 
See \cite{[DW07]} for details.

In the orientable case $s^2\equiv1$ {\it mod} $(d)$,
and so $s\equiv\pm1$ {\it mod} $(p)$, for each odd prime factor of $d$.
(If $d$ is even then $s\equiv\pm1$ {\it mod} $(2p)$.)
Elementary considerations now show that it suffices to assume
that $d=2^{r+1}$, $2^rp$ or $pp'$, for some $r\geq2$ and odd primes $p<p'$.
However this question is open even for the smallest cases (with $d=8$ or 12).

It can be shown that when $(d,q,s)=(5,1,2)$ or $(8,1,3)$
the mapping torus $M(g)$ is a {\it simple\/} $PD_4$-complex.
(This negates an earlier hope.)

When $\pi/F\cong D$ we have $\pi\cong G*_FH$, and we saw earlier that if
$G$ and $H$ are $\mathbb{S}^3$-groups and $F$ is noncyclic $\pi$ 
is the fundamental group of a $\mathbb{S}^3\times\mathbb{E}^1$-manifold.
However if $F$ is cyclic but neither $G$ nor $H$ is cyclic there may be
no geometric manifold realizing $\pi$.
If the double covers of $G\backslash S^3$ and $H\backslash S^3$ 
are homotopy equivalent then $\pi$ is realised by the union 
of two twisted $I$-bundles via a homotopy equivalence, 
which is a finite $PD_4 $-complex with $\chi=0$.
For instance, the spherical space forms corresponding to
$G=Q(40)$ and $H=Q(8)\times(Z/5Z)$ are doubly covered by
$L(20,1)$ and $L(20,9)$, respectively, 
which are homotopy equivalent but not homeomorphic.
The spherical space forms corresponding to $G=Q(24)$ and $H=Q(8)\times(Z/3Z)$ 
are doubly covered by $L(12,1)$ and $L(12,5)$, respectively, 
which are not homotopy equivalent.

In each case it remains possible that some extensions of $\mathbb{Z}$ or $D$
by normal $\mathbb{S}^3$-subgroups may be realized 
by manifolds with nonlinear $k$-invariants.

\section{$T$- and $Kb$-bundles over $RP^2$ with $\partial\not=0$}

Let $p:E\to RP^2$ be a bundle with fibre $T$ or $Kb$.
Then $\pi=\pi_1(E)$ is an extension of $Z/2Z$ by $G/\partial\mathbb{Z}$, 
where $G$ is the fundamental group of the fibre and $\partial$ is the
connecting homomorphism. If $\partial\not=0$ then $\pi$ has two ends,
$F$ is cyclic and central in $G/\partial\mathbb{Z}$ and $\pi$ 
acts on it by inversion, 
since $\pi$ acts nontrivially on $\mathbb{Z}=\pi_2 (RP^2)$.

If the fibre is $T$ then $\pi$ has a presentation of the form
\begin{equation*}
\langle t,u,v\mid uv=vu,\medspace u^n=1,\medspace tut^{-1}=u^{-1}\! ,
\medspace tvt^{-1}=u^a v^\epsilon\! ,\medspace t^2 =u^b v^c\rangle,
\end{equation*}
where $n>0$ and $\epsilon=\pm 1$. Either
\begin{enumerate}
\item $F$ is cyclic, 
$\pi\cong (Z/nZ)\rtimes_{-1}\mathbb{Z}$ and $\pi/F\cong\mathbb{Z}$; or

\item $F=\langle s,u\mid s^2=u^m\!,\medspace sus^{-1}=u^{-1} \rangle$;
or (if $\epsilon=-1$)   

\item $F$ is cyclic, $\pi=\langle s,t,u\mid s^2=t^2=u^b\! ,
\medspace sus^{-1} =tut^{-1} =u^{-1}\rangle$ 
and $\pi/F\cong D$.
\end{enumerate}
 
In case (2) $F$ cannot be dihedral.
If $m$ is odd $F\cong A(m,2)$ while if 
$m=2^r k$ with $r\geq 1 $ and $k$ odd $F\cong Q(2^{r+2} k)$. 
On replacing $v$ by $u^{[a/2]} v$, if necessary, we may arrange that $a=0$, 
in which case $\pi\cong F\times\mathbb{Z}$, or $a=1$, in which case
\begin{equation*}
\pi=\langle t,u,v\mid t^2=u^m\! ,\medspace tut^{-1}=u^{-1},
\medspace vtv^{-1}=tu,\medspace uv=vu\rangle,
\end{equation*} 
so $\pi/F\cong\mathbb{Z}$.
                   
If the fibre is $Kb$ then $\pi$ has a presentation of the form
\begin{equation*}
\langle t,u,w\mid uwu^{-1} =w^{-1}\! ,\medspace u^n =1,\medspace 
tut^{-1}=u^{-1}\! ,\medspace twt^{-1}=u^a w^\epsilon,
\medspace t^2 =u^b w^c\rangle,
\end{equation*}
where $n>0$ is even (since $\mathrm{Im}(\partial)\leq\zeta\pi_1(Kb)$) 
and $\epsilon=\pm 1$. 
On replacing $t$ by $ut$, if necessary, we may assume that $\epsilon=1$.
Moreover, $tw^2t^{-1}=w^{\pm2}$ since $w^2$ generates the commutator 
subgroup of $G/\partial\mathbb{Z}$, so $a$ is even and $2a\equiv 0~mod~(n)$,
$t^2 u=ut^2 $ implies that $c=0$, and $t.t^2.t^{-1}=t^2$ implies that 
$2b\equiv 0~mod~(n)$. 
As $F$ is generated by $t$ and $u^2$, and cannot be dihedral, 
we must have $n=2b$.
Moreover $b$ must be even, as $w$ has infinite order and $t^2w=wt^2$.
Therefore

\smallskip
(4) $F\cong Q(8k)$, $\pi/F\cong D$ and

$\pi=\langle t,u,w\mid uwu^{-1} =w^{-1}\! ,\medspace tut^{-1}=u^{-1}\! ,
\medspace tw=u^a wt, \medspace t^2 =u^{2k}\rangle$.

\smallskip
\noindent In all cases $\pi$ has a subgroup of index at most 2 which is 
isomorphic to $F\times\mathbb{Z}$.

Each of these groups is the fundamental group of such a bundle space.
(This may be seen by using the description of such bundle spaces given
in \S5 of Chapter 5.)
Orientable 4-manifolds which fibre over $RP^2$ with fibre $T$ and
$\partial\not=0$ are mapping tori of involutions of $\mathbb{S}^3$-manifolds,
and if $F$ is not cyclic two such bundle 
spaces with the same group are diffeomorphic \cite{[Ue91]}.

\begin{theorem}
Let $M$ be a closed orientable $4$-manifold 
with fundamental group $\pi$.
Then $M$ is homotopy equivalent to an $\mathbb{S}^3\times\mathbb{E}^1$-manifold
which fibres over $RP^2$ if and only $\chi(M)=0$ and $\pi$ is of 
type $(1)$ or $(2)$ above.
\end{theorem}

\begin{proof} If $M$ is an orientable 
$\mathbb{S}^3\times\mathbb{E}^1$-manifold then $\chi(M)=0$ 
and $\pi/F\cong\mathbb{Z}$,
by Theorem 11.1 and Lemma 3.14. 
Moreover $\pi$ must be of type (1) or (2) if $M$ fibres over
$RP^2$, and so the conditions are necessary. 

Suppose that they hold. 
Then $\widetilde M\cong\mathbb{R}^4 \setminus\{ 0\} $ 
and the homotopy type of $M$ is determined by $\pi$ and $k(M)$, 
by Theorem 11.1.
If $F\cong Z/nZ$ then $M_F=\widetilde M/F$ is homotopy equivalent to some 
lens space $L(n,s)$.
As the involution of $Z/nZ$ which inverts a generator
can be realized by an isometry of $L(n,s)$, $M$ is homotopy equivalent to
an $\mathbb{S}^3\times\mathbb{E}^1 $-manifold which fibres over $S^1 $. 

If $F\cong Q(2^{r+2} k)$ or $A(m,2)$ then $F\times\mathbb{Z}$ 
can only act freely and properly on $\mathbb{R}^4 \setminus\{ 0\}$ 
with the ``linear" $k$-invariant
\cite{[HM86]}.
Therefore $M_F$ is homotopy equivalent to a spherical space form $S^3/F$.
The class in $Out(Q(2^{r+2} k))$ represented by the automorphism which sends 
the generator $t$ to $tu$ and fixes $u$ is induced by conjugation 
in $Q(2^{r+3} k)$ and so can be realized by a (fixed point free) isometry 
$\theta$ of $S^3 /Q(2^{r+2} k)$.
Hence $M$ is homotopy equivalent to a bundle space 
$(S^3 /Q(2^{r+2} k))\times S^1$ or 
$(S^3 /Q(2^{r+2} k))\times_\theta S^1$ if $F\cong Q(2^{r+2} k)$. 
A similar conclusion holds when $F\cong A(m,2)$ as the corresponding 
automorphism is induced by conjugation in $Q(2^3 d)$.

With the results of \cite{[Ue91]} it follows in all cases that $M$ 
is homotopy equivalent to the total space of a torus bundle over $RP^2 $. 
\end{proof}
 
Theorem 11.8 does not assume that there is a homomorphism 
${u:\pi\to Z/2Z}$ such that $u^* (x)^3=0$ (as in \S5 of Chapter 5).
If $F$ is cyclic or $A(m,2)$ this condition is a purely algebraic 
consequence of the other hypotheses.
For let $C$ be a cyclic normal subgroup of maximal order in $F$.
(There is an unique such subgroup, except when $F=Q(8)$.)
The centralizer $C_\pi (C)$ has index 2 in $\pi$ and so there is a 
homomorphism
$u:\pi\to Z/2Z$ with kernel $C_\pi (C)$.

When $F$ is cyclic $u$ factors through $\mathbb{Z}$ 
and so the induced map on cohomology
factors through $H^3 (\mathbb{Z};\mathbb{Z}^u)=0$.

When $F\cong A(m,2)$ the 2-Sylow subgroup is cyclic of order 4,
and the inclusion of $Z/4Z$ into $\tau$ induces isomorphisms 
on cohomology with 2-local coefficients. 
In particular, $H^q (F;\mathbb{Z}_{(2)}^u )=0$ 
or $Z/2Z$ according as $q$ is even 
or odd. It follows easily that the restriction from 
$H^3 (\pi;\mathbb{Z}_{(2)}^u )$ to
$H^3 (Z/4Z;\mathbb{Z}_{(2)}^u )$ is an isomorphism.
Let $y$ be the image of $u^* (x)$ in $H^1 (Z/4Z;\mathbb{Z}_{(2)}^u)=Z/2Z$.
Then $y^2 $ is an element of order 2 in 
$H^2(Z/4Z;\mathbb{Z}_{(2)}^u\otimes\mathbb{Z}_{(2)}^u)=
H^2(Z/4Z;\mathbb{Z}_{(2)})\cong Z/4Z$,
and so $y^2=2z$ for some $z\in H^2 (Z/4Z;\mathbb{Z}_{(2)} )$.
But then $y^3=2yz=0$ in $H^3(Z/4Z;\mathbb{Z}_{(2)}^u)=Z/2Z$, 
and so $u^*(x)^3$ has image 0 in $H^3 (\pi;\mathbb{Z}_{(2)}^u)=Z/2Z$.
Since $x$ is a 2-torsion class this implies that $u^*(x)^3 =0$.      
                               
Is there a similar argument when $F$ is a generalized quaternionic group?

If $M$ is nonorientable, $\chi(M)=0$ and has fundamental group $\pi$ of type 
(1) or (2) then $M$ is homotopy equivalent to the 
mapping torus of the orientation reversing 
self homeomorphism of $S^3 $ or of $RP^3 $, and does not fibre over $RP^2 $.
If $\pi$ is of type (3) or (4) then the 2-fold covering space 
with fundamental group $F\times\mathbb{Z}$ is homotopy equivalent 
to a product $L(n,s)\times S^1 $. 
However we do not know which $k$-invariants give
total spaces of bundles over $RP^2 $.

\section{Some remarks on the homeomorphism types}

In this brief section we shall assume that $M$ is orientable
(except in the final sentence) and that 
$\pi\cong F\rtimes_\theta\mathbb{Z}$.
In contrast to the situation for the other geometries, the Whitehead groups
of fundamental groups of $\mathbb{S}^3\times\mathbb{E}^1$-manifolds are usually
nontrivial.
Computation of $Wh(\pi)$ is difficult as the $Nil$ groups occuring in the 
Waldhausen exact sequence relating $Wh(\pi)$ to the algebraic 
$K$-theory of $F$ seem intractable.

We can however compute the relevant surgery obstruction groups 
modulo 2-torsion and show that the structure sets are usually infinite.
There is a Mayer-Vietoris sequence 
$L_5^s(F)\to L_5^s(\pi)\to L_4^u(F)\to L_4^s(F)$, 
where the superscript $u$ signifies that
the torsion must lie in a certain subgroup of $Wh(F)$ \cite{[Ca73]}.
The right hand map is (essentially) $\theta_*-1$.
Now $L_5^s(F)$ is a finite 2-group and $L_4^u(F)\sim L_4^s(F)\sim\mathbb{Z}^R$ 
{\it mod\/} 2-torsion, 
where $R$ is the set of irreducible real representations of $F$ 
(see \cite[Chapter 13A]{[Wl]}).
The latter correspond to the conjugacy classes of $F$, up to inversion.
(See \cite[\S12.4]{[Se]}.)
In particular, if $\pi\cong F\times\mathbb{Z}$ then 
$L_5^s(\pi)\sim\mathbb{Z}^R$ {\it mod\/}
2-torsion, and so has rank at least 2 if $F\not=1$.
As $[\Sigma M,G/TOP]\cong\mathbb{Z}$ {\it mod\/} 2-torsion 
and the group of self homotopy equivalences of such a manifold is finite, 
by Theorem 11.2, there are infinitely many distinct topological 
4-manifolds simple homotopy equivalent to $M$.

For instance, as $Wh(\mathbb{Z}\oplus(Z/2Z))=0$ \cite{[Kw86]} and 
$L_5 (\mathbb{Z}\oplus(Z/2Z),+)\cong\mathbb{Z}^2 $
\cite[Theorem 13A.8]{[Wl]}, 
the set $S_{TOP} (RP^3 \times S^1 )$ is infinite. 
Although all of the manifolds in this homotopy type are doubly covered by 
$S^3 \times S^1$ only $RP^3 \times S^1 $ is itself geometric.
Similar estimates hold for the other manifolds covered by 
${S^3\times\mathbb{R}}$ (if $\pi\not\cong\mathbb{Z}$).

An explicit parametrization of the set of homeomorphism classes
of manifolds homotopy equivalent to $RP^4\#{RP^4}$ may be found
in \cite{[BDK07]}.

%% file: m5-12.tex
\chapter{Geometries with compact models}

There are three geometries with compact models, namely $\mathbb{S}^4$, 
$\mathbb{CP}^2$ and ${\mathbb{S}^2\times\mathbb{S}^2}$.
The first two of these are easily dealt with, as there is only one other geometric manifold,
namely $RP^4$, and for each of the two projective spaces there is one other 
(nonsmoothable) manifold of the same homotopy type.
There are eight $\mathbb{S}^2\times\mathbb{S}^2$-manifolds, 
seven of which are total spaces of bundles
with base and fibre each $S^2$ or $RP^2$.
We shall consider also the two other such bundle spaces
covered by $S^2\tilde\times S^2$, 
although they are not geometric.

The universal covering space $\widetilde M$ of a closed 4-manifold $M$ 
is homeomorphic to $S^2\times S^2$ if and only if $\pi=\pi_1(M)$ is finite, 
$\chi(M)|\pi|=4$ and $w_2(\widetilde M)=0$. 
(The condition $w_2 (\widetilde M)=0$ may be restated entirely in terms of $M$, 
but at somewhat greater length.) 
If these conditions hold and $\pi$ is cyclic then $M$ is homotopy equivalent 
to an $\mathbb{S}^2\times\mathbb{S}^2$-manifold, 
except when $\pi=Z/2Z$ and $M$ is nonorientable, 
in which case there is one other (non-geometric) homotopy type.
The $\mathbb{F}_2$-cohomology ring, Stiefel-Whitney classes and $k$-invariants 
must agree with those of bundle spaces when $\pi\cong (Z/2Z)^2$.
However there is again one other (non-geometric) homotopy type.
If $\chi(M)|\pi|=4$ and $w_2 (\widetilde M)\not=0$ then either $\pi=1$,
in which case $M\simeq S^2\tilde\times S^2$ or $CP^2\sharp CP^2$, 
or $M$ is nonorientable and $\pi=Z/2Z$;                    
in the latter case $M\simeq RP^4\sharp CP^2 $, the nontrivial $RP^2$-bundle 
over $S^2$, and $\widetilde M\simeq S^2\tilde\times S^2$.       

The number of homeomorphism classes within each homotopy type is one or two 
if $\pi=Z/2Z$ and $M$ is orientable, two if $\pi=Z/2Z$, $M$ is nonorientable 
and $w_2 (\widetilde M)=0$, 
four if $\pi=Z/2Z$ and $w_2(\widetilde M)\not=0$,  
four if $\pi\cong Z/4Z$, and at most eight if $\pi\cong (Z/2Z)^2$.
In the final case we do not know whether there are enough 
fake self homotopy equivalences to account for all the normal invariants 
with trivial surgery obstruction. 
However, in (at least) nine of the 13 cases a PL 4-manifold with 
the same homotopy type as a geometric manifold 
or $S^2 \tilde\times S^2 $ is homeomorphic to it.
(In seven of these cases the homotopy type is determined by the Euler characteristic, 
fundamental group and Stiefel-Whitney classes.)
Each nonorientable manifold has a fake twin, with the same homotopy type
but opposite Kirby-Siebenmann invariant.

For the full details of some of the arguments when $\pi\cong{Z/2Z}$ 
we refer to the papers \cite{[KKR92],[HKT94]} and \cite{[Te97]},
and when $\pi\cong{Z/4Z}$ we refer to \cite{[HH21]}.

\section{The geometries $\mathbb{S}^4$ and $\mathbb{CP}^2$}

The unique element of $Isom(\mathbb{S}^4)=O(5)$ of order 2 
which acts freely on $S^4$ is $-I$.
Therefore $S^4$ and $RP^4$ are the only $\mathbb{S}^4$-manifolds.
The manifold $S^4$ is determined up to homeomorphism 
by the conditions $\chi(S^4)=2$ and $\pi_1(S^4)=1$ \cite{[FQ]}.
                                                                                
\begin{lemma}
A closed $4$-manifold $M$ is homotopy equivalent to $RP^4$ 
if and only if $\chi(M)=1$ and $\pi_1 (M)=Z/2Z$.
\end{lemma}

\begin{proof} 
The conditions are clearly necessary. Suppose that they hold. 
Then $\widetilde M\simeq S^4 $ and $w_1 (M)=w_1 (RP^4)=w$, say, 
since any orientation preserving self
homeomorphism of $\widetilde M$ has Lefshetz number 2.
Since $RP^\infty=K(Z/2Z,1)$ may be obtained from $RP^4$ 
by adjoining cells of dimension at least 5
we may assume $c_M=c_{RP^4}f$, where $f:M\to RP^4 $.
Since $c_{RP^4} $ and $c_M$ are each 4-connected 
$f$ induces isomorphisms on homology with coefficients $Z/2Z$.
Considering the exact sequence of homology corresponding to the short
exact sequence of coefficients 
\[
0\to\mathbb{Z}^w\to\mathbb{Z}^w\to Z/2Z\to 0,
\]
we see that $f$ has odd degree.
By modifying $f$ on a 4-cell $D^4\subset M$ we may arrange that 
$f$ has degree 1, and the lemma then follows from Theorem 3.2. 
\end{proof} 

This lemma may also be proven by comparison of the $k$-invariants 
of $M$ and $RP^4$, 
as in \cite[Theorem 4.3]{[Wl67]}.

The surgery obstruction homomorphism is determined by an Arf invariant 
and maps $[RP^4;G/TOP]$ onto $Z/2Z$ \cite[Theorems 13.A.1 and 13.B.5]{[Wl]}, 
and hence $S_{TOP}(RP^4)$ has two elements.
(See the discussion of nonorientable manifolds with fundamental group 
$Z/2Z$ in \S6 below for more details.)
As every self homotopy equivalence of $RP^4 $ is homotopic to the identity
\cite{[Ol53]} there is one fake $RP^4$. 
The fake $RP^4$ is denoted $*RP^4$ and is not smoothable \cite{[Ru84]}.
  
There is a similar characterization of the homotopy type of 
the complex projective plane.

\begin{lemma}
A closed $4$-manifold $M$ is homotopy equivalent to $CP^2 $ 
if and only if $\chi(M)=3$ and $\pi_1 (M)=1$.
\end{lemma}
                                           
\begin{proof} 
The conditions are clearly necessary. Suppose that they hold.
Then $H^2 (M;\mathbb{Z})$ is infinite cyclic and so 
there is a map $f_M:M\to CP^\infty=K(Z,2)$ 
which induces an isomorphism on $H^2 $. 
Since $CP^\infty$ may be obtained from $CP^2$ by 
adjoining cells of dimension at least 6 
we may assume $f_M =f_{CP^2} g$, where $g:M\to CP^2 $ and 
$f_{CP^2} :CP^2 \to CP^\infty$ is the natural inclusion. 
As $H^4 (M;\mathbb{Z})$ is generated by $H^2 (M;\mathbb{Z})$, 
by Poincar\'e duality,
$g$ induces an isomorphism on cohomology and hence is a homotopy equivalence. 
\end{proof} 

In this case the surgery obstruction homomorphism is determined by the difference of signatures 
and maps $[CP^2 ;G/TOP]$ onto $\mathbb{Z}$. 
The structure set $S_{TOP} (CP^2 )$ again has two elements.
Since $[CP^2,CP^2]\cong [CP^2,CP^\infty]\cong H^2 (CP^2;\mathbb{Z})$, by obstruction theory,
there are two homotopy classes of self homotopy equivalences, 
represented by the identity and by complex conjugation.
Thus every self homotopy equivalence of $CP^2$ is homotopic to a homeomorphism,
and so there is one fake $CP^2 $. 
The fake $CP^2 $ is also known as the Chern manifold $Ch$ or $*CP^2 $, 
and is not smoothable \cite{[FQ]}.
Neither of these manifolds admits a nontrivial fixed point free action, 
as any self map of $CP^2$ or $*CP^2$ has nonzero Lefshetz number, 
and so $CP^2$ is the only $\mathbb{CP}^2 $-manifold.

\section{The geometry $\mathbb{S}^2\times\mathbb{S}^2$}

The manifold $S^2\times S^2$ is determined up to homotopy equivalence 
by the conditions $\chi(S^2 \times S^2)=4$, 
$\pi_1(S^2\times S^2)=1$ and $w_2(S^2\times S^2)=0$,
by Theorem 5.19. 
These conditions in fact determine $S^2\times S^2$ up to homeomorphism 
\cite{[FQ]}.
Hence if $M$ is an $\mathbb{S}^2\times\mathbb{S}^2$-manifold 
its fundamental group $\pi$ is finite, 
$\chi(M)|\pi|=4$ and $w_2 (\widetilde M)=0$. 

The isometry group of $\mathbb{S}^2\times\mathbb{S}^2$ 
is a semidirect product $(O(3)\times{O(3)})\rtimes(Z/2Z)$.
The $Z/2Z$ subgroup is generated by the involution $\tau$
which switches the factors ($\tau (x,y)=(y,x)$), 
and acts on $O(3)\times O(3)$ by $\tau(A,B)\tau=(B,A)$ for $A,B\in O(3)$.
In particular, $(\tau(A,B))^2=id$ if and only if $AB=I$, 
and so such an involution fixes $(x,Ax)$, for any $x\in S^2 $. 
Thus there are no free $Z/2Z$-actions in which the factors are switched.
The element $(A,B)$ generates a free $Z/2Z$-action if and only if $A^2 =B^2 =I$
and at least one of $A,B$ acts freely, i.e.\ if $A$ or $B=-I$.   
After conjugation with $\tau$ if necessary we may assume that $B=-I$, 
and so there are four conjugacy classes in 
$Isom(\mathbb{S}^2 \times\mathbb{S}^2 )$ of free $Z/2Z$-actions. 
(These may be distinguished by the multiplicity  
of 1 as an eigenvalue of $A$.) 
In each case projection onto the second factor induces a fibre bundle
projection from the orbit space to $RP^2 $, with fibre $S^2 $.
If $A\not=I$ the other projection induces an orbifold bundle 
with general fibre $S^2$ over $\mathbb{D}$,
$S(2,2)$ or $RP^2$.
                                       
If the involutions $(A,B)$ and $(C,D)$ generate a free $(Z/2Z)^2 $-action 
$(AC,BD)$ is also a free involution. 
By the above paragraph, one element of each of these ordered pairs must be $-I$.
It follows that (after conjugation with $\tau$ if necessary) 
the $(Z/2Z)^2 $-actions are generated by pairs $\{(A,-I),(-I,I)\}$, 
where $A^2 =I$. 
Since $A$ and $-A$ give rise to the same subgroup, 
there are two free $(Z/2Z)^2 $-actions.
The orbit spaces are the total spaces of $RP^2 $-bundles over $RP^2 $.

If $(\tau(A,B))^4 =id$ then $(BA,AB)$ is a fixed point free involution 
and so $BA=AB=-I$.
Since $(A,I)\tau(A,-A^{-1} )(A,I)^{-1} =\tau(I,-I)$ every free $Z/4Z$-action 
is conjugate to the one generated by $\tau(I,-I)$. 
The manifold $S^2\times{S^2}/\langle\tau(I,-I)\rangle$
does not fibre over a surface.
It is however the union of the tangent disc bundle of $RP^2$ 
(the image of the diagonal of ${S^2\times{S^2}}$)  
with the mapping cylinder of the double cover of $L(8,1)$.
(See \S12.9 below.)

In the next section we shall see that these eight geometric manifolds 
may be distinguished by their fundamental group and Stiefel-Whitney classes.
Note that if $F$ is a finite group then $q(F)\geq 2/|F|>0$, 
while $q^{SG} (F)\geq 2$.
Thus $S^4 $, $RP^4 $ and the geometric manifolds with $|\pi|=4$ 
have minimal Euler characteristic for
their fundamental groups (i.e., $\chi(M)=q(\pi)$), 
while $S^2 \times S^2 /(-I,-I)$ has
minimal Euler characteristic among 
$PD_4^+$-complexes realizing $Z/2Z$. 

\section{Bundle spaces}

There are two $S^2$-bundles over $S^2$, since $\pi_1 (SO(3))=Z/2Z$.
The total space $S^2 \tilde \times S^2 $ of the nontrivial $S^2 $-bundle
over $S^2 $ is determined up to homotopy equivalence by the conditions 
$\chi(S^2 \tilde \times S^2)=4$, $\pi_1 (S^2 \tilde \times S^2)=1$,
$w_2 (S^2 \tilde \times S^2)\not=0$ and $\sigma (S^2 \tilde\times S^2 )=0$,
by Theorem 5.19. 
The bundle space is homeomorphic to the connected sum $CP^2 \sharp -CP^2 $.
However there is one fake $S^2 \tilde \times S^2$, 
which is homeomorphic to $CP^2 \sharp -*CP^2$ and is not smoothable 
\cite{[FQ]}.
The manifolds $CP^2\sharp CP^2$ and $CP^2\sharp *CP^2$ also have
$\pi_1=0$ and $\chi=4$. 
It is easily seen that any self homotopy equivalence of
either of these manifolds has nonzero Lefshetz number,
and so they do not properly cover any other 4-manifold.

Since the Kirby-Siebenmann obstruction of a closed 4-manifold 
is natural with respect to covering maps and dies on passage 
to 2-fold coverings,                         
the non-smoothable manifold $CP^2\sharp -*CP^2$ admits no nontrivial 
free involution.
The following lemma implies that $S^2\tilde\times S^2$ admits no orientation 
preserving free involution, and hence no free action of $Z/4Z$ or $(Z/2Z)^2 $.        
\begin{lemma}
Let $M$ be a closed $4$-manifold with fundamental 
group $\pi=Z/2Z$ and universal covering space $\widetilde M$. Then 
\begin{enumerate}
\item $w_2 (\widetilde M)=0$ if and only if $w_2 (M)=u^2 $ for some 
$u\in H^1(M;\mathbb{F}_2)$; and

\item if $M$ is orientable and $\chi(M)=2$ then $w_2(\widetilde M)=0$ 
and so $\widetilde M\cong S^2\times S^2$.
\end{enumerate}
\end{lemma}

\begin{proof} 
The Cartan-Leray cohomology spectral sequence 
(with coefficients $\mathbb{F}_2$) 
for the projection $p:\widetilde M\to M$ gives an exact sequence 
\[
0\to H^2 (\pi;\mathbb{F}_2)\to H^2 (M;\mathbb{F}_2)\to H^2 (\widetilde M;\mathbb{F}_2),
\]
in which the right hand map is induced by $p$ and has image 
in the subgroup fixed under the action of $\pi$.  
Hence $w_2 (\widetilde M)=p^* w_2 (M)$ is 0 if and only if 
$w_2 (M)$ is in the image of $H^2(\pi;\mathbb{F}_2)$.
Since $\pi=Z/2Z$ this is so if and only if $w_2 (M)=u^2$ 
for some $u\in H^1(M;\mathbb{F}_2)$. 

Suppose that $M$ is orientable and $\chi(M)=2$. 
Then $H^2 (\pi;\mathbb{Z})=H^2 (M;\mathbb{Z})=Z/2Z$.
Let $x$ generate $H^2 (M;\mathbb{Z})$ and let $\bar x$ be its image under 
reduction modulo (2) in $H^2 (M;\mathbb{F}_2)$. 
Then $\bar x\cup\bar x=0$ in $H^4 (M;\mathbb{F}_2)$ 
since $x\cup x=0$ in $H^4 (M;\mathbb{Z})$.
Moreover as $M$ is orientable $w_2(M)=v_2(M)$ and so 
$w_2 (M)\cup \bar x=\bar x\cup\bar x=0$.
Since the cup product pairing on $H^2 (M;\mathbb{F}_2)\cong (Z/2Z)^2 $ 
is nondegenerate it follows that $w_2 (M)=\bar x$ or 0. 
Hence $w_2 (\widetilde M)$ is the reduction of $p^* x$ or is 0.
The integral analogue of the above exact sequence implies that the natural map 
from $H^2 (\pi;\mathbb{Z})$ to $H^2 (M;\mathbb{Z})$ is an isomorphism and so $p^*(H^2 (M;\mathbb{Z}))=0$. 
Hence $w_2 (\widetilde M)=0$ and so $\widetilde M\cong S^2 \times S^2 $. 
\end{proof} 

There are two $S^2$-bundles over $Mb$, since $\pi_1(BO(3))=Z/2Z$.
Each restricts to a trivial bundle over $\partial Mb$.
A map from $\partial Mb$ to $O(3)$ extends across $Mb$ 
if and only if it is homotopic to a constant map,
since $\pi_1(O(3))=Z/2Z$, and so there are four $S^2$-bundles over 
$RP^2=Mb\cup D^2$. (See also Theorem 5.10.)

The orbit space $M=(S^2 \times S^2)/\langle(A,-I)\rangle$ fibres over $RP^2$,
and is orientable if and only if $det(A)=-1$.
If $A$ has a fixed point $P\in S^2 $ then the image of $\{ P\}\times S^2 $ 
in $M$ is a section which represents a nonzero class 
in $H_2 (M;\mathbb{F}_2)$.
If $A=I$ or is a reflection across a plane the fixed point set 
has dimension $>0$ 
and so this section has self intersection 0. 
As the fibre $S^2 $ intersects the section in one point 
and has self intersection 0 it follows that $v_2 (M)=0$ 
and so $w_2 (M)=w_1 (M)^2 $ in these two cases.        
If $A$ is the half-turn about the $z$-axis let $n(x,y)=\frac1r(x,y,1)$, 
where $r=\sqrt{x^2+y^2+1}$, for $(x,y)\in\mathbb{R}^2$.
Let $\sigma_t[\pm{s}]=[n(tx,ty),s]\in{M}$,
for $s=(x,y,z)\in{S^2}$ and $|t|$ small.
Then $\sigma_t$ is an isotopy of sections with self intersection 1.
Finally, if $A=-I$ then the image of the diagonal $\{ (x,x)|x\in S^2 \}$ 
is a section with self intersection 1. 
Thus in these two cases $v_2 (M)\not=0$.
Therefore, by part (1) of the lemma, $w_2 (M)$ is the square 
of the nonzero element of $H^1 (M;\mathbb{F}_2)$ if $A=-I$ and is 0 if $A$ is a rotation.
Thus these bundle spaces may be distinguished by their Stiefel-Whitney classes, 
and every $S^2 $-bundle over $RP^2 $ is geometric. 

The group $E(RP^2 )$ is connected and the natural map 
from $SO(3)$ to $E(RP^2 )$ induces an isomorphism on $\pi_1 $,
by Lemma 5.15.
Hence there are two $RP^2 $-bundles over $S^2 $, up to fibre homotopy equivalence.
The total space of the nontrivial $RP^2 $-bundle over $S^2 $ is the quotient of
$S^2 \tilde \times S^2 $ by the bundle involution which is the antipodal map on each fibre.
If we observe that $S^2 \tilde \times S^2 \cong CP^2\sharp -CP^2$ is the union 
of two copies of the $D^2 $-bundle which is the mapping cone 
of the Hopf fibration and that this involution 
interchanges the hemispheres we see that this space is homeomorphic
to $RP^4\sharp CP^2 $.

There are two $RP^2 $-bundles over $RP^2 $.
(The total spaces of each of the latter bundles 
have fundamental group $(Z/2Z)^2 $, 
since $w_1 :\pi\to\pi_1 (RP^2)=Z/2Z$ restricts nontrivially to the fibre,
and so is a splitting homomorphism for the homomorphism 
induced by the inclusion of the fibre.)
They may be distinguished by their orientation double covers. 
Each is geometric, and $RP^2\tilde\times{RP^2}$
is also the total space of $S^2$-orbifold bundles
over $\mathbb{D}(2)$ and $RP^2(2)$.

\section{Cohomology and Stiefel-Whitney classes}

We shall show that if $M$ is a closed connected 4-manifold 
with fundamental group $\pi$ such that $\chi(M)|\pi|=4$ 
then $H^*(M;\mathbb{F}_2)$ is isomorphic 
to the cohomology ring of one of the above bundle spaces, 
as a module over the Steenrod algebra $\mathcal A_2$.
(In other words, 
there is an isomorphism which preserves Stiefel-Whitney classes.) 
This is an exercise in Poincar\'e duality and the Wu formulae.

The classifying map induces an isomorphism 
$H^1(\pi;\mathbb{F}_2)\cong H^1(M;\mathbb{F}_2)$ and a
monomorphism $H^2(\pi;\mathbb{F}_2)\to H^2(M;\mathbb{F}_2)$.
If $\pi=1$ then $M$ is homotopy equivalent to $S^2\times S^2$, 
$S^2\tilde\times S^2$ or $CP^2\sharp CP^2$, and the result is clear.

\smallskip
\noindent $\pi=Z/2Z$.\qua In this case $\beta_2(M;\mathbb{F}_2)=2$.
Let $x$ generate $H^1(M;\mathbb{F}_2)$. Then $x^2\not=0$, so $H^2(M;\mathbb{F}_2)$ has a basis $\{ x^2,u\}$.
If $x^4=0$ then $x^2u\not=0$, by Poincar\'e duality, and so $H^3(M;\mathbb{F}_2)$ is generated by $xu$.
Hence $x^3=0$, for otherwise $x^3=xu$ and $x^4=x^2u\not=0$.
Therefore $v_2(M)=0$ or $x^2$, and clearly $v_1(M)=0$ or $x$.
Since $x$ restricts to 0 in $\widetilde M$ we must have 
$w_2(\widetilde M)=v_2(\widetilde M)=0$.
(The four possibilities are realized by the four $S^2$-bundles over $RP^2$.)

If $x^4\not=0$ then we may assume that $x^2u=0$ and that $H^3(M;\mathbb{F}_2)$ is generated by $x^3$.
In this case $xu=0$. Since $Sq^1(x^3)=x^4$ we have $v_1(M)=x$, and $v_2(M)=u+x^2$.
In this case $w_2(\widetilde M)\not=0$, since $w_2(M)$ is not a square. 
(This possibility is realized by the nontrivial $RP^2$-bundle over $S^2$.)

In each case, $xSq^1u=Sq^1(xu)+x^2u=0$, and so $Sq^1u=0$.

\smallskip
\noindent $\pi\cong(Z/2Z)^2$.\qua In this case 
$\beta_2(M;\mathbb{F}_2)=3$ and $w_1(M)\not=0$.
Fix a basis $\{ x,y\}$ for $H^1(M;\mathbb{F}_2)$. 
Then $\{ x^2,xy,y^2\}$ is a basis for $H^2(M;\mathbb{F}_2)$,
since $H^2(\pi;\mathbb{F}_2)$ and $H^2(M;\mathbb{F}_2)$ both have dimension 3.
Since the ring $H^*(M;\mathbb{F}_2)$ is generated by $H^1(M;\mathbb{F}_2)$,  
it determines the image of $[M]$ in $H_4(\pi;\mathbb{F}_2)$.
 
If $x^3$, $y^3$ and $(x+y)^3$ are not all distinct, 
we may assume that $x^3=y^3$.
Then $x^4=Sq^1(x^3)=Sq^1(y^3)=y^4$.  
Hence $x^4=y^4=0$ and $x^2y^2\not=0$, by the nondegeneracy of cup product on
$H^2(M;\mathbb{F}_2)$. 
Hence $x^3=y^3=0$ and so $H^3(M;\mathbb{F}_2)$ is generated by $\{ x^2y,xy^2\}$.
Now $Sq^1(x^2y)=x^2y^2$ and $Sq^1(xy^2)=x^2y^2$, so $v_1(M)=x+y$.
Also $Sq^2(x^2)=0=x^2xy$, $Sq^2(y^2)=0=y^2xy$ and $Sq^2(xy)=x^2y^2$, so $v_2(M)=xy$.
Since the restrictions of $x$ and $y$ to the orientation cover $M^+$ agree we
have $w_2(M^+)=x^2\not=0$.                                                
(This possibility is realized by $RP^2\times RP^2$.)
                                                        
If $x^3$, $y^3$ and $(x+y)^3$ are all distinct
then we may assume that (say) $y^3$ and $(x+y)^3$ generate $H^3(M;\mathbb{F}_2)$.
If $x^3\not=0$ then $x^3=y^3+(x+y)^3=x^3+x^2y+xy^2$ and so $x^2y=xy^2$.
But then we must have $x^4=y^4=0$, by the nondegeneracy of cup product on $H^2(M;\mathbb{F}_2)$. 
Hence $Sq^1(y^3)=y^4=0$ and $Sq^1((x+y)^3)=(x+y)^4=x^4+y^4=0$, and so $v_1(M)=0$,
which is impossible, as $M$ is nonorientable.
Therefore $x^3=0$ and so $x^2y^2\not=0$.
After replacing $y$ by $x+y$, if necessary, we may assume $xy^3=0$ (and hence $y^4\not=0$).
Poincar\'e duality and the Wu relations give $v_1(M)=x+y$, $v_2(M)=xy+x^2$
and hence $w_2(M^+)=0$. 
(This possibility is realized by $RP^2\tilde\times{RP^2}$.)

\smallskip
\noindent $\pi=Z/4Z$.\qua In this case 
$\beta_2(M;\mathbb{F}_2)=1$ and $w_1(M)\not=0$. 
Similar arguments give
$H^*(M;\mathbb{F}_2)\cong\mathbb{F}_2[w,x,y]/(w^2,wx,x^2+wy,xy,x^3,y^2)$,
where $w=w_1(M)$, $x$ has degree 2 and $y$ has degree 3.
Hence $Sq^1x=0$, $Sq^1y=wy$ and $v_2(M)=x$.
 
In all cases, 
a class $x\in H^1(M;\mathbb{F}_2)$ such that $x^3=0$ 
may be realized by a map from $M$ to $K(Z/2Z,1)=RP^\infty$ which 
factors through $P_2(RP^2)$,
since $k_1(RP^2)$ generates 
$H^3(Z/2Z;\pi_2(RP^2))\cong{H^3(Z/2Z;\mathbb{F}_2)}=Z/2Z$. 
However there are such 4-manifolds which do not fibre over $RP^2$.

\section{The action of $\pi$ on $\pi_2 (M)$}

Let $M$ be a closed 4-manifold with finite fundamental group $\pi$ and 
orientation character $w=w_1(M)$.
The intersection form $S(\widetilde M)$ on 
$\Pi=\pi_2(M)=H_2 (\widetilde M;\mathbb{Z})$
is unimodular and symmetric, and $\pi$ acts $w$-isometrically
(that is, $S(ga,gb)=w(g)S(a,b)$ for all $g\in\pi$ and $a$, $b\in\Pi$).

The two inclusions of $S^2 $ as factors of $S^2 \times S^2 $ determine the 
standard basis for $\pi_2 (S^2 \times S^2 )$.                                                                         
Let $J=\left(\begin{smallmatrix}
0&1\\
1&0
\end{smallmatrix}\right)$
be the matrix of the intersection 
form $\bullet$ on $\pi_2 (S^2 \times S^2 )$, with respect to this basis.
The group $Aut(\pm\bullet)$ of automorphisms of $\pi_2(S^2\times S^2)$ which 
preserve this intersection form up to sign is the dihedral group of order eight,
and is generated by the diagonal matrices and $J$ or 
$K=
\left(\begin{smallmatrix}
0&1\\
-1&0
\end{smallmatrix}\right)$.
The subgroup of strict isometries has order four, 
and is generated by $-I$ and $J$.
(Note that the isometry $J$ is induced by the involution $\tau$.)
                                                                                
Let $f$ be a self homeomorphism of $S^2 \times S^2 $ and let $f_* $ be the 
induced automorphism of $\pi_2(S^2\times S^2)$. 
The Lefshetz number of $f$ is $2+trace(f_*)$ if $f$ is orientation preserving 
and $trace(f_* )$ if $f$ is orientation reversing. 
As any self homotopy equivalence which induces the identity on $\pi_2 $ 
has nonzero Lefshetz number the natural representation of a group $\pi$ 
of fixed point free self homeomorphisms of $S^2\times S^2$ into 
$Aut(\pm\bullet)$ is faithful. 

Suppose first that $f$ is a free involution, so $f_*^2=I$.
If $f$ is orientation preserving then $trace(f_*)=-2$ so $f_* =-I$.
If $f$ is orientation reversing then $trace(f_*)=0$, 
so $f_*=\pm JK=\pm 
\left(\begin{smallmatrix}
1&0\\
0&-1
\end{smallmatrix}\right)$.
Note that if $f'=\tau f\tau$ then $f'_* =-f_* $, so after conjugation by
$\tau$, if necessary, we may assume that $f_* =JK$. 

If $f$ generates a free $Z/4Z$-action the induced automorphism must be $\pm K$.
Note that if $f'=\tau f\tau$ then $f'_* =-f_* $, so after conjugation by
$\tau$, if necessary, we may assume that $f_* =K$.            
       
Since the orbit space of a fixed point free action of $(Z/2Z)^2 $ 
on $S^2 \times S^2 $ 
has Euler characteristic 1 it is nonorientable, 
and so the action is generated by two 
commuting involutions, one of which is orientation preserving and one of which is not. 
Since the orientation preserving involution must act via $-I$ and the orientation reversing
involutions must act via $\pm JK$ the action of $(Z/2Z)^2 $ is essentially unique.

The standard inclusions of $S^2=CP^1 $ into the summands of 
$CP^2\sharp -CP^2\cong S^2\tilde\times S^2$
determine a basis for $\pi_2 (S^2 \tilde\times S^2 )\cong\mathbb{Z}^2 $. 
Let $\tilde J=
\left(\begin{smallmatrix}
1&0\\
0&-1
\end{smallmatrix}\right)$ be the matrix of the 
intersection form $\tilde\bullet$ on 
$\pi_2 (S^2\tilde\times S^2)$ with respect to this basis.
The group $Aut(\pm\tilde\bullet)$ of automorphisms of 
$\pi_2 (S^2 \tilde\times S^2)$ which preserve this intersection form up to sign 
is the dihedral group of order eight, 
and is also generated by the diagonal matrices and 
$J=\left(\begin{smallmatrix}
0&1\\
1&0
\end{smallmatrix}\right)$.
The subgroup of strict isometries has order four, 
and consists of the diagonal matrices.  
A nontrivial group of fixed point free self homeomorphisms of 
$S^2 \tilde\times S^2 $ must have order 2,
since $S^2 \tilde\times S^2 $ admits no fixed point free orientation 
preserving involution, by Lemma 12.3.
If $f$ is an orientation reversing free involution of $S^2 \tilde\times S^2 $ 
then $f_*=\pm J$. 
Since the involution of $CP^2$ given by complex conjugation is orientation 
preserving it is isotopic to a selfhomeomorphism $c$ which fixes a 4-disc.
Let $g=c\sharp id_{CP^2}$.
Then $g_*=
\left(\begin{smallmatrix}
-1&0\\
0&1
\end{smallmatrix}\right)$ and so $g_*Jg_*^{-1}=-J$.
Thus after conjugating $f$ by $g$, if necessary, we may assume that $f_*=J$.

All self homeomorphisms of $CP^2\sharp CP^2$ preserve the sign of the 
intersection form, and thus are orientation preserving. 
With part (2) of Lemma 12.3, 
this implies that no manifold in this homotopy type 
admits a free involution.

\section{Homotopy type}

The {\it quadratic $2$-type} of $M$ is the quadruple 
$[\pi,\pi_2 (M),k_1 (M),S(\widetilde M)]$.
Two such quadruples $[\pi,\Pi,\kappa,S]$ and $[\pi',\Pi',\kappa',S']$ 
with $\pi$ a finite group, $\Pi$ a finitely generated, 
$\mathbb{Z}$-torsion-free $\mathbb{Z}[\pi]$-module, 
$\kappa\in H^3(\pi;\Pi)$ and $S:\Pi\times\Pi\to \mathbb{Z}$ a unimodular symmetric bilinear pairing
on which $\pi$ acts $\pm$-isometrically are {\it equivalent} if there is an isomorphism 
$\alpha:\pi\to\pi'$ and an (anti)isometry 
$\beta:(\Pi,S)\to(\Pi',(\pm)S')$ which is $\alpha$-equivariant 
(i.e., such that $\beta(gm)=\alpha(g)\beta(m)$ for all $g\in\pi$ and $m\in\Pi$)
and $\beta_* \kappa=\alpha^* \kappa'$ in $H^3 (\pi,\alpha^* \Pi' )$.   
Such a quadratic 2-type determines homomorphisms 
$w:\pi\to\mathbb{Z}^\times=Z/2Z$ (if $\Pi\not=0$)
and $v:\Pi\to Z/2Z$ by the equations $S(ga,gb)=w(g)S(a,b)$ and
$v(a)\equiv S(a,a)~mod~(2)$, for all $g\in\pi$ and $a$, $b\in\Pi$.  
(These correspond to the orientation character $w_1(M)$ and the Wu class 
$v_2(\widetilde M)=w_2(\widetilde M)$, of course.)

Let $\gamma:A\to\Gamma(A)$ be the universal quadratic functor of Whitehead.
Then the pairing $S$ may be identified with an indivisible element of
$\Gamma(Hom_{\mathbb{Z}}(\Pi,\mathbb{Z}))$.
Via duality, this corresponds to an element $\widehat S$ of $\Gamma(\Pi)$, 
and the subgroup generated by the image of $\widehat S$ is a 
$\mathbb{Z}[\pi]$-submodule. 
Hence $\pi_3=\Gamma(\Pi)/\langle \widehat S\rangle$ 
is again a finitely generated, 
$\mathbb{Z}$-torsion-free $\mathbb{Z}[\pi]$-module.
Let $B$ be the Postnikov 2-stage corresponding to the algebraic 2-type 
$[\pi,\Pi,\kappa]$.
A $PD_4$-{\it polarization} of the quadratic 2-type $q=[\pi,\Pi,\kappa,S]$ is a 
3-connected map $f:X\to B$, where $X$ is a $PD_4$-complex, 
$w_1(X)=w\pi_1(f)$ and $\tilde f_*(\widehat S_{\widetilde X})=\widehat S$ in 
$\Gamma(\Pi)$.
Let $S_4^{PD}(q)$ be the set of equivalence classes of $PD_4$-polarizations 
of $q$, where $f:X\to B\sim g:Y\to B$ if there is a map $h:X\to Y$ 
such that $f\simeq gh$.

\begin{theorem}
{\rm[Te]}\qua There is an effective, transitive action of               
the torsion subgroup of 
$\mathbb{Z}^w\otimes_{\mathbb{Z}[\pi]}\Gamma(\Pi)$ 
on $S_4^{PD}(q)$.
\end{theorem}

\begin{proof} 
(We shall only sketch the proof.) 
Let $f:X\to B$ be a fixed $PD_4$-polarization of $q$.
We may assume that $X=K\cup_g e^4$, where $K=X^{[3]}$ is the 3-skeleton and $g\in\pi_3(K)$
is the attaching map.
Given an element $\alpha$ in $\Gamma(\Pi)$ whose image in
$\mathbb{Z}^w\otimes_{\mathbb{Z}[\pi]}\Gamma(\Pi)$ 
lies in the torsion subgroup, 
let $X_\alpha=K\cup_{g+\alpha} e^4$.
Since $\pi_3(B)=0$ the map $f|_K$ extends to a map $f_\alpha:X_\alpha\to B$,
which is again a $PD_4$-polarization of $q$.
The equivalence class of $f_\alpha$ depends only on the image of $\alpha$ in
$\mathbb{Z}^w\otimes_{\mathbb{Z}[\pi]}\Gamma(\Pi)$.
Conversely, if $g:Y\to B$ is another $PD_4$-polarization of $q$ then
$f_*[X]-g_*[Y]$ lies in the image of
$Tors(\mathbb{Z}^w\otimes_{\mathbb{Z}[\pi]}\Gamma(\Pi))$ 
in $H_4(B;\mathbb{Z}^w)$.
See \S2 of \cite{[Te]} for the full details.
\end{proof} 

\begin{cor}
If $X$ and $Y$ are $PD_4$-complexes with the same quadratic $2$-type
then each may be obtained by adding a single $4$-cell to $X^{[3]}=Y^{[3]}$.
\qed
\end{cor}

If $w=0$ and the Sylow 2-subgroup of $\pi$ has cohomological period
dividing 4 then 
$Tors(\mathbb{Z}^w\otimes_{\mathbb{Z}[\pi]}\Gamma(\Pi))=0$ \cite{[Ba88]}.
In particular, if $M$ is orientable and $\pi$ is finite cyclic 
then the equivalence class of the quadratic 2-type determines 
the homotopy type \cite{[HK88]}.
Thus in all cases considered here the quadratic 2-type determines 
the homotopy type of the orientation cover. 

The group $Aut(B)=Aut([\pi,\Pi,\kappa])$ acts on $S_4^{PD}(q)$ and the orbits 
of this action correspond to the homotopy types of $PD_4$-complexes $X$ 
admitting such polarizations $f$.
When $q$ is the quadratic 2-type of $RP^2\times RP^2$ this action 
is nontrivial.
(See below in this section. Compare also Theorem 10.5.)

The next lemma shall enable us to determine the possible $k$-invariants.

\begin{lemma}
Let $M$ be a closed $4$-manifold with fundamental group
$\pi=Z/2Z$ and universal covering space $S^2 \times S^2 $. 
Then the first $k$-invariant of $M$
is a nonzero element of $H^3 (\pi;\pi_2 (M))$.
\end{lemma}

\begin{proof} 
The first $k$-invariant is the primary obstruction 
to the existence of a cross-section to the classifying map 
$c_M:M\to K(Z/2Z,1)=RP^\infty$ and is the only obstruction 
to the existence of such a cross-section for $c_{P_2(M)}$.
The only nonzero differentials in the Cartan-Leray cohomology spectral sequence 
(with coefficients $Z/2Z$) for the projection $p:\widetilde M\to M$ are at the 
$E_3 ^{**} $ level.
By the results of Section 4, $\pi$ acts trivially on 
$H^2 (\widetilde M;\mathbb{F}_2)$, since $\widetilde M=S^2 \times S^2 $.
Therefore $E_3^{22} =E_2^{22}\cong (Z/2Z)^2 $ and $E_3^{50} =E_2^{50}=Z/2Z$.
Hence $E_\infty^{22}\not=0$, so $E_\infty^{22} $ maps onto 
$H^4 (M;\mathbb{F}_2)=Z/2Z$ and 
$d_3^{12} :H^1 (\pi;H^2 (\widetilde M;\mathbb{F}_2))\to H^4 (\pi;\mathbb{F}_2)$ 
must be onto.
But in this region the spectral sequence is identical with the corresponding
spectral sequence for $P_2 (M)$. 
It follows that the image of $H^4 (\pi;\mathbb{F}_2)=Z/2Z$ in 
$H^4 (P_2 (M);\mathbb{F}_2)$
is 0, and so $c_{P_2(M)}$ does not admit a cross-section. Thus $k_1 (M)\not=0$. 
\end{proof} 

If $\pi=Z/2Z$ and $M$ is orientable then $\pi$ acts via $-I$ on 
$\mathbb{Z}^2$ and the
$k$-invariant is a nonzero element of $H^3(Z/2Z;\pi_2(M))=(Z/2Z)^2$.
The isometry which transposes the standard generators of 
$\mathbb{Z}^2$ is $\pi$-linear, 
and so there are just two equivalence classes of quadratic 2-types to consider.
The $k$-invariant which is invariant under transposition is realised by 
$(S^2\times S^2)/\langle(-I,-I)\rangle$, while the other $k$-invariant is realized by 
the other orientable $S^2$-bundle space.
Thus $M$ must be homotopy equivalent to one of these spaces.              
                                                                     
If $\pi=Z/2Z$, $M$ is nonorientable and $w_2(\widetilde M)=0$ then 
$H^3(\pi;\pi_2(M))=Z/2Z$ and there is only one quadratic 2-type to consider.
There are four equivalence classes of $PD_4$-polarizations,
as 
$Tors(\mathbb{Z}^w\otimes_{\mathbb{Z}[\pi]}\Gamma(\Pi))\cong (Z/2Z)^2$.
The corresponding $PD_4$-complexes have the form $K\cup_f e^4$, 
where $K={(S^2\times RP^2)\setminus{int}D^4}$ is the 3-skeleton 
of $S^2\times{RP^2}$ and $f\in\pi_3(K)$. 
Two choices for $f$ give total spaces of $S^2$-bundles over $RP^2$. 
A third choice gives $RP^4\sharp_{S^1} RP^4=E\cup_\tau{E}$, 
where $E=S^2\times{D^2}/(s,d)\sim(-s,-d)$ and $\tau$ 
is the twist map of $\partial{E}=S^2\tilde\times{S^1}$ \cite{[KKR92]}.
This is the total space of an orbifold bundle with general fibre $S^2$ 
over $S(2,2)$, but is not geometric.
The fourth homotopy type has nontrivial Browder-Livesay invariant,
and so is not realizable by a closed manifold \cite{[HM78]}.
The product $S^2\times{RP^2}$ is characterized by the additional 
conditions that $w_2(M)=w_1(M)^2\not=0$ (i.e., $v_2(M)=0$) and that there 
is an element $u\in H^2 (M;\mathbb{Z})$ which generates 
an infinite cyclic direct summand and is such that $u\cup u=0$. 
(See Theorem 5.19.)
The nontrivial nonorientable $S^2$-bundle space has $w_2(M)=0$.
The manifold $RP^4\sharp_{S^1} RP^4 $ also has $w_2(M)=0$, 
but it may be distinguished from the bundle space by the $Z/4Z$-valued
quadratic function on $\pi_2 (M)\otimes (Z/2Z)$ introduced in 
\cite{[KKR92]}.

If $\pi=Z/2Z$ and $w_2 (\widetilde M)\not=0$ then $H^3 (\pi_1;\pi_2 (M))=0$,
and the quadratic 2-type is unique. 
(Note that the argument of Lemma 12.5 breaks down here 
because $E_\infty^{22} =0$.)        
There are two equivalence classes of $PD_4$-polarizations,
as $Tors(\mathbb{Z}^w\otimes_{\mathbb{Z}[\pi]}\Gamma(\Pi))=Z/2Z$.
They are each of the form $K\cup_f e^4$, 
where $K=(RP^4\sharp CP^2)\setminus{int}D^4$ is the 3-skeleton 
of $RP^4\sharp CP^2$ and $f\in\pi_3(K)$.
The bundle space $RP^4\sharp CP^2$ is characterized by the additional condition 
that there is an element $u\in H^2 (M;\mathbb{Z})$ which generates an infinite 
cyclic direct summand and such that $u\cup u=0$. (See Theorem 5.19.)                   
In \cite{[HKT94]} it is shown that any closed 4-manifold $M$ with $\pi=Z/2Z$, 
$\chi(M)=2$ and $w_2(\widetilde M)\not=0$ is homotopy equivalent to 
$RP^4\sharp CP^2$.  

If $\pi\cong Z/4Z$ then 
$H^3(\pi;\pi_2(M))\cong\mathbb{Z}^2/(I-K)\mathbb{Z}^2=Z/2Z$, 
since $\Sigma_{k=1}^{k=4} f_*^k=\Sigma_{k=1}^{k=4}K^k=0$.
The $k$-invariant is nonzero, since it restricts to the $k$-invariant of 
the orientation double cover.
In this case $Tors(\mathbb{Z}^w\otimes_{\mathbb{Z}[\pi]}\Gamma(\Pi))=0$ 
and so $M$ is homotopy equivalent to 
$(S^2\times S^2)/\langle\tau(I,-I)\rangle$.

Finally, let $\pi\cong(Z/2Z)^2 $ be the diagonal subgroup of 
$Aut(\pm\bullet)<GL(2,\mathbb{Z})$,
and let $\alpha$ be the automorphism induced by conjugation by $J$.
The standard generators of $\pi_2 (M)=\mathbb{Z}^2 $ 
generate complementary $\pi$-submodules, 
so that $\pi_2 (M)$ is the direct sum 
$\tilde{\mathbb{Z}}\oplus\alpha^*\tilde{\mathbb{Z}}$ 
of two infinite cyclic modules.
The isometry $\beta=J$ which transposes the factors is $\alpha$-equivariant, 
and $\pi$ and $V=\{\pm I\} $ act nontrivially on each summand. 
If $\rho$ is the kernel of the action of $\pi$ on $\tilde{\mathbb{Z}}$ 
then $\alpha(\rho)$ is the kernel 
of the action on $\alpha^* \tilde{\mathbb{Z}}$, and $\rho\cap\alpha(\rho)=1$.
Let $j_V:V\to\pi$ be the inclusion.
As the projection of $\pi=\rho\oplus V$ onto $V$ is compatible with the action, 
$H^*(j_V;\tilde{\mathbb{Z}})$ is a split epimorphism and so
$H^* (V;\tilde{\mathbb{Z}})$ is a direct summand of 
$H^* (\pi;\tilde{\mathbb{Z}})$.
This implies in particular that the differentials in the LHSSS
$H^p (V;H^q (\rho;\tilde{\mathbb{Z}}))\Rightarrow 
H^{p+q} (\pi;\tilde{\mathbb{Z}})$ 
which end on the row $q=0$ are all 0.
Hence $H^3(\pi;\tilde{\mathbb{Z}})\cong 
H^1(V;\mathbb{F}_2)\oplus H^3(V;\tilde{\mathbb{Z}})\cong (Z/2Z)^2$.
Similarly $H^3 (\pi;\alpha^* \tilde{\mathbb{Z}})\cong (Z/2Z)^2 $, 
and so $H^3 (\pi;\pi_2 (M))\cong (Z/2Z)^4 $.
The $k$-invariant must restrict to the $k$-invariant of each double cover, 
which must be nonzero, by Lemma 12.5.
Let $K_V$, $K_\rho$ and $K_{\alpha(\rho)}$ be the kernels of the restriction
homomorphisms from $H^3(\pi;\pi_2(M))$ to $H^3(V;\pi_2(M))$, $H^3(\rho;\pi_2(M))$
and $H^3(\alpha(\rho);\pi_2(M))$, respectively.
Now $H^3 (\rho;\tilde{\mathbb{Z}})=
H^3 (\alpha(\rho);\alpha^* \tilde{\mathbb{Z}})=0$,
$H^3 (\rho;\alpha^* \tilde{\mathbb{Z}})=
H^3 (\alpha(\rho);\tilde{\mathbb{Z}})=Z/2Z$ and
$H^3 (V;\tilde{\mathbb{Z}})=H^3 (V;\alpha^* \tilde{\mathbb{Z}})=Z/2Z$.
Since the restrictions are epimorphisms $|K_V|=4$ and 
$|K_\rho|=|K_{\alpha(\rho)}|=8$.
It is easily seen that ${|K_\rho\cap K_{\alpha(\rho)}|}=4$.
Moreover 
Ker$(H^3(j_V;\tilde{\mathbb{Z}}))\cong H^1(V;H^2(\rho;\tilde{\mathbb{Z}}))\cong 
H^1(V;H^2(\rho;\mathbb{F}_2))$ 
restricts nontrivially to 
$H^3(\alpha(\rho);\tilde{\mathbb{Z}})\cong H^3(\alpha(\rho);\mathbb{F}_2)$,
as can be seen by reduction modulo (2),
and similarly $\mathrm{Ker}(H^3(j_V;\alpha^*\tilde{\mathbb{Z}}))$ 
restricts nontrivially to $H^3(\rho;\alpha^*\tilde{\mathbb{Z}})$.
Hence $|K_V\cap K_\rho|=|K_V\cap K_{\alpha(\rho)}|=2$ and 
$K_V\cap K_\rho\cap K_{\alpha(\rho)}=0$.
Thus $|K_V\cup K_\rho\cup K_{\alpha(\rho)}|=8+8+4-4-2-2+1=13$ and so
there are at most three possible $k$-invariants.
Moreover the automorphism $\alpha$ and the isometry $\beta=J$ act on the 
$k$-invariants by transposing the factors. 
The $k$-invariant of $RP^2\times RP^2$ is invariant under this transposition,
while that of $RP^2\tilde\times{RP^2}$ is not,
for the $k$-invariant of its orientation cover is not invariant.
Thus there are two equivalence classes of quadratic 2-types to be considered.
Since
$Tors(\mathbb{Z}^w\otimes_{\mathbb{Z}[\pi]}\Gamma(\Pi))\cong (Z/2Z)^2$
each has four equivalence classes of $PD_4$-polarizations.
In each case the quadratic 2-type determines the 
$\mathbb{F}_2$-cohomology ring, 
since it determines the orientation cover (see \S4).
The canonical involution of the direct product interchanges two of these
polarizations in the $RP^2\times RP^2$ case, 
and so there are seven homotopy types of $PD_4$-complexes to consider.

In \cite{[HH21]} the arguments of \cite{[HM78]} are extended to show that $RP^2\times{RP^2}$ is the only closed 4-manifold in its quadratic 2-type, 
while there are two closed 4-manifolds in the quadratic 2-type of
$RP^2\tilde\times{RP^2}$.
Thus there are three homotopy types of closed 4-manifolds $M$ 
with $\pi\cong(Z/2Z)^2$ and $\chi(M)=1$.
Two of these are represented by the bundle spaces.
The third homotopy type is represented by the quotient of $RP^4\sharp_{S^1}RP^4=E\cup_\tau{E}$ by the free involution 
which sends $[s,d]$ in one copy of $E$ to $[-s,d]$.

\section{Surgery}

We may assume that $M$ is a proper quotient of $S^2\times S^2$ or of
$S^2\tilde\times S^2$, so $|\pi|\chi(M)=4$ and $\pi\not=1$.
In the present context every homotopy equivalence is simple since $Wh(\pi)=0$ 
for all groups $\pi$ of order $\leq 4$ \cite{[Hg40]}.

Suppose first that $\pi=Z/2Z$. 
Then $H^1(M;\mathbb{F}_2)=Z/2Z$ and $\chi(M)=2$, 
so $H^2 (M;\mathbb{F}_2)\cong (Z/2Z)^2$. 
The $\mathbb{F}_2$-Hurewicz homomorphism from $\pi_2 (M)$ 
to $H_2 (M;\mathbb{F}_2)$ has cokernel $H_2 (\pi;\mathbb{F}_2)=Z/2Z$. 
Hence there is a map $\beta:S^2 \to M$ such that $\beta_* [S^2]\not=0$ in 
$H_2 (M;\mathbb{F}_2)$.        
If moreover $w_2 (\widetilde M)=0$ then $\beta^* w_2 (M)=0$, 
since $\beta$ factors through $\widetilde M$.
Then there is a self homotopy equivalence $f_\beta$ of $M$ with nontrivial 
normal invariant in $[M;G/TOP]$, by Lemma 6.5.
Note also that $M$ is homotopy equivalent to a PL 4-manifold (see \S 6 above).

If $M$ is orientable $[M;G/TOP]\cong\mathbb{Z}\oplus (Z/2Z)^2 $.
The surgery obstruction groups are $L_5 (Z/2Z,+)=0$ and 
$L_4 (Z/2Z,+)\cong\mathbb{Z}^2 $,                            
where the surgery obstructions are determined by the signature 
and the signature of the double cover \cite[Theorem 13.A.1]{[Wl]}.
Hence it follows from the surgery exact sequence that 
$S_{TOP} (M)\cong (Z/2Z)^2$.
Since $w_2 (\widetilde M)=0$ (by Lemma 12.3) there is 
a self homotopy equivalence $f_\beta$ of $M$ with nontrivial 
normal invariant and so there are at most two 
homeomorphism classes within the homotopy type of $M$. 
If $w_2(M)=0$ then $KS(M)=0$, since $\sigma(M)=0$,
and so $M$ is homeomorphic to the bundle space \cite{[Te97]}.
On the other hand if $w_2(M)\not=0$ there is an 
$\alpha\in H^2(M;\mathbb{F}_2)$ such that $\alpha^2\not=0$,
and $\alpha=kerv(\hat{f})$ for some homotopy equivalence $f:N\to{M}$.
We then have $KS(N)=f^*(KS(M)+\alpha^2)$ (see \S2 of Chapter 6 above),
and so $KS(N)\not=KS(M)$.
Thus there are three homeomorphism classes of orientable closed 4-manifolds 
with $\pi=Z/2Z$ and $\chi=2$.
One of these is a fake $(S^2\times S^2)/\langle(-I,-I)\rangle$ and is not smoothable.

Suppose now that $M$ is non-orientable, with orientable double cover $M^+$.
The natural maps from $L_4(1)$ to $L_4(Z/2Z,-)$ and 
$L_4((Z/2Z)^2,-)$ are trivial, while $L_4(Z/4Z,-)=0$ 
\cite{[Wl76]}. 
Thus we may change the value of $KS(M)$ at will,
by surgering the normal map $M\sharp{E_8}\to M\sharp{S^4}\cong{M}$.
In all cases $KS(M^+)=0$, and so $M^+$ is the total space 
of an $S^2$-bundle over $S^2$ or $RP^2$.

Nonorientable closed 4-manifolds with fundamental group $Z/2Z$ have been
classified in \cite{[HKT94]}. 
The surgery obstruction groups are $L_5 (Z/2Z,-)=0$ and $L_4 (Z/2Z,-)=Z/2Z$,
and $\sigma_4(\hat g)=c(\hat g)$ for any normal map 
$\hat g:M\to{G/TOP}$ \cite[Theorem 13.A.1]{[Wl]}. 
Therefore $\sigma_4(\hat g)=(w_1(M)^2\cup \hat g^*(k_2))[M]$ 
\cite[Theorem 13.B.5]{[Wl]}.
(See also \S2 of Chapter 6 above.)
As $w_1(M)$ is not the reduction of a class in $H^1(M;\mathbb{Z}/4\mathbb{Z})$ 
its square is nonzero and so there is an element $\hat g^*(k_2)$
in $H^2(M;\mathbb{F}_2)$ such that this cup product is nonzero.
Hence $S_{TOP}(M)\cong (Z/2Z)^2$, since $[M;G/TOP]\cong (Z/2Z)^3 $.
There are two homeomorphism types within each homotopy type if 
$w_2 (\widetilde M)=0$; 
if $w_2 (\widetilde M)\not=0$ 
there are four corresponding homeomorphism types \cite{[HKT94]}:
$RP^4\sharp CP^2$ (the nontrivial $RP^2$-bundle over $S^2$), 
$RP^4\sharp *CP^2$, $*RP^4\sharp CP^2$, 
which have nontrivial Kirby-Siebenmann invariant,
and $(*RP^4)\sharp*CP^2$, which is smoothable \cite{[RS97]}.
Thus there are ten homeomorphism classes of nonorientable closed 4-manifolds
with $\pi=Z/2Z$ and $\chi=2$.

The image of $[M;G/PL]$ in $[M;G/TOP]$ is a subgroup of index 2 
(see \cite[\S15]{[Si71]}). 
It follows that if $M$ is the total space of an $S^2 $-bundle over $RP^2 $ 
any homotopy equivalence $f:N\to M$ where $N$ is also PL is homotopic 
to a homeomorphism. 
(For then $S_{TOP}(M)\cong (Z/2Z)^2$, 
and the nontrivial element of the image 
of $S_{PL} (M)$ is represented by a self homotopy equivalence of $M$. 
The case $M=S^2 \times RP^2$ was treated in \cite{[Ma79]}. 
See also \cite{[Te97]} for the cases with $\pi=Z/2Z$ and $w_1(M)=0$.)           
This is also true if $M=S^4$, $RP^4$, $CP^2$, $S^2\times S^2$ or 
$S^2\tilde\times S^2$.

If $\pi\cong Z/4Z$ or $(Z/2Z)^2 $ then $\chi(M)=1$, 
and so $M$ is nonorientable.
As the $\mathbb{F}_2$-Hurewicz homomorphism is 0 in these cases Lemma 6.5
does not apply to give any exotic self homotopy equivalences.

If $\pi\cong (Z/2Z)^2$ then $[M;G/TOP]\cong (Z/2Z)^4 $ and the surgery 
obstruction groups are $L_5 ((Z/2Z)^2,-)=0$ and $L_4 ((Z/2Z)^2,-)=Z/2Z$ 
\cite[Theorem 3.5.1 ]{[Wl76]}.
Since $w_1(M)$ is a split epimorphism $L_4(w_1(M))$ is an isomorphism,
so the surgery obstruction is detected by the Kervaire-Arf invariant.
As $w_1(M)^2\not=0$ we find that $S_{TOP} (M)\cong(Z/2Z)^3$.
Thus there are between six and 24 homeomorphism classes 
of closed 4-manifolds with $\pi\cong (Z/2Z)^2$ and $\chi=1$, 
of which half are not stably smoothable.

If $\pi\cong Z/4Z$ then $[M;G/TOP]\cong(Z/2Z)^2$ 
and the surgery obstruction groups $L_4 (Z/4Z,-)$ 
and $L_5 (Z/4Z,-)$ are both 0 \cite[Theorem 3.4.5]{[Wl76]}. 
Hence $S_{TOP}(M)\cong(Z/2Z)^2$, and so
there are at most four homeomorphism classes 
of closed 4-manifolds with $\pi\cong Z/4Z$ and $\chi=1$.
The original version of \cite{[HH21]} 
relied on an argument for the assertion (*) 
{\sl all self homotopy equivalences of $M$ are homotopic to self homeomorphisms},
which was later felt to have a gap (in the identification of the $\mathbb{Z}[\pi]$-module structure for an entry in a spectral sequence).
This argument has been replaced by one based on first using Kreck's modified surgerey to determine a stable homeomorphism classification, 
and then using surgery to show that each stable homeomorphism class
contains a representative with Euler characteristic 1.
The conclusion is that the assertion (*) is correct, and that there are 4 homeomorphism types of TOP 4-manifolds in the homotopy type of $M$,
two of which are stably smoothable.

In the final section we give a candidate for another PL 4-manifold 
in this homotopy type.
(This is from \cite{[HH21]}.)

\section{A smooth fake version of $S^2\times{S^2}/\langle\tau(I,-I)\rangle$?}

Let $M=S^2\times{S^2}/\langle\tau(I,-I)\rangle$, as in the previous section,
and let $M^+$ be its orientable double cover.
Let $\Delta=\{(s,s)\mid{s}\in{S^2}\}$ be the diagonal in $S^2\times{S^2}$.
We may isotope $\Delta$ to a nearby sphere 
which meets $\Delta$ transversely in two points,
by rotating the first factor, 
and so $\Delta$ has self-intersection $\pm2$.
The diagonal is invariant under $(-I,-I)$, 
and so $\delta=\Delta/\langle\sigma^2\rangle\cong{RP^2}$ embeds in 
$M^+=S^2\times{S^2}/\langle(-I,-I)\rangle$ 
with an orientable regular neighbourhood.
Since $\sigma(\Delta)\cap\Delta=\emptyset$ this also embeds in $M$.
We shall see that the complementary region also has a simple description.

We shall view $S^2$ here as the unit sphere in the space
of purely imaginary quaternions.
The standard inner product on this space is given by
$v\bullet{w}=\mathfrak{Re}(v\bar{w})$,
for $v,w$ purely imaginary quaternions.
We shall also identify $S^3$ with the unit quaternions $\mathbb{H}_1$.
Let $C_x=\{(s,t)\in{S^2\times{S^2}}\mid{s\bullet{t}=x}\}$, 
for all $x\in[-1,1]$.
Then $C_1=\Delta$ and $C_{-1}=\tau(I,-I)(\Delta)$,
while $C_x\cong{C_0}$ for all $|x|<1$.
The map
$p:S^3\to{C_0}$ given by $p(q)=(q\mathbf{i}q^{-1},q\mathbf{j}q^{-1})$
for all $q\in{S^3}$ is a 2-fold covering projection,
and so $C_0\cong{RP^3}$.

It is easily seen that $N=\cup_{x\geq\varepsilon}C_x$ and $\tau(I,-I)(N)$ 
are regular neighbourhoods of $\Delta$ and $\tau(I,-I)(\Delta)$, 
respectively, while $C=\cup_{x\in[-\varepsilon,\varepsilon]}C_x\cong
{C_0\times[-\varepsilon,\varepsilon]}$.
In particular, $N$ and $\tau(I,-I)(N)$ are each homeomorphic to the total space 
of the unit disc bundle in $T_{S^2}$, 
and $\partial{N}\cong{C_0}\cong{RP^3}$. 
The subsets $C_x$ are invariant under $(-I,-I)$.
Hence $N(\delta)=N/\langle(-I,-I)\rangle$ is the total space 
of the tangent disc bundle of $RP^2$.
In particular, $\partial{N(\delta)}\cong{L(4,1)}$ and
$\delta$ represents the nonzero element of $H_2(M;\mathbb{F}_2)$, 
since it has self-intersection 1 in $\mathbb{F}_2$.
(It is not hard to show that any embedded surface representing
the nonzero element of $H_2(M;\mathbb{F}_2)$ is non-orientable
but lifts to $M^+$, and so has an orientable regular neighbourhood.) 

We also see that 
$C/\langle(-I,-I)\rangle\cong{L(4,1)}\times[-\varepsilon,\varepsilon]$.
The self map of ${S^3}$ defined by right multiplication
by $\frac1{\sqrt{2}}(\mathbf{1}+\mathbf{k})$ lifts $\tau(I,-I)$,
since $p(q.\frac1{\sqrt{2}}(\mathbf{1}+\mathbf{k}))=\tau(I,-I)(p(q))$.
Hence $C_0/\langle\tau(I,-I)\rangle=
S^3/\langle\frac1{\sqrt{2}}(\mathbf{1}+\mathbf{k})\rangle=L(8,1)$,
and so $MC=C/\langle\sigma\rangle$ is the mapping cylinder of 
the double cover $L(4,1)\to{L(8,1)}$.
Since $S^2\times{S^2}=N\cup{C}\cup\tau(I,-I)(N)$ it follows that 
$M=N(\delta)\cup{MC}$.

This construction suggests a candidate for another smooth 4-manifold 
in the same (simple) homotopy type.
Let $M'=N(\delta)\cup{MC'}$, where $MC'$ is the mapping cylinder of 
the double cover $L(4,1)\to{L(8,5)}$.
Then $\pi_1(M')\cong{Z/4Z}$ and $\chi(M')=1$, and so 
there is a homotopy equivalence $h:M'\simeq{M}$.
In this case, $kerv(h)$ is a complete invariant.
The difficulty is in finding an easily analyzed explicit choice for $h$,
for which we can compute $kerv(h)$.

Is there a computable homeomorphism (or diffeomorphism) invariant 
that can be applied here?
Most readily computable invariants are invariants of homotopy type.
These manifolds have $w_2\not=0$ but $w_1^2=0$ and $w_1w_2=0$, 
so admit $Pin^c$-structures.
Is it possible to distinguish them by consideration of the associated Arf invariants?

If $M$ and $M'$ are not homeomorphic then every closed 4-manifold 
with $\pi\cong{Z/4Z}$ and $\chi=1$ is homeomorphic to one of $M$, 
$M'$, $*M$ or $*M'$.

%% file: m5-13.tex
\chapter{Geometric decompositions of bundle spaces}

We begin by considering which closed 4-manifolds with geometries 
of euclidean factor type are mapping tori of homeomorphisms of 3-manifolds.
We also show that (as an easy consequence of the Kodaira classification of 
surfaces) a complex surface is diffeomorphic to a mapping torus 
if and only if its Euler characteristic is 0 and its fundamental group 
maps onto $\mathbb{Z}$ with finitely generated kernel, 
and we determine the relevant 3-manifolds and diffeomorphisms.
In \S2 we consider when an aspherical 4-manifold which is the total space
of a surface bundle is geometric or admits a geometric decomposition.
If the base and fibre are hyperbolic the only known examples 
are virtually products.
In \S3 we shall give some examples of torus bundles over closed surfaces 
which are not geometric,
some of which admit geometric decompositions of type $\mathbb{F}^4$ 
and some of which do not.
In \S4 we apply some of our earlier results to the characterization 
of certain complex surfaces.
In particular, we show that a complex surfaces fibres smoothly over
an aspherical orientable 2-manifold if and only if 
it is homotopy equivalent to the total space of a surface bundle.
In the final two sections we consider first $S^1$-bundles over geometric 
3-manifolds and then the existence of symplectic structures on geometric
4-manifolds.

\section{Mapping tori}

In \S3-5 of Chapter 8 and \S3 of Chapter 9 we used 3-manifold theory to 
characterize mapping tori of homeomorphisms of geometric 3-manifolds which 
have product geometries. 
Here we shall consider instead which 4-manifolds with product 
geometries or complex structures are mapping tori.

\begin{theorem} 
Let $M$ be a closed geometric $4$-manifold with 
$\chi(M)=0$ and such that $\pi=\pi_1 (M)$ is an extension of $\mathbb{Z}$ 
by a finitely generated normal subgroup $K$.
Then $K$ is the fundamental group of a geometric $3$-manifold.
\end{theorem}

\begin{proof} 
Since $\chi(M)=0$ the geometry must be either an 
infrasolvmanifold geometry or a product geometry $\mathbb{X}^3\times\mathbb{E}^1$, 
where $\mathbb{X}^3$ is one of the 3-dimensional geometries 
$\mathbb{S}^3$, $\mathbb{S}^2\times\mathbb{E}^1$, $\mathbb{H}^3$, 
$\mathbb{H}^2\times\mathbb{E}^1$ or $\widetilde{\mathbb{SL}}$.
If $M$ is an infrasolvmanifold then $\pi$ is torsion-free and virtually 
poly-$Z$ of Hirsch length 4, 
so $K$ is torsion-free and virtually poly-$Z$ of Hirsch length 3, 
and the result is clear.

If $\mathbb{X}^3=\mathbb{S}^3$ then $\pi$ is a discrete cocompact subgroup 
of $O(4)\times E(1)$. 
Since $O(4)$ is compact the image of $\pi$ in $E(1)$ is infinite,
and thus has no nontrivial finite normal subgroup.
Therefore $K$ is a subgroup of $O(4)\times\{1\}$.
Since $\pi$ acts freely on $S^3\times\mathbb{R}$ the subgroup $K$ 
acts freely on $S^3$,
and so $K$ is the fundamental group of an $\mathbb{S}^3$-manifold.
If $\mathbb{X}^3=\mathbb{S}^2\times\mathbb{E}^1$ 
then $\pi$ is virtually $\mathbb{Z}^2$.
Hence $K$ has two ends, and so $K\cong\mathbb{Z}$, 
$\mathbb{Z}\oplus(Z/2Z)$ or $D$, by Corollary 4.5.2.
Thus $K$ is the fundamental group of an $\mathbb{S}^2\times\mathbb{E}^1$-manifold.
              
In the remaining cases $\mathbb{X}^3$ is of aspherical type.
The key point here is that a discrete cocompact subgroup of the Lie group
$Isom(\mathbb{X}^3\times\mathbb{E}^1)$ must meet the radical of this group in a 
lattice subgroup.
Suppose first that $\mathbb{X}^3=\mathbb{H}^3$.
After passing to a subgroup of finite index if necessary, we may assume
that $\pi\cong H\times\mathbb{Z}<PSL(2,\mathbb{C})\times\mathbb{R}$, 
where $H$ is a discrete cocompact subgroup of $PSL(2,\mathbb{C})$.
If $K\cap(\{1\}\times\mathbb{R})=1$ then $K$ is commensurate with $H$,
and hence is the fundamental group of an $X$-manifold.
Otherwise the subgroup generated by $K\cap H=K\cap PSL(2,\mathbb{C})$ and 
$K\cap(\{1\}\times\mathbb{R})$ has finite index in $K$ and is isomorphic to 
$(K\cap H)\times\mathbb{Z}$.
Since $K$ is finitely generated so is $K\cap H$, and hence it is finitely
presentable, since $H$ is a 3-manifold group.
Therefore $K\cap H$ is a $PD_2 $-group and so $K$ is the fundamental group of a 
$\mathbb{H}^2\times\mathbb{E}^1$-manifold.

If $\mathbb{X}^3=\mathbb{H}^2\times\mathbb{E}^1$ then we may assume that 
$\pi\cong H\times\mathbb{Z}^2<PSL(2,\mathbb{R})\times\mathbb{R}^2$,
where $H$ is a discrete cocompact subgroup of $PSL(2,\mathbb{R})$.
Since such groups do not admit nontrivial maps to $\mathbb{Z}$ 
with finitely generated kernel $K\cap H$ must be commensurate with $H$, 
and we again see that $K$
is the fundamental group of an $\mathbb{H}^2\times\mathbb{E}^1$-manifold.

A similar argument applies if $\mathbb{X}^3=\widetilde{\mathbb{SL}}$. 
We may assume that $\pi\cong H\times\mathbb{Z}$ where $H$ 
is a discrete cocompact subgroup of $Isom(\widetilde{\mathbb{SL}})$.
Since such groups $H$ do not admit nontrivial maps to $\mathbb{Z}$ 
with finitely generated 
kernel $K$ must be commensurate with $H$ and so is the fundamental group of 
a $\widetilde{\mathbb{SL}}$-manifold.
\end{proof} 

\begin{cor}
Suppose that $M$ has a product geometry $\mathbb{X}\times\mathbb{E}^1$.
If $\mathbb{X}^3=\mathbb{E}^3$, $\mathbb{S}^3$, $\mathbb{S}^2\times\mathbb{E}^1$, $\widetilde{\mathbb{SL}}$ 
or $\mathbb{H}^2\times\mathbb{E}^1$ then $M$ is the mapping torus of an isometry of an 
$\mathbb{X}^3$-manifold with fundamental group $K$.
If $\mathbb{X}^3=\mathbb{N}il^3$ or $\mathbb{S}ol^3$ then $K$ is the fundamental group of 
an $\mathbb{X}^3$-manifold or of a $\mathbb{E}^3$-manifold. 
If $\mathbb{X}^3=\mathbb{H}^3$ then $K$ is the fundamental group of a
$\mathbb{H}^3$- or $\mathbb{H}^2\times\mathbb{E}^1$-manifold.
\end{cor}

\begin{proof}
In all cases $\pi$ is a semidirect product $K\rtimes_\theta\mathbb{Z}$ 
and may be realised by the mapping torus 
of a self homeomorphism of a closed 3-manifold with fundamental group $K$.
If this manifold is an $\mathbb{X}^3$-manifold then 
the outer automorphism class of $\theta$
is finite (see Chapters 8 and 9) and $\theta$ may then be realized by an isometry of an $\mathbb{X}^3$-manifold.
Infrasolvmanifolds are determined up to diffeomorphism by their fundamental 
groups \cite{[Ba04]}, as are $\widetilde{\mathbb{SL}}\times\mathbb{E}^1$- 
and $\mathbb{H}^2\times\mathbb{E}^2$-manifolds \cite{[Vo77]}.
This is also true of $\mathbb{S}^2\times\mathbb{E}^2$- and 
$\mathbb{S}^3\times\mathbb{E}^1$-manifolds, provided $K$ is not finite cyclic, 
and when $K$ is cyclic every such $\mathbb{S}^3\times\mathbb{E}^1$-manifold 
is a mapping torus of an isometry of a suitable lens space \cite{[Oh90]}.
Thus if $M$ is an $\mathbb{X}^3\times\mathbb{E}^1$-manifold and $K$ is the 
fundamental group of an $\mathbb{X}^3$-manifold $M$ is the mapping torus of 
an isometry of an $\mathbb{X}^3$-manifold with fundamental group $K$. 
\end{proof} 

There are (orientable) $\mathbb{N}il^3\times\mathbb{E}^1$- and 
$\mathbb{S}ol^3\times\mathbb{E}^1$-manifolds 
which are mapping tori of self homeomorphisms of flat 3-manifolds, 
but which are not mapping tori of self homeomorphisms of $\mathbb{N}il^3$- 
or $\mathbb{S}ol^3$-manifolds. 
(See Chapter 8.) 
There are analogous examples when $\mathbb{X}^3=\mathbb{H}^3$. 
(See \S4 of Chapter 9.)

We may now improve upon the characterization of mapping tori up to homotopy
equivalence from Chapter 4.

\begin{theorem} 
Let $M$ be a closed $4$-manifold with fundamental 
group $\pi$.
Then $M$ is homotopy equivalent to the mapping torus $M(\Theta)$ of a self 
homeomorphism of a closed $3$-manifold with one of the geometries 
$\mathbb{E}^3$, $\mathbb{N}il^3$, $\mathbb{S}ol^3$, $\mathbb{H}^2\times\mathbb{E}^1$, 
$\widetilde{\mathbb{SL}}$ or $\mathbb{S}^2\times\mathbb{E}^1$ if and only if
\begin{enumerate}
\item $\chi(M)=0$;

\item $\pi$ is an extension of $\mathbb{Z}$ by 
a finitely generated normal subgroup $K$;
 and

\item $K$ has a nontrivial torsion-free abelian normal subgroup $A$.
\end{enumerate}                                                                             
\noindent If $\pi$ is torsion-free $M$ is $s$-cobordant to $M(\Theta)$,
while if moreover $\pi$ is solvable $M$ is homeomorphic to $M(\Theta)$.

\end{theorem}

\begin{proof} 
The conditions are clearly necessary.
Since $K$ has an infinite abelian normal subgroup it has one or two ends.
If $K$ has one end then $M$ is aspherical and so $K$ is a $PD_3$-group, 
by Theorem 4.1 and the Addendum to Theorem 4.5.
Condition (3) then implies that $M'$ is homotopy equivalent to a 
closed 3-manifold with one of the first five of the geometries listed above, 
by Theorem 2.14.
If $K$ has two ends then $M'$ is homotopy equivalent to $S^2\times S^1$, 
$S^2\tilde\times S^1$, $RP^2\times S^1$ or $RP^3\sharp RP^3$, 
by Corollary 4.5.2.

In all cases $K$ is isomorphic to the fundamental group of a closed
3-manifold $N$ which is either Seifert fibred or a $\mathbb{S}ol^3$-manifold, 
and the outer automorphism class $[\theta]$ determined by the extension 
may be realised by a self homeomorphism $\Theta$ of $N$.                        
The manifold $M$ is homotopy equivalent to the mapping torus $M(\Theta)$.
Since $Wh(\pi)=0$, by Theorems 6.1 and 6.3, 
any such homotopy equivalence is simple.

If $K$ is torsion-free and solvable then $\pi$ is virtually poly-$Z$, 
and so $M$ is homeomorphic to $M(\Theta)$, by Theorem 6.11.
Otherwise $N$ is a closed $\mathbb{H}^2\times\mathbb{E}^1$- or 
$\widetilde{\mathbb{SL}}$-manifold,
and so $Wh(\pi\times\mathbb{Z}^n)=0$ for all $n\geq0$, by Theorem 6.3.
The surgery obstruction homorphisms $\sigma_i^N$ are isomorphism 
for all large $i$ \cite{[Ro11]}.
Therefore $M$ is $s$-cobordant to $M(\Theta)$, by Theorem 6.8.
\end{proof} 

Mapping tori of self homeomorphisms of $\mathbb{H}^3$- and 
$\mathbb{S}^3$-manifolds satisfy conditions (1) and (2). 
In the hyperbolic case there is the additional condition

\smallskip
(3-H) {\sl $K$ has one end and no noncyclic abelian subgroup.}

\smallskip
\noindent 
An aspherical closed 3-manifold is hyperbolic if its fundamental
group has no noncyclic abelian subgroup \cite{[B-P]}.
If every $PD_3 $-group is a 3-manifold group 
and the group rings of such hyperbolic groups are regular coherent, 
then Theorem 13.2 extends to show that 
a closed 4-manifold $M$ with fundamental group $\pi$ is $s$-cobordant
to the mapping torus of a self homeomorphism of a hyperbolic 3-manifold 
if and only these three conditions hold.

In the spherical case the appropriate additional conditions are

\smallskip 
(3-S) {\sl $K$ is a fixed point free finite subgroup of $SO(4)$ and
(if $K$ is not cyclic) 
the characteristic automorphism of $K$ determining $\pi$ is realized by an 
isometry of $S^3/K$}; and

\smallskip
(4-S) {\sl the first nontrivial $k$-invariant of $M$ is ``linear".}

\smallskip
\noindent The list of fixed point free finite subgroups of $SO(4)$ is well known.
(See Chapter 11.) If $K$ is cyclic or $Q\times Z/p^jZ$ for some odd prime $p$
or $T^*_k$ then the second part of (3-S) and (4-S) are redundant,
but the general picture is not yet clear \cite{[HM86]}.

The classification of complex surfaces leads easily to 
a complete characterization of the 3-manifolds and diffeomorphisms 
such that the corresponding mapping tori admit complex structures.
(Since $\chi(M)=0$ for any mapping torus $M$ we do not need to enter 
the imperfectly charted realm of surfaces of general type.)

\begin{theorem}
Let $N$ be a closed orientable $3$-manifold with $\pi_1 (N)=\nu$
and let $\theta:N\to N$ be an orientation preserving self diffeomorphism. 
Then the mapping torus $M(\theta)$ admits a complex structure if and only 
if one of the following holds:
\begin{enumerate}
\item $N=S^3/G$ where $G$ is a fixed point free finite subgroup of $U(2)$ and 
the monodromy is as described in $\cite{[Kt75]}$;

\item $N=S^2\times S^1$ (with no restriction on $\theta$);

\item $N=S^1\times S^1\times S^1$ and the image of $\theta$ in $SL(3,\mathbb{Z})$ 
either has finite order or satisfies the equation $(\theta^2-I)^2=0$;

\item $N$ is the flat 3-manifold with holonomy of order $2$, 
$\theta$ induces the identity on $\nu/\nu'$ and $|tr(\theta|_{\nu'})|\leq2$;

\item $N$ is one of the flat 3-manifolds with holonomy cyclic of order $3$, 
$4$ or $6$ and $\theta$ induces the identity on $H_1 (N;\mathbb{Q})$;

\item $N$ is a $\mathbb{N}il^3$-manifold and either the image of $\theta$ 
in $Out(\nu)$ has finite order
or $M(\theta)$ is a $\mathbb{S}ol^4_1$-manifold;

\item $N$ is a $\mathbb{H}^2\times\mathbb{E}^1$- or 
$\widetilde{\mathbb{SL}}$-manifold, $\zeta\nu\cong\mathbb{Z}$  
and the image of $\theta$ in $Out(\nu)$ has finite order.
\end{enumerate}
\end{theorem}

\begin{proof} The mapping tori of these diffeomorphisms 
admit 4-dimensional geometries, 
and it is easy to read off which admit complex structures from 
\cite{[Wl86]}.
In cases (3), (4) and (5) note that a complex surface 
is K\"ahler if and only if its first Betti number is even, 
and so the parity of this Betti number should be invariant under passage 
to finite covers.
(See \cite[Proposition 4.4]{[Wl86]}.)

The necessity of these conditions follows from examining 
the list of minimal complex surfaces $X$ with $\chi(X)=0$ 
in \cite[Table 10]{[BHPV]}, 
together with Bogomolov's theorem on class $VII_0$ surfaces
\cite{[Tl94]} and Perelman's work on geometrization. 
(Fundamental group considerations exclude blowups of ruled surfaces of genus $>1$. All other non-minimal surfaces have $\chi>0$.)
\end{proof} 

In particular, $N$ must be Seifert fibred and most orientable Seifert fibred 
3-manifolds (excepting only the orientable
$\mathbb{H}^2\times\mathbb{E}^1$-manifolds with nonorientable base
orbifolds, $RP^3\sharp RP^3$ and the Hantzsche-Wendt flat 3-manifold) occur. 
Moreover, in most cases 
(with exceptions as in (3), (4) and (6)) the image of $\theta$ in 
$Out(\nu)$ must have finite order. 
Some of the resulting 4-manifolds arise as mapping tori 
in several distinct ways.
The corresponding result for complex surfaces of the form $N\times S^1$ 
for which the obvious smooth $S^1$-action is holomorphic was given 
in \cite{[GG95]}.
In \cite{[EO94]} it is shown that if $\beta_1(N)=0$ then $N\times S^1$ 
admits a complex structure if and only if $N$ is Seifert fibred, 
and the possible complex structures on such products are determined.
Conversely, the following result is very satisfactory from the 
4-dimensional point of view.

\begin{theorem} 
Let $X$ be a complex surface. Then $X$ is diffeomorphic to
the mapping torus of a self diffeomorphism of a closed $3$-manifold 
f and only if $\chi(X)=0$ and $\pi=\pi_1(X)$ is an extension of $\mathbb{Z}$ 
by a finitely generated normal subgroup.
\end{theorem}

\begin{proof} The conditions are clearly necessary. 
Sufficiency of these conditions again
follows from the classification of complex surfaces, as in Theorem 13.3. 
\end{proof} 

\section{Surface bundles and geometries}

Let $p:E\to B$ be a bundle with base $B$ and fibre $F$ 
aspherical closed surfaces.
Then $p$ is determined up to bundle isomorphism by the 
epimorphism $p_*:\pi=\pi_1(E)\to\pi_1(B)$.
If $\chi(B)=\chi(F)=0$ then $E$ has geometry
$\mathbb{E}^4$, $\mathbb{N}il^3\times\mathbb{E}^1$, $\mathbb{N}il^4$ or
$\mathbb{S}ol^3\times\mathbb{E}^1$, by Ue's Theorem.
When the fibre is $Kb$ the geometry must be $\mathbb{E}^4$ or 
$\mathbb{N}il^3\times\mathbb{E}^1$, for then $\pi$ has a normal chain 
$\zeta\pi_1(Kb)\cong\mathbb{Z}<\sqrt{\pi_1(Kb)}\cong\mathbb{Z}^2$, 
so $\zeta\sqrt\pi$ has rank at least 2.
Hence a $\mathbb{S}ol^3\times\mathbb{E}^1$- or $\mathbb{N}il^4$-manifold $M$ 
is the total space of a $T$-bundle over $T$ if and only if $\beta_1(\pi)=2$.
If $\chi(F)=0$ but $\chi(B)<0$ then $E$ need not be geometric.
(See Chapter 7 and \S3 below.)
We shall assume henceforth that $F$ is hyperbolic, i.e.\ that $\chi(F)<0$. 
Then $\zeta\pi_1(F)=1$ and so the characteristic homomorphism 
$\theta:\pi_1(B)\to Out(\pi_1(F))$ determines $\pi$ up to isomorphism, 
by Theorem 5.2. 

\begin{theorem} 
Let $B$ and $F$ be closed surfaces with 
$\chi(B)=0$ and $\chi(F)<0$.
Let $E$ be the total space of the $F$-bundle over $B$ corresponding to a 
homomorphism $\theta:\pi_1(B)\to Out(\pi_1(F))$. 
Then 
\begin{enumerate}
\item $E$ admits the geometry $\mathbb{H}^2\times\mathbb{E}^2$ 
if and only if $\theta$ has finite image;

\item $E$ admits the geometry $\mathbb{H}^3\times\mathbb{E}^1$ 
if and only if $\mathrm{Ker}(\theta)\cong\mathbb{Z}$ and $\mathrm{Im}(\theta)$ 
 contains the class of a pseudo-Anasov homeomorphism of $F$;

\item otherwise $E$ is not geometric.
\end{enumerate} 
If $\mathrm{Ker}(\theta)\not=1$ then $E$ virtually 
has a geometric decomposition.                                                    
\end{theorem} 

\begin{proof} 
Let $\pi=\pi_1(E)$. 
Since $E$ is aspherical, $\chi(E)=0$ and $\pi$ is not solvable the only 
possible geometries are $\mathbb{H}^2\times\mathbb{E}^2$, $\mathbb{H}^3\times\mathbb{E}^1$ 
and $\widetilde{\mathbb{SL}}\times\mathbb{E}^1$.
If $E$ has a proper geometric decomposition the pieces must all have $\chi=0$,
and the only other geometry that may arise is $\mathbb{F}^4$.
In all cases the fundamental group of each piece has a nontrivial abelian normal subgroup.

If $\mathrm{Ker}(\theta)\not=1$ then $E$ is virtually 
a cartesian product $N\times S^1$, where $N$ 
is the mapping torus of a self diffeomorphism $\psi$ of $F$ 
whose isotopy class in $\pi_0(Diff(F))\cong Out(\pi_1(F))$ 
generates a subgroup of finite index in $\mathrm{Im}(\theta)$.
Since $N$ is a Haken 3-manifold it has a geometric decomposition, 
and hence so does $E$.
The mapping torus $N$ is an $\mathbb{H}^3$-manifold if and only if 
$\psi$ is pseudo-Anasov.
In that case the action of $\pi_1(N)\cong\pi_1(F)\rtimes_\psi\mathbb{Z}$ 
on $H^3$ extends to an embedding 
$p:\pi/\sqrt\pi\to Isom(\mathbb{H}^3)$, by Mostow rigidity.                                          
Since $\sqrt\pi\not=1$ we may also find a homomorphism $\lambda:\pi\to D<Isom(E^1)$ 
such that $\lambda(\sqrt\pi)\cong\mathbb{Z}$. 
Then $\mathrm{Ker}(\lambda)$ is an extension of $\mathbb{Z}$ by $F$ 
and is commensurate with $\pi_1(N)$, so
is the fundamental group of a Haken $\mathbb{H}^3$-manifold, $\widehat N$ say.
Together these homomorphisms determine a free cocompact action of $\pi$ on $H^3\times E^1$.
If $\lambda(\pi)\cong\mathbb{Z}$ then $M=\pi\backslash (H^3\times E^1)$ 
is the mapping torus of a self homeomorphism of $\widehat N$; 
otherwise it is the union of two twisted $I$-bundles over $\widehat N$.
In either case it follows from standard 3-manifold theory that 
since $E$ has a similar structure $E$ and $M$ are diffeomorphic.

If $\theta$ has finite image then $\pi/C_\pi(\pi_1(F))$ is a finite extension
of $\pi_1(F)$ and so acts properly and cocompactly on $\mathbb{H}^2$.
We may therefore construct an $\mathbb{H}^2\times\mathbb{E}^2$-manifold
with group $\pi$ and which fibres over $B$ as in Theorems 7.3 and 9.9. 
Since such bundles are determined up to diffeomorphism by their 
fundamental groups $E$ admits this geometry.

If $E$ admits the geometry $\mathbb{H}^2\times\mathbb{E}^2$ then 
$\sqrt\pi=\pi\cap Rad(Isom(\mathbb{H}^2\times\mathbb{E}^2))=
\pi\cap(\{1\}\times\mathbb{R}^2)\cong\mathbb{Z}^2$ 
\cite[Proposition 8.27]{[Rg]}. Hence $\theta$ has finite image.

If $E$ admits the geometry $\mathbb{H}^3\times\mathbb{E}^1$ 
then $\sqrt\pi=\pi\cap(\{1\}\times\mathbb{R})\cong\mathbb{Z}$ 
\cite[Proposition 8.27]{[Rg]}.
Hence $\mathrm{Ker}(\theta)\cong\mathbb{Z}$ and $E$ is finitely covered 
by a cartesian product $N\times S^1$,
where $N$ is a hyperbolic 3-manifold which is also an $F$-bundle over $S^1$. 
The geometric monodromy of the latter bundle is a pseudo-Anasov diffeomorphism
of $F$ whose isotopy class is in $\mathrm{Im}(\theta)$.

If $\rho$ is the group of a 
$\widetilde{\mathbb{SL}}\times\mathbb{E}^1$-manifold 
then $\sqrt\rho\cong\mathbb{Z}^2$ and $\sqrt\rho\cap K'\not=1$ 
for all subgroups $K$ of finite index, 
and so $E$ cannot admit this geometry. 
\end{proof} 

It follows from the Virtual Fibration Theorem that
every $\mathbb{H}^3\times\mathbb{E}^1$-manifold is 
finitely covered by such a bundle space \cite{[Ag13]}.
Let $X$ be the exterior of the figure eight knot,
and let $\phi$ be the self homeomorphism of 
$\partial{X\times{S^1}}=S^1\times{S^1}\times{S^1}$
which preserves the longitude of the knot and swaps the meridian and 
the third factor.
Then $M=X\times{S^1}\cup_\phi{X\times{S^1}}$ 
fibres over $T$ with fibre $T\sharp{T}$ and 
$\theta$ is injective.
It is not geometric, but is the union of two pieces of type
$\mathbb{H}^3\times\mathbb{E}^1$.

We shall assume henceforth that $B$ is also hyperbolic.
Then $\chi(E)>0$ and $\pi_1(E)$ has no solvable subgroups of Hirsch length 3.
Hence the only possible geometries on $E$ are $\mathbb{H}^2\times\mathbb{H}^2$, 
$\mathbb{H}^4$ and $\mathbb{H}^2(\mathbb{C})$.                           
(These are the least well understood geometries, 
and little is known about the possible lattices.)
                         
\begin{theorem} 
Let $B$ and $F$ be closed hyperbolic surfaces, 
and let $E$ be the total space of the $F$-bundle over $B$ corresponding 
to a homomorphism $\theta:\pi_1(B)\to Out(\pi_1(F))$. 
Then the following are equivalent:
\begin{enumerate}
\item $E$ admits the geometry $\mathbb{H}^2\times\mathbb{H}^2$;

\item $E$ is finitely covered by a cartesian product of surfaces;

\item $\theta$ has finite image.             

\end{enumerate}     
\end{theorem}

\begin{proof} 
Let $\pi=\pi_1(E)$ and $\phi=\pi_1(F)$. 
If $E$ admits the geometry $\mathbb{H}^2\times\mathbb{H}^2$ 
it is virtually a cartesian product,
by Corollary 9.9.1, and so (1) implies (2). 
                                                            
If $\pi$ is virtually a direct product of $PD_2$-groups then
$[\pi:\phi{C_\pi(\phi)}]<\infty$, by Theorem 5.4.
Therefore the image of $\theta$ is finite and so (2) implies (3). 
              
If $\theta$ has finite image then $\mathrm{Ker}(\theta)\not=1$ 
and $\pi/C_\pi(\phi)$ is a finite extension of $\phi$. 
Hence there is a homomorphism $p:\pi\to Isom(\mathbb{H}^2)$ with 
kernel $C_\pi(\phi)$ and with image a discrete cocompact subgroup. 
Let $q:\pi\to \pi_1(B)< Isom(\mathbb{H}^2)$.
Then $(p,q)$ embeds $\pi$ as a discrete cocompact subgroup of $Isom(\mathbb{H}^2\times\mathbb{H}^2)$,
and the closed 4-manifold $M=\pi\backslash (H^2\times H^2)$ 
clearly fibres over $B$.
Such bundles are determined up to diffeomorphism by the corresponding
extensions of fundamental groups, by Theorem 5.2.
Therefore $E$ admits the geometry $\mathbb{H}^2\times\mathbb{H}^2$ 
and so (3) implies (1). 
\end{proof} 

If $F$ is orientable and of genus $g$ its mapping class group 
$M_g=Out(\pi_1(F))$ has only finitely many conjugacy classes 
of finite groups. (See Chapter 7 of \cite{[FM]}.)
In Corollary 13.7.2 we shall show that no such bundle space $E$
is homotopy equivalent to a $\mathbb{H}^2(\mathbb{C})$-manifold. 
U.Hamenst\"adt has announced that if $\pi_1(E)$ is word-hyperbolic
then $\sigma(E)\not=0$,
and so $E$ cannot admit the geometry $\mathbb{H}^4$ \cite{[Ha13]}.
Thus only finitely many orientable bundle spaces
with given Euler characteristic $>0$ are geometric.
(This follows also from \cite{[Wa72]}.)


If $\mathrm{Im}(\theta)$ contains the outer automorphism class 
determined by a Dehn twist on $F$ then $E$ admits no metric 
of nonpositive curvature \cite{[KL96]}.
For any given base and fibre genera there are only finitely 
many extensions of $PD_2^+$-groups by $PD_2^+$-groups 
which contain no noncyclic abelian subgroups \cite{[Bo09]}.
For such groups $\theta$ must be injective.
There are bundle spaces with $\theta$ injective, 
since the $PD_2^+$-group of genus 4 embeds in $M_2$
\cite{[Re06]}.
Are there infinitely many such with given base and fibre?

If $B$ and $E$ are orientable and $F$ has genus $g$ then
$\sigma(E)=-\theta^* \tau\cap[B]$, 
where $\tau\in H^2(M_g;\mathbb{Z})$ is induced from a universal class 
in $H^2(Sp_{2g}(\mathbb{Z});\mathbb{Z})$ via the natural representation 
of $M_g$ as symplectic isometries of the intersection form on 
$H_1(F;\mathbb{Z})\cong\mathbb{Z}^{2g}$.
(\cite{[Me73]} -- see \cite{[En98],[Ho01]} for alternative approaches.)
In particular, if $g=2$ then $\sigma(E)=0$, since $H^2(M_2;\mathbb{Q})=0$.

Every closed orientable $\mathbb{H}^2\times\mathbb{H}^2$-manifold 
has a 2-fold cover which is a complex surface, and has signature 0.
Conversely, if a bundle space $E$ is a complex surface 
and $p$ is a holomorphic submersion then $\sigma(E)=0$ 
implies that the fibres are isomorphic,
and so $E$ is an $\mathbb{H}^2\times\mathbb{H}^2$-manifold 
\cite{[Ko99]}. 
This is also so if $p$ is a holomorphic fibre bundle.
(See \cite[\S V.6]{[BHPV]}.)
Any holomorphic submersion with base of genus at most 1 or
fibre of genus at most 2 is a holomorphic fibre bundle \cite{[Ks68]}. 
There are such holomorphic submersions in which $\sigma(E)\not=0$ 
and so which are not virtually products. (See \cite[\S V.14]{[BHPV]}.) 
The image of $\theta$ must contain the outer automorphism class
determined by a pseudo-Anasov homeomorphism and not be virtually abelian 
\cite{[Sh97]}.

If a bundle space $E$ with $\chi(E)>0$ has a proper geometric decomposition 
the pieces are reducible $\mathbb{H}^2\times\mathbb{H}^2$-manifolds,
the cusps are $\mathbb{H}^2\times\mathbb{E}^1$-manifolds
and the inclusions of the cusps induce injections on $\pi_1$. 
There are non-geometric examples.
One with $B$ and $F$ of genus 2 and $\mathrm{Im}(\theta)\cong{D}$ 
is given in \cite{[Hi11]}.

\section{Geometric decompositions of Seifert fibred 4-manifolds}

Most Seifert fibred 4-manifolds with hyperbolic base orbifold
have geometric decompositions with all pieces of type 
$\mathbb{H}^2\times\mathbb{E}^2$,
since $\mathrm{Im}(\theta)\cong\pi\mathcal{G}$, 
where $\mathcal{G}$ is a finite graph of finite groups,
and we may decompose the base correspondingly.
In particular, if the general fibre is $Kb$ the manifold is geometric,
by Theorem 9.5, since $Out(\mathbb{Z}\rtimes_{-1}\mathbb{Z})$ is finite.
However torus bundles over surfaces need not be geometric,
as we shall show.
We also give a Seifert fibred 4-manifold 
with no geometric decomposition.
(This can happen only if the base orbifold has 
at least two cone points of order 2, 
at which the action has one eigenvalue $-1$.)

We show first that there is no closed 4-manifold with the geometry 
$\mathbb{F}^4$.
If $G=Isom(\mathbb{F}^4)$ and $\Gamma<G$ is an $\mathbb{F}^4$-lattice 
then $\Gamma\cap{Rad(G)}$ is a lattice in $Rad(G)\cong\mathbb{R}^2$, 
and $\Gamma/\Gamma\cap{Rad(G)}$ is a discrete cocompact subgroup
of $G/Rad(G)$ \cite[Proposition 8.27]{[Rg]}.
Hence $\sqrt\Gamma=\Gamma\cap{Rad(G)}\cong\mathbb{Z}^2$ and
$\Gamma/\sqrt\Gamma$ is a subgroup of finite index in $GL(2,\mathbb{Z})$.
Therefore $v.c.d.\Gamma=3$ and so $\Gamma\backslash F^4$ is not a closed
4-manifold.
As observed in Chapter 7, such quotients are Seifert fibred,
and the base is a punctured hyperbolic orbifold.
Thus if $M$ is a compact manifold with interior $\Gamma\backslash{F^4}$
the double $DM=M\cup_\partial M$ is Seifert fibred but is not geometric, 
since $\sqrt\pi\cong\mathbb{Z}^2$ but $[\pi:C_\pi (\sqrt\pi)]$ is infinite. 

The orientable surface of genus 2 can be represented as a double in 
two distinct ways;
we shall give corresponding examples of nongeometric torus bundles
which admit geometric decompositions of type $\mathbb{F}^4$.

{\bf 1.} Let $F(2)$ be the free group of rank two and let 
$\gamma:F(2)\to SL(2,\mathbb{Z})$
have image the commutator subgroup $SL(2,\mathbb{Z})'$, 
which is freely generated by 
$X=\left(\begin{smallmatrix}
2&1\\
1&1
\end{smallmatrix}\right)$ 
and
$Y=\left(\begin{smallmatrix}
1&1\\
1&2
\end{smallmatrix}\right)$.
The natural surjection from $SL(2,\mathbb{Z})$ to $PSL(2,\mathbb{Z})$ 
induces an isomorphism of commutator subgroups.
(See \S2 of Chapter 1.)
The parabolic subgroup $PSL(2,\mathbb{Z})'\cap Stab(0)$ is generated 
by the image of 
$XY^{-1}X^{-1}Y=\left(\begin{smallmatrix}
-1&0\\
-6&-1
\end{smallmatrix}\right)$.
Hence $[Stab(0):PSL(2,\mathbb{Z})'\cap Stab(0)]=6
=[PSL(2,\mathbb{Z}):PSL(2,\mathbb{Z})']$, 
and so $PSL(2,\mathbb{Z})'$ has a single cusp,
represented by 0.
The quotient space $PSL(2,\mathbb{Z})'\backslash H^2$ is 
the once-punctured torus. 
Let $N\subset PSL(2,\mathbb{Z})'\backslash H^2$ be 
the complement of an open horocyclic neighbourhood of the cusp.
The double $DN$ is the orientable surface of genus 2.
The semidirect product $\Gamma=\mathbb{Z}^2\rtimes_\gamma F(2)$ is a lattice 
in $Isom(\mathbb{F}^4)$,
and the double of the bounded manifold with interior $\Gamma\backslash F^4$
is a torus bundle over $DN$.
          
{\bf 2.} Let $\delta:F(2)\to SL(2,\mathbb{Z})$ have image the subgroup 
which is freely generated by 
$U=\left(\begin{smallmatrix}
1&0\\
2&1
\end{smallmatrix}\right)$
and 
$V=\left(\begin{smallmatrix}
1&2\\
0&1
\end{smallmatrix}\right)$.
Let $\bar\delta:F(2)\to PSL(2,\mathbb{Z})$ be the composed map. 
Then $\bar\delta$ is injective and $[PSL(2,\mathbb{Z}):\bar\delta(F(2))]=6$.
(Note that $\delta(F(2))$ and $-I$ together generate 
the level 2 congruence subgroup.)
Moreover $[Stab(0):\bar\delta(F(2))\cap Stab(0)]=2$. 
Hence $\bar\delta(F(2))$ has three cusps, represented by 0, $\infty$ and 1, 
and $\bar\delta(F(2))\backslash H^2$ is the thrice-punctured sphere.
The corresponding parabolic subgroups are generated by $U$, 
$V$ and $VU^{-1}$, respectively.
Doubling the complement $N$ of disjoint horocyclic neighbourhoods 
of the cusps in $\bar\delta(F(2))\backslash H^2$ 
again gives an orientable surface of genus 2.
The presentation for $\pi_1 (DN)$ derived from this construction is
\begin{equation*}
\langle U,V,U_1,V_1,s,t\mid s^{-1}Us=U_1,t^{-1}Vt=V_1,
VU^{-1}=V_1U_1^{-1}\rangle,
\end{equation*}
which simplifies to the usual presentation 
$\langle U,V,s,t\mid s^{-1}V^{-1}sV=t^{-1}U^{-1}tU\rangle$.
The semidirect product $\Delta=\mathbb{Z}^2\rtimes_\delta F(2)$ 
is a lattice in $Isom(\mathbb{F}^4)$,
and doubling again gives a torus bundle over $DN$.
                 
{\bf 3.} 
Let $\pi$ be the group with presentation
\begin{gather*}
\langle {a,b,v,w,x,y,z}\mid{ab=ba},~v^2=w^2=x^2=a,~vb=bv,~
wbw^{-1}=xbx^{-1}=b^{-1}\!,\\
y^2=z^2=ab,~yay^{-1}=zaz^{-1}=ab^2,~yby^{-1}=zbz^{-1}=b^{-1},~vwx=yz
\rangle.
\end{gather*}
Then $\pi=\pi_1(E)$, where $E$ is a Seifert manifold with base 
the hyperbolic orbifold $B=S^2(2,2,2,2,2)$ and regular fibre $T$.
The action $\theta$ of $\pi_1^{orb}(B)$ on $\sqrt\pi=\langle{a,b}\rangle$
is generated by
$\theta(u)=\left(\smallmatrix1&0\\
0&-1\endsmallmatrix\right)$ and
$\theta(x)=\left(\smallmatrix1&0\\
2&-1\endsmallmatrix\right)$.
Hence $\mathrm{Im}(\theta)\cong{D}$ is infinite and
has infinite index in $GL(2,\mathbb{Z})$.
Thus $E$ is not geometric and has no pieces of type $\mathbb{F}^4$.
On the other hand $B$ has no proper decomposition into hyperbolic pieces,
and so $E$ has no geometric decomposition at all.

\newpage
\section{Complex surfaces and fibrations}

It is an easy consequence of the classification of surfaces that a minimal
compact complex surface $S$ is ruled over a curve $C$ of genus $\geq2$
if and only if $\pi_1(S)\cong\pi_1(C)$ and $\chi(S)=2\chi(C)$. 
(See \cite[Table 10]{[BHPV]}.)
We shall give a similar characterization of the complex surfaces which
admit holomorphic submersions to complex curves of genus $\geq2$, 
and more generally of quotients of such surfaces by free actions 
of finite groups.
However we shall use the classification only to handle 
the cases of non-K\"ahler surfaces.

\begin{theorem} 
Let $S$ be a complex surface.
Then $S$ has a finite covering space which admits a holomorphic submersion onto a complex curve, 
with base and fibre of genus $\geq2$, if and only if $\pi=\pi_1(S)$ has normal 
subgroups $K<\hat\pi$ such that $K$ and $\hat\pi/K$ are $PD_2^+$-groups, 
$[\pi:\hat\pi]<\infty$ and $[\pi:\hat\pi]\chi(S)=\chi(K)\chi(\hat\pi/K)>0$.
\end{theorem}

\begin{proof} 
The conditions are clearly necessary.
Suppose that they hold.
Then $S$ is aspherical, by Theorem 5.2.
In particular, $\pi$ is torsion-free and $\pi_2 (S)=0$, so $S$ is minimal.
Let $\widehat{S}$ be the finite covering space corresponding to $\hat\pi$.
Then $\chi(\widehat{S})>0$ and $\beta_1(\widehat{S})\geq 4$.
If $\beta_1(\widehat{S})$ were odd then $\widehat{S}$ 
would be minimally properly elliptic, 
by the classification of surfaces.
But then $\widehat{S}$ would have a singular fibre,
since $\chi(\widehat{S})>0$.
Hence $\pi_1(\widehat{S})\cong\pi^{orb}(B)$,
where $B$ is the base orbifold of an elliptic fibration
of a deformation of $S$ \cite[Proposition 2]{[Ue86]}.
Therefore $\beta_1(\widehat{S})$ is even,
and hence $\widehat{S}$ and $S$ are K\"ahler.
(See \cite[Chapter IV, Theorem 3.1]{[BHPV]}.)

After enlarging $K$ if necessary we may assume that $\pi/K$ 
has no nontrivial finite normal subgroup.
It is then isomorphic to a discrete cocompact group of isometries 
of $\mathbb{H}^2$.
Since $S$ is K\"ahler and $\beta_1^{(2)}(\pi/K)\not=0$
there is a properly discontinuous holomorphic action of $\pi/K$ on 
$\mathbb{H}^2$ and a $\pi/K$-equivariant holomorphic map 
from the covering space $S_K$ to $\mathbb{H}^2$, 
with connected fibres \cite[Theorems 4.1 and 4.2]{[ABR92]}. 
Let $B$ and $\widehat B$ be the complex curves $\mathbb{H}^2/(\pi/K)$ and 
$\mathbb{H}^2/(\hat\pi/K)$, respectively,
and let $h:S\to B$ and $\hat h:\widehat S\to\widehat B$ be the induced maps.
The quotient map from $\mathbb{H}^2$ to $\hat B$ is a covering projection, 
since $\hat\pi/K$ is torsion-free, 
and so $\pi_1(\hat h)$ is an epimorphism with 
kernel $K$.

The map $h$ is a submersion away from the preimage of a finite subset 
$D\subset B$.
Let $F$ be the general fibre and $F_d$ the fibre over $d\in D$.
Fix small disjoint discs $\Delta_d\subset B$ about each point of $D$, 
and let $B^*=B\setminus\cup_{d\in D}\Delta_d$,
$S^*=h^{-1}(B^*)$ and $S_d=h^{-1}(\Delta_d)$.
Since $h|_{S^*}$ is a submersion $\pi_1(S^*)$ is an extension 
of $\pi_1(B^*)$ by $\pi_1(F)$.
The inclusion of $\partial S_d$ into $S_d\setminus{F_d}$ 
is a homotopy equivalence.
Since $F_d$ has real codimension 2 in $S_d$, 
the inclusion of $S_d\setminus{F_d}$ into $S_d$ is 2-connected.
Hence $\pi_1(\partial S_d)$ maps onto $\pi_1(S_d)$. 
                               
Let $m_d=[\pi_1(F_d)]:Im(\pi_1(F))]$. 
After blowing up $S^*$ at singular points of $F_d$ we may assume that $F_d$  
has only normal crossings. 
We may then pull $h|_{S_d}$ back over a suitable branched covering of $\Delta_d$
to obtain a singular fibre $\widetilde F_d$ with no multiple components 
and only normal crossing singularities.
In that case $\widetilde F_d$  is obtained from $F$ by shrinking vanishing cycles, 
and so $\pi_1(F)$ maps onto $\pi_1(\widetilde F_d)$.
Since blowing up a point on a curve does not change the fundamental group 
it follows from \cite[\S III.9]{[BHPV]} 
that in general $m_d$ is finite.

We may regard $B$ as an orbifold with cone singularities of order $m_d$ 
at $d\in D$.
By the Van Kampen theorem (applied to the space $S$ and the orbifold $B$)
the image of $\pi_1(F)$ in $\pi$ is a normal subgroup
and $h$ induces an isomorphism from $\pi/\pi_1(F)$ to $\pi_1^{orb}(B)$.
Therefore the kernel of the canonical map from $\pi_1^{orb}(B)$ to $\pi_1(B)$
is isomorphic to $K/$Im$(\pi_1(F))$.
But this is a finitely generated normal subgroup of infinite index in 
$\pi_1^{orb}(B)$,
and so must be trivial.
Hence $\pi_1(F)$ maps onto $K$, and so $\chi(F)\leq\chi(K)$. 

Let $\widehat D$ be the preimage of $D$ in $\widehat B$. 
The general fibre of $\hat h$ is again $F$. 
Let $\widehat F_d$ denote the fibre over $d\in\widehat D$.
Then $\chi(\widehat S)=\chi(F)\chi(B)+\Sigma_{d\in\widehat D} 
(\chi(\widehat F_d)-\chi(F))$ 
and $\chi(\widehat F_d)\geq\chi(F)$ \cite[Proposition III.11.4]{[BHPV]}. 
Moreover $\chi(\widehat F_d)>\chi(F)$ unless 
$\chi(\widehat F_d)=\chi(F)=0$ \cite[Remark III.11.5]{[BHPV]}.
Since $\chi(\widehat B)=\chi(\hat\pi/K)<0$, 
$\chi(\widehat S)=\chi(K)\chi(\hat\pi/K)$
and $\chi(F)\leq\chi(K)$ it follows that $\chi(F)=\chi(K)<0$ and 
$\chi(\widehat F_d)=\chi(F)$ for all $d\in\widehat D$. 
Therefore $\widehat F_d\cong F$ for all $d\in\widehat D$ and so 
$\hat h$ is a holomorphic submersion.
\end{proof} 

Similar results have been found independently by Kapovich and Kotschick
\cite{[Ka98],[Ko99]}.
Kapovich assumes instead that $K$ is $FP_2$ and $S$ is aspherical.
As these hypotheses imply that $K$ is a $PD_2$-group, by Theorem 1.19,
the above theorem applies.

We may construct examples of such surfaces as follows. 
Let $n>1$ and $C_1$ and $C_2$ be two curves such that $Z/nZ$ acts freely on $C_1$ 
and with isolated fixed points on $C_2$.
The quotient of $C_1\times C_2$ by the diagonal action is a complex surface
$S$ and projection from $C_1\times C_2$ to $C_2$ induces a 
holomorphic mapping from $S$ onto $C_2/(Z/nZ)$ with critical values 
corresponding to the fixed points.

\begin{cor}
The surface $S$ admits such a holomorphic submersion onto a complex curve 
if and only if $\pi/K$ is a $PD_2^+$-group for some such $K$.
\qed
\end{cor}

\begin{cor}
No bundle space $E$ is homotopy equivalent to a closed 
$\mathbb{H}^2(\mathbb{C})$-manifold.
\end{cor}

\begin{proof} 
Since $\mathbb{H}^2(\mathbb{C})$-manifolds have 2-fold coverings which are 
complex surfaces, we may assume that $E$ is homotopy equivalent to
a complex surface $S$.
By the theorem, $S$ admits a holomorphic submersion onto a complex curve.
But then $\chi(S)>3\sigma(S)$ \cite{[Li96]}, and so $S$ cannot be a
$\mathbb{H}^2(\mathbb{C})$-manifold.
\end{proof} 

The relevance of Liu's work was observed by Kapovich,
who has also found a cocompact $\mathbb{H}^2(\mathbb{C})$-lattice which is 
an extension of a $PD_2^+$-group by a finitely generated normal
subgroup, but which is not almost coherent \cite{[Ka98]}.
                                                 
Similar arguments may be used to show that a K\"ahler surface $S$ is a minimal 
properly elliptic surface with no singular fibres if and only if $\chi(S)=0$ 
and $\pi=\pi_1(S)$ has a normal subgroup $A\cong\mathbb{Z}^2$ 
such that $\pi/A$ is virtually torsion-free and indicable, 
but is not virtually abelian. 
(This holds also in the non-K\"ahler case as a consequence 
of the classification of surfaces.)
Moreover, if $S$ is not ruled it is a complex torus, a hyperelliptic surface, 
an Inoue surface, a Kodaira surface or a minimal elliptic surface if and only 
if $\chi(S)=0$ and $\pi_1(S)$ has a normal 
subgroup $A$ which is poly-$Z$ and not cyclic, 
and such that $\pi/A$ is infinite and virtually torsion-free indicable.
(See \cite[Theorem X.5]{[H2]}.)

We may combine Theorem 13.7 with some observations deriving from the 
classification of surfaces for our second result.

\begin{theorem}
Let $S$ be a complex surface such that $\pi=\pi_1(S)\not=1$.
If $S$ is homotopy equivalent to the total space $E$ of a bundle over 
a closed orientable $2$-manifold then $S$ is diffeomorphic to $E$.
\end{theorem}                        
    
\begin{proof} 
Let $B$ and $F$ be the base and fibre of the bundle, respectively. 
Suppose first that $\chi(F)=2$.                               
Then $\chi(B)\leq0$, for otherwise $S$ would be simply-connected.
Hence $\pi_2(S)$ is generated by an embedded $S^2$ with self-intersection 0,
and so $S$ is minimal. 
Therefore $S$ is ruled over a curve diffeomorphic to $B$, 
by the classification of surfaces. 
                                                              
Suppose next that $\chi(B)=2$. 
If $\chi(F)=0$ and $\pi\not\cong\mathbb{Z}^2$ 
then $\pi\cong{\mathbb{Z}\oplus(Z/nZ)}$ for some $n>0$.
Then $S$ is a Hopf surface and so is determined up to diffeomorphism 
by its homotopy type \cite[Theorem 12]{[Kt75]}. 
If $\chi(F)=0$ and $\pi\cong\mathbb{Z}^2$ or if $\chi(F)<0$ then 
$S$ is homotopy equivalent to $S^2\times F$, 
so $\chi(S)<0$, $w_1(S)=w_2(S)=0$ and $S$ is ruled over a curve 
diffeomorphic to $F$.
Hence $E$ and $S$ are diffeomorphic to $S^2\times F$.

In the remaining cases $E$ and $F$ are both aspherical.
If $\chi(F)=0$ and $\chi(B)\leq 0$ then $\chi(S)=0$ and $\pi$ has one end.
Therefore $S$ is a complex torus, a hyperelliptic surface, 
an Inoue surface, a Kodaira surface or a minimal properly elliptic surface.
(This uses Bogomolov's theorem on class $VII_0$ surfaces 
\cite{[Tl94]}.)
The Inoue surfaces are mapping tori of self-diffeomorphisms of $S^1\times S^1\times S^1$,
and their fundamental groups are not extensions 
of $\mathbb{Z}^2$ by $\mathbb{Z}^2$,
so $S$ cannot be an Inoue surface.
As the other surfaces are Seifert fibred 4-manifolds $E$ and $S$ are diffeomorphic \cite{[Ue91]}.
                     
If $\chi(F)<0$ and $\chi(B)=0$ then $S$ is a minimal properly elliptic surface.
Let $A$ be the normal subgroup of the general fibre in an elliptic fibration.
Then $A\cap\pi_1(F)=1$ (since $\pi_1(F)$ has no nontrivial abelian
normal subgroup) and so $[\pi:A.\pi_1(F)]<\infty$.
Therefore $E$ is finitely covered by a cartesian product $T\times F$, 
and so is Seifert fibred.
Hence $E$ and $S$ are diffeomorphic, by \cite{[Ue91]}.

The remaining case ($\chi(B)<0$ and $\chi(F)<0$)
is an immediate consequence of Theorem 13.7, since such bundles 
are determined by the corresponding extensions of fundamental groups. 
(See Theorem 5.2.)
\end{proof} 

A 1-connected 4-manifold which fibres over a 2-manifold is 
homeomorphic to ${CP^1\times{CP^1}}$ or ${CP^2\sharp\overline{CP^2}}$.
(See Chapter 12.)
Is there such a complex surface of general type?     
(No surface of general type is diffeomorphic to
${CP^1\times{CP^1}}$ or ${CP^2\sharp\overline{CP^2}}$ 
\cite{[Qi93]}.)

\begin{cor}
If moreover the base has genus $0$ or $1$ or the 
fibre has genus $2$ then $S$ is finitely covered by a cartesian product.
\end{cor}

\begin{proof} A holomorphic submersion with fibre of genus 2 is 
the projection of a holomorphic fibre bundle and hence $S$ is virtually a product \cite{[Ks68]}.
\end{proof}

Up to deformation there are only finitely many algebraic surfaces 
with given Euler characteristic $>0$ which admit holomorphic submersions 
onto curves \cite{[Pa68]}. 
By the argument of the first part of Theorem 13.1 this remains true
without the hypothesis of algebraicity, for any such complex surface must be 
K\"ahler, and K\"ahler surfaces are deformations of algebraic surfaces. 
(See \cite[Theorem 4.3]{[Wl86]}.)
Thus the class of bundles realized by complex surfaces is very restricted.
Which extensions of $PD_2^+$-groups by $PD_2^+$-groups 
are realized by complex surfaces (i.e., not necessarily aspherical)?

The equivalence of the conditions 
``$S$ is ruled over a complex curve of genus $\geq 2$",
``$\pi=\pi_1 (S)$ is a $PD_2^+ $-group and $\chi(S)=2\chi(\pi)<0$" and
``$\pi_2 (S)\cong\mathbb{Z}$, $\pi$ acts trivially on $\pi_2 (S)$ 
and $\chi(S)<0$" 
also follows by an argument similar to that used in Theorems 13.7 and 13.8.
(See \cite[Theorem X.6]{[H2]}.)

If $\pi_2 (S)\cong\mathbb{Z}$ and $\chi(S)=0$ then 
$\pi$ is virtually $\mathbb{Z}^2$.
The finite covering space with fundamental group $\mathbb{Z}^2$ is K\"ahler, 
and therefore so is $S$.
Since $\beta_1 (S)>0$ and is even, we must have $\pi\cong\mathbb{Z}^2$,
and so $S$ is either ruled over an elliptic curve or is a minimal properly elliptic surface,
by the classification of complex surfaces.
In the latter case the base of the elliptic fibration is $CP^1$, 
there are no singular
fibres and there are at most 3 multiple fibres. (See \cite{[Ue91]}.) 
Thus $S$ may be obtained from a cartesian product $CP^1 \times E$ 
by logarithmic transformations.
(See \cite[\S V.13]{[BHPV]}.) Must $S$ in fact be ruled?

If $\pi_2 (S)\cong\mathbb{Z}$ and $\chi(S)>0$ then $\pi=1$, by Theorem 10.1.
Hence $S\simeq CP^2$ and so $S$ is analytically isomorphic to $CP^2$,
by a result of Yau.
(See \cite[Chapter 5, Theorem 1.1]{[BHPV]}.)

\section{$S^1$-Actions and foliations by circles}

For each of the geometries $\mathbb{X}^4=\mathbb{S}^3\times\mathbb{E}^1,
\mathbb{H}^3\times\mathbb{E}^1,\mathbb{S}^2\times\mathbb{E}^2,
\mathbb{H}^2\times\mathbb{E}^2,{\widetilde{\mathbb{SL}}\times\mathbb{E}^1}$ 
or $\mathbb{E}^4$ the radical of $Isom(\mathbb{X}^4)_o$
is a vector subgroup $V=\mathbb{R}^d$, where $d=1,1,2,2,2$ or 4, respectively.
Similarly, if $\mathbb{X}^4=\mathbb{S}ol^4_1,
\mathbb{N}il^3\times\mathbb{E}^1,\mathbb{N}il^4$ or
$\mathbb{S}ol^3\times\mathbb{E}^1$ the nilradical 
or its commutator subgroup is a vector subgroup, 
of dimension $1,2, 2$ or 3, respectively.
In each case the model space has a corresponding foliation by 
copies of $\mathbb{R}^d$, which is invariant under the isometry group.

If $\mathbb{X}^4\not=\mathbb{E}^4$ or $\mathbb{N}il^4$
then $V$ is central in $Isom(\mathbb{X}^4)_o$.
For $\mathbb{N}il^4$ we may replace $V$ by 
the centre $\zeta{Nil^4}=\mathbb{R}$,
while if $\mathbb{X}^4=\mathbb{E}^4$
then $V$ is a cocompact subgroup of $Isom(\mathbb{E}^4)$.
In each case,
if $\pi$ is a lattice in $Isom(\mathbb{X}^4)$ 
then $\pi\cap{V}\cong\mathbb{Z}^d$.
The corresponding closed geometric 4-manifolds have natural foliations 
with leaves flat $d$-manifolds, and have finite coverings which 
admit an action by the torus $\mathbb{R}^d/\mathbb{Z}^d$ 
with all orbits of maximal dimension $d$.
These actions lift to principal actions (i.e., without exceptional orbits) 
on suitable finite covering spaces.
(This does not hold for all actions.
For instance, $S^3$ admits non-principal $S^1$-actions without fixed points.) 

If a closed manifold $M$ is the total space of an $S^1$-bundle
then $\chi(M)=0$, and if it is also aspherical then 
$\pi_1(M)$ has an infinite cyclic normal subgroup.
As lattices in $Isom(\mathbb{S}ol^4_{m,n})$ or $Isom(\mathbb{S}ol^4_0)$ 
do not have such subgroups, it follows that 
no other closed geometric 4-manifold is finitely covered 
by the total space of an $S^1$-bundle.
Is every geometric 4-manifold $M$ with $\chi(M)=0$ nevertheless 
foliated by circles?

If a complex surface is foliated by circles it 
admits one of the above geometries, and so it must be Hopf, 
hyperelliptic, Inoue of type $S^\pm_{N\dots }$, Kodaira, 
minimal properly elliptic, ruled over an elliptic curve or a torus.  
With the exception of some algebraic minimal properly elliptic surfaces
and the ruled surfaces over elliptic curves with $w_2\not=0$
all such surfaces admit $S^1$-actions without fixed points.

Conversely, the total space $E$ of an $S^1$-orbifold bundle $\xi$ over a 
geometric 3-orbifold is geometric, except when the base $B$ has geometry 
$\mathbb{H}^3$ or $\widetilde{\mathbb{SL}}$ and the characteristic class 
$c(\xi)$ has infinite order.                                
More generally, $E$ has a (proper) geometric decomposition if and only if $B$
is a $\widetilde{\mathbb{SL}}$-orbifold and $c(\xi)$ has finite order or $B$
has a (proper) geometric decomposition and the restrictions of $c(\xi)$ to 
the hyperbolic pieces of $B$ each have finite order.

Total spaces of circle bundles over aspherical Seifert fibred 3-manifolds,
$\mathbb{S}ol^3 $-manifolds or $\mathbb{H}^3$-manifolds have characterizations
refining Theorem 4.12 and parallel to those of Theorem 13.2.

\begin{theorem}
Let $M$ be a closed $4$-manifold with fundamental group $\pi$.
Then: 
\begin{enumerate}
\item $M$ is $s$-cobordant to the total space $E$ of an 
$S^1$-bundle over an aspherical closed Seifert fibred $3$-manifold 
or a $\mathbb{S}ol^3 $-manifold if and only if $\chi (M)=0$ and $\pi$ 
has normal subgroups $A<B$ such that $A\cong\mathbb{Z}$, 
$\pi/A$ has finite cohomological dimension
and $B/A$ is abelian.
If $h(B/A)>1$ then $M$ is homeomorphic to $E$.

\item $M$ is $s$-cobordant to the total space of an 
$S^1$-bundle over the mapping torus of a self homeomorphism of
an aspherical surface if and only if $\chi(M)=0$ and $\pi$ has normal 
subgroups $A<B$ such that $A\cong\mathbb{Z}$, $\pi/A$ is torsion free, 
$B$ is finitely generated and $\pi/B\cong\mathbb{Z}$.

\item $M$ is $s$-cobordant to the total space of an $S^1$-bundle over a closed $\mathbb{H}^3$-manifold if and only if $\chi(M)=0$,
$\sqrt\pi\cong\mathbb{Z}$, $\pi/\sqrt\pi$ is torsion free
and has no noncyclic abelian subgroup,
and $\pi$ has a finitely generated normal subgroup $B$ such that $\sqrt\pi<B$
and $e(\pi/B)=2$.
\end{enumerate}
\end{theorem}
               
\begin{proof} (1)\qua The conditions are clearly necessary. 
If they hold then $h(\sqrt\pi)\geq h(B/A)+1\geq2$, and so $M$ is aspherical.
Since $\pi/A$ has finite cohomological dimension it is a $PD_3 $-group,
by Theorem 4.12.
Hence it is the fundamental group of a closed Seifert fibred 3-manifold
or $\mathbb{S}ol^3$-manifold, $N$ say, by Theorem 2.14.
If $h(\sqrt\pi)=2$ then $\sqrt\pi\cong\mathbb{Z}^2$, by Theorem 9.2.
Hence $B/A\cong\mathbb{Z}$ and $[\sqrt\pi:B]<\infty$.
Since conjugation in $\pi$ preserves the chain of normal subgroups
$A<B\leq\sqrt\pi$ it follows that $C_\pi(\sqrt\pi)$ has finite index in $\pi$.
Hence $M$ is $s$-cobordant to an $\mathbb{H}^2\times\mathbb{E}^2$-
or $\widetilde{\mathbb{S}}\times\mathbb{E}^1$-manifold $E$, 
by Theorems 9.5, 9.6 and 9.12.
The Seifert fibration of $E$ factors through an
an $S^1$-bundle over $N$.
If $h(\sqrt\pi)>2$ then $\pi$ is virtually poly-$Z$. 
Hence $N$ is a $\mathbb{E}^3$-, $\mathbb{N}il^3$- or 
$\mathbb{S}ol^3 $-manifold and $M$ is homeomorphic to $E$, by Theorem 6.11.

(2)\qua The conditions are again necessary.
If they hold then $B/A$ is infinite, so $B$ has one end and
hence is a $PD_3$-group, by Theorem 4.5.
Since $B/A$ is torsion-free it is a $PD_2$-group, by Bowditch's Theorem,
and so $\pi/A$ is the fundamental group of a mapping torus, $N$ say. 
As $Wh(\pi)=0$, by Theorem 6.4, $M$ is simple homotopy equivalent to the 
total space $E$ of an $S^1$-bundle over $N$.
Since $\pi\times\mathbb{Z}$ is square root closed accessible
$M\times S^1$ is homeomorphic to $E\times S^1$ \cite{[Ca73]},
and so $M$ is $s$-cobordant to $E$.

(3)\qua The conditions are necessary, by the Virtual Fibration Theorem 
\cite{[Ag13]}.
If they hold then, $\pi$ has a subgroup $\bar\pi$ of index $\leq2$
such that $B<\bar\pi$ and $\bar\pi/B\cong\mathbb{Z}$.
Hence $B/\sqrt\pi$ is a $PD_2$-group, as in (2), and $\bar\pi/\sqrt\pi$
is the fundamental group of a closed $\mathbb{H}^3$-manifold.
Hence $\pi/\sqrt\pi$ is also the fundamental group of a closed
$\mathbb{H}^3$-manifold, $N$ say \cite{[Zn86]}. 
As $Wh(\pi)=0$, by Theorem 6.4, $M$ is simple homotopy equivalent to the 
total space $E$ of an $S^1$-bundle over $N$.
Hence $M$ is $s$-cobordant to $E$, by Theorem 10.7 of \cite{[FJ89]}.
\end{proof} 

Simple homotopy equivalence implies $s$-cobordism
for such bundles over aspherical 3-manifolds,
using \cite{[Ro11]}.
However we do not yet have good intrinsic characterizations of the fundamental 
groups of such 3-manifolds.

\section{Symplectic structures}

Let $M$ be a closed orientable 4-manifold.
If $M$ fibres over an orientable surface $B$ and the image 
of the fibre in $H_2(M;\mathbb{R})$ is nonzero
then $M$ has a symplectic structure \cite{[Th76]}.
The tangent bundle along the fibres is an $SO(2)$-bundle on $M$
which restricts to the tangent bundle of each fibre,
and so the homological condition is automatic unless the fibre is a torus.
This condition is also necessary if the base $B$ has genus $g>1$ \cite{[Wc05]},
but if $B$ is also a torus every such $M$ has 
a symplectic structure \cite{[Ge92]}.
Theorem 4.9 of \cite{[Wc05]} implies that any such bundle space $M$
is symplectic if and only if it is virtually symplectic.
If $N$ is an orientable 3-manifold then $N\times{S^1}$ admits 
a symplectic structure if and only if $N$ is a mapping torus \cite{[FV08]}. 

If $M$ admits one of the geometries
$\mathbb{CP}^2$, $\mathbb{S}^2\times\mathbb{S}^2$, 
$\mathbb{S}^2\times\mathbb{E}^2$, $\mathbb{S}^2\times\mathbb{H}^2$,
$\mathbb{H}^2\times\mathbb{E}^2$, $\mathbb{H}^2\times\mathbb{H}^2$ or
$\mathbb{H}^2(\mathbb{C})$ 
then it has a 2-fold cover which is K\"ahler, and therefore symplectic.
If $M$ admits one of the geometries $\mathbb{E}^4$, $\mathbb{N}il^4$, 
$\mathbb{N}il^3\times\mathbb{E}^1$ or $\mathbb{S} ol^3\times\mathbb{E}^1$ 
then it is symplectic if and only if $\beta_2(M)\not=0$.
For this condition is clearly necessary,
while if $\beta_2(M)>0$ then $\beta_1(M)\geq2$ 
and so $M$ fibres over a torus.
In particular,
an $\mathbb{E}^4$-manifold is symplectic if and only if it
is one of the eight flat K\"ahler surfaces,
while a $\mathbb{N}il^4$-manifold is symplectic if and only if it is a
nilmanifold, i.e., $\sqrt\pi=\pi$.
Every such infrasolvmanifold is virtually symplectic.

As any closed orientable manifold with one of the geometries
$\mathbb{S}^4$, $\mathbb{S}^3\times\mathbb{E}^1$,
$\mathbb{S} ol^4_{m,n}$ (with $m\not=n$), 
$\mathbb{S}ol^4_0$ or $\mathbb{S}ol^4_1$ has $\beta_2=0$, 
no such manifold is symplectic.
A closed $\widetilde{\mathbb{SL}}\times\mathbb{E}^1$-manifold 
is virtually a product $N\times{S^1}$, 
but the 3-manifold factor is not a mapping torus,
and so no such manifold is symplectic \cite{[Et01],[FV08]}.

All $\mathbb{H}^3\times\mathbb{E}^1$-manifolds are virtually symplectic,
since they are virtually products, by Theorem 9.3, 
and $\mathbb{H}^3$-manifolds are virtually mapping tori \cite{[Ag13]}.

The issue is less clear for the geometry $\mathbb{H}^4$. 
Symplectic 4-manifolds with index 0 have Euler characteristic 
divisible by 4, by Corollary 10.1.10 of \cite{[GS]},
and so covering spaces of odd degree of the Davis 120-cell
space are nonsymplectic.

Let $p:S\to{B}$ be a Seifert fibration with hyperbolic base,
action $\alpha$ and Euler class $e^\mathbb{Q}(p)$.
Then it follows from \cite{[Wc05]} that $S$ is virtually symplectic 
if and only if either $\mathrm{Im}(\alpha)$ is finite and $e^\mathbb{Q}(p)=0$, 
or $\mathrm{Im}(\alpha)$ is an infinite unipotent group 
which leaves $\pm{e^\mathbb{Q}(p)}$ invariant, 
or $tr(\alpha(g))>2$ for some $g\in\pi_1^{orb}(B)$.
Which Seifert fibred 4-manifolds are symplectic?

If $M$ is symplectic and has a nontrivial $S^1$-action 
then either $M$ is rational or ruled, 
or the action is fixed-point free, 
with orbit space a mapping torus \cite{[Bo14]}.

Which orientable bundle spaces over nonorientable base surfaces are symplectic?
Which 4-dimensional mapping tori and $S^1$-bundle spaces are symplectic?
The question of which $S^1$-bundle spaces are virtually symplectic 
is largely settled in \cite{[BF15]}.

%% file: m5-14.tex
\part{2-Knots}

\chapter{Knots and links}

In this chapter we introduce the basic notions and constructions of knot theory.
Many of these apply equally well in all dimensions, and for the most part
we have framed our definitions in such generality, although our
main concern is with 2-knots (embeddings of $S^2$ in $S^4$).
In particular, we show how the classification of higher dimensional knots 
may be reduced (essentially) to the classification of certain closed manifolds,
and we give Kervaire's characterization of high dimensional knot groups.
In the final sections we comment briefly on links 
and the groups of links, homology spheres and their groups.

\section{Knots}

The standard orientation of $\mathbb{R}^n$ induces 
an orientation on the unit $n$-disc
$D^n=\{ (x_1,\dots x_n)\in\mathbb{R}^n\mid \Sigma x_i^2\leq1\}$ 
and hence on its boundary $S^{n-1}=\partial D^n$, 
by the convention ``outward normal first".
We shall assume that standard discs and spheres have such orientations.
Qualifications shall usually be omitted when there is no risk of ambiguity.
In particular, we shall often abbreviate
$X(K)$, $M(K)$ and $\pi K$ (defined below) as $X$, $M$ and $\pi$, respectively.
                                                       
An $n$-{\it knot} is a locally flat embedding $K:S^n\to S^{n+2}$. 
(We shall also use the terms ``classical knot" when $n=1$, 
``higher dimensional knot"
when $n\geq2$ and ``high dimensional knot" when $n\geq3$.)
It is determined up to (ambient) isotopy by its image $K(S^n)$,
considered as an oriented codimension 2 submanifold of $S^{n+2}$,
and so we may let $K$ also denote this submanifold.
Let $r_n$ be an orientation reversing self homeomorphism of $S^n$.
Then $K$ is {\it invertible, $+$amphicheiral} or {\it $-$amphicheiral} 
if it is isotopic to $K\rho=Kr_n$, $rK=r_{n+2}K$ or $-K=rK\rho$, respectively.
An $n$-knot is {\it trivial} if it is isotopic to the composite of equatorial inclusions
$S^n\subset S^{n+1}\subset S^{n+2}$.
                                     
Every knot has a product neighbourhood: there is an embedding $j:S^n\times D^2$ 
onto a closed neighbourhood $N$ of $K$, such that $j(S^n\times\{0\})=K$ and 
$\partial N$ is bicollared in $S^{n+2}$ \cite{[KS75],[FQ]}. 
We may assume that $j$ is orientation preserving.
If $n\geq2$ it is then unique up to isotopy {\sl rel} $S^n\times\{0\}$.
The {\it exterior} of $K$ is the compact $(n+2)$-manifold
$X(K)=S^{n+2}\setminus{int N}$ with boundary 
$\partial X(K)\cong S^n\times S^1$, 
and is well defined up to homeomorphism.
It inherits an orientation from $S^{n+2}$.
An $n$-knot $K$ is trivial if and only if $X(K)\simeq S^1$; this
follows from Dehn's Lemma if $n=1$, is due to Freedman if $n=2$ 
(\cite{[FQ]} -- see Corollary 17.1.1 below) 
and is an easy consequence of the $s$-cobordism theorem if $n\geq3$.

The {\it knot group} is $\pi K=\pi_1(X(K))$.
An oriented simple closed curve isotopic to the oriented boundary of 
a transverse disc $\{ j\}\times S^1$ is called a {\it meridian} for $K$, 
and we shall also use this term
to denote the corresponding elements of $\pi$.
If $\mu$ is a meridian for $K$, represented by a simple closed curve on
$\partial X$ then $X\cup_\mu D^2$ is a deformation retract 
of $S^{n+2}\setminus\{*\}$ and so is contractible. 
Hence $\pi$ is generated by the conjugacy class of its meridians.

Assume for the remainder of this section that $n\geq2$. 
The group of pseudoisotopy classes of self homeomorphisms of $S^n\times S^1$ 
is $(Z/2Z)^3$, generated by reflections in either factor 
and by the map $\tau$ given by
$\tau(x,y)=(\rho(y)(x),y)$ for all $x$ in $S^n$ and $y$ in $S^1$, 
where $\rho:S^1\to SO(n+1)$ is
an essential map \cite{[Gl62],[Br67],[Kt69]}.
As any self homeomorphism of $S^n\times S^1$ extends across $D^{n+1}\times S^1$ 
the {\it knot manifold} $M(K)=X(K)\cup (D^{n+1}\times S^1)$
obtained from $S^{n+2}$ by surgery on $K$ is well defined, 
and it inherits an orientation from $S^{n+2}$ via $X$.
Moreover $\pi_1(M(K))\cong\pi K$ and $\chi(M(K))=0$.
Conversely, suppose that $M$ is a closed orientable 4-manifold with $\chi(M)=0$
and $\pi_1 (M)$ is generated by the conjugacy class of a single element.
(Note that each conjugacy class in $\pi$ corresponds to 
an unique isotopy class of oriented simple closed curves in $M$.)
Surgery on a loop in $M$ representing such an element gives a 1-connected 4-manifold $\Sigma$
with $\chi(\Sigma)=2$ which is thus homeomorphic to $S^4 $ and which contains
an embedded 2-sphere as the cocore of the surgery.
We shall in fact study 2-knots through such 4-manifolds,
as it is simpler to consider closed manifolds rather than pairs.

There is however an ambiguity when we attempt to recover $K$ from $M=M(K)$.
The cocore $\gamma=\{0\}\times S^1\subset D^{n+1}\times S^1\subset M$
of the original surgery is well defined up to isotopy by the conjugacy class of a meridian in 
$\pi K=\pi_1(M)$.
(In fact the orientation of $\gamma$ is irrelevant for what follows.)
Its normal bundle is trivial, so $\gamma$ has a product neighbourhood, $P$ say,
and we may assume that $M\setminus{intP}=X(K)$.
But there are two essentially distinct ways of identifying $\partial X$ with 
$S^n\times S^1=\partial(S^n\times D^2)$, modulo self homeomorphisms of
$S^n\times S^1$ that extend across $S^n\times D^2$.
If we reverse the original construction of $M$ we recover
$(S^{n+2},K)=(X\cup_j S^n\times D^2,S^n\times\{0\})$.
If however we identify $S^n\times S^1$ with $\partial X$ by means of $j\tau$
we obtain a new pair 
\begin{equation*}
(\Sigma,K^*)=(X\cup_{j\tau}S^n\times D^2,S^n\times\{0\}).
\end{equation*}
It is easily seen that $\Sigma\simeq S^{n+2}$, and hence $\Sigma\cong S^{n+2}$.
We may assume that the homeomorphism is orientation preserving.
Thus we obtain a new $n$-knot $K^*$, which we shall call the {\it Gluck reconstruction} of $K$.
The knot $K$ is {\it reflexive} if it is determined as an unoriented submanifold by its exterior,
i.e., if $K^*$ is isotopic to $K$, $rK$, $K\rho$ or $-K$.

If there is an orientation preserving homeomorphism from $X(K_1)$ to
$X(K)$ then $K_1$ is isotopic to $K$, $K^*$, $K\rho$ or $K^*\rho$.
If the homeomorphism also preserves the homology class of the meridians then $K_1$
is isotopic to $K$ or to $K^*$.
Thus $K$ is determined up to an ambiguity of order at most 2 by $M(K)$ together 
with the conjugacy class of a meridian.

A {\it Seifert hypersurface} for $K$ is a locally flat, oriented codimension 1 
submanifold $V$ of $S^{n+2}$ with (oriented) boundary $K$.
By a standard argument these always exist.
(Using obstruction theory it may be shown that the projection $pr_2j^{-1}:\partial X\to S^n\times S^1\to S^1$
extends to a map $p:X\to S^1$ \cite{[Ke65]}.
By topological transversality we may assume that $p^{-1}(1)$ is a bicollared, proper codimension 1
submanifold of $X$. The union $p^{-1}(1)\cup j(S^n\times [0,1])$
is then a Seifert hypersurface for $K$.)
We shall say that $V$ is {\it minimal} if the natural homomorphism
from $\pi_1(V)$ to $\pi K$ is a monomorphism.

In general there is no canonical choice of Seifert surface.
However there is one important special case.
An $n$-knot $K$ is {\it fibred} if there is such a map $p:X\to S^1$
which is the projection of a fibre bundle. (Clearly $K^*$ is then fibred also.)
The exterior is then the mapping torus of a self homeomorphism $\theta$ of the 
fibre $F$ of $p$, called the {\it (geometric) monodromy} of the bundle.
Such a map $p$ extends to a fibre bundle projection $q:M(K)\to S^1$,
with fibre $\widehat F=F\cup D^{n+1}$, called the {\it closed fibre} of $K$.
Conversely, if $M(K)$ fibres over $S^1$ then the cocore $\gamma$ is homotopic
(and thus isotopic) to a cross-section of the bundle projection, and so
$K$ is fibred.
If the monodromy is represented by a self-homeomorphism of finite order
then it has nonempty fixed point set, and the closed monodromy 
$\widehat\theta$ has finite order.
However the results of \cite{[Hn]} and \cite{[La]} may be used 
to show that the closed monodromy of the spun trefoil knot $\sigma3_1$ 
has finite order, 
but as $\pi_1(F)\cong F(2)$ has no automorphism of order 6 
\cite{[Me74]}
there is no representative of finite order with nonempty fixed point set.

\section{Covering spaces}

Let $K$ be an $n$-knot. Then 
$H_1(X(K);\mathbb{Z})\cong\mathbb{Z}$ and $H_i(X(K);\mathbb{Z})=0$ if $i>1$, 
by Alexander duality.
The meridians are all homologous and generate $\pi/\pi'=H_1(X;\mathbb{Z})$, 
and so determine a canonical isomorphism with $\mathbb{Z}$. 
Moreover $H_2(\pi;\mathbb{Z})=0$, since it is a quotient of $H_2(X;\mathbb{Z})=0$.

We shall let $X'(K)$ and $M'(K)$ denote the covering spaces corresponding to the commutator subgroup.
(The cover $X'/X$ is also known as the {\it infinite cyclic cover} of the knot.)
Since $\pi/\pi'=\mathbb{Z}$ the (co)homology groups of $X'$ 
are modules over the group ring $\mathbb{Z}[\mathbb{Z}]$, 
which may be identified with the ring of integral Laurent polynomials
$\Lambda=\mathbb{Z}[t,t^{-1}]$. 
If $A$ is a $\Lambda$-module, let $zA$ be the $\mathbb{Z}$-torsion submodule,
and let $e^iA=Ext_\Lambda^i(A,\Lambda)$.
                                                
Since $\Lambda$ is noetherian the (co)homology of a finitely generated free
$\Lambda$-chain complex is finitely generated.
The Wang sequence for the projection of $X'$ onto $X$ may be identified with
the long exact sequence of homology corresponding to the exact sequence of 
coefficients 
\begin{equation*}
0\to \Lambda\to\Lambda\to\mathbb{Z}\to 0.
\end{equation*}
Since $X$ has the homology of a circle it follows easily that
multiplication by $t-1$ induces automorphisms of the modules
$H_i(X;\Lambda)$ for $i>0$.  Hence these homology modules are all
finitely generated torsion $\Lambda$-modules.  It follows that
$Hom_\Lambda(H_i(X;\Lambda),\Lambda)$ is 0 for all $i$, and the UCSS
collapses to a collection of short exact sequences
\begin{equation*}
0\to e^2H_{i-2}(X;\Lambda)\to H^i(X;\Lambda)\to e^1H_{i-1}(X;\Lambda)\to 0.
\end{equation*}
The infinite cyclic covering spaces $X'$ and $M'$ behave homologically 
much like $(n+1)$-manifolds,
at least if we use field coefficients \cite{[Mi68],[Ba80]}.
If $H_i(X;\Lambda)=0$ for $1\leq i\leq (n+1)/2$ then $X'$ is acyclic; 
thus if also $\pi=\mathbb{Z}$ then $X\simeq S^1$ and so $K$ is trivial.
All the classifications of high dimensional knots to date 
assume that $\pi=\mathbb{Z}$ and that $X'$ is highly connected. 

When $n=1$ or 2 knots with $\pi=\mathbb{Z}$ are trivial, 
and it is more profitable to work 
with the universal cover $\widetilde X$ (or $\widetilde M$).
In the classical case $\widetilde X$ is contractible \cite{[Pa57]}. 
In higher dimensions $X$ is aspherical only when 
the knot is trivial \cite{[DV73]}.
Nevertheless the closed 4-manifolds $M(K)$ obtained by surgery on 2-knots are 
often aspherical.
(This asphericity is an additional reason for choosing to work with $M(K)$ 
rather than $X(K)$.)

\section{Sums, factorization and satellites}

The {\it sum} of two knots $K_1$ and $K_2$ may be defined (up to isotopy) 
as the $n$-knot $K_1\sharp K_2$ obtained as follows.
Let $D^n(\pm)$ denote the upper and lower hemispheres of $S^n$.
We may isotope $K_1$ and $K_2$ so that each $K_i(D^n(\pm))$ contained in
$D^{n+2}(\pm)$, $K_1(D^n(+))$ is a trivial $n$-disc in $D^{n+2}(+)$,
$K_2(D^n(-))$ is a trivial $n$-disc in $D^{n+2}(-)$ and $K_1|_{S^{n-1}}=K_2|_{S^{n-1}}$
(as the oriented boundaries of the images of $D^n(-)$).
Then we let $K_1\sharp K_2=K_1|_{D^n(-)}\cup K_2|_{D^n(+)}$.
By van Kampen's theorem $\pi (K_1\sharp K_2)=\pi K_1*_{\mathbb{Z}}\pi K_2$
where the amalgamating subgroup is generated by a meridian in each knot group.
It is not hard to see that $X'(K_1\sharp K_2)\simeq X'(K_1)\vee X'(K_2)$,
and so $\pi(K_1\sharp K_2)'\cong \pi(K_1)'*\pi(K_2)'$.

The knot $K$ is {\it irreducible} if it is not the sum of two nontrivial knots.
Every knot has a finite factorization into irreducible knots 
\cite{[DF87]}. 
(For 1- and 2-knots whose groups have finitely generated commutator subgroups
this follows easily from the Grushko-Neumann theorem 
on factorizations of groups as free products.)
In the classical case the factorization is essentially unique, 
but if $n\geq3$ there are $n$-knots with several distinct 
such factorizations \cite{[BHK81]}.
Almost nothing is known about the factorization of 2-knots.

If $K_1$ and $K_2$ are fibred then so is their sum, and the closed fibre of 
$K_1\sharp K_2$
is the connected sum of the closed fibres of $K_1$ and $K_2$.
However in the absence of an adequate criterion for a 2-knot to fibre, 
we do not know whether every summand of a fibred 2-knot is fibred.
In view of the unique factorization theorem for oriented 3-manifolds 
we might hope that there would be a similar theorem for fibred 2-knots.
However the closed fibre of an irreducible 2-knot need not be an irreducible 3-manifold.
(For instance, the Artin spin of a trefoil knot is an irreducible fibred 2-knot,
but its closed fibre is $(S^2\times S^1)\sharp(S^2\times S^1)$.)

A more general method of combining two knots is the process of forming satellites.
Although this process arose in the classical case, where it is intimately connected with
the notion of torus decomposition, 
we shall describe only the higher-dimensional
version of \cite{[Kn83]}.                                        
Let $K_1$ and $K_2$ be $n$-knots (with $n\geq2$) and let $\gamma$ be
a simple closed curve in $X(K_1)$, with a product neighbourhood $U$.
Then there is a homeomomorphism $h$ which carries 
$S^{n+2}\setminus{int}U\cong{S^n\times{D^2}}$
onto a product neighbourhood of $K_2$.
The knot $\Sigma(K_2;K_1,\gamma)=hK_1$ is
called the {\it satellite} of $K_1$ about $K_2$ relative to $\gamma$.
We also call $K_2$ a {\it companion} of $hK_1$.
If either $\gamma=1$ or $K_2$ is trivial then $\Sigma(K_2;K_1,\gamma)=K_1$.
If $\gamma$ is a meridian for $K_1$ then $\Sigma(K_2;K_1,\gamma)=K_1\sharp K_2$.
If $\gamma$ has finite order in $\pi K_1$ let $q$ be that order; 
otherwise let $q=0$.
Let $w$ be a meridian in $\pi K_2$. 
Then $\pi\Sigma(K_2;K_1,\gamma)\cong 
(\pi K_2/\langle\langle w^q\rangle\rangle)*_{Z/qZ}\pi K_1$, 
where $w$ is identified with $\gamma$ in $\pi K_1$,
by Van Kampen's theorem.

\section{Spinning, twist spinning and deform spinning}

The first nontrivial examples of higher dimensional knots
were given by Artin \cite{[Ar25]}.
We may paraphrase his original idea as follows.
As the half space 
$\mathbb{R}^3_+=\{ (w,x,y,z)\in\mathbb{R}^4\mid w=0,z\geq0\}$
is spun about the axis $A=\{(0,x,y,0)\}$ it sweeps out $\mathbb{R}^4$,
and any arc in $\mathbb{R}^3_+$ with endpoints on $A$ sweeps out a 2-sphere.

This construction has been extended, 
first to twist-spinning and roll-spinning \cite{[Fo66],[Ze65]}, 
and then more generally to deform spinning \cite{[Li79]}.
Let $g$ an orientation preserving self-homeomorphism 
of $S^{n+2}$ which is the identity on the $n$-knot $K$. 
If $g$ does not twist the normal bundle of $K$ in $S^{n+2}$
then the sections of $M(g)$ determined by points of $K$
have canonical ``constant" framings,
and surgery on such a section in the pair $(M(g),K\times{S^1})$ gives
an $(n+1)$-knot, called the {\it deform spin\/} of $K$ determined by $g$.
The deform spin is {\it untwisted\/} if $g$ preserves 
a Seifert hypersurface for $K$.
If $g$ is the identity this gives the {\it Artin spin} $\sigma{K}$, 
and $\pi\sigma{K}=\pi{K}$.

Twist spins are defined by maps supported in a collar of $\partial{X}={K\times{S^1}}$.
(If $n=1$ we use the 0-framing.)
Let $r$ be an integer.
The self-map  $t_r$ of $S^{n+2}$ defined 
by $t_r(k,z,x)=(k, e^{2\pi{irx}}z,x)$ on $\partial{X}\times[0,1]$ 
and the identity elsewhere gives the $r$-{\it twist spin} $\tau_rK$.
Clearly $\tau_0K=\sigma{K}$.
The group of $\tau_r K$ is obtained from $\pi K$ by adjoining the relation
making the $r^{th}$ power of (any) meridian central.
Zeeman discovered the remarkable fact that if $r\not=0$ 
then $\tau_r K$ is fibred, 
with closed fibre the $r$-fold cyclic branched cover of $S^{n+2}$, 
branched over $K$,
and monodromy the canonical generator of the group of covering transformations \cite{[Ze65]}.
Hence $\tau_1K$ is always trivial.
More generally,
if $g$ is an {\it un}twisted deformation of $K$ and $r\not=0$
then the knot determined by $t_rg$ is fibred \cite{[Li79]}.
(See also \cite{[GK78],[Mo83],[Mo84]} and \cite{[Pl84']}.)
Twist spins of $-$amphicheiral knots are $-$amphicheiral, while twist spinning
interchanges invertibility and $+$amphicheirality \cite{[Li85]}.

If $K$ is a classical knot the factors of the closed fibre of $\tau_r K$ are 
the cyclic branched covers of the prime factors of $K$, 
and are Haken, hyperbolic or Seifert fibred.
With some exceptions for small values of $r$, the factors are aspherical, 
and $S^2\times S^1$ is never a factor \cite{[Pl84]}.          
If $r>1$ and $K$ is nontrivial then $\tau_r K$ is nontrivial, 
by the Smith Conjecture. 
If $K$ is a deform spun 2-knot then the order ideal of
$H_1(\pi{K};\mathbb{Q}\Lambda)$ is invariant under the
involution $t\mapsto{t^{-1}}$ \cite{[BM09]}.
  
\section{Ribbon and slice knots}

Two $n$-knots $K_0$ and $K_1$ are {\it concordant} if there is a 
locally flat embedding
$\mathcal{K}:S^n\times [0,1]\to S^{n+2}\times [0,1]$ 
such that $K_i=\mathcal{K}\cap{S^n}\times\{i\}$ for $i=0,1$.
They are {\it $s$-concordant} if there is a concordance whose exterior is an 
$s$-cobordism ({\it rel} $\partial$) from $X(K_0)$ to $X(K_1)$.
(If $n>2$ this is equivalent to ambient isotopy, by the $s$-cobordism theorem.) 

An $n$-knot $K$ is a {\it slice} knot if it is concordant to the unknot;
equivalently, if it bounds a properly embedded $(n+1)$-disc $\Delta$ in $D^{n+3}$. 
Such a disc is called a {\it slice disc} for $K$.
Doubling the pair $(D^{n+3},\Delta)$ gives an $(n+1)$-knot which meets 
the equatorial $S^{n+2}$ of $S^{n+3}$ transversally in $K$; 
if the $(n+1)$-knot can be chosen to be trivial then $K$ is {\it doubly slice}.
All even-dimensional knots are slice \cite{[Ke65]}, 
but not all slice knots are doubly 
slice, and no adequate criterion is yet known.
The sum $K\sharp -K$ is a slice of $\tau_1K$ and so is doubly slice 
\cite{[Su71]}.
Twist spins of doubly slice knots are doubly slice.
                         
An $n$-knot $K$ is a {\it ribbon} knot if it is the boundary of an immersed 
$(n+1)$-disc $\Delta$ in $S^{n+2}$ whose only singularities are transverse 
double points, the double point sets being a disjoint union of discs. 
Given such a ``ribbon" $(n+1)$-disc $\Delta$ in $S^{n+2}$ the cartesian product
$\Delta\times D^p\subset S^{n+2}\times D^p\subset S^{n+2+p}$ determines a 
ribbon $(n+1+p)$-disc in $S^{n+2+p}$.
All higher dimensional ribbon knots derive from ribbon 1-knots by this process 
\cite{[Yn77]}.
As the $p$-disc has an orientation reversing involution, 
this easily imples that all 
ribbon $n$-knots with $n\geq2$ are $-$amphicheiral.
The Artin spin of a 1-knot is a ribbon 2-knot. 
Each ribbon 2-knot has a Seifert hypersurface which is a once-punctured 
connected sum of copies of $S^1\times S^2$ \cite{[Yn69]}.
Hence such knots are reflexive.
(See \cite{[Su76]} for more on geometric properties of such knots.)
                                      
An $n$-knot $K$ is a {\it homotopy ribbon} knot if it is a slice knot with 
a slice disc whose exterior $W$ has a handlebody 
decomposition consisting of 0-, 1- and 2-handles.
The dual decomposition of $W$ relative to $\partial W=M(K)$ has only handles
of index $\geq{n+1}$, and so $(W,M)$ is $n$-connected.
(The definition of ``homotopically ribbon" for 1-knots
used in Problem 4.22 of \cite{[GK]} requires only that 
this latter condition be satisfied.)
More generally, we shall say that $K$ is $\pi_1$-{\it slice} if the inclusion 
of $X(K)$ into the exterior of some slice disc induces an isomorphism on
fundamental groups.

Every ribbon knot is homotopy ribbon and hence slice \cite{[Hi79]},
while if $n\geq2$ every homotopy ribbon $n$-knot is $\pi_1$-slice. 
Nontrivial classical knots are never $\pi_1$-slice, 
since the longitude of a slice knot is nullhomotopic 
in the exterior of a slice disc.    
It is an open question whether every classical slice knot is ribbon.
However in higher dimensions ``slice" does not even imply ``homotopy ribbon".
(The simplest example is $\tau_23_1$ - see below.)

Most of the conditions considered here depend only on the $h$-cobordism 
class of $M(K)$.
An $n$-knot $K$ is slice if and only if $M=\partial{W}$, where $W$ is an
homology $S^1\times{D^{n+2}}$ and the image of $\pi{K}$ normally generates
$\pi_1(W)$, and it is $\pi_1$-slice if and only if we may assume also that the
inclusion of $M$ into $\partial{W}$ induces an isomorphism on $\pi_1$.
The knot $K$ is doubly slice if and only if  
$M$ embeds in $S^1\times{S^{n+2}}$ via a map which induces an isomorphism on
first homology.                                                                  

\section{The Kervaire conditions}

A group $G$ has {\it weight} 1 if it has an element whose conjugates generate $G$.
Such an element is called a {\it weight element} for $G$, 
and its conjugacy class is called a {\it weight class} for $G$.
If $G$ is solvable then it has weight 1 if and only if $G/G'$ is cyclic, 
for a solvable group with trivial abelianization must be trivial. 

If $\pi$ is the group of an $n$-knot $K$ then 
{\sl
\begin{enumerate}
\item $\pi$ is finitely presentable;

\item $\pi$ is of weight 1;

\item $H_1(\pi;\mathbb{Z})=\pi/\pi'\cong\mathbb{Z}$; and

\item $H_2(\pi;\mathbb{Z})=0$.
\end{enumerate}}
Kervaire showed that any group satisfying these conditions is 
an $n$-knot group, for every $n\geq3$ \cite{[Ke65]}. 
These conditions are also necessary when $n=1$ or 2, 
but are then no longer sufficient,
and there are as yet no corresponding characterizations 
for 1- and 2-knot groups.
If (4) is replaced by the stronger condition that $\mathrm{def}(\pi)=1$ then 
$\pi$ is a 2-knot group, but this condition is not necessary \cite{[Ke65]}.
(See \S 9 of this chapter, \S4 of Chapter 15 and \S4 of Chapter 16 for
examples with deficiency $\leq0$.)
Gonzalez-Acu\~na has given a characterization of
2-knot groups as groups admitting certain presentations \cite{[GA94]}.
(Note also that if $\pi$ is a high dimensional knot group then $q(\pi)\geq0$, 
and $q(\pi)=0$ if and only if $\pi$ is a 2-knot group.) 

Every knot group has a Wirtinger presentation, i.e., one in which
the relations are all of the form $x_j=w_jx_0w_j^{-1}$,
where $\{ x_i, 0\leq i\leq n\}$ is the generating set \cite{[Yj70]}.
If $K$ is a nontrivial 1-knot then $\pi{K}$ has 
a Wirtinger presentation of deficiency 1.
A group has such a presentation if and only if it has weight 1 and has a 
deficiency 1 presentation $P$ such that the presentation of the trivial 
group obtained by adjoining the relation killing a weight element is 
AC-equivalent to the empty presentation \cite{[Yo82']}.
Any such group is the group of a 2-knot which is a smooth embedding
in the standard smooth structure on $S^4$ \cite{[Le78]}.
The group of a nontrivial 1-knot $K$ has one end \cite{[Pa57]}, 
so $X(K)$ is aspherical, 
and $X(K)$ collapses to a finite 2-complex, so $g.d.\pi{K}=2$. 
If $\pi$ is an $n$-knot group then $g.d.\pi=2$ if and only if 
$c.d.\pi=2$ and $\mathrm{def}(\pi)=1$, by Theorem 2.8.

Since the group of a homotopy ribbon $n$-knot (with $n\geq2$) 
is the fundamental group of a $(n+3)$-manifold $W$ with $\chi(W)=0$ 
and which can be built with 0-, 1- and 2-handles only, 
such groups also have deficiency 1.
Conversely, if a finitely presentable group $\pi$ has weight 1 
and deficiency 1 then we may use such a presentation to construct 
a 5-dimensional handlebody $W=D^5\cup\{ h_i^1\}\cup\{ h_j^2\}$ 
with $\pi_1(\partial W)=\pi_1(W)\cong\pi$ and $\chi(W)=0$.
Adjoining another 2-handle $h$ along a loop representing a weight class for 
$\pi_1(\partial W)$ gives a homotopy 5-ball $B$ with 1-connected boundary.
Thus $\partial B\cong S^4$, and the boundary of the cocore of the 2-handle $h$
is clearly a homotopy ribbon 2-knot with group $\pi$.
(In fact any group of weight 1 with a Wirtinger presentation of deficiency 1
is the group of a ribbon $n$-knot, for each $n\geq2$ 
\cite{[Yj69]} -- see \cite{[H3]}.)

The deficiency may be estimated in terms of the minimum number 
of generators of the $\Lambda$-module $e^2(\pi'/\pi'')$. 
Using this observation, it may be shown that if $K$ is the sum 
of $m+1$ copies of $\tau_23_1$ then $\mathrm{def}(\pi K)=-m$ 
\cite{[Le78]}.
There are irreducible 2-knots whose groups have deficiency $-m$, 
for each $m\geq0$ \cite{[Kn83]}.

A knot group $\pi$ has two ends if and only if $\pi'$ is finite.
We shall determine all such 2-knots in \S1 of Chapter 15.
Nontrivial torsion-free knot groups have one end 
\cite{[Kl93]}.
There are also 2-knot groups with infinitely many ends \cite{[GM78]}.
The simplest is perhaps the group with presentation 
\begin{equation*}
\langle a,b,t\mid a^3=1,\medspace aba^{-1}=b^2,\medspace tat^{-1}=a^2\rangle.
\end{equation*} 
The first two relations imply that $b^7=1$,
and so this group is an HNN extension of  
$\langle{a,b}\rangle\cong{Z/7Z\rtimes_2Z/3Z}$, 
with associated subgroups both $\langle{a}\rangle\cong{Z/3Z}$.
It is also the group of a satellite of $\tau_23_1$ 
with companion Fox's Example 10.
(The smallest possible HNN base for such a 2-knot group is $Q(16)$ \cite{[BH17]}.)

\section{Weight elements, classes and orbits}

Two $2$-knots $K$ and $K_1$ have homeomorphic exteriors if and only if there is 
a homeomorphism from $M(K_1)$ to $M(K)$ which carries the conjugacy class of a 
meridian of $K_1$ to that of $K$ (up to inversion).
In fact if $M$ is any closed orientable 4-manifold with $\chi(M)=0$ and with
$\pi=\pi_1(M)$ of weight 1 then surgery on a weight class gives a 2-knot with group $\pi$.
Moreover, if $t$ and $u$ are two weight elements and $f$ 
is a self homeomorphism of $M$ such that $u$ is conjugate 
to $f_*(t^{\pm1})$ then surgeries on $t$ and $u$ lead to
knots whose exteriors are homeomorphic 
(via the restriction of a self homeomorphism of $M$ isotopic to $f$).
Thus the natural invariant to distinguish between knots with isomorphic groups 
is not the weight class, but rather the orbit of the weight class under the
action of self homeomorphisms of $M$.
In particular, the orbit of a weight element under $Aut(\pi)$ is a 
well defined invariant, which we shall call the {\it weight orbit}.
If every automorphism of $\pi$ is realized by a self homeomorphism of $M$
then the homeomorphism class of $M$ and the weight orbit together form a
complete invariant for the (unoriented) knot, up to Gluck reconstruction.
(This is the case if $M$ is an infrasolvmanifold.)
                     
For oriented knots we need a refinement of this notion.
If $w$ is a weight element for $\pi$ then we shall call the set
$\{\alpha(w)\mid \alpha\in Aut(\pi),~\alpha(w)\equiv w~mod~\pi'\}$ 
a {\it strict} weight orbit for $\pi$.
A strict weight orbit determines a transverse orientation for the 
corresponding knot (and its Gluck reconstruction).
An orientation for the ambient sphere is determined by an orientation 
for $M(K)$.
If $K$ is invertible or $+$amphicheiral then there is 
a self homeomorphism of $M$ which is orientation preserving 
or reversing (respectively) and which reverses the 
transverse orientation of the knot, 
i.e., carries the strict weight orbit to its inverse.
Similarly, if $K$ is $-$amphicheiral there is an orientation reversing 
self homeomorphism of $M$ which preserves the strict weight orbit.

\begin{theorem}
Let $G$ be a group of weight $1$ and with $G/G'\cong\mathbb{Z}$. 
Let $t$ be an element of $G$ whose image generates $G/G'$ and let $c_t$ be the 
automorphism of $G'$ induced by conjugation by $t$. Then
\begin{enumerate}
\item $t$ is a weight element if and only if $c_t$ is meridianal;

\item two weight elements $t$, $u$ are in the same weight class if and only if 
there is an inner automorphism $c_g$ of $G'$ such that $c_u=c_gc_tc_g^{-1}$;

\item two weight elements $t$, $u$ are in the same strict weight orbit 
if and only if there is an automorphism $d$ of $G'$ 
such that $c_u=dc_td^{-1}$ and $dc_td^{-1}c_t^{-1}$ is 
an inner automorphism;

\item if $t$ and $u$ are weight elements then $u$ is conjugate to $(g''t)^{\pm1}$ for some $g''$ 
in $G''$.
\end{enumerate}
\end{theorem}

\begin{proof} The verification of (1-3) is routine.
If $t$ and $u$ are weight elements then, up to inversion, 
$u$ must equal $g't$ for some $g'$
in $G'$. Since multiplication by $t-1$ is invertible on $G'/G''$ 
we have $g'=khth^{-1}t^{-1}$
for some $h$ in $G'$ and $k$ in $G''$. Let $g''=h^{-1}kh$. 
Then $u=g't=hg''th^{-1}$.
\end{proof}

An immediate consequence of this theorem is that if $t$ and $u$ are in the same
strict weight orbit then $c_t$ and $c_u$ have the same order. 
Moreover if $C$ is the centralizer of $c_t$ in $Aut(G')$ then the strict weight orbit 
of $t$ contains at most $[Aut(G'):C.Inn(G')]\leq |Out(G')|$ weight classes.
In general there may be infinitely many weight orbits 
\cite{[Pl83']}.
However if $\pi$ is metabelian the weight class (and hence the weight orbit) 
is unique up to inversion, by part (4) of the theorem.

\section{The commutator subgroup}

It shall be useful to reformulate the Kervaire conditions in terms of the 
automorphism of the commutator subgroup induced by conjugation by a meridian.
An automorphism $\phi$ of a group $G$ is {\it meridianal} if 
$\langle\langle  g^{-1}\phi(g)\mid g\in G\rangle\rangle_G=G$.
If $H$ is a characteristic subgroup of $G$ and $\phi$ is meridianal
the induced automorphism of $G/H$ is then also meridianal.
In particular, $H_1(\phi)-1$ maps $H_1(G;\mathbb{Z})=G/G'$ onto itself.
If $G$ is solvable an automorphism satisfying the latter condition is 
meridianal, for a solvable perfect group is trivial.

It is easy to see that no group $G$ with $G/G'\cong\mathbb{Z}$ can have 
$G'\cong\mathbb{Z}$ or $D$.
It follows that the commutator subgroup of a knot group never has two ends. 

\begin{theorem}
{\rm[HK78],[Le78]}\qua 
A finitely presentable group $\pi$ is a high dimensional knot group 
if and only if $\pi\cong\pi'\rtimes_\theta\mathbb{Z}$ 
for some meridianal automorphism $\theta$ of $\pi'$ 
such that $H_2(\theta)-1$ 
is an automorphism of $H_2(\pi';\mathbb{Z})$.
\qed
\end{theorem}                                  

If $\pi$ is a knot group then $\pi'/\pi''$ is a finitely generated 
$\Lambda$-module.
Levine and Weber have made explicit the conditions under which 
a finitely generated $\Lambda$-module may be the commutator subgroup 
of a metabelian high dimensional knot group \cite{[LW78]}.
Leaving aside the $\Lambda$-module structure, 
Hausmann and Kervaire have characterized the finitely generated abelian groups 
$A$ that may be commutator subgroups of high dimensional knot groups 
\cite{[HK78]}.  
``Most" can occur; there are mild restrictions on 2- and 3-torsion, 
and if $A$ is infinite it must have rank at least 3.
We shall show that the abelian groups which are commutator subgroups of
2-knot groups are $\mathbb{Z}^3$, $\mathbb{Z}[\frac12 ]$ 
(the additive group of dyadic rationals)
and the cyclic groups of odd order.
(See Theorems 15.4 and 15.10.)
The commutator subgroup of a nontrivial classical knot group is never abelian.
                                                  
Hausmann and Kervaire also showed that any finitely generated abelian group
could be the centre of a high dimensional knot group \cite{[HK78']}. 
We shall show that the centre of a 2-knot group is either $\mathbb{Z}^2$,
torsion-free of rank 1, $\mathbb{Z}\oplus(Z/2Z)$ or is finite. 
(See Theorem 15.11.
In all known cases the centre is $\mathbb{Z}^2$, 
$\mathbb{Z}\oplus(Z/2Z)$, $\mathbb{Z}$, $Z/2Z$ or 1.)
A classical knot group has nontrivial centre if and only if 
the knot is a torus knot \cite{[BZ]}; the centre is then $\mathbb{Z}$.
    
Silver has given examples of high dimensional knot groups $\pi$ 
with $\pi'$ finitely generated but not finitely presentable 
\cite{[Si91]}. 
He has also shown that there are embeddings $j:T\to S^4$ such that 
$\pi_1(S^4\setminus{j(T)})'$ is finitely generated but not finitely presentable 
\cite{[Si97]}.
However no such 2-knot groups are known. 
If the commutator subgroup is finitely generated 
then it is the unique HNN base \cite{[Si96]}.
Thus knots with such groups have no minimal Seifert hypersurfaces.

The first examples of high dimensional knot groups which are not 2-knot 
groups made use of Poincar\'e duality with coefficients $\Lambda$.
Farber \cite{[Fa77]} and Levine \cite{[Le77]} independently found  
the following theorem.

\begin{theorem}
[Farber, Levine] 
Let $K$ be a $2$-knot and $A=H_1(M(K);\Lambda)$.
Then $H_2(M(K);\Lambda)\cong \overline{e^1A}$, and there is a nondegenerate 
$\mathbb{Z}$-bilinear pairing $[~,~]:zA\times zA\to\mathbb{Q}/\mathbb{Z}$ such that 
$[t\alpha,t\beta]=[\alpha,\beta]$ for all $\alpha$ and $\beta$ in $zA$.
\qed
\end{theorem}

Most of this theorem follows easily from Poincar\'e duality with coefficients $\Lambda$,
but some care is needed in order to establish the symmetry of the pairing.
When $K$ is a fibred 2-knot, with closed fibre $\widehat F$, 
the Farber-Levine pairing is just the standard linking pairing on the torsion 
subgroup of $H_1(\widehat F;\mathbb{Z})$, 
together with the automorphism induced by the monodromy.                          
In particular, Farber observed that although the group $\pi$ 
with presentation 
\begin{equation*}
\langle a,t\mid tat^{-1}=a^2\!,\medspace a^5=1\rangle
\end{equation*} 
is a high dimensional knot group, if $\ell$ is any nondegenerate 
$\mathbb{Z}$-bilinear pairing on $\pi'\cong Z/5Z$ with values in 
$\mathbb{Q}/\mathbb{Z}$ then $\ell(t\alpha,t\beta)=-\ell(\alpha,\beta)$ 
for all $\alpha$, $\beta$ in $\pi'$, 
and so $\pi$ is not a 2-knot group. 

\begin{cor}
{\rm[Le78]}\,
$H_2(\pi';\mathbb{Z})$ is a quotient of 
$\overline{Hom_\Lambda(\pi'/\pi'',\mathbb{Q}(t)/\Lambda)}$.
\qed
\end{cor}

Every orientation preserving meridianal automorphism 
of a torsion-free 3-manifold group is realizable by a fibred 2-knot.

\begin{theorem}
Let $N$ be a closed orientable $3$-manifold 
such that $\nu=\pi_1(N)$ is torsion-free.
If $K$ is a $2$-knot such that $(\pi K)'\cong\nu$ then $M(K)$ is 
homotopy equivalent to the mapping torus of a self homeomorphism of $N$.
If $\theta$ is a meridianal automorphism of $\nu$ then 
$\pi=\nu\rtimes_\theta\mathbb{Z}$ is the group of a fibred $2$-knot 
with fibre $N$ if and only if $\theta_*(c_{N*}[N])=c_{N*}[N]$.
\end{theorem}
                     
\begin{proof}
The first assertion follows from Corollary 4.5.4.

The classifying maps for the fundamental groups induce 
a commuting diagram involving the Wang sequences of $M(K)$ and $\pi$ 
from which the necessity of the orientation condition follows easily. 
(It is vacuous if $\nu$ is a free group.)

Let $N=P\sharp{R}$ where $P$ is a connected sum of $r$
copies of $S^1\times S^2$ and the summands of $R$ 
are aspherical.
If $\theta_*(c_{N*}[N])=c_{N*}[N]$ then $\theta$ may be realized by an 
orientation preserving self homotopy equivalence $g$ of $N$ 
\cite{[Sw74]}.
We may assume that $g$ is a connected sum of homotopy equivalences
between the irreducible factors of $R$ and a self homotopy
equivalence of $P$, by the Splitting Theorem of 
\cite{[HL74]}.
The factors of $R$ are either Haken, hyperbolic or Seifert-fibred,
by the Geometrization Conjecture (see \cite{[B-P]}),
and homotopy equivalences between such manifolds are homotopic 
to homeomorphisms, by \cite{[Hm]}, 
Mostow rigidity and \cite{[Sc83]}, respectively.
A similar result holds for $P=\sharp^r(S^1\times S^2)$, by \cite{[La]}.
Thus we may assume that $g$ is a self homeomorphism of $N$.
Surgery on a weight class in the mapping torus of $g$ gives 
a fibred 2-knot with closed fibre $N$ and group $\pi$.
\end{proof}

If $N$ is hyperbolic, Seifert fibred or if its prime factors 
are Haken or $S^1\times S^2$ then
the mapping torus is determined up to homeomorphism among fibred 4-manifolds 
by its homotopy type, 
since homotopy implies isotopy in each case, 
by Mostow rigidity, \cite{[Sc85],[BO91]} and \cite{[HL74]}, respectively.

Yoshikawa has shown that a finitely generated abelian group 
is the base of some HNN extension which is a high dimensional knot group 
if and only if it satisfies the restrictions on torsion of \cite{[HK78]}, 
while if a knot group has a non-finitely 
generated abelian base then it is metabelian.
Moreover a 2-knot group $\pi$ which is an HNN 
extension with abelian base is either metabelian or has base 
$\mathbb{Z}\oplus(Z/\beta Z)$ for some odd $\beta\geq1$ \cite{[Yo86],[Yo92]}.
We shall show that in the latter case $\beta$ must be 1, 
and so $\pi$ has a deficiency 1 presentation 
$\langle t,x\mid tx^nt^{-1}=x^{n+1}\rangle$.
(See Theorem 15.14.)
No nontrivial classical knot group is an HNN extension with abelian base.
(This is implicit in Yoshikawa's work, and can also be deduced from the
facts that classical knot groups have cohomological dimension $\leq2$
and symmetric Alexander polynomial.)

\section{Deficiency and geometric dimension}

J.H.C.Whitehead raised the question 
``is every subcomplex of an aspherical 2-complex also aspherical?"
This is so if the fundamental group of the subcomplex is a 1-relator group
\cite{[Go81]} or is locally indicable \cite{[Ho82]} 
or has no nontrivial superperfect normal subgroup \cite{[Dy87]}. 
Whitehead's question has interesting connections with knot theory.
(For instance, the exterior of a ribbon $n$-knot or of a ribbon concordance 
between classical knots is homotopy equivalent to such a 2-complex.
The asphericity of such ribbon exteriors has been raised in 
\cite{[Co83]} and \cite{[Go81]}.) 

If the answer to Whitehead's question is YES, 
then a high dimensional knot group 
has geometric dimension at most 2 if and only if it has deficiency 1 
(in which case it is a 2-knot group).
For let $G$ be a group of weight 1 and with $G/G'\cong\mathbb{Z}$.
If $C(P)$ is the 2-complex corresponding to a presentation of deficiency 1 
then the 2-complex obtained by adjoining a 2-cell to $C(P)$ along a loop 
representing a weight element for $G$ is 1-connected and has 
Euler characteristic 1, and so is contractible.          
The converse follows from Theorem 2.8.                          
On the other hand a positive answer in general implies
that there is a group $G$ such that $c.d.G=2$ and $g.d.G=3$ 
\cite{[BB97]}.

If the answer is NO then either there is a finite nonaspherical 2-complex $X$ 
such that $X\cup_f D^2$ is contractible for some $f:S^1\to X$ 
or there is an infinite ascending chain of nonaspherical 2-complexes 
whose union is contractible \cite{[Ho83]}.
In the finite case $\chi(X)=0$ and so $\pi=\pi_1(X)$ has deficiency 1; 
moreover, 
$\pi$ has weight 1 since it is normally generated by the conjugacy class 
represented by $f$.
Such groups are 2-knot groups. 
Since $X$ is not aspherical $\beta^{(2)}_1(\pi)\not=0$, by Theorem 2.4,
and so $\pi'$ cannot be finitely generated, by Lemma 2.1.
    
A group is called {\it knot-like} if it has abelianization 
$\mathbb{Z}$ and deficiency 1.
If the commutator subgroup of a classical knot group 
is finitely generated then it is free.                                                                                               
Using the result of Corollary 2.5.1 above and the fact that the Novikov
completions of $\mathbb{Z}[G]$ with respect to epimorphisms 
from $G$ onto $\mathbb{Z}$ are weakly finite, 
Kochloukova has shown that this holds more generally for all
knot-like groups \cite{[Ko06]}.
(See Corollary 4.3.1 above.) 
This answers an old question of Rapaport, who           
established this in the 2-generator, 1-relator case \cite{[Rp60]}.

In particular, if the group of a fibred 2-knot has a presentation of 
deficiency 1 then its commutator subgroup is free. 
Any 2-knot with such a group is $s$-concordant to a fibred homotopy ribbon knot.
(See \S6 of Chapter 17.) 
As $S^2\times S^1$ is never a factor of the closed fibre of a nontrivial 
twist spin $\tau_r K$ \cite{[Pl84]}, 
it follows that if $r>1$ and $K$ is nontrivial
then $\mathrm{def}(\pi\tau_r K)\leq 0$ and $\tau_r K$ is not a homotopy ribbon 
2-knot.

If a knot group has a 2-generator 1-relator Wirtinger presentation 
it is an HNN extension with free base and associated subgroups \cite{[Yo88]}.
This paper also gives an example $\pi$ with $g.d.\pi=2$
and a deficiency 1 Wirtinger presentation which also has a 
2-generator 1-relator presentation but which is not such an HNN extension
(and so has no 2-generator 1-relator {\it Wirtinger} presentation).


\begin{lemma}
If $G$ is a group with $\mathrm{def}(G)=1$ 
and $e(G)=2$ then $G\cong\mathbb{Z}$.
\end{lemma}      

\begin{proof} 
The group $G$ has an infinite cyclic subgroup $A$ of finite index,
since $e(G)=2$.
Let $C$ be the finite 2-complex corresponding to a presentation of deficiency 1 for $G$,
and let $D$ be the covering space corresponding to $A$.
Then $D$ is a finite 2-complex with $\pi_1(D)=A\cong\mathbb{Z}$ and 
$\chi(D)=[\pi:A]\chi(C)=0$.
Since $H_2(D;\mathbb{Z}[A])=H_2(\widetilde D;\mathbb{Z})$ is a submodule 
of a free $\mathbb{Z}[A]$-module and is of rank $\chi(D)=0$ it is 0.
Hence $\widetilde D$ is contractible, and so $G$ must be torsion-free and 
hence abelian.
\end{proof}

This lemma is also a consequence of Theorem 2.5.
It follows immediately that $\mathrm{def}(\pi\tau_23_1)=0$, 
since $\pi\tau_23_1\cong (Z/3Z)\rtimes_{-1}\mathbb{Z}$.
Moreover, if $K$ is a classical knot such that $\pi'$ finitely generated but
nontrivial then $H^1(\pi;\mathbb{Z}[\pi])=0$, and so $X(K)$ is aspherical, 
by Poincar\'e duality.

\begin{theorem}
Let $K$ be a $2$-knot.
Then $\pi=\pi{K}\cong\mathbb{Z}$ if and only if $\mathrm{def}(\pi)=1$ 
and $\pi_2(M(K))=0$.
\end{theorem}

\begin{proof} 
The conditions are necessary, by Theorem 11.1.
If they hold then $\beta^{(2)}_j(M)=\beta^{(2)}_j(\pi)$ 
for $j\leq2$ \cite[Theorem 6.54]{[Lu]}, 
and so $0=\chi(M)=\beta^{(2)}_2(\pi)-2\beta^{(2)}_1(\pi)$.
Now $\beta^{(2)}_1(\pi)-\beta^{(2)}_2(\pi)\geq \mathrm{def}(\pi)-1=0$, 
by Corollary 2.4.1.
Therefore $\beta^{(2)}_1(\pi)=\beta^{(2)}_2(\pi)=0$ and so $g.d.\pi\leq2$, 
by the same Corollary.
In particular, the manifold $M$ is not aspherical.
Hence $H^1(\pi;\mathbb{Z}[\pi])\cong H_3(M;\mathbb{Z}[\pi])\not=0$.
Since $\pi$ is torsion-free it is indecomposable as a free product 
\cite{[Kl93]}.
Therefore $e(\pi)=2$ and so $\pi\cong\mathbb{Z}$, by Lemma 14.5.
\end{proof}

In fact $K$ must be trivial (\cite{[FQ]} - see Corollary 17.1.1).
A simpler argument is used in \cite{[H1]} to show that 
if $\mathrm{def}(\pi)=1$ then 
$\pi_2(M)$ maps onto $H_2(M;\Lambda)$, which is nonzero if $\pi'\not=\pi''$.

\section{Asphericity}

The outstanding property of the exterior of a classical knot 
is that it is aspherical.
Swarup extended the classical Dehn's lemma criterion for unknotting to show 
that if $K$ is an $n$-knot such that the natural inclusion of $S^n$ 
(as a factor of $\partial X(K)$) into $X(K)$ is null homotopic then 
$X(K)\simeq S^1$, provided $\pi K$ is accessible \cite{[Sw75]}.
Since it is now known that finitely presentable groups are accessible 
\cite{[DD]},
it follows that the exterior of a higher dimensional knot is aspherical if 
and only if the knot is trivial.
Nevertheless, we shall see that the closed 4-manifolds $M(K)$ 
obtained by surgery on 2-knots are often aspherical.

\begin{theorem}
Let $K$ be a $2$-knot.
Then $M(K)$ is aspherical if and only if $\pi=\pi{K}$ is a $PD_4$-group,
which must then be orientable.
\end{theorem}

\begin{proof} 
The condition is clearly necessary.
Suppose that it holds.
Let $M^+$ be the covering space associated to $\pi^+=\mathrm{Ker}(w_1(\pi))$.
Then $[\pi:\pi^+]\leq2$, so $\pi'<\pi^+$.
Since $\pi/\pi'\cong Z$ and $t-1$ acts invertibly on $H_1(\pi';\mathbb{Z})$ 
it follows that $\beta_1(\pi^+)=1$.
Hence $\beta_2(M^+)=0$, since $M^+$ is orientable and $\chi(M^+)=0$.
Hence $\beta_2(\pi^+)$ is also 0, so $\chi(\pi^+)=0$,
by Poincar\'e duality for $\pi^+$.
Therefore $\chi(\pi)=0$ and so $M$ must be aspherical,
by Corollary 3.5.1.
\end{proof}

We may use this theorem to give more examples of high dimensional knot groups 
which are not 2-knot groups.
Let $A\in GL(3,\mathbb{Z})$ be such that ${\det(A)=-1}$, 
$\det(A-I)=\pm1$ and $\det(A+I)=\pm1$.
The characteristic polynomial of $A$ must be either $f_1(X)=X^3-X^2-2X+1$, 
$f_2(X)=X^3-X^2+1$, $f_3(X)=X^3f_1(X^{-1})$ or $f_4(X)=X^3f_2(X^{-1})$.
(There are only two conjugacy classes of such matrices, up to inversion,
for it may be shown that the rings $\mathbb{Z}[X]/(f_i(X))$ are 
principal ideal domains.)
The group $\mathbb{Z}^3\rtimes_A\mathbb{Z}$ satifies the Kervaire conditions, 
and is a $PD_4$-group.
However it cannot be a 2-knot group, since it is nonorientable.
(Such matrices have been used to construct fake $RP^4$s 
\cite{[CS76']}.)
                                                                                  
Is every (torsion-free) 2-knot group $\pi$ with $H^s(\pi;\mathbb{Z}[\pi])=0$ 
for $s\leq2$ a $PD_4$-group?
Is every 3-knot group which is also a $PD_4^+$-group a 2-knot group?
(Note that by Theorem 3.6 such a group cannot have deficiency 1.)

We show next that knots with such groups cannot be nontrivial satellites.

\begin{theorem}
Let $K=\Sigma(K_2;K_1,\gamma)$ be a satellite $2$-knot.
If $\pi=\pi{K}$ is a $PD_4$-group then $K=K_1$ or $K_2$.
In particular, $K$ is irreducible.
\end{theorem}

\begin{proof} 
Let $q$ be the order of $\gamma$ in $\pi K_1$.
Then $\pi\cong \pi K_1*_CB$, where  
$B=\pi K_2/\langle\langle w^q\rangle\rangle$, and $C$ is cyclic.
Since $\pi$ is torsion-free $q=0$ or 1.
Suppose that $K\not= K_1$. 
Then $q=0$, so $C\cong\mathbb{Z}$, while $B\not=C$.
If $\pi K_1\not=C$ then $\pi K_1$ and $B$ have infinite index in $\pi$,
and so $c.d.\pi K_1\leq 3$ and $c.d.B\leq 3$,
by Strebel's Theorem.
A Mayer-Vietoris argument then gives $4=c.d.\pi\leq3$,
which is impossible.
Therefore $K_1$ is trivial and so $K=K_2$.
\end{proof}
    
\section{Links}

A {\it $\mu$-component $n$-link} is a locally flat embedding $L:\mu S^n \to S^{n+2} $.
The {\it exterior} of $L$ is $X(L)=S^{n+2}\setminus{int}N(L)$, 
where $N(L)\cong \mu{S^n}\times D^2 $
is a regular neighbourhood of the image of $L$, and the {\it group} of $L$ is $\pi L=\pi_1 (X(L))$.
Let $M(L)=X(L)\cup\mu{D^{n+1}}\times S^1 $ be the closed manifold obtained by 
surgery on $L$.

An $n$-link $L$ is {\it trivial} if it bounds a collection of $\mu$ disjoint 
locally flat 2-discs in $S^n$.
It is {\it split} if it is isotopic to one which is the union of nonempty 
sublinks $L_1 $ and $L_2 $ whose images lie in disjoint discs in $S^{n+2} $,
in which case we write $L=L_1 \amalg L_2$, and it is a {\it boundary} link if 
it bounds a collection of $\mu$ disjoint hypersurfaces in $S^{n+2}$.
Clearly a trivial link is split, and a split link is a boundary link; 
neither implication can be reversed if $\mu>1$.
Knots are boundary links, and many arguments about knots that depend 
on Seifert hypersurfaces extend readily to boundary links.
The definitions of slice and ribbon knots and $s$-concordance extend naturally
to links.

A 1-link is trivial if and only if its group is free, and is
split if and only if its group is a nontrivial free product,
by the Loop Theorem and Sphere Theorem, respectively. 
(See \cite[Chapter 1]{[H3]}.)
Guti\'errez has shown that if $n\geq 4$ an $n$-link $L$ is trivial 
if and only if $\pi L$ is freely generated by meridians 
and the homotopy groups $\pi_j (X(L))$ are all 0,
for $2\leq j\leq (n+1)/2$ \cite{[Gu72]}. 
His argument applies also when $n=3$.
While the fundamental group condition is necessary when $n=2$,
we cannot yet use surgery to show that it is a complete criterion for
triviality of 2-links with more than one component.
We shall settle for a weaker result.
                             
\begin{theorem}
Let $M$ be a closed $4$-manifold with $\pi_1 (M)$ 
free of rank $r$ and $\chi(M)=2(1-r)$. 
If $M$ is orientable it is $s$-cobordant to $\sharp^r (S^1 \times S^3 )$, 
while if it is nonorientable it is $s$-cobordant to 
$(S^1\tilde\times S^3 )\sharp (\sharp^{r-1} (S^1 \times S^3))$.
\end{theorem}

\begin{proof}
We may assume without loss of generality that $\pi_1 (M)$ has a free basis 
$\{ x_1,\dots,x_r\}$
such that $x_i$ is an orientation preserving loop for all $i>1$,
and we shall use $c_{M*} $ to identify $\pi_1 (M)$ with $F(r)$.
Let $N=\sharp^r (S^1 \times S^3)$ if $M$ is orientable and let 
$N=(S^1\tilde\times S^3 )\sharp (\sharp^{r-1} (S^1 \times S^3))$ otherwise.
(Note that $w_1 (N)=w_1 (M)$ as homomorphisms from $F(r)$ to $\{\pm1\}$.)
Since $c.d.\pi_1 (M)\leq 2$ and $\chi(M)=2\chi(\pi_1 (M))$ we have
$\pi_2 (M)\cong\overline{ H^2 (F(r);\mathbb{Z}[F(r)])}$, by Theorem 3.12.
Hence $\pi_2 (M)=0$ and so 
$\pi_3 (M)\cong H_3 (\widetilde M;\mathbb{Z})\cong
{E}=\overline{H^1 (F(r);\mathbb{Z}[F(r)])}$,
by the Hurewicz theorem and Poincar\'e duality. 
Similarly, we have $\pi_2 (N)=0$ and $\pi_3 (N)\cong{E}$. 

Let $c_M=g_M h_M $ be the factorization of $c_M$ through $P_3 (M)$,
the third stage of the Postnikov tower for $M$.
Thus $\pi_i (h_M)$ is an isomorphism if $i\leq 3$ and $\pi_j (P_3 (M))=0$ if $j>3$.
As $K(F(r),1)=\vee^r S^1$ each of the fibrations $g_M $ and $g_N$ clearly have cross-sections
and so there is a homotopy equivalence $k:P_3 (M)\to P_3 (N)$ such that 
$g_M=g_N k$.
(See \cite[\S5.2]{[Ba]}.)
We may assume that $k$ is cellular.
Since $P_3 (M)=M\cup\{ cells~of~dimension\geq5\} $ it follows that $kh_M=h_N f$
for some map $f:M\to N$.
Clearly $\pi_i (f)$ is an isomorphism for $i\leq 3$.
Since the universal covers $\widetilde M$ and $\widetilde N$ are 2-connected 
open 4-manifolds the induced map $\tilde f:\widetilde M\to\widetilde N$ 
is an homology isomorphism, and so is a homotopy equivalence.
Hence $f$ is itself a homotopy equivalence. 
As $Wh(F(r))=0$ any such homotopy equivalence is simple.

If $M$ is orientable $[M,G/TOP]\cong\mathbb{Z}$, 
since $H^2(M;Z/2Z)=0$.
As the surgery obstruction in $L_4(F(r))\cong\mathbb{Z}$ 
is given by a signature difference,
it is a bijection, and so the normal invariant of $f$ is trivial.
Hence there is a normal cobordism $F:P\to N\times I$ with $F|\partial_- P=f$ 
and $F|\partial_+ P=id_N$. 
There is another normal cobordism $F':P'\to N\times I$ from $id_N$ to itself 
with surgery obstruction $\sigma_5 (P',F')=-\sigma_5 (P,F)$ in $L_5 (F(r))$, 
by Theorem 6.7 and Lemma 6.9.
The union of these two normal cobordisms along $\partial_+ P=\partial_- P'$
is a normal cobordism from $f$ to $id_N $ with surgery obstruction 0, and so
we may obtain an $s$-cobordism $W$ by 5-dimensional surgery (rel $\partial$). 

A similar argument applies in the nonorientable case.
The surgery obstruction is then a bijection from $[N;G/TOP]$ to 
$L_4 (F(r),-)=Z/2Z$,
so $f$ is normally cobordant to $id_N$, 
while $L_5 (\mathbb{Z},-)=0$, 
so $L_5(F(r),-)\cong{L_5(F(r-1))}$
and the argument of \cite{[FQ]} still applies. 
\end{proof}

\begin{cor}
Let $L$ be a $\mu$-component $2$-link
such that $\pi L$ is freely generated by $\mu$ meridians. 
Then $L$ is $s$-concordant to the trivial $\mu$-component link.
\end{cor}

\begin{proof}
Since $M(L)$ is orientable, $\chi(M(L))=2(1-\mu)$
and $\pi_1 (M(L))\cong\pi L=F(\mu)$, there is an $s$-cobordism $W$
with $\partial W=M(L)\cup M(\mu)$, by Theorem 14.9.
Moreover it is clear from the proof of that theorem that we may assume that the
elements of the meridianal basis for $\pi L$ are freely homotopic to loops
representing the standard basis for $\pi_1 (M(\mu))$.
We may realise such homotopies by $\mu$ disjoint embeddings of annuli 
running from meridians for $L$ to such standard loops in $M(\mu)$.
Surgery on these annuli (i.e., replacing $D^3\times S^1\times [0,1]$
by $S^2\times D^2\times [0,1]$) then gives an $s$-concordance 
from $L$ to the trivial $\mu$-component link.   
\end{proof}

A similar strategy may be used to give an alternative proof 
of the higher dimensional unlinking theorem of \cite{[Gu72]} 
which applies uniformly for $n\geq 3$.
The hypothesis that $\pi L$ be freely generated by meridians cannot be dropped entirely \cite{[Po71]}.
On the other hand, if $L$ is a 2-link whose longitudes are all null homotopic
then the pair $(X(L),\partial X(L))$ is homotopy equivalent to
the pair
$(\sharp^\mu S^1\times D^3,\partial(\sharp^\mu S^1\times D^3))$ 
\cite{[Sw77]}, 
and hence the Corollary applies.

There is as yet no satisfactory splitting criterion for higher-dimensional links.
However we can give a stable version for 2-links.

\begin{theorem}
Let $M$ be a closed $4$-manifold such that $\pi=\pi_1 (M)$
is isomorphic to a nontrivial free product $G*H$.
Then $M$ is stably homeomorphic to a connected sum $M_G\sharp M_H $ 
with $\pi_1 (M_G)\cong G$
and $\pi_1 (M_H )\cong H$.
\end{theorem}

\begin{proof}
Let $K=K_G\cup [-1,1]\cup K_H /(*_G \sim -1,+1\sim *_H )$, where
$K_G $ and $K_H $ are $K(G,1)$- and $K(H,1)$-spaces with basepoints $*_G $ and $*_H $ 
(respectively).
Then $K$ is a $K(\pi,1)$-space and so there is a map $f:M\to K$ which induces an isomorphism of
fundamental groups. 
We may assume that $f$ is transverse to $0\in [-1,1]$,
so $V=f^{-1} (0)$ is a submanifold of $M$ 
with a product neighbourhood $V\times [-\epsilon,\epsilon]$.
We may also assume that $V$ is connected, by the arc-chasing argument 
of Stallings' proof of Kneser's conjecture. 
(See \cite[page 67]{[Hm]}.)
Let $j:V\to M$ be the inclusion.
Since $fj$ is a constant map and $\pi_1 (f)$ is an isomorphism $\pi_1 (j)$ is the trivial
homomorphism, and so $j^* w_1 (M)=0$.
Hence $V$ is orientable and so there is a framed link $L\subset V$ 
such that surgery on $L$ in $V$ gives $S^3 $ \cite{[Li62]}.
The framings of the components of $L$ in $V$ extend to framings in $M$.
Let $W=M\times [0,1]\cup_{L\times D^2\times [-\epsilon,\epsilon]\times\{1\}} 
(\mu D^2 \times D^2\times [-\epsilon,\epsilon])$,
where $\mu$ is the number of components of $L$.
Note that if $w_2 (M)=0$ then we may choose the framed link $L$ so that $w_2 (W)=0$ also \cite{[Kp79]}.
Then $\partial W=M\cup\widehat M$, where $\widehat M$ is the result of surgery 
on $L$ in $M$.
The map $f$ extends to a map $F:W\to K$ such that $\pi_1 (F|_{\widehat M})$ 
is an isomorphism and $(F|_{\widehat M})^{-1} (0)\cong S^3 $.
Hence $\widehat M$ is a connected sum as in the statement.
Since the components of $L$ are null-homotopic in $M$ they may be isotoped into disjoint discs, 
and so $\widehat M\cong M\sharp(\sharp^\mu S^2\times S^2 )$.
This proves the theorem. 
\end{proof}

Note that if $V$ is a homotopy 3-sphere then $M$ is a connected sum, 
for $V\times\mathbb{R}$ is then homeomorphic to $S^3 \times\mathbb{R}$, by 1-connected surgery.

\begin{theorem}
Let $L$ be a $\mu$-component $2$-link with sublinks 
$L_1$ and $L_2=L\setminus{L_1}$ such that there is 
an isomorphism from $\pi L$ to $\pi L_1 *\pi L_2 $ which is compatible 
with the homomorphisms determined by the inclusions of $X(L)$ 
into $X(L_1 )$ and $X(L_2 )$.
Then $X(L)$ is stably homeomorphic to $X(L_1 \amalg L_2 )$.
\end{theorem} 

\begin{proof}
By Theorem 14.10, $M(L)\sharp(\sharp^a S^2 \times S^2)\cong N\sharp P$, 
where $\pi_1 (N)\cong \pi L_1 $ and $\pi_1 (P)\cong \pi L_2 $. 
On undoing the surgeries on the components of $L_1 $ and $L_2 $, respectively, 
we see that $M(L_2)\sharp(\sharp^a S^2 \times S^2)\cong N\sharp\bar P$,
and $M(L_1)\sharp(\sharp^a S^2 \times S^2)\cong \bar N\sharp P$,
where $\bar N$ and $\bar P$ are simply connected. 
Since undoing the surgeries on all the components of $L$ gives
$\sharp^a S^2 \times S^2\cong \bar N\sharp\bar P$, 
$\bar N$ and $\bar P$ are each connected sums of copies of $S^2\times{S^2}$, 
so $N$ and $P$ are stably homeomorphic 
to $M(L_1 )$ and $M(L_2 )$, respectively. The result now follows easily. 
\end{proof}

Similar arguments may be used to show that, firstly,
if $L$ is a 2-link such that $c.d.\pi L\leq 2$ and there is an isomorphism
$\theta:\pi L\to \pi L_1*\pi L_2 $ which is compatible with the natural maps 
to the factors then there is a map 
$f_o :M(L)_o =M(L)\setminus{int}D^4 \to M(L_1 )\sharp M(L_2)$
such that $\pi_1 (f_o)=\theta$ and $\pi_2 (f_o )$ is an isomorphism;
and secondly, if moreover $f_o $ extends to a homotopy equivalence
$f:M(L)\to M(L_1 )\sharp M(L_2 )$ and the factors of $\pi L$ are either 
classical link groups or are square root closed accessible then $L$ is 
$s$-concordant to the split link $L_1 \amalg L_2 $.
(The surgery arguments rely on \cite{[AFR97]} and \cite{[Ca73]}, respectively.)
However we do not know how to bridge the gap between the algebraic hypothesis
and obtaining a homotopy equivalence.

\section{Link groups}

If $\pi$ is the group of a $\mu$-component $n$-link $L$ then 

{\sl
\begin{enumerate}
\item $\pi$ is finitely presentable;

\item $\pi$ is of weight $\mu$;

\item $H_1(\pi;\mathbb{Z})=\pi/\pi'\cong\mathbb{Z}^\mu$; and

\item (if $n>1$) $H_2(\pi;\mathbb{Z})=0$.
\end{enumerate}}
Conversely, any group satisfying these conditions is the group of an $n$-link,
for every $n\geq3$ \cite{[Ke65']}. 
(Note that $q(\pi)\geq2(1-\mu)$, with equality if and only if $\pi$ is the group of a 2-link.)
If (4) is replaced by the stronger condition that 
$\mathrm{def}(\pi)=\mu$         
(and $\pi$ has a deficiency $\mu$ Wirtinger presentation) 
then $\pi$ is the group of a (ribbon) 2-link which is a sublink of a (ribbon) 
link whose group is a free group.
(See \cite[Chapter 1]{[H3]}.)
The group of a classical link satisfies (4) if and only if the link 
splits completely as a union of knots in disjoint balls.
If subcomplexes of aspherical 2-complexes are aspherical then 
a higher-dimensional link group group has geometric dimension at most 2 
if and only if it has deficiency $\mu$ (in which case it is a 2-link group).

A link $L$ is a boundary link if and only if there is an epimorphism from
$\pi(L)$ to the free group $F(\mu)$ which carries a set of meridians to a free basis.
If the latter condition is dropped $L$ is said to be an homology boundary link.
Although sublinks of boundary links are clearly boundary links,
the corresponding result is not true for homology boundary links.
It is an attractive conjecture that every even-dimensional link is a slice link.    
This has been verified under additional hypotheses on the link group.
For a 2-link $L$ it suffices that there be a homomorphism $\phi:\pi L\to G$ 
where $G$ is a high-dimensional link group such that 
$H_3(G;\mathbb{F}_2)=H_4(G;\mathbb{Z})=0$ and
where the normal closure of the image of $\phi$ is $G$ \cite{[Co84]}.
In particular, sublinks of homology boundary 2-links are slice links.

A choice of (based) meridians for the components of a link $L$ determines a homomorphism 
$f:F(\mu)\to \pi L$ which induces an isomorphism on abelianization.
If $L$ is a higher dimensional link $H_2(\pi L;\mathbb{Z})=H_2(F(\mu);\mathbb{Z})=0$ 
and hence $f$ induces isomorphisms on all the nilpotent quotients 
$F(\mu)/F(\mu)_{[n]}\cong\pi L/(\pi L)_{[n]}$, and a monomorphism 
$F(\mu)\to\pi L/(\pi L)_{[\omega]}= \pi L/\cap_{n\geq1}(\pi L)_{[n]}$ 
\cite{[St65]}. 
(In particular, if $\mu\geq2$ then $\pi L$ contains a nonabelian free subgroup.)
The latter map is an isomorphism if and only if $L$ is a homology boundary link.
In that case the homology groups of the covering space $X(L)^\omega$ corresponding
to $\pi L/(\pi L)_{[\omega]}$ are modules over 
$\mathbb{Z}[\pi L/(\pi L)_{[\omega]}]\cong\mathbb{Z}[F(\mu)]$,
which is a coherent ring of global dimension 2.
Poincar\'e duality and the UCSS then give rise to an isomorphism
$e^2e^2 (\pi L/(\pi L)_{[\omega]}) \cong
\overline{e^2(\pi L/(\pi L)_{[\omega]})}$,
where $e^i(M)=Ext^i_{\mathbb{Z}[F(\mu)]}(M,\mathbb{Z}[F(\mu)])$, 
which is the analogue of the Farber-Levine pairing for 2-knots.                  
The argument of \cite{[HK78']} may be adapted to show that 
every finitely generated 
abelian group is the centre of the group of some $\mu$-component
boundary $n$-link, for any $\mu\geq1$ and $n\geq3$.
However the centre of the group of a 2-link with more than one component 
must be finite.
(All known examples have trivial centre.)

\begin{theorem}
Let $L$ be a $\mu$-component $2$-link.
If $\mu>1$ then 
\begin{enumerate}
\item $\pi{L}$ has no infinite amenable normal subgroup;
             
\item $\pi{L}$ is not an ascending HNN extension over a finitely 
generated base.
\end{enumerate}
\end{theorem}

\begin{proof}
Since $\chi(M(L))=2(1-\mu)$ the $L^2$-Euler characteristic formula gives
$\beta_1^{(2)}(\pi{L})\geq\mu-1$.
Therefore $\beta_1^{(2)}(\pi{L})\not=0$ if $\mu>1$, and so 
the result follows from Lemma 2.1 and Corollary 2.3.1.
\end{proof}

In particular, the exterior of a 2-link with more than one component never fibres over $S^1$. 
(This is true of all higher dimensional links: 
see \cite[Theorem 5.12]{[H3]}.)
Moreover a 2-link group has finite centre and is never amenable.
In contrast, we shall see that there are many 2-knot groups which have infinite centre or are solvable.
                                                                          
The exterior of a classical link is aspherical if and only the link is unsplittable,
while the exterior of a higher dimensional link with more 
than one component is never aspherical \cite{[Ec76]}.
Is $M(L)$ ever aspherical?                                                   

\section{Homology spheres}

A closed connected $n$-manifold $M$ is an {\it homology $n$-sphere\/} 
if $H_q(M;\mathbb{Z})=0$ for $0<q<n$. 
In particular, it is orientable and so $H_n(M;\mathbb{Z})\cong\mathbb{Z}$.
If $\pi$ is the group of an homology $n$-sphere then 
{\sl
\begin{enumerate}
\item $\pi$ is finitely presentable;

\item $\pi$ is perfect, i.e., $\pi=\pi'$; and

\item $H_2(\pi;\mathbb{Z})=0$.
\end{enumerate}}
A group satisfying the latter two conditions is said to be {\it superperfect}.
Every finitely presentable superperfect group
is the group of an homology $n$-sphere, for every $n\geq5$ 
\cite{[Ke69]},
but in low dimensions more stringent conditions hold.
As any closed 3-manifold has a handlebody structure with one 0-handle and equal
numbers of 1- and 2-handles, homology 3-sphere groups have deficiency 0.
Every perfect group with a presentation of deficiency 0 is superperfect, 
and is an homology 4-sphere group \cite{[Ke69]}.
However none of the implications ``$G$ is an homology 3-sphere group"
$\Rightarrow$ ``$G$ is finitely presentable, perfect and $\mathrm{def}(G)=0$" 
$\Rightarrow$ ``$G$ is an homology 4-sphere group"
$\Rightarrow$ ``$G$ is finitely presentable and superperfect" can be reversed,
as we shall now show.

Although the finite groups $SL(2,\mathbb{F}_p)$ are perfect and have
deficiency 0 for each prime $p\geq5$ \cite{CR80}, 
the binary icosahedral group $I^*=SL(2,\mathbb{F}_5)$ 
is the only nontrivial finite perfect group with cohomological period 4, 
and thus is the only finite homology 3-sphere group.

Let $G=\langle x,s\mid x^3=1, sxs^{-1}=x^{-1}\rangle$ be the group of
$\tau_23_1$ and let 
\[H=\langle a,b,c,d\mid bab^{-1}=a^2, cbc^{-1}=b^2, dcd^{-1}=c^2,
ada^{-1}=d^2\rangle\] 
be the Higman group \cite{[Hg51]}.
Then $H$ is perfect and $\mathrm{def}(H)=0$, 
so there is an homology 4-sphere $\Sigma$ with group $H$. 
Surgery on a loop representing $sa^{-1}$ in 
$\Sigma\sharp M(\tau_23_1)$ gives an homology
4-sphere with group $\pi=(G*H)/\langle\langle sa^{-1}\rangle\rangle$.
Then $\pi$ is the semidirect product $\rho\rtimes H$,
where $\rho=\langle\langle G'\rangle\rangle_\pi$ is the normal closure of the
image of $G'$ in $\pi$.
The obvious presentation for this group has deficiency -1.
We shall show that this is best possible.

Let $\Gamma=\mathbb{Z}[H]$.
Since $H$ has cohomological dimension 2 \cite{[DV73']}
the augmentation ideal $I=\mathrm{Ker}(\varepsilon:\Gamma\to\mathbb{Z})$
has a short free resolution
\begin{equation*}
C_*:0\to\Gamma^4\to\Gamma^4\to I\to 0.
\end{equation*}
Let $B=H_1(\pi;\Gamma)\cong\rho/\rho'$.
Then $B\cong\Gamma/\Gamma(3,a+1)$ as a left $\Gamma$-module 
and there is an exact sequence 
\begin{equation*}
0\to B\to A\to I\to0,
\end{equation*}
in which $A=H_1(\pi,1;\Gamma)$ is a relative homology group \cite{[Cr61]}.
Since $B\cong\Gamma\otimes_\Lambda(\Lambda/\Lambda(3,a+1))$,
where $\Lambda=\mathbb{Z}[a,a^{-1}]$, there is a free resolution 
\begin{equation*}
\begin{CD}
0\to\Gamma@>(3,a+1)>>\Gamma^2@>\left(\smallmatrix a+1\\
-3\endsmallmatrix\right)
>>\Gamma\to B\to 0.
\end{CD}
\end{equation*}

Suppose that $\pi$ has deficiency 0.
Evaluating the Jacobian matrix associated to an optimal presentation for $\pi$ 
via the natural epimorphism from $\mathbb{Z}[\pi]$ to $\Gamma$ gives
a presentation matrix for $A$ as a module \cite{[Cr61],[Fo62]}.
Thus there is an exact sequence 
\begin{equation*}
D_*:\dots\to\Gamma^n\to\Gamma^n\to A\to0.
\end{equation*}
A mapping cone construction leads to an exact sequence  
of the form 
\begin{equation*}
D_1\to C_1\oplus D_0\to B\oplus C_0\to0
\end{equation*}
and hence to a presentation of deficiency $0$ for $B$ of the form
\begin{equation*}
D_1\oplus C_0\to C_1\oplus D_0\to B.
\end{equation*}
Hence there is a free resolution
\begin{equation*}
0\to L\to\Gamma^p\to\Gamma^p\to B\to0.
\end{equation*}
Schanuel's Lemma gives an isomorphism $\Gamma^{1+p+1}\cong L\oplus\Gamma^{p+2}$,
on comparing these two resolutions of $B$.
Since $\Gamma$ is weakly finite the endomorphism of $\Gamma^{p+2}$ given by
projection onto the second summand is an automorphism.
Hence $L=0$ and so $B$ has a short free resolution.
In particular, $Tor_2^\Gamma(R,B)=0$ for any right $\Gamma$-module $R$.
But it is easily verified that if $\overline{B}\cong \Gamma/(3,a+1)\Gamma$
is the conjugate right $\Gamma$-module then
$Tor_2^\Gamma(\overline{B},B)\not=0$.
Thus our assumption was wrong, and $\mathrm{def}(\pi)=-1<0$.
 
Our example of an homology 4-sphere group with negative deficiency is 
``very infinite" in the sense that the Higman 
group $H$ has no finite quotients, and therefore no finite-dimensional 
representations over any field \cite{[Hg51]}.
Livingston has constructed examples with arbitrarily large negative deficiency,
which are extensions of $I^*$ by groups which map onto $\mathbb{Z}$.
His argument uses only homological algebra for infinite cyclic covers 
\cite{[Li05]}.

Let $I^*=SL(2,\mathbb{F}_5)$ act diagonally on $(\mathbb{F}_5^2)^k$,
and let $G_k$ be the universal central extension of the perfect group 
$(\mathbb{F}_5^2)^k\rtimes{I^*}$.
In \cite{[HW85]} it is shown that for $k$ large, 
this superperfect group is not the group of an homology 4-sphere.
In particular, it has negative deficiency.
It seems unlikely that deficiency 0 is a necessary condition
for a finite perfect group to be an homology 4-sphere group.

Kervaire's criteria may be extended further to the groups of links 
in homology spheres.
Unfortunately, the condition $\chi(M)=0$ is central to most of our arguments,
and is satisfied only when the link has one component.

We may also use the Higman group $H$ to construct fibred knots with group of 
cohomological dimension $d$, for any finite $d\geq1$.
If $\Sigma$ is an homology 4-sphere with $\pi_1(\Sigma)\cong{H}$
then surgery on a loop in $\Sigma\times{S^1}$ representing 
$(a,1)\in{H\times\mathbb{Z}}$ gives a homotopy 5-sphere, 
and so $H\times\mathbb{Z}$ is the group of a fibred 3-knot.
Every superperfect group is the group of an homology 5-sphere, 
and a similar construction shows that $H^k\times\mathbb{Z}$ 
is the group of a fibred 4-knot, for all $k\geq0$.
Similarly, $\pi=H^{k}\times\pi3_1$ is a high-dimensional knot group
with $\pi'$ finitely presentable and $c.d.\pi=2k+2$, for all $k\geq0$.

On the other hand,
if $K$ is a 2-knot with group $\pi=\pi{K}$ such that $\pi'$ is
finitely generated then $M(K)'$ is a $PD_3$-space, by Theorem 4.5.
Hence $\pi'$ has a factorization as a free product of $PD_3$-groups
and indecomposable virtually free groups,
by the theorems of Turaev and Crisp.
In particular, $v.c.d.\pi'=0$, 1 or 3, and so $v.c.d.\pi=1$, 2 or 4.
Thus $H^k\times\mathbb{Z}$ is not a 2-knot group, if $k\geq1$.
 
These observations suggest several questions:
\begin{enumerate}
\item are there any 2-knot groups $\pi$ with $c.d.\pi=3$?
\item what are the groups of fibred $n$-knots?
\item in particular, is $H^k\times\pi3_1$ realized by a fibred 3-knot, 
if $k\geq2$?
\end{enumerate}

%% file: m5-15.tex
\chapter{Restrained normal subgroups}

It is plausible that if a 2-knot group $\pi=\pi K$ has an
infinite restrained normal subgroup $N$ then either $\pi'$ is finite or 
$\pi\cong\Phi=Z*_2$ (the group of Fox's Example 10)
or $M(K)$ is aspherical and $\sqrt\pi\not=1$ or $N\cong\mathbb{Z}$ 
and $\pi/N$ has infinitely many ends.
In this chapter we shall give some evidence in this direction.  
We begin by determining the  2-knot groups $\pi$ with $\pi'$ finite.
In \S2 we show that if $\pi$ has one end and $\beta_1^{(2)}(\pi)=0$
then it has no non-trivial finite normal subgroup.
We then give several characterizations of $\Phi$,
which plays a somewhat exceptional role.
In \S4 we establish the above tetrachotomy when $N$ is nilpotent
and $h(N)\geq2$ or $N$ is abelian of rank 1.
In \S5 we assume that $N$ is almost coherent,
locally virtually indicable and restrained.
In \S6 we show that virtually torsion-free solvable 2-knot groups
with $\pi'$ infinite and other than $\Phi$ are virtually poly-$Z$,
and give a version of the Tits alternative for 2-knot groups.
In \S7 we complete Yoshikawa's study of the
2-knot groups which are HNN extensions over abelian bases.
It is not known whether any 2-knot groups has an infinite
locally finite normal subgroup.
We consider this issue briefly in \S8.

\section{Finite commutator subgroup}

It is a well known consequence of the asphericity 
of the exteriors of classical knots 
that classical knot groups are torsion-free.
The first examples of higher dimensional knots whose groups have nontrivial 
torsion were given by Mazur \cite{[Ma62]} and Fox \cite{[Fo62]}. 
These examples are 2-knots whose groups have finite commutator subgroup. 
We shall show that if $\pi$ is such a group then $\pi'$ must be a CK group, 
and that the images of meridianal automorphisms in $Out(\pi')$ 
are conjugate, up to inversion.
Hence $\pi$ is determined by $\pi'$.
In each case there is just one 2-knot group with given finite
commutator subgroup.
Many of these groups can be realized by twist spinning classical knots.
Zeeman introduced twist spinning in order to study Mazur's example;
Fox used hyperplane cross sections, but his examples
(with $\pi'\cong Z/3Z$) are in fact twist spins \cite{[Kn83']}.

\begin{lemma} 
An automorphism of $Q(8)$ is meridianal if and 
only if it is conjugate to $\sigma$.
\end{lemma}
		      
\begin{proof} Since $Q(8)$ is solvable an automorphism is meridianal 
if and only if the induced automorphism of $Q(8)/Q(8)'$ is meridianal. 
It is easily verified that all such elements of 
$Aut(Q(8))\cong(Z/2Z)^2\rtimes SL(2,\mathbb{F}_2)$ are conjugate to $\sigma$.
\end{proof}

\begin{lemma} 
All nontrivial automorphisms of $I^*$ are meridianal.
Moreover each automorphism is conjugate to its inverse.                          
The nontrivial outer automorphism class of $I^*$ cannot be realised by a 
$2$-knot group.
\end{lemma}
                                                             
\begin{proof} 
The first assertion is immediate,
since $\zeta{I^*}$ is the only nontrivial proper normal subgroup of $I^*$.
Each automorphism is conjugate to its inverse, since $Aut(I^*)\cong S_5$.
Considering the Wang sequence for the map ${:M(K)'\to{M(K)}}$ 
shows that the meridianal automorphism induces the identity on 
$H_3(\pi';\mathbb{Z})$, 
and so the nontrivial outer automorphism class cannot occur,
by Lemma 11.4.
\end{proof}

The elements of order 2 in $A_5\cong Inn(I^*)$ are all conjugate, 
as are the elements of order 3.
There are two conjugacy classes of elements of order 5.
(Thus there are four weight classes in $\mathbb{Z}\times{I^*}$.)

\begin{lemma} 
An automorphism of $T^*_k$ is meridianal if and only if it is conjugate
to $\rho^{3^{k-1}}$ or $\rho^{3^{k-1}}\eta$.
These have the same image in $Out(T^*_k)$.
\end{lemma}

\begin{proof} 
Since $T^*_k$ is solvable an automorphism is meridianal if and only if the 
induced automorphism of
$T^*_k/(T^*_k)'$ is meridianal.
Any such automorphism is conjugate to either $\rho^{2j+1}$ or to $\rho^{2j+1}\eta$ for some
$0\leq j<3^{k-1}$.
(Note that 3 divides $2^{2j}-1$ but does not divide $2^{2j+1}-1$.)
However among them only those with $2j+1=3^{k-1}$ satisfy the 
isometry condition of Theorem 14.3.
\end{proof}

\begin{theorem} 
Let $K$ be a $2$-knot with group $\pi=\pi K$.
If $\pi'$ is finite then $\pi'\cong P\times (Z/nZ)$ where $P=1$, $Q(8)$, 
$I^*$ or $T^*_k$, and $(n,2|P|)=1$, 
and the meridianal automorphism sends $x$ and $y$ in $Q(8)$ to $y$ and $xy$,
is conjugation by a noncentral element on $I^*$,
sends $x$, $y$ and $z$ in $T^*_k$ to $y^{-1}$, $x^{-1}$ and $z^{-1}$, 
and is $-1$ on the cyclic factor.
\end{theorem}

\begin{proof} 
Since $\chi(M(K))=0$ and $\pi$ has two ends $\pi'$ has cohomological period 
dividing 4, by Theorem 11.1, and so is among the groups listed in
\S2 of Chapter 11. 
As the meridianal automorphism of $\pi'$ induces a meridianal automorphism 
on the quotient by any characteristic subgroup, 
we may eliminate immediately the groups $O^*(k)$, $A(m,e)$ and $Z/2mZ$
and their products with $Z/nZ$ since these all have abelianization 
cyclic of even order.
If $k>1$ the subgroup generated by $x$ in $Q(8k)$ is a characteristic subgroup 
of index 2.
Since $Q(2^na)$ is a quotient of $Q(2^na,b,c)$ by a characteristic subgroup 
(of order $bc$) this eliminates this class also.
Thus there remain only the above groups.

Automorphisms of a group $G=H\times J$ such that $(|H|,|J|)=1$ 
correspond to pairs of automorphisms $\phi_H$ and $\phi_J$ of $H$ and $J$, 
respectively,
and $\phi$ is meridianal if and only if $\phi_H$ and $\phi_J$ are.
Multiplication by $s$ induces a meridianal automorphism of
$Z/mZ$ if and only if $(s-1,m)=(s,m)=1$. 
If $Z/mZ$ is a direct factor of $\pi'$ then it
is a direct summand of $\pi'/\pi''=H_1(M(K);\Lambda)$ and so 
$s^2\equiv1$ modulo $(m)$,
by Theorem 14.3. Hence we must have $s\equiv-1$ modulo $(m)$.
The theorem now follows from Lemmas 15.1 -- 15.3.
\end{proof}

Finite cyclic groups are realized by the 2-twist spins of 2-bridge knots,
while $Q(8)$, $T^*_1$ and $I^*$ are realized by $\tau_33_1$, 
$\tau_43_1$ and $\tau_53_1$, respectively.
As the groups of 2-bridge knots have 2-generator 1-relator presentations,
the groups of these twist spins have 2-generator presentations of deficiency 0.
In particular, $\pi\tau_r3_1$ has the presentation
$\langle a,t\mid tat^{-1}=at^2at^{-2},\medspace t^ra=at^r\rangle$.

Let $F$ be the binary polyhedral group with presentation
\[
\langle{y,z}\mid(yz)^2=y^r=z^s\rangle,
\]
where  $(r,s)=(2,2)$, (3,3) or (3,5). 
Then $F\cong{Q(8)}$, $T^*_1$ or $I^*$, respectively, and
the image of $y^r$ is a central element of order 2. 
Let $G$ be the group with presentation 
\[
\langle{u,w}\mid(u^nw^n)^2=u^{rn},~u^r=w^s\rangle,
\]
where  $(n,2rs)=1$, and let $a=u^{2r}=w^{2s}$. 
If $mn\equiv1$ {\it mod} $(2rs)$ then the rules $f(y)=u^n$, $f(z)=w^n$,
$g(u)=y^m$ and $g(w)=z^m$ define homomorphisms $f:F\to{G}$ 
and $g:G\to{F}$ such that $gf=id_F$.
The image of $a$ is central in $G$, of order $n$, 
and generates $G/F\cong{Z/nZ}$.
Hence $G\cong{F}\times{Z/nZ}$.
Since these groups have balanced 2-generator presentations 
the semidirect products $\pi\cong{G}\rtimes\mathbb{Z}$ 
have 3-generator presentations of deficiency $-1$
(one new generator for the meridian $t$ and two new relations 
determining $tut^{-1}$ and $twt^{-1}$).

The groups with $\pi'\cong Q(8)\times(Z/nZ)$  (for $n$ odd) are realized by fibred 2-knots \cite{[Yo80]}, 
but if $n>1$ no such group can be realized by a twist spin. 
(See \S3 of Chapter 16.)
An extension of the twist spin construction may be used to realize
such groups by smooth fibred knots in the standard $S^4$, 
if $n=3$, 5, 11, 13, 19, 21 or 27 \cite{[Kn88],[Tr90]}. 
Is this so in general?  
The other groups are realized by the 2-twist spins of certain pretzel knots 
\cite{[Yo80]}. 

If $\pi'\cong T_k^*\times(Z/nZ)$ (with $(n,6)=1$) 
then $\pi$ has a presentation
\[
\langle s,x,y,z\mid {x^2=(xy)^2=y^2,}
\medspace z^\alpha=1,\medspace zxz^{-1}=y,\medspace zyz^{-1}=xy,
\]
\[
sxs^{-1}=y^{-1},\medspace sys^{-1}=x^{-1},\medspace szs^{-1}=z^{-1}\rangle,
\]
where $\alpha=3^kn$.
This is equivalent to the presentation
\[
\langle s,x,y,z\mid {z^\alpha=1,}\medspace zxz^{-1}=y,\medspace  
zyz^{-1}=xy, \medspace sxs^{-1}=y^{-1}\!,\medspace  szs^{-1}=z^{-1}\rangle.
\]
For conjugating $zxz^{-1}=y$ by $s$ gives $z^{-1}y^{-1}z=sys^{-1}$,
so $sys^{-1}=x^{-1}$,
while conjugating $zyz^{-1}=xy$ by $s$ gives $x=yxy$,
so $x^2=(xy)^2$, and conjugating this by $s$ gives $y^2=(xy)^2$.
On replacing $s$ by $t=xzs$  we obtain the presentation
\[
\langle t,x,y,z\mid z^\alpha=1,\medspace zxz^{-1}=y,\medspace
zyz^{-1}=xy,\medspace txt^{-1}=xy,\medspace tzt^{-1}=yz^{-1}\rangle.
\]
We may use the second and final relations to eliminate the generators 
$x$ and $y$ to obtain a 2-generator presentation of deficiency $-1$.
(When $n=k=1$ we may relate this to the above presentation for $\pi\tau_43_1$ 
by setting $a=zx^2$.)

Are there similar 2-generator presentations of deficiency $-1$ if ${\pi'\cong{F}\times(Z/nZ)}$ with $F\cong{Q(8)}$ or $I^*$ and $n>1$?
It seems unlikely that any of these groups have presentations of deficiency 0.

If $P=1$ or $Q(8)$ the weight class is unique up to inversion,
while $T^*_k$ and $I^*$ have 2 and 4 weight orbits, respectively,
by Theorem 14.1.
If $\pi'=T^*_1$ or $I^*$ each weight orbit is realized by a branched 
twist spun torus knot \cite{[PS87]}.

The group $\pi\tau_53_1\cong\mathbb{Z}\times I^*=
\mathbb{Z}\times SL(2,\mathbb{F}_5)$ 
is the common member of two families          
of high dimensional knot groups which are not otherwise 2-knot groups.
If $p$ is a prime greater than 3 then $SL(2,\mathbb{F}_p)$ 
is a finite superperfect group.
Let $e_p=\left(\smallmatrix 1&1\\0&1\endsmallmatrix\right)$.
Then $(1,e_p)$ is a weight element for $\mathbb{Z}\times SL(2,\mathbb{F}_p)$.
Similarly, $(I^*)^m$ is superperfect and $(1, e_5,\dots ,e_5)$
is a weight element for $G=\mathbb{Z}\times (I^*)^m$, for any $m\geq0$.
However $SL(2,\mathbb{F}_p)$ has cohomological period $p-1$ 
(see Corollary 1.27 of \cite{[DM85]}), while $\zeta (I^*)^m\cong(Z/2Z)^m$ 
and so $(I^*)^m$ does not have periodic cohomology if $m>1$.

Let $G(n)\cong{Q(8n)*_{Q(8)}}$ be the HNN extension with presentation
\[\langle{t,x,y}\mid x^{2n}=y^2,~yxy^{-1}=x^{-1},~ 
tyt^{-1}=yx^{-1},~ty^{-1}x^nt^{-1}=x^n\rangle.\]
Then $G(n)$ is a 2-knot group with an element of order $4n$,
for all $n\geq1$ \cite{[Kn89]},
and $\zeta{G(n)}=\langle{y^2}\rangle=Z/2Z$, if $n\geq2$.
The group $G(2)$ has the smallest finite base among knot groups $A*_C$ 
with infinitely many ends \cite{[BH17]}.

\section{Finite normal subgroups}

The results of this section perhaps belong more naturally in Chapter 3,
but are first used in this Chapter.
(The strategy may extend to the case of locally-finite
normal subgroups, but the argument used in \cite{[Hi17]} now seems flawed.)

\begin{theorem}
Let $X$ be a finite $PD_4$-complex with fundamental group 
$\pi$ and such that $\chi(X)=0$.
If $\beta_1^{(2)}(\pi)=0$ then $X$ is strongly minimal.
\end{theorem}

\begin{proof}
Since $X$ is a finite complex 
the $L^2$-Euler characteristic formula holds,
and so $\chi(X)=\beta_2^{(2)}(X)-2\beta_1^{(2)}(X)$.
Hence $\beta_2^{(2)}(X)=0$ also.
Since $\beta_1^{(2)}(X)=\beta_1^{(2)}(\pi)=0$,
$\beta_2^{(2)}(X)\geq\beta_2^{(2)}(\pi)\geq0$ and $\chi(X)=0$,
it follows that $\beta_2^{(2)}(X)=\beta_2^{(2)}(\pi)$.
Hence $H^2(c_X;\mathbb{Z}[\pi])$ is an isomorphism, 
by part (3) of Theorem 3.4.
\end{proof}

\begin{theorem}
Let $X$ be an orientable strongly minimal $PD_4$-complex.
If $\pi=\pi_1(X)$ has one end then it has no nontrivial finite normal subgroup.
\end{theorem}

\begin{proof}
Since $\pi$ has one end,  
$H_s(X;\mathbb{Z}[\pi])=0$ for $s\not=0$ or 2.
Poincar\'e duality and $c_X$ give an isomorphism 
$\Pi=H_2(X;\mathbb{Z}[\pi])\cong\overline{H^2(X;\mathbb{Z}[\pi])}$.
Since $X$ is strongly minimal, this in turn is isomorphic to $\overline{H^2(\pi;\mathbb{Z}[\pi])}$.
(See \S3.5 and \S10.8.)
The action of the covering group on the comology modules
$H^j(\pi;\mathbb{Z}[\pi])$ is induced by right multiplication 
on the coefficients.

Suppose that $\pi$ has a non-trivial finite normal subgroup $T$.
The LHS spectral sequence for $\pi$ as an extension of $\pi/T$
by $T$ collapses to give
$H^j(\pi;\mathbb{Z}[\pi])\cong{H^j}(\pi/T;\mathbb{Z}[\pi/T]))$,
and so $T$ acts trivially on $H^j(\pi;\mathbb{Z}[\pi])$, for all $j$.
Let $g\in{T}$ have order $p>1$, 
and let $C=\langle{g}\rangle\cong{Z/pZ}$.
Then $H_{i+3}(C;\mathbb{Z})\cong{H_i(C;\Pi)}$, 
for all $i\geq2$, by Lemma 2.10,
and so there is an exact sequence 
\[
0\to{Z/pZ}\to\Pi\to\Pi\to0,
\]
since $C$ has cohomological period 2 and acts trivially on $\Pi$. 
But $\Pi\cong\overline{H^2(\pi;\mathbb{Z}[\pi])}$ is torsion-free,
by Proposition 13.7.1 of \cite{[Ge]},
since $\pi$ is finitely presentable.
Hence there is no such element $g$ or subgroup $T$.
\end{proof}

\begin{cor}
If $\pi$ is a $2$-knot group with $\pi'$ infinite and $\beta_1^{(2)}(\pi)=0$
then $\pi$ has no non-trivial finite normal subgroup.
\end{cor}

\begin{proof}
If $e(\pi)=\infty$ then $\beta_1^{(2)}(\pi)>0$,
by an easy variation on Theorem 1.35(5) of \cite{[Lu]}.
Hence if  $\beta_1^{(2)}(\pi)=0$ and $\pi'$ is infinite then $e(\pi)=1$.
Thus the theorem applies to $X=M(K)$, for any 2-knot with group $\pi$.
\end{proof}

\section{The group $\Phi$}

The group $\Phi=Z*_2$ is the group of Example 10 of Fox.
This group is an ascending HNN extension with base $\mathbb{Z}$,
is metabelian, and has commutator subgroup isomorphic to the dyadic rationals. 
Since the 2-complex corresponding to the presentation 
$\langle a,t\mid tat^{-1}=a^2\rangle$ is aspherical, $g.d.\Phi=2$.

Example 10 of Fox is the boundary of the ribbon $D^3$ in $S^4$ obtained 
by ``thickening" a suitable immersed ribbon $D^2$ in $S^3$ for the 
stevedore's knot $6_2$ \cite{[Fo62]}.
Such a ribbon disc may be constructed by applying the method of 
\cite[\S1.7]{[H3]} to the equivalent presentation 
$\langle t,u,v\mid vuv^{-1}=t,\medspace tut^{-1}=v\rangle$ 
for $\Phi$ (where $u=ta$ and $v=t^2at^{-1}$).

The following alternative characterizations of $\Phi$ shall be useful.

\begin{theorem} 
Let $\pi$ be a $2$-knot group with 
maximal locally finite normal subgroup $T$.
Then $\pi/T\cong\Phi$ if and only if $\pi$ is elementary amenable and $h(\pi)=2$.
Moreover the following are equivalent:
\begin{enumerate}
\item $\pi$ has an abelian normal subgroup $A$ of rank $1$ 
such that $e(\pi/A)=2$;

\item $\pi$ is elementary amenable, $h(\pi)=2$ and $\pi$ 
has an abelian normal subgroup $A$ of rank $1$; 
                                
\item  $\pi$ is almost coherent, elementary amenable and $h(\pi)=2$; 

\item $\pi\cong\Phi$.
\end{enumerate}
\end{theorem}

\begin{proof} 
Since $\pi$ is finitely presentable and has infinite cyclic abelianization
it is an HNN extension $\pi\cong{B}*_\phi$ with base $B$ 
a finitely generated subgroup of $\pi'$, by Theorem 1.13. 
Since $\pi$ is elementary amenable the extension must be ascending.
Since $h(\pi'/T)=1$ and $\pi'/T$ has no nontrivial locally-finite normal subgroup
$[\pi'/T:\sqrt{\pi'/T}]\leq2$. 
The meridianal automorphism of $\pi'$ induces a meridianal automorphism on 
$(\pi'/T)/\sqrt{\pi'/T}$ and so $\pi'/T=\sqrt{\pi'/T}$.
Hence $\pi'/T$ is a torsion-free rank 1 abelian group.
Let $J=B/B\cap T$. Then $h(J)=1$ and $J\leq \pi'/T$ so $J\cong\mathbb{Z}$.
Now $\phi$ induces a monomorphism $\psi:J\to J$ and $\pi/T\cong J*_\psi$.
Since $\pi/\pi'\cong\mathbb{Z}$ we must have $J*_\psi\cong\Phi$. 

If (1) holds then $\pi$ is elementary amenable and $h(\pi)=2$.
Suppose (2) holds.  
We may assume without loss of generality that $A$ is 
the normal closure of an element of
infinite order, and so $\pi/A$ is finitely presentable.
Since $\pi/A$ is elementary amenable and $h(\pi/A)=1$, 
it is virtually $\mathbb{Z}$.
Therefore $\pi$ is virtually an HNN extension with base a finitely generated subgroup
of $A$, and so is coherent.
If (3) holds then $\pi\cong\Phi$, by Corollary 3.17.1.
Since $\Phi$ clearly satisfies conditions (1-3) this proves the theorem.
\end{proof}

\begin{cor}
If $T$ is finite and $\pi/T\cong\Phi$ then 
$T=1$ and $\pi\cong\Phi$.
\qed
\end{cor}

The group $\Phi$ is the only virtually torsion-free solvable 2-knot group 
which is not virtually poly-$Z$.
(See \S15.6 below.)
It is also the only solvable 2-knot group other than $\mathbb{Z}$ 
with deficiency 1 (see \S16.4 below) or in the class $\mathcal{X}$.

\begin{theorem}
Let $\pi$ be a $2$-knot group in $\mathcal{X}$ and with a nontrivial 
abelian normal subgroup $E$. 
Then either $\pi\cong\Phi$ or $E\leq\zeta\pi$, $\zeta\pi\cap\pi'=1$
and  $\pi$ is an iterated free product 
of (one or more) torus knot groups,
with amalgamation over copies of $\mathbb{Z}$.
\end{theorem}

\begin{proof} 
If $\pi$ is solvable then $\pi\cong Z*_m$,
for some $m\not=0$, by Corollary 2.6.1.
Since $\pi/\pi'\cong\mathbb{Z}$ we must have $m=2$ and so $\pi\cong\Phi$.

Otherwise $E\cong\mathbb{Z}$, 
$[\pi:C_\pi(E)]\leq 2$ and $C_\pi(E)'$ is free, 
of rank $>1$, by Theorem 2.7.
Hence $E\cap C_\pi(E)'=1$ and so $E$ maps injectively to $H=\pi/C_\pi(E)'$.
As $H$ has an abelian normal subgroup of index at most 2 
and $H/H'\cong\mathbb{Z}$, we must have $H\cong\mathbb{Z}$.
Hence  $C_\pi(E)=\pi$, and so $\pi'$ is free.
The further structure of $\pi$ is then due to Strebel \cite{[St76]}.
\end{proof}

For instance, the group $\pi$ with presentation 
$\langle{x,y,z}\mid{x^2=y^3,~y^2=z^3}\rangle$ is a 2-knot group, 
since $\pi=\langle\langle{x^{-1}z^2}\rangle\rangle$ 
and $\mathrm{def}(\pi)=1$. 
An Alexander polynomial calculation 
shows that $\pi$ is not a torus knot group,
although $\zeta\pi=\langle{y^6}\rangle\cong\mathbb{Z}$.

\section{Nilpotent normal subgroups}

The class of groups with infinite nilpotent normal subgroups 
includes the groups of torus knots and twist spins, $\Phi$, 
and all 2-knot groups with finite commutator subgroup.
If there is such a subgroup of Hirsch length $h>1$ 
the knot manifold is aspherical; 
this case is considered further in Chapter 16.

\begin{lemma} Let $G$ be a group with normal subgroups $T<N$
such that $T$ is a torsion group, $N$ is nilpotent, 
$N/T\cong\mathbb{Z}$ and $G/N$ is finitely generated. 
Then $G$ has an infinite cyclic normal subgroup $C\leq{N}$.
\end{lemma}

\begin{proof} 
Let $z\in{N}$ have infinite order, and
let $P=\langle\langle{z}\rangle\rangle$.
Then $P$ is nilpotent and $P/P\cap{T}\cong\mathbb{Z}$.
We may choose $g_0,\dots,g_n\in{G}$ such that their images generate $G/N$,
and such that $g_0zg_0^{-1}=z^\varepsilon{t_0}$ and $g_izg_i^{-1}=zt_i$, 
where $\varepsilon=\pm1$ and $t_i$ has finite order,
for all $0\leq{i}\leq{n}$.
The conjugates of $t_0,\dots,t_n$ generate the torsion subgroup of
$P/P'$, and so it has bounded exponent, $e$ say. 

Let $Q=\langle\langle{z^e}\rangle\rangle$.
If $P$ is abelian then we may take $C=Q$.
Otherwise, we may replace $P$ by $Q$,
which has strictly smaller nilpotency class, and the result
follows by induction on the nilpotency class.
\end{proof}

Since all 2-knot groups with finite commutator subgroup are known, 
we do not need to consider them in the next theorem.

\begin{theorem} 
Let $K$ be a $2$-knot with group $\pi=\pi K$.
If $\pi'$ is infinite and $\pi$ has a nilpotent normal subgroup $N$ 
with $h=h(N)\geq1$ then $h\leq4$ and
\begin{enumerate}
\item if $\pi'$ is finitely generated then so is $N$;

\item if $h=1$ and $N\not\leq\pi'$ or if $N$ is finitely generated 
then $N\cong\mathbb{Z}$ and either $M(K)$ is aspherical or $e(\pi/N)=\infty$;

\item if $h=2$ then $N\cong\mathbb{Z}^2$ and $M(K)$ is aspherical;

\item if $h=3$ or $4$ then $M(K)$ is homeomorphic to an
infrasolvmanifold;

\item if $N$ is not finitely generated then $h=1$ and
either $\pi\cong\Phi$ or $M(K)$ is aspherical, and so $N$ is 
torsion-free and abelian.
\end{enumerate}
In all cases $\sqrt\pi$ is nilpotent and torsion free.
\end{theorem}

\begin{proof}
Since $N$ and $\pi'$ are infinite, $e(\pi)=1$ and $\beta_1^{(2)}(\pi)=0$.
If $N\cap\pi'$ is finite then $N$ is finitely generated.
If $\pi'$ is finitely generated and $N\cap\pi'$ is infinite
then $\pi'$ has one end, since it cannot be virtually $\mathbb{Z}$.
Hence $M(K)$ is aspherical and $\pi'$ is a $PD_3$-group, by Theorem 4.5.
Then $N\cap\pi'$ is finitely generated, by Theorem 2.14, and so (1) follows.

Suppose that $h=1$. 
If $N\not\leq\pi'$ then $T=\pi'\cap{N}$ is the
torsion subgroup of $N$ and $N/T\cong\mathbb{Z}$.
There is an infinite cyclic  subgroup $S\leq{N}$ 
which is normal in $\pi$, by Lemma 15.9.
Since $S\cap\pi'=1$ and $Aut(S)$ is abelian,
$S$ is in fact central, and $\pi'S$ has finite index in $\pi$.
Hence $\pi'S\cong\pi'\times{S}$, so $\pi'$ is finitely presentable, 
and then $N$ is finitely generated, by (1).
If $N$ is finitely generated then its torsion subgroup is finite.
(See \cite[5.2.7 and 5.2.17]{[Ro]}.)
Hence $N$ is torsion-free, by Corollary 15.6.1, and so $N\cong\mathbb{Z}$.
Since $\pi'$ is infinite and no knot group is virtually $\mathbb{Z}^2$,
we have $e(\pi/N)\not=0$ or 2.
If $e(\pi/N)=1$ then $H^2(\pi;\mathbb{Z}[\pi])=0$, by Theorem 1.17, 
and so $M(K)$ is aspherical, by Corollary 3.5.2.
Otherwise $e(\pi/N)=\infty$.
This proves (2).

A Mayer-Vietoris argument as in Lemma 14.8 shows that if $e(\pi/N)=\infty$
then $\pi$ is not a $PD_4$-group.

Cases (3) and (4) follow immediately from Theorems 9.2 and 8.1, respectively.

If  $N$ is not finitely generated then $h=1$ and $N\leq\pi'$, by the above
results, and so $\pi'$ and $\pi/N$ are infinite.
If $\pi'/N$ is finite then $\pi'$ is virtually nilpotent and so 
$\pi$ is an ascending HNN extension 
with base $B$ a finitely generated virtually nilpotent group such that $h(B)=1$.
It follows easily that $\pi \cong\Phi$.

Suppose that $\pi'/N$ is infinite.
If $\pi'/N$ is finitely generated then $e(\pi/N)=1$,
and we set $G_1=\pi'$.
Otherwise,
$\pi'=\cup_n{G_n}$ is a strictly increasing union 
of subgroups $G_n$ which contain $N$,
and such that $G_n/N$ is finitely generated.
If every $G_n/N$ is finite then $\pi$ is elementary amenable 
and $h(\pi)=2$, so $\pi\cong\Phi$, by Theorem 15.7.
Hence we may assume that $G_n/N$ is infinite, for all $n$.
Since $G_n$ is finitely generated, it has an infinite cyclic normal subgroup $C_n$, by Lemma 15.9.
Let $J_n=G_n/C_n$. 
Then $G_n$ is an extension of $J_n$ by $\mathbb{Z}$,
and $J_n$ is in turn an extension of $G_n/N$ by $N/C_n$.
Since $G_n/N$ is finitely generated $J_n$ is an increasing union
$J_n=\cup_{m\geq1}{J_{n,m}}$ of finitely generated 
infinite subgroups which are extensions of $G_n/N$ by torsion groups.
If these torsion subgroups are all finite then 
$J_n$ is not finitely generated, since $N/C_n$ is infinite,
and ${[J_{n,m+1}:J_{n,m}]}$ is finite for all $m$.
Moreover, $H^0(J_{n,m};\mathbb{Z}[J_{n,m}])=0$, for all $m$,
since $J_{n,m}$ is infinite.
Therefore $H^s(J_n;F)=0$ for $s\leq1$ 
and any free $\mathbb{Z}[J_n]$-module $F$,
by the Gildenhuys-Strebel Theorem (with $r=1$).
Otherwise, $J_{n,m}$ has one end for $m$ large,
and we may apply Theorem 1.15 to obtain the same conclusion.
An LHSSS argument now shows that
$H^s(G_n;F)=0$ for $s\leq2$ and any free $\mathbb{Z}[G_n]$-module $F$.
Another application of Theorem 1.15 shows that
$H^s(\pi';F)=0$ for $s\leq2$ and any free $\mathbb{Z}[\pi']$-module $F$,
and then another LHSSS argument gives
$H^s(\pi;\mathbb{Z}[\pi])=0$ for $s\leq2$.
Hence $M(K)$ is aspherical, and so $N$ is torsion-free.
Thus (5) holds.

The final assertion is clear.
\end{proof}

\begin{cor} 
If $\pi'$ is infinite then $N$ is torsion-free.
\qed
\end{cor}

If $M(K)$ is aspherical and $N$ is a (torsion-free) 
abelian normal subgroup of rank 1 in $\pi$ then $N\leq\pi'$ 
and $Aut(N)$ is abelian,
so $N$ is a $\Lambda$-module.
It is not known whether $N$ must be finitely generated as a module.
(In fact, it is not known whether central subgroups of $PD_n$-groups
must be finitely generated as groups.)
Some further constraints on such groups may be found in \cite{[Hi18a]}.

\begin{cor}
If $\pi'$ is infinite and $\sqrt\pi$ is not nilpotent then 
$h(\sqrt\pi)\leq1$, the torsion subgroup of $\sqrt\pi$ is infinite,
and $e(\pi/\sqrt\pi)=1$ or $\infty$.
\end{cor}
          
\begin{proof} 
If $h(\sqrt\pi)=2$ then $\sqrt\pi\cong\mathbb{Z}^2$, by Theorem 9.2,
while if $h>2$ then $\pi$ is virtually poly-$Z$, by Theorem 8.1.
Finitely generated nilpotent groups are polycyclic.
If $e(\pi/\sqrt\pi)=0$ or 2 then $\pi$ is an ascending HNN extension
on a virtually nilpotent base, and we may apply Theorem 15.7.
Thus if $\sqrt\pi$ is not nilpotent then
$h(\sqrt\pi)\leq1$ and $e(\pi/\sqrt\pi)=1$ or $\infty$.
Moreover, $e(\pi)=1$ and $\beta_1^{(2)}=0$, and 
the torsion subgroup of $\sqrt\pi$ is non-trivial, 
since $\sqrt\pi$ is non-abelian.
Hence the torsion subgroup of $\sqrt\pi$ is infinite, 
by Corollary 15.6.1. 
\end{proof}

Is $\sqrt\pi$ always nilpotent?

\begin{cor}
If $\pi'$ is finitely generated and infinite then $N$ is torsion-free 
poly-$Z$ and either $\pi'\cap{N}=1$ or $M(K)'$ is homotopy equivalent 
to an aspherical Seifert fibred 3-manifold or 
$M(K)$ is homeomorphic to an infrasolvmanifold.
If moreover $h(\pi'\cap{N})=1$ then $\zeta\pi'\not=1$.
\end{cor}

\begin{proof} 
If $\pi'\cap{N}$ is torsion-free and $h(\pi'\cap{N})=1$
then $Aut(\pi'\cap{N})$ is abelian.
Hence $\pi'\cap{N}\leq\zeta\pi'$.
The rest of the corollary is clear from the theorem,
with Theorems 8.1 and 9.2.
\end{proof}

Twists spins of classical knots give examples with $N$ abelian,
of rank 1, 2 or 4, 
and $\pi'$ finite, $M$ aspherical or  $e(\pi/N)=\infty$.
If $\pi'\cap{N}=1$ and $\pi'$ is torsion free then
$\pi'$ is a free product of $PD_3^+$-groups and free groups.
Is every 2-knot $K$ such that $\zeta\pi\not\leq\pi'$ 
and $\pi$ is torsion-free $s$-concordant to a fibred knot?

We may construct examples of 2-knot groups $\pi$
such that $\zeta\pi'\not=1$ as follows.
Let $N$ be a closed 3-manifold such that $\nu=\pi_1(N)$ has weight 1 and 
$\nu/\nu'\cong\mathbb{Z}$, and let $w=w_1(N)$. 
Then $H^2(N;\mathbb{Z}^w)$ is infinite cyclic.
Let $M_e$ be the total space of the $S^1$-bundle over $N$ 
with Euler class $e\in H^2(N;\mathbb{Z}^w)$.
Then $M_e$ is orientable, and $\pi_1(M_e)$ has weight 1 
if $e=\pm 1$ or if $w\not=0$ and $e$ is odd.
Surgery on a weight class in $M_e$ gives $S^4$, 
so $M_e\cong M(K)$ for some 2-knot $K$. 

In particular, we may take $N$ to be the result of 0-framed surgery 
on a 1-knot.
If the 1-knot is $3_1$ or $4_1$ (i.e., is fibred of genus 1) 
then the resulting 2-knot group has commutator subgroup $\Gamma_1$. 
If instead we assume that the 1-knot is {\it not} fibred
then $N$ is not fibred \cite{[Ga87]} and so we get a 2-knot group $\pi$
with $\zeta\pi'\cong\mathbb{Z}$ but $\pi'$ not finitely generated.

For examples with $w\not=0$ we may take one of the groups
\[
\langle{t,z,a_1,\dots,a_n,b_1,\dots,b_n}\mid\Pi[a_i,b_i]=z^q,~tzt^{-1}=z^{-1},
~ta_it^{-1}=b_{n-i},\]
\[
tb_it^{-1}=a_{n-i}b_{n-i} ~\forall{i}\rangle,
\]
where $q$ and $n$ are odd.
When $n=1$ we get all the other 2-knot groups with 
commutator subgroup $\Gamma_q$. (See Theorem 16.14.)
When $n>1$ these groups have $\sqrt\pi\cong\mathbb{Z}<\pi'$,
with image $Z/qZ$ in $\pi'/\pi''$.

\begin{theorem}
Let $\pi$ be a $2$-knot group such that $\zeta\pi$ is infinite.
Then either
\begin{enumerate}
\item$e(\pi)=1$ and $\zeta\pi\cong\mathbb{Z}^2$,
or $\zeta\pi$ is torsion-free of rank $1$;  or 
\item$e(\pi)=2$ and $\zeta\pi\cong\mathbb{Z}$ or $\mathbb{Z}\oplus{Z/2Z}$.
\end{enumerate}
\end{theorem}

\begin{proof}
If $\pi$ has one end and $\zeta\pi$ is infinite
then $\pi$ has no nontrivial finite normal subgroup,
by Corollary 15.6.1.
Hence $\zeta\pi$ is torsion free.
If $h(\zeta\pi)>2$ then $\pi$ is torsion-free and
virtually poly-$Z$, and $h(\pi)=4$.
In particular, $\pi'$ is infinite, 
and so $[\pi:\zeta\pi]=\infty$ \cite[10.1.4]{[Ro]}.
Hence $h(\zeta\pi)=3$.
But then $\pi/\zeta\pi$ has two ends and cyclic abelianization,
so surjects onto $\mathbb{Z}$.
Hence $\zeta\pi\leq\pi'$, and $[\pi':\zeta\pi]<\infty$.
Therefore $\pi''$ is finite \cite[10.1.4]{[Ro]}, and so $\pi''=1$.
Since the meridianal automorphism fixes the centre, 
this is incompatible with $\pi/\pi'\cong\mathbb{Z}$.
Therefore either $h(\zeta\pi)=2$, in which case $\zeta\pi\cong\mathbb{Z}^2$, 
or $\zeta\pi$ is torsion-free of rank $\leq1$.

If $\pi$ has two ends then $\pi'$ is finite, 
and so $\zeta\pi$ is finitely generated and of rank 1.
It follows from the classification of such 2-knot groups that
$\zeta\pi\cong\mathbb{Z}\oplus{Z/2Z}$ if $\pi'$ has even order, 
and otherwise $\zeta\pi\cong\mathbb{Z}$.
\end{proof}

In all known cases the centre of a 2-knot group is cyclic, 
$\mathbb{Z}\oplus(Z/2Z)$ or $\mathbb{Z}^2$.

The knots $\tau_04_1$, the trivial knot, $\tau_33_1$ and $\tau_63_1$ 
are fibred and their groups have centres 1, $\mathbb{Z}$, 
$\mathbb{Z}\oplus(Z/2Z)$ and $\mathbb{Z}^2$, 
respectively.
Let $G=\pi\tau_63_1$ and let $\gamma$ be a loop in $X(\tau_63_1)$ which
represents a generator of $G'\cap\zeta{G}\cong\mathbb{Z}$.
If $K$ is the satellite knot $\Sigma(\tau_23_1;\tau_63_1,\gamma)$ 
and $\pi=\pi{K}$ then $\zeta\pi\cong\mathbb{Z}$ 
and $\pi/\zeta\pi$ has infinitely many ends \cite{[Yo82]}.
This paper also gives an example with a minimal Seifert hypersurface 
and such that $\sqrt\pi=\zeta\pi=Z/2Z$.
Yoshikawa's construction may be extended as follows. 
Let $q>0$ be odd and let $k_q$ be a 2-bridge knot such that 
the 2-fold branched cyclic cover of $S^3$,
branched over $k_q$ is a lens space $L(3q,r)$, 
for some $r$ relatively prime to $q$.
Let $K_1=\tau_2k_q$ be the 2-twist spin of $k_q$, and let $K_2=\tau_3k$ 
be the 3-twist spin of a nontrivial knot $k$. 
Let $\gamma$ be a simple closed curve 
in $X(K_1)$ with image $[\gamma]\in\pi{K_1}$ of order $3q$, 
and let $w$ be a meridian for $K_2$.
Then $w^3$ is central in $\pi{K_2}$.
The group of the satellite of $K_1$ about $K_2$ relative to $\gamma$ is the generalized free
product
\[
\pi=\pi{K_2}/\langle\langle{w^{3q}}\rangle\rangle*_{w=[\gamma]}\pi{K_1}.
\]
Hence $\sqrt\pi=\langle{w^3}\rangle\cong{Z/qZ}$, while $\zeta\pi=1$.
If we use a 2-knot $K_1$ with group $(Q(8)\times{Z/3qZ})\rtimes_\theta\mathbb{Z}$
instead and choose $\gamma$ so that $[\gamma]$ has order $6q$ then we obtain
examples with $\sqrt\pi\cong{Z/2qZ}$ and $\zeta\pi=Z/2Z$.

\section{Almost coherent, restrained and locally virtually indicable}

In this section we use the fact that knot groups are HNN extensions 
over finitely generated bases, by the Bieri-Strebel Theorem,
to show that the tetrachotomy of the introduction holds,
under mild coherence hypotheses on $\pi{K}$ or $N$.

\begin{theorem} 
Let $K$ be a $2$-knot whose group $\pi=\pi K$
is an ascending HNN extension over an $FP_2$ base $B$ with finitely many ends.
Then either $\pi'$ is finite or $\pi\cong\Phi$ or $M(K)$ is aspherical.
\end{theorem}

\begin{proof}                       
This follows from Theorem 3.17, since a group with abelianization 
$\mathbb{Z}$ cannot be virtually $\mathbb{Z}^2$.
\end{proof}

\begin{cor}
If $B$ is $FP_3$ and has one end then 
$\pi'=B$ and is a $PD_3^+$-group.
\end{cor}

\begin{proof} This follows from \cite[Lemma 3.4]{[BG85]}, 
as in Theorem 2.13.
\end{proof}

Is $M(K)$ aspherical whenever $B$ is finitely generated and one-ended?
Does the corollary hold if we assume only that $B$ is $FP_2$ and has one end?

\begin{cor}
If $\pi$ is an ascending HNN extension over an $FP_2$ base $B$ 
and has an infinite restrained normal subgroup $A$ then 
either $\pi'$ is finite or $\pi\cong\Phi$ or $M(K)$ is aspherical 
or $\pi'\cap A=1$ and $\pi/A$ has infinitely many ends.
\end{cor}

\begin{proof} 
If $B$ is finite or $A\cap{B}$ is infinite then $B$ has 
finitely many ends (cf.\ Corollary 1.15.1) and Theorem 15.12 applies.
Therefore we may assume that $B$ has infinitely many ends 
and $A\cap{B}$ is finite.
But then $A\not\leq\pi'$, so $\pi$ is virtually $\pi'\times\mathbb{Z}$.
Hence $\pi'=B$ and $M(K)'$ is a $PD_3$-complex.
In particular, $\pi'\cap A=1$ and $\pi/A$ has infinitely many ends.
\end{proof}

If $K$ is the $r$-twist spin of an irreducible 1-knot then the $r^{th}$ power
of a meridian is central in $\pi$ and either $\pi'$ is finite or $M(K)$ 
is aspherical. (See \S3 of Chapter 16.)
The final possibility is realized by Artin spins of nontrivial torus knots.

\begin{theorem} 
Let $K$ be a $2$-knot whose group $\pi=\pi K$ is an HNN extension 
with $FP_2$ base $B$ and associated subgroups $I$ and $\phi(I)=J$.
If $\pi$ has a restrained normal subgroup $N$ which 
is not locally finite and $\beta_1^{(2)}(\pi)=0$ then either 
$\pi'$ is finite or $\pi\cong\Phi$ or $M(K)$ is aspherical or 
$N$ is locally virtually $\mathbb{Z}$ and $\pi/N$ has infinitely many ends.
\end{theorem}

\begin{proof}
If $\pi'\cap N$ is locally finite then it follows from Britton's lemma
(on normal forms in HNN extensions) that either $B\cap N=I\cap N$ or
$B\cap N= J\cap N$. 
Moreover $N\not\leq\pi'$ (since $N$ is not locally finite), 
and so $\pi'/\pi'\cap N$ is finitely generated.
Hence $B/B\cap N\cong I/I\cap N\cong J/J\cap N$.
Thus either $B=I$ or $B=J$ and so the HNN extension is ascending.
If $B$ has finitely many ends we may apply Theorem 15.12.
Otherwise $B\cap N$ is finite, 
so $\pi'\cap N=B\cap N$ and $N$ is virtually $\mathbb{Z}$.
Hence $\pi/N$ is commensurable with $B/B\cap N$, and $e(\pi/N)=\infty$.

If $\pi'\cap N$ is locally virtually $\mathbb{Z}$ and $\pi/\pi'\cap N$ 
has two ends then $\pi$ is an ascending HNN extension on a two-ended base.
Hence $\pi$ is almost coherent, elementary amenable and $h(\pi)=2$, 
so $\pi\cong\Phi$, by Theorem 15.7. 
Otherwise we may assume that either $\pi/\pi'\cap N$ has one end or
$\pi'\cap N$ has a finitely generated, one-ended subgroup.
In either case $H^s(\pi;\mathbb{Z}[\pi])=0$ for $s\leq2$, by Theorem 1.18,
and so $M(K)$ is aspherical, by Theorem 3.5.                 
\end{proof}

Note that $\beta_1^{(2)}(\pi)=0$ if $N$ is amenable. 
Every knot group is an HNN extension with finitely generated base 
and associated subgroups, by Theorem 1.13, 
and in all known cases these subgroups are $FP_2$.

\begin{theorem} 
Let $K$ be a $2$-knot such that $\pi=\pi K$
has an almost coherent, locally virtually indicable, 
restrained normal subgroup $E$ which is not locally finite.
Then either $\pi'$ is finite or $\pi\cong\Phi$ or $M(K)$ is aspherical
or $E$ is locally virtually $\mathbb{Z}$ and $\pi/E$ has one 
or infinitely many ends.
\end{theorem} 

\begin{proof}
Let $F$ be a finitely generated subgroup of $E$.
Since $F$ is $FP_2$ and virtually indicable it has a subgroup of finite index 
which is an HNN extension over a finitely generated base, by Theorem 1.13.
Since $F$ is restrained the HNN extension is ascending, 
and so $\beta_1^{(2)}(F)=0$, by Lemma 2.1.
Hence $\beta_1^{(2)}(E)=0$ and so $\beta_1^{(2)}(\pi)=0$ 
\cite[Theorem 7.2.(3)]{[Lu]}.
 
If $E$ is locally virtually $\mathbb{Z}$
then $E$ is elementary amenable and $h(E)=1$.
If $\pi/E$ is finite then $\pi'$ is finite.
If $\pi/E$ has two ends then $\pi$ is almost coherent,
elementary amenable and $h(\pi)=2$, and so $\pi\cong\Phi$, by Theorem 15.7.
If $E$ has a finitely generated, one-ended subgroup
and $\pi$ is not elementary amenable of Hirsch length 2
then $H^s(\pi;\mathbb{Z}[\pi])=0$ for $s\leq2$, by Theorem 1.17.
Hence $M(K)$ is aspherical, by Theorem 3.5.

The remaining possibility is that $E$ is locally virtually $\mathbb{Z}$ 
and $\pi/E$ has one or infinitely many ends.
\end{proof}

Does this theorem hold without any coherence hypothesis? 
Note that the other hypotheses hold if $E$ is elementary amenable and 
$h(E)\geq2$.
If $E$ is abelian and $\pi/E$ has one end then $M(K)$ is aspherical,
by Theorems 1.17 and 3.5.
If $M(K)$ is aspherical then $c.d.E\cap\pi'\leq3$.
Slightly stronger hypotheses on $E$ then imply that $\pi$ has
a nontrivial torsion-free abelian normal subgroup.

\begin{theorem} 
Let $N$ be a group which is either elementary amenable or locally $FP_3$, 
virtually indicable and restrained. 
If $c.d.N\leq 3$ then $N$ is virtually solvable.
\end{theorem}

\begin{proof} 
Suppose first that $N$ is locally $FP_3$ and virtually indicable,
and let $E$ be a finitely generated subgroup of $N$ 
which maps onto $\mathbb{Z}$.
Then $E$ is an ascending HNN extension $B*_\phi$ 
with $FP_3$ base $B$ and associated subgroups.
If $c.d.B=3$ then $H^3(B;\mathbb{Z}[E])\cong
{H^3(B;\mathbb{Z}[B])\otimes_{\mathbb{Z}[B]}\mathbb{Z}[E]}\not=0$
and the homomorphism $H^3(B;\mathbb{Z}[E])\to{H^3(B;\mathbb{Z}[E])}$ 
in the Mayer-Vietoris sequence for the HNN extension is not onto
\cite[Lemma 3.4 and Remark 3.5]{[BG85]}.
But then $H^4(E;\mathbb{Z}[E])\not=0$, contrary to $c.d.N\leq3$.
Therefore $c.d.B\leq2$, and so $B$ is elementary amenable, by Theorem 2.7.
Hence $N$ is elementary amenable, and so is virtually solvable by Theorem 1.11. 
\end{proof}

In particular, $\zeta\sqrt N$ is a nontrivial, torsion-free abelian characteristic subgroup of $N$. 
A similar argument shows that if $N$ is locally $FP_n$, virtually indicable,
restrained and $c.d.N\leq n$ then $N$ is virtually solvable.

The final result of this section extends a result of \cite{[Yo97]},
where it was assumed that the centre is finitely generated.

\begin{theorem} 
Let $K$ be a $2$-knot with a minimal Seifert 
hypersurface, and let $\pi=\pi{K}$.
Then $\zeta\pi\cong 1$, $Z/2Z$, $\mathbb{Z}$, 
$\mathbb{Z}\oplus(Z/2Z)$ or $\mathbb{Z}^2$.
\end{theorem}

\begin{proof} 
Since $K$ has a minimal Seifert 
hypersurface, $V$ say, $\pi$ is an HNN extension
$\pi=HNN(B;\phi:I\cong J)$, with $B$ finitely presentable and
$I\cong J\cong\pi_1(V)$.
If the HNN extension is ascending then the result follows from 
Corollary 15.12.1, since $\pi'$ is then a $PD_3$-group or is finite.
If the HNN extension is not ascending then $\pi'$ is the increasing union 
of iterated generalized free products of conjugates of $B$, 
amalgamated over conjugates of $I$ and $J$,
and so $\zeta\pi\cap\pi'\leq{I\cap{J}}$.
Since $I$ is a 3-manifold group, 
$\zeta{I}$ is finitely generated,
and therefore so are $\zeta\pi\cap\pi'$ and $\zeta\pi$.
If $\pi'$ is finite then $\zeta\pi\cong\mathbb{Z}$ or 
$\mathbb{Z}\oplus(Z/2Z)$, by inspection of the groups in Theorem 15.4.
Otherwise $\zeta\pi$ is torsion free, by Theorem 15.10, and so $\zeta\pi=1$, $\mathbb{Z}$ or $\mathbb{Z}^2$, by Theorem 15.11.
\end{proof}

\section{A Tits alternative for 2-knot groups}

In this section we shall show that virtually torsion-free
solvable 2-knot groups other than $\Phi$ and with 
infinite commutator subgroup are polycyclic. 
These groups shall be determined explicitly in \S4 of Chapter 16.
We shall show also that a version of the Tits alternative holds 
for 2-knot groups.
The class of groups considered next probably includes all 
elementary amenable 2-knot groups and perhaps even all
restrained 2-knot groups (giving the full Tits alternative).

\begin{theorem} 
Let $\pi$ be a $2$-knot group. 
Then the following are equivalent:
\begin{enumerate}
\item $\pi$ is restrained, locally $FP_3$ and locally virtually 
indicable;

\item $\pi$ is restrained and is an ascending HNN extension $B*_\phi$, 
where $B$ is $FP_3$  
and virtually indicable;

\item $\pi$ is elementary amenable and has a torsion-free normal
subgroup $N\not=1$;

\item $\pi$ is elementary amenable and has an abelian 
normal subgroup $A$ of rank $\geq1$;

\item $\pi$ is elementary amenable and is an ascending HNN extension 
$B*_\phi$, where $B$ is $FP_2$;

\item $\pi'$ is finite or $\pi\cong\Phi$
or $\pi$ is torsion-free virtually poly-$Z$ and $h(\pi)=4$.
\end{enumerate}
\end{theorem}                    

\begin{proof} 
Condition (1) implies (2) by Theorem 1.13,
since an HNN extension (such as a knot group) is restrained if and only if 
it is ascending and the base is restrained.
If (2) holds then either $\pi'$ is finite or $\pi\cong\Phi$ or 
$\pi'=B$ and is a $PD_3$-group, by Theorem 15.12 and Corollary 15.12.2.
In the latter case $B$ has a subgroup of finite index 
which maps onto $\mathbb{Z}^2$,
since it is virtually indicable and admits a meridianal automorphism.
Hence $B$ is virtually poly-$Z$, by Corollary 2.13.1 
(together with the remark following it).
Hence (2) implies (6). 
If (3) holds then either $H^s(N;\mathbb{Z}[\pi])=0$ for all $s\geq0$ or
$N$ is virtually solvable \cite[Proposition 3]{[Kr93']}.
Hence either $\pi$ is torsion-free virtually poly-$Z$ and $h(\pi)=4$,
by Theorem 8.1, or (4) holds.
If (4) holds and $\pi'$ is infinite then $h(\pi)>1$,
and either $\pi\cong\Phi$, 
by Theorem 15.7 (if $h(\pi)=2$),
or $M(K)$ is aspherical, 
by Theorems 1.17 and 8.1 (if $h(\pi)\geq3$). 
Hence (4) implies (6). 
Condition (5) implies (6)  by Theorem 15.12. 
On the other hand (6) implies (1-5).
\end{proof}

In particular, if  a 2-knot $K$ has a minimal Seifert hypersurface 
and $\pi=\pi{K}$ is restrained then these conditions hold,
by the Tits alternative for 3-manifold groups \cite{[EJ73]}.

\section{Abelian HNN bases}

We shall complete Yoshikawa's study of 2-knot groups which are
HNN extensions with abelian base.
The first four paragraphs of the following proof outline the arguments of
\cite{[Yo86],[Yo92]}. (Our contribution is the final paragraph, 
which shows that there is no torsion when the base has rank 1.)
Let $BS(p,q)$ be the {\it Baumslag-Solitar group\/} with presentation
$\langle{t,x}\mid {tx^pt^{-1}=x^q}\rangle$, where $p,q\not=0$.

\begin{theorem} 
Let $\pi$ be a $2$-knot group which is an HNN extension with abelian base.
Then either $\pi$ is metabelian or $\pi\cong{BS(n,n+1)}$ for some $n>1$.
\end{theorem}

\begin{proof} 
Suppose that $\pi=A*_\phi=HNN(A;\phi:B\to C)$ where $A$ is abelian.
Let $j$ and $j_C$ be the inclusions of $B$ and $C$ into $A$,
and let $\tilde\phi=j_C\phi$.
Then $\tilde\phi-j:B\to A$ is an isomorphism, by the Mayer-Vietoris sequence 
for homology with coefficients $\mathbb{Z}$ for the HNN extension.
Hence $rank(A)=rank(B)=r$, say, and the torsion subgroups $TA$, $TB$ and $TC$
of $A$, $B$ and $C$ coincide.

Suppose first that $A$ is not finitely generated.
Since $\pi$ is finitely presentable and $\pi/\pi'\cong\mathbb{Z}$
it is also an HNN extension with finitely generated base 
and associated subgroups, by the Bieri-Strebel Theorem (1.13).
We may assume the base is a subgroup of $A$.
Considerations of normal forms imply that the latter HNN structure 
must be ascending, and so $\pi$ is metabelian \cite{[Yo92]}.

Assume now that $A$ is finitely generated.
Then $TA=TB=TC$, 
since $TB$ and $TC\leq{TA}$ and they have the same finite order.
Therefore the image of $TA$ in $\pi$ is a finite normal subgroup $N$.
If $r=0$ then clearly $B=A$ and so $\pi$ is metabelian.
If $r>0$ then $\pi$ has one end, by a Mayer-Vietoris argument 
with coefficients $\mathbb{Z}[\pi]$,
and $\beta_i^{(2)}(\pi)=0$, for all $i$ \cite[Theorem 7.2]{[Lu]}.
But then $TA$ must be trivial, by Theorem 15.2.
Let $M_j=|det(j)|$ and $M_\phi=|det(\phi)|$. 
Applying the Mayer-Vietoris sequence for homology with coefficients
$\Lambda$, we find that $t\tilde\phi-j$ is injective and 
$\pi'/\pi''\cong H_1(\pi;\Lambda)$ has rank $r$ as an abelian group.
Now $H_2(A;\mathbb{Z})\cong A\wedge A$ 
(see \cite[page 334]{[Ro]}) and so 
$H_2(\pi;\Lambda)\cong\mathrm{Cok}(t\wedge_2\tilde\phi-\wedge_2j)$
has rank $\binom{r}2$.
Let $\delta_i(t)=\Delta_0(H_i(\pi;\Lambda))$, for $i=1$ and 2.
Then $\delta_1(t)=det(t\phi-j)$ and $\delta_2(t)=det(t\phi\wedge\phi-j\wedge{j})$.
Moreover $\delta_2(t^{-1})$ divides $\delta_1(t)$, by Theorem 14.3.
In particular, $\binom r2\leq r$, and so $r\leq3$.
If $r=2$ then $\binom r2=1$ and $\delta_2(t)=\pm(tM_\phi-M_j)$.
Comparing coefficients of the terms of highest and lowest degree in
$\delta_1(t)$ and $\delta_2(t^{-1})$, we see that $M_j=M_\phi$,
so $\delta_2(1)\equiv0$ mod $(2)$, 
which is impossible since $|\delta_1(1)|=1$.
If $r=3$ a similar comparison of coefficients in $\delta_1(t)$ and 
$\delta_2(t^{-1})$ shows that $M_j^3$ divides $M_\phi$ and $M_\phi^3$ 
divides $M_j$, so $M_j=M_\phi=1$.
Hence $\phi$ is an isomorphism, and so $\pi$ is metabelian.

There remains the case $r=1$, when $A,B$ and $C\cong\mathbb{Z}$.
Let $x$ be a generator for $A$. Then $B$ and $C$ are  generated by 
$x^n$ and $x^m$ for some exponents $m,n>0$, and $\phi(x^n)=x^m$.
We may assume that $m\geq{n}$. 
It then follows from the considerations of the 
first paragraph that $m=n+1$, and so $\pi\cong{BS(n,n+1)}$.
Since $\pi$ is metabelian if $n=1$ this completes the proof.
\end{proof}

Every knot group in $\mathcal{X}$ 
is either a quotient of a torus knot group or a quotient of $BS(n,n+1)$, 
for some $n\geq1$ \cite{[DM18]}.
We may write $\pi=\pi\mathcal{G}$, 
where $(\mathcal{G},\Gamma)$ is a finite graph of groups 
with all vertex and edge groups infinite cyclic.
Then $\beta_1(G)\geq\beta_1(\Gamma)$, 
and so either $\Gamma$ is a tree and $\zeta\pi\cong\mathbb{Z}$ 
or $\pi$ is an HNN extension with base such a group.
A key element in \cite{[DM18]} is the use of \cite{[Ho02]} to show that 
since $\pi$ has weight 1 the vertices of $\Gamma$ have valence $\leq2$.
In each case $\pi$ has a 2-generator presentation.
Does it have a 1-relator presentation?

The 2-knot group $\pi$ with presentation 
$\langle{x,y,z}\mid{x^2=y^3,~y^2=z^3}\rangle$ 
has $\zeta\pi\cong\mathbb{Z}$, 
and is a quotient of the group of the (4,9)-torus knot,
since $x^4=z^9$ in $\pi$.
However $\pi$ is not a torus knot group, as observed in \S2.
(This group also has the 1-relator presentation
$\langle{x,z}\mid{z^6=x^2z^{-3}x^2}\rangle$.) 

The 2-knot group $G$ with presentation $\langle{a,b,t}\mid{a^2=b^3,~ta^3t^{-1}=b^4}\rangle$
is a quotient of $BS(8,9)$, and has the same Alexander polynomial.
However $\zeta{G'}\not=1$, and so $G\not\cong{BS(8,9)}$.
(It is again a 1-relator group, 
but is not an HNN extension with free base, 
and so does not have a 1-relator Wirtinger presentation 
\cite[Theorems 1 and 2]{[Yo88]}.)

\section{Locally-finite normal subgroups}

If Theorem 15.6 could be extended to exclude non-trivial {\it locally}-finite normal subgroups $T$ then 
it would follow that (for a 2-knot group  $\pi$)
\begin{enumerate}
\item$\sqrt\pi$ is nilpotent;
\item{}if $\sqrt\pi$ has nontrivial torsion then either $\pi'$ is finite or $\sqrt\pi$ is finite and $\beta_1^{(2)}(\pi)>0$;
\item{}if $\pi$ is elementary amenable then it is virtually torsion-free solvable,
as in part (6) of Theorem 15.17.
\item{}if $\pi$ has an infinite elementary amenable normal subgroup $E$
then either $\pi$ is elementary amenable or  $E$ is torsion-free abelian and $c.d.E\leq2$.
\end{enumerate}
Moreover, the hypothesis ``$H^2(\pi;\mathbb{Z}[\pi])=0$" in Theorem 9.1 
could be dropped.
(We may assume that $T=\langle\langle{x}\rangle\rangle_\pi$ is the normal closure
of an element $x$ of prime order $p>1$, since $\pi$ has no finite normal subgroup,
by Corollary 15.6.1.)

The argument for the key Lemma 3.2 of \cite{[Hi17]}
(which asserts that any such $T$ should act trivially on $H^2(\pi;\mathbb{Z}[\pi])$)
now seems inadequate.
There is a near counter-example; 
it fails for $G=A\rtimes\mathbb{Z}$, 
where $A\cong(Z/pZ)^\infty$ is the elementary abelian 
group of exponent $p$ on generators $a_n$, 
and $\mathbb{Z}$ acts by translating the indices.
(This group is not $FP_2$.)
The lemma remains true when $T$ is an increasing union of finite 
normal subgroups,
but this is only enough to establish the result on centres:
either $\pi'$ is finite (of even order) or $\zeta\pi$ is torsion-free.
Nevertheless there seems to be considerable evidence suggesting that the
above assertions should hold for all 2-knot groups.
(This would establish our tetrachotomy for the cases with $N$ elementary amenable.)

In earlier versions of this section we cited at length 
some results from \cite[Chapter 6]{[H1]},
which considered a class of groups between solvable and elementary amenable,
and showed that if $\pi$ is a 2-knot group in this class then $\pi$ 
has a maximal locally finite normal subgroup $T$ such that either $\pi/T$ is virtually poly-$Z$ of Hirsch length 4
or $\pi/T\cong\Phi$. 
We  retain only the following lemma.

\begin{lemma} 
Let $G$ be an $FP_2$ group with a torsion normal subgroup $T$ such
that $\mathbb{Z}[G/T]$ is coherent.
Then $T/T'$ has finite exponent.
In particular, if $G$ is solvable then $T=1$ if and only if
$H_1(T;\mathbb{F}_p)=0$ for all primes $p$.
\end{lemma}

\begin{proof} 
Let $C_*$ be a free $\mathbb{Z}[G]$-resolution of the augmentation module 
$\mathbb{Z}$ which is finitely generated in degrees $\leq2$.
Then $T/T'=H_1(\mathbb{Z}[G/T]\otimes_G C_*)$,
and so it is finitely presentable as a $\mathbb{Z}[G/T]$-module,
since $\mathbb{Z}[G/T]$ is coherent.
If $T/T'$ is generated by elements $t_i$ of order $e_i$ then
$\Pi e_i$ is a finite exponent for $T/T'$.
If $G$ is solvable then so is $T$, 
and $T=1$ if and only if $T/T'=1$. 
Since $T/T'$ has finite exponent $T/T'=1$ 
if and only if $H_1(T;\mathbb{F}_p)=0$ for all primes $p$.
\end{proof}

If $G/T$ is virtually $Z*_m$ or virtually poly-$Z$ then 
$\mathbb{Z}[G/T]$ is coherent \cite{[BS79]}.
Note also that $\mathbb{F}_p[Z*_m]$ is a coherent Ore domain 
of global dimension 2,
while if $J$ is a torsion-free virtually poly-$Z$ group then
$\mathbb{F}_p[J]$ is a noetherian Ore domain of global dimension $h(J)$.
(See \cite[\S 4.4 and \S 13.3]{[Pa]}.)

%% file: m5-16.tex
\chapter{Abelian normal subgroups of rank $\geq2$}

If $K$ is a 2-knot such that $h(\sqrt{\pi K})=2$ then 
$\sqrt{\pi K}\cong\mathbb{Z}^2$,
by Theorems 9.2 and 15.10.
The most familiar examples of such knots are the branched twist spins
of torus knots, but there are others.
In all cases, the knot manifolds are $s$-cobordant to
$\mathbb{SL}\times\mathbb{E}^1$-manifolds.
The groups of the branched twist spins of torus knots 
are the 3-knot groups which are $PD_4^+$-groups 
and have centre of rank 2, with some power of a weight element being central.
There are related characterizations of the groups of branched twist spins 
of other prime 1-knots.

If $h(\sqrt{\pi K})>2$ then $M(K)$ is homeomorphic to an infrasolvmanifold,
so $K$ is determined up to Gluck reconstruction and
change of orientations by $\pi{K}$ and a weight orbit.
We find the groups and their strict weight orbits,
and give optimal presentations for all such groups.
Two are virtually $\mathbb{Z}^4$.
Whether $K$ is amphicheiral or invertible is known, 
and in most cases whether it is reflexive is also known. 
(See Chapter 18 and \cite{[Hi11']}.)

\section{The Brieskorn manifolds $M(p,q,r)$}

Let $M(p,q,r)=\{(u,v,w)\in \mathbb{C}^3\mid u^p+v^q+w^r=0\}\cap S^5$.
Thus $M(p,q,r)$ is a {\it Brieskorn $3$-manifold} (the link of an isolated 
singularity of the intersection of $n$ algebraic hypersurfaces in 
$\mathbb{C}^{n+2}$, for some $n\geq1$).
It is clear that $M(p,q,r)$ is unchanged by a 
permutation of $\{ p,\medspace q,\medspace r\}$.

Let $s= hcf\{ pq,pr,qr\}$.
Then $M(p,q,r)$ admits an effective $S^1$-action 
given by $z(u,v,w)=(z^{qr/s}u,z^{pr/s}v,z^{pq/s}w)$ for $z\in S^1$ and 
$(u,v,w)\in M(p,q,r)$, and so is Seifert fibred.
More precisely, let $\ell=lcm\{ p, q,r\}$, $p'=lcm\{p,r\}$, 
$q'=lcm\{q,r\}$ and $r'=lcm\{ p,q\}$, 
$s_1=qr/q'$, $s_2=pr/p'$ and $s_3=pq/r'$ and
$t_1=\ell/q'$, $t_2=\ell/p'$ 
and $t_3=\ell/r'$.
Let $g=(2+(pqr/\ell)-s_1-s_2-s_3)/2$.
Then $M(p,q,r)=M(g;s_1(t_1,\beta_1),s_2(t_2,\beta_2),s_3(t_3,\beta_3))$, 
in the notation of \cite{[NR78]}, 
where the coefficients $\beta_i$ are determined {\it mod} $t_i$ 
by the equation
\begin{equation*}
e=-(qr\beta_1+pr\beta_2+pq\beta_3)/\ell)=-pqr/\ell^2
\end{equation*}
for the generalized Euler number.
(See \cite{[NR78]}.)
If $p^{-1}+q^{-1}+r^{-1}\leq1$ the Seifert fibration is essentially unique. 
(See \cite[Theorem 3.8]{[Sc83']}.)
In most cases the triple $\{ p,q,r\}$ is determined by
the Seifert structure of $M(p,q,r)$. 
(Note however that, for example, $M(2,9,18)\cong M(3,5,15)$ \cite{[Mi75]}.)

The map $f:M(p,q,r)\to \mathbb{CP}^1$ given by $f(u,v,w)=[u^p:v^q]$ is 
constant on the orbits of the $S^1$-action, 
and the exceptional fibres are those above $0$, $-1$ and $\infty$ in 
$\mathbb{CP}^1$.
In particular, if $p$, $q$ and $r$ are pairwise relatively prime $f$ is the
orbit map and $M(p,q,r)$ is Seifert fibred over the orbifold $S^2(p,q,r)$.
The involution $c$ of $M(p,q,r)$ induced by complex conjugation in $\mathbb{C}^3$ 
is orientation preserving and is compatible with $f$ and complex conjugation 
in $\mathbb{CP}^1$.
Let $A(u,v,w)=(u,v,e^{2\pi i/r}w)$ and $h(u,v,w)=(u,v)/(|u|^2+|v|^2)$,
for $(u,v,w)$ $\in M(p,q,r)$.
The $Z/rZ$-action generated by $A$ commutes with the $S^1$-action,
and these actions agree on their subgroups of order $r/s$.
The projection to the orbit space $M(p,q,r)/\langle A\rangle$
may be identified with the map $h:M(p,q,r)\to S^3$,
which is an $r$-fold cyclic branched covering, 
branched over the $(p,q)$-torus link.
(See \cite[Lemma 1.1]{[Mi75]}.)

The 3-manifold $M(p,q,r)$ is a $\mathbb{S}^3$-manifold if and only if 
$p^{-1}+q^{-1}+r^{-1}>1$.
The triples $(2,2,r)$ give lens spaces.
The other triples with $p^{-1}+q^{-1}+r^{-1}>1$ 
are permutations of $(2,3,3)$, $(2,3,4)$ or $(2,3,5)$,
and give the three CK 3-manifolds with fundamental groups 
$Q(8)$, $T_1^*$ and $I^*$.
The manifolds $M(2,3,6)$, $M(3,3,3)$ and $M(2,4,4)$ are 
$\mathbb{N}il^3$-manifolds;
in all other cases $M(p,q,r)$ is a $\widetilde{\mathbb{SL}}$-manifold 
(in fact, a coset space of $\widetilde{SL}$ \cite{[Mi75]}), 
and $\sqrt{\pi_1(M(p,q,r))}\cong\mathbb{Z}$.
If $p$, $q$ and $r$ are pairwise relatively prime $M(p,q,r)$
is a $\mathbb{Z}$-homology 3-sphere.

\section{Rank 2 subgroups}

In this section we shall show that an abelian normal subgroup of rank 2 
in a 2-knot group is free abelian and not contained in the 
commutator subgroup.

\begin{lemma} 
Let $\nu$ be the fundamental group of a closed
$\mathbb{H}^2\times\mathbb{E}^1$-, $\mathbb{S}ol^3$- or 
$\mathbb{S}^2\times\mathbb{E}^1$-manifold. 
Then $\nu$ admits no meridianal automorphism.
\end{lemma}
                
\begin{proof}
The fundamental group of a closed $\mathbb{S}ol^3$- or 
$\mathbb{S}^2\times\mathbb{E}^1$-manifold has a 
characteristic subgroup with quotient having two ends.
If $\nu$ is a lattice in $Isom^+(\mathbb{H}^2\times\mathbb{E}^1)$ 
then $\sqrt\nu\cong\mathbb{Z}$,
and either $\sqrt\nu=\zeta\nu$ and is not contained in $\nu'$ 
or $C_\nu(\sqrt\nu)$ is a characteristic subgroup of index 2 in $\nu$.
In none of these cases can $\nu$ admit a meridianal automorphism.
\end{proof}

\begin{theorem}
Let $\pi$ be a finitely presentable group.
Then the following are equivalent:
\begin{enumerate}
\item$\pi$ is a $PD_4^+$-group of type $FF$ and weight $1$, 
and with an abelian normal subgroup $A$ of rank $2$;
\item$\pi=\pi{K}$ where $K$ is a $2$-knot and $\pi$ has
an abelian normal subgroup $A$ of rank $2$;
\item $\pi=\pi{K}$ where $K$ is a $2$-knot such that $M(K)$ 
is Seifert fibred over an aspherical $2$-orbifold $B$.
\end{enumerate}
If so, then $A\cong\mathbb{Z}^2$, $\pi'\cap A\cong\mathbb{Z}$,
$\pi'\cap A\leq\zeta\pi'\cap I(\pi')$, $[\pi:C_\pi(A)]\leq 2$ 
and $\pi'=\pi_1(N)$, where $N$ is a 
$\mathbb{N}il^3$- or $\widetilde{\mathbb{SL}}$-manifold.
If $\pi$ is solvable then $M(K)$ is
homeomorphic to a $\mathbb{N}il^3\times\mathbb{E}^1$-manifold.
If $\pi$ is not solvable then $M(K)$ 
is $s$-cobordant to a $\widetilde{\mathbb{SL}}\times\mathbb{E}^1$-manifold
which fibres over $S^1$ with fibre $N$, 
and closed monodromy of finite order.
\end{theorem}

\begin{proof}
If (1) holds then $\chi(\pi)=0$ \cite{[Ro84]}, 
and so $\beta_1(\pi)>0$, by Lemma 3.14.
Hence $H_1(\pi;\mathbb{Z})\cong\mathbb{Z}$,
since $\pi$ has weight 1. 
If $h(\sqrt\pi)=2$ then $A\cong\mathbb{Z}^2$, by Theorem 9.2.
Otherwise, $\pi$ is virtually poly-$Z$, and so $A\cong\mathbb{Z}^2$ again.
We may assume that $A$ is maximal among such subgroups.
Then $\pi/A\cong\beta=\pi_1^{orb}(B)$, for some aspherical 2-orbifold $B$,
and $\pi=\pi_1(M)$, where $M$ is an orientable 4-manifold
which is Seifert fibred over $B$, by Corollary 7.3.1.
Surgery on a loop in $M$ representing a normal generator of $\pi$
gives a 2-knot $K$ with $M(K)\cong{M}$ and $\pi{K}\cong\pi$.
Thus (2) and (3) hold. Conversely, these each imply (1).

Since no 2-orbifold group has infinite cyclic abelianization,
$A\not\leq\pi'$, and so $\pi'\cap{A}\cong\mathbb{Z}$.
If $\tau$ is the meridianal automorphism of $\pi'/I(\pi')$ then 
$\tau-1$ is invertible, 
and so $\tau$ cannot have $\pm1$ as an eigenvalue.
Hence $\pi'\cap A\leq I(\pi')$. 
(In particular, $\pi'$ is not abelian.)
The image of $\pi/C_\pi(A)$ in $Aut(A)\cong GL(2,\mathbb{Z})$ is triangular,
since $\pi'\cap A\cong\mathbb{Z}$ is normal in $\pi$.
Therefore as $\pi/C_\pi(A)$ has finite cyclic abelianization
it must have order at most 2.
Thus $[\pi:C_\pi(A)]\leq2$, so $\pi'<C_\pi(A)$ and $\pi'\cap A<\zeta\pi'$.
The subgroup $H$ generated by $\pi'\cup A$ has finite index in $\pi$ 
and so is also a $PD_4^+$-group.
Since $A$ is central in $H$ and maps onto $H/\pi'$ we have 
$H\cong\pi'\times\mathbb{Z}$.
Hence $\pi'$ is a $PD_3^+$-group with nontrivial centre.
As the nonabelian flat 3-manifold groups either admit no meridianal
automorphism or have trivial centre,
$\pi'=\pi_1(N)$ for some $\mathbb{N}il^3$- or 
$\widetilde{\mathbb{SL}}$-manifold $N$,
by Theorem 2.14 and Lemma 16.1.

The manifold $M(K)$ is $s$-cobordant to the mapping torus $M(\Theta)$ 
of a self homeomorphism of $N$, by Theorem 13.2.
If $N$ is a $\mathbb{N}il^3$-manifold $M(K)$ is homeomorphic to $M(\Theta)$,
by Theorem 8.1, 
and $M(K)$ must be a $\mathbb{N}il^3\times\mathbb{E}^1$-manifold, 
since $\mathbb{S}ol_1^4$-lattices do not have rank 2 abelian normal subgroups, 
while $\mathbb{N}il^4$-lattices cannot have abelianization $\mathbb{Z}$, 
as they have characteristic
rank 2 subgroups contained in their commutator subgroups.
Since $[\pi:C_\pi(A)]\leq2$ and $A\not\leq\pi'$ the meridianal
outer automorphism class has finite order.
Therefore if $N$ is a $\widetilde{\mathbb{SL}}$-manifold
then $M(\Theta)$ is a $\widetilde{\mathbb{SL}}\times\mathbb{E}^1$-manifold,
by Theorem 9.8.
\end{proof}

The space underlying the base $B$ is $S^2$, $RP^2$ or $D^2$,
since $\beta/\beta'$ is cyclic. 
Since $\pi'$ is infinite, so is $\beta$,
and $A$ is central if and only if $B$ is orientable.
If $B=S^2(a_1,\dots,a_r)$ then $r\geq3$,
and if $r=3$ then $\beta$ has weight 1
if and only if h.c.f.$(a_1,a_2,a_3)=1$.
In \cite{[HH13],[Hi20]}  the work of \cite{[Ho02]} is extended 
to show that if $r>3$ and $\beta$ has weight 1 then 
at most two disjoint pairs each have a common factor.
No known  examples have $r>4$.
If $B=RP^2(b_1,\dots,b_r)$ then $r=2$ or 3,
the $b_i$ are pairwise relatively prime, 
and the action $\alpha$ is via 
$\left(\begin{smallmatrix}
0&1\\ 1&0
\end{smallmatrix}\right)$. 
If $r=3$ then either $b_1=2$ or $\{b_1,b_2,b_3\}=\{3,4,5\}$.
When $r=2$ the group $\beta$ has weight 1, 
but this is undecided  when $r=3$.
If $B=\mathbb{D}(c_1,\dots,c_r,\overline{d_1},\dots,\overline{d_s})$
then $r\leq2$, $2r+s\geq3$, 
the $c_i$ are odd and relatively prime, at most one  $d_j$ is even
and $\alpha$ acts via 
$\left(\begin{smallmatrix}
1&0\\ 0&-1
\end{smallmatrix}\right).$
(However, $\pi^{orb}\mathbb{D}(3,\overline2)$ is finite.)
The weight of $\beta$ is 1 if $r\leq1$ or if $r=2$ and $s=0$, 
but is not known otherwise.

It is not hard to see that if $M$ is an orientable Seifert fibred 4-manifold 
then $\pi=\pi_1(M)$ has weight 1 if and only if $\pi/\pi'\cong\mathbb{Z}$ 
and $\pi/\sqrt\pi$ has weight 1.
Thus each of the bases allowed by the above considerations 
is realized by a 2-knot. 

\begin{theorem} 
Let $\pi$ be a $2$-knot group such that $\zeta\pi$ has rank $r>0$. 
If $\zeta\pi$ has nontrivial torsion it is finitely generated and $r=1$.
If $r>1$ then $\zeta\pi\cong\mathbb{Z}^2$,
$\zeta\pi'=\pi'\cap\zeta\pi\cong\mathbb{Z}$, and $\zeta\pi'\leq\pi''$.
\end{theorem}

\begin{proof} 
The first assertion follows from Theorem 15.10.
If $\zeta\pi$ had rank greater than 2 then $\pi'\cap\zeta\pi$ would
contain an abelian normal subgroup of rank 2, contrary to Theorem 16.2.
Therefore $\zeta\pi\cong\mathbb{Z}^2$ and $\pi'\cap\zeta\pi\cong\mathbb{Z}$.
Moreover $\pi'\cap\zeta\pi\leq\pi''$, since $\pi/\pi'\cong\mathbb{Z}$.
In particular $\pi'$ is nonabelian and $\pi''$ has nontrivial centre.
Hence $\pi'$ is the fundamental group of a $\mathbb{N}il^3$- or 
$\widetilde{\mathbb{SL}}$-manifold, by Theorem 16.2,
and so $\zeta\pi'\cong\mathbb{Z}$.
It follows easily that $\pi'\cap\zeta\pi=\zeta\pi'$.
\end{proof}

\section{Twist spins of torus knots}

The commutator subgroup of the group of the $r$-twist spin of a classical knot 
$K$ is the fundamental group of the $r$-fold cyclic branched cover of $S^3$, 
branched over $K$ \cite{[Ze65]}.
The $r$-fold cyclic branched cover of a sum of knots is the connected sum of 
the $r$-fold cyclic branched covers of the factors, 
and is irreducible if and only if the knot is prime.
Moreover the cyclic branched covers of a prime knot are either aspherical or 
finitely covered by $S^3$;
in particular no summand has free fundamental group \cite{[Pl84]}.
The cyclic branched covers of prime knots with nontrivial companions are 
Haken 3-manifolds \cite{[GL84]}.
The $r$-fold cyclic branched cover of the $(p,q)$-torus knot $k_{p,q}$ 
is the Brieskorn manifold $M(p,q,r)$ \cite{[Mi75]}.
The $r$-fold cyclic branched cover of a simple nontorus knot is a hyperbolic 
3-manifold if $r\geq3$, excepting only the 3-fold cyclic branched cover of the 
figure-eight knot, which is the Hanztsche-Wendt flat 3-manifold \cite{[Du83]}.
(In particular, there are only four $r$-fold cyclic branched covers of 
nontrivial knots for any $r>2$ which have finite fundamental group.)

\begin{theorem} 
Let $M$ be the $r$-fold cyclic branched cover 
of $S^3$, branched over a knot $K$, 
and suppose that $r>2$ and that $\sqrt{\pi_1(M)}$ is infinite.
Then $K$ is uniquely determined by $M$ and $r$, and either $K$ is 
a torus knot or $K\cong 4_1$ and $r=3$.
\end{theorem}
                                           
\begin{proof} 
As the connected summands of $M$ are the cyclic branched covers 
of the factors of $K$, any homotopy sphere summand must be standard,
by the proof of the Smith conjecture.
Therefore $M$ is aspherical, and is either Seifert fibred or is a
$\mathbb{S}ol^3$-manifold, by Theorem 2.14.
It must in fact be a $\mathbb{E}^3$-, $\mathbb{N}il^3$- or 
$\widetilde{\mathbb{SL}}$-manifold, by Lemma 16.1.
If there is a Seifert fibration which is preserved by the automorphisms of the 
branched cover the fixed circle (the branch set of $M$) must be a fibre of the 
fibration (since $r>2$) which therefore passes to a Seifert fibration of $X(K)$.
Thus $K$ must be a $(p,q)$-torus knot, for some relatively prime integers
$p$ and $q$ \cite{[BZ]}.
These integers may be determined arithmetically from $r$ and the formulae 
for the Seifert invariants of $M(p,q,r)$ given in \S1.
Otherwise $M$ is flat \cite{[MS86]} and so $K\cong 4_1$ and $r=3$, 
by \cite{[Du83]}.
\end{proof}

All the knots whose 2-fold branched covers are Seifert fibred 
are torus knots or Montesinos knots. 
(This class includes the 2-bridge knots and pretzel knots, 
and was first described in \cite{[Mo73]}.) 
The number of distinct knots whose 2-fold branched cover is a given Seifert 
fibred 3-manifold can be arbitrarily large \cite{[Be84]}. 
Moreover for each $r\geq2$ there are distinct simple 
1-knots whose $r$-fold cyclic branched covers are homeomorphic 
\cite{[Sa81],[Ko86]}.

If $K$ is a fibred 2-knot with monodromy of finite order $r$ and if $(r,s)=1$
then the $s$-fold cyclic branched cover of $S^4$, branched over $K$ is again a 
4-sphere and so the branch set gives a new 2-knot, 
which we shall call the $s$-fold cyclic branched cover of $K$.
This new knot is again fibred, with the same fibre and monodromy the $s^{th}$ 
power of that of $K$ \cite{[Pa78],[Pl86]}. 
If $K$ is a classical knot we shall let $\tau_{r,s}K$ denote the
$s$-fold cyclic branched cover of the $r$-twist spin of $K$.
We shall call such knots {\it branched twist spins}, for brevity.

Using properties of $S^1$-actions on smooth homotopy 4-spheres, 
Plotnick obtains the following result \cite{[Pl86]}.                      

\begin{thm}
[Plotnick] 
A $2$-knot is fibred with periodic monodromy if and only if it is 
(up to Gluck reconstruction) a branched twist spin of a knot 
in a homotopy $3$-sphere.
\qed
\end{thm}

Here ``periodic monodromy" means that the fibration of the exterior 
of the knot has a characteristic map of finite order. 
It is not in general sufficient that the 
{\it closed} monodromy be represented by a map of finite order.
(For instance, if $K$ is a fibred 2-knot with $\pi'\cong Q(8)\times (Z/nZ)$ 
for some $n>1$ then the meridianal automorphism of $\pi'$ has order 6,
and so it follows from the observations above that $K$ is not a twist spin.)

The homotopy 3-sphere must be standard, by Perelman's work (see \cite{[B-P]}).
In our application in the next theorem we are able to show this directly.

\begin{theorem} 
A group $G$ which is not virtually solvable is 
the group of a branched twist spin of a torus knot if and only if it is a 
$3$-knot group and a $PD_4^+$-group with centre of rank $2$, 
some nonzero power of a weight element being central.
\end{theorem}

\begin{proof}              
If $K$ is a cyclic branched cover of $\tau_rk_{p,q}$
then $M(K)$ fibres over $S^1$ with fibre $M(p,q,r)$ and monodromy 
of order $r$,
and so the $r^{th}$ power of a meridian is central.
Moreover the monodromy commutes with the natural $S^1$-action on $M(p,q,r)$
(see \cite[Lemma 1.1]{[Mi75]}) and hence preserves a Seifert fibration. 
Hence the meridian commutes with $\zeta\pi_1(M(p,q,r))$, 
which is therefore also central in $G$.
Since $\pi_1(M(p,q,r))$ is either a $PD_3^+$-group with infinite centre 
or is finite, 
the necessity of the conditions is evident.

Conversely, if $G$ is such a group then $G'$ is the fundamental group of a 
Seifert fibred 3-manifold, $N$ say, by Theorem 2.14.
Moreover $N$ is ``sufficiently complicated" in the sense of \cite{[Zn79]}, 
since $G'$ is not virtually solvable.
Let $t$ be an element of $G$ whose normal closure is the whole group,
and such that $t^n$ is central for some $n>0$.
Let $\theta$ be the automorphism of $G'$ determined by $t$,
and let $m$ be the order of the outer automorphism class $[\theta]\in Out(G')$.
By \cite[Corollary 3.3]{[Zn79]} there is a fibre preserving 
self homeomorphism $\tau$ of $N$ inducing $[\theta]$ 
such that the group of homeomorphisms 
of $\widetilde N\cong R^3$ generated by the covering group
$G'$ together with the lifts of $\tau$ is an extension of $Z/mZ$ by $G'$,
and which is a quotient of the semidirect product 
$\hat G=G/\langle\langle t^n\rangle\rangle\cong G'\rtimes_\theta(Z/nZ)$.
Since the self homeomorphism of $\widetilde N$ corresponding to the image of 
$t$ has finite order it has a connected 1-dimensional fixed point set,
by Smith theory. 
The image $P$ of a fixed point in $N$ determines a cross-section 
$\gamma=\{P\}\times S^1$ of the mapping torus $M(\tau)$.
Surgery on $\gamma$ in $M(\tau)$ gives a 2-knot with group $G$ which is fibred 
with monodromy (of the fibration of the exterior $X$) of finite order.
We may then apply Plotnick's Theorem to conclude that the 2-knot
is a branched twist spin of a knot in a homotopy 3-sphere.
Since the monodromy respects the Seifert fibration and leaves the centre of 
$G'$ invariant, the branch set must be a fibre, 
and the orbit manifold a Seifert fibred homotopy 3-sphere.
Therefore the orbit knot is a torus knot in $S^3$, 
and the 2-knot is a branched twist spin of a torus knot.
\end{proof}

If $K$ is a 2-knot with group as in Theorem 16.5 then $M(K)$ is 
aspherical, and so is homotopy equivalent to some such knot manifold.
(See also Theorem 16.2.)

Can we avoid the appeal to Plotnick's Theorem in the above argument?
There is a stronger result for fibred 2-knots.
The next theorem is a version of \cite[Proposition 6.1]{[Pl86]}, 
starting from more algebraic hypotheses.

\begin{theorem} 
Let $K$ be a fibred $2$-knot such that $\pi{K}$ has centre of rank $2$, 
some power of a weight element being central.
Then $M(K)$ is homeomorphic to $M(K_1)$, where
$K_1$ is some branched twist spin of a torus knot.
\end{theorem}
                      
\begin{proof} 
Let $F$ be the closed fibre and $\phi:F\to F$ the characteristic map.
Then $F$ is a Seifert fibred manifold, as above.
Now the automorphism of $F$ constructed as in Theorem 16.5 induces the same
outer automorphism of $\pi_1(F)$ as $\phi$, and so these maps must be homotopic.
Therefore they are in fact isotopic \cite{[Sc85],[BO91]}. 
The theorem now follows.
\end{proof}

If $(p,q)=(r,s)=1$ then $M(\tau_{r,s}k_{p,q})$ 
is Seifert fibred over $S^2(p,q,r)$. 
(This can be derived from \S16.1.)
More generally, we have the following lemma.

\begin{lemma}
Let $B$ be an aspherical orientable $2$-orbifold. 
If $\beta=\pi_1^{orb}(B)$ is normally generated by an element 
of finite order $r$ then $B=S^2(m,n,r)$, where $(m,n)=1$.
\end{lemma}

\begin{proof}
Since $\beta$ is infinite and $\beta/\beta'$ is cyclic, $B=S^2(\alpha_1,\dots,\alpha_k)$, 
for some $\alpha_i>1$ and $k\geq3$.
Hence $\beta$ has the presentation
\[
\langle
{x_1,\dots,x_k}\mid {x_i^{\alpha_i}=1}~\forall~1\leq{i}\leq{k},~\Pi{x_i}=1
\rangle.
\]
Every such group is nontrivial, 
and elements of finite order are conjugate to powers of
a cone point generator $x_i$ \cite[Theorem 4.8.1]{[ZVC]}.
If $x_k^s$ is a weight element then so is $x_k$, 
and so
$\langle{x_1,\dots,x_{k-1}}\mid 
{x_i^{\alpha_i}=1}~\forall~1\leq{i}\leq{k-1},
~\Pi{x_i}=1\rangle$ is trivial.
Hence $k-1=2$ and $(\alpha_1,\alpha_2)=1$.
Similarly, $(s,\alpha_3)=1$, so $r=\alpha_3$.
Setting $m=\alpha_1$ and $n=\alpha_2$, 
we see that $B=S^2(m,n,r)$ and $(m,n)=1$.
\end{proof}

The central extension of $\pi_1^{orb}(S^2(6,10,15))$ by $\mathbb{Z}^2$
with presentation
\[
\langle{u,v,x,y,z}\mid {u,v}~central,~x^6=u^5,~y^{10}=uv,~z^{15}=uv^{-1},~xyz=u\rangle
\]
is torsion-free, has abelianization $\mathbb{Z}$,
and is normally generated by $yz^{-1}$.
The corresponding Seifert 4-manifold with base $S^2(6,10,15)$
is a knot manifold.

The manifold obtained by 0-framed surgery on the reef knot $3_1\#-3_1$
is the Seifert fibred 3-manifold $N=M(0;(2,1),(2,-1),(3,1),(3,-1))$,
and the total space of the $S^1$-bundle 
over $N$ with Euler class a generator of $H^2(N;\mathbb{Z})$ 
is a knot manifold which is Seifert fibred over $S^2(2,2,3,3)$.
(Using $k_{p,q}$ instead of $3_1$ gives a knot manifold with
Seifert base $S^2(p,p,q,q)$.)

In each of these cases, it follows from Lemma 16.7 that
no power of any weight element is central.
Thus no such knot is a branched twist spin, and 
the final hypothesis of Theorem 16.5
is not a consequence of the others.

If $p$, $q$ and $r$ are pairwise relatively prime then 
$M(p,q,r)$ is an homology sphere and $\pi=\pi\tau_rk_{p,q}$
has a central element which maps to a generator of $\pi/\pi'$.
Hence $\pi\cong\pi'\times\mathbb{Z}$ and $\pi'$ has weight 1.
Moreover if $t$ is a generator for the $\mathbb{Z}$ summand 
then an element $h$ of $\pi'$ is a weight element for $\pi'$ 
if and only if $ht$ is a weight element for $\pi$.
This correspondance also gives a bijection between conjugacy classes of such 
weight elements.
If we exclude the case $(2,3,5)$ then $\pi'$ has infinitely many 
distinct weight orbits,
and moreover there are weight elements such that no power is central 
\cite{[Pl83]}.
Therefore we may obtain many 2-knots whose groups are as in Theorem 16.5 
but which are not themselves branched twist spins by surgery on weight 
elements in $M(p,q,r)\times S^1$.

We may apply arguments similar to those of Theorems 16.5 and 16.6
in attempting to understand twist spins of other knots.
As only the existence of homeomorphisms of finite order and 
``homotopy implies isotopy" require different justifications, 
we shall not give proofs for the following assertions.

Let $G$ be a 3-knot group such that $G'$ is the fundamental group of an
aspherical 3-manifold $N$ and in which some nonzero power of a weight element 
is central.
Then $N$ is Seifert fibred, hyperbolic or Haken, 
by Perelman's work \cite{[B-P]}.
If $N$ is hyperbolic we may use Mostow rigidity to show that $G$ is the group 
of some branched twist spin $K$ of a simple non-torus knot.
Moreover, if $K_1$ is another fibred 2-knot with group $G$
and hyperbolic fibre then $M(K_1)$ is homeomorphic to $M(K)$.
In particular the simple knot and the order of the twist are
determined by $G$.
Similarly if $N$ is Haken, but neither hyperbolic nor Seifert fibred, 
then we may use \cite{[Zn82]} to show that $G$ is the group 
of some branched twist spin of a prime non-torus knot.
Moreover, 
the prime knot and the order of the twist are determined by $G$ 
\cite{[Zn86]},
since all finite group actions on $N$ are geometric \cite{[B-P]}.

\section{Solvable $PD_4$-groups}

If $K$ is a 2-knot such that $h(\sqrt{\pi{K}})>2$ then $\pi{K}$ 
is virtually poly-$Z$ and $h(\pi)=4$, by Theorem 8.1.
In this section we shall determine all such groups and their weight orbits.

\begin{lemma} 
Let $G$ be torsion-free and virtually poly-$Z$ 
with $h(G)=4$ and $G/G'\cong\mathbb{Z}$. 
Then $G'\cong\mathbb{Z}^3$ or $G_6$
or $\sqrt{G'}\cong\Gamma_q$ (for some $q>0$) and $G'/\sqrt{G'}\cong Z/3Z$ 
or $1$.
\end{lemma}

\begin{proof} Let $H=G/\sqrt{G'}$. 
Then $H/H'\cong\mathbb{Z}$ and $h(H')\leq1$, since $\sqrt{G'}=G'\cap\sqrt{G}$ 
and $h(G'\cap\sqrt{G})\geq h(\sqrt{G})-1\geq2$, by Theorem 1.4.
Now $H'$ cannot have two ends, since $H/H'\cong\mathbb{Z}$,
and so $H'=G'/\sqrt{G'}$ is finite.

If $\sqrt{G'}\cong\mathbb{Z}^3$ then $G'\cong\mathbb{Z}^3$ or $G_6$, 
since these are the only flat 3-manifold groups which admit 
meridianal automorphisms.

If $\sqrt{G'}\cong\Gamma_q$ for some $q>0$ then 
$\zeta\sqrt{G'}\cong\mathbb{Z}$ 
is normal in $G$ and so is central in $G'$.
Using the known structure of automorphisms of $\Gamma_q$, it follows that
the finite group $G'/\sqrt{G'}$ must act on 
$\sqrt{G'}/\zeta\sqrt{G'}\cong\mathbb{Z}^2$ 
via $SL(2,\mathbb{Z})$ and so must be cyclic. 
Moreover it must have odd order, and hence be 1 or $Z/3Z$, 
since ${G/\sqrt{G'}}$ has infinite cyclic abelianization.
\end{proof}

There is a fibred 2-knot $K$ with $\pi{K}\cong{G}$ if and only if $G$
is orientable, by Theorems 14.4 and 14.7.

\begin{theorem} 
Let $\pi$ be a $2$-knot group with $\pi'\cong\mathbb{Z}^3$,
and let $\,C$ be the image of the meridianal automorphism in $SL(3,\mathbb{Z})$.
Then $\Delta_C(t)=\det(tI-C)$ is irreducible, $|\Delta_C(1)|=1$
and $\pi'$ is isomorphic to an ideal in the domain $R=\Lambda/(\Delta_C(t))$.
Two such groups are isomorphic if and only if the polynomials are equal 
(after inverting $t$, if necessary) and the ideal classes then agree.
There are finitely many ideal classes for each such polynomial and each class 
(equivalently, each such matrix) is realized by some $2$-knot group.
Moreover $\sqrt\pi=\pi'$ and $\zeta\pi=1$.
Each such group $\pi$ has two strict weight orbits.
\end{theorem}
                      
\begin{proof} 
Let $t$ be a weight element for $\pi$ and let $C$ be the matrix 
of the action of $t$ by conjugation on $\pi'$, with respect to some basis.
Then $\det(C-I)=\pm1$, since $t-1$ acts invertibly.
Moreover if $K$ is a 2-knot with group $\pi$ then
$M(K)$ is orientable and aspherical, so $\det(C)=+1$. 
Conversely, surgery on the mapping torus of the self homeomorphism of 
$S^1\times S^1\times S^1$
determined by such a matrix $C$ gives a 2-knot with group 
$\mathbb{Z}^3\rtimes_C\mathbb{Z}$.

The Alexander polynomial of $K$ is the characteristic polynomial 
$\Delta_K(t)=\det(tI-C)$ which has the form $t^3-at^2 +bt-1$, 
for some $a$ and $b=a\pm1$.
It is irreducible, since it does not vanish at $\pm1$. 
Since $\pi'$ is annihilated by $\Delta_K(t)$ it is an $R$-module;
moreover as it is torsion-free it embeds in $\mathbb{Q}\otimes\pi'$,
which is a vector space over the field of fractions $\mathbb{Q}\otimes R$.
Since $\pi'$ is finitely generated and $\pi'$ and $R$ each have rank 3 
as abelian groups it follows that $\pi'$ is isomorphic to an
ideal in $R$.
Moreover the characteristic polynomial of $C$ cannot be cyclotomic 
and so no power of $t$ can commute with any nontrivial element of $\pi'$.
Hence $\sqrt\pi=\pi'$ and $\zeta\pi=1$.

By Lemma 1.1 two such semidirect products are isomorphic if and only if
the matrices are conjugate up to inversion.
The conjugacy classes of matrices in $SL(3,\mathbb{Z})$ with given 
irreducible characteristic polynomial $\Delta(t)$
correspond to the ideal classes of 
$\Lambda/(\Delta(t))$, by Theorem 1.4.
Therefore $\pi$ is determined by the ideal class of $\pi'$,
and there are finitely many such 2-knot groups with given
Alexander polynomial.

Since $\pi''=1$ there are just two strict weight orbits,
by Theorem 14.1.
\end{proof}

We shall call 2-knots with such groups {\it Cappell-Shaneson} 2-knots.

In \cite{[AR84]} matrix calculations are used to show that any matrix $C$ 
as in Theorem 16.9 is conjugate to one with first row $(0,0,1)$.
Given this, 
it is easy to see that the corresponding Cappell-Shaneson 2-knot
group has a presentation
\[
\langle t,x,y,z\mid xy=yx,\medspace xz=zx,\medspace txt^{-1}=z,
\medspace tyt^{-1}=x^my^nz^p\! ,\medspace tzt^{-1}=x^qy^rz^s\rangle.
\]
Since $p$ and $s$ must be relatively prime these relations 
imply $yz=zy$. 
We may reduce the number of generators and relations on setting $z=txt^{-1}$.
The next lemma gives a more conceptual exposition of part of this result.

\begin{lemma} 
Let $\Delta_a(t)=t^3-at^2+(a-1)t-1$ for some $a\in\mathbb{Z}$, 
and let $M$ be an ideal in the domain $R=\Lambda/(\Delta_a(t))$.
Then $M$ can be generated by $2$ elements as an $R$-module. 
\end{lemma}

\begin{proof} 
Let $D=a(a-2)(a-3)(a-5)-23$ be the discriminant of $\Delta_a(t)$.
If a prime $p$ does not divide $D$ then 
$\Delta_a(t)$ has no repeated roots {\it mod} $(p)$.
If $p$ divides $D$ then $p>5$, 
and there are integers $\alpha_p$, $\beta_p$ 
such that $\Delta_a(t)\equiv$ $(t-\alpha_p)^2(t-\beta_p)$ {\it mod} $(p)$.
Let $K_p=\{ m\in M\mid (t-\beta_p)m\in pM\}$.

If $\beta_p\equiv\alpha_p$ {\it mod} $(p)$ then
$\alpha_p^3\equiv-\Delta_a(0)=1$ and $(1-\alpha_p)^3\equiv\Delta_a(1)=-1$ 
{\it mod} $(p)$.
Together these congruences imply that $p=7$ and
$\alpha_p\equiv2$ {\it mod} $(p)$.
Fix a basis for $M\cong\mathbb{Z}^3$, and let
$A$ be the matrix of multiplication by $t$ on $M$.
If $M/7M\cong(R/(7,t-2))^3$ then $A=2I+7B$ for some $\mathbb{Z}$-matrix $B$.
On expanding $\det(A)=\det(2I+7B)$ and $\det(A-I)=\det(I+7B)$ {\it mod} 
$(49)$,
we find that $\det(A)=1$ and $\det(A-I_3)=\pm1$ cannot both hold.
(Compare Lemma A2 of \cite{[AR84]}.)
Thus $M/7M\cong{R/(7,(t-2)^3)}$ or $R/(7,(t-2)^2)\oplus{R/(7,t-2)}$,
and $K_7$ has index at least 7 in $M$. 
For all primes $p\not=7$ dividing $D$ the index of $K_p$ is $p^2$.
Since 
\begin{equation*}
\frac17+\Sigma_{p|D,p\not=7} \frac1{p^2}<\frac17+\int_2^\infty\frac1{t^2}dt<1,
\end{equation*}
the image of $\cup_{p|D}K_p$ in $M/\cap_{p|D}K_p$ is a proper subset,
and so $M-\cup_{p|D}K_p$ is nonempty.
Let $x$ be an element of $M-\cup_{p|D}K_p$ 
which is not $\mathbb{Z}$-divisible in $M$. 
Then $N=M/Rx$ is finite, 
and is generated by at most two elements as an abelian group, 
since $M\cong\mathbb{Z}^3$.
For each prime $p$ the $\Lambda/p\Lambda$-module $M/pM$ is an extension of 
$N/pN$ by the submodule $X_p$ generated by the image of $x$ 
and its order ideal is generated by the image of $\Delta_a(t)$ in the P.I.D. 
$\Lambda/p\Lambda\cong\mathbb{F}_p[t,t^{-1}]$.
If $p$ does not divide $D$ the image of $\Delta_a(t)$ in $\Lambda/p\Lambda$ is
square free.
If $p$ divides $D$ the order ideal of $X_p$ is divisible by $t-\alpha_p$,
since $(t-\beta_p)x\not\in{pM}$.
In all cases the order ideal of $N/pN$ is square free and so $N/pN$ is cyclic
as a module.
By the Chinese Remainder Theorem there is an element $y\in M$ whose image is
a generator of $N/pN$, for each prime $p$ dividing the order of $N$.
The image of $y$ in $N$ generates $N$ as a module, by Nakayama's Lemma.
\end{proof}

The cited result of \cite{[AR84]} is equivalent to showing that 
$M$ has an element $x$ such that $M/Rx$ is cyclic as an abelian group,
from which it follows that $M$ is generated as an abelian group by $x$, 
$tx$ and some third element $y$.

\begin{lemma} 
Let $\pi$ be a finitely presentable group such that $\pi/\pi'\cong\mathbb{Z}$,
and let $R=\Lambda$ or $\Lambda/p\Lambda$ for some prime $p\geq2$. Then
\begin{enumerate}
\item if $\pi$ can be generated by $2$ elements $H_1(\pi;R)$ is cyclic as an 
$R$-module;

\item if $\mathrm{def}(\pi)=0$ then $H_2(\pi;R)$ is cyclic as an $R$-module.
\end{enumerate}
\end{lemma}

\begin{proof} 
If $\pi$ is generated by $t$ and $x$, say, 
we may assume that the image of $t$ generates $\pi/\pi'$ and that $x\in\pi'$.
Then $\pi'$ is generated by the conjugates $\{t^kxt^{-k}|k\in\mathbb{Z}\}$,
and $H_1(\pi;R)=R\otimes_\Lambda(\pi'/\pi'')$ is generated by the image of 
$x$.

Let $X$ be the 2-complex determined by a deficiency 0 presentation for $\pi$. 
Then $H_0(X;R)=R/(t-1)$ and $H_1(X;R)$ are $R$-torsion modules, 
and $H_2(X;R)$ is the kernel of a map between finitely generated 
free $R$-modules.
Hence $H_2(X;R)\cong{R}$, as it has rank 1 and $R$ is a noetherian UFD.
Since $H_2(\pi;R)$ is a quotient of $H_2(X;R)$,
by Hopf's Theorem,
it is cyclic as an $R$-module.
\end{proof}

If $M(K)$ is aspherical $H_2(\pi;R)\cong\overline{Ext^1_R(H_1(\pi;R),R)}$,
by Poincar\'e duality and the UCSS. 
In particular, if $R=\Lambda/p\Lambda$ for some prime $p\geq2$ then
$H_1(\pi;R)$ and $\overline{H_2(\pi;R)}$ are non-canonically isomorphic.

\begin{theorem} 
Let $\pi=\mathbb{Z}^3\rtimes_C\mathbb{Z}$ be the group of a 
Cappell-Shaneson $2$-knot, and let $\Delta(t)=\det(tI-C)$. 
Then $\pi$ has a $3$ generator presentation of deficiency $-2$.
Moreover the following are equivalent.
\begin{enumerate}
\item $\pi$ has a $2$ generator presentation of deficiency $0$;

\item $\pi$ is generated by $2$ elements;

\item $\mathrm{def}(\pi)=0$;

\item $\pi'$ is cyclic as a $\Lambda$-module.
\end{enumerate}
\end{theorem}
                                                              
\begin{proof} 
The first assertion follows immediately from Lemma 16.10.
Condition (1) implies (2) and (3), since $\mathrm{def}(\pi)\leq0$, 
by Theorem 2.5, while (2) implies (4), by Lemma 16.11.
If $\mathrm{def}(\pi)=0$ then $H_2(\pi;\Lambda)$ is cyclic as a $\Lambda$-module, 
by Lemma 16.11.
Since $\pi'=H_1(\pi;\Lambda)\cong\overline{H^3(\pi;\Lambda)}
\cong\overline{Ext^1_\Lambda(H_2(\pi;\Lambda),\Lambda)}$, 
by Poincar\'e duality and the UCSS, 
it is also cyclic and so (3) also implies (4).
If $x$ generates $\pi'$ as a $\Lambda$-module it is easy to see 
that $\pi$ has a presentation 
\begin{equation*}
\langle t,x\mid xtxt^{-1}=txt^{-1}x,\medspace 
t^3xt^{-3}=t^2x^at^{-2}tx^bt^{-1}x\rangle,
\end{equation*}
for some integers $a,b$, and so (1) holds.
\end{proof}
In fact Theorem A.3 of \cite{[AR84]} implies that any such group 
has a 3 generator presentation of deficiency -1, 
as remarked before Lemma 16.10.
Since $\Delta(t)$ is irreducible the $\Lambda$-module $\pi'$ is determined  
by the Steinitz-Fox-Smythe row ideal $(t-n,m+np)$ in the domain 
$\Lambda/(\Delta(t))$. 
(See \cite[Chapter 3]{[H3]}.)
Thus $\pi'$ is cyclic if and only if this ideal is principal. 
In particular, this is not so for the concluding example of \cite{[AR84]},
which gives rise to the group with presentation
\[
\langle t,x,y,z\mid xz=zx,\medspace yz=zy,\medspace 
txt^{-1}=y^{-5}z^{-8}\! ,\medspace tyt^{-1}=y^2z^3\! ,\medspace 
tzt^{-1}=xz^{-7}\rangle.
\]
Since $(\Delta_a(t))\not=(\Delta_a(t^{-1}))$ 
(as ideals of $\mathbb{Q}\Lambda$), for all $a$,
no Cappell-Shaneson 2-knot is a deform spin,
by the criterion of \cite{[BM09]}.

Let $G(+)$ and $G(-)$ be the extensions of $\mathbb{Z}$ by $G_6$ 
with presentations
\begin{equation*}
\langle t,x,y\mid xy^2x^{-1}y^2=1, txt^{-1}=(xy)^{\mp 1}, tyt^{-1}=x^{\pm 1} 
\rangle.
\end{equation*}
In each case, using the final relation to eliminate the generator $x$ gives
a 2-generator presentation of deficiency 0, which is optimal, by Theorem 2.5.
         
\begin{theorem} 
Let $\pi$ be a $2$-knot group with $\pi'\cong G_6$.
Then $\pi\cong G(+)$ or $G(-)$.
In each case $\pi$ is virtually $\mathbb{Z}^4$, 
$\pi'\cap\zeta\pi=1$ and $\zeta\pi\cong\mathbb{Z}$.
Every strict weight orbit representing a given generator $t$ 
for $\pi/\pi'$ contains an unique element of the form $x^{2n}t$, 
and every such element is a weight element.
\end{theorem}
                      
\begin{proof} Since $Out(G_6)$ is finite $\pi$ is virtually 
$G_6\times\mathbb{Z}$ and hence is virtually $\mathbb{Z}^4$. 
The groups $G(+)$ and $G(-)$ are the only orientable flat
4-manifold groups with $\pi/\pi'\cong\mathbb{Z}$.
The next assertion ($\pi'\cap\zeta\pi=1$) follows as $\zeta G_6=1$.
It is easily seen that $\zeta G(+)$ and $\zeta G(-)$ are generated by 
the images of $t^3$ and $t^6x^{-2}y^2(xy)^{-2}$, respectively,
and so in each case $\zeta\pi\cong\mathbb{Z}$.

If $t\in\pi$ represents a generator of $\pi/\pi'$
it is a weight element, since $\pi$ is solvable.
We shall use the notation of \S2 of Chapter 8 for automorphisms of $G_6$.

Suppose first that $\pi=G(+)$ and $c_t=ja$.
If $u$ is another weight element with the same image in $\pi/\pi'$
then we may assume that $u=gt$ for some $g\in\pi''=G_6'$,
by Theorem 14.1.
Suppose that $g=x^{2m}y^{2n}z^{2p}$.
Let $\lambda(g)=m+n-p$ and $w=x^{2n}y^{2p}$.
Then $w^{-1}gtw=x^{2\lambda(g)}t$.
On the other hand, if $\psi\in{Aut(G_6)}$ then 
$\psi{c_{gt}}\psi^{-1}=c_{ht}$ for some $h\in{G_6'}$ 
if and only if the images of $\psi$ and $ja$ in $Aut(G_6)/G_6'$ commute.
If so, $\psi$ is in the subgroup generated by $\{def^{-1},jb,ce\}$.
It is easily verified that $\lambda(h)=\lambda(g)$
for any such $\psi$. 
(It suffices to check this for the generators of $Aut(G_6)$.)
Thus $x^{2n}t$ is a weight element representing $[ja]$,
for all $n\in\mathbb{Z}$, and $x^{2m}t$ and $x^{2n}t$ are in 
the same strict weight orbit if and only if $m=n$.

The proof is very similar when $\pi=G(-)$ and $c_t=jb$.
The main change is that we should define the homomorphism $\lambda$
by $\lambda(x^{2m}y^{2n}z^{2p})=m-n+p$. 
\end{proof}

The group $G(+)$ is the group of the 3-twist spin of the figure eight knot
($G(+)\cong\pi\tau_3 4_1$.)
Although $G(-)$ is the group of a fibred 2-knot, by Theorem 14.4, 
it is not the group of any twist spin, by Theorem 16.4.
We can also see this directly.

\begin{cor}
No $2$-knot with group $G(-)$ is a twist-spin.
\end{cor}

\begin{proof}
If $G(-)=\pi\tau_rK$ for some 1-knot $K$ then
the $r^{th}$ power of a meridian is central in $G(-)$.
The power $(x^{2n}t)^r$ is central if and only if $(d^{2n}jb)^r=1$
in $Aut(G_6)$.
But $(d^{2n}jb)^3=d^{2n}f^{2n}e^{-2n}(jb)^3=(de^{-1}f)^{2n+1}$.
Therefore $d^{2n}jb$ has infinite order,
and so $G(-)$ is not the group of a twist-spin.
\end{proof}

\begin{theorem} 
Let $\pi$ be a $2$-knot group with 
$\pi'\cong\Gamma_q$ for some $q>0$, and let $\theta$ be the image of the
meridianal automorphism in $Out(\Gamma_q)$.
Then either $q=1$ and $\theta$ is conjugate to
$[\left(\smallmatrix 1&-1\\
1&0\endsmallmatrix\right),0]$ or
$[\left(\smallmatrix 1&1\\
1&2\endsmallmatrix\right),0]$,
or $q$ is odd and $\theta$ is conjugate to
$[\left(\smallmatrix 1&1\\
1&0\endsmallmatrix\right),0]$ or its inverse.
Each such group $\pi$ has two strict weight orbits.
\end{theorem}

\begin{proof} If $(A,\mu)$ is a meridianal automorphism of $\Gamma_q$
the induced automorphisms of $\Gamma_q/\zeta\Gamma_q\cong Z^2$ and 
$tors(\Gamma_q/\Gamma_q')\cong Z/qZ$ are also meridianal, and so
$\det(A-I)=\pm1$ and $\det(A)-1$ is a unit {\it mod} $(q)$.
Therefore $q$ must be odd and $\det(A)=-1$ if $q>1$, and
the characteristic polynomial $\Delta_A(X)$ of $A$ must be 
$X^2-X+1$, $X^2-3X+1$, $X^2-X-1$ or $X^2+X-1$. 
Since the corresponding rings $\mathbb{Z}[X]/(\Delta_A(X))$ 
are isomorphic to $\mathbb{Z}[(1+\sqrt{-3})/2]$ or $\mathbb{Z}[(1+\sqrt5)/2]$),
which are PIDs, $A$ is conjugate to one of
$\left(\smallmatrix 1&-1\\
1&0\endsmallmatrix\right)$, 
$\left(\smallmatrix 1&1\\
1&2\endsmallmatrix\right)$,
$\left(\smallmatrix 1&1\\
1&0\endsmallmatrix\right)$,
or
$\left(\smallmatrix 1&1\\
1&0\endsmallmatrix\right)^{-1}=
\left(\smallmatrix 0&1\\
1&-1\endsmallmatrix\right)$,
by Theorem 1.4.
Now 
$[A,\mu][A,0][A,\mu]^{-1}=[A,\mu(I-\det(A)A)^{-1}]$
in $Out(\Gamma_q)$.
(See \S7 of Chapter 8.)
As in each case $I-\det(A)A$ is invertible, it follows that
$\theta$ is conjugate to $[A,0]$ or to $[A^{-1},0]=[A,0]^{-1}$.
Since $\pi''\leq\zeta\pi'$ the final observation follows from Theorem 14.1.
\end{proof}

The groups $\Gamma_q$ are discrete cocompact subgroups of the Lie group 
$Nil^3$ and the coset spaces are $S^1$-bundles over the torus.
Every automorphism of $\Gamma_q$ is orientation preserving and each of 
the groups allowed by Theorem 16.14 is the group of some fibred 2-knot, 
by Theorem 14.4.
The automorphism $[\left(\smallmatrix 1&-1\\1&0\endsmallmatrix\right), 0]$
is realised by $\tau_63_1$,
and $M(\tau_63_1)$ is Seifert fibred over $S^2(2,3,6)$.
In all the other cases $\theta$ has infinite order 
and the group is not the group of any twist spin. 
If $q>1$ no such 2-knot is deform spun,
by the criterion of \cite{[BM09]},
since $(t^2-t-1)\not=(t^{-2}-t^{-1}-1)$ (as ideals of $\mathbb{Q}\Lambda$).

The 2-knot groups with commutator subgroup $\Gamma_1$ have presentations
\begin{equation*}
\langle t,x,y\mid xyxy^{-1}=yxy^{-1}x,\medspace txt^{-1}=xy,\medspace 
tyt^{-1}=w \rangle,
\end{equation*}
where $w=x^{-1}$, $xy^2$ or $x$ (respectively), while those with
commutator subgroup $\Gamma_q$ with $q>1$ have presentations
\begin{equation*}
\langle t,u,v,z\mid uvu^{-1}v^{-1}=z^q\! ,\medspace tut^{-1}=v,\medspace
tvt^{-1}=zuv,\medspace tzt^{-1}=z^{-1} \rangle.
\end{equation*}
(Note that as $[v,u]=t[u,v]t^{-1}=[v,zuv]=[v,z]z[v,u]z^{-1}=[v,z][v,u]$, 
we have $vz=zv$ and hence $uz=zu$ also.)
These are easily seen to have 2 generator 
presentations of deficiency 0 also.                                             
                         
The other $\mathbb{N}il^3$-manifolds which arise as the closed fibres 
of fibred 2-knots are Seifert fibred over $S^2(3,3,3)$. 
(In all other cases the fundamental group has no meridianal automorphism.) 
These are 2-fold branched covers of $(S^3,k(e,\eta))$,
where $k(e,\eta)=\mathfrak{m}(e;(3,1),(3,1),(3,\eta))$ is a Montesinos link,
for some $e\in\mathbb{Z}$ and $\eta=\pm1$ \cite{[Mo73],[BZ]}.
If $e$ is even this link is a knot, and is
invertible, but not amphicheiral. (See \cite[\S 12E]{[BZ]}.)
Let $\pi(e,\eta)=\pi\tau_2k(e,\eta)$.

\begin{theorem} 
Let $\pi$ be a $2$-knot group such that $\sqrt{\pi'}\cong \Gamma_q$ 
(for some $q\geq1$) and $\pi'/\sqrt{\pi'}\cong Z/3Z$.
Then $q$ is odd and either
\begin{enumerate}
\item$\pi=\pi(e,\eta)$, 
where $e$ is even and $\eta=1$ or $-1$, and has the presentation

$\langle{t,x,z}\mid{x^3}=(x^{3e-1}z^{-1})^3=z^{3\eta},~
txt^{-1}=x^{-1},~tzt^{-1}=z^{-1}\rangle$; \quad{or}

\item$\pi$ has the presentation 

$\langle{t,x}\mid{x^3}=t^{-1}x^{-3}t,
~x^{3e-1}t^{-1}xt=xtx^{-1}t^{-1}x^{-1}\rangle$,
where $e$ is even.
\end{enumerate}
Every strict weight orbit representing a given generator $t$ 
for $\pi/\pi'$ contains an unique element of the form $u^nt$,
where $u=z^{-1}x$ in case $(1)$ and $u=t^{-1}xtx$ in case $(2)$, 
and every such element is a weight element.
\end{theorem}

\begin{proof} 
It follows easily from Lemma 16.8 that $\zeta\sqrt{\pi'}=\zeta\pi'$ and
$G=\pi'/\zeta\pi'$ is isomorphic to $\pi_1^{orb}(S^2(3,3,3))\cong{Z^2}\rtimes_{-B}(Z/3Z)$, 
where $-B=
\left(\smallmatrix 0&-1\\
1&-1\endsmallmatrix\right)$.
As $\pi'$ is a torsion-free central extension of $G$ by $\mathbb{Z}$
it has a presentation
\begin{equation*}
\langle 
h,x,y,z\mid x^3=y^3=z^{3\eta}=h,\medspace xyz=h^e\rangle
\end{equation*}
for some $\eta=\pm1$ and $e\in\mathbb{Z}$.
The image of $h$ in $\pi'$ generates $\zeta\pi'$,
and the images of $u=z^{-1}x$ and $v=zuz^{-1}$ 
in $G=\pi'/\langle h\rangle$ 
form a basis for the translation subgroup $T(G)\cong \mathbb{Z}^2$.
Hence $\sqrt{\pi'}$ is generated by $u$, $v$ and $h$.
Since 
\[
vuv^{-1}u^{-1}=xz^{-2}xzx^{-2}z=xz^{-3}.zxzx.x^{-3}z
=x.x^{-3}z^{-3}(y^{-1}h^e)^3x^{-1}=h^{3e-2-\eta},
\] 
we must have $q=|3e-2-\eta|$.
Since $\pi'$ admits a meridianal automorphism and 
$\pi'/(\pi')^2\cong Z/(2,e-1)$,
$e$ must be even, and so $q$ is odd.
Moreover, $\eta=1$ if and only if 3 divides $q$.
In terms of the new generators, $\pi'$ has the presentation
\[
\langle{u,v,z}\mid{zuz^{-1}=v},~zvz^{-1}=u^{-1}v^{-1}z^{9e-6\eta-3},~
vuv^{-1}u^{-1}=z^{9e-6-3\eta}\rangle.
\]
Since $\pi'/\pi''$ is finite, $Hom(\pi',\zeta\pi')=0$,
and so the natural homomorphism from $Aut(\pi')$ to $Aut(G)$ is injective.
An automorphism $\phi$ of $\pi'$ must preserve 
characteristic subgroups such as $\zeta\pi'=\langle{z^3}\rangle$ and 
$\sqrt{\pi'}=\langle{u,v,z^3}\rangle$.
Let $K$ be the subgroup of $Aut(\pi')$ consisting of automorphisms 
which induce the identity on the subquotients 
$\pi'/\sqrt{\pi'}\cong{Z/3Z}$ and 
$\sqrt{\pi'}/\zeta\pi'\cong\mathbb{Z}^2$.
Automorphisms in $K$ also fix the centre, and are of the form $k_{m,n}$,
where
\[{k_{m,n}(u)=uz^{3\eta{s}}},~{k_{m,n}(v)=vz^{3\eta{t}}}\quad
\mathrm{and}\quad{k_{m,n}(z)=z^{3\eta{p}+1}u^mv^n},\]
for $(m,n)\in\mathbb{Z}^2$. 
These formulae define an automorphism if and only if
\[s-t=-n(3e-2-\eta),\quad
s+2t=m(3e-2-\eta)\quad\mathrm{and}
\]
\[
6p=(m+n)((m+n-1)(3e-2-\eta)+2(\eta-1)).
\]
In particular, the parameters $m$ and $n$ determine $p$, $s$ and $t$.
Conjugation by $u$ and $v$ give
$c_u=k_{-2,-1}$ and $c_v=k_{1,-1}$, respectively.
Let $k=k_{1,0}$.
If $\eta=1$ these equations have solutions
$p,s,t\in\mathbb{Z}$ for all $m,n\in\mathbb{Z}$,
and $K=\langle{k},c_u\rangle$.
If $\eta=-1$ we must have $m+n\equiv0$ {\it mod} $(3)$,
and $K=\langle{c_u},c_v\rangle$.
In each case, $K\cong\mathbb{Z}^2$.
We may define further automorphisms $b$ and $r$ by setting
\[{b(u)=v^{-1}z^{3\eta{e}-3}},\quad 
b(v)=uvz^{3\eta(e-1)}\quad\mathrm{and}\quad {b(z)=z};\quad\mathrm{and}\]
\[{r(u)=v^{-1}},\quad{r(v)=u^{-1}}\quad\mathrm{and}\quad{r(z)=z^{-1}}.\]
Then $b^6=r^2=(br)^2=1$ and $rk=kr$, while $b^4=c_z$
is conjugation by $z$.
If $\eta=-1$ then $Out(\pi')\cong(Z/2Z)^2$ is
generated by the images of $b$ and $r$,
while if $\eta=1$ then $Out(\pi')\cong{S_3}\times{Z/2Z}$
is generated by the images of $b,k$ and $r$.
Since $bk^nr(u)\equiv{u}$ {\it mod} $\pi''$,
such automorphisms are not meridianal.
Thus if $\eta=+1$ there are two classes of meridianal automorphisms 
(up to conjugacy and inversion in $Out(\pi')$),
represented by $r$ and $kr$, 
while if $\eta=-1$ there is only one, represented by $r$.
It is easily seen that $\pi(e,\eta)\cong\pi'\rtimes_r\mathbb{Z}$.
If $\eta=1$ and the meridianal automorphism is $kr$ then 
$kr(x)=x^{-1}zx^{2-3e}$ and $kr(z)=x^{-1}$, 
which leads to the presentation (2).

In each case,
$\pi$ has an automorphism which is the identity on $\pi'$ 
and sends $t$ to $th$, and 
the argument of Theorem 16.13 may be adapted 
to prove the final assertion.
\end{proof}

If $\pi\cong\pi(e,\eta)$ then 
$H_1(\pi;\Lambda/3\Lambda)\cong{H_2(\pi;\Lambda/3\Lambda)}
\cong(\Lambda/(3,t+1))^2$, 
and so the presentation in (1) is optimal, by Lemma 16.11.
The subgroup $A=\langle{t^2,x^3}\rangle\cong\mathbb{Z}^2$ 
is normal but not central.
The quotient $\pi(e,\eta)/A$ is 
$\pi_1^{orb}(\mathbb{D}^2(\overline{3},\overline{3},\overline{3}))$,
and $M(\tau_2k(e,\eta))$ is Seifert fibred over
$\mathbb{D}^2(\overline{3},\overline{3},\overline{3})$.

The presentation in (2) is also optimal, by Theorem 3.6.
Since $k^3=c_u^{-1}c_v$ and $(kr)^6=k^6$ in $Aut(\pi')$,
$t^6u^2v^{-2}$ is central in $\pi$.
The subgroup $A=\langle{t^6u^2v^{-2},x^3}\rangle\cong\mathbb{Z}^2$ 
is again normal but not central.
The quotient $\pi/A$ is $\pi_1^{orb}(\mathbb{D}^2(3,\overline{3}))$,
and the knot manifold is Seifert fibred over $\mathbb{D}^2(3,\overline{3})$.

The knot $k(0,-1)=9_{46}$ is doubly slice,
and hence so are all of its twist spins.
Since $M(\tau_2k(0,-1))$ is a $\mathbb{N}il^3\times\mathbb{E}^1$-manifold
it is determined up to homeomorphism by its group.
Hence every knot with group $\pi(0,-1)$ must be doubly slice.
No other non-trivial 2-knot with torsion-free, 
elementary amenable group is doubly slice.
(The Alexander polynomials of Fox's knot, the Cappell-Shaneson 2-knots 
and the knots of Theorem 16.14 are irreducible and not constant, 
while the Farber-Levine pairings of the others are not hyperbolic.)

%% file: m5-17.tex
\chapter{Knot manifolds and geometries}

In this chapter we shall attempt to characterize certain 2-knots in terms 
of algebraic invariants.
As every 2-knot $K$ may be recovered (up to orientations and Gluck 
reconstruction) from $M(K)$ together with the orbit of a weight class in 
$\pi=\pi K$ under the action of self homeomorphisms of $M$,
we need to characterize $M(K)$ up to homeomorphism. 
After some general remarks on the algebraic 2-type in \S1,
and on surgery in \S2, 
we shall concentrate on three special cases: when $M(K)$ is aspherical,
when $\pi'$ is finite and when $g.d.\pi=2$.

When $\pi$ is torsion-free and virtually poly-$Z$ the surgery obstructions 
vanish, and when it is poly-$Z$ the weight class is unique.
The surgery obstruction groups are notoriously difficult to compute
if $\pi$ has torsion.
However we can show that there are infinitely many distinct 2-knots $K$ 
such that $M(K)$ is simple homotopy equivalent to $M(\tau_23_1)$;
among these knots only $\tau_23_1$ has a minimal Seifert hypersurface.
If $\pi=\Phi$ the homotopy type of $M(K)$ determines the exterior of the knot;
the difficulty here is in finding a homotopy equivalence from $M(K)$ 
to a standard model. 

In the final sections we shall consider which knot manifolds are homeomorphic 
to geometric 4-manifolds or complex surfaces.
If $M(K)$ is geometric then either $K$ is a Cappell-Shaneson knot or the 
geometry must be one of $\mathbb{E}^4$, $\mathbb{N}il^3\times\mathbb{E}^1$, 
$\mathbb{S}ol^4_1$, $\widetilde{\mathbb{SL}}\times\mathbb{E}^1$,
$\mathbb{H}^3\times\mathbb{E}^1$ or $\mathbb{S}^3\times\mathbb{E}^1$.
If $M(K)$ is homeomorphic to a complex surface then
either $K$ is a branched twist spin of a torus knot or $M(K)$ admits one of the
geometries $\mathbb{N}il^3\times\mathbb{E}^1$, $\mathbb{S}ol^4_0$ or 
$\widetilde{\mathbb{SL}}\times\mathbb{E}^1$.

\section{Homotopy classification of $M(K)$}

Let $K$ and $K_1$ be 2-knots and suppose that $\alpha:\pi=\pi K\to \pi K_1$
and $\beta:\pi_2(M)\to\pi_2(M_1)$ determine an isomorphism of the algebraic
2-types of $M=M(K)$ and $M_1=M(K_1)$. 
Since the infinite cyclic covers $M'$ and $M_1'$ are homotopy equivalent 
to 3-complexes there is a map $h:M'\to M_1'$ such that
$\pi_1(h)=\alpha|_\pi$ and $\pi_2(h)=\beta$.
If $\pi=\pi K$ has one end then $\pi_3(M)\cong \Gamma(\pi_2(M))$ and so 
$h$ is a homotopy equivalence.
Let $t$ and $t_1=\alpha(t)$ be corresponding generators of
$Aut(M'/M)$ and $Aut(M_1'/M_1)$, respectively.
Then $h^{-1}t_1^{-1}ht$ is a self homotopy equivalence of $M'$ which fixes 
the algebraic 2-type.
If this is homotopic to $id_{M'}$ then $M$ and $M_1$ are homotopy equivalent,
since up to homotopy they are the mapping tori of $t$ and $t_1$, respectively.
Thus the homotopy classification of such knot manifolds may be 
largely reduced to determining the obstructions to homotoping a self-map 
of a 3-complex to the identity.

We may use a similar idea to approach this problem in another way.
Under the same hypotheses on $K$ and $K_1$ there is a map 
$f_o:M\setminus{intD^4}\to M_1$ inducing 
isomorphisms of the algebraic 2-types.          
If $\pi$ has one end $\pi_3(f_o)$ is an epimorphism, 
and so $f_o$ is 3-connected.
If there is an extension $f:M\to M_1$ then it is a homotopy equivalence, as
it induces isomorphisms on the homology of the universal covering spaces.

If $c.d.\pi\leq2$ the algebraic 2-type is determined by $\pi$, 
for then $\pi_2(M)\cong\overline{H^2(\pi;\mathbb{Z}[\pi])}$, 
by Theorem 3.12, 
and the $k$-invariant is 0.
In particular, if $\pi'$ is free of rank $r$ then $M(K)$ 
is homotopy equivalent to the mapping torus of a self-homeomorphism 
of $\sharp^rS^1\times S^2$,                                               
by Corollary 4.5.1.   
If $\pi=\Phi$ then $K$ is one of Fox's examples.
(See \S17.6 below.) 

The related problem of determining the homotopy type 
of the exterior of a 2-knot
has been considered in \cite{[Lo81], [Pl83]} and \cite{[PS85]}.
The examples considered in \cite{[Pl83]} do not test the adequacy of 
the algebraic 2-type for the present problem, 
as in each case either $\pi'$ is finite or $M(K)$ is aspherical.
The examples of \cite{[PS85]} probably show that in general $M(K)$ is 
not determined by $\pi$ and $\pi_2(M(K))$ alone.

\section{Surgery}
The natural transformations $I_G:G\to L^s_5(G)$ defined in Chapter 6 
clearly factor through $G/G'$.
If $\alpha:G\to\mathbb{Z}$ induces an isomorphism on abelianization 
then the homomorphism $\hat I_G=I_G\alpha^{-1}I_{\mathbb{Z}}^{-1}$ 
is a canonical splitting for $L_5(\alpha)$.

\begin{theorem} 
Let $K$ be a $2$-knot. 
If $L^s_5(\pi{K})\cong\mathbb{Z}$ and $N$ is simple homotopy equivalent 
to $M(K)$ then $N$ is $s$-cobordant to $M(K)$.
\end{theorem}

\begin{proof} 
Since $M=M(K)$ is orientable and 
$[M,G/TOP]\cong H^4(M;\mathbb{Z})\cong\mathbb{Z}$
the surgery obstruction map $\sigma_4:[M(K),G/TOP]\to L_4^s(\pi{K})$ 
is injective, by Theorem 6.6.
The image of $L_5(\mathbb{Z})$ under $\hat I_{\pi{K}}$ acts trivially on $S_{TOP}^s(M(K))$, by Theorem 6.7.
Hence there is a normal cobordism with obstruction 0
from any simple homotopy equivalence $f:N\to M$ to $id_M$.
\end{proof}

\begin{cor}
[Freedman] 
A $2$-knot $K$ is trivial if and only if $\pi{K}\cong\mathbb{Z}$.
\end{cor}
              
\begin{proof} 
The condition is clearly necessary.
Conversely, 
if $\pi{K}\cong\mathbb{Z}$ then $M(K)$ is homeomorphic to $S^3\times S^1$, 
by Theorem 6.11.
Since the meridian is unique up to inversion 
and the unknot is clearly reflexive the result follows.
\end{proof}

Theorem 17.1 applies if $\pi$ is a classical knot group \cite{[AFR97]},
or if $c.d.\pi=2$ and $\pi$ is square root closed accessible \cite{[Ca73]}.

Surgery on an $s$-concordance ${\mathcal K}$ from $K_0$ to $K_1$ gives an 
$s$-cobordism from $M(K_0)$ to $M(K_1)$ in which the meridians are conjugate.
Conversely, if $M(K)$ and $M(K_1)$ are $s$-cobordant via such an $s$-cobordism 
then $K_1$ is $s$-concordant to $K$ or $K^*$.
In particular, if $K$ is reflexive then $K$ and $K_1$ are $s$-concordant.

The next lemma follows immediately from Perelman's work (see \cite{[B-P]}).
In the original version of this book, 
we used a 4-dimensional surgery argument. 
We retain the statement so that the numbering of the other
results not be changed.

\begin{lemma}
Let $K$ be a $2$-knot. 
Then $K$ has a Seifert hypersurface 
which contains no fake $3$-cells.
\qed
\end{lemma}
                     
\section{The aspherical cases}

Whenever the group of a 2-knot $K$ contains a sufficiently large abelian 
normal subgroup $M(K)$ is aspherical. 
This holds for most twist spins of prime knots.

\begin{theorem} 
Let $K$ be a $2$-knot with group $\pi=\pi{K}$. 
If $\sqrt\pi$ is abelian of rank $1$ and $e(\pi/\sqrt\pi)=1$ 
or if $h(\sqrt\pi)\geq2$ then $\widetilde M(K)$ 
is homeomorphic to $\mathbb{R}^4$.
\end{theorem}
                      
\begin{proof}                                       
If $\sqrt\pi$ is abelian of rank 1
and $\pi/\sqrt\pi$ has one end then $M$ is aspherical,
by Theorem 15.11, 
and $\pi$ is 1-connected at $\infty$ \cite[Theorem 1]{[Mi87]}. 
If $h(\sqrt\pi)=2$ then $\sqrt\pi\cong\mathbb{Z}^2$ and 
$M$ is $s$-cobordant to the mapping torus of a self homeomorphism 
of a $\widetilde{\mathbb{SL}}$-manifold,
by Theorem 16.2.
If $h(\sqrt\pi)\geq3$ then $M$ is homeomorphic to an infrasolvmanifold, 
by Theorem 8.1.
In all cases, $\widetilde M$ is contractible and 1-connected at $\infty$, 
and so is homeomorphic to $\mathbb{R}^4$ \cite{[FQ]}.
\end{proof}

Is there a 2-knot $K$ with $\widetilde M(K)$ contractible but not
1-connected at $\infty$?

\begin{theorem} 
Let $K$ be a $2$-knot such that $\pi=\pi K$ is
torsion-free and virtually poly-$Z$.
Then $K$ is determined up to Gluck reconstruction by $\pi$ 
together with a generator of $H_4(\pi;\mathbb{Z})$ and the strict weight orbit 
of a meridian.
\end{theorem}

\begin{proof} 
If $\pi\cong\mathbb{Z}$ then $K$ must be trivial, and so we may assume
that $\pi$ is torsion-free and virtually poly-$Z$ of Hirsch length 4. 
Hence $M(K)$ is aspherical and is determined up to homeomorphism by $\pi$, and
every automorphism of $\pi$ may be realized by a self homeomorphism of $M(K)$, 
by Theorem 6.11.
Since $M(K)$ is aspherical its orientation is determined by a generator 
of $H_4(\pi;\mathbb{Z})$.
\end{proof}

This theorem applies in particular to the Cappell-Shaneson 2-knots, which have
an unique strict weight orbit, up to inversion.
(A similar argument applies to Cappell-Shaneson $n$-knots with $n>2$, provided
we assume also that $\pi_i(X(K))=0$ for $2\leq i\leq (n+1)/2$.)

\begin{theorem}
Let $K$ be a $2$-knot with group $\pi=\pi{K}$.
Then $K$ is $s$-concordant to a fibred knot with closed fibre 
a $\widetilde{\mathbb{SL}}$-manifold if and only if $\pi'$ is finitely
generated, $\zeta\pi'\cong\mathbb{Z}$ and $\pi$ is not virtually solvable.
The fibred knot is determined up to Gluck reconstruction by 
$\pi$ together with a generator of
$H_4(\pi;\mathbb{Z})$ and the strict weight orbit of a meridian.
\end{theorem}

\begin{proof} The conditions are clearly necessary.
If they hold then $M(K)$ is aspherical, by Theorem 15.11, 
so every automorphism of $\pi$ is induced by a self homotopy equivalence 
of $M(K)$. 
Moreover $\pi'$ is a $PD_3$-group, by Theorem 4.5.(3).
As $\zeta\pi'\cong\mathbb{Z}$ and $\pi$ is not virtually solvable,
$\pi'$ is the fundamental group of a $\widetilde{\mathbb{SL}}$-manifold, 
by Lemma 16.1.
Therefore $M(K)$ is determined up to $s$-cobordism by $\pi$, 
by Theorem 13.2.
The rest is standard.
\end{proof}

Branched twist spins of torus knots are perhaps the most important 
examples of such knots, but there are others.
(See Chapter 16.)
Is every 2-knot $K$ such that $\pi=\pi K$ is a $PD_4^+$-group determined up to
$s$-concordance and Gluck reconstruction by $\pi$ together 
with a generator of $H_4(\pi;\mathbb{Z})$ and a strict weight orbit?
Is $K$ $s$-concordant to a fibred knot with aspherical closed fibre
if and only if $\pi'$ is finitely generated and has one end?
(This follows from \cite{[Ro11]} and Theorem 6.8 if $Wh(\pi)=0$.)

\section{Quasifibres and minimal Seifert hypersurfaces}

Let $M$ be a closed 4-manifold with fundamental group $\pi$.
If $f:M\to S^1$ is a map which is transverse to $p\in S^1$ then
$\widehat V=f^{-1}(p)$ is a codimension 1 submanifold with a product 
neighbourhood $N\cong \widehat V\times [-1,1]$.
If moreover the induced homomorphism $f_*:\pi\to\mathbb{Z}$ 
is an epimorphism and each of the inclusions 
$j_\pm:\widehat V\cong
\widehat V\times\{\pm 1\}\subset W=M\setminus{V}\times (-1,1)$
induces monomorphisms on fundamental groups then we
shall say that $\widehat V$ is a {\it quasifibre} for $f$.
The group $\pi$ is then an HNN extension with base $\pi_1(W)$ and
associated subgroups $j_{\pm *}(\pi_1(\widehat V)$,
by Van Kampen's Theorem.
Every fibre of a bundle projection is a quasifibre.
We may use the notion of quasifibre to interpolate between the homotopy
fibration theorem of Chapter 4 and a TOP fibration theorem.
(See also Theorem 6.12 and Theorem 17.7.)

\begin{theorem} 
Let $M$ be a closed $4$-manifold with $\chi(M)=0$ and such that 
$\pi=\pi_1(M)$ is an extension of $\mathbb{Z}$ by a finitely
generated normal subgroup $\nu$.
If there is a map $f:M\to S^1$ inducing an epimorphism with kernel $\nu$ and 
which has a quasifibre $\widehat V$ then the infinite cyclic covering space
$M_\nu$ associated with $\nu$ is homotopy equivalent to $\widehat V$.
\end{theorem}
                     
\begin{proof}
As $\nu$ is finitely generated the monomorphisms $j_{\pm*}$ must be 
isomorphisms.
Therefore $\nu$ is finitely presentable, and so $M_\nu$ is a $PD_3$-complex,
by Theorem 4.5.
Now $M_\nu\cong W\times\mathbb{Z}/\sim$,
where $(j_+(v),n)\sim(j_-(v),n+1)$ for all $v\in\widehat V$ 
and $n\in\mathbb{Z}$.
Let $\tilde j(v)$ be the image of $(j_+(v),0)$ in $M_\nu$.
Then $\pi_1(\tilde j)$ is an isomorphism.
A Mayer-Vietoris argument shows that $\tilde j$ has degree 1,
and so $\tilde j$ is a homotopy equivalence.
\end{proof}

One could use duality instead to show that
$H_s=H_s(W,\partial_{\pm} W;\mathbb{Z}[\pi])= 0$ 
for $s\not=2$, while $H_2$ is a stably free $\mathbb{Z}[\pi]$-module, 
of rank $\chi(W,\partial_{\pm} W)=0$.
Since $\mathbb{Z}[\pi]$ is weakly finite this module is 0,
and so $W$ is an $h$-cobordism.

\begin{cor}
Let $K$ be a $2$-knot with group $\pi=\pi{K}$.
If $\pi'$ is finitely generated and $K$ has a minimal Seifert 
hypersurface $V$ such that every self homotopy 
equivalence of $\widehat V$ is homotopic to a homeomorphism 
then $M(K)$ is homotopy equivalent to $M(K_1)$, 
where $M(K_1)$ is a fibred $2$-knot with fibre $V$.
\end{cor}
                     
\begin{proof}
Let $j_+^{-1}:M(K)'\to\widehat V$ be a homotopy 
inverse to the homotopy equivalence $j_+$, and let $\theta$ be a self 
homeomorphism of $\widehat V$ homotopic to $j_+^{-1}j_-$. Then
$j_+\theta j_+^{-1}$ is homotopic to a generator of $Aut(M(K)'/M(K))$, 
and so the mapping torus of $\theta$ is homotopy equivalent to $M(K)$.
Surgery on this mapping torus gives such a knot $K_1$.
\end{proof}

If a Seifert hypersurface $V$ for a 2-knot has fundamental group $\mathbb{Z}$
then the Mayer-Vietoris sequence for $H_*(M(K);\Lambda)$ gives 
$H_1(X')\cong\Lambda/(ta_+-a_-)$, 
where $a_\pm :H_1(V)\to H_1(S^4\setminus{V})$. 
Since $H_1(X)=\mathbb{Z}$ we must have $a_+-a_-=\pm 1$. 
If $a_+a_-\not=0$ then $V$ is minimal. 
However one of $a_+$ or $a_-$ could be 0, 
in which case $V$ may not be minimal.
The group $\Phi$ is realized by ribbon knots 
with such minimal Seifert hypersurfaces 
(homeomorphic to $S^2\times S^1\setminus{int{D^3}}$) \cite{[Fo62]}. 
Thus minimality does not imply that $\pi'$ is finitely generated.
                                           
If a 2-knot $K$ has a minimal Seifert hypersurface then $\pi{K}$
is an HNN extension with finitely presentable base 
and associated subgroups.
It remains an open question whether every 2-knot group is of this type.
(There are high dimensional knot groups which are not so 
\cite{[Si91], [Si96]}.)
Yoshikawa has shown that there are ribbon 2-knots 
whose groups are HNN extensions with base a torus knot group 
and associated subgroups $\mathbb{Z}$ but which cannot be expressed as 
HNN extensions with base a free group \cite{[Yo88]}.

An argument of Trace implies that if $V$ is a Seifert hypersurface for a
fibred $n$-knot $K$ then there is a degree-1 map from 
$\widehat V=V\cup D^{n+1}$ to the closed fibre $\widehat F$ \cite{[Tr86]}.
For the embedding of $V$ in $X$ extends to an embedding of $\widehat V$ in $M$,
which lifts to an embedding in $M'$.
Since the image of $[\widehat V]$ in $H_{n+1}(M;\mathbb{Z})$ is Poincar\'e dual 
to a generator of $H^1(M;\mathbb{Z})=Hom(\pi,\mathbb{Z})=[M,S^1]$ 
its image in $H_{n+1}(M';\mathbb{Z})\cong\mathbb{Z}$ is a generator. 
In particular, if $K$ is a fibred 2-knot and $\widehat F$ has a summand which 
is aspherical or whose fundamental group is a nontrivial finite group 
then $\pi_1(V)$ cannot be free.
Similarly, as the Gromov norm of a 3-manifold does not 
increase under degree 1 maps, if $\widehat F$ is a $\mathbb{H}^3$-manifold 
then $\widehat V$ cannot be a graph manifold \cite{[Ru90]}.
Rubermann observes also that the ``Seifert volume" of \cite{[BG84]} 
may be used instead to show that if $\widehat F$ is a
$\widetilde{\mathbb{SL}}$-manifold then $\widehat V$ 
must have nonzero Seifert volume. 
(Connected sums of $\mathbb{E}^3$-, $\mathbb{S}^3$-,
$\mathbb{N}il^3$-, $\mathbb{S}ol^3$-, $\mathbb{S}^2\times\mathbb{E}^1$-
or $\mathbb{H}^2\times\mathbb{E}^1$-manifolds all
have Seifert volume 0 \cite{[BG84]}.)

\section{The spherical cases}

Let $\pi$ be a 2-knot group with commutator subgroup
$\pi'\cong P\times (Z/(2r+1)Z)$, 
where $P=1$, $Q(8)$, $T_k^*$ or $I^*$.  
The meridianal automorphism induces the identity on the set of irreducible 
real representations of $\pi'$, except when $P=Q(8)$.
(It permutes the three nontrivial 1-dimensional representations
when $\pi'\cong Q(8)$, and similarly when $\pi'\cong Q(8)\times(Z/nZ)$.)
It then follows as in Chapter 11 that $L^s_5(\pi)$ has rank $r+1$, $3(r+1)$, 
$3^{k-1}(5+7r)$ or $9(r+1)$, respectively.                                       Hence if $\pi'\not=1$ then there are infinitely many 
distinct 2-knots with group $\pi$,
since the group of self homotopy equivalences of $M(K)$ is finite.

The simplest nontrivial such  group is
$\pi=\pi\tau_23_1\cong(Z/3Z)\rtimes_{-1}\mathbb{Z}$.
If $K$ is any 2-knot with this group then $M(K)$ is homotopy equivalent to 
$M(\tau_2 3_1)$.
Since $Wh(Z/3Z)=0$ \cite{[Hg40]} and $L_5(Z/3Z)=0$ \cite{[Ba75]} 
we have $L_5^s(\pi)\cong L_4(\pi')\cong\mathbb{Z}^2$,
but we do not know whether $Wh(\pi)=0$. 

\begin{theorem} 
Let $K$ be a $2$-knot with group $\pi=\pi{K}$ such that $\pi'\cong{Z/3Z}$, 
and which has a minimal Seifert hypersurface. 
Then $K=\pm\tau_23_1$.
\end{theorem}

\begin{proof}
Let $V$ be a minimal Seifert hypersurface for $K$.
Let $\widehat V=V\cup D^3$ and $W=M(K)\setminus{V}\times (-1,1)$.
Then $W$ is an $h$-cobordism from $\widehat V$ to itself. (See the remark
following Theorem 6.)
Therefore $W\cong\widehat V\times I$, by surgery over $Z/3Z$. 
(Note that $Wh(Z/3Z)=L_5(Z/3Z)=0$.)
Hence $M$ fibres over $S^1$.
The closed fibre must be the lens space $L(3,1)$,
by Perelman's work (see \cite{[B-P]}).
Since $\tau_23_1$ is invertible and reflexive, $K$ must be $\pm\tau_23_1$.
\end{proof}

As none of the other 2-knots with this group has a minimal Seifert surface,
they are all further counter-examples to the most natural 
4-dimensional analogue of Farrell's fibration theorem.
We do not know whether any of these knots (other than $\tau_2 3_1$) is 
PL in some PL structure on $S^4$.
       
Let $F$ be an $\mathbb{S}^3$-group,
and let $W=(W;j_\pm)$ be an $h$-cobordism with 
homeomorphisms $j_\pm:N\to \partial_\pm W$, where $N=S^3/F$.
Then $W$ is an $s$-cobordism \cite{[KS92]}.
The set of such $s$-cobordisms from $N$ to itself
is a finite abelian group with respect to stacking of cobordisms.
All such $s$-cobordisms are products if $F$ is cyclic, 
but there are nontrivial examples if $F\cong Q(8)\times (Z/pZ)$, 
for any odd prime $p$ \cite{[KS95]}.
If $\phi$ is a self-homeomorphism of $N$ the closed 4-manifold $Z_\phi$ 
obtained by identifying the ends of $W$ via $j_+\phi j_-^{-1}$ is homotopy
equivalent to $M(\phi)$.
However if $Z_\phi$ is a mapping torus of a self-homeomorphism of $N$ then
$W$ is trivial.
In particular, if $\phi$ induces a meridianal automorphism of $F$ then
$Z_\phi\cong M(K)$ for an exotic 2-knot $K$ with 
$\pi'\cong F$ and which has a minimal Seifert hypersurface, 
but which is not fibred with geometric fibre.

\section{Finite geometric dimension 2}

Knot groups with finite 2-dimensional Eilenberg-Mac Lane complexes
have deficiency 1, by Theorem 2.8, and so are 2-knot groups.
This class includes all classical knot groups, 
all knot groups with free commutator subgroup and all knot groups
in the class ${\mathcal X}$ (such as those of Theorems 15.7 and 15.18).

\begin{theorem} 
Let $K$ be a $2$-knot.
If $\pi{K}$ is a $1$-knot group or an ${\mathcal X}$-group 
then $M(K)$ is determined up to $s$-cobordism by its homotopy type.
\end{theorem}

\begin{proof} 
This is an immediate consequence of Lemma 6.9,
if $\pi$ is an ${\mathcal X}$-group.
If $\pi$ is a nontrivial classical knot group it follows from Theorem 17.1, 
since $Wh(\pi{K})=0$ \cite{[Wd78]} 
and $L^s_5(\pi{K})\cong\mathbb{Z}$ \cite{[AFR97]}.
\end{proof}

\begin{cor}
Any $2$-knot $K$ with group $\Phi$ is ambient isotopic to one of Fox's
examples.
\end{cor}

\begin{proof} 
There is an unique homotopy type of closed orientable 4-manifolds $M$ 
with $\pi_1(M)=\Phi$ and $\chi(M)=0$ \cite{[Hi09],[Hi20]}.
Since $\Phi$ is metabelian, $s$-cobordism implies homeomorphism and 
there is an unique weight class up to inversion,
so $X(K)$ is determined by the homotopy type of $M(K)$.
Since Examples 10 and 11 of \cite{[Fo62]} are ribbon knots they are 
$-$amphicheiral and are determined up to reflection by their exterior.
\end{proof}

Fox's examples are mirror images.
Since their Alexander polynomials are asymmetric,
they are not invertible. (Thus they are not isotopic.)
Nor are they deform spins, 
by the criterion of \cite{[BM09]}.

\begin{theorem} 
A $2$-knot $K$ with group $\pi=\pi{K}$ is $s$-concordant 
to a fibred knot with closed fibre $\sharp^r(S^1\times S^2)$ 
if and only if $\mathrm{def}(\pi)=1$ and $\pi'$ is finitely generated.
Moreover any such fibred $2$-knot is reflexive and homotopy ribbon.
\end{theorem}
 
\begin{proof} 
The conditions are clearly necessary. 
If they hold then $\pi'\cong F(r)$, for some $r\geq0$, 
by Corollary 4.3.1.
Then $M=M(K)$ is homotopy equivalent to a PL 4-manifold $N$ which fibres
over $S^1$ with fibre $\sharp^r(S^1\times S^2)$, 
by Corollary 4.5.1. 
Moreover $Wh(\pi)=0$, by Lemma 6.3, 
and $\pi$ is square root closed accessible,
so $I_\pi$ is an isomorphism, by Lemma 6.9, 
so there is an $s$-cobordism $W$ from $M$ to $N$, by Theorem 17.1.
We may embed an annulus $A=S^1\times [0,1]$ in $W$ so that 
$M\cap A=S^1\times\{0\}$ is a meridian for $K$ and $N\cap A=S^1\times\{1\}$.
Surgery on $A$ in $W$ then gives an $s$-concordance from $K$ to such a 
fibred knot $K_1$, 
which is reflexive \cite{[Gl62]} and homotopy ribbon \cite{[Co83]}.              \end{proof}

Based self homeomorphisms of $N=\sharp^r(S^1\times S^2)$
which induce the same automorphism of $\pi_1(N)=F(r)$
are conjugate up to isotopy if their images in $E_0(N)$ are conjugate
\cite{[HL74]}.
Their mapping tori are then orientation-preserving homeomorphic.
(Compare Corollary 4.5.)
The group $E_0(N)$ is a semidirect product $(Z/2Z)^r\rtimes{Aut(F(r))}$, 
with the natural action \cite{[Hn]}.
It follows easily that every fibred 2-knot with $\pi'$ free
is determined (among such knots) 
by its group together with the weight orbit of a meridian.
(However, $\pi 3_1$ has infinitely many weight orbits \cite{[Su85]}.)
Is every such group the group of a ribbon knot?

If $K=\sigma k$ is the Artin spin of a fibred 1-knot then $M(K)$ fibres over 
$S^1$ with fibre $\sharp^r (S^2\times S^1)$.
However not all such fibred 2-knots arise in this way.
It follows easily from Lemma 1.1 and the fact that 
$Out(F(2))\cong{GL}(2,\mathbb{Z})$ 
that there are just three knot groups $G\cong{F}(2)\rtimes\mathbb{Z}$,
namely $\pi 3_1$ (the trefoil knot group), $\pi 4_1$ 
(the figure eight knot group) and the group with presentation
\begin{equation*}
\langle x,y,t\mid txt^{-1}=y,\, tyt^{-1}=xy\rangle.
\end{equation*}                                                     
(Two of the four presentations given in \cite{[Rp60]} 
present isomorphic groups.)
The Alexander polynomial of the latter example is $t^2-t-1$, 
which is not symmetric,
and so this is not a classical knot group.
(See also \cite{[AY81], [Rt83]}.)

\begin{theorem} 
Let $K$ be a $2$-knot with group $\pi=\pi K$. Then
\begin{enumerate}
\item{if $c.d.\pi\leq2$ then $K$ is $\pi_1$-slice;}
\item{if $K$ is $\pi_1$-slice then $H_3(c_{M'(K)};\mathbb{Z})=0$};
\item{if $\pi'$ is finitely generated then $K$ is $\pi_1$-slice if and only if $\pi'$ is free};
\end{enumerate}
\end{theorem}
                     
\begin{proof} 
Let $M=M(K)$. If $c.d.\pi\leq2$ then 
$H_p(\pi;\Omega_{4-p}^{Spin})=0$ for $p>0$, 
since $\Omega_3^{Spin}=0$ and $H_2(\pi;Z/2Z)=0$,
and so signature gives an isomorphism $\Omega_4^{Spin}(K(\pi,1))\cong\mathbb{Z}$.
Hence $M=\partial{W}$ for some $Spin$ 5-manifold $W$ with a map $f$ to $K(\pi,1)$
such that $f|_M=c_M$. 
We may modify $W$ by elementary surgeries to make $\pi_1(f)$ an isomorphism 
and $H_2(W;\mathbb{Z})=0$.
Attaching a 2-handle to $W$ along a meridian for $K$ then gives 
a contractible 5-manifold with 1-connected boundary,
and so $K$ is $\pi_1$-slice. 

Let $R$ be an open regular neighbourhood in $D^5$ 
of a $\pi_1$-slice disc $\Delta$.
Since $c_M$ factors through $D^5\setminus{R}$ 
the first assertion follows from the exact 
sequence of homology (with coefficients $\Lambda$) 
for the pair $(D^5\setminus{R},M)$.

If $\pi'$ is finitely generated then $M'$ is a $PD_3$-space, by Theorem 4.5.
The image of $[M']$ in $H_3(\pi';\mathbb{Z})$ 
determines a projective homotopy equivalence 
of modules $C^2/\partial^1(C^1)\simeq{A}(\pi')$, 
by the argument of \cite[Theorem 4]{[Tu90]}.
(This does not need Turaev's assumption that $\pi'$ be finitely presentable.)
If this image is 0 then $id_{{A}(\pi')}\sim0$, 
so ${A}(\pi')$ is projective and $c.d.\pi'\leq1$. 
Therefore $\pi'$ is free.
The converse follows from Theorem 17.9, or from part (1) of this theorem.
\end{proof}

If $\pi'$ is free must $K$ be homotopy ribbon?
This would follow from ``homotopy connectivity implies geometric connectivity",
but our situation is just beyond the range of known results.
Is $M(K)$ determined up to $s$-cobordism by its group whenever $g.d.\pi\leq2$?
If $g.d.\pi\leq2$ then $\mathrm{def}(\pi)=1$,
and so $\pi$ is the group of a homotopy ribbon $2$-knot.
(See \S6 of Chapter 14.)
The conditions  ``$g.d.\pi\leq2$" and ``$\mathrm{def}(\pi)=1$"
are equivalent if the Whitehead Conjecture holds.
(See \S9 of Chapter 14.)
Are these conditions also equivalent to ``homotopy ribbon"?
(See \cite{[Hi08']} for further connections between these properties.)

\section{Geometric 2-knot manifolds}

The 2-knots $K$ for which $M(K)$ is homeomorphic to an infrasolvmanifold are
essentially known. 
There are three other geometries which may be realized by such knot manifolds.
All known examples are fibred, 
and most are derived from twist spins of classical knots.
However there are examples (for instance, those with $\pi'\cong
Q(8)\times(Z/nZ)$ for some $n>1$) which cannot be constructed from twist spins.
The remaining geometries may be eliminated very easily; only
$\mathbb{H}^2\times\mathbb{E}^2$ and $\mathbb{S}^2\times\mathbb{E}^2$
require a little argument.
                                          
\begin{theorem} 
Let $K$ be a $2$-knot with group $\pi=\pi K$. 
If $M(K)$ admits a geometry then the geometry is one of
$\mathbb{E}^4$, $\mathbb{N}il^3\times\mathbb{E}^1$, $\mathbb{S}ol^4_0$, $\mathbb{S}ol^4_1$,
$\mathbb{S}ol^4_{m,n}$ (for certain $m\not=n$ only), $\mathbb{S}^3\times\mathbb{E}^1$,
$\mathbb{H}^3\times\mathbb{E}^1$ or $\widetilde{\mathbb{SL}}\times\mathbb{E}^1$.
All these geometries occur.
\end{theorem}

\begin{proof} 
The knot manifold $M(K)$ is homeomorphic to an infrasolvmanifold 
if and only if $h(\sqrt\pi)\geq3$, by Theorem 8.1.
It is then determined up to homeomorphism by $\pi$.
We may then use the observations of \S10 of Chapter 8 to show
that $M(K)$ admits a geometry of solvable Lie type.
By Lemma 16.7 and Theorems 16.12 and 16.14 $\pi$ must be either 
$G(+)$ or $G(-)$, 
$\pi(e,\eta)$ for some even $b$ and $\epsilon=\pm 1$ or 
$\pi'\cong\mathbb{Z}^3$ or $\Gamma_q$ for some odd $q$.
We may identify the geometry on looking more closely at the
meridianal automorphism.

If $\pi\cong G(+)$ or $G(-)$ then $M(K)$ 
admits the geometry $\mathbb{E}^4 $.
If $\pi\cong\pi(e,\eta)$ then $M(K)$ is the mapping torus 
of an involution of a $\mathbb{N}il^3$-manifold, 
and so admits the geometry $\mathbb{N}il^3\times\mathbb{E}^1$.
If $\pi'\cong\mathbb{Z}^3 $ then $M(K)$ is homeomorphic to a $\mathbb{S}ol_{m,n}^4$- 
or $\mathbb{S}ol_0^4$-manifold.
More precisely, we may assume (up to change of orientations) that
the Alexander polynomial of $K$ is $t^3 -(m-1)t^2+mt-1$ for some integer $m$.
If $m\geq 6$ all the roots of this cubic are positive and 
the geometry is $\mathbb{S}ol^4_{m-1,m} $. 
If $0\leq m\leq 5$ two of the roots are complex conjugates and 
the geometry is $\mathbb{S}ol^4_0 $.
If $m<0$ two of the roots are negative and $\pi$ has a subgroup of index 2 
which is a discrete cocompact subgroup of $Sol^4_{m',n'}$,
where $m'=m^2 -2m+2$ and $n'=m^2-4m+1$, 
so the geometry is $\mathbb{S}ol^4_{m',n'}$.

If $\pi'\cong\Gamma_q $ and the image of the meridianal automorphism 
in $Out(\Gamma_q)$ has finite order then
$q=1$ and $K=\tau_63_1$ or $(\tau_63_1)^*=\tau_{6,5}3_1$. 
In this case $M(K)$ admits the geometry $\mathbb{N}il^3\times\mathbb{E}^1$. 
Otherwise (if $\pi'\cong\Gamma_q $ and the image of the 
meridianal automorphism in $Out(\Gamma_q)$ has infinite order) 
$M(K)$ admits the geometry $\mathbb{S}ol^4_1$.

If $K$ is a branched $r$-twist spin of the $(p,q)$-torus knot then $M(K)$
is a ${\mathbb{S}^3\times\mathbb{E}^1}$-manifold if $p^{-1}+q^{-1}+r^{-1}>1$,
and is a $\widetilde{\mathbb{SL}}\times \mathbb{E}^1$-manifold if
$p^{-1}+q^{-1}+r^{-1}<1$.
(The case $p^{-1}+q^{-1}+r^{-1}=1$ gives the
$\mathbb{N}il^3\times\mathbb{E}^1$-manifold $M(\tau_63_1)$.)
The manifolds obtained from 2-twist spins of 2-bridge knots and certain other 
``small" simple knots also have geometry $\mathbb{S}^3\times\mathbb{E}^1$.
Branched $r$-twist spins of simple (nontorus) knots with $r>2$ give 
$\mathbb{H}^3\times\mathbb{E}^1$-manifolds, 
excepting $M(\tau_34_1)\cong M(\tau_{3,2}4_1)$, which 
is the $\mathbb{E}^4$-manifold with group $G(+)$.

Every orientable $\mathbb{H}^2\times\mathbb{E}^2$-manifold 
is double covered by a K\"ahler surface \cite{[Wl86]}.
Since the unique double cover of a 2-knot manifold $M(K)$ has first Betti 
number 1 no such manifold can be an $\mathbb{H}^2\times\mathbb{E}^2$-manifold. 
(If $K$ is fibred we could instead exclude this geometry by Lemma 16.1.)
Since $\pi$ is infinite and $\chi(M(K))=0$ we may exclude
the geometries $\mathbb{S}^4$, $\mathbb{CP}^2$ and $\mathbb{S}^2\times\mathbb{S}^2$,
and $\mathbb{H}^4$, $\mathbb{H}^2(\mathbb{C})$, 
$\mathbb{H}^2\times\mathbb{H}^2$ 
and $\mathbb{S}^2\times\mathbb{H}^2$, respectively.
The geometry $\mathbb{S}^2\times\mathbb{E}^2$ may be excluded by Theorem 
10.10 or Lemma 16.1 (no group with two ends admits a meridianal automorphism),
while $\mathbb{F}^4$ is not realized by any closed 4-manifold. 
\end{proof}

In particular, no knot manifold is a $\mathbb{N}il^4$-manifold 
or a $\mathbb{S}ol^3\times\mathbb{E}^1$-manifold, and many 
$\mathbb{S}ol^4_{m,n}$-geometries do not arise in this way.
The knot manifolds which are infrasolvmanifolds or have geometry 
$\mathbb{S}^3\times\mathbb{E}^1$ are essentially known, by 
Theorems 8.1, 11.1, 15.4 and \S4 of Chapter 16.
The knot is uniquely determined up to Gluck reconstruction and change of
orientations if $\pi'\cong\mathbb{Z}^3$ 
(see Theorem 17.4 and the subsequent remarks above), 
$\Gamma_q$ (see \S3 of Chapter 18) or $Q(8)\times (Z/nZ)$ 
(since the weight class is then unique up to inversion). 
There are only six knots whose knot manifold admits
the geometry $\mathbb{S}ol^4_0$,
for $\mathbb{Z}[X]/(\Delta_a(X))$ is a PID if $0\leq{a}\leq5$.
(See the tables of \cite{[AR84]}.)
If the knot is fibred with closed fibre a lens space 
it is a 2-twist spin of a 2-bridge knot \cite{[Tr89]}.
The other knot groups corresponding to infrasolvmanifolds have infinitely
many weight orbits.

\begin{cor}
If $M(K)$ admits a geometry then it fibres over $S^1$,
and if $\pi$ is not solvable the closed monodromy has finite order.
\end{cor}

\begin{proof}
This is clear if $M(K)$ is an $\mathbb{S}^3\times\mathbb{E}^1$-manifold
or an infrasolvmanifold, and follows from Corollary 13.1.1 
and Theorem 16.2 if the geometry is 
$\widetilde{\mathbb{SL}}\times\mathbb{E}^1$.

If $M(K)$ is a $\mathbb{H}^3\times\mathbb{E}^1$-manifold 
we refine the argument of Theorem 9.3.
Since $\pi/\pi'\cong\mathbb{Z}$ and $\sqrt\pi=\pi\cap(\{1\}\times\mathbb{R})\not=1$
we may assume $\pi\leq{Isom(\mathbb{H}^3)\times\mathbb{R}}$,
and so $\pi'\leq{Isom(\mathbb{H}^3)\times\{1\}}$. 
Hence $\pi'$ is the fundamental group of a closed $\mathbb{H}^3$-manifold,
$N$ say, and $M(K)$ is the quotient of $N\times\mathbb{R}$ 
by the action induced by a meridian.
Thus $M(K)$ is a mapping torus, and so fibres over $S^1$.
\end{proof}

If the geometry is $\mathbb{H}^3\times\mathbb{E}^1$
is $M(K)\cong M(K_1)$ for some branched twist spin of a simple non-torus knot?
(See \S3 of Chapter 16.)
 
\begin{cor}
If $M(K)$ is Seifert fibred it is a 
$\widetilde{\mathbb{SL}}\times\mathbb{E}^1$-, 
$\mathbb{N}il^3\times\mathbb{E}^1$- or 
$\mathbb{S}^3\times\mathbb{E}^1$-manifold. 
\end{cor}

\begin{proof}
This follows from Ue's Theorem, Theorem 16.2 and Theorem 17.11.
\end{proof}

Are any 2-knot manifolds $M$ total spaces of orbifold bundles
with hyperbolic general fibre?
(The base $B$ must be flat, since $\chi(M)=0$, 
and $\pi_1^{orb}(B)$ must have cyclic abelianization.
Hence $B=S^2(2,3,6)$, $\mathbb{D}^2(\overline{3},\overline{3},\overline{3})$,
or $\mathbb{D}^2(3,\overline{3})$.)

It may be shown that if $k$ is a nontrivial 1-knot and $r\geq2$ then 
$M(\tau_rk)$ is geometric if and only if $k$ is simple, and has
a proper geometric decomposition if and only if $k$ is prime but not simple.
(The geometries of the pieces are then
$\mathbb{H}^2\times\mathbb{E}^2$,
$\widetilde{\mathbb{SL}}\times\mathbb{E}^1$
or $\mathbb{H}^3\times\mathbb{E}^1$.)
This follows from the fact that the $r$-fold cyclic branched cover of 
$(S^3,k)$ admits an equivariant JSJ decomposition, 
and has finitely generated $\pi_2$ if and only if $k$ is prime.

\section{Complex surfaces and 2-knot manifolds}

If a complex surface $S$ is homeomorphic to a 2-knot manifold $M(K)$ then $S$ 
is minimal, since $\beta_2 (S)=0$, and has Kodaira dimension $\kappa (S)=1$, 
0 or $-\infty$, since $\beta_1 (S)=1$ is odd. 
If $\kappa(S)=1$ or 0 then $S$ is elliptic and admits a compatible geometric 
structure, of type $\widetilde{\mathbb{SL}}\times\mathbb{E}^1$ or 
$\mathbb{N}il^3\times\mathbb{E}^1$, respectively \cite{[Ue90],[Ue91],[Wl86]}.
The only complex surfaces with $\kappa(S)=-\infty$, $\beta_1 (S)=1$ and 
$\beta_2 (S)=0$ are Inoue surfaces, which are not elliptic, 
but admit compatible geometries of 
type $\mathbb{S}ol^4_0$ or $\mathbb{S}ol^4_1$, and Hopf surfaces \cite{[Tl94]}.
An elliptic surface with Euler characteristic 0 has no exceptional fibres 
other than multiple tori. 

If $M(K)$ has a complex structure compatible with a geometry then the geometry 
is one of $\mathbb{S}ol^4_0$, $\mathbb{S}ol^4_1$, $\mathbb{N}il^3\times\mathbb{E}^1$,
$\mathbb{S}^3\times\mathbb{E}^1 $ or 
$\widetilde{\mathbb{SL}}\times\mathbb{E}^1 $,
by Theorem 4.5 of \cite{[Wl86]}.
Conversely, if $M(K)$ admits one of the first three of these geometries 
then it is homeomorphic to an Inoue surface of type $S_M$, 
an Inoue surface of type 
$\vrule width 0pt depth 7pt S^{(+)}_{N,p,q,r;t}$ 
or $S^{(-)}_{N,p,q,r} $, 
or an elliptic surface of Kodaira dimension 0, respectively.
(See \cite{[In74],[EO94]} and \cite[Chapter V]{[BHPV]}.)

\begin{lemma}
Let $K$ be a branched $r$-twist spin of the $(p,q)$-torus knot.
Then $M(K)$ is homeomorphic to an elliptic surface.
\end{lemma}
                      
\begin{proof} We shall adapt the argument of \cite[Lemma 1.1]{[Mi75]}. 
(See also \cite{[Ne83]}.)                             
Let $V_0 =
\{ (z_1,z_2,z_3)\in\mathbb{C}^3\setminus\{0\}|z_1^p+z_2^q+z_3^r =0\}$, 
and define an action of $\mathbb{C}^\times$ on $V_0 $ by 
$u.v=(u^{qr}z_1,u^{pr}z_2,u^{pq}z_3 )$ for all $u$ in $\mathbb{C}^\times$ 
and $v=(z_1,z_2,z_3)$ in $V_0 $.
Define functions $m:V_0 \to\mathbb{R}^+ $ and $n:V_0 \to m^{-1} (1)$ 
by 
\[m(v)=(|z_1 |^p +|z_2 |^q +|z_3 |^r )^{1/pqr} 
\quad\mathrm{and}\quad{n(v)=m(v)^{-1}v},\quad\forall~v\in{V_0}.
\]
The map $(m,n):V_0 \to m^{-1} (1)\times\mathbb{R}^+$ 
is an $\mathbb{R}^+$-equivariant homeomorphism, 
and so $m^{-1} (1)$ is homeomorphic to $V_0 /\mathbb{R}^+ $.
Therefore there is a homeomorphism from $m^{-1} (1)$ to $M(p,q,r)$, 
under which the action of the group of $r^{th} $ roots of unity on
$m^{-1} (1)=V_0 /\mathbb{R}^+$ corresponds to the group of covering 
homeomorphisms of $M(p,q,r)$ as the branched cyclic cover of $S^3 $, 
branched over the $(p,q)$-torus knot \cite{[Mi75]}.
The manifold $M(K)$ is the mapping torus of some generator 
of this group of self homeomorphisms of $M(p,q,r)$.
Let $\omega$ be the corresponding primitive $r${th} root of unity.
If $t>1$ then $t\omega$ generates a subgroup $\Omega$ 
of $\mathbb{C}^\times $ which acts freely and holomorphically 
on $V_0 $, and the quotient $V_0 /\Omega$ is an
elliptic surface over the curve $V_0 /\Omega$.
Moreover $V_0/\Omega$ is homeomorphic to the mapping torus of the self 
homeomorphism of $m^{-1} (1)$ which maps $v$ to 
$m(t\omega.v)^{-1}.t\omega.v=\omega m(t.v)^{-1} t.v$.
Since this map is isotopic to the map sending $v$ to $\omega.v$ this mapping 
torus is homeomorphic to $M(K)$.
\end{proof}

The Kodaira dimension of the elliptic surface in the above lemma is 1, 0 or 
$-\infty$ according as $p^{-1} +q^{-1} +r^{-1} $ is $<1 $, 1 or $>1$. 
In the next theorem we shall settle the case
of elliptic surfaces with $\kappa= -\infty$.

\begin{theorem} 
Let $K$ be a $2$-knot. Then $M(K)$ is homeomorphic
to a Hopf surface if and only if $K$ or its Gluck
reconstruction is a branched $r$-twist spin of the 
$(p,q)$-torus knot for some $p,q$ and $r$ such that 
$p^{-1} +q^{-1} +r^{-1} >1$.
\end{theorem}

\begin{proof} 
If $K=\tau_{r,s}k_{p,q}$ then $M(K)$ is homeomorphic to an elliptic surface, 
by Lemma 17.13, 
and the surface must be a Hopf surface if $p^{-1} +q^{-1} +r^{-1} >1$.

If $M(K)$ is homeomorphic to a Hopf surface then $\pi$ has two ends, 
and there is a monomorphism $h:\pi=\pi K\to GL(2,\mathbb{C})$
onto a subgroup which contains a contraction $c$ 
(Kodaira - see \cite{[Kt75]}).
Hence $\pi'$ is finite and $h(\pi')=h(\pi)\cap SL(2,\mathbb{C})$,
since $\det(c)\not=1$ and $\pi/\pi'\cong\mathbb{Z}$.
Finite subgroups of $SL(2,\mathbb{C})$ are conjugate to subgroups
of $SU(2)=S^3$, and so are cyclic, binary dihedral
or isomorphic to $T_1^*$, $O^*_1$ or $I^*$.
Therefore $\pi\cong\pi\tau_2k_{2,n}$,
$\pi\tau_33_1$, $\pi\tau_43_1$ or $\pi\tau_53_1$,
by Theorem 15.12 and the subsequent remarks.
Hopf surfaces with $\pi\cong\mathbb{Z}$ or $\pi$ nonabelian
are determined up to diffeomorphism by their 
fundamental groups \cite[Theorem 12]{[Kt75]}.
Therefore $M(K)$ is homeomorphic to the manifold 
of the corresponding torus knot.
If $\pi'$ is cyclic there is an unique weight orbit.
The weight orbits of $\tau_43_1$ are realized by $\tau_2k_{3,4}$ and
$\tau_43_1$, while the weight orbits of $T^*_1$ are realized by 
$\tau_2k_{3,5}$, $\tau_3k_{2,5}$, $\tau_53_1$ and $\tau_{5,2}3_1$ 
\cite{[PS87]}.
Therefore $K$ agrees up to Gluck reconstruction with a branched
twist spin of a torus knot.
\end{proof}

The Gluck reconstruction of a branched twist spin of a classical knot is
another branched twist spin of that knot \cite[\S6]{[Pl84']}.

Elliptic surfaces with $\beta_1=1$ and $\kappa=0$ are
$\mathbb{N}il^3\times\mathbb{E}^1 $-manifolds, and so a knot manifold $M(K)$ 
is homeomorphic to such an elliptic surface if and only if $\pi K$ is 
virtually poly-$Z$ and $\zeta\pi K\cong\mathbb{Z}^2 $. 
For minimal properly elliptic surfaces (those with $\kappa=1$) 
we must settle for a characterization up to $s$-cobordism.

\begin{theorem} 
Let $K$ be a $2$-knot with group $\pi=\pi K$. 
Then $M(K)$ is $s$-cobordant to a minimal properly elliptic surface 
if and only if 
$\zeta\pi\cong\mathbb{Z}^2 $ and $\pi'$ is not virtually poly-$Z$.
\end{theorem}

\begin{proof} If $M(K)$ is a minimal properly elliptic surface then 
it admits a compatible geometry of type $\widetilde{\mathbb{SL}}\times\mathbb{E}^1$ and $\pi$ is 
isomorphic to a discrete cocompact subgroup of 
$Isom_o (\widetilde{\mathbb{SL}})\times\mathbb{R}$, 
the maximal connected subgroup of 
$Isom_o (\widetilde{\mathbb{SL}}\times\mathbb{E}^1)$, 
for the other components consist of orientation reversing 
or antiholomorphic isometries.
(See \cite[Theorem 3.3]{[Wl86]}.)
Since $\pi$ meets
$\zeta{Isom_o (\widetilde{\mathbb{SL}})\times\mathbb{R})}\cong\mathbb{R}^2 $
in a lattice subgroup $\zeta\pi\cong\mathbb{Z}^2 $ and projects nontrivially
onto the second factor,
$\pi'=\pi\cap Isom_o (\widetilde{\mathbb{SL}})$ 
and is the fundamental group of a $\widetilde{\mathbb{SL}}$-manifold. 
Thus the conditions are necessary.

Suppose that they hold. 
Then $M(K)$ is $s$-cobordant to a $\widetilde{\mathbb{SL}}\times\mathbb{E}^1$-manifold
which is the mapping torus $M(\Theta)$ of a 
self homeomorphism of a $\widetilde{\mathbb{SL}}$-manifold,
by Theorem 16.2.
As $\Theta$ must be orientation preserving and induce 
the identity on $\zeta\pi'\cong\mathbb{Z}$ the group $\pi$ is contained in 
$Isom_o (\widetilde{\mathbb{SL}})\times\mathbb{R}$.
Hence $M(\Theta)$ has a compatible structure as 
an elliptic surface \cite[Theorem 3.3]{[Wl86]}.
\end{proof}

An elliptic surface with Euler characteristic 0 is a Seifert fibred 4-manifold,
and so is determined up to diffeomorphism by its fundamental group if the base orbifold is euclidean
or hyperbolic \cite{[Ue90],[Ue91]}.
Using this result (instead of \cite{[Kt75]}) together with 
Theorem 16.6 and Lemma 17.12, 
it may be shown that if $M(K)$ is homeomorphic to a minimal properly elliptic
surface and some power of a weight element is central in $\pi K$ then $M(K)$ 
is homeomorphic to $M(K_1 )$, 
where $K_1 $ is some branched twist spin of a torus knot.
However in general there may be infinitely many algebraically distinct weight
classes in $\pi K$ and we cannot conclude that $K$ is itself such a branched
twist spin.

%% file: m5-18.tex
\chapter{Reflexivity}

The most familiar invariants of knots are derived from the knot complements, 
and so it is natural to ask whether every knot is determined by its complement.
This has been confirmed for classical knots \cite{[GL89]}. 
Given a higher dimensional knot there is at most one other knot
(up to change of orientations) with homeomorphic exterior.
The first examples of non-reflexive 2-knots were given by 
Cappell and Shaneson \cite{[CS76]};
these are fibred with closed fibre $\mathbb{R}^3/\mathbb{Z}^3$. 
Gordon gave a different family of examples \cite{[Go76]}, 
and Plotnick extended his work 
to show that no fibred 2-knot with monodromy of odd order is reflexive.
It is plausible that this may be so whenever the order is greater than 2, 
but this is at present unknown.
                                       
We shall decide the questions of amphicheirality, 
invertibility and reflexivity for most fibred 2-knots with closed fibre 
an $\mathbb{E}^3$- or $\mathbb{N}il^3$-manifold,
excepting only the knots with groups as in part (ii) of Theorem 16.15, 
for which reflexivity remains open.
The other geometrically fibred 2-knots have closed fibre
an $\mathbb{S}^3$-, $\mathbb{H}^3$- or $\widetilde{\mathbb{SL}}$-manifold.
Branched twist spins $\tau_{r,s}k$ of simple 1-knots $k$ 
form an important subclass.
We shall show that such branched twist spins are reflexive if and only if $r=2$.
If, moreover, $k$ is a torus knot then $\tau_{r,s}k$ is $+$amphicheiral but
is not invertible.

This chapter is partly based on joint work with Plotnick and Wilson 
(in \cite{[HP88]} and \cite{[HW89]}, respectively).

\section{Sections of the mapping torus}

Let $\theta$ be a self-homeomorphism of a closed 3-manifold $F$,
with mapping torus $M(\theta)=F\times_\theta{S^1}$,
and canonical projection $p_\theta:M(\theta)\to{S^1}$,
given by $p_\theta([x,s])=e^{2\pi{is}}$ for all $[x,s]\in{M(\theta)}$. 
Then $M(\theta)$ is orientable if and only if $\theta$ is orientation-preserving.
If $\theta'=h\theta{h^{-1}}$ for some self-homeomorphism $h$ of $F$
then $[x,s]\mapsto[h(x),s]$ defines a homeomorphism 
$m(h):M(\theta)\to{M(\theta')}$ such that $p_{\theta'}m(h)=p_\theta$.
Similarly, if $\theta'$ is isotopic to $\theta$ then
$M(\theta')\cong{M(\theta)}$.

If $P\in{F}$ is fixed by $\theta$ then the image 
of $P\times[0,1]$ in $M(\theta)$ is a section of $p_\theta$.
In particular, if the fixed point set of $\theta$ is connected 
there is a canonical isotopy class of sections.
(Two sections are isotopic if and only if 
they represent conjugate elements of $\pi$.)
In general, we may isotope $\theta$ to have a fixed point $P$.
Let $t\in\pi_1(M(\theta))$ correspond to the constant section of
$M(\theta)$, and let $u=gt$ with $g\in\pi_1(F)$.
Let $\gamma:[0,1]\to{F}$ be a loop representing $g$.
There is an isotopy $h_s$ from $h_0=id_F$ to $h=h_1$ which drags $P$
around $\gamma$, so that $h_s(P)=\gamma(s)$ for all $0\leq{s}\leq1$.
Then $H([f,s])=[(h_s)^{-1}(f),s]$ defines a homeomorphism
$M(\theta)\cong{M(h^{-1}\theta)}$.
Under this homeomorphism the constant section of
$p_{h^{-1}\theta}$ corresponds to the section of $p_\theta$
given by $m_u(t)=[\gamma(t),t]$, which represents $u$.
If $F$ is a geometric 3-manifold we may assume that $\gamma$
is a geodesic path.

Suppose henceforth that $F$ is orientable,
$\theta$ is orientation-preserving and fixes a basepoint $P$, 
and induces a meridianal automorphism $\theta_*$ of $\nu=\pi_1(F)$. 
The loop sending $[u]=e^{2\pi iu}$ to $[P,u]$, for all $0\leq u\leq1$, 
is the {\it canonical} cross-section of the mapping torus,
and the corresponding element 
$t\in\pi=\pi_1(M)=\nu\rtimes_{\theta_*}\mathbb{Z}$ 
is a weight element for $\pi$.
Surgery on the image $C=\{P\}\times{S^1}$ of this section gives a 2-knot,
with exterior the complement of an open regular neighbourhood $R$ of $C$.
Choose an embedding $J:D^3\times S^1\to M$ onto $R$.
Let $M_o=M\setminus{intR}$ and let $j=J|_{\partial D^3\times S^1}$.                                                       
Then $\Sigma=M_o\cup_j S^2\times D^2$ and 
$\Sigma_\tau=M_o\cup_{j\tau} S^2\times D^2$
are homotopy 4-spheres and the images of $S^2\times\{0\}$ represent
2-knots $K$ and $K^*$ with group $\pi$ and exterior $M_o$.

Let $\widetilde{F}$ be the universal covering space of $F$, and let 
$\tilde\theta$ be the lift of $\theta$ which fixes some chosen basepoint.
Let $\widehat{M}=\widetilde{N}\times_{\tilde\theta}S^1$ be the (irregular) 
covering space corresponding to the subgroup of $\pi$ generated by $t$. 
This covering space shall serve as a natural model 
for a regular neighbourhood of $C$ in our geometric arguments below.

\section{Reflexivity for fibred 2-knots}

Let $K$ be an $n$-knot with exterior $X$ and group $\pi$.
If it is reflexive there is such a self-homeomorphism which changes
the framing of the normal bundle.
This restricts to a self-homeomorphism of $X$ which (up to changes of orientations) 
is the nontrivial twist $\tau$ on $\partial X\cong{S^n}\times{S^1}$.
(See \S1 of Chapter 14).

If  $K$ is invertible or $\pm$amphicheiral there is a self-homeomorphism $h$ 
of $(S^{n+2},K)$ which changes the orientations appropriately, 
but does not twist the normal bundle of $K(S^n)\subset{S^{n+2}}$. 
Thus if $K$ is $-$amphicheiral there is such an $h$ which 
reverses the orientation of $M(K)$ and for which the induced automorphism
$h'_*$ of $\pi'$ commutes with the meridianal automorphism $c_t$, 
while if $K$ is invertible or $+$ampicheiral 
there is a self-homeomorphism $h$ such that $h'_*c_th'_*=c_t^{-1}$ and which 
preserves or reverses the orientation.
We shall say that $K$ is {\it strongly\/} $\pm$amphicheiral or invertible if 
there is such an $h$ which is an involution.

Suppose now that $n=2$.
Every self-homeomorphism $f$ of $X$ extends ``radially" 
to a self-homeomorphism $h$ of $M(K)=X\cup{D^3\times{S^1}}$ 
which maps the cocore $C=\{0\}\times{S^1}$ to itself.
If $f$ preserves both orientations or reverses both orientations 
then it fixes the meridian, and we may assume that $h|_C=id_C$.
If $f$ reverses the meridian, 
we may still assume that it fixes a point on $C$.
We take such a fixed point as the basepoint for $M(K)$.

Now $\partial X\cong{S^2\times_AS^1}$, where $A$ is
the restriction of the monodromy $\theta$ to 
${\partial (F\setminus{int D^3})}\cong{S^2}$.
(There is an unique isotopy class of homeomorphisms $A$
which are compatible with the orientations of the spheres $S^1$, 
$S^2$ and $X\subset{S^4}$.)
Roughly speaking, the local situation -- 
the behaviour of $f$ and $A$ on $D^3\times S^1$ --
determines the global situation.
Assume that $f$ is a fibre preserving self homeomorphism 
of $D^3\times_A S^1$ which induces a linear map $B$ on each fibre $D^3$. 
If $A$ has infinite order, the question as to when $f$ ``changes the framing", i.e., induces $\tau$ on $\partial D^3\times_A S^1$ is delicate. 
(See \S2 and \S3 below).
But if $A$ has finite order we have the following easy result.
                                           
\begin{lemma} 
Let $A$ in $SO(3)$ be a rotation of order $r\geq2$ and
let $B$ in $O(3)$ be such that $BAB^{-1}=A^{\pm1}$, 
so that $B$ induces a diffeomorphism $f_B$
of $D^3\times_A S^1$. 
If $f_B$ changes the framing then $r=2$.
\end{lemma}
                       
\begin{proof} 
We may choose coordinates for $\mathbb{R}^3$ so that $A=\rho_{s/r}$, 
where $\rho_u$ is the matrix of rotation through $2\pi u$ radians 
about the $z$-axis in $\mathbb{R}^3$, and $0<s<r$.
Let $\rho:D^3\times_A S^1\to D^3\times S^1$ be the diffeomorphism given by
$\rho([x,u])=(\rho_{-su/r},e^{2\pi{i}u})$, 
for all $x\in D^3$ and $0\leq u\leq1$.

If $BA=AB$ then $f_B([x,u]) =[Bx,u]$.
If moreover $r\geq3$ then $B=\rho_v$ for some $v$, 
and so $\rho f_B\rho^{-1}(x,e^{2\pi{i}u})=(Bx,e^{2\pi{i}u})$
does not change the framing.
But if $r=2$ then $A=diag[-1,-1,1]$ and there is more choice for $B$.
In particular, $B=diag[1,-1,1]$ acts dihedrally: 
$\rho_{-u}B\rho_u=\rho_{-2u}B$,
and so $\rho_{-u}f_B\rho_u(x,e^{2\pi{i}u})=(\rho_{-u}x,e^{2\pi{i}u})$, 
i.e., $\rho_{-u}f_B\rho_u$ is the twist $\tau$.
                                    
If $BAB^{-1}=B^{-1}$ then $f_B([x,u])=[Bx,1-u]$. 
In this case             
$\rho f_B\rho^{-1}(x,e^{2\pi{i}u})=
(\rho_{-s(1-u)/r}B\rho_{su/r}x,e^{-2\pi{i}u})$.
If $r\geq3$ then $B$ must act as a reflection in the first two coordinates, 
so $\rho f_B\rho^{-1}(x,e^{2\pi{i}u})=(\rho_{-s/r}Bx,e^{-2\pi{i}u})$ 
does not change the framing.
But if $r=2$ we may take $B=I$, and then
$\rho f_B\rho^{-1}(x,e^{2\pi{i}u})=
(\rho_{(u-1)/2}\rho_{u/2}x,e^{-2\pi{i}u})=
(\rho_{(u-{\frac12})}x,e^{-2\pi{i}u})$, 
which after reversing the $S^1$ factor is just $\tau$.
\end{proof}

We can sometimes show that reflexivity depends only on the knot manifold,
and not the weight orbit.

\begin{lemma}
Let $K$ be a fibred $2$-knot. 
If there is a self homeomorphism $h$ of $X(K)$
which is the identity on one fibre and such that $h|_{\partial{X}}=\tau$ 
then all knots $\tilde{K}$ with $M(\tilde{K})\cong{M(K)}$ are reflexive.
In particular, this is so if the monodromy of $K$ has order $2$. 
\end{lemma}

\begin{proof}
We may extend $h$ to a self-homeomorphism $\hat{h}$ of $M(K)$
which fixes the surgery cocore $C\cong{S^1}$.
After an isotopy of $h$, we may assume that it is the identity on a product
neighbourhood $N=\hat{F}\times[-\epsilon,\epsilon]$ of the closed fibre.
Since any weight element for $\pi$ may be represented by a section
$\gamma$ of the bundle which coincides with $C$ outside $N$,
we may use $h$ to change the framing of the normal bundle 
of $\gamma$ for any such knot.
Hence every such knot is reflexive.

If the monodromy $\theta$ has order 2 the diffeomorphism $h$ of $\hat{F}\times_\theta{S^1}$ 
given by $h([x,s])=[x,1-s]$ which 
``turns the bundle upside down" changes the framing of the normal bundle and fixes one fibre.
\end{proof}

This explains why $r=2$ is special. 
The reflexivity of 2-twist spins is due to Litherland. 
See the footnote to \cite{[Go76]} and also \cite{[Mo83],[Pl84']}. 
                              
The hypotheses in the next lemma seem very stringent, 
but are satisfied in our applications below,
where we shall seek a homeomorphism $\widetilde{F}\cong\mathbb{R}^3$ 
which gives convenient representations of the maps in question, 
and then use an isotopy from the identity to
$\widetilde\theta$ to identify $M(\widetilde\theta)$
with $\mathbb{R}^3\times{S^1}$.
                              
\begin{lemma} 
Suppose that $\widetilde{F}\cong\mathbb{R}^3$ and that
$h$ is an orientation preserving self-homeomorphism of $M$ which
which fixes $C$ pointwise.
If $h_*|_{\pi'}$ is induced by a basepoint preserving
self-homeomorphism $\omega$ of $F$ which commutes with $\theta$
and if there is an isotopy $\gamma$ from $id_{\widetilde{F}}$ to 
$\tilde\theta$ which commutes with the basepoint preserving lift 
$\tilde\omega$ then $h$ does not change the framing.
\end{lemma}

\begin{proof}                                                                    Let $h$ be an orientation preserving self homeomorphism of $M$ 
which fixes $C$ pointwise,
and let $R$ be a regular neighbourhood of $C$.
Suppose that $h$ changes the framing.
We may assume that $h|_R$ is a bundle automorphism
which agrees with the radial extension of $\tau$ from 
$\partial{R}=S^2\times S^1$ to $R$. 
Since $h$ fixes the meridian, $h_*\theta_*=\theta_*h_*$.
Let $\omega$ be a basepoint preserving self diffeomorphism of $F$ 
which induces $h_*|_\nu$ and commutes with $\theta$. 
Then we may define a self diffeomorphism $h_\omega$ of $M$ by
$h_\omega([x,s])=[\omega(x),s]$ for all $[x,s]$ in $M=F\times_\theta S^1$.

Since $h_{\omega*}=h_*$ and $M$ is aspherical, 
$h$ and $h_\omega$ are homotopic.
Therefore the lifts $\hat h$ and $\hat h_\omega$ 
to basepoint preserving maps of $\widehat{M}$ are properly homotopic. 
Let $\tilde\omega$ be the lift of $\omega$ to a basepoint preserving
map of $\widetilde{F}$. 
Note that $\tilde\omega$ is orientation preserving, and so is
isotopic to $id_{\widetilde{F}}$.

Given an isotopy $\gamma$ from $\gamma(0)=id_{\widetilde{F}}$ to 
$\gamma(1)=\tilde\theta$ we may define a diffeomorphism 
$\rho_\gamma:\widetilde{F}\times{S^1}\to\widehat{M}$ by                 
$\rho_\gamma(x,e^{2\pi it})=[\gamma(t)(x),t]$.
Now $\rho_\gamma^{-1}\hat h_\omega\rho_\gamma(l,[u])=
(\gamma(u)^{-1}\tilde\omega\gamma(u)(l),[u])$.
Thus if $\gamma(t)\tilde\omega=\tilde\omega\gamma(t)$ 
for all $t$ then                             
$\rho_\gamma^{-1}\hat h_\omega\rho_\gamma=\tilde\omega\times id_{S^1}$,
and so $\hat h$ is properly homotopic to $id_{\widehat M}$. 

Since the radial extension of $\tau$ and $\rho_\gamma^{-1}\hat h\rho_\gamma$ 
agree on $D^3\times S^1$ 
they are properly homotopic on $\mathbb{R}^3\times S^1$ 
and so $\tau$ is properly homotopic to the identity.
Now $\tau$ extends uniquely to a self diffeomorphism $\tau$ of $S^3\times S^1$, 
and any such proper homotopy extends to a homotopy from $\tau$ to the identity.
Let $p$ be the projection of $S^3\times S^1$ onto $S^3$. 
The suspension of $p\tau$,
restricted to the top cell of $\Sigma(S^3\times S^1)=S^2\vee S^4\vee S^5$ 
is the nontrivial element
of $\pi_5(S^4)$, whereas the corresponding restriction 
of the suspension of $p$ is trivial. 
(See \cite{[CS76],[Go76]}).
This is contradicts $p\tau\sim{p}$.
Therefore $h$ cannot change the framing.
\end{proof}

Note that in general there is no isotopy from $id_F$ to $\theta$.

We may use a similar argument to give a sufficient condition 
for knots constructed from mapping tori to be $-$amphicheiral. 
As we shall not use this result below we shall only
sketch a proof.                                                  
                                           
\begin{lemma} 
Let $F$ be a closed orientable $3$-manifold with universal cover 
$\widetilde{F}\cong\mathbb{R}^3$.
Suppose now that there is an orientation reversing self diffeomorphism 
$\psi:F\to{F}$ which commutes with $\theta$ and which fixes $P$.
If there is a path $\gamma$ from $I$ to $\Theta=D\theta(P)$ which commutes 
with $\Psi=D\psi(P)$ then each of $K$ and $K^*$ is $-$amphicheiral.
\end{lemma}

\begin{proof} The map $\psi$ induces an orientation reversing 
self diffeomorphism of $M$ which fixes $C$ pointwise. 
We may use such a path $\gamma$ to define a diffeomorphism 
$\rho_\gamma:\widetilde{F}\times S^1\to \widetilde{M}$.
We may then verify that $\rho_\gamma^{-1}\hat h\rho_\gamma$ 
is isotopic to $\Psi\times id_{S^1}$,
and so $\rho_\gamma^{-1}\hat h\rho_\gamma|_{\partial D^3\times S^1}$ 
extends across $S^2\times D^2$.
\end{proof}

\section{Cappell-Shaneson knots}

Let $A\in SL(3,\mathbb{Z})$ be such that $\det(A-I)=\pm 1$.
The knots $K$ obtained by surgery on the canonical cross-section of 
the mapping torus of the corresponding
automorphism of the 3-torus $T^3$
are the Cappell-Shaneson knots determined by $A$.
On replacing $A$ by $A^{-1}$
(which corresponds to changing the orientation of the knot), 
if necessary, we may assume that $\det(A-I)=+1$.
Inversion in each fibre of $M(K)$ fixes a circle, 
and passes to an orientation reversing involution of $(S^4,K)$.
Hence each such knot is strongly $-$amphicheiral.
However, it is not invertible, since $\Delta_K(t)=\det(tI-A)$ is not symmetric.

Cappell and Shaneson showed that if $A$ has no negative eigenvalue 
then the associated knots are not reflexive.
In a footnote they observed that if $A$ has negative eigenvalues
the two knots obtained from $A$ are equivalent if and only if 
there is a matrix $B$ in $GL(3,\mathbb{Z})$ such that $AB=BA$ 
and the restriction of $B$ to the 
negative eigenspace of $A$ has negative determinant. 
We shall translate this criterion into one 
involving algebraic numbers, 
and 
show that up to change of orientations 
there is just one reflexive Cappell-Shaneson 2-knot.

\begin{theorem}
Let $A\in SL(3,\mathbb{Z})$ satisfy $\det(A-I)=1$.
If $A$ has trace $-1$ then the corresponding Cappell-Shaneson knot is reflexive,
and is determined up to change of orientations among all $2$-knots 
with metabelian group by its Alexander polynomial $t^3+t^2-2t-1$.
If the trace of $A$ is not $-1$ then the corresponding Cappell-Shaneson 
knots are not reflexive.
\end{theorem}

\begin{proof} Let $a$ be the trace of $A$. 
Then the characteristic polynomial of $A$ 
is $\Delta_a(t)=t^3-at^2+(a-1)t-1=t(t-1)(t-a+1)-1$.
It is easy to see that $\Delta_a$ is irreducible; 
indeed, it is irreducible modulo $(2)$.
Since the leading coefficient of $\Delta_a$ is positive and $\Delta_a(1)<0$ 
there is at least one positive eigenvalue.
If $a>5$ all three eigenvalues are positive (since $\Delta_a(0)=-1$, 
$\Delta_a(\frac12)=(2a-11)/8>0$ and $\Delta_a(1)=-1$).
If $0\leq a\leq5$ there is a pair of complex eigenvalues.

Thus if $a\geq 0$ there are no negative eigenvalues, and so
$\gamma(s)=sA+(1-s)I$ (for $0\leq{s}\leq1$) defines an isotopy 
from $I$ to $A$ in $GL(3,\mathbb{R})$.
Let $h$ be a self homeomorphism of $(M,C)$ such that $h(*)=*$. 
We may assume that $h$ is orientation preserving and preserves
the meridian.
Since $M$ is aspherical $h$ is homotopic to a map $h_B$, 
where $B\in SL(3,\mathbb{Z})$ commutes with $A$.
Applying Lemma 18.3 with $\tilde\theta=A$ and $\tilde\omega=B$, 
we see that $h$ must preserve the framing and so $K$ is not reflexive.

We may assume henceforth that $a<0$.
There are then three real roots $\lambda_i$, for $1\leq i\leq3$, 
such that $a-1<\lambda_3<a<\lambda_2<0<1<\lambda_1<2$.
Note that the products  $\lambda_i(\lambda_i-1)$ are all positive, 
for $1\leq i\leq3$.

Since the eigenvalues of $A$ are real and distinct there is a matrix
$P$ in $GL(3,\mathbb{R})$ such that $\tilde A=PAP^{-1}$ 
is the diagonal matrix $diag[\lambda_1,\lambda_2,\lambda_3]$.
If $B$ in $GL(3,\mathbb{Z})$ commutes with $A$ then 
$\tilde B=PBP^{-1}$ commutes with $\tilde A$ 
and hence is also diagonal (as the $\lambda_i$ are distinct).
On replacing $B$ by $-B$ if necessary we may assume that $\det(B)=+1$.
Suppose that $\tilde B=diag[\beta_1,\beta_2,\beta_3]$.
We may isotope $PAP^{-1}$ linearly to $diag[1,-1,-1]$.
If $\beta_2\beta_3>0$ for all such $B$ then $PBP^{-1}$ is isotopic to $I$ 
through block diagonal matrices and we may again conclude that the 
knot is not reflexive.
On the other hand if there is such a $B$ with $\beta_2\beta_3<0$ 
then the knot is reflexive.
Since $\det(B)=+1$ an equivalent criterion for reflexivity is that $\beta_1<0$.

The roots of $\Delta_{-1}(t)=t^3+t^2-2t-1$ are the Galois conjugates
of $\zeta_7+\zeta_7^{-1}$, where $\zeta_7$ is a primitive 7th root of unity.
The discriminant of $\Delta_{-1}(t)$ is 49,
which is not properly divisible by a perfect square,
and so $\mathbb{Z}[t]/(\Delta_{-1}(t))$ is the full ring of integers in
the maximal totally real subfield of the cyclotomic field $\mathbb{Q}(\zeta_7)$.
As this ring has class number 1 it is a PID.
Hence any two matrices in $SL(3,\mathbb{Z})$ 
with this characteristic polynomial are conjugate, by Theorem 1.4.
Therefore the knot group is unique and determines $K$ 
up to Gluck reconstruction and change of orientations, by Theorem 17.5.
Since $B=-A-I$ has determinant 1 and $\beta_1=-\lambda_1-1<0$, 
the corresponding knot is reflexive.

Suppose now that $a<-1$.
Let $F$ be the field $\mathbb{Q}[t]/(\Delta_a(t))$ and let $\lambda$ 
be the image of $X$ in $F$.
We may view $\mathbb{Q}^3$ as a $\mathbb{Q}[t]$-module 
and hence as a 1-dimensional $F$-vector space via the action of $A$.                                                            
If $B$ commutes with $A$ then it induces an automorphism of this vector space 
which preserves a lattice and so determines a unit $u(B)$ in $O_F$,
the ring of integers in $F$.
Moreover $\det(B)=N_{F/\mathbb{Q}}u(B)$.
If $\sigma$ is the embedding of $F$ in $R$ which sends $\lambda$ to 
$\lambda_1$ and $P$ and $B$ are as above we must have $\sigma(u(B))=\beta_1$.

Let $U=O_F^\times$ be the group of all units in $O_F$, and let $U^\nu$, 
$U^\sigma$, $U^+$ and $U^2$ be the subgroups of units of norm 1, units whose 
image under $\sigma$ is positive, totally positive units and squares, 
respectively.
Then $U\cong\mathbb{Z}^2\times\{\pm1\}$, 
since $F$ is a totally real cubic number field,
and so $[U:U^2]=8$.
The unit $-1$ has norm $-1$, and $\lambda$ is a unit of norm 1 in $U^\sigma$ 
which is not totally positive.
Hence $[U:U^\nu]=[U^\nu\cap U^\sigma:U^+]=2$.
It is now easy to see that there is a unit of norm 1 that is not in $U^\sigma$ 
(i.e., $U^\nu\not= U^\nu\cap U^\sigma$) if and only if
every totally positive unit is a square (i.e., $U^+=U^2$).

The image of $t(t-1)$ in $F$ is $\lambda(\lambda-1)$, 
which is totally positive and is a unit (since $t(t-1)(t-a+1)=1+\Delta_a(t)$).
Suppose that it is a square in $F$.
Then $\phi=\lambda-(a-1)$ is a square (since 
$\lambda(\lambda-1)(\lambda-(a-1))=1)$.    
The minimal polynomial of $\phi$ is $g(Y)=Y^3+(2a-3)Y^2+(a^2-3a+2)Y-1$.
If $\phi=\psi^2$ for some $\psi$ in $F$ then $\psi$ is a root of $h(Z)=g(Z^2)$
and so the minimal polynomial of $\psi$ divides $h$.
This polynomial has degree 3 also, since $\mathbb{Q}(\psi)=F$, 
and so $h(Z)=p(Z)q(Z)$ for some
polynomials $p(Z)=Z^3+rZ^2+sZ+1$ and $q(Z)=Z^3+r'Z^2+s'Z-1$ with integer coefficients.
Since the coefficients of $Z$ and $Z^5$ in $h$ are 0 we must have $r'=-r$ and $s'=-s$.
Comparing the coefficients of $Z^2$ and $Z^4$ then gives the equations
$2s-r^2=2a-3$ and $s^2-2r=a^2-3a+2$. Eliminating $s$ we find that
$r(r^3+(4a-6)r-8)=-1$ and so $1/r$ is an integer. 
Hence $r=\pm1$ and so $a=-1$ or 3, contrary to hypothesis.                                                            
Thus there is no such matrix $B$ and so the Cappell-Shaneson knots 
corresponding to $A$ are not reflexive.
\end{proof}

\section{The Hantzsche-Wendt flat 3-manifold}

In \S2 of Chapter 8 we determined the automorphisms 
of the flat 3-manifold groups from a purely algebraic point of view.
In this chapter we use instead the geometric approach,
to study fibred 2-knots with closed fibre
the {\it Hantzsche-Wendt} flat 3-manifold,
$HW={G_6}\backslash\mathbb{R}^3$.

The group of affine motions of 3-space is 
$Aff(3)=\mathbb{R}^3\rtimes{GL(3,\mathbb{R})}$.
The action is given by $(v,A)(x)=Ax+v$, for all $x\in\mathbb{R}^3$.
Therefore $(v,A)(w,B)=(v+Aw,AB)$.
Let $\{e_1,e_2,e_3\}$ be the standard basis of $\mathbb{R}^3$,
and let $X,Y,Z\in{GL(3,\mathbb{Z})}$ be the diagonal matrices
$X=diag[1,-1,-1]$, $Y=diag[-1,1,-1]$ and $Z=diag[-1,-1,1]$.
Let $x=(\frac12e_1,X)$, $y=(\frac12(e_2-e_3),Y)$ and 
$z={(\frac12(e_1-e_2+e_3),Z)}$.
Then the subgroup of $Aff(3)$ generated by $x$ and $y$ is $G_6$.
The translation subgroup $T=G_6\cap\mathbb{R}^3$ is free abelian,
with basis $\{x^2,y^2,z^2\}$.
The holonomy group $H=\{I,X,Y,Z\}\cong(Z/2Z)^2$
is the image of $G_6$ in $GL(3,\mathbb{R})$.
(We may clearly take $\{1,x,y,z\}$ as
coset representatives for $H$ in $G_6$.
The commutator subgroup $G_6'$ is free abelian, 
with basis $\{x^4,y^4,x^2y^2z^{-2}\}$.
Thus $2T<G_6'<T$, $T/G_6'\cong(Z/2Z)^2$ and $G_6'/2T\cong{Z/2Z}$.

Automorphisms of the fundamental group of a flat $n$-manifold 
are induced by conjugation in $Aff(n)$,
by a theorem of Bieberbach.
In particular, $Aut(G_6)\cong{N/C}$ and $Out(G_6)\cong{N/CG_6}$,
where $C=C_{Aff(3)}(G_6)$ and $N=N_{Aff(3)}(G_6)$.

If $(v,A)\in{Aff(3)}$ commutes with all elements of $G_6$
then $AB=BA $ for all $B\in{H}$, so $A$ is diagonal, and 
$v+Aw=w+Bv$ for all $(w,B)\in{G_6}$.
Taking $B=I$, we see that $Aw=w$ for all $w\in\mathbb{Z}^3$, so $A=I$,
and then $v=Bv$ for all $B\in{H}$, so $v=0$.
Thus $C=1$, and so $Aut(G_6)\cong{N}$.

If $(v,A)\in{N}$ then $A\in{N_{GL(3,\mathbb{R})}(H)}$ 
and $A$ preserves $T=\mathbb{Z}^3$, so $A\in{N_{GL(3,\mathbb{Z})}(H)}$.
Therefore $W=AXA^{-1}$ is in $H$.
Hence $WA=AX$ and so $WAe_1=Ae_1$ is 
up to sign the unique basis vector fixed by $W$.
Applying the same argument to $AYA^{-1}$ and $AZA^{-1}$,
we see that $N_{GL(3,\mathbb{R})}(H)$ is 
the group of ``signed permutation matrices",
generated by the diagonal matrices and permutation matrices.
Let
\[ P=
\left(
\begin{matrix}
0& 1 & 0\\
1 & 0 &0\\
0& 0& -1
\end{matrix}
\right)
\quad \mathrm{and}\quad 
J=
\left(
\begin{matrix}
0& 1 &0\\
0&0&-1\\
1& 0& 0
\end{matrix}
\right).
\]
If $A$ is a diagonal matrix in $GL(3,\mathbb{Z})$ then $(0,A)\in{N}$.
Thus $\tilde{a}=(0,-X)$, $\tilde{b}=(0,-Y)$ and $\tilde{c}=(0,-Z)$ are in $N$.
It is easily seen that $N\cap\mathbb{R}^3=\frac12\mathbb{Z}^3$,
with basis $\tilde{d}=(\frac12{e_1},I)$, $\tilde{e}=(\frac12{e_2},I)$ and 
$\tilde{f}=(\frac12{e_3},I)$.
It is also easily verified that $\tilde{i}=(-\frac14e_3,P)$ 
and $\tilde{j}=(\frac14(e_1-e_2),J)$ are in $N$,
and that $N$ is generated by 
$\{\tilde{a},\tilde{b},\tilde{c},\tilde{d},\tilde{e},
\tilde{f},\tilde{i},\tilde{j}\}$.
These generators correspond to the generators given in \S2 of Chapter 8,
and we shall henceforth drop the tildes.

The natural action of $N$ on $\mathbb{R}^3$ is isometric,
since ${N_{GL(3,\mathbb{R})}(H)<O(3)}$,
and so $N/G_6$ acts isometrically on the orbit space $HW$.
In fact every isometry of $HW$ lifts to an affine transformation of
$\mathbb{R}^3$ which normalizes $G_6$, and so $Isom(HW)\cong{Out(G_6)}$.
The isometries which preserve the orientation
are represented by pairs $(v,A)$ with $\det(A)=1$. 

Since $G_6$ is solvable and $H_1(G_6)\cong(Z/4Z)^2$,
an automorphism $(v,A)$ of $G_6$ is meridianal 
if and only if its image in $Aut(G_6/T)\cong{GL(2,2)}$ has order 3.
Thus its image in $Out(G_6)$ is conjugate to $[j]$, $[j]^{-1}$, 
$[ja]$ or $[jb]$.
The latter pair are orientation-preserving 
and each is conjugate to its inverse (via $[i]$).
However $(ja)^3=1$ while $(jb)^3=de^{-1}f$, so $[jb]^3=[ab]\not=1$.
Thus $[ja]$ is not conjugate to $[jb]^\pm$, and the knot groups
$G(+)=G_6\rtimes_{[ja]}\mathbb{Z}$ and 
$G(-)=G_6\rtimes_{[jb]}\mathbb{Z}$ are distinct.
The corresponding knot manifolds are the mapping tori 
of the isometries of $HW$ determined by $[ja]$ and $[jb]$, 
and are flat 4-manifolds.

\section{2-knots with group $G(\pm)$}

We shall use the geometry of flat 3-manifolds to
show that 2-knots with group $G(+)$ or $G(-)$ are not reflexive,
and that among these knots only $\tau_34_1, \tau_34_1^*$ 
and the knots obtained by surgery on the section of $M([jb])$ 
defined by $\gamma|_{[0,1]}$ admit orientation-changing symmetries.

\begin{theorem}
Let $K$ be a $2$-knot with group $\pi\cong{G(+)}$ or $G(-)$.
Then $K$ is not reflexive.
\end{theorem}

\begin{proof}
The knot manifold $M=M(K)$ is homeomorphic to
the flat 4-manifold $\mathbb{R}^4/\pi$,
by Theorem 8.1 and the discussion in Chapter 16.
The weight orbit of $K$ may be represented by a geodesic 
simple closed curve $C$ through the basepoint $P$ of $M$.
Let $\gamma$ be the image of $C$ in $\pi$.
Let $h$ be a self-homeomorphism of $M$ which fixes $C$ pointwise.
Then $h$ is based-homotopic to an affine diffeomorphism 
$\alpha$, and then $\alpha_*(\gamma)=h_*(\gamma)=\gamma$.
Let $\widehat{M}\cong\mathbb{R}^3\times{S^1}$ be the covering space 
corresponding to the subgroup $\langle\gamma\rangle\cong\mathbb{Z}$,
and fix a lift $\widehat{C}$.
A homotopy from $h$ to $\alpha$ lifts to a proper homotopy
between the lifts $\widehat{h}$ and $\widehat{\alpha}$ to 
self-homeomorphisms fixing $\widehat{C}$.
Now the behaviour at $\infty$ of these maps is determined by
the behaviour near the fixed point sets, as in Lemma 18.3.
Since the affine diffeomorphism $\widehat{\alpha}$ does not change
the framing of the normal to $\widehat{C}$ it follows that
$\widehat{h}$ and $h$ do not change the normal framings either.
\end{proof}

The orthogonal matrix $-JX$ is a rotation though $\frac{2\pi}3$
about the axis in the direction $e_1+e_2-e_3$.
The fixed point set of the isometry $[ja]$ of $HW$ is
the image of the line $\lambda(s)=s(e_1+e_2-e_3)-\frac14e_2$.
The knots corresponding to the canonical section are $\tau_34_1$ 
and its Gluck reconstruction $\tau_34_1^*$.

\begin{lemma}
If $n=0$ then $C_{Aut(G_6)}(ja)$ and $N_{Aut(G_6)}(\langle{ja}\rangle)$
are generated by $\{ja,def^{-1},abce\}$ and $\{ja,ice,abce\}$,
respectively.
The subgroup which preserves the orientation of $\mathbb{R}^3$ 
is generated by $\{ja,ice\}$.

If $n\not=0$ then $N_{Aut(G_6)}(\langle{d^{2n}ja}\rangle)=
C_{Aut(G_6)}(d^{2n}ja)$, and is generated by 

\noindent$\{d^{2n}ja,def^{-1}\}$.
This subgroup acts orientably on $\mathbb{R}^3$.
\end{lemma}

\begin{proof}
This is straightforward.
(Note that $abce=j^3$ and $def^{-1}=(ice)^2$.)
\end{proof}

\begin{lemma}
The mapping torus $M([ja])$ has an orientation reversing involution 
which fixes a canonical section pointwise,
and an orientation reversing involution which fixes a
canonical section setwise but reverses its orientation.
There is no orientation preserving involution of $M$ which reverses 
the orientation of any section.
\end{lemma}

\begin{proof}
Let $\omega=abcd^{-1}f=abce(ice)^{-2}$\!,
and let $p=\lambda(\frac14)=\frac14(e_1-e_3)$.
Then $\omega=(2p,-I_3)$, $\omega^2=1$,
$\omega{ja}=ja\omega$ and $\omega(p)=ja(p)=p$.
Hence $\Omega=m([\omega])$ is an orientation reversing involution 
of $M([ja])$ which fixes the canonical section determined by 
the image of $p$ in $HW$.

Let $\Psi([f,s])=[[iab](f),1-s]$ for all $[f,s]\in{M([ja])}$.
This is well-defined, since $(iab){ja}(iab)^{-1}=(ja)^{-1}$, 
and is an involution, since $(iab)^2=1$.
It is clearly orientation reversing, 
and since $iab(\lambda(\frac18))=\lambda(\frac18)$ 
it reverses the section determined by the image 
of $\lambda(\frac18)$ in $HW$.

On the other hand, 
$\langle{ja},ice\rangle\cong{Z/3Z}\rtimes_{-1}\mathbb{Z}$,
and the elements of finite order in this group do not invert $ja$.
\end{proof}

\begin{theorem}
Let $K$ be a $2$-knot with group $G(+)$ and weight element $u=x^{2n}t$,
where $t$ is the canonical section.
If $n=0$ then $K$ is strongly $\pm$amphicheiral,
but is not strongly invertible.
If $n\not=0$ then $K$ is neither amphicheiral nor invertible.
\end{theorem}

\begin{proof}
Suppose first that $n=0$.
Since $-JX$ has order 3 it is conjugate in $GL(3,\mathbb{R})$
to a block diagonal matrix $\Lambda(-JX)\Lambda^{-1}=\left(\smallmatrix 1&0\\
0&R(\frac{2\pi}3)\endsmallmatrix\right)$,
where $R(\theta)\in{GL(2,\mathbb{R})}$ is rotation through $\theta$.
Let $R_s=R(\frac{2\pi}3s)$ and 
$\xi(s)=((I_3-A_s)p,A_s)$, where 
$A_s= \Lambda^{-1}\left(\smallmatrix 1&0\\
0&R_s\endsmallmatrix\right)\Lambda$, for $s\in\mathbb{R}$.
Then $\xi$ is a 1-parameter subgroup of $Aff(3)$,
such that $\xi(s)(p)=p$ and $\xi(s)\omega=\omega\xi(s)$ for all $s$.
In particular, $\xi|_{[0,1]}$ is a path from $\xi(0)=1$ 
to $\xi(1)=ja$ in $Aff(3)$.
Let $\Xi:\mathbb{R}^3\times{S^1}\to{M(ja)}$ be the homeomorphism given by
$\Xi(v,e^{2\pi{is}})=[\xi(v),s]$ for all $(v,s)\in \mathbb{R}^3\times[0,1]$.
Then $\Xi^{-1}\Omega\Xi=\omega\times{id_{S^1}}$ and so
$\Omega$ does not change the framing.
Therefore $K$ is strongly $-$amphicheiral.
 
Similarly, if we let $\zeta(s)=((I_3-A_s)\lambda(\frac18),A_s)$
then $\zeta(s)(\lambda(\frac18))=\lambda(\frac18)$ and
$iab\zeta(s)iab=\zeta(s)^{-1}$, for all $s\in\mathbb{R}$,
and $\zeta|_{[0,1]}$ is a path from 1 to $ja$ in $Aff(3)$.
Let $Z:\mathbb{R}^3\times{S^1}\to{M(ja)}$ be the homeomorphism given by
$Z(v,e^{2\pi{is}})=[\zeta(s)(v),s]$ for all $(v,s)\in \mathbb{R}^3\times{S^1}$.
Then 
\[Z^{-1}\Psi{Z}(v,z)=(\zeta(1-s)^{-1}iab\zeta(s)(v),z^{-1})
=(\zeta(1-s)^{-1}\zeta(s)^{-1}iab(v),z^{-1})\]
\[=((ja)^{-1}iab(v),z^{-1})\]
for all $(v,z)\in \mathbb{R}^3\times{S^1}$.
Hence $\Psi$ does not change the framing, 
and so $K$ is strongly $+$amphicheiral.
However it is not strongly invertible, by Lemma 18.15.

If $n\not=0$ every self-homeomorphism $h$ of $M(K)$
preserves the orientation and fixes the meridian, by Lemma 18.7,
and so $K$ is neither amphicheiral nor invertible.
\end{proof}

A similar analysis applies when the knot group is $G(-)$,
i.e., when the meridianal automorphism is
$jb=(\frac14(e_1-e_2),-JY)$.
(The orthogonal matrix $-JY$ is now a rotation though $\frac{2\pi}3$
about the axis in the direction ${e_1-e_2+e_3}$.)
All 2-knots with group $G(-)$ are fibred,
and the characteristic map $[jb]$ has finite order, 
but none of these knots are twist-spins, by Corollary 16.12.1.

\begin{lemma}
If $n=0$ then $C_{Aut(G_6)}(jb)$ and $N_{Aut(G_6)}(\langle{jb}\rangle)$
are generated by $\{jb\}$ and $\{jb,i\}$,
respectively.
If $n\not=0$ then 
$N_{Aut(G_6)}(\langle{d^{2n}jb}\rangle)=C_{Aut(G_6)}(d^{2n}jb)$ 
and is generated by $\{d^{2n}jb,de^{-1}f\}$.
These subgroups act orientably on $\mathbb{R}^3$.
\hfill $\Box$
\end{lemma}

The isometry $[jb]$ has no fixed points in $G_6\backslash\mathbb{R}^3$.
We shall defined a preferred section as follows.
Let $\gamma(s)=\frac{2s-1}8(e_1-e_2)-\frac18e_3$, for $s\in\mathbb{R}$.
Then $\gamma(1)=jb(\gamma(0))$, 
and so $\gamma|_{[0,1]}$ defines a section of $p_{[jb]}$.
We shall let the image of $(\gamma(0),0)$ be the basepoint for $M([jb])$.

\begin{theorem}
Let $K$ be a $2$-knot with group $G(-)$ and weight element $u=x^{2n}t$,
where $t$ is the canonical section.
If $n=0$ then $K$ is strongly $+$amphicheiral but not invertible.
If $n\not=0$ then $K$ is neither amphicheiral nor invertible.
\end{theorem}

\begin{proof}
Suppose first that $n=0$.
Since $i(\gamma(s))=\gamma(1-s)$ for all $s\in\mathbb{R}$
the section defined by $\gamma|_{[0,1]}$ is fixed setwise 
and reversed by the orientation reversing involution 
${[f,s]\mapsto[[i](f),1-s]}$.
Let $B_s$ be a 1-parameter subgroup of $O(3)$ such that $B_1=jb$.
Then we may define a path from 1 to $jb$ in $Aff(3)$ by setting
$\zeta(s)={((I_3-B_s)\gamma(s),B_s)}$ for $s\in\mathbb{R}$.
We see that $\zeta(0)=1$, $\zeta(1)=jb$,
$\zeta(s)(\gamma(s))=\gamma(s)$ 
and $i\zeta(s)i=\zeta(s)^{-1}$ for all $0\leq{s}\leq1$.
As in Theorem 18.9 it follows that the involution does not change the framing
and so $K$ is strongly $+$amphicheiral.

The other assertions follow from Lemma 18.10, as in Theorem 18.9.
\end{proof}

\section{$\mathbb{N}il^3$-fibred knots}

The 2-knot groups $\pi$ with $\pi'$ a $\mathbb{N}il^3$-lattice
were determined in Theorems 14 and 15 of Chapter 16.
The group $Nil=Nil^3$ is a subgroup of $SL(3,\mathbb{R})$
and is diffeomorphic to $R^3$, with multiplication
given by 
\[
[r,s,t][r',s',t']=[r+r',s+s',rs'+t+t'].
\]
(See Chapter 7). 
The kernel of the natural homomorphism from $Aut_{Lie}(Nil)$ to 
$Aut_{Lie}(R^2)=GL(2,\mathbb{R})$ induced by abelianization 
($Nil/Nil'\cong\mathbb{R}^2$) is isomorphic to 
$Hom_{Lie}(Nil,\zeta Nil)\cong\mathbb{R}^2$.
The set underlying the group $Aut_{Lie}(Nil)$ is the cartesian product
$GL(2,\mathbb{R})\times\mathbb{R}^2$, with 
$(A,\mu)= (\left(\smallmatrix a& c\\
          b& d\endsmallmatrix\right),(m_1,m_2))$ 
acting via $(A,\mu)([r,s,t])=$
\begin{equation*}
[ar+cs,br+ds,m_1r+m_2s+(ad-bc)t+bcrs+\frac{ab}2r(r-1)+\frac{cd}2s(s-1)].
\end{equation*}
The Jacobian of such an automorphism is $(ad-bc)^2$, 
and so it is orientation preserving.
The product of $(A,\mu)$ with
$(B,\nu)= (\left(\smallmatrix g&j\\
          h&k\endsmallmatrix\right),(n_1,n_2))$ is
\[
(A,\mu)\circ(B,\nu)=(AB,\mu B+\det(A)\nu+\frac12\eta(A,B)),
\]
where
\[
\eta(A,B)= 
(abg(1-g)+cdh(1-h)-2bcgh,abj(1-j)+cdk(1-k)-2bcjk).
\]
In particular, $Aut_{Lie}(Nil)$ is {\it not} a semidirect product of 
$GL(2,\mathbb{R})$ with $\mathbb{R}^2$.
For each $q>0$ in $\mathbb{Z}$ the stabilizer of $\Gamma_q$ 
in $Aut_{Lie}(Nil)$ is ${GL(2,\mathbb{Z})\times(q^{-1}\mathbb{Z}^2)}$, 
which is easily verified to be $Aut(\Gamma_q)$.
(See \S7 of Chapter 8).
Thus every automorphism of $\Gamma_q$ extends to an automorphism of $Nil$.
(This is a special case of a theorem of Malcev on embeddings of torsion-free 
nilpotent groups in 1-connected nilpotent Lie groups - see \cite{[Rg]}).

Let the identity element $[0,0,0]$ and its images in $N_q=Nil/\Gamma_q$
be the basepoints for $Nil$ and for these coset spaces.
The extension of each automorphism of $\Gamma_q$ to $Nil$
induces a basepoint and orientation preserving self homeomorphism of $N_q$.

If $K$ is a 2-knot with group $\pi=\pi K$ and $\pi'\cong\Gamma_q$ 
then $M=M(K)$ is 
homeomorphic to the mapping torus of such a self homeomorphism of $N_q$.
(In fact, such mapping tori are determined up to diffeomorphism by their fundamental groups).
Up to conjugacy and involution there are just three classes of meridianal automorphisms 
of $\Gamma_1$ and one of $\Gamma_q$, for each odd $q>1$. (See Theorem 16.13).
Since $\pi''\leq\zeta\pi'$ it is easily seen that $\pi$ has just two strict
weight orbits.
Hence $K$ is determined up to Gluck reconstruction and changes of orientation 
by $\pi$ alone, by Theorem 17.4.
(Instead of appealing to 4-dimensional surgery to realize automorphisms of
$\pi$ by basepoint and orientation preserving self homeomorphisms of $M$
we may use the $S^1$-action on $N_q$ to construct such a self homeomorphism
which in addition preserves the fibration over $S^1$).
We shall show that the knots with $\pi'\cong\Gamma_1$ and whose characteristic
polynomials are $X^2-X+1$ and $X^2-3X+1$ are not reflexive, 
while all other 2-knots with $\pi'\cong\Gamma_q$ for some odd $q\geq1$ are reflexive.

\begin{lemma} 
Let $K$ be a fibred $2$-knot with closed fibre $N_1$ and 
Alexander polynomial $X^2-3X+1$. Then $K$ is $+$amphicheiral.
\end{lemma}

\begin{proof}  
Let $\Theta=(A,(0,0))$ be the automorphism of $\Gamma_1$ with 
$A=\left(\smallmatrix 1&1\\
          1&2\endsmallmatrix\right)$.
Then $\Theta$ induces a basepoint and orientation preserving 
self diffeomorphism $\theta$ of $N_1$.
Let $M=N_1\times_\theta S^1$ and let $C$ be the canonical section.                         
A basepoint and orientation preserving self diffeomorphism $\psi$ of $N_1$ 
such that $\psi\theta\psi^{-1}=\theta^{-1}$ induces a self diffeomorphism
of $M$ which reverses the orientations of $M$ and $C$. 
If moreover it does not twist the normal bundle of $C$ 
then each of the 2-knots $K$ and $K^*$ obtained by surgery on $C$ 
is $+$amphicheiral.
We may check the normal bundle condition by using an isotopy from $\Theta$ 
to $id_{Nil}$ to identify $\widehat M$ with $Nil\times S^1$.                                   
Thus we seek an automorphism $\Psi=(B,\mu)$ of $\Gamma_1$ 
such that $\Psi\Theta_t\Psi^{-1}=\Theta_t^{-1}$, 
or equivalently $\Theta_t\Psi\Theta_t=\Psi$, 
for some isotopy $\Theta_t$ from $\Theta_0=id_{Nil}$ 
to $\Theta_1=\Theta$.

Let $P=\left(\smallmatrix 0&-1\\
          1&0\endsmallmatrix\right)$.
Then $PAP^{-1}=A^{-1}$, or $APA=P$.
It may be checked that the equation $\Theta(P,\mu)\Theta=(P,\mu)$ reduces to a linear equation for 
$\mu$ with unique solution $\mu=-(2,3)$.
Let $\Psi=(P,-(2,3))$ and let $h$ be the induced diffeomorphism of $M$.

As the eigenvalues of $A$ are both positive it lies on a 1-parameter subgroup,
determined by 
$L=\ln(A)=m\left(\smallmatrix 1&-2\\
          -2&-1\endsmallmatrix\right)$,
where $m=(\ln((3+\sqrt5)/2))/\sqrt5$.
Now $PLP^{-1}=-L$ and so $P\exp(tL)P^{-1}=\exp(-tL)=(\exp(tL)^{-1}$, 
for all $t$.
We seek an isotopy $\Theta_t=(\exp(tL),v_t)$ from $id_{Nil}$ to $\Theta$ such
that $\Theta_t\Psi\Theta_t=\Psi$ for all $t$.
It is easily seen that this imposes a linear condition on $v_t$ which has
an unique solution, and moreover $v_0=v_1=(0,0)$.

Now $\rho^{-1}h\rho(x,u)=(\Theta_{1-u}\Psi\Theta_u(x),1-u)=
(\Psi\Theta_{1-u}\Theta_u,1-u)$.
The loop $u\mapsto\Theta_{1-u}\Theta_u$
is freely contractible in $Aut_{Lie}(Nil)$,
since $\exp((1-u)L)\exp(uL)=\exp(L)$. 
It follows easily that $h$ does not change the framing of $C$.
\end{proof}

Instead of using the one-parameter subgroup determined by $L=\ln(A)$ 
we may use the polynomial isotopy given by
$A_t=\left(\smallmatrix 1&t\\
t&1+t^2\endsmallmatrix\right)$,
for $0\leq t\leq1$.
A similar argument could be used for the polynomial $X^2-X+1$,
which is realized by $\tau_63_1$ and its Gluck reconstruction.
On the other hand, the polynomial $X^2+X-1$ is not symmetric, 
and so the corresponding knots are not $+$amphicheiral. 

\begin{theorem} 
Let $K$ be a fibred $2$-knot with closed fibre $N_q$.
\begin{enumerate}
\item If the fibre is $N_1$ and the monodromy has characteristic polynomial
$X^2-X+1$ or $X^2-3X+1$ then $K$ is not reflexive;

\item
If the fibre is $N_q$ ($q$ odd) and the monodromy has characteristic polynomial
$X^2\pm X-1$ then $K$ is reflexive.
\end{enumerate}
\end{theorem}

\begin{proof} 
As $\tau_63_1$ is shown to be not reflexive in \S7 below, 
we shall concentrate on the knots with polynomial $X^2-3X+1$, 
and then comment on how our argument may be modified to handle the other cases.
                                    
Let $\Theta$, $\theta$ and $M=N_1\times_\theta S^1$ be as in Lemma 18.12,
and let $\widehat M=Nil\times_\Theta S^1$ be as in \S1.                          We shall take $[0,0,0,0]$ as the basepoint of $\widehat M$ 
and its image in $M$ as the basepoint there.

Suppose that $\Omega=(B,\nu)$ is an automorphism of $\Gamma_1$ which commutes with $\Theta$.
Since the eigenvalues of $A$ are both positive, $A(u)=uA+(1-u)I$ is 
invertible and $A(u)B=BA(u)$, for all $0\leq u\leq1$.
We seek a path of the form $\gamma(u)=(A(u),\mu(u))$ with commutes with $\Omega$.
Equating the second elements of the ordered pairs $\gamma(u)\Omega$ and $\Omega\gamma(u)$, 
we find that $\mu(u)(B-\det(B)I)$ is uniquely determined.
If $\det(B)$ is an eigenvalue of $B$ then there is 
a corresponding eigenvector $\xi$ in $\mathbb{Z}^2$.
Then $BA\xi=AB\xi=\det(B)A\xi$, so $A\xi$ is also an eigenvector of $B$.
Since the eigenvalues of $A$ are irrational we must have $B=\det(B)I$ 
and so $B=I$.
But then $\Omega\Theta=(A,\nu A)$ and $\Theta\Omega=(A,\nu)$, 
so $\nu(A-I)=0$ and hence $\nu=0$.
Therefore $\Omega=id_{Nil}$ and there is no difficulty in finding such a path.
Thus we may assume that $B-\det(B)I$ is invertible, and then $\mu(u)$ is uniquely determined.
Moreover, by the uniqueness, when $A(u)=A$ or $I$ we must have $\mu(u)=(0,0)$.
Thus $\gamma$ is an isotopy from $\gamma(0)=id_{Nil}$ 
to $\gamma(1)=\Theta$ (through diffeomorphisms of $Nil$)
and so determines a diffeomorphism $\rho_\gamma$ from 
$\mathbb{R}^3\times S^1$ to $\widehat M$ via 
$\rho_\gamma(r,s,t,[u])=[\gamma(u)([r,s,t]),u]$.

A homeomorphism $f$ from $\Sigma$ to $\Sigma_\tau$ carrying $K$ to $K_\tau$ 
(as unoriented submanifolds) extends to a self homeomorphism $h$ of $M$
which leaves $C$ invariant, but changes the framing. 
We may assume that $h$ preserves the orientations of $M$ and $C$, by Lemma 18.12. 
But then $h$ must preserve the framing, by Lemma 18.3.
Hence there is no such homeomorphism and such knots are not reflexive.

If $\pi\cong\pi\tau_63_1$ then we may assume that
the meridianal automorphism is
$\Theta=(\left(\smallmatrix 1&-1\\
          1&0\endsmallmatrix\right),(0,0))$.
As an automorphism of $Nil$, $\Theta$ fixes the centre pointwise, and it has order 6.
Moreover $(\left(\smallmatrix 0&1\\
          1&0\endsmallmatrix\right),(0,0)$ is an involution of $Nil$
which conjugates $\Theta$ to its inverse, and so $M$ admits an orientation reversing involution.
It can easily be seen that any automorphism of $\Gamma_1$ which commutes with $\Theta$ 
is a power of $\Theta$, and the rest of the argument is similar.

If the monodromy has characteristic polynomial $X^2\pm X-1$ 
we may assume that the meridianal automorphism is $\Theta=(D,(0,0))$,
where $D=\left(\smallmatrix 1&1\\
          1&0\endsmallmatrix\right)$ or its inverse.
As $\Omega=(-I,(-1,1))$ commutes with $\Theta$ (in either case) it
determines a self homeomorphism $h_\omega$ of $M=N_q\times_\theta S^1$ 
which leaves the meridianal circle $\{0\}\times S^1$ pointwise fixed.
The action of $h_\omega$ on the normal bundle may be 
detected by the induced action on $\widehat M$.
In each case there is an isotopy from $\Theta$ to 
$\Upsilon=(\left(\smallmatrix 1&0\\
          0&-1\endsmallmatrix\right),(1,0))$
which commutes with $\Omega$, and so
we may replace $\widehat M$ by the mapping torus $Nil\times_\Upsilon S^1$.
(Note also that $\Upsilon$ and $\Omega$ act linearly
under the standard identification of $Nil$ with $\mathbb{R}^3$).

Let $R(u)\in SO(2)$ be rotation through $\pi u$ radians,
and let $v(u)=\left(\smallmatrix 0\\
          u\endsmallmatrix\right)$, for $0\leq u\leq1$.
Then $\gamma(u)=\left(\smallmatrix 1&0\\
          v(u)&R(u)\endsmallmatrix\right)$
defines a path $\gamma$ in $SL(3,\mathbb{R})$ from $\gamma(0)=id_{Nil}$ 
to $\gamma(1)=\Upsilon$ which we may use to identify the mapping torus 
of $\Upsilon$ with $\mathbb{R}^3\times S^1$.
In the ``new coordinates" $h_\omega$ acts by sending $(r,s,t,e^{2\pi iu})$
to $(\gamma(u)^{-1}\Omega\gamma(u)(r,s,t),e^{2\pi iu})$.
The loop sending $e^{2\pi iu}$ in $S^1$ to $\gamma(u)^{-1}\Omega\gamma(u)$
in $SL(3,\mathbb{R})$ is freely homotopic to the loop 
$\gamma_1(u)^{-1}\Omega_1\gamma_1(u)$,
where $\gamma_1(u)=\left(\smallmatrix 1&0\\
          0&R(u)\endsmallmatrix\right)$
and $\Omega_1=diag[-1,-1,1]$.
These loops are essential in $SL(3,\mathbb{R})$,
since 
$\Omega_1\gamma_1(u)^{-1}\Omega_1\gamma_1(u)=
\left(\smallmatrix 1&0\\
          0&R(2u)\endsmallmatrix\right)$.
Thus $h_\omega$ induces the twist $\tau$ on the normal bundle of the meridian, 
and so the knot is equivalent to its Gluck reconstruction.
\end{proof}

We may refine one aspect of Litherland's observations 
on symmetries of twist spins \cite{[Li85]} as follows.
It is convenient to use here the formulation of twist spinning given in
\cite{[Ze65]}.
Let $k$ be a 1-knot which meets a small 3-ball $B$ in an unknotted arc
$b=B\cap{k}$,
and let $\rho_\phi$ be rotation of $\overline{S^3\setminus{B}}$ 
about an axis with endpoints $\partial{b}$. Then $\tau_rk$ is the 2-knot with image
\[
\cup_{0\leq\phi\leq2\pi}
(\rho_{r\phi}(k\setminus{b})\times\{e^{i\phi}\})\cup(S^0\times{D^2})
~\subset~S^4=(\overline{S^3\setminus{B}}\times{S^1})\cup(S^2\times{D^2})
\]
If $k$ is strongly invertible there is an involution $h$ of $S^3$
which rotates $S^3$ about an unknotted axis $A$ and 
such that $h(K)=K$ and $A\cap{K}=S^0$.
Let $P\in{A}\cap{K}$ and let $B$ be a small ball about $P$ such that $h(B)=B$.
Define an involution $H$ of $S^4$ by
$H(s,\theta)=(\rho_{-r\theta}h\rho_{-r\theta}(s),-\theta)$ 
for $(s,\theta)\in(S^3\setminus{B})\times{S^1}$ and 
$H(s,z)=(h(s),\bar{z})$ for $(s,z)\in{S^2\times{D^2}}$.
Then $H(\tau_rk)=\tau_rk$ and so $\tau_rk$ is strongly $+$amphicheiral.

Since the trefoil knot $3_1$ is strongly invertible $\tau_63_1$ is strongly 
$+$amphicheiral.
The involution of $X(\tau_63_1)$ extends to an involution of $M(\tau_63_1)$
which fixes the canonical section $C$ pointwise and 
does not change the framing of the normal bundle,
and hence $(\tau_63_1)^*$ is also $+$amphicheiral. 

\begin{theorem}
Let $K$ be a $2$-knot with group $\pi(e,\eta)$.
Then $K$ is reflexive.
If $K=\tau_2k(e,\eta)$ it is strongly $+$amphicheiral,
but no other $2$-knot with this group is $+$amphicheiral.
\end{theorem}

\begin{proof}
The first assertion follows from Lemma 18.2,
while $\tau_2k(e,\eta)$ is strongly $+$amphicheiral since
the Montesinos knot $k(e,\eta)$ is strongly invertible.

Let $t$ be the image of the canonical section in $\pi$.
Every strict weight orbit representing 
the image of $t$ in $\pi/\pi'$ 
contains an unique element of the form $u^nt$, by Theorem 16.15.
Let $\alpha$ be the automorphism of $\pi$ given by 
$\alpha(t)=t^{-1}$, $\alpha(x)=x$ and $\alpha(y)=y$.
If there were an automorphism $\theta$ such that $\theta(u^nt)=(u^nt)^{-1}$
then $\alpha{c_t}\theta(u^nt)=u^{-n}t$.
Hence if $n\not=0$ then $K$ is not $+$amphicheiral.
\end{proof}

The other 2-knot groups $\pi$ with $\pi'$ a $\mathbb{N}il^3$-lattice 
are as in Theorem 16.15.
The geometric argument used in \cite{[Hi11']} for the knots with group
$\pi(e,\eta)$ appears to break down for the 2-knots with groups as in part (2) of Theorem 16.15, and the questions of reflexivity and
$+$amphicheirality are open for these.
Automorphisms of $\mathbb{N}il^3$-lattices preserve orientation,
and so no fibred 2-knot with closed fibre a $\mathbb{N}il^3$-manifold is $-$amphicheiral or invertible.

It has been shown that for many of the Cappell-Shaneson knots 
at least one of the (possibly two) corresponding smooth homotopy 4-spheres 
is the standard $S^4$ \cite{[AR84]}.
Can a similar study be made in the $Nil$ cases?

\section{Other geometrically fibred knots}

We shall assume henceforth throughout this section that $k$ is a prime simple 
1-knot, i.e., that $k$ is either a torus knot or a hyperbolic knot. 

\begin{lemma}
Let $A$ and $B$ be automorphisms of a group $\pi$ such that $AB=BA$, 
$A(h)=h$ for all $h$ in $\zeta\pi$ and the images of $A^i$ and $B$ in 
$Aut(\pi/\zeta\pi)$ are equal.
Let $[A]$ denote the induced automorphism of $\pi/\pi'$.
If $I-[A]$ is invertible in $End(\pi/\pi')$ then $B=A^i$ in $Aut(\pi)$.
\end{lemma}

\begin{proof} 
There is a homomorphism $\epsilon:\pi\to\zeta\pi$ such that 
$BA^{-i}(x)=x\epsilon(x)$ for all $x$ in $\pi$.
Moreover $\epsilon A=\epsilon$, since $BA=AB$.
Equivalently, $[\epsilon](I-[A])=0$, where $[\epsilon]:\pi/\pi'\to\zeta\pi$ is 
induced by $\epsilon$.
If $I-[A]$ is invertible in $End(\pi/\pi')$ then $[\epsilon]=0$ and so
$B=A^i$.
\end{proof}

Let $p=ap'$, $q=bq'$ and $r=p'q'c$, where $(p,q)=(a,r)=(b,r)=1$.
Let $A$ denote both the canonical generator of the $Z/rZ$ action on 
the Brieskorn manifold $M(p,q,r)$ given by $A(u,v,w)=(u,v,e^{2\pi i/r}w)$
and its effect on $\pi_1(M(p,q,r))$.
The image of the Seifert fibration of $M(p,q,r)$ under the projection
to the orbit space $M(p,q,r)/\langle A\rangle\cong S^3$ is the Seifert 
fibration of $S^3$ with one fibre of multiplicity $p$ 
and one of multiplicity $q$.
The quotient of $M(p,q,r)$ by the subgroup generated by $A^{p'q'}$ 
may be identified with $M(p,q,p'q')$.
We may display the factorization of these actions as follows:
\begin{equation*}
\begin{CD}
M(p,q,r)@> /S^1>> P^2(p,q,r) \\
@VVV                   @VVV     \\         
M(p,q,p'q')@> /S^1>> P^2(p,q,p'q') \\
@VVV                   @VVV     \\         
(S^3,(p,q))@> /S^1>> S^2
\end{CD}
\end{equation*}

Sitting above the fibre in $S^3$ of multiplicity $p$ in both $M$'s we find 
$q'$ fibres of multiplicity $a$,
and above the fibre of multiplicity $q$ we find $p'$ fibres of multiplicity $b$.
But above the branch set, a principal fibre in $S^3$, we have one fibre of 
multiplicity $c$ in $M(p,q,r)$, but a principal fibre in $M(p,q,p'q')$.
(Note that the orbit spaces $P^2(p,q,r)$ and $P^2(p,q,p'q')$ are in fact
homeomorphic.)

We have the following characterization of the centralizer of $A$ in $Aut(\pi)$.

\begin{theorem}
Assume that $p^{-1}+q^{-1}+r^{-1}\leq1$, and let
$A$ be the automorphism of $\pi=\pi_1(M(p,q,r))$ of order $r$ induced by the 
canonical generator of the branched covering transformations.
If $B$ in $Aut(\pi)$ commutes with $A$ then $B=A^i$ for some $0\leq i<r$.
\end{theorem}

\begin{proof}                                                                    The 3-manifold $M=M(p,q,r)$ is aspherical, 
with universal cover $\mathbb{R}^3$,
and $\pi$ is a central extension of $Q(p,q,r)$ by $\mathbb{Z}$.
Here $Q=Q(p,q,r)$ is a discrete planar group with signature
$((1-p')(1-q')/2;a,\dots,a,b,\dots,b,c)$ 
(where there are $q'$ entries $a$ and $p'$ entries $b$).
Note that $Q$ is Fuchsian, 
except for $Q(2,3,6)\cong\mathbb{Z}^2$.
(In general, $Q(p,q,pq)$ is a $PD_2^+$-group of genus $(1-p)(1-q)/2$).

There is a natural homomorphism from $Aut(\pi)$ to $Aut(Q)=Aut(\pi/\zeta\pi)$.
The strategy shall be to show first that $B=A^i$ in $Aut(Q)$ and then lift to $Aut(\pi)$.
The proof in $Aut(Q)$ falls naturally into three cases.

{\sl Case 1}. $r=c$. In this case $M$ is a homology 3-sphere, 
fibred over $S^2$ with three exceptional fibres of multiplicity $p$, $q$ and $r$.
Thus $Q\cong\Delta(p,q,r)=\langle q_1,q_2,q_3\mid 
q_1^p=q_2^q=q_3^r=q_1q_2q_3=1\rangle$,
the group of orientation preserving symmetries of a tesselation of $H^2$ 
by triangles with angles $\pi/p$, $\pi/q$ and $\pi/r$.
Since $Z_r$ is contained in $S^1$, $A$ is inner.
(In fact it is not hard to see that the image of $A$ in $Aut(Q)$ 
is conjugation by $q_3^{-1}$.
See \S3 of \cite{[Pl83]}).                        

It is well known that the automorphisms of a triangle group correspond
to symmetries of the tessellation (see Chapters V and VI of \cite{[ZVC]}). 
Since $p$, $q$ and $r$ are pairwise relatively prime there are no self 
symmetries of the $(p,q,r)$ triangle. 
So, fixing a triangle $T$, all symmetries take $T$ to another triangle.
Those that preserve orientation correspond to elements of $Q$ acting by inner 
automorphisms, 
and there is one nontrivial outerautomorphism, $R$ say, given by reflection 
in one of the sides of $T$. We can assume $R(q_3)=q_3^{-1}$.

Let $B$ in $Aut(Q)$ commute with $A$.
If $B$ is conjugation by $b$ in $Q$ then $BA=AB$ is equivalent to $bq_3=q_3b$,
since $Q$ is centreless.
If $B$ is $R$ followed by conjugation by $b$ then $bq_3=q_3^{-1}b$.
But since $\langle q_3\rangle=Z_r$ in $Q$ is generated 
by an elliptic element the normalizer of $\langle q_3\rangle$
in $PSL(2,\mathbb{R})$ consists of elliptic elements 
with the same fixed point as $q_3$.
Hence the normalizer of $\langle q_3\rangle$ in $Q$ 
is just $\langle q_3\rangle$.
Since $r>2$, $q_3\not=q_3^{-1}$ and so we must have $bq_3=q_3b$, 
$b=q_3^i$ and $B=A^i$.
(Note that if $r=2$ then $R$ commutes with $A$ in $Aut(Q)$).

{\sl Case 2}. $r=p'q'$ so that $Z_r\cap S^1=1$.                      
The map from $S^2(p,q,p'q')$ to $S^2$ is branched over three points in $S^2$.
Over the point corresponding to the fibre of multiplicity $p$ in $S^3$
the map is $p'$-fold branched; it is $q'$-fold branched                                         
over the point corresponding to the fibre of multiplicity $q$ in $S^3$,
and it is $p'q'$-fold branched over the point $*$ corresponding
to the branching locus of $M$ over $S^3$.
                                        
Represent $S^2$ as a hyperbolic orbifold $H^2/\Delta(p,q,p'q')$.
(If $(p,q,r)=(2,3,6)$ we use instead the flat orbifold 
$\mathbb{R}^2/\Delta(2,3,6)$).
Lift this to an orbifold structure on $S^2(p,q,p'q')$, 
thereby representing $Q=Q(p,q,p'q')$ into $PSL(2,\mathbb{R})$.
Lifting the $Z_{p'q'}$-action to $H^2$ gives 
an action of the semidirect product $Q\rtimes Z_{p'q'}$ on $H^2$, 
with $Z_{p'q'}$ acting as rotations about a point 
$\tilde*$ of $H^2$ lying above $*$.
Since the map from $H^2$ to $S^2(p,q,p'q')$ is unbranched at $\tilde*$ (equivalently,
$Z_r\cap S^1=1$), $Q\cap Z_{p'q'}=1$. Thus $Q\rtimes Z_{p'q'}$
acts effectively on $H^2$, with quotient $S^2$ and three branch points, 
of orders $p$, $q$ and $p'q'$.

In other words, $Q\rtimes Z_{p'q'}$ is isomorphic to $\Delta(p,q,p'q')$.
The automorphism $A$ extends naturally to an automorphism of $\Delta$, 
namely conjugation by an element of order $p'q'$, 
and $B$ also extends to $Aut(\Delta)$, since $BA=AB$.

We claim $B=A^i$ in $Aut(\Delta)$.
We cannot directly apply the argument in Case 1, 
since $p'q'$ is not prime to $pq$.
We argue as follows.
In the notation of Case 1, $A$ is conjugation by $q_3^{-1}$.
Since $BA=AB$,
$B(q_3)=q_3^{-1}B(q_3)q_3$, which forces $B(q_3)=q_3^j$.
Now  $q_3^{-1}B(q_2)q_3 =AB(q_2)=B(q_3^{-1})B(q_2)B(q_3)=q_3^{-j}B(q_2)q_3^j$,
or $B(q_2)=q_3^{1-j}B(q_3)q_3^{j-1}$.
But $B(q_2)$ is not a power of $q_3$, so $q_3^{1-j}=1$, 
or $j\equiv 1$ modulo $(r)$.
Thus $B(q_3)=q_3$.
This means that the symmetry of the tessellation that realizes $B$ shares 
the same fixed point as $A$, 
so $B$ is in the dihedral group fixing that point, 
and now the proof is as before.
                    
{\sl Case 3}. $r=p'q'c$ (the general case). 
We have $Z_{p'q'c}$ contained in $Aut(\pi)$, but $Z_{p'q'c}\cap S^1=Z_c$,
so that $Z_c$ is the kernel of the composition 
\begin{equation*}
Z_r\to Out(\pi)\to Out(Q).
\end{equation*}
Let $\bar Q$ be the extension corresponding to the abstract kernel 
$Z_{p'q'}\to Out(Q)$.
(The extension is unique since $\zeta Q=1$).
Then $\bar Q$ is a quotient of the semidirect product $Q(p,q,r)\rtimes (Z/rZ)$
by a cyclic normal subgroup of order $c$.

Geometrically, this corresponds to the following.
The map from $S^2(p,q,r)$ to $S^2$ is branched as in Case 2, 
over three points with branching indices $p$, $q$ and $p'q'$.
This time, represent $S^2$ as $H^2/\Delta(p,q,p'q')$.
Lift to an orbifold structure on $S^2(p,q,r)$ with one cone point of order $c$.
Lifting an elliptic element of order $r$ in $\Delta(p,q,r)$
to the universal orbifold cover of $S^2(p,q,r)$ gives $Z_r$ contained in $Aut(Q(p,q,r))$
defining the semidirect product.
But $Q(p,q,r)\cap Z_r=Z_c$, so the action is ineffective.
Projecting to $Z_{p'q'}$ and taking the extension $\bar Q$ kills the 
ineffective part of the action.
Note that $Q(p,q,r)$ and $Z_r$ inject into $\bar Q$.

As in Case 2, $\bar Q\cong\Delta(p,q,r)$,
$A$ extends to conjugation by an element of order $r$ in $\bar Q$, and 
$B$ extends to an automorphism of $Q(p,q,r)\rtimes Z_r$, since $BA=AB$.
Now $(q_3,p'q')$ in $Q(p,q,r)\rtimes Z_r$ normally generates the kernel of
$Q(p,q,r)\rtimes Z_r\to\bar Q$, where $q_3$ is a rotation of order $c$ 
with the same fixed point as the generator of $Z_r$.
In other words, 
$A$ in $Aut(Q(p,q,r))$ is such that $A^{p'q'}$ is conjugation by $q_3$.
Since $BA^{p'q'}=A^{p'q'}B$ the argument in Case 2 shows that $B(q_3)=q_3$.
So $B$ also gives an automorphism of $\bar Q$, and now the argument 
of Case 2 finishes the proof.

We have shown that $B=A^i$ in $Aut(Q)$.
Since $A$ in $Aut(\pi)$ is the monodromy of a fibred knot in $S^4$
(or, more directly, since $A$ is induced by a branched cover 
of a knot in a homology sphere), $I-[A]$ is invertible.
Thus the Theorem now follows from Lemma 18.15.
\end{proof}

\begin{theorem}
Let $k$ be a prime simple knot in $S^3$.
Let $0<s<r$, $(r,s)=1$ and $r>2$. Then $\tau_{r,s}k$ is not reflexive.
\end{theorem}

\begin{proof} 
We shall consider separately the three cases 
(a) $k$ a torus knot and the branched cover aspherical;
(b) $k$ a torus knot and the branched cover spherical;
and (c) $k$ a hyperbolic knot.
                                                       
{\it Aspherical branched covers of torus knots.}                                              
Let $K=\tau_{r,s}(k_{p,q})$ where $r>2$ and $M(p,q,r)$ is aspherical.
Then $X(K)=(M(p,q,r)\setminus{int D^3})\times_{A^s}S^1$,                
$M=M(K)=M(p,q,r)\times_{A^s}S^1$ and
$\pi=\pi K\cong \pi_1(M(p,q,r))\rtimes_{A^s}\mathbb{Z}$.
If $K$ is reflexive there is a homeomorphism $f$ of $X$ which changes
the framing on $\partial X$.
Now $k_{p,q}$ is strongly invertible -- there is an involution 
of $(S^3,k_{p,q})$ fixing two points of the knot and reversing the meridian.
This lifts to an involution of $M(p,q,r)$ fixing two points of the branch set and conjugating $A^s$ to 
$A^{-s}$, thus inducing a diffeomorphism of $X(K)$ which reverses the meridian.
By Lemma 18.1 this preserves the framing, 
so we can assume that $f$ preserves the meridian of $K$.
It must also preserve the orientation,
by the remark following Theorem 2.14.
Since $M(p,q,r)$ is a $\widetilde{\mathbb{SL}}$-manifold
$\widetilde{M(p,q,r)}\cong\mathbb{R}^3$ and automorphisms of $\pi_1(M(p,q,r))$ 
are induced by fibre-preserving self-diffeomorphisms \cite{[Sc83']}. 
The remaining hypothesis of Lemma 18.3 is satisfied, by Theorem 18.16.
Therefore there is no such self homeomorphism $f$, and $K$ is not reflexive.

\smallskip
{\it Spherical branched covers of torus knots.}        
We now adapt the previous argument to the spherical cases.
The analogue of Theorem 18.16 is valid, except for $(2,5,3)$.
We sketch the proofs.

$(2,3,3)$: $M(2,3,3=S^3/Q(8)$. The image in $Aut(Q(8)/\zeta Q(8))\cong S_3$ 
of the automorphism $A$ induced by the 3-fold cover of the trefoil knot has order 3
and so generates its own centralizer.

$(2,3,4)$: $M(2,3,4)=S^3/T_1^*$. 
In this case the image of $A$ in $Aut(T_1^*)\cong S_4$ must be a 4-cycle,
and generates its own centralizer.

$(2,3,5)$: $M(2,3,5)=S^3/I^*$. 
In this case the image of $A$ in $Aut(I^*)\cong S_5$
must be a 5-cycle, and generates its own centralizer. 
                   
$(2,5,3)$: We again have $I^*$, but in this case $A^3=I$, 
say $A=(123)(4)(5)$.
Suppose $BA=AB$.
If $B$ fixes 4 and 5 then it is a power of $A$. 
But $B$ may transpose 4 and 5, and then $B=A^iC$, 
where $C=(1)(2)(3)(45)$ represents the nontrivial outer automorphism class 
of $I^*$.

Now let $K=\tau_{r,s}(k_{p,q})$ as usual, 
with $(p,q,r)$ one of the above four triples,
and let $M=M(p,q,r)\times_{A^s} S^1$.
As earlier, if $K$ is reflexive we have a homeomorphism $f$ which preserves 
the meridian $t$ and changes the framing on $D^3\times_{A^s} S^1$.
Let $\widehat M=S^3\times_{\hat A^s}S^1$ be the cover of $M$ corresponding 
to the meridian subgroup, where $\hat A$ is a rotation about an axis.
Let $f$ be a basepoint preserving self homotopy equivalence of $M$ such that 
$f_*(t)=t$ in $\pi$.
Let $B$ in $Aut(\pi_1(M(p,q,r))$ be induced by $f_*$, so $BA^s=A^sB$.
The discussion above shows that $B=A^{si}$ except possibly for $(2,5,3)$.
But if $B$ represented the outer automorphism of $I^*$ then after lifting to 
infinite cyclic covers we would have a homotopy equivalence of $S^3/I^*$ 
inducing $C$, contradicting Lemma 11.4.
So we have an obvious fibre preserving diffeomorphism $f_B$ of $M$.

The proof that $\hat f_B$ is homotopic to $id_{\widehat M}$ is exactly as in 
the aspherical case.
To see that $\hat f_B$ is homotopic to $\hat f$ (the lift of $f$ to a 
basepoint preserving proper self homotopy equivalence of $\widehat M$)
we investigate whether $f_B$ is homotopic to $f$.
Since $\pi_2(M)=0$ we can homotope $f_B$ to $f$ on the 2-skeleton of $M$.
On the 3-skeleton we meet an obstruction in 
$H^3(M;\pi_3)\cong{H^3(M;\mathbb{Z})}\cong\mathbb{Z}$,
since $M$ has the homology of $S^3\times S^1$.
But this obstruction is detected on the top cell of $M(p,q,r)$ and just 
measures the difference of the degrees of $f$ and $f_B$ on the infinite 
cyclic covers \cite{[Ol53]}.
Since both $f$ and $f_B$ are orientation preserving homotopy equivalences 
this obstruction vanishes.
On the 4-skeleton we have an obstruction in $H^4(M;\pi_4)=Z/2Z$, 
which may not vanish.
But this obstruction is killed when we lift to $\widehat M$, 
since the map from $\widehat M$ to $M$ has even degree, 
proving that $\hat f_B\simeq \hat f$.

We now use radial homotopies on $S^3\times S^1$ to finish, as before.
                             
\smallskip
{\it Branched covers of hyperbolic knots.}        
Let $k$ be hyperbolic.
Excluding $N_3(4_1)$ 
(the 3-fold cyclic branched cover of the figure eight knot),
$N=N_r(k)$ is a closed hyperbolic 3-manifold, 
with $\langle\alpha\rangle\cong Z/rZ$ acting by isometries.
As usual, we assume there is a homeomorphism $f$ of $M=M(\tau_{r,s}(k))$ 
which changes the framing on $D^3\times_{A^s} S^1$.
As in the aspherical torus knot case, it shall suffice to show that the lift 
$\hat f$ on $\widehat M$ is properly homotopic to a map of 
$(\mathbb{R}^3\times S^1,D^3\times S^1)$ that does not change the framing 
on $D^3\times S^1$.

Letting $B=f_*$ on $\nu=\pi_1(N)$, we have $BA^sB^{-1}=A^{\pm s}$, 
depending on whether
$f_*(t)=t^{\pm1}$ in $\pi=\nu\rtimes_{A^s}\mathbb{Z}$.
There is an unique isometry $\beta$ of $N$ realizing the class of $B$ in 
$Out(\nu)$, by Mostow rigidity, and $\beta\alpha^s\beta^{-1}=\alpha^{\pm s}$.
Hence there is an induced self diffeomorphism $f_\beta$ of 
$M=N\times_{\alpha^s} S^1$.
Note that $f_*=(f_\beta)_*$ in $Out(\pi)$, so $f$ is homotopic to $f_\beta$.
We cannot claim that $\beta$ fixes the basepoint of $N$, but 
$\beta$ preserves the closed geodesic fixed by $\alpha^s$.

Now $\widehat M=H^3\times_{\hat\alpha^s}S^1$ where $\hat\alpha^s$ is an 
elliptic rotation about an axis $L$, and $\hat f_\beta$ is fibrewise an 
isometry $\hat\beta$ preserving $L$.
We can write $H^3=\mathbb{R}^2\times L$ (non-metrically!)
by considering the family of hyperplanes perpendicular to $L$,
and then $\hat\beta$ is just an element of $O(2)\times E(1)$ and
$\hat\alpha^s$ is an element of $SO(2)\times\{1\}$.
The proof of Lemma 18.1, with trivial modifications, shows that, 
after picking coordinates and ignoring orientations, 
$\hat f_\beta$ is the identity.
This completes the proof of the theorem.
\end{proof}

The manifolds $M(p,q,r)$ with $p^{-1}+q^{-1}+r^{-1}<1$
are coset spaces of $\widetilde{SL}$ \cite{[Mi75]}.
Conversely, let $K$ be a 2-knot obtained by surgery on the canonical cross-section of
$N\times_\theta S^1$, where $N$ is such a coset space.
If $\theta$ is induced by an automorphism of $\widetilde{SL}$ which 
normalizes $\nu=\pi_1(N)$ then it has finite order, since 
$N_{\widetilde{SL}}(\nu)/\nu\cong 
N_{PSL(2,\mathbb{R})}(\nu/\zeta\nu)/(\nu/\zeta\nu)$.
Thus if $\theta$ has infinite order we cannot expect to use such 
geometric arguments to analyze the question of reflexivity. 

We note finally that since torus knots are strongly invertible, 
their twist spins are strongly $+$amphicheiral.
(See the paragraph after Theorem 18.13.)
However, since automorphisms of $\widetilde{\mathbb{SL}}$-lattices 
preserve orientation,
no such knot is $-$amphicheiral or invertible.

\section{Smooth 2-knots in $S^4$}

In this book we have chosen to work on the level of 
homotopy equivalence or TOP $s$-cobordism, 
and to set aside all considerations of differential topology
(with some minor exceptions in Chapters 8, 12 and 13),
even though geometric manifolds, total spaces of fibre bundles 
and all the knot manifolds considered above have natural smooth structures.
In the case of 2-knots, although spinning, twist-spinning and 
related constructions give smooth knots in the standard $S^4$, 
the other major constructions (including Gluck reconstruction)
give only smooth knots in smooth homotopy 4-spheres.
(They typically involve first constructing
the knot manifold and then performing surgery on a meridianal loop
to obtain a 1-connected 4-manifold $\Sigma$ with $\chi(\Sigma)=0$,
which is then homeomorphic to $S^4$ by Freedman's theorem.)

We shall summarize here what is presently known about the smooth
homotopy 4-spheres deriving from 2-knot constructions.
Following Levine \cite{[Le78]}, we shall say that a 2-knot is {\it ordinary} 
if it is TOP equivalent to a smooth knot in $S^4$,
which shall henceforth denote the 4-sphere with its standard smoothing.

1). In \S6 of Chapter 14 we noted that if a knot group $\pi$ has a presentation 
of deficiency 1 such that the balanced presentation for the trivial group
obtained by adjoining a relator representing a meridian for $\pi$ 
is Andrews-Curtis equivalent to the empty presentation then $\pi$ 
is the group of an ordinary 2-knot \cite{[Le78], [Yo82']}. 
(However, it is widely believed that there may be balanced presentations 
of the trivial group which are not Andrews-Curtis equivalent to the empty presentation.)

2).  Cappell-Shaneson 2-knots have been the most popular source 
of candidates for potential counter-examples to the Smooth Poincar\'e 
Conjecture in dimension 4.
The  Cappell-Shaneson 2-knots with group $\pi=\mathbb{Z}^3\rtimes_A\mathbb{Z}$
are parametrized up to Gluck reconstruction (and reflection) 
by the ideal class monoid of the ring $\Lambda/(\Delta_A)$,
by Theorems 16.9 and  17.5.
If $A\in{SL(3,\mathbb{Z})}$ corresponds to the principal ideal class then
surgery with the untwisted framing always gives $S^4$ [AR84]. 
It was shown more recently that (for the principal ideal class)
surgery with the twisted framing (i.e., the Gluck reconstruction
of the result of surgery with the untwisted framing) also gives $S^4$ \cite{[Ak10]}.
The arguments involved delicate study of handle decompositions.

Gompf has introduced a new approach,
which avoids direct analysis of handlebody structures,
using instead logarithmic transformations of ``fishtail neighbourhoods" 
to relate mapping tori with different monodromies.
Applied to Cappell-Shaneson 2-knots,
it follows that certain infinite families of matrices
with distinct traces (and thus not conjugate) give rise to diffeomorphic 2-spheres. 
We may assume that the defining matrices are in standard form
\[ A=
\left(
\begin{matrix}
0& a & b\\
0 & c &d\\
1& 0& n-c
\end{matrix}
\right),
\]
with $\det(A-I)=1$ and $\Delta_A(t)=f_n(t)$.
Then $b=(c-1)(n-c-1)$ and $ad=f_n(c)$, 
so $A$ is determined by the triple $(c,d,n)$.
(The case $c=d=1$ corresponds to the principal ideal class.)
If $A$ and $A'$ are two such matrices and $(c',d',n')=(c,d,n+kd)$ for some $k\in\mathbb{Z}$
then the homotopy 4-spheres obtained via surgery with the untwisted (respectively, twisted)
framing on the mapping tori of $A$ and $A'$ are diffeomorphic \cite{[Go10]}.
Hence all Cappell-Shaneson 2-knots corresponding to principal ideal classes 
are ordinary,
since we may reduce to the case $n=-1$, 
which is reflexive by Theorem 18.5,
and appeal to \cite{[AR84]}.

The equivalence relation on such matrices generated by similarity and 
this condition on the entries is called {\it Gompf equivalence}.
If two such matrices are each equivalent to a matrix with trace $n$  
in the interval $-64\leq{n}\leq69$ then they are Gompf equivalent.
Hence all 2-knots corresponding to such matrices are ordinary,
since the ring $\Lambda/(f_{-1})$ is a PID.
Moreover,  the endomorphism $\Phi_n$ of $\mathbb{Z}[t]$ given by
$\Phi_n(t)=t^2+(n-4)t+1$ induces an isomorphism of rings
$\varphi_n:\Lambda/(f_n)\cong\Lambda/(f_{5-n})$,
with inverse $\varphi_{5-n}$, 
thereby explaining a symmetry observed in calculations by Gompf.
The induced isomorphism of ideal class monoids is compatible 
with Gompf equivalence.
It is an open question whether all $A\in{SL(3,\mathbb{Z})}$ 
such that $\det(A-I)=1$ are Gompf equivalent \cite{[KY17]}.

3). Branched twist spins and their Gluck reconstructions are ordinary 
\cite{[Go76],[Pa78]}.
(A fibred 2-knot with geometric monodromy of finite order 
is a branched twist spin, 
by \cite{[Pa78], [Pl86]} together with the 3-dimensional Poincar\'e conjecture, 
since proven by Perelman.)

4). Most 2-knot groups $\pi$ with $\pi'$ finite are realized by twist spins,
but only a handful of those with $\pi'\cong{Q(8)\times\mathbb{Z}/n\mathbb{Z}}$
have been realized by ordinary knots \cite{[Kn88],[Tr90]}.
(See Chapter 15 above.
These groups are not realized by fibred 2-knots with 
geometric monodromy of finite order.
See \S3 of Chapter 16.)

5). It is not known whether all 2-knots with torsion-free 
virtually poly-$Z$ knot groups (other than twist spins) are ordinary.
The main theorem of \cite{[Go10]} may apply also when $\pi'$ is non-abelian,
but this has not been tried thus far.

This section would be redundant should the Smooth Poincar\'e Conjecture 
in dimension 4 be proven!

%% file: m5-bib.tex
\bibliographystyle{amsalpha}

%% file: m5-i.tex
\begin{theindex}

{\sl Expressions beginning with Greek 

characters and non-alphabetic symbols are
listed at the end of this index}.
\bigskip

\item{ $(A,\beta,C)$\quad (isometry of $\mathbb{S}^2\times\mathbb{E}^2$), 204}

\item{ $A(m,e)$\quad (metacyclic group 

\subitem of order $2^e m$), 222}

\item{ $A(\pi)$\quad (augmentation ideal 

\subitem of $\mathbb{Z}[\pi]$), 36}

\item{ action  (of a Seifert fibration), 145}

\item{ admits a geometry, 132}
               
\item{ algebraic 2-type

\subitem ($[\pi,\pi_2 (M),k_1 (M)]$), 26}

\item{ algebraic mapping torus, 74}

\item{ almost coherent, 16}

\item{ almost complex structure,  150}

\item{ almost finitely presentable ($FP_2$), 14}

\item{ almost linear $k$-invariant, 225}

\item{ amenable group, 8}
                                          
\item{ amphicheiral knot, 
271}       

\item{ Artin spin of a knot\quad ($\sigma K$), 
276}

\item{ ascendant, 5}

\item{ aspherical (orbifold), 138}

\item{ automorphisms of $\Gamma_q$, 168}

\indexspace
             
\item{ $B_1 - B_4$\quad (nonorientable flat 

\subitem  3-manifold groups), 154}

\item{ $BE(X)$\quad (classifying space), 89} 

\item{ $BS(p,q)$ \quad (Baumslag-Solitar group), 308}
 
\item{ bad (orbifold), 138}

\item{ Bass Conjectures, 14}

\item{ Baumslag-Solitar group\quad ($BS(p,q)$), 308}
            
\item{ Bieri's Theorem

\subitem (Theorem 8.8 of [Bi]), 20}

\item{ Bieri-Strebel Theorem [BS78], 14}

\item{ Bogomolov's Theorem, 149}

\item{ boundary link, 
287}

\item{ Bowditch's Theorem, 21}

\item{ branched twist spin, 
315}

\item{ Brieskorn manifold\quad 

\subitem ($M(p,q,r)$), 
311}

\item{ Brown-Geoghegan Theorem  

\subitem [BG85], 18}

\indexspace
             
\item{ $c(\hat g)$\quad (Kervaire-Arf invariant 

\subitem of $\hat g:M\to G/TOP$), 117}

\item{ $c_M :M\to K(\pi_1 (M),1)$\quad 

\subitem (classifying map), 26}

\item{ $\mathbb{CP}^2 $\quad (geometry of complex 

\subitem projective plane), 234}

\item{ $Ch=*CP^2$\quad (the fake 

\subitem complex projective plane), 235}

\item{ $C_G (H)$\quad (centralizer 

\subitem of a subgroup), 3}

\item{ $Cl$\quad (Waldhausen's class 

\subitem of groups),  112}
                
\item{ canonical cross-section, 
342}
                  
\item{ Cappell-Shaneson knot, 
320}

\item{ Cartan-Leray spectral sequence, 26}

\item{ centre of a group $G$\quad ($\zeta G$), 3}

\item{ characteristic subgroup, 3}

\item{ class $VII$ (complex surface), 149}

\item{ classifying map\quad 

\subitem ($c_M :M\to K(\pi_1 (M),1)$), 26}

\item{ closed fibre, 
273}
               
\item{ closed manifold, 26}

\item{ codimension-2 Kervaire invariant, 117}
             
\item{ coherent group, 15}

\item{ coherent ring, 15}

\item{ cohomology intersection pairing, 66}

\item{ coinduced (module), 21}
   
\item{ commutator subgroup 

\subitem of a group $G$ ($G'$), 3}

\item{ companion, 
275}
    
\item{ complex surface, 148, 
261}

\item{ complex torus, 149}

\item{ conjugate $\bar M$ of a module $M$, 13}

\item{ connecting homomorphism 

\subitem $\partial:\pi_2(B)\to\pi_1(F)$, 89}

\item{ Crisp's Theorem [Cr00], 34}

\item{ cusp, 139}    

\indexspace
             
\item{ $D$\quad (infinite dihedral group 

\subitem $(Z/2Z)*(Z/2Z)$), 16}

\item{ deficiency\quad  ($\mathrm{def}(P)$, $\mathrm{def}(\pi)$), 28}

\item{ deform spin of a knot, 
276}

\item{ $dim_{{\mathcal N}(\pi)}(M)$ (von Neumann 

\subitem dimension of $M$), 24}

\item{ doubly slice knot, 
277}

\indexspace

\item{ $e(G)$ \quad (number of ends of the group  

\subitem $G$, $=0$, 1, 2 or $\infty$), 16}

\item{ $e^\mathbb{Q}(\pi)$ \quad (Euler class of an aspherical 

\subitem Seifert fibration), 145}

\item{ $\mathbb{E}^n $\quad (flat geometry), 134}                          

\item{ $E(n)$\quad (isometry group of $\mathbb{E}^n$), 134}

\item{ $E(X)$, $E_0 (X)$\quad (space of self 

\subitem homotopy equivalences), 89}

\item{ $EG$\quad (class of elementary amenable 

\subitem groups), 9}

\item{ $ev$, $ev^{(2)}$ \quad (evaluation homomorphisms), 

\subitem 49}
 
\item{ elementary amenable, 9}
 
\item{ elliptic surface, 149, 
261, 336}

\item{ ends (and $H^1(G;\mathbb{Z}[G])$), 16}

\item{ equivariant (co)homology, 25}

\item{ Euler characteristic formula,

\subitem $\chi(X)=\Sigma(-1)^i\beta_i^{(2)}(X)$, 27}
 
\item{ Euler class (of a Seifert fibration), 145}

\item{ evaluation homomorphism ($ev$, $ev^{(2)}$), 

\subitem 49}
                    
\item{ extension of groups, 4}

\item{ exterior of a knot\quad ($X(K)$, $X$), 
269}

\indexspace

\item{ $f_\alpha$\quad (self homotopy equivalence 

\subitem of a closed 4-manifold), 118}
                 
\item{ $f_M:M\to P_2 (M)$\quad (second map 

\subitem of Postnikov tower), 26}

\item{ $FF$, $FP$, $FP_n$\quad  (finiteness 

\subitem conditions), 14}

\item{ $F(r)$\quad (free group), 3}

\item{ $\mathbb{F}^4$\quad (geometry of $T_{\mathbb{H}^2}$), 133, 
258}

\item{ Farrell's Theorem  [Fa74], 18}

\item{ fibration theorem, 124}
                                   
\item{ fibred knot, 
271}

\item{ finite $k$-skeleton, 23}
                 
\item{ finite $PD_n$-complex, 33}                         
                 
\item{ finite $PD_n$-space, 33}                         

\item{ finitely dominated

\subitem (chain complex), 23}

\item{ flat manifold, 134}                   

\item{ flat $n$-manifold group, 134}

\item{ foliation by circles, 
263}

\item{ F\o lner exhaustion, 9}

\item{ fundamental triple 

\subitem (of a $PD_3$-complex), 34}

\indexspace
                                                     
\item{ $g.d.$\quad (geometric dimension), 28}

\item{ $G_1 - G_6$\quad (orientable flat 

\subitem 3-manifold groups), 153}

\item{ $G(\pm)$\quad  (flat 2-knot groups), 
320}

\item{ generalized Eilenberg-Mac Lane

\subitem complex, 214}

\item{ geometric decomposition, 139}

\item{ geometric dimension 

\subitem of a group ($g.d.$), 28}

\item{ geometry, 132}

\item{ Gildenhuys-Strebel Theorem 

\subitem [GS81], 17}

\item{ Gluck reconstruction of an

\subitem $S^2$-orbifold bundle, 139}

\item{ Gluck reconstruction of a knot $K$\quad 

\subitem ($K^*$), 
270}

\item{ Gompf equivalence,  364}

\item{ good (orbifold), 138}

\item{ graph manifold, 114}

\item{ Gromov's Theorem \cite[\S8.A]{[Gr]}, 28}

\indexspace

\item{ $\mathbb{H}^2\times\mathbb{H}^2$\quad (semisimple product

\subitem geometry), 188}
                                            
\item{ $\mathbb{H}^4$, $\mathbb{H}^2(\mathbb{C})$\quad (rank 1 geometries), 
 192}

\item{ $\mathbb{H}^2\times\mathbb{E}^2$\quad (product geometry), 182}
                                           
\item{ $\mathbb{H}^3\times\mathbb{E}^1$\quad (product geometry), 185}

\item{ $H_i (X;R[G/H])$, $H^i (X;R[G/H])$\quad 

\subitem (equivariant (co)homology), 25}

\item{ $h(G)$\quad (Hirsch length 

\subitem of a group $G$), 10}

\item{ Haken 3-manifold, 114}

\item{ Hantzsche-Wendt flat 3-manifold 

\subitem group\quad ($G_6$), 154}                        

\item{ Hendrik's Theorem [Hn], 34}

\item{ Hilbert ${\mathcal N}(\pi)$-module, 23}

\item{ Hilbert ${\mathcal N}(\pi)$-complex, 24}

\item{ Hirsch length  \quad ($h(G)$), 4, 10}

\item{ Hirsch-Plotkin radical \quad  ($\sqrt G$), 7}

\item{ homology 4-sphere, 
290}

\item{ holonomy group, 134}

\item{ homotopy ribbon knot, 
275}

\item{ Hopf surface, 149, 
336}

\item{ hyperelliptic surface, 149}

\indexspace

\item{ $I(G)=\{ g\in G\mid \exists n>0,g^n\in G'\}$, 3}

\item{ $I^*$\quad  (binary icosahedral group), 222}

\item{ $I_\pi$\quad (homomorphism 

\subitem from $H_1(\pi;\mathbb{Z})$ to $L_1^s(\pi)$), 119}

\item{ $I_\pi^+$\quad (homomorphism 

\subitem from Ker$(w)$ to $L_1^s(\pi,w)$), 119}

\item{ indicable group, 3}           

\item{ infinite cyclic covering space\quad 

\subitem ($E_\nu$, $X'(K)$, $M'(K)$), 77, 
271}

\item{ infinite dihedral group\quad 

\subitem ($D=(Z/2Z)*(Z/2Z)$), 16}
                                                                 
\item{ infranilmanifold, 134}

\item{ infrasolvmanifold, 135, 177}

\item{ Inoue surface, 149, 
336}

\item{ intersection pairing, 66}

\item{ invertible knot, 
269}

\item{ irreducible knot, 
273}

\item{ $Isom(\mathbb{X})$, 132}

\indexspace

\item{ $J_+ (F)$ (kernel of action of $Out(F)$

\subitem  on  $H_3(F;\mathbb{Z})$), 219}

\item{ $J(F)$\quad  (automorphisms of $F$ inducing

\subitem  $\pm1$ on $H_3(F;\mathbb{Z})$), 220}

\item{ Johnson's trichotomy

\subitem (surface bundle groups), 92, 145}

\indexspace

\item{ $k(e,\eta)$ \quad (Montesinos knot), 
322}

\item{ $k_1 (M)$\quad (first $k$-invariant), 26} 

\item{ Kaplansky rank

\subitem ($\kappa(P)=dim_\mathbb{Q}\,\mathbb{Q}\otimes_\pi{P}$), 14}

\item{ $Kb$\quad (Klein bottle), 89}

\item{ $k_{p,q}$\quad ($(p,q)$-torus knot), 
312}

\item{ $kerv(\hat g)$\quad (codimension-2 Kervaire 

\subitem invariant of $\hat g$), 117}

\item{ Kervaire-Arf invariant, 117}

\item{ knot, 
269}
                   
\item{ knot group\quad ($\pi K$), 
270}

\item{ knot-like group, 
282}

\item{ knot manifold ($M(K)$), 
270}

\item{ Kodaira surface, 149}

\item{ $KS(M)$ (Kirby-Siebenmann invariant), 117}

\indexspace

\item{ $\ell P$\quad (locally $P$), 3}

\item{ $\ell^2(\pi)$ ($L^2$-completion of $\mathbb{C}[\pi]$), 23}

\item{ $L^2$-Betti number, 24, 26}              

\item{ $L^2$-homology, 24}              

\item{ lattice, 132}

\item{ linear $k$-invariant, 225}

\item{ link, 
284}

\item{ link group, 
284,288}

\item{ LHSSS\quad (Lyndon-Hochschild-Serre 

\subitem spectral sequence), 16}

\item{ locally $P$\quad   ($\ell P$), 3}

\item{ locally finite, 3}

\item{ L\"uck's Theorem  [L\"u94], 27}

\indexspace

\item{ $Mb$\quad (M\"obius band), 106}

\item{ $M(K)$\quad (closed manifold 

\subitem arising from a knot $K$), 
270}

\item{ $M(f)$\quad (mapping torus of a self

\subitem  homotopy equivalence $f$), 77}

\item{ $M(p,q,r)$ \quad (Brieskorn manifold), 
309}

\item{ mapping torus, 76, 186, 249} 

\item{ maximal finite normal subgroup 

\subitem (of a group with two ends), 16}

\item{ Mayer-Vietoris sequence 

\subitem of Waldhausen,  112}

\item{ Melvin's Theorem [Me84], 100}

\item{ meridian, 
270}

\item{ meridianal automorphism, 
279}

\item{ minimal complex surface, 148}

\item{ minimal model 

\subitem (for a $PD_4$-complex), 213}

\item{ minimal Seifert hypersurface, 
271}

\item{ monodromy, 
271}

\item{ morphism of Hilbert ${\mathcal N}(\pi)$-module, 24}

\item{ Mostow orbifold bundle, 139}
                       
\item{ Mostow rigidity, 192}

\indexspace

\item{ $n$-dimensional geometry, 132}

\item{ $\mathbb{N}il^3$
(nilpotent Lie geometry),  134}

\item{ $\mathbb{N}il^3\times\mathbb{E}^1$
(nilpotent Lie geometry), 

\subitem 135, 164}

\item{ $\mathbb{N}il^4$\quad (nilpotent Lie geometry),  

\subitem 135, 164}

\item{ ${\mathcal N}(\pi)$\quad  (von Neumann algebra),  23}

\item{ $n$-knot, 
269}

\item{ $N_G (H)$\quad (normalizer of a subgroup), 3}

\item{ nilradical, 135}

\item{ normal closure of $S$ in $G$ 

\subitem ($\langle\langle S\rangle\rangle_G$),  3}

\item{ nonsingular (of $\lambda_X$), 67}

\item{ Novikov ring, 73}

\indexspace

\item{ $Out(G)$\quad (group of outer 

\subitem automorphism classes), 3}

\item{ $O^*_1$\quad (binary octahedral group), 221}
                  
\item{ $O^*_k$\quad (extended binary octahedral 

\subitem group), 221}

\item{ orbifold bundle, 138}

\item{ordinary 2-knot, 361}

\item{ orientable $PD_n$-group

\subitem ($PD_n^+$-group), 21}

\item{ outer automorphism group, 3}

\indexspace

\item{ $P$ ($=PSL(2,\mathbb{R})$), 188}

\item{ $P_2 (M)$\quad (second stage 

\subitem of Postnikov tower), 26}

\item{ $PD_3 $-complex\quad (3-dimensional

\subitem Poincar\'e duality complex), 33}

\item{ $PD_3$-group, 37}

\item{ $PD_n $-complex, 33}
                                         
\item{ $PD_n^{(+)}$-group, 21}

\item{ $PD_n $-space, 33}

\item{ $PD_4$-polarization, 241}

\item{ $P^+$ (orientable cover), 33}

\item{ piece (of a geometric 

\subitem decomposition), 139}

\item{ Plotnick's Theorem [Pl86], 
314}

\item{ Poincar\'e duality, 32}

\item{ poly-, 4}

\item{ problem of the four exponentials, 137}

\item{ proper geometric decomposition, 139}

\item{ proper graph manifold, 114}

\indexspace
                       
\item{ $q(\pi)$, $q^{SG}(\pi)$ \quad (minimal 

\subitem Euler characteristic), 57}

\item{ $Q(2^na,b,c)$\quad (generalized 

\subitem quaternionic group), 223}
                                                                  
\item{ $Q(8k)$\quad (quaternionic group 

\subitem of order $8k$), 221}

\item{ quadratic 2-type 

\subitem ($[\pi,\pi_2 (M),k_1 (M),S(\widetilde M)]$), 241}
                           
\item{ quasifibre, 
329}

\item{ quaternion group\quad ($Q(8)$), 221}

\indexspace

\item{ rational surface, 149}

\item{ reduced intersection pairing  

\subitem ($\lambda_X$), 67}

\item{ reduced $L^2$-homology, 24}

\item{ reducible\quad ($\mathbb{H}^2\times\mathbb{H}^2$-manifold), 188}

\item{ reflexive knot, 270}

\item{ regular coherent ring, 15}

\item{ regular noetherian ring, 15}

\item{ restrained (group), 10}

\item{ ribbon knot, 275}

\item{ ruled surface, 149} 

\indexspace

\item{ $SG$\quad (class generated by groups 

\subitem of subexponential growth), 9}

\item{ $S^1$-actions, 263}

\item{ $\mathbb{S}^3$-group, 225}

\item{ $S_4^{PD}(P)$ \quad (polarized

\subitem  $PD_4$-complexes), 200}

\item{ $S^s_{TOP} (M)$\quad ($s$-cobordism 

\subitem structure set), 116}

\item{ $\mathbb{S}^4$\quad (spherical geometry), 234}

\item{ $\mathbb{S}^2\times\mathbb{S}^2$\quad (compact 

\subitem product geometry), 235}

\item{ $\mathbb{S}ol^4_{m,n}$, $\mathbb{S}ol^3\times\mathbb{E}^1$, 
\quad (solvable 

\subitem Lie geometries), 136, 164}

\item{ $\mathbb{S}ol^4_0$\quad (solvable Lie geometry), 

\subitem 137, 164}

\item{ $\mathbb{S}ol^4_1$\quad (solvable Lie geometry), 

\subitem 137, 165}

\item{ $\mathbb{S}^3\times\mathbb{E}^1$  (2-ended spherical-euclidean 

\subitem product geometry), 224}

\item{ $\mathbb{S}^2\times\mathbb{E}^2$ (1-ended spherical-euclidean 

\subitem product geometry), 203, 208}

\item{ $\mathbb{S}^2\times\mathbb{H}^2$\quad  (spherical-hyperbolic

\subitem product geometry), 203}

\item{ $\widetilde{\mathbb{SL}}\times\mathbb{E}^1$, 182}
 
\item{ safe extension, 24}
                          
\item{ satellite, 273}

\item{ $s$-concordant, 274}
                                            
\item{ Seifert fibred (4-manifold), 145}

\item{ Seifert hypersurface, 271}

\item{ semidirect product\quad ($G\!\rtimes_\theta\mathbb{Z}$), 4}

\item{ slice knot, 274}

\item{ solvable Lie type, 132, 176}

\item{ spin (Artin) of a knot\quad ($\sigma K$), 274}

\item{ split link, 285}

\item{square-root closed, 121}

\item{ $s$-rigid, 117}

\item{ stably homeomorphic, 122}

\item{ strict weight orbit, 278}
                 
\item{ Strebel's Theorem  [St77], 21}                                                                
\item{ strongly minimal, 213}

\item{ subnormal, 5}

\item{ sum of knots\quad ($K_1\sharp K_2$), 272}

\item{ superperfect, 290}

\item{ surface bundles, 89, 254}

\item{ surgery exact sequence, 116}

\item{ Swan complex, 219}

\item{ symplectic structure, 150, 265}

\indexspace
                                                                              
\item{ $T$\quad (torus), 89}        
                   
\item{ $T^*_1$\quad  (binary tetrahedral group), 221}
                   
\item{ $T^*_k$\quad  (extended binary tetrahedral 

\subitem group), 221}

\item{ $T(\pi)$ \quad (translation subgroup), 

\subitem 134, 136}

\item{ Tits alternative, 29, 39, 306}

\item{ translation subgroup ($T(\pi)$), 134, 136}

\item{ triangular (solvable Lie group), 136} 

\item{ trivial knot, 269}

\item{ trivial link, 285}

\item{ Turaev's Theorem  [Tu90], 34}

\item{ twist spin of a knot\quad ($\tau_r K$), 274}

\item{ type I, II, III  (Johnson's trichotomy 

\subitem for surface bundle groups), 92}

\item{ type R\quad (solvable Lie group), 136}

\indexspace

\item{ UCSS\quad (universal coefficient 

\subitem spectral sequence), 26}

\item{ Ue's Theorem, 146}

\item{ unreduced $L^2$-homology, 24}

\indexspace

\item{ $vP$\quad  (virtually $P$), 4}

\item{ virtual bundle group, 143}

\item{ virtually (qualifying a property 

\subitem of a group or space), 4}

\item{ von Neumann dimension of a 

\subitem Hilbert module ($dim_{{\mathcal N}(\pi)} M$), 24}

\indexspace

\item{ Waldhausen's Mayer-Vietoris 

\subitem sequence for $K$-theory, 112}

\item{ Weak Bass Conjecture, 14}

\item{ weak isomorphism, 24}

\item{ weak $PD_r$-group, 70}

\item{ weakly exact, 24}

\item{ weakly finite (ring), 15}
                                                                                
\item{ weight (class, element), 276}

\item{ weight orbit, 277}

\item{ Whitehead quadratic functor\quad 

\subitem ($\Gamma(-)$), 241}

\item{ Whitehead's question, 281}

\indexspace

\item{ $\mathcal X$\quad (class of groups), 31}

\item{ $X(K)$\quad (knot exterior), 269} 

\item{ $X_H $\quad (covering space with

\subitem fundamental group $H$), 25}

\item{ $\mathbb{X}$-manifold, 132}

\indexspace
     
\item{ $\mathbb{Z}^w$\quad ($w$-twisted integers), 13}
                                                                     
\item{ $Z*_m$\quad (group with presentation 

\subitem $\langle a,t\mid tat^{-1}=a^m\rangle$), 29}

\item{ $\mathbb{Z}\rtimes_{-1}\mathbb{Z}$\quad (fundamental group 

\subitem of Klein bottle, $\cong Z*_{-1}$), 29}

\bigskip
\centerline{{\sl Greek characters}}
\bigskip

\item{ $\alpha$-twisted endomorphism, 73}

\item{ $\beta_i(-)$\quad (Betti number), 25}

\item{ $\beta^{(2)}_i(-)$\quad ($L^2$-Betti number), 26}

\item{ $\beta^u$\quad ($u$-twisted Bockstein), 197}

\item{ $\Gamma(-)$\quad (Whitehead quadratic 

\subitem functor), 241}

\item{ $\Gamma_q$\quad (nilpotent group), 7}

\item{ $\Delta_a(X)=X^3-aX^2+(a-1)X-1$, 317}

\item{ $\zeta G$\quad (centre of a group), 3}

\item{ $\zeta_2G$\quad ($\zeta_2G/\zeta G=\zeta(G/\zeta G)$), 8}

\item{ $\eta_G$\quad (cohomology class), 70}

\item{ $\kappa(P)$\quad (Kaplansky rank), 14}

\item{ $\lambda_X$\quad (reduced intersection

\subitem pairing), 67}

\item{ $\Lambda=\mathbb{Z}[\mathbb{Z}]\cong\mathbb{Z}[t,t^{-1}]$\quad 

\subitem (Laurent polynomial ring), 6}
                    
\item{ $\pi K$\quad (knot group), 270}

\item{ $\pi_1$-slice, 275}

\item{ $\pi(e,\eta )$\quad (group of 2-twist spin 

\subitem of Montesinos knot), 323}

\item{$\pi^{orb}(B)$, orbifold fundamental 

\subitem group, 138}

\item{ $[\pi,m]_f$-complex, 32}

\item{ $\pi^+$ ($\mathrm{Ker}(w)$, 33}

\item{ $\sigma K$\quad  (Artin spin of $K$), 274}

\item{ $\tau$\quad (the twist of $S^2\times S^1$), 3}
                           
\item{ $\tau_r K$\quad ($r$-twist spin 

\subitem of a knot $K$), 274}

\item{ $\tau_{r,s} K$\quad (branched twist spin 

\subitem of a knot $K$), 313}

\item{ $\Phi$\quad ($\cong Z*_2$, 2-knot group), 297}  
                        
\item{ $\chi(\pi)$\quad (Euler characteristic 

\subitem of $vFP$ group $\pi$), 14}

\newpage
\centerline{{\sl Non-alphabetic symbols}}
\bigskip 

\item{ boundary $\partial:\pi_2 (B)\to\pi_1 (F)$\quad 

\subitem (connecting homomorphism), 89}

\item{ dagger $\dagger$: $L^\dagger=\overline{Hom_{\mathbb{Z}[\pi]}(L,{\mathbb{Z}[\pi]})}$ 
 
\subitem the conjugate dual module, 67}

\item{ double angle brackets $\langle\langle{\ }\rangle\rangle$:
$\langle\langle S\rangle\rangle_G$\quad 

\subitem (normal closure of $S$ in $G$), 3}
   
\item{ overbar $\bar{\ }$: anti-involution $\bar g=w(g)g^{-1}$, 

\subitem conjugate module $\overline M$, 13}  
 
\item{ plus superscript ${}^+$: orientation

\subitem preserving subgroup, cover, 33} 

\item{ prime $'$: commutator subgroup $G'$, 

\subitem maximal abelian cover $X'$, 3, 271}            

\item{ semidirect product: $G\rtimes_\theta\mathbb{Z}$, 4}
                                                                                 \item{ sharp $\sharp$: sum of knots $K_1\sharp K_2$, 272}
                      
\item{ star $^*$: $K^*$ (Gluck reconstruction

\subitem of a knot $K$), 270}

\item{ surd $\sqrt{\ }$: $\sqrt G$ (Hirsch-Plotkin 

\subitem radical of a group $G$), 7}

\item{ tilde $\widetilde{\ }$: $\widetilde X$ (universal cover), 25}

\end{theindex}